\documentclass[12pt]{book}

\usepackage{epigraph}
\usepackage{titling}
\usepackage{Preamble}

\newcommand{\subtitle}[1]{%
	\posttitle{%
		\par\end{center}
	\begin{center}\large#1\end{center}
	\vskip0.5em}%
}

\title{\Huge\textsf{Consequence  Relations}}
\subtitle{\Large An Introduction to the Tarski-Lindenbaum Method}

\author{\Large A. Citkin and A. Muravitsky}

\date{}

\addtocounter{chapter}{0}

\begin{document}

\maketitle

\setlength{\epigraphwidth}{.8\textwidth}
\thispagestyle{empty}

\vspace{5in}
 
 \epigraph{\large\emph{From one perspective, all intellectual and artistic pursuits are efforts to understand the world, including ourselves and our relation to the rest of the world. If successful, they would not `leave everything as it is’ but bring about some, however slight, desired change of the world in the form of modifications of our surroundings or our consciousness.}}{\large Hao Wang, Beyond Analytic Philosophy (1981)}

\frontmatter

\chapter*{Foreword}
When Tarski introduced the concept of a consequence operation in his 1930 paper ‘Fundamentale Begriffe der Methodologie der deduktiven Wissenschaften I’, readers could have been forgiven if they thought that it was just ‘abstract nonsense’ reflecting an insistence, characteristic of the Lvov-Warsaw school of mathematics, on always striving for the highest possible level of generality. In effect, Tarski took the already quite abstract notion of a closure operation that had been introduced into topology by Kuratowski, and noticed that if we omit two of its conditions (that the closure of the empty set is empty and that closure distributes over unions), then we have an important instantiation in logic: the set of all logical consequences of a set of propositions satisfies the remaining conditions for closure operations. Yet, despite this abstractness, in the following decades the notion became a versatile instrument in the logician’s tool-box. 

In the process, it underwent some subtle refinements, of which we mention two. Although in his initial intuitive explanations Tarski used the language of relations saying, for example, that “The set of all consequences of the set $A$ of sentences is here denoted by the symbol ‘$\textbf{Cn}(A)$’”, both the official definition and its formal development were carried out in the language of operations. It was not until the 1970s, following remarks of Scott, that it became standard practice to treat the two as trivially inter-definable and, in some accounts, to focus principally on the relational presentation. 

Again, in his 1930 paper Tarski quickly moved from consequence operations themselves to the ‘deductive or closed systems’ that they determine. Given a consequence operation \textbf{Cn}, a deductive system with respect to \textbf{Cn} is any set $A$ of propositions such that $A = \textbf{Cn}(A)$. Since $A$ is not in general closed under substitution then, so long as it is fixed, substitution plays little part in the investigation. However, if we direct attention to the consequence operations themselves – or even to those deductive systems that are of the kind $\textbf{Cn}(\varnothing)$, which is typical for the theorem-sets of formal logic --- then the picture changes. As {\L}o\'{s} and Suzko emphasized in 1958, in that context it becomes natural to consider the case that \textbf{Cn} itself is closed under substitution, i.e. that $\sigma(\textbf{Cn}(A))\subseteq\textbf{Cn}(\sigma(A)))$ for any suitably defined substitution function $\sigma$ on the underlying language. The term ‘structural’ was introduced for this condition on consequence operations or relations (not to be confused with a quite different sense of the same term that is current in proof-theory in the tradition of Gentzen).

By the 1960s the first books making free use of the notion of structural consequence relations began to appear. The present writer vividly remembers coming across Rasiowa {\&} Sikorski’s celebrated volume \emph{The Mathematics of Metamathematics} as a graduate student in 1963, the year of its publication, and the profound influence that it had on him. In the 1970s Rasiowa followed up with her \emph{Algebraic Approach to Non-Classical Logic}, which widened the selection of non-classical logics under consideration. In both volumes, consequence relations were used in association with another basic concept of algebra, also going back to Tarski and Lindenbaum as far as its application to logic is concerned – quotient structures. Under suitable conditions, satisfied in classical logic and many of its non-classical variants, we can define an equivalence relation on formulae of the language that is well-behaved with respect to all its logical connectives, thus forming a quotient structure of equivalence classes of formulae that is equipped with algebraic operations corresponding to the connectives; this is now known as a Lindenbaum (or Lindenbaum-Tarski) algebra for that logic. The construction allows one to apply to logic powerful techniques of universal algebra, and to lift Tarski’s consequence operations/relations to the quotient level. 

However, it seems fair to say that for Rasiowa and Sikorski, center stage was occupied by quotient structures and their algebraic manipulations, while consequence operations were rather auxiliary. In particular, the {\L}o\'{s}-Suzko condition of structurality was not considered at all. It was in 1988 that W\'{o}jcicki’s \emph{Theory of Logical Calculi} gave consequence a central place in the story, reflected in the subtitle Basic Theory of Consequence Operations, with particular attention accorded to those consequence operations that are generated by ‘logical matrices’, especially those matrices with a single ‘designated element’ (cf. the masterly review by Bull in \emph{Studia Logica} 1991 50:623-629). 

The present text of Citkin {\&} Muravitsky continues in the tradition of these authors, taking into account advances that have been obtained since the publication of Wójcicki’s study but hitherto available only in the journals. One such advance, due to Muravitsky in a paper of 2014, is the concept of a ‘unital’ consequence operation (or relation). This is a very broad class of structural consequence relations over the language of a propositional logic (or, indeed, first-order logic) that can be used to define Lindenbaum-Tarski quotient algebras. It includes several other such classes in the literature, notably those of ‘implicative’ and ‘Fregean’ consequence relations. The quotient structures can then be used to obtain results about the logic, particularly negative results to the effect that certain kinds of formulae are not theorems of certain logical systems, or are not consequences of certain sets of formulae. 

Advanced undergraduate students in mathematics, and graduate students in philosophy, with a solid course in classical logic as well as some exposure to a bit of model theory and basic ideas of computability, should be able to understand and benefit from this text. May it serve them as well as did Rasiowa and Sikorski in my own youth and W\'{o}jcicki in following years.\\

\vspace{.2in}
\begin{tabular}{lllllllllllllllllll}
	&&&&&&&&&&&&&&&&&\text{David Makinson}\\
	&&&&&&&&&&&&&&&&&\text{London, UK}\\
\end{tabular}

\tableofcontents

\mainmatter

\chapter[Introduction]{Introduction}
\label{chapter:introduction}

\section{Overview of the key concepts}\label{section:key-concepts}
This book is devoted to the study of the concept of a \emph{consequence relation} in formal languages. The idea of ​​an argument that allows one to draw conclusions from predetermined premises goes back to Aristotle. French historian J-P. Vernant writes about this turning point in the development of argumentation in ancient Greek culture,
\begin{quote}
	``Historically, rhetoric and sophistry, by analysing the forms of discourse as the means of winning the contest in the assembly and the tribunal, opened the way for Aristotle's inquiries, which in turn defined the rules of proof along with the technique of persuasion, and thus laid down a logic of the verifiably true, a matter of theoretical understanding, as opposed to the logic of the apparent or probable, which presided over the hazardous debates on practical questions.'' (Cf.~\cite{vernant1982}, chapter 4.) 
\end{quote}

Based on his observations of the public debate and pedagogical practice of the time, Aristotle called an argument dialectical if it was presented as an argument that begins with predetermined, possibly conditional, premises and leads to a definite conclusion.\footnote{W. Kneale and M. Kneale write, ``In its earliest sense the word `dialectic' is the name of the method of argument which is characteristic of metaphysics. It is derived from the verb \textgreek{dial'ege{\textsigma}{\texttheta}ai}, which means `discuss', and, as we have already seen, Aristotle thinks of a dialectic premiss as one chosen by a disputant in an argument.'' Cf.~\cite{kneales1962}, chapter 1, section 3.}
On the contrary, the conclusion of a demonstrative, or didactic, argument was made from self-evident judgments, regarded as axioms. W. Kneale and M. Kneale explain the difference between the two types of reasoning as follows.
\begin{quote}
	``In demonstration, we start from true premisses and arrive with necessity at a true conclusion: in other words, we have proof. In dialectical argument, on the other hand, the premisses are not known to be true, and there is no necessity that the conclusion be true. If there is an approach to truth through dialectic, it must be more indirect.'' (\cite{kneales1962}, chapter 1, section 1)
\end{quote}

The notion of a consequence relation is mainly associated with dialectical reasoning, in which demonstrative reasoning is considered as a special case.\\

Most languages that we will deal with, although not all of them, fall into the category of propositional;\footnote{In justification, we remind the reader of the judgment of {\L}ukasiewicz: ``The logic of propositions is the basis of all logical and mathematical systems.'' Cf.~\cite{lukasiewicz1934}.} however, sometimes we will use first-order language in our discussion. The concept of \emph{logical consequence} in formal languages was introduced by A. Tarski in the 1930s; cf.~\cite{tarski1930a,tarski1930b,tarski1936b}. Carnap's earlier attempts to define logical consequence were, as Tarski rightly notes, ``connected rather closely with the particular properties of the formalized language which was chosen as the subject of investigation.'' (Cf.~\cite{tarski1936}; quoted from~\cite{tarski1983}, p.414.) According to Carnap, a sentence $\alpha$ follows logically from a set $X$ of sentences if $X$ and the negation of $\alpha$ together generate a contradictory class.\footnote{Regarding this definition, Tarski writes, ``The decisive element of the above definition obviously is the concept `contradictory'. Carnap's definition of this concept is too complicated and special to be reproduced here without long and troublesome explanations.'' (quoted from~\cite{tarski1983}, p. 414.)} 

Tarski defines logical consequence in terms of \emph{consequence operators}, not in terms of consequence relations; however, as was explicitly stated by D. Scott~\cite{scott1974}, section 1, (see also~\cite{scott1971}, section II, and~\cite{shoesmith-smiley2008}, pp. 19--21), the two concepts can be used interchangeably. Thus, we have two notions for expressing the same conception of logical consequence. This prompts us to introduce a unifying term, \emph{abstract logic}, so that fixing an abstract logic $\mathcal{S}$, we can use both the consequence relation $\vdashS$ and the consequence operator $\textbf{Cn}_{\mathcal{S}}$ in the same context.

Tarski regarded the notion of consequence operator as primitive, governed by axioms. In his exposition, this also concerned a non-empty set, S, (of \emph{meaningful sentences}) of cardinality less than or equal to $\aleph_0$, which is fixed for consideration. From our point of view, Tarski’s axioms for the consequence operator still give you space to develop, to a large extent (see, e.g.,~\cite{tarski1930b}), the concept of logical consequence, but the absence of any structure in the set of meaningful sentences is not so good. That is why, from the very beginning, we distinguish atomic and compound
\emph{well-formed sentences}, or \emph{formulas} (as we prefer to call Tarski's ``meaningful sentences'') and impose structural characteristics on the set S of all formulas so that S is not merely an infinite set but becomes a \emph{formula algebra} $\FormAl$, where $\Lan$ is a \emph{formal language} that is fixed over consideration. This leads us directly to the \emph{Lindenbaum method}, one of the main lines of our exposition, to which Chapter~\ref{chapter:consequence} and all subsequent chapters are devoted.

One of the properties of the Tarski consequence operator is that the set of consequences can increase, but never decreases with an increase in the set of premises. Such consequence operators are called \emph{monotone}, in contrast to \emph{non-monotone} consequence operators. Accordingly, the abstract logic and the consequence relation both associated with a monotone consequence operator are also called monotone, for monotonicity is one of the conditions of the consequence relation. In this book, we deal only with \emph{monotone abstract logics}.

Another property, \emph{finitarity}, marked by Tarski makes an important distinction between \emph{finitary} and \emph{non-finitary} abstract logics. On the other hand, one important property, \emph{structurality}, was missing in Tarski's axiomatization of logical consequence. Only in 1958, the concept of a \emph{structural consequence operator} was introduced in~\cite{los-suszko1958}. Accordingly, the abstract logic and the consequence operator both associated with a structural consequence operator are also called \emph{structural}. In this book, we deal with structural and non-structural abstract logics.

Given an abstract logic $\mathcal{S}$ and a set $X\cup\lbrace\alpha\rbrace$ of formulas, the two mutually exclusive results are possible --- $X\vdashS\alpha$ or $X\not\vdashS\alpha$. Tarski's intended interpretation of the former is expressed in~\cite{tarski1930b} as follows.
\begin{quote}
	``Let $A$ be an arbitrary set of sentences of a particular discipline. With the help of certain operations, the so-called \emph{rules of inference}, new sentences are derived from the set $A$, called the \emph{consequences of the set} $A$.''
	(quoted from~\cite{tarski1983}, p. 63)
\end{quote}

This view of Tarski on how to define the relation of logical consequence, namely in a purely syntactic form, was probably obtained from his observation of how all deductive sciences could be constructed if their scientific languages were completely formalized.\footnote{We would like to draw the reader's attention to the title of~\cite{tarski1930b}.}
In contrast to this view, J. {\L}ukasiewicz and A. Tarski showed in~\cite{lukasiewicz-tarski1930} a different way of doing this.
Namely, they gave the definition of \emph{logical matrix} and formulated the notion of the \emph{system generated by a logical matrix}. Moreover, which is important for the purposes of this book, they first formulated the \emph{Lindenbaum theorem} in print, which laid the foundation for the Lindenbaum method. 

If $\textbf{Cn}_{\mathcal{S}}$  is the consequence operator of an abstract logic $\mathcal{S}$, then, by definition, $\textbf{Cn}_{\mathcal{S}}(\varnothing)$ consists of all $\mathcal{S}$-\emph{theorems}. A logical matrix, whose algebra is the formula algebra and whose \emph{logical filter} (that is the set of \emph{designated elements}) is $\textbf{Cn}_{\mathcal{S}}(\varnothing)$, is a \emph{Lindenbaum matrix}. The Lindenbaum theorem reads: Given an abstract logic $\mathcal{S}$, the Lindenbaum matrix is an \emph{adequate logical matrix} for the set of $\mathcal{S}$-theorems, assuming that the last set is closed under \emph{formula substitution}. Thus,  we see that Lindenbaum matrix is a universal tool for separating $\mathcal{S}$-theorems from those formulas which are not $\mathcal{S}$-theorems. The Lindenbaum method has evolved in the further development of this idea.

The first step in this direction is to allow any set $\textbf{Cn}_{\mathcal{S}}(X)=\set{\alpha}{X\vdashS\alpha}$, the $\mathcal{S}$-\emph{theory} generated by $X$, to be a logical filter (and it is here that we need the structurality of $\textbf{Cn}_{\mathcal{S}}$). The updated Lindenbaum matrix in this way is denoted by $\LinS[X]$. It turns out that $\LinS[X]$ is an adequate matrix for all and only those formulas $\alpha$ that $X\vdashS\alpha$. Thus, given an structural abstract logic $\mathcal{S}$, the set of all matrices $\LinS[X]$, known as a \emph{bundle}, gives one a means to solve the problem: either $X\vdashS\alpha$ or $X\not\vdashS\alpha$, for any set $X\cup\lbrace\alpha\rbrace$ of formulas.

The Lindenbaum method opens up two perspectives in the field of abstract logic. We will consider them separately.

First, an immediate consequence of Lindenbaum’s theorem is that any structural abstract logic is determined by a class of logical matrices, which provides another alternative for defining a logical consequence, and not just by inference rules. Such a logical consequence is called a \emph{matrix consequence}. A special cases of a matrix consequence is a \emph{single-matrix consequence} when an abstract logic is determined by a single matrix. The criterion for the latter case was discovered by J. {\L}o\'{s} and R. Suszko in terms of the concept of \emph{uniformity}; cf.~\cite{los-suszko1958}. However, as R. W\'{o}jcicki noticed in~\cite{wojcicki1970}, {\L}o\'{s} and Suszko's argument is valid only for finitary abstract logics. Further, W\'{o}jcicki managed to establish in~\cite{wojcicki1988} a criterion for a single-matrix consequence in terms of uniformity and the concept of \emph{couniformity}. 

What is striking in the formulations of both criteria is essential that the cardinality of the set of variables of the object language $\Lan$, $\Var_{\Lan}$, must be greater than or equal to $\aleph_{0}$. It seems unlikely that the validity of both theorems depends on the cardinality of the set of variables. Therefore, it would be interesting to drop the condition that $\card{\Var_{\Lan}}\ge\aleph_{0}$, or to show that without this condition the theorems are invalid.

This brings us back to the restriction imposed by Tarski on the cardinality of the set of formulas. This restriction can be satisfied even if the set of propositional variables is finite while keeping the set of formulas infinite.

The second perspective of Lindenbaum method is related to the second alternative in `$X\vdashS\alpha$ or $X\not\vdashS\alpha$'. Namely, while `$X\vdashS\alpha$' can be established using inference rules, other means will be required to establish `$X\not\vdashS\alpha$'. A. Kuznetsov proposed in~\cite{kuznetsov1979} a broad conception of \emph{separating tools} which also covers our problem. His approach he explained in the following words:
\begin{quote}
	``[...] let us consider those of them [objects], the use of which is based on a special relation $R$ between formulas and these objects, which is preserved under the rules of inference [postulated] in a given calculus [...] If such a relation $R$ is fixed and an object $\alpha$ stands in the relation $R$ to each of given formulas (axioms, hypotheses), but a formula $A$ does not, we say that $\alpha$ separates (module $R$) the formula $A$ from the given formulas. [...] Objects of this kind, with such a relation $R$ to formulas and so agreed with the rules of inference of the calculus under consideration, we call \textit{separating tools}.'' (our translation; comp.~\cite{mur2014a})
\end{quote}

Although in case of structural abstract logic $\mathcal{S}$, a Lindenbaum matrix $\LinS[X]$ is a separating tool to decide whether $X\vdashS\alpha$, however, the algebra of this matrix can be very complex; therefore, the direct application of $\LinS[X]$ can be difficult to implement, if at all possible. The way out of this difficulty may lie in narrowing the class of structural abstract logics in the hope of obtaining more manageable Lindenbaum matrices.

This perspective prompted A. Muravitsky to introduce in~\cite{mur2014a} a class of \emph{unital abstract logics}, that is, the structural abstract logics $\aLog$ with the property: for every set $X$ of formulas, the logical filter $\textbf{Cn}_{\mathcal{S}}(X)$ of the Lindenbaum matrix $\LinS[X]$ is a congruence class with respect to the congruence on the algebra of $\LinS[X]$, generated by the set $\textbf{Cn}_{\mathcal{S}}(X)\times\textbf{Cn}_{\mathcal{S}}(X)$. As shown in~\cite{mur2014a}, the class of unital logics properly contains the class of \emph{implicative logics} introduced by H. Rasiowa~\cite{ras74}, among which we find classical propositional logic, intuitionistic propositional logic, modal logic S4 and many others, and the class of \emph{Fregean logics} introduced in~\cite{czelakowski-pigozzi2004a}.

The quintessence of the Lindenbaum method is the concept of a \emph{Lindenbaum-Tarski algebra}. 
Although the Lindenbaum–Tarski algebra can be defined for any abstract logic $\mathcal{S}$ and for each $\mathcal{S}$-theory, as, for example, in~\cite{citkin-mur2013} and \cite{citkin-mur2016}, we apply this definition exclusively for unital logic, since for some specifications of unital logic, such as \emph{implicational unital logic}, the Lindenbaum-Tarski algebra exhibits good universal properties, at best being a \emph{free algebra} over some \emph{varieties} and \emph{quasi-varieties}. It is in this capacity that the Lindenbaum-Tarski algebra is useful for finding separating tools.

The discussion above applies to propositional languages. Next, we move on to formal languages that do not contain quantifiers but contain the equality sign `$\approx$', interpreted as an equality relation, that is, as the diagonal of the Cartesian square of the carrier of an algebra. Replacing the logical filter of the logical matrix with the diagonal of the Cartesian square of the carrier of the latter, we obtain an E-\emph{matrix}. Adapting the notion of a matrix consequence to this new type of logical matrix leads to the concept of an \emph{equational consequence}. If not all E-matrices are used, but only those that belong to a class $\E$ of E-matrices, we have an $\E$-\emph{consequence}, which is denoted by $\vdashE$. Similar to the notion of a matrix consequence, the notion of an E-consequence is based on the notion of \emph{valuation} of equalities in an E-matrix. The equational consequence is the formalization of, perhaps, one of the oldest forms of mathematical reasoning, especially in algebra.

With the introduction of the concepts of an E-matrix and equational consequence, the concept of a \emph{Lindenbaum-Tarski matrix for $\E$ relative to $X$}, symbolically $\LTE[X]$, plays the same role for the equational consequence as the Lindenbaum-Tarski algebra does for the matrix consequence. Yet, given a set $X\cup\lbrace\epsilon\rbrace$ of equalities with the condition that the variables occurring in the equality $\epsilon$ are contained in the set of the variables occurring in the equalities of $X$, a more effective tool to determine whether $X\vdashE\epsilon$ is a special E-matrix related to the set $X$, which is called a \emph{Mal'cev matrix} and denoted by $\ME[X]$. In the context of separating tools, the algebra of any E-matrix valuating the equalities of $X$ is a \emph{homomorphic image} of the algebra of $\ME[X]$; and if the class $\A[\E]$ of all algebras of the matrices of $\E$ is a quasi-variety, then the algebra of each $\ME[X]$ belongs to $\A[\E]$. 

We arrive at a different type of consequence relation associated with a class $\E$ of E-matrices if, instead of relying on the notion of valuation, we base the new definition on the notion of \emph{validation} of equalities in an E-matrix.
Given a nonempty class $\E$ of E-matrices, this new consequence relation is called \emph{equational} L-\emph{consequence} associated with $\E$, or $\EL$-\emph{consequence} for short. The abstract logic associated with $\EL$-consequence is denoted by $\ELog$. If the class $\A[\E]$ is a variety, the set $\textbf{Cn}_{\ELog}(\varnothing)$ is called the \emph{equational logic} (relative to $\E$). In particular, if $\hat{\E}$ is the class of all E-matrices of a certain signature, $\textbf{Cn}_{\hat{\E}_{\textsf{L}}}(\varnothing)$ is known as the \emph{Birkhoff equational logic}.\\

The treatment of consequence relations changes when the formal language of the subject matter changes. This is evident in the transition from sentential languages to quantifier-free languages with equality. Even more radical changes are required when we deal with first-order languages, which include generalized quantifiers. The consequence relation for such languages is called  $\Q$-\textit{consequence}. The most difficult task in formulating the notion of $\Q$-consequence is the definition of a new semantics. Thus, $\Q$-\textit{structures} and $\Q$-\textit{matrices} replace algebras and logical matrices for sentential languages, as well as algebras and E-matrices for languages ​​without quantifiers with equality. The most radical part of the new semantics is the presence of infinite operations that interpret the quantifiers. Moreover, these operations may be partial. Nevertheless, there is a similarity with the matrix consequence, and the notion of a \emph{Lindenbaum $\Q$-matrix} plays a role in the $\Q$-consequence relations analogous to one the notion of a Lindenbaum matrix plays in the matrix consequence relations.\\

Effective decidability (or computability) problems associated with matrix consequence, require new concepts specific to the field. In particular, we introduce the relation of \emph{formula indistinguishability} associated with a given algebra. Since each relation of formula indistinguishability is an equivalence on a formula algebra, it is important to be able to construct effectively a \emph{complete set of representatives} of all equivalence classes with respect to this relation. In addition, since computability concerns limited resources, a key player here is the concept of a \emph{restricted Lindenbaum matrix}.\\

Thus, the book is devoted to the study of the field of application of the method, which arose from the concept of the \emph{Lindenbaum matrix} by A. Lindenbaum and the \emph{Lindenbaum theorem}, within the framework of the concept of a \emph{consequence relation} by A. Tarski and in the context of the conception of \emph{separating tools} by A. Kuznetsov. The unifying term \emph{Tarski-Lindenbaum method} is intended to refer to the first two headings as the key topics of this study. Our implementation of the Tarski-Lindenbaum method aims to emphasize the role of the conception of separating tools.

\section{Overview of the contents}
\paragraph{Chapter 2: Preliminaries}contains concepts and facts from the areas of set theory, topology, universal algebra, model theory, and computability theory that will be used in this book. We expect the reader to refer to this chapter as needed.

\paragraph{Chapter 3: Formal Languages}introduces the type of formal language that is used in Chapters 2--6 and 10 of the book. This sentential schematic language is illustrated by a number of examples. Also in this chapter, we define the operation of uniform (or simultaneous) substitution and the semantics for this language by introducing the concept of a logical matrix. The latter concept is also illustrated by a number of examples.

\paragraph{Chapter 4: Logical Consequence}explains the central concept of this book --- the concept of consequence relation, as well as its counterpart, that of consequence operator. We prove the main properties of these concepts. In addition, we show several ways in which consequence relations can be defined and modified. We pay special attention to the determination of consequence relations using logical matrices and inference rules.The first way leads to the Lindenbaum method and the Lindenbaum theorem. The second method is very natural, and we consider several of its options. In particular, we discuss in some detail the definition of consequence relations defined by modus rules. We also discuss in this chapter the idea of completeness.

\paragraph{Chapter 5: Matrix Consequence}is devoted to matrix consequence. We especially focus on the single-matrix consequence. Namely, we prove the {\L}o\'{s}-Suszko-W\'{o}jcicki theorem, the W\'{o}jcicki theorem, and the Shoesmith-Smiley theorem; all these theorems are about when a structural abstract logic (a unifying term for the structural consequence relation and structural consequence operator) has an adequate logical matrix. We conclude this chapter with two discussions --- on finitary matrix consequence and on the conception of separating tools.

\paragraph{Chapter 6: Unital Abstract Logics}introduces the concept of unital abstract logic. We illustrate this concept with a number of well-known examples of logical consequence, as well as some classes of logical consequence, such as implicative logics and Fregean logics. The quintessence of this chapter is the concept of a Lindenbaum-Tarski algebra. Also in this chapter, we are developing a technique that allows us to consider the Lindenbaum-Tarski algebra as a separating tool, or at least as a starting point for creating such a means.
As illustrative examples, we give a detailed description of the Lindenbaum-Tarski algebras with one and two variables of the classical propositional calculus and the Lindenbaum-Tarski algebra with one variable of the intuitionistic propositional calculus, the so-called Rieger-Nishimura algebra. We also include some applications of Lindenbaum-Tarski algebras.

\paragraph{Chapter 7: Equational Consequence}develops a theory of the consequence relation for a schematic language of terms with equality. 
This consequence relation, called the equational consequence, is determined in a semantic way by means of E-matrices, as well as syntactically by inference rules, and their equivalence is established as a completeness theorem. We also define the concept of the Lindenbaum-Tarski matrix and the concept of the Mal'cev matrix for the equational consequence. Further, we prove an analog of the Mal'cev first and second theorems, as well as analogs of the Dick and Tietze theorems for the equational consequence. The equational consequence based on implication logic is illustrated by two examples: the class of Boolean algebras with equality and the class of Heyting algebras with equality.

\paragraph{Chapter 8: Equational L-Consequence} While the equational consequence of Chapter 7 is a direct application of the concept of matrix consequence to E-matrices, the concept of equational L-consequence is an extension of Birkhoff's ``equational logic''; namely, in contrast to the equational logic, we assume that the equational L-consequence can be defined for any nonempty set of E-matrices.  Also in this chapter, we show the connection between the equational consequence and the equational L-consequence and illustrate the latter concept with two examples: the class of Boolean algebras with equality and the class of Heyting algebras with equality. In addition, we define the Lindenbaum-Tarski matrices for the equational L-consequence and demonstrate the usefulness of this concept.

\paragraph{Chapter 9: $\Q$-Consequence} The topic of this chapter is the treatment of a consequence relation in the framework of predicate languages with generalized quantifiers, which is called in this book $\Q$-consequence. We discuss two patterns of $\Q$-consequence: one is defined by $\Q$-structures and the other, for specified languages, through deductive systems. Since the former is largely related to matrix consequence, we employ the Lindenbaum method for its characterization. Then, this method is used to consider models of the three kinds of $\Q$-consequence given deductively. Completeness results are given for three $\Q$-consequence relations, but only for a limited case,  which, however, includes the subcase of all formal theorems.

\paragraph{Chapter 10: Decidability} This chapter discusses the effective decidability of problems related to finite logical matrices and finite atlases. Among these questions is the problem of the triviality of the logic of a given finite matrix, the problem of weak equivalence of two given finite matrices, and also the question of whether two given finite atlases determine the same consequence relation.

\chapter{Preliminaries}

\section{Preliminaries from set theory}\label{section:perlimenaries-set-theory}
In this book, we adhere to a na\v{\i}ve version of the Zermelo-Fraenkel axiomatic system with the Axiom of Choice, \textit{ZFC}. In this context and essentially throughout the book, the understanding of the word `na\v{\i}ve' differs from its use in colloquial speech as `simple' and `inexperienced'.
Rather, the language of set theory, which includes membership $\in$, inclusion $\subseteq$, proper inclusion $\subset$ and some other relations, including functions, as well as Cartesian products based on them, is employed to use the results of a strict axiomatic theory, namely \textit{ZFC}\index{ZFC}. The words `class', `family', etc. will be used as synonyms of the word `set'. 

In addition to the notation above, we use the following concepts and their notation.
\begin{itemize}
	\item The set with no elements, the \textit{\textbf{empty set}}, is denoted by $\emptyset$.
	\item We use the set-roster notation to define a set whose elements can be explicitly presented; for example, a set whose elements are $a,b,c$, will be written as $\lbrace a,b,c\rbrace$.
	\item A set that can be defined by a predicate, say $P(x)$, applied to the elements of a given set $A$, will be written by the set-builder notation
	as $\set{x\in A}{P(x)}$; given $A$ and $P(X)$, the existence of that set in \textit{ZF} is guaranteed by the Axiom of Comprehension.
	\item Given a set $X$, the \textit{\textbf{power set}} of $X$ is defined as follows:
	\[
	\mathcal{P}(X):=\set{Y}{Y\subseteq X};
	\]
	the existence of $\mathcal{P}(X)$ is guaranteed in \textit{ZF} by the Axiom of Power Set. 
	\item Given sets $X$ and $Y$, `$X\Subset Y$' denotes that $X$ is a finite (perhaps empty) subset of $Y$. 
	\item Given sets $X$ and $Y$, $X\cap Y$, $X\cup Y$, $X\setminus Y$ stand for the \textit{\textbf{union}}, \textit{\textbf{intersection}} of $X$ and $Y$, and the \textit{\textbf{difference}} by subtracting $Y$ from $X$,  respectively; the first two operations are generalized to any family $\lbrace X_i\rbrace_{i\in I}$ and denoted by
	$\bigcap_{i\in I}\lbrace X_i\rbrace_{i\in I}$ and $\bigcup_{i\in I}\lbrace X_i\rbrace_{i\in I}$, respectively.
	\item Given entities $x$ and $y$ (in \textit{ZFC}, $x$ and $y$ are sets, but we can treat them as elements of any nature), $(x,y)$ denotes the set $\lbrace x,\lbrace x,y\rbrace\rbrace$ which is called an \textit{\textbf{ordered pair}}; the characteristic property of this notion is the following:
	\[
	(x_1,y_1)=(x_2,y_2)~\Longleftrightarrow~
	x_1=x_2~\text{and}~y_1=y_2.
	\]
	\item The \emph{\textbf{Cartesian}} (or \emph{\textbf{direct}}) \emph{\textbf{product of two sets}}\index{Cartesian product}\index{direct product} (not necessarily distinct) $X$ and $Y$ (in this order) is the set
	\[
	X\times Y:=\set{(x,y)}{x\in X~\text{and}~y\in Y}.
	\]
	\item $R$ is a \emph{\textbf{binary relation on}}\index{binary relation} a set $X$ if $R\subseteq X\times X$; if $R$ is coincident with the equality relation on $X$, it is denoted by $\Delta_X$ and called the \emph{\textbf{diagonal}} of $X\times X$, that is
	\[
	\Delta_X:=\set{(x,x)}{x\in X};
	\] 
	if a set $X$ is fixed over a context, we write simply $\Delta$;
	\textit{partially ordered} and \textit{linear relations} will be used explicitly, while \textit{well-ordered relations} on ordinal and cardinal numbers only implicitly.
	\item Given relations $R$ and $S$ on a set $X$, the \textit{\textbf{composition}} of $R$ and $S$, symbolically $R\circ S$, (in this order) is the following binary relation on $X$:
	\[
	R\circ S:=\set{(x,y)}{\text{there exists $z\in X$ such that $(x,z)\in R$ and $(z,y)\in S$}}.
	\]
	\item $f$ is a \emph{\textbf{function}} (or \textit{\textbf{map}})\index{function}  from a set $X$ \textit{\textbf{into}} (or \textbf{\textbf{to}}) a set $Y$, in symbols $f:X\longrightarrow Y$, if $f\subseteq X\times Y$ and for every $x\in X$, there exists a unique $y\in Y$ such that $(x,y)\in f$; the set of all functions from $X$ to $Y$ is denoted by $Y^{X}$; instead of `$(x,y)\in f$', we write `$f(x)=y$' or `$f:x\mapsto y$';
	a function $f:X\longrightarrow Y$ is \emph{\textbf{one-to-one}}, or \emph{\textbf{injective}}, if $f(x_1)=f(x_2)$ implies $x_1=x_2$; $f$ is \emph{\textbf{onto}}, or $\emph{\textbf{surjective}}$, if for any $y\in Y$, there is an $x\in X$ such that $f(x)=y$. Note that composition (as defined above) of two functions is a function.
	\item The \emph{\textbf{Cartesian product of a family}} $\lbrace X_i\rbrace_{i\in I}$ is denoted by $\prod_{i\in I}\lbrace X_i\rbrace$, or simply by $\prod_{i\in I} X_i$, and is defined as the set of all such functions $f$ from $I$ into $\bigcup_{i\in I}\lbrace X_i\rbrace$ such that $f(i)\in X_i$; as customary, we call the elements of $\prod_{i\in I} X_i$  \textit{\textbf{sequences}}; if $x\in\prod_{i\in I} X_i$, we write $x_i$ rather than $x(i)$.
	\item If $f\in Y^{X}$ and $A\subseteq X$, then we denote
	\[
	f(A):=\set{y\in Y}{\text{there is $x\in A$ suth that $f(x)=y$}};
	\]
	and if $B\subseteq Y$, we denote
	\[
	f^{-1}(B):=\set{x\in X}{f(x)\in B}.
	\]
	\item Given functions $f\in Y^{X}$ and $g\in Z^{Y}$, a composite function $g\circ f:X\longrightarrow Z$ (in this order) is defined as usual
	\[
	g\circ f(x):=g(f(x)).
	\]
\end{itemize}

\textit{Ordinals} are sets that are generated from the empty set $\emptyset$ by means of the operations of \textit{successor}, $X^{\prime}:=X\cup\lbrace X\rbrace$, and union. All ordinals can be well-ordered by $\in$; in relation to ordinals, this relation is denoted by $<$. The finite ordinals are:
$0:=\emptyset$, $1:=0\cup\lbrace 0\rbrace$, etc. 

The \textit{Axiom of Infinity}\index{axiom!of infinity} states that the collection consisting of $0,1,\ldots$ is a set. It is denoted by $\mathbb{N}$; that is,
\[
\mathbb{N}:=\lbrace 0,1,\ldots\rbrace.
\]

Arranging the elements of $\mathbb{N}$ by $<$, we obtain a well-ordered set of type $\omega$, which is the least infinite ordinal.

\begin{center}
	\textbf{The Axiom of Choice}
\end{center}
All the concepts mentioned above can be defined in \textit{ZF}. However, they will not be enough to develop the reasoning that we need in this book. We need the Axiom of Choice.

The \emph{Axiom of Choice}\index{axiom!of choice} is the following statement.
\[
\begin{array}{l}
	\textit{For any nonempty collection $X$ of sets, there is a function $f$}\\
	\textit{from $X$ to $\bigcup X$ such that for any nonempty $Y\in X$, $f(Y)\in Y$}.
\end{array}\tag{AC}
\]

The function $f$ in the above statement is called a \emph{choice function} for $X$.

The following statements proved to be equivalent to (AC) in \textit{ZF}.
\begin{itemize}
	\item[(a)] Every set can be well-ordered.
	\item[(b)] The Cartesian product of a nonempty set of nonempty sets is nonempty.
	\item[(c)] In any nonempty partially ordered set, if any chain of this set has the greatest element, then this set contains a maximal element.  (\textit{Zorn's Lemma})
	\item[($\text{c}^{\ast}$)] In any nonempty partially ordered set, if any chain of this set has the greatest element, then any element of this set is less than or equal to a maximal element. (\textit{Zorn's Lemma})
\end{itemize}

Sets $X$ and $Y$ are said to have the \textit{same cardinality} if there is an injective function from $X$ onto $Y$; and the \textit{cardinality of $X$ is less} than the cardinality of $Y$ if there is an injective function from $X$ into $Y$, but not vice versa.

Employing (AC), we define the cardinal of a set as follows.
Given a set $X$, its \textit{\textbf{cardinal}}, symbolically $\card{X}$, is the least ordinal $\kappa$ such that $\kappa$ and $X$ have the same cardinality. 
That the cardinal of a set always exists follows from
the property (a) above, as well as from the facts that every well-ordered set is isomorphic to some ordinal and that all ordinals are well-ordered. The first property is equivalent to (AC), but the last two properties are proven in \textit{ZF}. Thus, two sets, say $X$ and $Y$, have the same cardinality if, and only if, $\card{X}=\card{Y}$. 

We denote:
\[
\aleph_0:=\card{\mathbb{N}}.
\]

If there is an injective function from $X$ into $Y$ but $\card{X}\neq\card{Y}$, we write `$\card{X}<\card{Y}$'.
As customary, we use the notation:
\[
\card{X}\le\card{Y}~\define~\card{X}<\card{Y}~\text{or}~\card{X}=\card{Y}.
\]

A set $X$ is said to be \textit{\textbf{infinite}} if there is no injective function from $X$ onto one of the elements of $\mathbb{N}$. If there is an injection from $X$ onto a proper subset of $X$, then $X$ is called  \textit{\textbf{Dedekind infinite}}. It can be proven in \textit{ZF} that a set is Dedekind infinite if, and only if, it contains a subset of the cardinality of $\mathbb{N}$.  Therefore, any Dedekind infinite set is infinite.
However, it can only be proven in \textit{ZFC} that a set is Dedekind infinite if, and only if, it is infinite. Thus one can prove in \textit{ZFC} that a set $X$ is infinite (or equivalently Dedekind infinite) if, and only if, $\card{X}\ge\aleph_0$.

\paragraph{References}
\begin{enumerate}
	\item~\cite{fraenkel-bar-hillel-levy1973}
	\item~\cite{halmos1974}
	\item~\cite{johnstone1987}
	\item~\cite{kuratowski-mostowski1976}
\end{enumerate}

\section{Preliminaries from topology}\label{section:topology}
In this book, we use some well-known facts from general topology.\\

A \textit{\textbf{topological space}} \index{topological space}(or simply \textit{\textbf{space}}) is a pair $(X,\Top)$ (perhaps both $X$ and $\Top$ with subscripts), where $\Top\subseteq\mathcal{P}(X)$ and such that 
\[
\begin{array}{cl}
	1^{\circ} &\varnothing\in\Top~\text{and}~X\in\Top;\\\\
	2^{\circ} &\text{for any $U_1,\ldots,U_n\in\Top$, $U_1\cap\ldots\cap U_n\in\Top$};\\\\
	3^{\circ} &\text{for any $\lbrace U_i\rbrace_{i\in I}\subseteq\Top$, $\bigcup\lbrace U_i\rbrace_{i\in I}\in\Top$}.
\end{array}
\]

In a topological space $(X,\Top)$, $X$ is called the carrier of the space and $\Top$ is the topology of the space. The elements of $\Top$ are called the \textit{\textbf{open sets}}\index{topological space!open sets} of the space. The complement of an open set is called \textit{\textbf{closed}}\index{topological space!closed sets}. The topology $\Top$ is called \textit{\textbf{discrete}} if $\Top=\mathcal{P}(X)$. In our applications, we will be assuming that $X\neq\varnothing$.
If $X$ is a finite set, any space $(X,\Top)$ is called \textit{\textbf{finite}.}

It is customary to apply the term \textit{space} (or \textit{topological space}) to $X$, meaning that $X$ is endowed with a family $\Top$ of subsets of $X$ satisfying the conditions $1^{\circ}$--$3^{\circ}$. It must be clear that $X$ can be endowed with different families of open sets. Slightly abusing notation, we write $X=(X,\Top)$.

A family $\lbrace U_i\rbrace_{i\in I}\subseteq\Top$ is called an \textit{\textbf{open cover}}\index{topological space!open cover} of a space $X=(X,\Top)$ if $\bigcup\lbrace U_i\rbrace_{i\in I}=X$. Any
subfamily $\lbrace U_i\rbrace_{i\in I_0}$ with $I_0\subseteq I$ is called a \textit{\textbf{subcover}} of the cover $\lbrace U_i\rbrace_{i\in I}$ if $\bigcup\lbrace U_i\rbrace_{i\in I_0}=X$. A subcover is called \textit{\textbf{finite}} if $I_0$ is a nonempty finite set.

A topological space is called \textit{\textbf{compact}}\index{topological space!compact} if each open cover has a finite subcover.

The following observation is obvious.
\begin{prop}\label{P:finite-discrete=compact}
	Any finite discrete topological space is compact.
\end{prop}

Let $X$ be a nonempty set. Suppose we have a family $\text{B}_{0}:=\lbrace Y_i\rbrace_{i\in I}\subseteq\mathcal{P}(X)$. Suppose that we want to define a topology $\Top$ on $X$ such that sets $Y_i$ are open in $X$ endowed with this topology. 

First we form all finite intersections $Y_{i_1}\cap\ldots\cap Y_{i_n}$ collecting them in a set $\text{B}$. Then, we define all unions $\bigcup_{j\in J}\set{Z_j}{Z_{j}\in\text{B}}$ collecting them in a set $\Top$, to which we add $\varnothing$ and $X$ (if they are not yet there). The sets of the resulting set $\Top$ we announce open and define the space $X:=(X,\Top)$. It is not difficult to check that the conditions $1^{\circ}$--$3^{\circ}$ are satisfied. The family $\text{B}_{0}$ is called a \textit{\textbf{subbase}} of the space $X$, and the set $\text{B}$ is a \textit{\textbf{base} of $X$}. It is obvious that all sets of the subbase $\text{B}_{0}$ are open in $X$. \\

Now we employ this way of introducing topology on a set to define the cartesian product of topological spaces\index{topological spaces!Cartesian product}. Assume that we have a collection $\lbrace X_i\rbrace_{i\in I}$ of topological spaces $X_{i}:=(X_{i},\Top_i)$, where $i\in I$. First, we define the cartesian product of the sets $X_i$; that is
\[
\text{P}:=\prod_{i\in I}X_i.
\]
An arbitrary element \textbf{x} of \text{P} will be denoted by $(x_i)_{i\in I}$; that is $\textbf{x}:=(x_i)_{i\in I}$. Given an $i\in I$, a \textit{\textbf{projection}} $\rho_{i}:\text{P}\longrightarrow X_i$ is defined by the mapping $\rho_{i}:\textbf{x}\mapsto x_i$. Now we define:
\[
\text{B}_{0}:=\set{\rho_{i}^{-1}(U)}{U\in\Top_i~\text{and}~i\in I}.
\]

Finally, $\text{B}_{0}$ is employed as a subbase to define a \textit{\textbf{product topology}}\index{topological spaces!product topology} (aka \textit{\textbf{Tychonoff topology)}}\index{topological spaces!Tychonoff topology} on \text{P}. Denoting the product topology on \text{P} by $\Top$, we call the topological space $(\text{P},\Top)$ a \textit{\textbf{cartesian product}} of the collection $\lbrace(X_i,\Top_i)\rbrace_{i\in I}$.

From the definition of the product topology, we see that its base \text{B} constituted by the sets of the form
\[
\rho_{i_1}^{-1}(U_{i_1})\cap\ldots\cap\rho_{i_k}^{-1}(U_{i_k}),
\]
where $U_{i_1}\in\Top_{i_1},\ldots,U_{i_k}\in\Top_{i_k}$. This means that the sets of this base can be represented as follows:
\begin{equation}\label{E:representation}
	\prod_{i\in I}Z_i,
\end{equation}
where for some $I_0\Subset I$, if $i\in I_0$, then $Z_i$ is an arbitrary set from $\Top_i$, and if $i\notin I_0$, then $Z_i=X_i$. In particular, if each $X_i$ in the cartesian product \text{P} is a finite discrete space, then in each representation~\eqref{E:representation}, there is a set $I_0\Subset I$ such that if $i\in I_0$, then $Z_i$ is an arbitrary subset of $X_i$, and if $i\notin I_0$, then $Z_i=X_i$.

In our application of cartesian products (Section~\ref{section:finitary-matrix-consequence}) the following proposition is essential.

\begin{prop}[Tychonoff theorem]\label{P:tychonoff}\index{Theorem!Tychonoff Theorem}
	The cartesian product of a collection of compact topological spaces is compact relative to the product topology.	
\end{prop}

\begin{cor}\label{C:product-finite-discrete-spaces}
	The cartesian product of a collection of discrete spaces is compact relative to the product topology.	
\end{cor}
\paragraph{References}
\begin{enumerate}
	\item~\cite{bourbaki1998}
	\item~\cite{kelley1975}
\end{enumerate}

\section{Preliminaries from algebra}\label{section:prelimenaries-algebra}

\subsection{Subalgebras, homomorphisms, direct products and subdirect products}
Let $\alg{A}=\langle\textsf{A};\Func,\Cons\rangle$ be an algebra of type $\Lan$, where \textsf{A}, called the \textit{\textbf{carrier}}, or the \textit{\textbf{universe}}, of \alg{A}, is a nonempty set, and $\Func$ and $\Cons$ are operations on \textsf{A} of arity $n>0$ and arity $n=0$, respectively. The operations of $\Func\cup\Cons$ are called the \textit{\textbf{fundamental operations}}\index{algebra!fundamental operation} of \alg{A}. We also denote the carrier of an algebra $\alg{A}$ by $|\alg{A}|$. We denote by $\Delta_{\textsf{A}}$ (or by $\Delta_{\alg{A}}$) the diagonal of $|\alg{A}|\times |\alg{A}|$; if an algebra \alg{A} is fixed over consideration, we write simply $\Delta$.

Assume that a nonempty set $\textsf{B}\subseteq\textsf{A}$ is such that for any $F\in\Func$  of arity $n$ and any $b_1,\ldots,b_n\in\textsf{B}$, $Fb_1\ldots b_n\in\textsf{B}$, and $\Cons\subseteq\textsf{B}$. Then the algebra $\alg{B}=\langle\textsf{B};\Func,\Cons\rangle$ is called a \textit{\textbf{subalgebra}} of the algebra \alg{A}. 

Given an algebra, the intersection of any set of subalgebras of this algebra, providing that the intersection of the carriers is nonempty, is a subalgebra of the algebra. 

Let \alg{A} be an algebra and $\varnothing\neq X\subseteq|\alg{A}|$. It is clear that the set of the subalgebras of \alg{A} that include $X$ is nonempty; for instance, \alg{A} itself is such a subalgebra (of itself).
The least (in the sense of $\subseteq$) subalgebra of \alg{A} that contains $X$ exists and is denoted by $\lbrack X\rbrack_{\alg{A}}$. If the set $\Cons\neq\varnothing$, then condition $X\neq\varnothing$ can be dropped, for then the algebra $\lbrack\varnothing \rbrack_{\alg{A}}$ will be the least subalgebra of $\alg{A}$ that contains $\Cons$. The algebra $\lbrack X\rbrack_{\alg{A}}$ is called the \textit{\textbf{subalgebra}} (of \alg{A}) \textit{\textbf{generated by the set}} $X$, where the elements of the latter set are called \textit{\textbf{generators}}.\index{algebra!generator}
\begin{prop}[\cite{cohn1981}, chapter II, proposition 5.1, \cite{gratzer2008}, {\S}9, lemma 3]\label{P:subalgebra-generated-by-X}
	Let {\em$\alg{A}=\langle\textsf{A};\Func,\Cons\rangle$} be an algebra of type $\Lan$. Also, let {\em$\varnothing\neq X\subseteq|\alg{A}|$}. Then
	the carrier of {\em$\lbrack X\rbrack_{\alg{A}}$} consists of the elements
	$\alpha[a_1,\ldots,a_n]$, where $\alpha$ is an $\Lan$-formula with $n$ variables {\em(}see Chapter~\ref{chapter:languages}{\em)} and $a_1,\ldots,a_n\in X$.
\end{prop}

Let \alg{A} and \alg{B} be algebras of type $\Lan$. A map $f:|\alg{A}|\longrightarrow|\alg{B}|$ is called a \textit{\textbf{homomorphism}}\index{algebra!homomorphism} if the following permutability conditions (for each $F\in\Func$ of arity $n$ and any $a_1,\ldots,a_n\in|\alg{A}|$) are satisfied:
\begin{equation}\label{E:homomorphism-1}
	f(Fa_1\ldots a_n)=F f(a_1)\ldots f(a_n),
\end{equation}
and for each $c\in\Cons$,
\begin{equation}\label{E:homomorphism-2}
	f(c)=c,
\end{equation}
where the left-hand occurrence of `$c$' in the above equality is an element of $|\alg{A}|$, and the right-hand occurrence of `$c$' is an element of $|\alg{B}|$.

From an algebraic point of view, it is convenient to regard constants as 0-ary operations, which we will do in what follows.

An `onto' homomorphism is called an \emph{\textbf{epimorphism}}. A one-one homomorphism is called an \emph{\textbf{embedding}}.\\

Let \alg{A} be an algebra of type $\Lan$. An equivalence on $|\alg{A}|$ is called a \textit{\textbf{congruence}}\index{algebra!congruence} if the following preservation conditions
(for each $F\in\Func_{\mathcal{L}}$ of arity $n$ and any $a_1,\ldots,a_n\in|\alg{A}|$) are satisfied:
\[
(a_1,b_1),\ldots,(a_n,b_n)\in\theta~\Longrightarrow~
(Fa_1\ldots a_n,Fb_1\ldots b_n)\in\theta.
\]

Given a congruence $\theta$ on an algebra $\alg{A}$ and $X\subseteq|\alg{A}|$, we denote:
\[
[X]\theta:=\set{x\in|\alg{A}|}
{(x,y)\in\theta,~\text{for some $y\in X$}}
\]
and call $[X]\theta$ the \textit{\textbf{congruence class on}}\index{algebra!congruence class} $\alg{A}$ \textit{\textbf{modulo a congruence}} $\theta$, \textit{\textbf{generated by the class}} $X\times X$. $X$ is called a \textit{\textbf{generating set}}.
If $X=\lbrace x\rbrace$, we write $[x]\theta$ (or $x/\theta$) instead of $[\lbrace x\rbrace]\theta$.
$\alg{A}\slash\theta(X)$ is a \textit{\textbf{quotient algebra relative to the congruence}}\index{algebra!quotient algebra} $\theta(X)$.\\

Let $\theta$ be a congruence on an algebra \alg{A}. It is well known that the map $x\mapsto x\slash\theta$ is a homomorphism. This homomorphism is called \textit{\textbf{natural}}; we denote it by $\mathfrak{n}_{\theta}$.\\

Let $f:X\longrightarrow Y$ be an arbitrary map of $X$ to $Y$. A relation
\[
\text{ker}(f):=\set{(x,y)\in X\times X}{f(x)=f(y)}
\]
is called the \textit{\textbf{kernel}} of $f$.
\begin{prop}
	Let {\em$f:\alg{A}\longrightarrow\alg{B}$} be a homomorphism. Then {\em$\text{ker}(f)$} is a congruence on {\em\alg{A}}.
\end{prop}

Let $\alg{A}$ be an algebra of type $\Lan$ and $\theta$ be a congruence on $\alg{A}$. Fix a set $X\subseteq|\alg{A}|$. We denote a special congruence on $\alg{A}$ as follows:
\[
\theta(X):=\bigcap\set{\theta\in\Congruence\,\alg{A}}{X\times X\subseteq\theta}.
\]

We note that
\begin{equation}\label{E:preliminaries-algebra-1}
	X\subseteq Y\subseteq|\alg{A}|~\Longrightarrow~\theta(Y)\subseteq\theta(X)
\end{equation}

\begin{prop}\label{P:homomorphic-image-is-subalgebra}
	The homomorphic image of an algebra {\em$\alg{A}$} under a homomorphism into an algebra {\em$\alg{B}$} is a subalgebra of {\em$\alg{B}$}.
\end{prop}

The following proposition is called the ``first isomorphism theorem'' in~\cite{cohn1981}, theorem II.3.7, and the ``Homomorphism Theorem'' in~\cite{burris-sankapp1981}, theorem 6.12, and in~\cite{gratzer2008}, {\S}11, theorem 1.\index{Theorem!Homomorphism Theorem}\index{algebra!Homomorhism Theorem}
\begin{prop}\label{P:isomorphism-first}
	Let {\em\alg{A}} and {\em\alg{B}} be algebras of type $\Lan$ and {\em$f:\alg{A}\longrightarrow\alg{B}$} be a homomorphism. Then the map {\em$g:a\slash\text{ker}(f)\mapsto f(a)$}, where {\em$a\in\alg{A}$}, is an isomorphism of {\em$\alg{A}\slash\text{ker}(f)$} onto the image of {\em\alg{A}} under $f$ such that {\em$f=g\circ \mathfrak{n}_{\text{ker}(f)}$}.
	{\em(See Figure~\ref{figure:isomorphism-theorem-first}.)}
\end{prop}
\begin{figure}
	\[
	\ctdiagram{
		\def\diagramunit{0.8pt}
		\ctinnermid
		\ctv 0,100:{\alg{A}}
		\ctv 100,100:{\alg{B}}
		\ctet 0,100,100,100:{f}
		\ctv 50,30:{\alg{A}\slash\text{ker}(f)}
		\ctel 0,100,50,30:{\mathfrak{n}_{\text{ker}(f)}}
		\cter 50,30,100,100:{g}
	}
	\]	
	\caption{Illustration of Proposition~\ref{P:isomorphism-first}}\label{figure:isomorphism-theorem-first}
\end{figure}

The following proposition is the first part of theorem II.3.11 (third isomorphism theorem) \index{Theorem!Third Isomorphism Theorem}\index{algebra!Third Isomorphism Theorem}in~\cite{cohn1981}.
\begin{prop}\label{P:congruence-thm-auxiliary}
	Let $\theta$ and $\phi$ be congruences on an algebra {\em$\alg{A}$} with $\theta\subseteq\phi$. Then the map
	$a\slash\theta\mapsto a\slash\phi$, where {\em$a\in|\alg{A}|$}, defines an epimorphism of {\em$\alg{A}\slash\theta$} onto {\em$\alg{A}\slash\phi$}.
\end{prop}

Let $\theta$ and $\phi$ be congruences on an algebra \alg{A} with $\theta\subseteq\phi$. We define:
\begin{equation}\label{E:definition-of-fraction-congruence}
	\phi\slash\theta:=\set{(a\slash\theta,b\slash\theta)}{(a,b)\in\phi}.
\end{equation}

It is well known that $\phi\slash\theta$ is a congruence on $\alg{A}\slash\theta$; cf.~\cite{burris-sankapp1981}, lemma 6.14.

The next proposition is called ``Second Isomorphism Theorem''\index{Theorem!Second Isomorphism Theorem}\index{algebra!Second Isomorphism Theorem} in~\cite{burris-sankapp1981}, theorem 6.15.
\begin{prop}\label{P:isomorphism-second}
	Let $\theta$ and $\phi$ be congruences on an algebra {\em\alg{A}} with $\theta\subseteq\phi$. Then the map $(a\slash\theta)\slash(\phi\slash\theta)\mapsto a\slash\phi$ is an isomorphism
	of {\em$(\alg{A}\slash\theta)(\phi\slash\theta)$} onto {\em$\alg{A}\slash\phi$}. \\
	{\em(See Figure~\ref{figure:isomorphism-second}.)}
\end{prop}
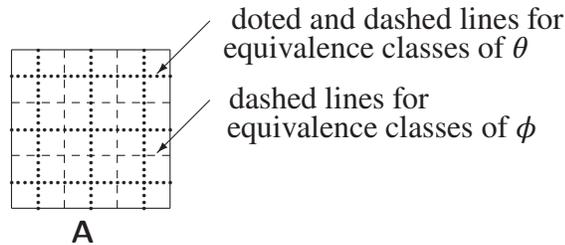
\begin{figure}[h!]
	\paragraph{Algebra}
	\[
	\ctdiagram{
		\ctet 75,71,55,51:{}
		\ctet 75,41,55,21:{}
		\ctnohead
		\cten 0,0,0,60:{}
		\cten 0,0,60,0:{}
		\cten 60,0,60,60:{}
		\cten 0,60,60,60:{}
		\ctnohead\ctdash
		\def\zzctdrawdashedge{\drawdashedge{3pt}{2.5pt}11}	
		\cten 20,0,20,60:{}
		\cten 40,0,40,60:{}
		\cten 0,20,60,20:{}
		\cten 0,40,60,40:{}
		\ctnohead\ctdot
		\def\zzctdrawdotedge{\drawdotedge{2.5pt}1}
		\cten 0,10,60,10:{}
		\cten 0,30,60,30:{}
		\cten 0,50,60,50:{}
		\cten 10,0,10,60:{}
		\cten 30,0,30,60:{}
		\cten 50,0,50,60:{}
		\ctv 146,71:{\text{doted and dashed lines for}}
		\ctv 138,60:{\text{equivalence classes of }\theta}
		\ctv 120,41:{\text{dashed lines for}}
		\ctv 140,30:{\text{equivalence classes of }\phi}
		\ctv 27,-10:{\alg{A}}
	}
	\]
	\paragraph{Quotients}
	\[		
	\begin{array}{ccc}
		\begin{array}{rcl}
			\ctdiagram{
				\ctnohead\ctdash
				\def\zzctdrawdashedge{\drawdashedge{3pt}{2.5pt}11}	
				\cten 0,0,0,14:{}
				\cten 0,0,14,0:{}
				\cten 14,0,14,14:{}
				\cten 0,14,14,14:{}
			} &
			\ctdiagram{
				\ctnohead\ctdash
				\def\zzctdrawdashedge{\drawdashedge{3pt}{2.5pt}11}	
				\cten 0,0,0,14:{}
				\cten 0,0,14,0:{}
				\cten 14,0,14,14:{}
				\cten 0,14,14,14:{}
			} &
			\ctdiagram{
				\ctnohead\ctdash
				\def\zzctdrawdashedge{\drawdashedge{3pt}{2.5pt}11}	
				\cten 0,0,0,14:{}
				\cten 0,0,14,0:{}
				\cten 14,0,14,14:{}
				\cten 0,14,14,14:{}
			} \\
			\ctdiagram{
				\ctnohead\ctdash
				\def\zzctdrawdashedge{\drawdashedge{3pt}{2.5pt}11}	
				\cten 0,0,0,14:{}
				\cten 0,0,14,0:{}
				\cten 14,0,14,14:{}
				\cten 0,14,14,14:{}
			} &
			\ctdiagram{
				\ctnohead\ctdash
				\def\zzctdrawdashedge{\drawdashedge{3pt}{2.5pt}11}	
				\cten 0,0,0,14:{}
				\cten 0,0,14,0:{}
				\cten 14,0,14,14:{}
				\cten 0,14,14,14:{}
			} &
			\ctdiagram{
				\ctnohead\ctdash
				\def\zzctdrawdashedge{\drawdashedge{3pt}{2.5pt}11}	
				\cten 0,0,0,14:{}
				\cten 0,0,14,0:{}
				\cten 14,0,14,14:{}
				\cten 0,14,14,14:{}
			} \\
			\ctdiagram{
				\ctnohead\ctdash
				\def\zzctdrawdashedge{\drawdashedge{3pt}{2.5pt}11}	
				\cten 0,0,0,14:{}
				\cten 0,0,14,0:{}
				\cten 14,0,14,14:{}
				\cten 0,14,14,14:{}
			} &
			\ctdiagram{
				\ctnohead\ctdash
				\def\zzctdrawdashedge{\drawdashedge{3pt}{2.5pt}11}	
				\cten 0,0,0,14:{}
				\cten 0,0,14,0:{}
				\cten 14,0,14,14:{}
				\cten 0,14,14,14:{}
			} &
			\ctdiagram{
				\ctnohead\ctdash
				\def\zzctdrawdashedge{\drawdashedge{3pt}{2.5pt}11}	
				\cten 0,0,0,14:{}
				\cten 0,0,14,0:{}
				\cten 14,0,14,14:{}
				\cten 0,14,14,14:{}
			} 
		\end{array}
		&&
		\begin{array}{rcl}
			\ctdiagram{
				\ctnohead
				\cten 0,0,0,4:{}
				\cten 0,0,4,0:{}
				\cten 4,0,4,4:{}
				\cten 0,4,4,4:{}
				\cten 6,0,6,4:{}
				\cten 6,0,10,0:{}
				\cten 10,0,10,4:{}
				\cten 6,4,10,4:{}
				\cten 0,6,0,10:{}
				\cten 0,6,4,6:{}
				\cten 4,6,4,10:{}
				\cten 0,10,4,10:{}
				\cten 6,6,6,10:{}
				\cten 6,6,10,6:{}
				\cten 10,6,10,10:{}
				\cten 6,10,10,10:{}
				\ctnohead\ctdash
				\def\zzctdrawdashedge{\drawdashedge{3pt}{2.5pt}11}
				\cten -2,-2,12,-2:{}
				\cten -2,-2,-2,12:{}
				\cten 12,-2,12,12:{}
				\cten -2,12,12,12:{}
			} 
			&
			\ctdiagram{
				\ctnohead
				\cten 0,0,0,4:{}
				\cten 0,0,4,0:{}
				\cten 4,0,4,4:{}
				\cten 0,4,4,4:{}
				\cten 6,0,6,4:{}
				\cten 6,0,10,0:{}
				\cten 10,0,10,4:{}
				\cten 6,4,10,4:{}
				\cten 0,6,0,10:{}
				\cten 0,6,4,6:{}
				\cten 4,6,4,10:{}
				\cten 0,10,4,10:{}
				\cten 6,6,6,10:{}
				\cten 6,6,10,6:{}
				\cten 10,6,10,10:{}
				\cten 6,10,10,10:{}
				\ctnohead\ctdash
				\def\zzctdrawdashedge{\drawdashedge{3pt}{2.5pt}11}
				\cten -2,-2,12,-2:{}
				\cten -2,-2,-2,12:{}
				\cten 12,-2,12,12:{}
				\cten -2,12,12,12:{}
			} 
			&	
			\ctdiagram{
				\ctnohead
				\cten 0,0,0,4:{}
				\cten 0,0,4,0:{}
				\cten 4,0,4,4:{}
				\cten 0,4,4,4:{}
				\cten 6,0,6,4:{}
				\cten 6,0,10,0:{}
				\cten 10,0,10,4:{}
				\cten 6,4,10,4:{}
				\cten 0,6,0,10:{}
				\cten 0,6,4,6:{}
				\cten 4,6,4,10:{}
				\cten 0,10,4,10:{}
				\cten 6,6,6,10:{}
				\cten 6,6,10,6:{}
				\cten 10,6,10,10:{}
				\cten 6,10,10,10:{}
				\ctnohead\ctdash
				\def\zzctdrawdashedge{\drawdashedge{3pt}{2.5pt}11}
				\cten -2,-2,12,-2:{}
				\cten -2,-2,-2,12:{}
				\cten 12,-2,12,12:{}
				\cten -2,12,12,12:{}
				\ctv 0,-3:{}
			} 
			\\
			\ctdiagram{
				\ctnohead
				\cten 0,0,0,4:{}
				\cten 0,0,4,0:{}
				\cten 4,0,4,4:{}
				\cten 0,4,4,4:{}
				\cten 6,0,6,4:{}
				\cten 6,0,10,0:{}
				\cten 10,0,10,4:{}
				\cten 6,4,10,4:{}
				\cten 0,6,0,10:{}
				\cten 0,6,4,6:{}
				\cten 4,6,4,10:{}
				\cten 0,10,4,10:{}
				\cten 6,6,6,10:{}
				\cten 6,6,10,6:{}
				\cten 10,6,10,10:{}
				\cten 6,10,10,10:{}
				\ctnohead\ctdash
				\def\zzctdrawdashedge{\drawdashedge{3pt}{2.5pt}11}
				\cten -2,-2,12,-2:{}
				\cten -2,-2,-2,12:{}
				\cten 12,-2,12,12:{}
				\cten -2,12,12,12:{}
			} 
			&
			\ctdiagram{
				\ctnohead
				\cten 0,0,0,4:{}
				\cten 0,0,4,0:{}
				\cten 4,0,4,4:{}
				\cten 0,4,4,4:{}
				\cten 6,0,6,4:{}
				\cten 6,0,10,0:{}
				\cten 10,0,10,4:{}
				\cten 6,4,10,4:{}
				\cten 0,6,0,10:{}
				\cten 0,6,4,6:{}
				\cten 4,6,4,10:{}
				\cten 0,10,4,10:{}
				\cten 6,6,6,10:{}
				\cten 6,6,10,6:{}
				\cten 10,6,10,10:{}
				\cten 6,10,10,10:{}
				\ctnohead\ctdash
				\def\zzctdrawdashedge{\drawdashedge{3pt}{2.5pt}11}
				\cten -2,-2,12,-2:{}
				\cten -2,-2,-2,12:{}
				\cten 12,-2,12,12:{}
				\cten -2,12,12,12:{}
			} 
			&	\ctdiagram{
				\ctnohead
				\cten 0,0,0,4:{}
				\cten 0,0,4,0:{}
				\cten 4,0,4,4:{}
				\cten 0,4,4,4:{}
				\cten 6,0,6,4:{}
				\cten 6,0,10,0:{}
				\cten 10,0,10,4:{}
				\cten 6,4,10,4:{}
				\cten 0,6,0,10:{}
				\cten 0,6,4,6:{}
				\cten 4,6,4,10:{}
				\cten 0,10,4,10:{}
				\cten 6,6,6,10:{}
				\cten 6,6,10,6:{}
				\cten 10,6,10,10:{}
				\cten 6,10,10,10:{}
				\ctnohead\ctdash
				\def\zzctdrawdashedge{\drawdashedge{3pt}{2.5pt}11}
				\cten -2,-2,12,-2:{}
				\cten -2,-2,-2,12:{}
				\cten 12,-2,12,12:{}
				\cten -2,12,12,12:{}
				\ctv 0,-3:{}
			} 
			\\
			\ctdiagram{
				\ctnohead
				\cten 0,0,0,4:{}
				\cten 0,0,4,0:{}
				\cten 4,0,4,4:{}
				\cten 0,4,4,4:{}
				\cten 6,0,6,4:{}
				\cten 6,0,10,0:{}
				\cten 10,0,10,4:{}
				\cten 6,4,10,4:{}
				\cten 0,6,0,10:{}
				\cten 0,6,4,6:{}
				\cten 4,6,4,10:{}
				\cten 0,10,4,10:{}
				\cten 6,6,6,10:{}
				\cten 6,6,10,6:{}
				\cten 10,6,10,10:{}
				\cten 6,10,10,10:{}
				\ctnohead\ctdash
				\def\zzctdrawdashedge{\drawdashedge{3pt}{2.5pt}11}
				\cten -2,-2,12,-2:{}
				\cten -2,-2,-2,12:{}
				\cten 12,-2,12,12:{}
				\cten -2,12,12,12:{}
			} 
			&
			\ctdiagram{
				\ctnohead
				\cten 0,0,0,4:{}
				\cten 0,0,4,0:{}
				\cten 4,0,4,4:{}
				\cten 0,4,4,4:{}
				\cten 6,0,6,4:{}
				\cten 6,0,10,0:{}
				\cten 10,0,10,4:{}
				\cten 6,4,10,4:{}
				\cten 0,6,0,10:{}
				\cten 0,6,4,6:{}
				\cten 4,6,4,10:{}
				\cten 0,10,4,10:{}
				\cten 6,6,6,10:{}
				\cten 6,6,10,6:{}
				\cten 10,6,10,10:{}
				\cten 6,10,10,10:{}
				\ctnohead\ctdash
				\def\zzctdrawdashedge{\drawdashedge{3pt}{2.5pt}11}
				\cten -2,-2,12,-2:{}
				\cten -2,-2,-2,12:{}
				\cten 12,-2,12,12:{}
				\cten -2,12,12,12:{}
			} 
			&		\ctdiagram{
				\ctnohead
				\cten 0,0,0,4:{}
				\cten 0,0,4,0:{}
				\cten 4,0,4,4:{}
				\cten 0,4,4,4:{}
				\cten 6,0,6,4:{}
				\cten 6,0,10,0:{}
				\cten 10,0,10,4:{}
				\cten 6,4,10,4:{}
				\cten 0,6,0,10:{}
				\cten 0,6,4,6:{}
				\cten 4,6,4,10:{}
				\cten 0,10,4,10:{}
				\cten 6,6,6,10:{}
				\cten 6,6,10,6:{}
				\cten 10,6,10,10:{}
				\cten 6,10,10,10:{}
				\ctnohead\ctdash
				\def\zzctdrawdashedge{\drawdashedge{3pt}{2.5pt}11}
				\cten -2,-2,12,-2:{}
				\cten -2,-2,-2,12:{}
				\cten 12,-2,12,12:{}
				\cten -2,12,12,12:{}
			} 
			\\
		\end{array}\\
		\alg{A}/\phi&&(\alg{A}/\theta)/(\phi/\theta)
	\end{array}
	\]	
	\caption{Illustration of Proposition~\ref{P:isomorphism-second}}\label{figure:isomorphism-second}	
\end{figure}
The following proposition could be regarded as a reverse of Proposition~\ref{P:isomorphism-second}.
\begin{prop}\label{P:congruence-extension}
	Let $\theta$ be a congruence on {\em$\alg{A}$}. Then for any congruence $\phi$ on {\em$\alg{A}\slash\theta$}, there is a congruence $\eta$ on {\em$\alg{A}$} such that $\theta\subseteq\eta$ and the map $(a\slash\theta)\slash\phi\mapsto a\slash\eta$ is an isomorphism of {\em$(\alg{A}\slash\theta)\slash\phi$} onto {\em$\alg{A}\slash\eta$}.
\end{prop}
\begin{proof}
	First of all, we note that the map $(a\slash\theta)\slash\phi\mapsto a\slash\eta$, where $a\in|\alg{A}|$, is defined correctly. For, for any congruence $\eta$ on \alg{A}  with $\theta\subseteq\eta$, if $a\slash\theta=b\slash\theta$, for some $a,b\in|\alg{A}|$, then $a\slash\eta=b\slash\eta$.
	
	Now we aim to define $\eta$. It can be done as follows.
	\[
	(a,b)\in\eta~\stackrel{\text{df}}{\Longleftrightarrow}~(a\slash\theta,b\slash\theta)\in\phi.
	\]
	
	It is easily seen that $\eta$ is a congruence on \alg{A} and $\theta\subseteq\eta$. Also, it is not difficult to see that the map
	$(a\slash\theta)\slash\phi\mapsto a\slash\eta$ is an isomorphism.
\end{proof}

Let $\lbrace\alg{A}_i\rbrace_{i\in I}$ be a family of algebras of type $\Lan$. We denote by
\[
\alg{A}:=\prod_{i\in I}\alg{A}_i
\]
the \textit{\textbf{direct product}}\index{algebra!direct product} of the algebras $\alg{A}_i$, which is an algebra of type $\Lan$, whose carrier $|\alg{A}|$ is the Cartesian product of  $\lbrace|\alg{A}_i|\rbrace_{i\in I}$ and each $n$-ary operation $f$ is defined pointwise: for any $x_1,\ldots,x_n\in|\alg{A}|$,
\[
f(x_1,\ldots,x_n):=(f((x_1)_{i},\ldots,(x_n)_i))_{i\in I},
\]
and each $\Lan$-constant $c:=(c_i)_{i\in I}$.

If all $\alg{A}_i$ are isomorphic to one another, we call $\alg{A}$ a \textit{\textbf{direct power}}\index{algebra!direct power}.
\begin{prop}[comp.~\cite{burris-sankapp1981}, chapter II,theorem 7.9]\label{P:direct-product-property}
	Given two sets {\em$\{\alg{A}_i\}_{i\in I}$} and {\em$\{\alg{A}_j\}_{j\in J}$} of algebras of type $\Lan$ with $I\cap J=\varnothing$, the algebras {\em$\prod_{i\in I}\alg{A}_{i}\times\prod_{j\in J}\alg{A}_{j}$} and {\em$\prod_{k\in I\cup J}\alg{A}_{k}$} are isomorphic with respect to the map
	$f(a,b)=(c_k)_{k\in I\cup J}$, where $c_k=a_k$, if $k\in I$, and $c_k=b_k$, if $k\in J$.
\end{prop}

A subalgebra $\alg{B}$ of the direct product $\alg{A}=\prod_{i\in I}\alg{A}_i$ is called a \textit{\textbf{subdirect product}}\index{algebra!subdirect product} of the family $\lbrace\alg{A}_i\rbrace_{i\in I}$ if each projection $p_i:x\mapsto x_i$ is an epimorphism of $\alg{B}$ onto $\alg{A}_i$. If an isomorphic image of $\alg{B}$ is a subalgebra of $\prod_{i\in I}\alg{A}_i$, we speak of a \textit{\textbf{subdirect embedding}}.

An algebra $\alg{A}$ is \textit{\textbf{subdirectly irreducible}}\index{algebra!subdirectly irreducible} if \alg{A} is trivial or $\Congruence\,\alg{A}\setminus\{\Delta\}$ has a least congruence. It is equivalent to say that \alg{A} is subdirectly irreducible if, and only if, for any subdirect embedding $f:\alg{A}\longrightarrow\prod_{i\in I}\alg{A}_i$, there is $i\in I$ such that $p_{i}\circ f:\alg{A}\longrightarrow\alg{A}_i$ is an isomorphism. Next we define a map $h: x\mapsto(f(x),g(x))$, for any $x\in|\alg{B}|$. It is clear that $h$ is a homomorphism from $\alg{B}$ into $\prod_{i\in I}\alg{A}_{i}\times
\alg{A}$.

\begin{prop}[Birkhoff]\label{P:Birkhoff's-theorem}\index{Theorem!Birkhoff Theorem}\index{algebra!Birkhoff Theorem}
	Every algebra of type $\Lan$ is isomorphic to a subdirect product of subdirectly irreducible algebras of type $\Lan$.
\end{prop}

We denote by $\mathcal{S}_{\mathcal{L}}$ the abstract class\footnote{The term \textit{abstract class} (of algebras) is defined in the next subsection.} of all subdirectly irreducible algebras of type $\Lan$. 

\paragraph{References}
\begin{enumerate}
	\item~\cite{burris-sankapp1981}
	\item~\cite{cohn1981}
	\item~\cite{gratzer2008} 
\end{enumerate}

\subsection{Class operators}

Given a class $\mathcal{K}$ of algebras of type $\Lan$, we employ the following class operators:\index{algebra!class operators}
\begin{itemize}
	\item $\alg{A}\in\Is\mathcal{K}$ if, and only if, $\alg{A}$ is an isomorphic image of some algebra of $\mathcal{K}$;\index{$\Is\mathcal{K}$}
	\item $\alg{A}\in\Su\mathcal{K}$ if, and only if, $\alg{A}$ is isomorphic to a subalgebra of some algebra of $\mathcal{K}$;\index{$\Su\mathcal{K}$}
	\item $\alg{A}\in\Ho\mathcal{K}$ if, and only if, $\alg{A}$ is a homomorphic image of some algebra of $\mathcal{K}$;\index{$Ho\mathcal{K}$}
	\item $\alg{A}\in\Pro\mathcal{K}$ if, and only if, $\alg{A}$ is isomorphic to the direct product of a nonempty set of algebras in $\mathcal{K}$; \index{$\Pro\mathcal{K}$}
	\item $\alg{A}\in\Pros\mathcal{K}$ if, and only if, $\alg{A}$ is a subdirect product of a nonempty set of algebras in $\mathcal{K}$;\index{$\Pros\mathcal{K}$}
	\item $\alg{A}\in\Pred\mathcal{K}$ if, and only if, $\alg{A}$ is isomorphic to a reduced product\footnote{The notion of reduced product is defined in Section~\ref{section:preliminaries-model-theory}.} of a nonempty set of algebras in $\mathcal{K}$;\index{$\Pred\mathcal{K}$}
	\item $\alg{A}\in\Pu\mathcal{K}$ if, and only if, $\alg{A}$ is isomorphic to an ultraproduct\footnote{The notion of ultraproduct is defined in Section~\ref{section:preliminaries-model-theory}.} of a nonempty set of algebras in $\mathcal{K}$.\index{$Pu\mathcal{K}$}
\end{itemize}

A class $\mathcal{K}$ is called \emph{abstract} if, along with each algebra in it, $\mathcal{K}$ contains its isomorphic images; in other words, $\mathbf{I}\mathcal{K}=\mathcal{K}$. We note that, applying any of the above operators to $\mathcal{K}$, we obtain an abstract class.

In the sequel, we will find useful another class operator. 
Proposition~\ref{P:Birkhoff's-theorem} induces the following definition.

Given a class $\mathcal{K}$ of algebras of type $\Lan$, 
\[
\begin{array}{rl}
	\alg{A}\in\Si\mathcal{K} \define &\text{there is $\alg{B}\in\mathcal{K}$ and there is $\mathcal{K}_0\subseteq\mathcal{S}_{\mathcal{L}}$ such that $\alg{A}\in\mathcal{K}_0$}\\ 
	&\text{and $\alg{B}$ is a subdirect product of all algebras of $\mathcal{K}_0$}.
\end{array}\index{$\Si\mathcal{K}$}
\] 

In terms of operators $\Pros$ and $\Si$, Proposition~\ref{P:Birkhoff's-theorem} can be reformulated as follows.
\begin{prop}\label{P:K-included-P_s-K_si}
	Given a class $\mathcal{K}$ of algebras of type $\Lan$, {\em$\mathcal{K}\subseteq\Pros\Si\mathcal{K}$}.
\end{prop}
\begin{proof}
	Let $\alg{A}\in\mathcal{K}$. In virtue of Proposition~\ref{P:Birkhoff's-theorem}, there is a set $\lbrace\alg{A}_i\rbrace_{i\in I}\subseteq\mathcal{S}_{\Lan}$ such that
	$\alg{A}$ is a subdirect product of $\lbrace\alg{A}_i\rbrace_{i\in I}$. By definition, $\lbrace\alg{A}_i\rbrace_{i\in I}\subseteq\Si\mathcal{K}$.	
\end{proof}

For any class operators $O_1$ and $O_2$, we denote
\[
O_1\le O_2\stackrel{\text{df}}{\Longleftrightarrow}\text{for any class $\mathcal{K}$, $O_{1}\mathcal{K}\subseteq O_{2}\mathcal{K}$}.
\]

\begin{prop}[\cite{burris-sankapp1981}, chapter II, lemma 9.2 ]\label{P:class-operators-inequalities}
	The following inequalities hold: $\Su\Ho\le\Ho\Su$, $\Pro\Su\le\Su\Pro$,
	$\Pro\Ho\le\Ho\Pro$, and $\Pros\le\Su\Pro$. Also, $\Ho\Ho=\Ho$, $\Su\Su=\Su$, and $\Pro\Pro=\Pro$.
\end{prop}

There is a simple characterization of the class operator $\Si$. In fact, it collects all subdirectly irreducible homomorphic images of any algebra of a given class $\mathcal{K}$.
\begin{prop}\label{P:characterisation-of-Si-operator}
	Let $\mathcal{K}$ be a class of algebras of type $\Lan$. Then {\em$\Si\mathcal{K}=(\Ho\mathcal{K})\cap\mathcal{S}_{\Lan}$}.
\end{prop}
\begin{proof}
	Suppose $\alg{A}\in\Si\mathcal{K}$. Then, by the definition of the $\Si$ operator and by the definition of subdirect products, the algebra $\alg{A}$ is subdirectly irreducible and is a homomorphic image of an algebra from $\mathcal{K}$.
	
	Now assume that $\alg{A}\in(\Ho\mathcal{K})\cap\mathcal{S}_{\Lan}$. This implies that there is an algebra $\alg{B}$ such that $\alg{A}$ is a homomorphic image with respect to an epimorphism $g:\alg{B}\longrightarrow\alg{A}$.
	Further, in virtue of Proposition~\ref{P:Birkhoff's-criterion-for varieties}, there is a set $\{\alg{A}\}_{i\in I}$ of subdirectly irreducible algebras and a subdirect embedding $f:\alg{B}\longrightarrow\prod_{i\in I}\alg{A}_i$.
	In virtue of Proposition~\ref{P:direct-product-property}, the map
	\[
	h:\prod_{i\in I}\alg{A}_{i}\times\alg{A}\rightarrow\prod\set{\alg{A}^{\prime}}{\alg{A}^{\prime}\in\{\alg{A}_i\}_{i\in I}\cup\{\alg{A}\}}
	\] 
	where $h(a,b)=(c)_{k\in I\cup\{j\}}$ with $j\notin I$, $c_{k}=c_{k}$, if $k\in I$, and
	$c_{k}=b$, if $k=j$, is an isomorphism. 
	
	Next we define a map:
	\[
	f^{\ast}: x\mapsto(f(x),g(x)),
	\]
	for any $x\in|\alg{B}|$. It is clearly seen that the map $h\circ f^{\ast}$ is a subdirect embedding of $\alg{B}$ into $\prod\set{\alg{A}^{\prime}}{\alg{A}^{\prime}\in\{\alg{A}_i\}_{i\in I}\cup\{\alg{A}\}}$. 
\end{proof}

\paragraph{References}
\begin{enumerate}
	\item~\cite{burris-sankapp1981}
	\item~\cite{gratzer2008} 	
\end{enumerate}

\subsection{Term algebra, varieties and quasi-varieties}\label{section:term-algebra}
Below we give algebraic characterizations of \emph{equational classes} and of  \emph{implicational classes} of algebras of the same type. These characterizations require formulas of first-order. We defer the definition of a first-order formula to Section~\ref{section:preliminaries-model-theory}. Here we use the first-order formulas of two special kinds. 

If we want to represent the signature of the algebras of type $\Lan$ symbolically, we introduce for each operation from $\Func$ and for each constant from $\Cons$ the symbols, usually denoted by the same letters, corresponding to those operations and constants. Thus we generate the set of function symbols $\Func_{\mathcal{L}}$ and the set of constant symbols $\Cons_\Lan$. In addition, we introduce a nonempty set $\Var_{\Lan}$ of (individual) variables.
Using these categories of symbols, we generate the set of \textit{\textbf{terms}}, by including in the latter all variables and constants, as well as the strings that are obtained according to the following rule: if $t_1,\ldots,t_n$ are terms and $F$ is a functional symbol of arity $n$, then $Ft_1\ldots t_n$ is a term. The last procedure can be regarded as an operation on the set of terms with the constant symbols as 0-ary operations. Thus we obtain an algebra of type $\Lan$, which is called a \textit{\textbf{term algebra}}\index{algebra!term algebra} and denoted by $\FormAl$.\index{$\FormAl$} The specificity of the term algebra of type $\Lan$ is that, given any algebra $\alg{A}$ of type $\Lan$ and an arbitrary map $f:\Var_{\Lan}\longrightarrow|\alg{A}|$, $f$ can be extended to a homomorphism
$\hat{f}:\FormAl\longrightarrow\alg{A}$. It is customarily to denote the initial map $f$ and its extension $\hat{f}$ by the letter $f$. The property just described is called a \textit{\textbf{universal mapping property}}\index{universal mapping property} of algebra $\FormAl$ \textit{\textbf{relative to}} the class of all algebras of type $\Lan$. Below we will see that some algebras have the universal property relative to special subclasses of the class of all algebras of type $\Lan$.

In this section, we describe first-order languages with equality and without additional predicate symbols. More information about first-order languages and their semantics can be found in the references at the end of Section~\ref{section:preliminaries-model-theory}. We use a unifying notation `$\FOL$'\index{$\FOL$} to denote first-order languages which are discussed in this subsection and in Section~\ref{section:preliminaries-model-theory}.
In Chapter~\ref{chapter:Q-consequence}, we will deal with first-order languages of a different grammar, called there $\Q$-languages.

Any $\FOL$ consists of a set $\Var$ of individual variables, logical constants `$\Land$', `$\Lor$, `$\Rarrow$', `$\Neg$', and quantifiers `$\Forall$' and `$\Exists$', as well as parentheses `$($' and 
`$)$'. 

The formal language $\Lan$ of terms is an integral part of $\FOL$. In this subsection, we use a simplified version of $\FOL$, which is not reflected in the notion of a term but is reflected in the notion of a formula.

If $t$ and $s$ are $\Lan$-terms, then $t\approx s$, where $\approx$ is regarded as the only (binary) predicate symbol (of this version of $\FOL$), is a first-order formula; we also call any such formula an \emph{\textbf{equality}}. If $A$ and $B$ are first-order formulas, then:
\begin{itemize}\label{FOL-forfula}
	\item[(1)] $(A\Land B)$ is a first-order formula;
	\item[(2)] $(A\Lor B)$ is a first-order formula;
	\item[(3)] $(A\Rarrow B)$ is a first-order formula;
	\item[(4)]  $\Neg A$ is a first-order formula;
	\item[(5)]  $\Forall x A$, where $x\in\VarL$, is a first-order formula;
	\item[(6)]  $\Exists x A$, where $x\in\VarL$, is a first-order formula.
\end{itemize}
Only equalities and the formulas that can be built according to (1)--(6) are \emph{first-order formulas} (or $\FOL$-\emph{formulas} or simply \emph{formulas}) of (this simplified version of) $\FOL$.\index{first-order formula}

A formula of the form (5) or (6) can be a constituent of a larger formula. Whether a formula of the form (5) or (6) is considered individually of as a constituent of a larger formula, $A$ is called the scope of the quantifier $\Forall$ or of the quantifier $\Exists$, respectively; and the variable $x$ following the quantifier, as well any occurrence of $x$ in $A$ is called \textit{bound}.\index{bound variable}
The occurrences of variables that are not in the scope of any quantifier with the equal variable are called \textit{free}\index{bound variable}. If an $\FOL$-formula is closed, it is called an $\FOL$-\textit{sentence}. A formula which does not contain quantifiers is called \textit{\textbf{quantifier-free}}. 

Let $A$ be a formula, all free variables of which constitute a set $\{x_1,\ldots,x_n\}$. The formula $\forall x_1\ldots\forall\ x_n A$ is called a \emph{universal closure}\index{universal closure} of $A$. It should be clear that if the formula $A$ has more than one free variable, it has more than one universal closures. However, from a semantic viewpoint, all these closures have equal truth values in all interpretations.\footnote{This is explained in Section~\ref{section:preliminaries-model-theory}. } Thus, at least semantically, they are indistinguishable. Therefore, we use the term \emph{the universal closure} of a formula $A$, which may have several instances. We denote the universal closure of $A$ by $\Forall\ldots\Forall A$. In particular, the universal closure $\Forall\ldots\Forall\,\alpha\approx\beta$ of an equality $\alpha\approx\beta$ is called an \textit{\textbf{identity}}.\index{identity} A formula of the form
\[
(\alpha_1\approx\beta_1\Land\ldots\Land\alpha_n\approx\beta_n)\Rarrow\alpha\approx\beta, \tag{$\ast$}
\]
where $n\ge 0$, is called a \textit{\textbf{quasi-equality}},\index{quasi-identity} and its universal closure
\[
\Forall\ldots\Forall ((\alpha_1\approx\beta_1\Land\ldots\Land\alpha_n\approx\beta_n)\Rarrow\alpha\approx\beta) 
\]
a \textit{\textbf{quasi-identity}}. In particular, each equality is a quasi-equality, and each identity is a quasi-identity, with $n=0$ in both cases.\\

Given an $\FOL$-equality $\alpha\approx\beta$ and an algebra $\alg{A}$ of type $\Lan$, we say that $\alpha\approx\beta$ is \textit{\textbf{valid}} in $\alg{A}$, symbolically
$\alg{A}\models\alpha\approx\beta$, if for any homomorphism $v:\FormAl\longrightarrow\alg{A}$, $v[\alpha]=v[\beta]$.

Similarly, a quasi-equality, say $(\ast)$, is \textit{\textbf{valid}} in $\alg{A}$, symbolically\index{valid quasi-identity}
\[
\alg{A}\models ((\alpha_1\approx\beta_1\Land\ldots\Land\alpha_n\approx\beta_n)\Rarrow\alpha\approx\beta),
\] if for any $v$,  $v[\alpha]=v[\beta]$, whenever all the equalities $v[\alpha_i] =v[\beta_i]$, where $1\le i\le n$, hold.
\begin{rem}
	{\em
		In Section~\ref{section:preliminaries-model-theory}, we will introduce a relation $\lr{\alg{A},\Delta}\models\Forall\ldots\Forall A$, where
		$\lr{\alg{A},\Delta}$ is an expansion of an algebra $\alg{A}$ and $\Forall\ldots\Forall A$ is an identity of quasi-identity. From our explanation there, it will be clear that 
		\[
		\alg{A}\models A~\Longleftrightarrow~\lr{\alg{A},\Delta}\models\Forall\ldots\Forall A.
		\]
		This is supposed to show a relationship between equalities and identities, on the one hand, and quasi-equalities and quasi-identities, on the other.
	}
\end{rem}

Let $\Sigma$ be a nonempty set of equalities or quasi-equalities. We write
\[
\alg{A}\models\Sigma
\]
if $\alg{A}\models A$, for all $A\in\Sigma$.\\

A nonempty class $\mathcal{K}$ of $\Lan$-algebras is an \textit{\textbf{equational class}}\index{equational class} if there is a set $\Sigma$ of identities such that
\[
\alg{A}\in\mathcal{K}~\Longleftrightarrow~\alg{A}\models\Sigma.\tag{$\ast\ast$}
\]

There is a very well-known and useful characterization of equational classes due G. Birkhoff.
\begin{prop}[Birkhoff, cf.~\cite{malcev1973}, {\S} 13, theorem 1] \label{P:Birkhoff's-criterion-for varieties}
	Let $\mathcal{K}$ be a nonempty class of algebras of type $\Lan$. Then $\mathcal{K}$ is an equational class if, and only if, the following conditions are satisfied:
	{\em\[
		\begin{array}{cl}
			(\text{a}) & \Su\mathcal{K}\subseteq\mathcal{K};\\
			(\text{b}) & \Ho\mathcal{K}\subseteq\mathcal{K};\\
			(\text{c}) &  \Pro\mathcal{K}\subseteq\mathcal{K}.
		\end{array}
		\]}
\end{prop}
A class $\mathcal{K}$ is called a \textit{\textbf{variety}}\index{variety} if the closures (a)--(c) of Proposition~\ref{P:Birkhoff's-criterion-for varieties} hold. Thus for a class to be equational and to be a variety are two equivalent conditions.

\begin{prop}[Tarski, cf.~\cite{malcev1973}, {\S} 13, theorems 2]
	The smallest variety containing a class $\mathcal{K}$ equals $\Ho\Su\Pro\mathcal{K}$.
\end{prop}

Given a class $\mathcal{K}$ of algebras of type $\Lan$, the variety $\Ho\Su\Pro\mathcal{K}$ is said to be \textit{\textbf{generated}} by $\mathcal{K}$, or one can also say that $\mathcal{K}$ \textit{\textbf{generates}} this variety.

Here is another useful characterization.
\begin{prop}[\cite{burris-sankapp1981}, chapter II, theorem 11.12]\label{P:HSP-K=HP_s-K}
	Let $\mathcal{K}$ be a class of algebras of type $\Lan$. Then
	$\Ho\Su\Pro\mathcal{K}=\Ho\Pros\mathcal{K}$	
\end{prop}

We conclude our discuss about equational classes with the following important property formulated in terms of variety.
\begin{prop}\label{P:variety-generated-si-algebras}
	Let $\mathcal{K}$ be a class of algebras of type $\Lan$. Then the variety generated by $\mathcal{K}$ coincides with the variety generated by {\em$\Si\mathcal{K}$}, that is, {\em$\Ho\Su\Pro\mathcal{K}=\Ho\Su\Pro\Si\mathcal{K}$}.
\end{prop}
\begin{proof}
	Our starting point is the inclusion $\mathcal{K}\subseteq\Pros\Si\mathcal{K}$ which is true in virtue of Proposition~\ref{P:K-included-P_s-K_si}. Then, we continue:
	\[
	\begin{array}{rl}
		\mathcal{K}\subseteq\Pros\Si\mathcal{K}\!\!&\Longrightarrow
		\mathcal{K}\subseteq\Su\Pro\Si\mathcal{K}\quad[\text{by Proposition~\ref{P:class-operators-inequalities}}]\\
		&\Longrightarrow\Pro\mathcal{K}\subseteq\Pro\Su\Pro\Si\mathcal{K}
		\subseteq\Su\Pro\Pro\Si\mathcal{K}=\Su\Pro\Si\mathcal{K}
		\quad[\text{by Proposition~\ref{P:class-operators-inequalities}}]\\
		&\Longrightarrow\Su\Pro\mathcal{K}\subseteq\Su\Su\Pro\Si\mathcal{K}
		=\Su\Pro\Si\mathcal{K}\quad[\text{by Proposition~\ref{P:class-operators-inequalities}}]\\
		&\Longrightarrow\Ho\Su\Pro\mathcal{K}\subseteq\Ho\Su\Pro\Si\mathcal{K}.
	\end{array}
	\]	
	
	To prove the other inclusion, we show that $\Si\mathcal{K}\subseteq\Ho\Su\Pro\mathcal{K}$. In proving the last inclusion, we have to remember that the last class is a variety and therefore, it is closed under the operator $\Ho$. (Proposition~\ref{P:Birkhoff's-criterion-for varieties})
	
	Assume $\alg{A}\in\Si\mathcal{K}$. Then there is an algebra $\alg{B}\in\mathcal{K}$ and a set $\lbrace\alg{A}_i\rbrace_{i\in I}\subseteq\mathcal{S}_{\Lan}$ such that $\alg{B}$ is a subdirect product of
	$\lbrace\alg{A}_i\rbrace_{i\in I}$ and $\alg{A}=\alg{A}_i$, for some $i\in I$.
	Since $\alg{A}$ is a homomorphic image of $\alg{B}$,
	$\alg{A}\in\Ho\Su\Pro\mathcal{K}$.
	
	Finally, applying the second part of Proposition~\ref{P:Birkhoff's-criterion-for varieties}, we conclude that
	$\Ho\Su\Pro\Si\mathcal{K}\subseteq\Ho\Su\Pro\mathcal{K}$.
\end{proof}

Now we turn to implicative classes, or quasi-varieties. 

A nonempty class $\mathcal{K}$ of $\Lan$-algebras is an \textit{\textbf{implicative class}},\index{implicative class} or \textit{\textbf{quasi-variety}}\index{quasi-variety} if there is a set $\Sigma$ of quasi-equalities such that the equivalence $(\ast\ast)$ holds.

The following characterization is due to A. Mal'cev.
\begin{prop}[\cite{malcev1973}, {\S} 11, theorem 2]\label{P:Mal'cev's-for-quasivarieties}
	Let $\mathcal{K}$ be a nonempty abstract class of algebras of type $\Lan$. Then $\mathcal{K}$ is a quasi-variety if, and only if, the following conditions are satisfied:
	{\em\[
		\begin{array}{cl}
			(\text{a}) & \text{a trivial (degenerate) $\Lan$-algebra belongs to $\mathcal{K}$};\\
			(\text{b}) & \Su\mathcal{K}\subseteq\mathcal{K};\\
			(\text{c}) &  \Pro\mathcal{K}\subseteq\mathcal{K};\\
			(\text{d}) &  \Pu\mathcal{K}\subseteq\mathcal{K}.
		\end{array}
		\]}
\end{prop}

For any nonempty class $\mathcal{K}$ of algebras of type $\Lan$ and any set $X\subseteq\Var_\Lan$, we define:
\[
\theta_{\mathcal{K}}(X):=\bigcap\set{\theta\in\Congruence\,\FormAl}{X^{2}\subseteq\theta~\text{and}~\FormAl\slash\theta\in\Is\Su\mathcal{K}}.
\]

This gives rise for the following definition of algebra $\FormAl\slash
\theta_{\mathcal{K}}(X)$ with the generators $\overline{X}:=\set{x\slash\theta_{\mathcal{K}}(X)}{x\in X}$. This algebra will also be denoted by 
$\mathfrak{F}_{\mathcal{K}}(\overline{X})$.

\begin{prop}[\cite{burris-sankapp1981}, theorem 10.10 and corollary 10.11]
	Given a nonempty class $\mathcal{K}$ of algebras of type $\Lan$, for any set $X\subseteq\Var_\Lan$, the algebra $\mathfrak{F}_{\mathcal{K}}(\overline{X})$ has the universal mapping property relative to the class $\mathcal{K}$. Moreover, given  an algebra {\em$\alg{A}\in\mathcal{K}$},
	for sufficiently large $X\subseteq\Var_{\Lan}$,
	{\em$\alg{A}\in\Ho\mathfrak{F}_{\mathcal{K}}(\overline{X})$}.
\end{prop}

The following proposition is a generalization of theorem 10.12 of~\cite{burris-sankapp1981}.
\begin{prop}\label{P:Birkhpff's-for-algebras}
	Let $\mathcal{K}$ be a nonempty class of algebras of type $\Lan$. Then the
	algebra $\FormAl\slash
	\theta_{\mathcal{K}}(X)$ belongs to the class $\Su\Pro\mathcal{K}$. Hence, if $\Su\mathcal{K}\subseteq\mathcal{K}$ and
	$\Pro\mathcal{K}\subseteq\mathcal{K}$, in particular if $\mathcal{K}$ is a quasi-variety, then  $\FormAl\slash
	\theta_{\mathcal{K}}(X)$ belongs to $\mathcal{K}$.
\end{prop}
\begin{proof}
	First we note that
	\[
	\FormAl\slash
	\theta_{\mathcal{K}}(X)\in\Pros\set{\FormAl\slash\theta}{X\subseteq\theta~\text{and}~\FormAl\slash\theta\in\Su\mathcal{K}};
	\]
	see~\cite{burris-sankapp1981}, {\S}8, exercise 11. This implies that
	\[
	\FormAl\slash
	\theta_{\mathcal{K}}(X)\in\Pros\Pro\mathcal{K}.
	\]
	Using Proposition~\ref{P:class-operators-inequalities}, we obtain that 
	$\FormAl\slash
	\theta_{\mathcal{K}}(X)\in\Su\Pro\mathcal{K}$.
\end{proof}

Given a variety (or quasi-variety) $\mathcal{K}$ of algebras of type $\Lan$ and $X\subseteq\Var_{\Lan}$, the algebra $\mathfrak{F}_{\mathcal{K}}(\overline{X})$\index{$\mathfrak{F}_{\mathcal{K}}(\overline{X})$} is called a \textit{\textbf{free algebra relative to}}\index{algebra!free} $\mathcal{K}$. It is well-known that, given $X\cup Y\subseteq\Var_{\Lan}$, if $\card{\overline{X}}=\card{\overline{Y}}$, then
the algebras $\mathfrak{F}_{\mathcal{K}}(\overline{X})$ and $\mathfrak{F}_{\mathcal{K}}(\overline{Y})$ are isomorphic. Then, if $\kappa=\card{\overline{X}}$, then the algebra $\mathfrak{F}_{\mathcal{K}}(\overline{X})$ is also denoted by $\mathfrak{F}_{\mathcal{K}}(\kappa)$ and call the\textit{\textbf{ free algebra of rank}} $\kappa$ relative to $\mathcal{K}$.
\paragraph{References}
\begin{enumerate}
	\item~\cite{burris-sankapp1981}
	\item~\cite{malcev1973}	
\end{enumerate}

\subsection{Hilbert algebras}
An algebra $\alg{A}=\langle\textsf{A};\rightarrow,\one\rangle$, where $\rightarrow$ is a binary operation and $\one$ (the \textit{unit}) a 0-ary operation (or a constant), is called a \textit{\textbf{Hilbert algebra}}\footnote{The notion of a Hilbert algebra was introduced in~\cite{diego1962}\index{Hilbert algebra}\index{algebra!Hilbert}; see also~\cite{diego1966}. In~\cite{ras74} Hilbert algebras are called \textit{positive implication algebras}. Diego~ also proved (see~\cite{diego1966}, theorem 3) that the Hilbert algebras form an equational class.} if for arbitrary elements $x,y,z\in\textsf{A}$, the following conditions are satisfied:
\[
\begin{array}{cl}
	(\text{d}_1) &x\rightarrow\one=\one,\\
	(\text{d}_2) &x\rightarrow(y\rightarrow x)=\one,\\
	(\text{d}_3) &(x\rightarrow(y\rightarrow z))\rightarrow((x\rightarrow y)\rightarrow(x\rightarrow z))=\one,\\
	(\text{d}_4) &(x\rightarrow y=\one~\text{and}~y\rightarrow x=\one)~\Longleftrightarrow~x=y.
\end{array}
\]

In any Hilbert algebra the following hold:
\[
x\rightarrow x=\one,
\]
and
\begin{equation}\label{E:diego-algebra-1}
	(x=\one~\text{and}~x\rightarrow y=\one)~\Longrightarrow~y=\one.
\end{equation}

The relation $\le$ defined as
\begin{equation}\label{E:diego-algebra-ordering}
	x\le y~\stackrel{\text{df}}{\Longleftrightarrow}~x\rightarrow y=\one
\end{equation}
is a partial ordering on \textsf{A} and, hence, $\one$ is the greatest element with respect to $\le$.
\paragraph{References}
\begin{enumerate}
	\item~\cite{diego1962}
	\item~\cite{ras74}	
\end{enumerate}

\subsection{Distributive lattices}
\label{section:distributive-lattices}

An algebra $\alg{A}=\langle\textsf{A};\wedge,\vee\rangle$ of type $\langle2,2\rangle$, where $\wedge$ (\textit{meet}) and $\vee$ (\textit{join}) are binary operations, is called a \textit{\textbf{distributive lattice}}\index{lattice!distributive } if the following equalities hold in $\alg{A}$ for arbitrary elements $x$, $y$ and $z$ of \textsf{A}:
\[
\begin{array}{cll}
	(\text{l}_1) &i)~~~x\wedge y=y\wedge x, &ii)~~~x\vee y=y\vee x,\\
	(\text{l}_2) &i)~~~x\wedge(y\wedge z)=(x\wedge y)\wedge z,
	&ii)~~~x\vee(y\vee z)=(x\vee y)\vee z,\\
	(\text{l}_3) &i)~~~(x\wedge y)\vee y=y, &ii)~~~x\wedge(x\vee y)=x,\\
	(\text{l}_4) &i)~~~x\wedge(y\vee z) =(x\wedge y)\vee(x\wedge z),
	&ii)~~~x\vee(y\wedge z) =(x\vee y)\wedge(x\vee z).\\
\end{array}
\]

In any distributive lattice the following identities hold:
\[
x\wedge x=x~~\text{and}~~x\vee x= x \tag{\textit{idempotent laws}}
\]

Indeed, for any $x,y\in\textsf{L}$,
\[
\begin{array}{rl}
	x\!\!\! &= x\wedge(x\vee y)\\&=(x\wedge x)\vee(x\wedge y)\quad[\mbox{according to ($\text{l}_{4}$--$i$)}]\\
	&= ((x\wedge y)\vee x)\wedge((x\wedge y)\vee x)\quad[\mbox{in virtue of ($\text{l}_{1}$--$ii$) and ($\text{l}_{4}$--$ii$)}]\\
	&= x\wedge x.\quad[\mbox{according to ($\text{l}_{3}$--$i$)}]
\end{array}
\]

The second idempotent law can be proven in a similar fashion. \\

If in the definition of a distributive lattice, we replace the laws $(\text{l}_4)$ with the idempotent laws, we receive the definition of a
\textit{\textbf{lattice}}. In this book, we will mostly deal with distributive lattices. However, some laws below are true for any lattice.\\

We observe that if $\alg{A}$ is a lattice, then for any $x,y\in\textsf{A}$,
\[
x\wedge y=x~\Longleftrightarrow~x\vee y=y.
\]

Indeed, assume that $x\wedge y=x$. Then we have: $x\vee y=(x\wedge y)\vee y=y$.

On the other hand, if $x\vee y=y$, then $x\wedge y=x\wedge(x\vee y)=x$.

The last equivalence induces the following definition. Given a lattice $\alg{A}$, we define a binary relation on $\textsf{A}$ as follows:
\begin{equation}\label{E:ordering-in-lattice}
	x\le y\stackrel{\text{df}}{\Longleftrightarrow} x\wedge y=x.
\end{equation}

According to the last equivalence, we have:
\[
x\le y~\Longleftrightarrow~x\vee y=y.
\]

Further, we notice that $\le$ is a partial ordering on \textsf{A} and, according to this ordering, $x\wedge y$ is the greatest lower bound and $x\vee y$ is the least upper bound of $\lbrace x,y\rbrace$, respectively. This, in particular, implies that the operations $\wedge$ and $\vee$ are monotone with respect to $\le$. Namely,
\[
x\le y \Longrightarrow (\text{$x\wedge z\le y\wedge z$ and $x\vee z\le y\vee z$}).
\]

A lattice \alg{A} is \textit{\textbf{complete}}\index{lattice!complete} if for any subset $X\subseteq\textsf{A}$, there are a greatest lower bound of $X$, symbolically 
$\bigwedge X$, and a least upper bound of $X$, symbolically $\bigvee X$,
with respect to $\le$.

\paragraph{References}
\begin{enumerate}
	\item~\cite{ras74}	
	\item~\cite{rs70}
\end{enumerate}
\subsection{Brouwerian semilattices}
An algebra $\alg{A}=\langle\textsf{A};\land,\rightarrow,\one\rangle$ is a \textit{\textbf{Brouwerian semilattice}}\index{semilattice!Brouwerian} (aka \textit{relatively pseudo-complemented semilattice}\index{semilattice!relatively pseudo-complemented} or \textit{implicative semilattice})\index{semilattice!implicative} if $\langle\textsf{A};\rightarrow,\one\rangle$ is a Hilbert algebra and $\langle\textsf{A};\land\rangle$ is a meet semilattice, that is, the properties
($\text{l}_{1}$--$i$) and ($\text{l}_{2}$--$i$) are satisfied; in addition, the following equalities hold:
\[
\begin{array}{cl}
	(\text{s}_1) &(x\land y)\rightarrow x=\one,\\
	(\text{s}_2) &(x\rightarrow y)\rightarrow((x\rightarrow z)\rightarrow(x\rightarrow(y\land z)))=\one.\\
\end{array}
\]

In any Brouwerian semilattice the equivalence 
\[
x\land y=x~\Longleftrightarrow~x\rightarrow y=\one.
\]
is true. This implies that the relations defined by~\eqref{E:diego-algebra-ordering} and~\eqref{E:ordering-in-lattice}  coincide. Thus we employ one symbol, $\le$, for this relation and use \eqref{E:diego-algebra-1} and \eqref{E:ordering-in-lattice} either as a defining equivalence or as a property.

The following properties hold in any Brouwerian semilattice:
\begin{equation}\label{E:relative-pseudo-complement-1}
	x\le y\rightarrow z~\Longleftrightarrow~x\land y\le z
\end{equation}
and
\begin{equation}\label{E:brouwerian-conjunction}
	x\land y=\one~\Longleftrightarrow~(x=\one~\text{and}~y=\one).
\end{equation}

\paragraph{References}
\begin{enumerate}
	\item~\cite{kohler1981}
	\item~\cite{nemitz1965}
\end{enumerate}

\subsection{Relatively pseudo-complemented lattices}
An algebra $\alg{A}=\langle\textsf{A};\land,\lor,\rightarrow,\one\rangle$ is a \textit{\textbf{relatively pseudo-complemented lattice}}\index{lattice!relatively pseudo-complemented lattice} if $\langle\textsf{A};\land,\lor\rangle$ is a distributive lattice and
$\langle\textsf{A};\land,\rightarrow,\one\rangle$ is a Bouwerian semilattice; in addition, the following equalities hold:
\[
\begin{array}{cl}
	(\text{p}_1) &x\rightarrow(x\lor y)=\one,\\
	(\text{p}_2) &(x\rightarrow z)\rightarrow((y\rightarrow z)\rightarrow((x\lor y)\rightarrow z))=\one.\\
\end{array}
\]

In fact, there is no need to require $\langle\textsf{A};\land,\lor\rangle$ to be distributive. The distributivity follows from~\eqref{E:relative-pseudo-complement-1}.

Also,  in virtue of~\eqref{E:relative-pseudo-complement-1}, in any relatively pseudo-complemented lattice,
\[
x\le y\rightarrow y;
\] 
that is, for any element $y$, $y\rightarrow y$ is the greatest element with respect to $\le$, that is, $y\rightarrow y=\one$.

\paragraph{References}
\begin{enumerate}
	\item~\cite{ras74}
	\item~\cite{rs70}	
\end{enumerate}

\subsection{Boolean algebras}\label{section:boolean-algebra}
Boolean algebras belong to the class of distributive lattices.

A \textit{\textbf{Boolean algebra}}\index{algebra!Boolean} is an algebra $\alg{A}=\langle\textsf{A};\wedge,\vee,\neg,\one\rangle$, where $\wedge$ and $\vee$ are the operations of a distributive lattice, that is, they satisfy  the equalities $(\text{l}_1)$--$(\text{l}_4)$  above, $\neg$ (\textit{complementation}) is a unary operation and $\one$ is a 0-ary operation, if $(\text{b}_1)$--$(\text{b}_2)$ below are satisfied in \alg{A} for arbitrary elements $x$, $y$ and $z$ of \textsf{A}:
\[
\begin{array}{cll}
	(\text{b}_1) &i)~~~ x\wedge\one=x, &ii)~~~x\vee\one=\one,\\
	(\text{b}_2) &i)~~~(x\wedge\neg x)\vee y=y,
	&ii)~~~(x\vee\neg x)\wedge y=y.\\
\end{array}
\]

If we apply the definition~\eqref{E:ordering-in-lattice} for a Boolean algebra $\alg{A}$, we realize that $\one$ is the greatest element and $\neg\one$ the least element in $\alg{A}$, respectively. Indeed, the former comes from 
($\text{b}_{1}$--$i$) and the latter, with the help of ($\text{b}_{2}$--$i$)
and ($\text{b}_{1}$--$i$), can be obtained as follows:
\[
x=(\one\wedge\neg\one)\vee x= \neg\one\vee x.
\]

Denoting
\begin{equation}\label{E:zero-definition}
	\zero:=\neg\one,
\end{equation}
we observe that for any $x\in\textsf{B}$, $\neg x$ is a unique element $y$ such that $x\wedge y=\zero$ and $x\vee y=\one$. The element $\neg x$ is called a \textit{complement} of $x$. Thus any Boolean algebra is a distributive lattice with a greatest element $\one$, a least element $\zero$ and a \textit{complementation}\index{complementation} $\neg x$. This implies that, in view of commutativity of $\wedge$ and $\vee$ (the equalities $(\text{l}_1)$), we can say the $x$ is the complement of $\neg x$, which immediately implies the equality
\begin{equation}\label{E:double-negation-equality}
	\neg\neg x=x.
\end{equation}

We also observe that
\begin{equation}\label{E:one-zero-in-boolean}
	\one=x\vee\neg x~~\text{and}~~\zero=x\wedge\neg x,
\end{equation}
for an arbitrary element $x\in\textsf{A}$.

We denote:
\[
x\rightarrow y:=\neg x\lor y.
\]
and
\[
x\leftrightarrow y:=(x\rightarrow y)\land(y\rightarrow x).
\]

The just defined implication $x\rightarrow y$ satisfies all the laws of the relatively complemented lattice. Thus a Boolean algebra can be defined as an expansion of a relatively complemented lattice with a new unary operation $\neg$ which satisfies ($\text{b}_1$)--($\text{b}_2$) and the identity $x\rightarrow y=\neg x\lor y$. This allows one to consider a Boolean algebra
in the signature $\langle\land,\lor,\rightarrow,\neg,\one\rangle$.

It is well known that given a Boolean algebra $\alg{A}=\langle\textsf{A};\wedge,\vee,\neg,\one\rangle$, the algebra
$\langle\textsf{A};\wedge,\vee,\rightarrow,\one\rangle$ is a relatively pseudo-complemented lattice; cf.~\cite{rs70}, theorem II.1.1, or~\cite{ras74}, theorem VI.1.2. Then, with the help of~\eqref{E:brouwerian-conjunction}, we immediately obtain:
\begin{equation}\label{E:relative-psedo-complement-3} 
	x\leftrightarrow y=\one~\Longleftrightarrow~x=y.
\end{equation}

It is obvious that a Boolean algebra is nontrivial if, and only if, $\zero\neq\one$. This implies that the simplest nontrivial Boolean algebras consists of two elements --- $\one$ and $\zero$.
Usually, this algebra is depicted by the following diagram.
\begin{figure}[!ht]
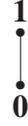
	
	\[
	\ctdiagram{
		\ctnohead
		\ctinnermid
		\ctel 0,0,0,20:{}
		\ctv 0,0:{\bullet}
		\ctv 0,20:{\bullet}
		\ctv 0,27:{\mathbf{1}}
		\ctv 0,-9:{\mathbf{0}}
	}
	\]\label{figure:two-element-boolean}
	\caption{A 2-element Boolean algebra}
	
\end{figure}

A simple observation gives us the following.
\begin{prop}\label{P:2-element-algebra}
	In any Boolean algebra {\em$\alg{A}=\langle\textsf{A};\wedge,\vee,\neg,\one\rangle$}, the algebra
	$\langle\lbrace\zero,\one\rbrace;\wedge,\vee,\neg,\one\rangle$ is a subalgebra of {\em\alg{A}}.
\end{prop}

\begin{prop}\label{P:finitely-generated-boolean}
	Let a Boolean algebra {\em$\alg{A}$} be generated by a nonempty set {\em$\textsf{A}_0$}. Then any element {\em$x\in\textsf{A}$} can be represented in each of the forms:
	\[
	x=\bigvee_{i=1}^{m}\bigwedge_{j=1}^{n_i}b_{ij},\tag{\textit{disjunctive normal decomposition}}
	\]
	or in the form:
	\[
	x=\bigwedge_{i=1}^{m}\bigvee_{j=1}^{n_i}b_{ij},\tag{\textit{conjunctive normal decomposition}}
	\]
	where for any $i$ and $j$ either {\em$b_{ij}\in\textsf{A}_0$}
	or {\em$\neg b_{ij}\in\textsf{A}_0$}. 
\end{prop}\index{disjunctive normal decomposition}\index{conjunctive normal decomposition}
\begin{proof}
	To prove this proposition, it suffices to notice that both decompositions belong to \alg{A}. Then, one shows that all elements that have a disjunctive normal decomposition constitute a subalgebra of $\alg{A}$; and the same is true for the elements that have a conjunctive normal decomposition; see~\cite{rs70}, chapter II, theorem 2.1, for detail.
\end{proof}

\begin{prop}[\cite{burris-sankapp1981}, corollary 1.9]\label{P:subdirectly-irredicible-Boolean}
	Up to isomorphism, there is only one subdirectly irreducible Boolean algebra, a 2-element Boolean algebra.
\end{prop}

As we will see below (Proposition~\ref{P:boolean-algebra-as-heyting}), the class of Boolean algebras can be regarded as a subclass of a larger class of Heyting algebras which are discussed in the next subsection. Therefore, many properties related to Heyting algebras are applicable to Boolean algebras.
For example, the properties that are discussed below in Proposition~\ref{P:Int-properties} for Heyting algebras (Section~\ref{section:heyting-algebra}) can be applied to Boolean algebras as well.
\paragraph{References}
\begin{enumerate}
	\item~\cite{burris-sankapp1981}
	\item~\cite{ras74}
	\item~\cite{rs70}	
\end{enumerate}

\subsection{Heyting algebras}\label{section:heyting-algebra}
A \textit{\textbf{Heyting algebra}} (aka \textit{pseudo-Boolean algebra})\index{algebra!Heyting}\index{algebra!pseudo-Boolean} is an algebra $\alg{A}=\langle\textsf{A};\wedge,\vee,\rightarrow,\neg,\one\rangle$ of type $\langle2,2,2,1,0\rangle$, where $\wedge$ and $\vee$, $\rightarrow$ (\textit{relative pseudo-complementation}) are binary operations, $\neg$ (\textit{pseudo-complementation}) is a unary operation and $\one$
(the unit) is a 0-ary operation, if, besides the equalities ($\text{l}_{1}$)--($\text{l}_{2}$) and ($\text{b}_{1}$) (Section~\ref{section:boolean-algebra}), the following equalities are satisfied for arbitrary elements $x$, $y$, $z$ of \textsf{H}:
\[
\begin{array}{cl}
	(\text{h}_1) &x\wedge(x\rightarrow y)=x\wedge y,\\
	(\text{h}_2) &(x\rightarrow y)\wedge y=y,\\
	(\text{h}_3) &(x\rightarrow y)\wedge(x\rightarrow z)=x\rightarrow(y\wedge z),\\
	(\text{h}_4) &x\wedge(y\rightarrow y)=x,\\
	(\text{h}_5) &\neg\one\vee y=y,\\
	(\text{h}_6) &\neg x=x\rightarrow\neg\one.\\
\end{array}
\]

It is customary to denote:
\[
x\leftrightarrow y:=(x\rightarrow y)\wedge(y\rightarrow x).
\]

As in the case of Boolean algebras, in Heyting algebras, in virtue of
$(\text{b}_1)$, the element $\one$ is a greatest element in Heyting algebra.
The property $(\text{h}_4)$ implies that
\begin{equation}\label{E:x-implies-x}
	x\rightarrow x=\one.
\end{equation}

Using the notation \eqref{E:zero-definition} and the property $(\text{h}_5)$, we conclude that $\zero$ is a least element in Heyting algebra. According to $(\text{h}_6)$, we have:
\[
\neg x=x\rightarrow\zero.
\]

Thus, the last identity and $(\text{h}_1)$ imply the following:
\begin{equation}\label{E:zero-in-heyting}
	x\wedge\neg x=\zero.
\end{equation}

The following property characterizes pseudo-complementation:
\begin{equation}\label{E:pseudo-complementation}
	x\le y\rightarrow z~\Longleftrightarrow~x\wedge y\le z.
\end{equation}

Indeed, assume first that $x\le y\rightarrow z$. Then, in view of ($\text{l}_1$--$i$), the monotonicity of $\wedge$ with respect to $\le$ and $(\text{h}_1)$, we have:
\[
x\wedge y\le y\wedge(y\rightarrow z)=y\wedge z\le z.
\]

Conversely, suppose that $x\wedge y\le z$, that is $x\wedge y= x\wedge y\wedge z$. Then, in virtue of $(\text{h}_2)$, $(\text{h}_4)$, $(\text{h}_3)$,
we obtain:
\[
\begin{array}{rl}
	x\le y\rightarrow x\!\!\!&=(y\rightarrow x)\wedge(y\rightarrow y)=y\rightarrow(x\wedge y)=y\rightarrow(x\wedge y\wedge z)\\
	&=(y\rightarrow(x\wedge y))\wedge(y\rightarrow z)\le y\rightarrow z.
\end{array}
\]

Using~\eqref{E:pseudo-complementation}, we receive  immediately:
\begin{equation}\label{E:less-than=implication}
	x\le y~\Longleftrightarrow~x\rightarrow y=\one;
\end{equation}
which in turn implies:
\[
x=y~\Longleftrightarrow~x\leftrightarrow y=\one.
\]
\begin{prop}\label{P:boolean-algebra-as-heyting}
	Let {\em$\alg{A}=\langle\textsf{A};\wedge,\vee,\rightarrow,\neg,\one\rangle$} be a Heyting algebra. Then its restriction {\em$\alg{B}=\langle\textsf{A};\wedge,\vee,\neg,\one\rangle$} is a Boolean algebra if the identity $x\rightarrow y=\neg x\vee y$ holds in {\em\alg{A}}.
	Conversely, if {\em$\alg{B}=\langle\textsf{B};\wedge,\vee,\neg,\one\rangle$} is a Boolean algebra. Then its expansion {\em$\alg{A}=\langle\textsf{B};\wedge,\vee,\rightarrow,\neg,\one\rangle$},
	where $x\rightarrow y:=\neg x\vee y$, is a Heyting algebra.
\end{prop}
\begin{proof}
	Let \alg{A} be a Heyting algebra in which $x\rightarrow y=\neg x\vee y$.
	We observe that, because of~\eqref{E:zero-in-heyting}, the identity ($\text{b}_2$--$i$) is true in \alg{A}. The identity ($\text{b}_2$--$ii$) is also true, for we have:
	\[
	(x\vee\neg x)\wedge y=(x\rightarrow x)\wedge y=y.
	\]
	
	Now assume that \alg{A} is an expansion of a Boolean algebra \alg{B}, where, by definition, $x\rightarrow y=\neg x\vee y$. We aim to show that all the properties $(\text{h}_1)$--$(\text{h}_6)$ are valid in \alg{A}. Indeed, we obtain:
	\[
	\begin{array}{l}
		x\wedge(x\rightarrow y)=x\wedge(\neg x\vee y)=(x\wedge\neg x)\vee(x\wedge y)=x\wedge y;\\
		(x\rightarrow y)\wedge y=(\neg x\vee y)\wedge y=y;\\
		(x\rightarrow y)\wedge(x\rightarrow z)=(\neg x\vee y)\wedge(\neg x\vee z)=
		\neg x\vee(y\wedge\neg x)\vee(\neg x\wedge z)\vee(y\wedge z)\\
		\qquad\qquad\qquad\qquad~=\neg x\vee(y\wedge z)=x\rightarrow(y\wedge z);\\
		x\wedge(y\rightarrow y)=x\wedge(\neg y\vee y)=x;\\
		\neg x=\neg x\vee\zero=\neg x\vee\neg\one=x\rightarrow\neg\one.
	\end{array}
	\]
\end{proof}

\begin{cor}\label{C:Cl-property}
	Let {\em$\alg{A}=\langle\textsf{A};\wedge,\vee,\neg,\one\rangle$} be a Boolean algebra. Then for an arbitrary {\em$x\in\textsf{A}$}, $\neg\neg x\rightarrow x=\one$, where the operation $x\rightarrow y$ is understood as $\neg x\vee y$.
\end{cor}
\begin{proof}
	Augmenting \alg{A} with $\rightarrow$ defined as above, we receive a Heyting algebra. Then we use~\eqref{E:x-implies-x} and~\eqref{E:double-negation-equality}.
\end{proof}

It should be clear that the simplest nontrivial Boolean algebra (Fig.~\ref{figure:two-element-boolean}) is also the simplest non-generate Heyting algebra.\\

We will be using Heyting algebras as a basic semantics for the intuitionistic propositional logic (Section~\ref{section:some-lindenbaum-algebras}). For this purpose, in the sequel, we will need the following propertied of Heyting algebras.
\begin{prop}\label{P:Int-properties}
	Let {\em$\alg{A}=\langle\textsf{A};\wedge,\vee,\rightarrow,\neg,\one\rangle$} be a Heyting algebra. For arbitrary elements $x$, $y$ and $z$ of {\em\textsf{H}} the following properties hold:
	{\em\[
		\begin{array}{cl}
			(\text{a}) &x\le y\rightarrow x,\\
			(\text{b}) &x\rightarrow y\le(x\rightarrow (y\rightarrow z))
			\rightarrow(x\rightarrow z),\\
			(\text{c}) &x\le y\rightarrow(x\wedge y),\\
			(\text{d}) &x\wedge y\le x,\\
			(\text{e}) &x\le x\vee y,\\
			(\text{f}) &x\rightarrow z\le(y\rightarrow z)\rightarrow((x\vee y)\rightarrow z),\\
			(\text{g}) &x\rightarrow y\le(x\rightarrow\neg y)\rightarrow\neg x,\\
			(\text{h}) &x\le\neg x\rightarrow y.
		\end{array}
		\]}
\end{prop}
\begin{proof}
	We obtain successively the following.
	
	(a) follows straightforwardly from $(\text{h}_2)$.
	
	To prove (b), using~\eqref{E:pseudo-complementation}, we have:
	\[
	\begin{array}{rl}
		(x\rightarrow y)\wedge(x\rightarrow (y\rightarrow z))\wedge x\!\!\!
		&=(x\rightarrow y)\wedge x\wedge(y\rightarrow z)\\
		&=x\wedge y\wedge (y\rightarrow z)=x\wedge y\wedge z\le z.
	\end{array}
	\]
	Then, we apply \eqref{E:pseudo-complementation} again twice to get (b).
	
	We derive (c) immediately from $x\wedge y\le x\wedge y$.
	
	(d) and (e) are obvious.
	
	To prove (f), we first have:
	\[
	\begin{array}{rl}
		(x\rightarrow z)\wedge(y\rightarrow z)\wedge(x\vee y)\!\!\!
		&=((x\rightarrow z)\wedge(y\rightarrow z)\wedge x)\vee
		((x\rightarrow z)\wedge(y\rightarrow z)\wedge y)\\
		&=(x\wedge z\wedge(z\rightarrow y))\vee(y\wedge z\wedge(x\rightarrow z))\\
		&=(x\wedge z)\vee(y\wedge z)=(x\vee y)\wedge z\le z.
	\end{array}
	\]
	Then, we apply~\eqref{E:pseudo-complementation} twice.
	
	To prove (g), applying~\eqref{E:pseudo-complementation} twice, we have:
	\[
	(x\rightarrow y)\wedge (x\rightarrow \neg y)\wedge x=x\wedge y\wedge\neg y=x\wedge\zero\le\zero,
	\]
	which implies that
	\[
	(x\rightarrow y)\wedge (x\rightarrow \neg y)\le x\rightarrow\zero=\neg x.
	\]
	It remains to apply~\eqref{E:pseudo-complementation}.
	
	Finally, we obtain (h) as follows:
	\[
	x\wedge\neg x=\zero\le y.
	\]
	Then, we apply~\eqref{E:pseudo-complementation}.
\end{proof}

Two obvious ways to obtain new Heyting algebras from old ones are the formation of subalgebras and homomorphic images. For instance, given a nontrivial Heyting algebra $\langle\textsf{A};\wedge,\vee,\rightarrow,\neg,\one\rangle$,
the algebra $\langle\lbrace\mathbf{0},\mathbf{1}\rbrace;\wedge,\vee,\rightarrow,\neg,\one\rangle$, where $\mathbf{0}:=\neg\mathbf{1}$, the smallest subalgebra of the first algebra.

The following method is very useful for obtaining homomorphic images of a given Heyting algebra. 

Let $\alg{A}=\langle\textsf{A};\wedge,\vee,\rightarrow,\neg,\one\rangle$ be a nontrivial Heyting algebra and $b\in\textsf{A}$. We consider the algebra
$\alg{A}_{b}=\langle\set{x\in\textsf{A}}{x\le b};\wedge_b,\vee_b,\rightarrow_b,\neg_b,\one_b\rangle$, where, by definition:
\[
\begin{array}{l}
	x \wedge_{b}y:=x\wedge y;\\
	x\vee_b{b}y:=x\vee y;\\
	x\rightarrow_{b}y:=(x\rightarrow)\wedge b;\\
	\neg_{b}x:=(\neg x)\wedge b;\\
	\mathbf{1}_{b}:=b.
\end{array}
\]

It turns out that $\alg{A}_{b}$ is a homomorphic image of $\alg{A}$ with respect to the homomorphic map $x\mapsto x\wedge b$. See~\cite{rs70}, chapter IV, {\S} 8, where the formation of $\alg{A}_b$ is called \emph{relativization}.

\paragraph{References}
\begin{enumerate}
	\item~\cite{ras74}
	\item~\cite{rs70}	
\end{enumerate}

\section{Preliminaries from model theory}\label{section:preliminaries-model-theory}
Let a first-order language $\FOL$ consist of a set $\Var$ of individual variables, logical constants `$\Land$', `$\Lor$, `$\Rarrow$', `$\Neg$' and quantifiers `$\Forall$' and `$\Exists$', as well as parentheses `$($' and 
`$)$'. Moreover, a signature of $\FOL$ includes a set $\Func_{\Lan}$ of functional symbols of arity greater than or equal to $1$, a set $\Cons_{\Lan}$ of individual constants which are treated as 0-ary functional symbols.

We continue to call the terms of $\FOL$ $\Lan$-\textit{\textbf{terms}}.

Any $\FOL$ must have a nonempty set $\mathcal{P}$ of predicate symbols; some have equality $\approx$, some don't. We will be assuming that the arity of each predicate symbol $P\in\mathcal{P}$ is greater or equal to 1.

The notion of a formula is similar to the one given in Section~\ref{section:term-algebra} with a change in the clause for atomic formulas:
if $P\in\mathcal{P}$ of arity $n$ and $t_1,\ldots,t_n$ are $\Lan$-terms,
then $P(t_1,\ldots,t_n)$ is an (\textit{atomic}) $\FOL$-\textit{formula}.
The other clauses are the same as (1)--(6) on p.~\pageref{FOL-forfula}. 
Thus defined formulas we call $\FOL$-\textit{\textbf{formulas}}, or simply \textit{\textbf{formulas}}; the latter is used in the contexts where confusion is unlikely.

If an $\FOL$-formula is closed, it is called an $\FOL$-\textit{\textbf{sentence}}.\index{sentence} As in Section~\ref{section:term-algebra}, a formula which does not contain quantifiers is called \textit{\textbf{quantifier-free}}.

Models of $\FOL$ ($\FOL$-\textit{\textbf{models}} for short, or simply \textit{\textbf{models}})\index{model} which we use in Chapters 4--8 of this book, where we apply some of the results of model theory, have one of the following two forms:
\begin{equation}\label{E:types-of-models}
	\begin{array}{lcl}
		(\text{a})\quad\langle\textsf{A};\Func_{\Lan},\Cons_{\Lan},D\rangle &\text{or}
		&(\text{b})\quad\langle\textsf{A};\Func_{\Lan},\Cons_{\Lan},\varDelta\rangle,
	\end{array}
\end{equation}
where $\Func_{\Lan}$ and $\Cons_{\Lan}$ are the sets of interpreted functional symbols and constants  of $\FOL$, $D$ is a unary predicate and $\varDelta$, the interpretation of $\approx$ in $\fA$, is the diagonal of $\textsf{A}\times\textsf{A}$;  usually, $\langle\textsf{A};\Func_{\Lan},\Cons_{\Lan}\rangle$ is called an algebra of the corresponding model. If \mat{M} is a model, its algebra is denoted by $\mathcal{A}[\mat{M}]$. For a class $\mathcal{M}$ of models, $\mathcal{A}[\mathcal{M}]$ denotes the class of their algebras. Given a model $\mat{M}$, we denote by $|\mat{M}|$ its carrier. 

The models of type (\ref{E:types-of-models}--b) we call $e$-\textit{\textbf{models}}. 

Further, let $v$ be a map of the set of individual variables into $|\mat{M}|$. We will be using (a well-known) satisfactory  relation of a first-order formula
$A$ in a model $\mat{M}$, saying that $\mat{M}$ \emph{satisfies} $A$ \emph{at} $v$, symbolically $(\mat{M}, v)\models A$ or $\mat{M}\models_{v} A$.\footnote{Other utterances of this relation are:
	`$v$ satisfies $\varphi$ in $\mat{M}$' (\cite{bell-slomson1969}, chapter 3, {\S} 2); `$\varphi$ is true of $v$ in $\mat{M}$' (\cite{hodges1993}, section 1.3); `$\mat{M}$ satisfies $\varphi$ with $v$' (\cite{enderton2001}, section 2.2); `$(\mat{M},\sigma)$ satisfies $\varphi$' or `$\mat{M}$ is a model for $\varphi$ with respect to $v$' (\cite{rautenber2010}, section 2.3); `$\varphi$ is satisfied by a sequence $v:\Var_{\Lan}\longrightarrow\textsf{A}$ in $\mat{M}$' or
	`a sequence $v:\Var_{\Lan}\longrightarrow\textsf{A}$ satisfies $\varphi$ in $\mat{M}$ (\cite{chang-keisler1990}, definition 1.3.12).} A model $\mat{M}$ \emph{\textbf{satisfies}} $\varphi$, or $\varphi$ \emph{\textbf{is satisfiable in}} $\mat{M}$, if there is such a \emph{valuation} $v$ of the individual variables in $\mat{M}$ that $(\mat{M},v)\models\varphi$.

We write $\mat{M}\models A$ if $\mat{M}\models_{v} A$, for any valuation $v$ in \mat{M}, and say that $\mat{M}$ \textit{\textbf{validates}} $A$. Given a model $\mat{M}$ and a quantifier-free formula $A$, the following equivalence must be clear:
\begin{equation}\label{E:equivalence-for-modeling}
	\mat{M}\models\Forall\ldots\Forall A~\Longleftrightarrow~\mat{M}\models A.
\end{equation}

It must be clear that any identity, say $\Forall\ldots\Forall\,(\alpha\approx\beta)$, and the quasi-identity
$\Forall\ldots\Forall\,(p\approx p\Rarrow\alpha\approx\beta)$ are validated by the same $\FOL$-models.

\begin{defn}
	Let $\Sigma$ be a set of {\em$\FOL$}-sentences and let $\mathcal{M}$ be a nonempty class of {\em$\FOL$}-models. We say that $\Sigma$ is \textbf{satisfiable} in
	$\mathcal{M}$ if there a model $\fA\in\mathcal{M}$, which satisfies all sentences of $\Sigma$.
\end{defn}

Let $\mathcal{M}$ be a nonempty class of $e$-models. The set of equalities validated by all models from $\mathcal{M}$ is called the \textit{\textbf{equational theory}}\index{equational theory} of $\mathcal{M}$; this equational theory will be denoted by $\EQ[\mathcal{M}]$.\index{$\EQ[\mathcal{M}]$}
Now, we define:
\[
\text{Mod}_{\EQ}(\mathcal{M}):=\set{\mat{M}}{\mat{M}\models\EQ[\mathcal{M}]}.
\]

It is clear that $\mathcal{A}[\text{Mod}_{\EQ}(\mathcal{M})]$ (of algebras) is an equational class (or variety) generated by the class $\mathcal{A}[\mathcal{M}]$.
\begin{prop}[\cite{burris-sankapp1981},~\cite{gratzer2008},~\cite{mckenzie-mcnulty-taylor2018}]\label{P:equational-theory}
	Let $\mathcal{M}$ be a nonempty class of models of the form {\em(\ref{E:types-of-models}--b)}.	Then {\em$\EQ[\mathcal{M}]=\EQ[\text{Mod}_{\EQ}(\mathcal{M})]$}. Furthermore,
	{\em$\mathcal{A}[\text{Mod}_{\EQ}(\mathcal{M})]=\Ho\Su\Pro\mathcal{A}[\mathcal{M}]=\Ho\Su\Pro\Si\mathcal{A}[\mathcal{M}]$}.
\end{prop}
\begin{rem}
	{\em The last conclusion of Proposition~\ref{P:equational-theory} follows from 
		Proposition~\ref{P:Birkhoff's-criterion-for varieties}, Proposition~\ref{P:class-operators-inequalities} and Proposition~\ref{P:variety-generated-si-algebras}.}
\end{rem}

Let $\mat{M}=\langle\textsf{A};\Func_{\Lan},\Cons_{\Lan},\varDelta\rangle$ be a model of type (\ref{E:types-of-models}--b). We define:
\[
\lbrace\mat{M}\rbrace_{\textit{si}}:=\set{\langle\textsf{B};\Func_{\Lan},\Cons_{\Lan},\varDelta\rangle}{\langle\textsf{B};\Func_{\Lan},\Cons_{\Lan}\rangle\in\mathcal{A}[\mat{M}]_{\textit{si}}}
\]
\begin{prop}\label{P:EQ[M]=EQ[{M}_si]}
	Given a model {\em\mat{M}} of type {\em(\ref{E:types-of-models}--b)}, 
	{\em$\EQ[\mat{M}]=\EQ[\lbrace\mat{M}\rbrace_{\textit{si}}]$}.
\end{prop}
\begin{proof}
	We prove that for any identity $\Forall\ldots\Forall\epsilon$,
	\[
	\mat{M}\models\Forall\ldots\Forall\epsilon~\Longleftrightarrow~\lbrace\mat{M}\rbrace_{\textit{si}}
	\models\Forall\ldots\Forall\epsilon.
	\]	
	
	Indeed, suppose $\mat{M}\models\Forall\ldots\Forall\epsilon$ and $\mat{N}\in\lbrace\mat{M}\rbrace_{\textit{si}}$. Since $\mathcal{A}[\mat{N}]$ is a homomorphic image of $\mathcal{A}[\mat{M}]$, in virtue of Proposition~\ref{P:preservation}--a, $\mat{N}\models\Forall\ldots\Forall\epsilon$.
	
	Conversely, assume that $\lbrace\mat{M}\rbrace_{\textit{si}}
	\models\Forall\ldots\Forall\epsilon$. By definition, there is a set $\lbrace\mat{M}_i\rbrace_{i\in I}\subseteq\lbrace\mat{M}\rbrace_{\textit{si}}$ 
	such that the algebra $\mathcal{A}[\mat{M}]$ is subdirectly embedded into
	$\prod_{i\in I}\mathcal{A}[\mat{M}_i]$. According to Proposition~\ref{P:preservation}--b, $\langle\prod_{i\in I}\mathcal{A}[\mat{M}_i],\varDelta\rangle\models\Forall\ldots\Forall\epsilon$ and hence, in virtue of Proposition~\ref{P:preservation}--c, $\mat{M}\models\Forall\ldots\Forall\epsilon$.
\end{proof}

In the sequel, we will use the following fact.
\begin{prop}\label{P:varieties-Boolean-algebras}
	There are only two varieties of Boolean algebras --- the trivial variety and the variety of all Boolean algebras.	
\end{prop}
\begin{proof}
	Let us consider a nontrivial variety. And let us take any nontrivial algebra in it. According to Proposition~\ref{P:2-element-algebra}, a 2-element algebra is a subalgebra of the latter. Hence, according to Proposition~\ref{P:Birkhoff's-criterion-for varieties}, this 2-element algebra belongs to the variety. By the same proposition, any direct power of the 2-element algebra belongs to the variety. However, any Boolean algebra is isomorphic to some direct power of a 2-element Boolean algebra (\cite{burris-sankapp1981}, corollary IV.1.12). Therefore, any Boolean algebra belongs to any nontrivial variety.
\end{proof}

Let $\lbrace\mat{M}_i\rbrace_{i\in I}$ be a family of $\FOL$-models. Also, let $\nabla$ be a filter over $I$. We fix the direct product
\[
\mat{M}:=\prod_{i\in I}\mat{M}_i
\]
and define:
\[
x\equiv_{\nabla}y\stackrel{\text{df}}{\Longleftrightarrow}\set{i\in I}{x_i=y_i}\in\nabla.
\]

It is well known that $\equiv_{\nabla}$ is congruence relation on the algebra $\alg{A}=\langle\textsf{A},\Func,\Cons\rangle$. We denote by $\alg{A}\slash\nabla$ the quotient algebra with respect to this congruence. To obtain a model, $\fA\slash\nabla$, we have to expand $\alg{A}\slash\nabla$ with an interpretation of $\approx$ in this algebra. Let us denote (temporarily) this interpretation by $\varDelta_D$. By definition, we have:
\[
(x\slash \nabla,y\slash\nabla)\in\varDelta_\nabla\stackrel{\text{df}}{\Longleftrightarrow}
\set{i\in I}{(x_i,y_i)\in\varDelta_i}\in\nabla.
\]

Since each $\varDelta_i$ is a diagonal, we observe that
\[
(x\slash\nabla,y\slash\nabla)\in\varDelta_\nabla\Longleftrightarrow x\slash \nabla=y\slash\nabla.
\]
This means that $\varDelta_\nabla$ is the diagonal of $|\fA\slash F|\times|\fA\slash F|$; we denote it simply by $\varDelta$, if no confusion arises.

We denote
\[
\prod_{\nabla}\mat{M}_i:=\langle\alg{A}\slash\nabla,\varDelta\rangle
\]
and call the last model a \textit{\textbf{reduced product}}\index{reduced product} of the family $\lbrace\fA_i\rbrace_{i\in I}$ with respect to the filter $\nabla$. If $\nabla$ be an ultrafilter over $I$, we call $\prod_{\nabla}\mat{M}_i$ an \textit{\textbf{ultraproduct}}\index{ultraproduct} with respect to the ultrafilter $\nabla$.

In case if each $\mat{M}_i=\langle\alg{A}_i,D_i\rangle$, where $D_i$ is a unary predicate on $|\fA_i|$, then, by definition, we have:
\begin{equation}\label{E:ultraproduct-df}
	\prod_{\nabla}\mat{M}_i:=\langle\alg{A}\slash\nabla,D\rangle,
\end{equation}
where
\[
x\slash\nabla\in D\stackrel{\text{df}}{\Longleftrightarrow}\set{i\in I}{x_i\in D_i}\in\nabla.
\]

Let $\mathcal{M}$ be a nonempty class of models of the form (\ref{E:types-of-models}--b). The set of quasi-identities validated by all models from $\mathcal{M}$ is called the \textit{\textbf{quasi-equational theory}} of $\mathcal{M}$; this quasi-equational theory will be denoted by $\QE[\mathcal{M}]$.
Now, we define:
\[
\text{Mod}_{\QE}(\mathcal{M}):=\set{\mat{M}}{\mat{M}\models\QE[\mathcal{M}]}.
\]

The class $\mathcal{A}[\text{Mod}_{\QE}(\mathcal{M})]$ (of algebras) is obviously an implicative class (or a quasi-variety) generated by the class $\mathcal{A}[\mathcal{M}]$.
\begin{prop}[\cite{burris-sankapp1981}, theorem V.2.25, \cite{gorbunov1998}, corollary 2.3.4]\label{P:quasi-equational-theory}
	Let $\mathcal{M}$ be a nonempty class of $e$-models.	Then {\em$\QE[\mathcal{M}]=\QE[\text{Mod}_{\QE}(\mathcal{M})]$}. Furthermore,
	{\em$\mathcal{A}[\text{Mod}_{\QE}(\mathcal{M})]=\Su\Pro\Pu\mathcal{A}[\mathcal{M}]$}.
\end{prop}
\begin{rem}
	{\em The last conclusion of Proposition~\ref{P:quasi-equational-theory} follows from Proposition~\ref{P:Mal'cev's-for-quasivarieties} and Proposition~\ref{P:class-operators-inequalities}.}
\end{rem}

In the sequel, we will need the following refined variant of compactness theorem.
\begin{prop}[\cite{bell-slomson1969},~\cite{chang-keisler1990}]\label{P:compactness-thm-refined}
	Let $\Sigma$ be a set of {\em$\FOL$}-sentences and let $\mathcal{M}$ be a nonempty class of {\em$\FOL$}-models closed under ultraproducts. Then $\Sigma$ is satisfiable in $\mathcal{M}$ if, and only if, every finite subset of $\Sigma$ is satisfiable in $\mathcal{M}$.
\end{prop}
\begin{proof}
	Suppose $\Sigma$ is satisfiable in $\mathcal{M}$; that is, there is a model $\mat{M}\in\mathcal{M}$, which satisfies all sentences of $\Sigma$. Then, obviously, $\fA$ satisfies any subset, in particular any finite subset, of $\Sigma$.
	
	Conversely, let $I$ denote the set of all finite subsets of $\Sigma$. By premise, for each $i\in I$, there is a model $\mat{M}_i$ which satisfies $i$.
	According to Corollary 4.1.11 of~\cite{chang-keisler1990}, there is an ultrafilter $\nabla$ over $I$ such that the ultraproduct $\prod_{\nabla}\mat{M}_i$ satisfies $\Sigma$. By premise, this model belongs to $\mathcal{M}$. Hence, $\Sigma$ is satisfiable in $\mathcal{M}$.
\end{proof}

A first-order formula $A$ is said to be \textit{\textbf{preserved}}:
\begin{itemize}
	\item \textit{\textbf{under homomorphisms}} if for any homomorphism $f$ of a model $\mat{M}$ into a model $\mat{N}$ and any valuation $v$ in $\mat{M}$,
	\[
	(\mat{M},v)\models\phi\Longrightarrow (\mat{N},v_f)\models\phi,
	\]
	where $v_f$ is a valuation in $\mat{N}$ such that $v_f[p]=f(v[p])$, for any individual variable $p$;
	\item \textit{\textbf{under direct products}} if for any class $\lbrace\mat{M}_i\rbrace_{i\in I}$ of $\FOL$-models and any valuation $v$ in
	$\mat{M}=\prod_{i\in I}\mat{M}_i$, whenever $(\mat{M}_i,v_i)\models\phi$, for every $i\in I$, where $v_i[p]$ is the projection of $v[p]$ on $\mat{M}_i$, then $(\mat{M},v)\models\phi$,	
\end{itemize}
\begin{prop}\label{P:preservation}
	{\em\[
		\begin{array}{cl}
			(\text{a}) &\textit{Any equality and any identity is preserved under homomorphisms}.\\
			(\text{b}) &\textit{Any equality and any identity, as well as any quasi-equality and any quasi-identity},\\ 
			&\text{are preserved under direct products}.\\
			(\text{c}) &\text{Both identities and quasi-idenities are preserved on subalgebras}.
		\end{array}
		\]}
\end{prop}
\begin{proof}
	(a) is part of the proof of theorem 3.2.4 of~\cite{chang-keisler1990}. (b) is covered by theorem 1 of~\cite{malcev1973}, chapter III, {\S} 7. For the proof of (c), see, e.g., \cite{malcev1973}, {\S} 7, theorem 1.
\end{proof}

\paragraph{References}
\begin{enumerate}
	\item~\cite{bell-slomson1969}
	\item~\cite{burris-sankapp1981}
	\item~\cite{chang-keisler1990}
	\item~\cite{enderton2001}
	\item~\cite{gorbunov1998}
	\item~\cite{gratzer2008}
	\item~\cite{hodges1993}
	\item~\cite{malcev1973}
	\item~\cite{mckenzie-mcnulty-taylor2018}
	\item~\cite{rautenber2010}	
\end{enumerate}

\section{Preliminaries from computability theory}\label{section:prelimanaries-computability}
In this section, we outline some ideas of computability and strive to prepare a reader unfamiliar with these ideas for reading Chapter~\ref{chapter:decidability}. Thus we intend to introduce the reader to the idea of an \emph{effective method} or \emph{algorithm}. From the very beginning, we need to make a few warnings.

First of all, we do not assume any aesthetic load in the term `effective method'. Secondly, although different levels of effectiveness can be distinguished in computability theory, in this book we deal mainly with the lowest. Even at this level, we will prove that certain problems have positive solutions in the sense that there is a desirable effective method or algorithm that \textit{\textbf{decides}} this problem. Such problems are called \textit{\textbf{decidable}}. The problems (from a certain category of problems) for which any algorithm does not produce a desirable result are called \textit{undecidable}.\index{problem!decidable}\index{problem!undecidable} This distinction requires a rigorous definition of the concept of an effective method. In this book, we however deal only with decidable problems. Thus, the purpose of this section is to discuss the concept of decidability from the point of view of an effective method without going into rigorous detail, that is, without giving an exact explication of an effective method as a computable function. The reader can find a precise definition of a computable function, decidability, and undecidability in the references at the end of this section.\\

Suppose a set $A$ is a subset of a set $U$. It is customary to think of $U$ as a `universal' set. Then the problem arises: Is there an effective method for any $x$ in $U$ to decide whether $x$ lies in $A$? For simplicity, let us assume that $U=\mathbb{N}$, where $\mathbb{N}$ is the set of nonnegative integers introduced in Section~\ref{section:perlimenaries-set-theory}
We would answer the above question in the affirmative, if we could find an effective function $\chi_A:\mathbb{N}\longrightarrow\lbrace 0,1\rbrace$ such that
\begin{equation}\label{E:characteristic-of-A}
	\chi_A(x)=\begin{cases}
		\begin{array}{cl}
			1 &\text{if $x\in A$}\\
			0 &\text{if $x\notin A$}.
		\end{array}
	\end{cases}
\end{equation}

The function $\chi_A$ is called the \textit{\textbf{characteristic function}}\index{characteristic function} of $A$. If we succeed in finding an effective (or algorithmic) $\chi_A$, we state that
\[
A=\set{x\in\mathbb{N}}{\chi_A(x)=1}
\]
and, hence, can claim that the set $A$ is decidable\index{decidable!set}.

The situations that we will meet in Chapter~\ref{chapter:decidability} are a little more complex. Roughly speaking, we will consider situations when $U\subset\mathbb{N}$ and is decidable and then ask the question whether a set $A\subseteq U$ is \textit{\textbf{decidable relative to}} $U$, that is, whether an algorithm exists for any $x\in U$ to decide if $x\in A$.

Let us define for the last decidability procedure a function $\chi_{A|U}:U\longrightarrow\lbrace 0,1\rbrace$:
\begin{equation}\label{E:characteristic-of-A|U}
	\chi_{A|U}(x):=\begin{cases}
		\begin{array}{cl}
			1 &\text{if $x\in A$}\\
			0 &\text{if $x\in U\setminus A$}.
		\end{array}
	\end{cases}
\end{equation}

The function $\chi_{A|U}$ is partial, since it is undefined on $\mathbb{N}\setminus U$. We agree that effective method can be implemented as a partial function, and later in Section~\ref{section:algorithmic-functions-on-N}, we will give the rationale for this view. However, even in this case, the set $A$ is decidable as soon as the characteristic functions $\chi_U$ and $\chi_{A|U}$ are realized as algorithms.

Indeed, we define a function $\chi^{\prime}_{A}:\mathbb{N}\longrightarrow\lbrace 0,1\rbrace$ as follows:
\[
\chi^{\prime}_{A}(x):=\begin{cases}
	\begin{array}{cl}
		1 &\text{if $\chi_U(x)=1$ and $\chi_{A|U}(x)=1$}\\
		0 &\text{if $\chi_U(x)=1$ and $\chi_{A|U}(x)=0$}\\
		0 &\text{if $\chi_U(x)=0$}.
	\end{array}
\end{cases}
\]

It is clear that for any $x\in\mathbb{N}$,
\[
x\in A~\Longleftrightarrow~\chi^{\prime}_{A}(x)=1.
\]	

Next in this section, we discuss the idea of ​​an effective method or algorithm that is implemented by a function, possibly partial, from $\underbrace{\mathbb{N}\times\ldots\times\mathbb{N}}_{n}$, where $n\ge 1$, or its proper subset into $\mathbb{N}$. Along with such functions, we also consider relations $R\subseteq\underbrace{\mathbb{N}\times\ldots\times\mathbb{N}}_{n}$ such that their characteristic function $\chi_R:\underbrace{\mathbb{N}\times\ldots\times\mathbb{N}}_{n}\longrightarrow\lbrace 0,1\rbrace$ is effective. We call such functions and relations \textit{\textbf{effective}} or \textit{\textbf{algorithmic}} if they can be defined in `a purely mechanical way'. (This phrase will be discussed in the following subsection.) After that, we discuss the idea of algorithmic functions and relations that are defined on other domains.

\subsection{Algorithmic functions and relations on $\mathbb{N}$}\label{section:algorithmic-functions-on-N}
It is not surprising that we start with functions of any finite number of variables defined on $\mathbb{N}$, taking into account the historical fact, which E. Beth expressed in words:
\begin{quote}
	``At any rate, the fascination inspired by the art of computation explains why, at a certain moment, a properly scientific interest in numbers developed.'' (\cite{beth1962}, section 25)
\end{quote}

What functions can be considered executable in accordance with an effective method? The last example of the last subsection suggests that to avoid the process of relative effectiveness \emph{ad infinitum}, it must be agreed that some functions are effective for granted. By agreement, the following functions form the \emph{basis} of effective functions:
\begin{itemize}
	\item the \emph{successor function} $\mathbf{s}:n\mapsto n+1$, where $n\in\mathbb{N}$; 
	\item the \emph{zero functions}: $\mathbf{o}^{n}:x\mapsto 0$, where $x\in\underbrace{\mathbb{N}\times\ldots\times\mathbb{N}}_{n}$ and $n\in\mathbb{N}$; in particular, $\mathbf{o}^{0}$ is just $0$;
	\item and the \emph{projections}: $\pi^{n}_{m}(x_1,\ldots,x_n)=x_m$, where each $x_i\in\mathbb{N}$ and $1\le m\le n$.
\end{itemize}

We argue that the basis functions can be performed in `a purely mechanical way'. Before moving on to this, we need consensus. The consensus is that any natural numbers $m$ and $n$ can be identified `purely mechanically' as different or the same. This needs to be explained.

To facilitate our discussion of the `purely mechanical method', let us assume that there is a device, a \emph{computer}, which is capable of performing such `mechanical' actions. Our first assumption about this device is that it takes a finite number of constructive objects as input and returns the constructive object as output. This is what V. Uspensky and A. Semenov write about the notion of a \emph{constructive object}:
\begin{quote}
	``This notion can be regarded not only as a basic notion but also as a primary one. Strictly speaking it belongs not to the theory of algorithms but to the introduction to it. As one old textbook says: ``A subject of an introduction to a science usually contains \emph{preliminary notions}, i.e. notions that can not be included into the science itself but are essential for it and are assumed by it''.
	
	All attempts to define the notion of a constructive object (including ours) either are very abstract and vague or define specific classes of constrictive objects. As all primitive notions it can be understood by examples.'' (\cite{uspensky-semenov1993}, section 1.0.0)
\end{quote}

Perhaps the main feature of the constructive object is that it can be represented by a word in a certain alphabet. For instance, using the decimal symbols '$0$', `$1$', \ldots, '$9$', \emph{digits}, any natural number can be written as a finite sequence of digits, providing that only one such a sequence starts with `$0$', that is the one which does not have more digits. The examples of such sequences are:
\[
0, ~1, ~2, ~31, ~502.
\]

Another way to depict natural numbers as words is the following. Any positive number is written as a sequence of occurrences of `$|$', and $0$ is depicted as an empty word, for which one can use any symbol different from `$|$', for example, `$\emptyword$'. Thus we obtain that $0$ is depicted by `$\emptyword$', $1$ by `$|$', $2$ by `$||$', etc. 

More general setting of depicting constructive objects will be discussed in the next subsection. 

Now it should be clear that the computer, using the two numbers $m$ and $n$ in one of the representations presented above, can decide `purely mechanically', that is, without resorting to intelligence in any reasonable sense of the word, whether they are represented by the same string of elementary characters.

We note that the equality relation `$=$' defined for the output of the same or different algorithms is a more general relation than the equality of two strings of symbols representing natural numbers in one of the above notations.
As we will see below, the relation $=$ between two algorithmic functions is also decidable in a `purely mechanical way'.

It is generally accepted that in order to be qualified as algorithmic, a method must have the following characteristics:
\begin{itemize}
	\item the algorithmic method is \textit{deterministic}; in other words, when implemented, the algorithm works like a function;
	\item the algorithmic method is independent of a single input, but can be applied to a data set;
	\item as a function, the algorithm may be partial; but if it is defined for input, it completes its work with a certain result.
\end{itemize}

A. A. Markov claims something similar as follows:
\begin{quote}
	``The following three features are characteristic of algorithms and determine their role in mathematics;\\
	a) the precision of the prescription, leaving no place to arbitrariness, and its universal comprehensibility --- the definiteness of the algorithm;\\
	b) the possibility of starting out with initial data, which may vary within given limits --- the generality of the algorithm;\\
	c) the orientation of the algorithm toward obtaining some desired result, which is indeed obtained in the end with proper initial data, --- the conclusiveness of the algorithm.'' (\cite{markov1962}, introduction; cp.~\cite{markov1960}, section 1.)
\end{quote}

Returning to the functions of the basis, we argue that the procedure for obtaining for a given natural number $n$ the number $\mathbf{s}(n)$ in any of the aforementioned representations can be regarded as purely mechanical, as well as the procedure for returning $0$ from any natural argument or from nothing. Also, the projections $\pi^{n}_{m}:\underbrace{\mathbb{N}\times\ldots\times\mathbb{N}}_{n}\longrightarrow\mathbb{N}$ can be performed purely mechanically, given that the parameters $m$ and $n$ are predefined.

Thus in this context, we can understand $\mathbb{N}$ as an effectively generated infinite list of \emph{numerals}:
\begin{equation}\label{E:sequence-numerals}
	0, ~\mathbf{s}(0),~\mathbf{s}(\mathbf{s}(0)),\ldots;
\end{equation}

We observe that, given a numeral $n$ of~\eqref{E:sequence-numerals}, which is not coincident with $0$, the computer can determine the (nearest) predecessor of $n$, that is, such a numeral $m$ that $\mathbf{s}(m)$ is coincident with $n$. This leads to the \textit{recursive definitions}\footnote{See about recursive definitions (or definitions by recursion) in~\cite{cutland1980}, section 2.4,~\cite{davis1958}, section 3.3, \cite{davis-sigal-weyuker1994}, section 3.2, or~\cite{malcev1970}, section 2.2} of effective functions. 

Recursion allows you to define a new function using one or two given functions by recursion equations. There are two types of recursion equations. The first type uses one given function, say $f$, and defines a new function of one variable, say $h$, by the following equations:
\[
\begin{cases}
	\begin{array}{l}
		h(0)=a\\
		h(\mathbf{s}(x))=f(x,h(x)),
	\end{array}
\end{cases}
\]
where $a$ is a given number.

The second type uses two given functions, say $f(x_1,\ldots,x_{n-1})$ and $g(x_1,\ldots,x_{n+1})$, and defines a function of $n$ variables with $n\ge 2$:
\[
\begin{cases}
	\begin{array}{l}
		h(x_1,\ldots,x_{n-1},0)=f(x_1,\ldots,x_{n-1})\\
		h(x_1,\ldots,x_{n-1},\mathbf{s}(x))=g(x_1,\ldots,x_{n-1},x,g(x_1,\ldots,x_{n-1},x))
	\end{array}
\end{cases}
\]

Using basis functions as initial functions, many other functions can be defined in a finite number of steps, using \textit{substitution} (or \textit{composition}) and suitable recursion equations. The functions obtained from the basis functions by a finite number of applications of substitution and recursion equations form the class of \textit{primitive recursive} functions.

There is a strong argument that a recursive definition whose recursion equations employ effective functions defines an effective function; cf.~\cite{cutland1980}, section 4.4,~\cite{davis-sigal-weyuker1994}, section 3.3,~\cite{davis1958}, section 3.1, or~\cite{malcev1970}, section 12.3.
In addition, this argument states that the basis functions are effective and, therefore, all primitive recursive functions are effective or algorithmic.
The argument, known as the \emph{Church thesis}, goes even further, claiming that a certain extension of the class of primitive recursive functions, among which some functions are partial, contains exactly all effective functions.\footnote{About the Church thesis, see, e.g.,~\cite{uspensky-semenov1993}, section 1.2.4, and references therein.}
In Section~\ref{section:enumerable-sets} we give an example and demonstrate the usefulness of partial effective functions.

The inclusion of partial functions is important. For illustration, let us consider two computable functions, $f(x)$ and $g(x)$, where $f(x)$ is partial and $g(x)$ is total. If however $f(x)$ is defined on all outputs of $g(x)$, then the composite function $f\circ g(x)$, obtained by substitution, $f(g(x))$, is a total computable function. \\

Below we demonstrate that some functions and predicates defined on $\mathbb{N}$ and on Cartesian products of $\mathbb{N}$ are effective, providing that the Church thesis is accepted.

For each function below, we first give a meaningful definition of the function, and then define it using recursion equations.
\[
\text{sg}(x):=\begin{cases}
	\begin{array}{cl}
		0 &\text{if $x$ is coincident with $0$}\\
		1 &\text{if otherwise};
	\end{array}
\end{cases}
\]

These are recursion equations for \text{sg}.
\[
\begin{cases}
	\begin{array}{l}
		\text{sg}(0)=\mathbf{o}^{0}\\
		\text{sg}(\mathbf{s}(x))=\mathbf{s}(\mathbf{o}^{1}(x)).
	\end{array}
\end{cases}
\]
\[
\text{pd}(x):=\begin{cases}
	\begin{array}{cl}
		0 &\text{if $x$ is coincident with $0$}\\
		\text{the predicessor of $x$} &\text{if otherwise}.
	\end{array}
\end{cases}
\]
\[
\begin{cases}
	\begin{array}{l}
		\text{pd}(0)=\mathbf{o}^{0}\\
		\text{pd}(\mathbf{s}(x))=\pi^{1}_{1}(x).
	\end{array}
\end{cases}
\]

Using the last function, we define:
\[
x\dotminus y:=\begin{cases}
	\begin{array}{cl}
		x &\text{if $y$ is coincident with $0$}\\
		\text{pd}(x\dotminus\text{pd}(y)) &\text{if otherwise};
	\end{array}
\end{cases}
\]
and now this function by recursion equations:
\[
\begin{cases}
	\begin{array}{l}
		x\dotminus 0=\pi_{1}^{1}(x)\\
		x\dotminus\mathbf{s}(y)=\text{pd}(x\dotminus y).
	\end{array}
\end{cases}
\]

Next we define for the numerals~\eqref{E:sequence-numerals} the following relations:
\[
x\le y~\define~1\dotminus\text{sg}(x\dotminus y)~\text{is coincident with}~1
\]
and
\[
x=y~\define~x\le y~\text{and}~y\le x.
\]

Taking this into account, we observe that the function $1\dotminus\text{sg}(x\dotminus y)$ is the characteristic function of $\le$ and, hence, $\le$ is a decidable binary relation on $\mathbb{N}$.

Let us temporarily denote the characteristic function of $\le$ by $\chi_1$; that is,
\[
\chi_1(x,y):=1\dotminus\text{sg}(x\dotminus y).
\]

Now we define:
\[
\chi_2(x,y):=\chi_1(x,y)\cdot\chi_1(y,x).
\]

It is obvious that $\chi_2$ is the characteristic function of $=$. We aim to show that $\chi_2$ is effective and, hence, $=$ is decidable.

Indeed, through two recursive definitions, we show that the ordinary operation
of addition, $x+y$, and that of multiplication, $x\cdot y$, are effective.

Thus we first define:
\[
x+y:=\begin{cases}
	\begin{array}{cl}
		x &\text{if $y=0$}\\
		\mathbf{s}(x+\text{pd}(y)) &\text{if otherwise}.
	\end{array}
\end{cases}
\]
This function can be given by the following recursion equations:
\[
\begin{cases}
	\begin{array}{l}
		x+0=\pi_{1}^{1}(x)\\
		x+\mathbf{s}(y)=\mathbf{s}(x+y).
	\end{array}
\end{cases}
\]
Finally, we define:
\[
x\cdot y:=\begin{cases}
	\begin{array}{cl}
		0 &\text{if $y=0$}\\
		(x\cdot\text{pd}(y))+x &\text{if otherwise}.
	\end{array}
\end{cases}
\]

The last function is given by the following recursion equations:
\[
\begin{cases}
	\begin{array}{l}
		x\cdot 0=\mathbf{o}^{1}(x)\\
		x\cdot\mathbf{s}(y)=(x\cdot y)+x.
	\end{array}
\end{cases}
\]

Recall that a (two-valued) \textit{predicate} is a function with a range of two values: \textit{true} and \textit{false}. We will be considering predicates that are defined on some universal set $U$ which will vary. A predicate $P(x)$ is called \textit{\textbf{decidable}} (relative to $U$) if the set
\[
A_P:=\set{x\in U}{P(x)~\text{is true}}
\]
is decidable.

As a consequence of the last definition we obtain that,
given effective functions $f$ of $m$ variables and $g$ of $k$ variables, the relation `$f(x)=g(y)$', where $x\in\underbrace{\mathbb{N}\times\ldots\times\mathbb{N}}_{m}$ and $y\in\underbrace{\mathbb{N}\times\ldots\times\mathbb{N}}_{k}$, is decidable; in particular, the equality relation `$x=y$' is decidable. It is not difficult to see that the nonequality relation `$x\neq y$' is also decidable.

This allows us not to distinguish, at least in theoretical considerations, the numerals of~\eqref{E:sequence-numerals} and the natural numbers of $\mathbb{N}$, and also not to distinguish the relations $\le$ and $=$ just defined respectively from the ordinary `less than or equal' relation and equality relation on $\mathbb{N}$.

In particular, this implies that any given nonempty finite set $A=\lbrace n_1,\ldots,n_k\rbrace$ of natural numbers is decidable; the characteristic function of $A$ can be defined, e.g., as follows:
\[
\chi_A(x):=\text{sg}(\chi_2(x,n_1)+\ldots+\chi_2(x,n_k)).
\]
The empty set is also decidable, since it can be defined by a decidable inconsistent condition; for instance,
\[
\emptyset=\set{x}{x\neq x}.
\]

Below we will also need the function $x^{y}$ which can be defined by the following recursion equations:
\[
\begin{cases}
	\begin{array}{l}
		x^{0}=\mathbf{s}(\mathbf{o}^{1}(x))\\
		x^{\mathbf{s}(y)}=x^{y}\cdot x.
	\end{array}
\end{cases}
\]

The following facts are well known.
\begin{itemize}
	\item The predicate $\text{Dv}(x,y)$ (`$x$ divides $y$') as defined below is decidable:
	\[
	\text{Dv}(x,y)=\begin{cases}
		\begin{array}{cl}
			\text{true} &\text{if $x\neq 0$ and $x$ divides $y$}\\
			\text{false} &\text{if otherwise}.
		\end{array}
	\end{cases}
	\]
	\item Let us define a predicate $\text{Pr}(x)$ on $\mathbb{N}$ as follows:
	\[
	\text{Pr}(x)~\text{is true}~\define~\text{$x$ is a prime number}.
	\]
	The predicate $\text{Pr}(x)$ is decidable. Cf., e.g.,~\cite{cutland1980}, section 2.4,~\cite{davis1958}, section 3.5, or~\cite{davis-sigal-weyuker1994}, section 3.6. 
	\item The predicate $\text{Pr}(x)$ is used to show that a function defined on $\mathbb{N}$ that, given $x\in\mathbb{N}$, returns the $x$th prime, $p_x$,  is algorithmic. See, e.g.,~\cite{cutland1980}, ibid,~\cite{davis1958}, ibid, or~\cite{davis-sigal-weyuker1994}, ibid.
\end{itemize}

\subsection{Word sets and functions}\label{section:word-sets}
Algorithmic functions can be defined not only for natural numbers or, more precisely, for numerals but also for constructive objects of a more general nature.

The elements of~\eqref{E:sequence-numerals} are examples of constructive objects. Another example are words in a given alphabet.

Let $\Sigma$ be a nonempty set of \textit{primitive symbols}. Taking a few symbols from $\Sigma$, say $a$, $b$ and $c$, and adding them together in any order as a string, we get a \textit{\textbf{word}}\index{word} over the \textit{\textbf{alphabet}} $\Sigma$. The set of all words over $\Sigma$ is denoted by $\Sigma^{\ast}$. These are examples of words over $\Sigma=\lbrace a,b,c\rbrace$:
\[
a,~b,~c,~ab,~ba,~cb,~bc,~ac,~ca.
\]
For convenience, we assume that the empty word, denoting it as before by $\emptyword$, is also contained in $\Sigma^{\ast}$.

As for the primitive symbols of alphabets, we adhere to the requirement, which was expressed by the authors of~\cite{fraenkel-bar-hillel-levy1973} as follows:
\begin{quote}
	``[$\ldots$] being a primitive symbol has to be an effective notion, enabling us to determine in a finite number of steps whether or not a given symbol is primitive.'' (\cite{fraenkel-bar-hillel-levy1973}, chapter V, {\S} 2.)
\end{quote}

As before, we can designate two sets $A$ and $U$, such that $A\subseteq U\subseteq \Sigma^{\ast}$ and raise the question, whether $A$ is decidable relative to $U$.

Now, imagination can give more examples of constructive objects. For example, let $\Sigma^{\ast}$ be the set of finite graphs, $U$ the set of finite connected graphs, and $A$ the set of finite connected trees.

The last example shows that $\Sigma$ can be infinite. But to ask questions about decidability, $\Sigma$ must be countable, that is, we must assume that
$\card{\Sigma}\le\aleph_0$. However, when $\card{\Sigma}=\aleph_0$, many questions, especially those of algorithmic character, can be reduced to similar questions in the framework of a finite alphabet. For example, let us take the infinite alphabet $\Sigma=\lbrace a_1, a_2,\ldots\rbrace$. On the other hand, let us consider the alphabet $\Delta=\lbrace a,|\rbrace$. We argue that every word over $\Sigma$ has its counterpart in $\Delta^{\ast}$. The last claim becomes clear when we notice the following correspondence: $a_1$ corresponds to $a|$, $a_2$ to $a||$, and so on. Therefore, in this section below, we restrict ourselves to finite alphabets.

In Chapter~\ref{chapter:decidability}, we will use as constructive object sentential formulas, finite lists of such formulas, as well as finite logical matrices. Despite the fact that it is impossible to give a rigorous definition of a constructive object (see Section~\ref{section:algorithmic-functions-on-N}), yet, practice shows that the existing descriptions can be reduced to the second example; that is, all existing examples are possible to consider as words over a suitable alphabet.
But even this is not more general than our first example --- numerals, as we are going to demonstrate below.\\

Let $\Sigma=\lbrace a_1,\ldots,a_n\rbrace$ be a finite alphabet. We define a function $f:i\mapsto 2i+1$, where $1\le i\le n$. Further, we define a function $@:\Sigma^{\ast}\longrightarrow\mathbb{N}$ as follows: 
\[
@(\emptyword)=3~\,\text{and}~\,@(a_{i_1}\ldots a_{i_k})=p_{0}^{f(i_1)}\cdot\ldots\cdot p_{k-1}^{f(i_k)},
\]
where $a_{i_1}\ldots a_{i_k}$ is an arbitrary nonempty word of $\Sigma^{\ast}$.

It should be clear that the map $@$ is injective, but not surjective. The map @ is an example of a \emph{translation} of the words of one alphabet, in this example $\Sigma$, to the words of another alphabet, in this example $\lbrace 0, 1,\ldots,9\rbrace$.\footnote{About the notion of a translation see~\cite{markov1962}, chapter 1, {\S} 6.} In our example, the translation is algorithmic. As the reader can notice, there is a reverse translation. Although the fundamental theorem of arithmetic guarantees that the reverse translation is algorithmic, it is only a partial function, since the map @ is not surjective.\\

Let $A\subseteq\Sigma^{\ast}$. $A$ is called \textit{\textbf{decidable}} if the set $@(A)$ is decidable. 

An intermediate consequence of this definition is that the set $\Sigma^{\ast}$ is decidable. If we denote $U:=@(\Sigma^{\ast})$, then, solving the decidability problem about a given set $A\subset\Sigma^{\ast}$, we can reduce it to the question, whether the set $@(A)$ is decidable relative to $U$.

\subsection{Effectively enumerable sets}\label{section:enumerable-sets}
For simplicity, we limit ourselves in this subsection with effective functions defined on subsets of $\mathbb{N}$.

We encounter a partial effective functions when we deal with semidecidable sets. A set $A\subseteq\mathbb{N}$ is \textit{semidecidable}, or \textit{\textbf{effectively denumerable}},\index{effectively enumerable} if there is an effective function function $\varphi_A$ such that
\[
\varphi_A(x)=\begin{cases}
	\begin{array}{cl}
		1 &\text{if $x\in A$}\\
		\text{undefined} &\text{if $x\notin A$}.
	\end{array}
\end{cases}
\]

The ground for the first name is obvious, and for the second --- the following (well-known) theorem.

\begin{quote}
	{\em A set $A$ is effectively enumerable if, and only if, $A=\varnothing$ or $A$ is the range of a total effective function.} 
\end{quote}
(Cf.~\cite{cutland1980}, chapter 7, section 2,~\cite{davis1958}, chapter 5, section 4, or~\cite{davis-sigal-weyuker1994}, chapter 4, section 4.)

If a set is decidable, it is also semidecidable, but the converse is not necessarily true. It is usually simpler to show the semidecidability of a set than to prove that the set is decidable. The following theorem is very useful for establishing the decidability of sets.
\begin{quote}
	(Post) {\em A set $A$ is decidable if, and only if, both $A$ and $\mathbb{N}\setminus A$ are effectively enumerable.} 
\end{quote}
(Cf. ibid.)

The last theorem is based on the previous one, and can be illustrated as follows.

Suppose both $A$ and $\mathbb{N}\setminus A$ are effectively enumerable and unequal to $\varnothing$. For convenience, we denote: $B:=\mathbb{N}\setminus A$. In virtue of the first theorem, there are two effective total functions, $\psi_A$ and $\psi_B$, such that $A=\psi_A(\mathbb{N})$ and $B=\psi_B(\mathbb{N})$. Then it is clear that any given $n\in\mathbb{N}$ is among the members of the following sequence:
\[
\psi_A(0),~\psi_B(0),~\psi_A(1),~\psi_B(1),\ldots
\]
If $n$ is a value of $\psi_A$, then $n\in A$; if $n$ is a value of $\psi_B$ in this sequence, then $n\notin A$. The described procedure is effective, since it is based on the effective functions $\psi_A$ and $\psi_B$ and on the decidable relation $=$.\\

Another decidability procedure uses a different criterion for effectively enumerable sets.
\begin{quote}
	{\em A set $A$ is effectively enumerable if, and only if, there is a total effective function $\phi_A(x,y)$ (with two variables) such that for any $n\in\mathbb{N}$, $n\in A$ if, and only if, the equation $\phi_A(x,n)=0$ has a solution.}
\end{quote}
(Cf.~\cite{cutland1980}, ibid,~\cite{malcev1970}, section 4.2.)

Indeed, if given sets $A$ and $B$ are as above and if we know functions $\psi_A$ and $\phi_B$, then for a given $n\in\mathbb{N}$, we effectively create the sequence:
\begin{equation}\label{E:enumerable-sets-2}
	\psi_A(0),~\phi_B(0,n),~\psi_A(1),~\phi_B(1,n),\ldots
\end{equation}\label{E:enumerable-sets}
Thus, if $n$ is a value of $\psi_A$, then $n\in A$; if one of the values of $\phi_B$ in this sequence takes $0$, then $n\notin A$. It is clear that one of these options should happen.\\

The situation is not too complicated if we have a universal decidable set $U\subset\mathbb{N}$ and all effective procedures are executed relative to $U$.
It is also not more difficult when in Chapter~\ref{chapter:decidability} we will work with words in a fixed alphabet, and not with numbers, where a function such as @ will be implied, but behind the scenes. Namely, the sequence~\eqref{E:enumerable-sets-2} will be implemented there on word sets.

\paragraph{References}
\begin{enumerate}
	\item~\cite{cutland1980}
	\item~\cite{davis1958}
	\item~\cite{davis-sigal-weyuker1994}
	\item~\cite{malcev1970}
	\item~\cite{markov1960}
	\item~\cite{markov1962}
	\item~\cite{uspensky-semenov1993}
\end{enumerate}

\chapter[Sentential Formal Languages]{Sentential Formal Languages and Their Interpretation}\label{chapter:languages}	
	
\section{Sentential Formal languages}\label{section:languages}
It is commonly accepted that formal logic is an enterprise, in the broad sense of the word, which is characterized, first of all, by the use, study and development  of formal languages. One of the grounds for this characterization is that formal language ensures the codification of accepted patterns of reasoning. From the outset, we have to admit that such a representation of the train of thought in argumentation confines it to a framework of the discrete rather than the continuous. From this viewpoint, an explicit argument is represented in a formal language by a sequence of discrete units.

It is also commonly accepted that logic, formal or otherwise, deals with the concept of truth. Each unit in a formally presented reasoning is assumed to bear truth-related content. Yet, the use of formal language for the codification of the units of truth-related contents does not make logic formal;
only the realization that the discrete units, called \textit{words}, of a formal language, constituting a deduction chain, are codes of \textit{forms of judgments} rather than ``formulations of factual, contingent knowledge'' (Carnap) does.\footnote{Compare with the notion of elementary propositional function in~\cite{whitehead-russell}, section A $\ast1$.} 

Thus we arrive at a twofold view on the units of formal language, which are used to form the codes of formally presented deductions. On the one hand, they are words of a given formal language, satisfying  some standard requirements; on the other, they are forms that do not possess any factual content but are intended to represent it. Moreover, the content that is meant to be represented by means of a formal unit is intended to generate an assertion, not command, exclamation, question and the like. Having this in mind, we will be calling these formal units  (\textit{sentential}) \textit{formulas} or \textit{terms}, depending on which role they play in our analysis.

Concerning formal deduction, or formal reasoning, we will make a distinction, and emphasize it, between
what is known in philosophy of logic as the \textit{theory of consequence}, or \textit{dialectic argument}, which is rooted in Aristotle's \textit{Analytica Priora}, where the central question is `What follows from what?', and what philosophy of logic calls the \textit{theory of demonstration}, which can be traced to \textit{Analytica Posteriora}, where the main question is `What can be obtained from known (accepted) premises?'. The formulas qualifying to be answers to the latter question will be called \textit{theses} or (formal) \textit{theorems}. We note that the class of such theses depends on the set of \textit{premises}, as well as on the means of the deduction in use. The theory of consequence will be analyzed in terms of a binary relation between the sets of formulas of a given formal language and the formulas of this language. As to ``the means of the deduction in use,'' it will be our main focus throughout the book.\\

After these preparatory remarks, we turn to definition of a formal language. Here we face an obstacle. Since each formal language depends on a particular set of symbols, in order to reach some generality, we are to deal with a schematic language which would contain all main characteristics of many particular languages that are in use nowadays; many other languages of this category can be defined in the future. All our formal languages in focus, with the exception of the ones we discuss in Chapters 7--9, are \textit{sentential}, that is, they are intended to serve for defining \textit{sentential logical systems}.\footnote{E. Beth rightly notes in~\cite{beth1965}, chapter IV, section: ``The sentential logic begins with the observation that not only single terms or concepts but also complete sentences can be considered as the content of a sentence.''}  However, for technical purposes, we will occasionally be employing predicate languages, which will be introduced as needed.

The schematic language we define in this section will not be the only one we are going to deal with. In order to advance Lindenbaum method, as we understand it, to the limits we are able to envisage, we will discuss it in the framework of some fragments of first order language as well.

A \textit{\textbf{sentential schematic language}}\index{language!sentential schematic} $\Lan$ is assumed to contain symbols of the following pairwise disjoint categories:\index{sentential language}
\begin{itemize}
	\item a nonempty set $\Var$ of \textit{sentential} (or \textit{propositional}) \textit{variables}; we will refer to the elements of $\Var$ as $p,q, r,\ldots$, using subscripts if needed; \index{variable!sentential}\index{variable!propositional}\index{propositional!variable}
	\item a set $\Func$  of  \textit{sentential} (or \textit{propositional}) \textit{connectives}\index{propositional!connective} is nonempty; we will refer to the elements of $\Func$ as $F_0,F_1,\ldots,F_{\gamma},\ldots$, where $\gamma$ is any ordinal; for each connective $F_i$, there is a natural number $\#(F_i)$ called the \textit{arity of the connective} $F_i$. We will be omitting index and write simply $F$, if confusion is unlikely; the cardinality of $\Func$ is not bounded from above; that is, $0<\card{\Func}$. The connectives of arity $0$, if any, are called \textit{sentential} (or \textit{propositional}) \textit{constants}\index{propositional!constant}. The constants form a set denoted by $\Cons$ (maybe empty). We will refer to the elements of $\Cons$ as $a,b,c,\ldots$, using subscripts if necessary.
\end{itemize} 

We repeat: the sets $\Var$, $\Cons$, and $\Func$, each consisting of symbols, are pairwise disjoint, and assume also that the elements of each of these sets are distinguishable in such a way that a concatenation of a number of these symbols can be read uniquely; that is, given a word $w$ of $\Lan$ of the length 1, one can decide whether $w\in\Var$ or $w\in\Cons$ or $w\in\Func$; if the word $w$ has the length greater than 1, then, we assume, it is possible to recognize each symbol of $w$ in a unique way.

Thus $\Lan$ gives space for defining many specific languages.
To denote a language defined by a specification within  $\Lan$ we will be using subscripts or other marks attached to $\Lan$. 

The most important concept within a language is that of \textit{formula}.\index{propositiona!formula} To refer to formulas formed by the means of $\Lan$, we use the letters $\alpha,\beta,\gamma,\ldots$ (possibly with subscripts or other marks), calling them $\Lan$-\textit{formulas}. Specifically, these symbols play a role of \textit{informal metavariables} for $\Lan$-formulas. For specified languages we reserve the right to employ other notations, about which the reader will be advised. For a special purpose (to deal with structural inference rules), we introduce \textit{formal metavariables}\index{metavariable} as well.

Given two sentential languages, $\Lan^{\prime}$ and $\Lan^{\prime\prime}$ with $\Var^{\prime}\subseteq\Var^{\prime\prime}$, $\Cons^{\prime}\subseteq\Cons^{\prime\prime}$ and $\Func^{\prime}\subseteq\Func^{\prime\prime}$, $\Lan^{\prime\prime}$ is called an \textit{\textbf{extension}}\index{language!extension} of $\Lan^{\prime}$ and the latter is a \textit{\textbf{restriction}}\index{language!restriction} of the former. An extension $\Lan^{\prime\prime}$ of a language $\Lan^{\prime}$ is called \textit{\textbf{primitive}}\index{language!primitive} if $\Cons^{\prime}=\Cons^{\prime\prime}$ and $\Func^{\prime}=\Func^{\prime\prime}$, that is, when $\Lan^{\prime\prime}$ is obtained only by adding new variables, if any, to $\Var^{\prime}$.

\begin{defn}[$\Lan$-formulas]\label{D:L-formulas}\index{$\Lan$-formula}
	The set {\em$\textbf{Fm}_{\mathcal{L}}$} of $\Lan$-\textbf{formulas} is formed inductively according to the following rules$:$\index{$\mathbf{Fm}_{\mathcal{L}}$} 
	{\em\[
		\begin{array}{cl}
		(\text{a}) &\Var\cup\Cons\subseteq\Forms_{\mathcal{L}};~\textit{if $\alpha\in\Var\cup\Cons$, $\alpha$ is called \textbf{atomic}};\\ 
		(\text{b}) & \textit{if }\alpha_1,\ldots,\alpha_n\in\FormsL\textit{ and $F\in\Func\textit{ with }\#(F)=n$, then }F\alpha_1\ldots\alpha_n\in\FormsL;\\
		(\text{c}) & \textit{the $\Lan$-formalas are only those $\Lan$-words which can be obtained}\\
		&\textit{successively according to the rules }(\text{a})\textit{ and }(\text{b}).\\
		\end{array}
		\]}
\end{defn}

When we deal with more than one language in one and the same context, say $\Lan^{\prime}$ and $\Lan^{\prime\prime}$, which are differ only by their sets of variables,
$\Var^{\prime}$ and $\Var^{\prime\prime}$ respectively, we denote their sets of
formulas by $\Forms_{\mathcal{L}_{\Var^{\prime}}}$ and $\Forms_{\mathcal{L}_{\Var^{\prime\prime}}}$, respectively.

We will find useful the notion of the \textit{\textbf{degree of}} a (given) \textit{\textbf{formula}} $\alpha$\index{degree of formula}, according to the following clauses:
\[
\begin{array}{cl}
(\text{a}) &\text{if $\alpha\in\Var\cup\Cons$, then its degree is $0$};\\
(\text{b}) &\text{if $\alpha=F\alpha_1\ldots\alpha_n$, then the degree of $\alpha$ is the sum of}\\ 
&\text{the degrees of all $\alpha_i$ plus $1$}.
\end{array}
\]
It is easy to see that the degree of a formula $\alpha$ equals the number of occurrences of the sentential connectives that are contained in $\alpha$.

As can be seen from the last definition, the degree of a formula is an estimate of its ``depth.'' Also, in the last definition, it is assumed that if $\alpha=F_{i}\alpha_1\ldots\alpha_n$, that is if $\alpha\notin\Var\cup\Cons$, it begins with some sentential connective ($F$ with $\#(F)=n$ in our case). But then, the words $\alpha_1,\ldots,\alpha_n$ must be recognizable as formulas. Thus it is important to note that for any given word $w$ of the language $\Lan$, one can decide whether $w$ is an $\Lan$-formula or not. One of the ways, by which it can be done, is the construction of a formula tree for the formula in question.\\

Technically, it is convenient to have the notion of metaformula associated with a given language $\Lan$.
\begin{defn}[$\Lan$-metaformulas]\label{D:L-metaformula}\index{$\Lan$-metaformula}\index{metaformula}\index{\em$\Mform$}
	Given a language $\Lan$, the set {\em$\Mform$} of $\Lan$-metaformulas is defined inductively by the clauses$:$
	{\em\[
		\begin{array}{cl}
		(\text{a}) &\Mvar\cup\Cons\subseteq\Mform;\textit{ where $\Mvar$ is a set of \textbf{metavariables} $\bm{\alpha},\bm{\beta},\bm{\gamma},\ldots$};\\
		& \textit{$($with or without subscripts$)$ with }\card{\Mvar}=\card{\FormsL}\textit{ and such that}\\
		&\Mvar\cap\Cons=\varnothing~\textit{and}~\Mvar\cap\Func=\varnothing;\\ 
		(\text{b}) & \textit{if $\phi_1,\ldots,\phi_n\in\Mform$ and $F\in\Func$ with $\#(F)=n$, then $F\phi_1\ldots\phi_n\in\Mform$};\\
		(\text{c}) & \textit{the $\Lan$-metaformalas are only those formal words which can be obtained}\\
		&\textit{successively according to the items }(\text{a})\textit{ and }(\text{b}).\\
		\end{array}
		\]}
\end{defn}

The notion of the \textit{\textbf{degree of metaformula}} is similar to that of formula.\index{degree of metaformula}

The intended interpretation of metavariables is $\Lan$-formulas.\\

Now we turn again to formulas. We will be using the following notation: Given an $\Lan$-formula $\alpha$,
\[
\Var(\alpha)~\textit{is the set of variables that occur in $\alpha$}.
\]

It will be assumed that the set $\Var$ is decidable and the sets $\Cons$ and $\Func$ are decidable when they are countable; that is, for any $\Lan$-word $w$ of the length 1, there is an effective procedure which decides whether $w\in\Var$, whether $w\in\Cons$, or whether $w\in\Func$. In this case, it is also effectively decidable whether any $\Lan$-word is an $\Lan$-formula or not.

The \textit{\textbf{cardinality of a language}} $\Lan$\index{language!cardinality}, symbolically $\card{\Lan}$, is $\card{\Forms_{\mathcal{L}}}$. We will be using both notations. It can be observed that for any language $\Lan$, $\card{\Lan}\ge\aleph_{0}$.

The notion of a \textit{\textbf{subformula of}}\index{subformula} a (given) \textit{\textbf{formula}} is defined inductively through the notion of the \textit{\textbf{set}} $\SubF{\alpha}$ \textit{\textbf{of all subformulas}} \textit{\textbf{of}} a (given) \textit{\textbf{formula}} $\alpha$ as follows:\index{$\SubF{\alpha}$}
\[
\begin{array}{cl}
(\text{a}) &\text{if $\alpha\in\Var\cup\Cons$, then $\SubF{\alpha}=\{\alpha\}$};\\
(\text{b}) &\text{if $\alpha=F\alpha_1\ldots\alpha_n$, then $\SubF{\alpha}=\{\alpha\}\cup\SubF{\alpha_1}\cup
	\ldots\cup\SubF{\alpha_n}$}.
\end{array}
\]
A subformula $\beta$ of a formula $\alpha$ is called \textit{\textbf{proper}}\index{proper subformula} if $\beta\in\SubF{\alpha}$ and $\beta\neq\alpha$. 

A useful view on a subformula of a formula is that the former is a subword of the latter and is a formula itself.

The following observation will implicitly be used in the sequel. 
\begin{prop}
	For each formula $\alpha$ that does not belong to $\Var\cup\Cons$, there is a functional symbol $F$ of arity $n$ and a unique sequence $\beta_1,\ldots,\beta_n$ of formulas such that $\alpha=F\beta_1\ldots\beta_n$.
\end{prop}
\noindent\emph{Proof}~can be found in~\cite{rosenbloom1950}, chapter IV, section 1.\\

The \textit{\textbf{formula tree}}\index{formula tree} for a formula $\alpha$ is a finite directed edge-weighted labeled tree where each node is labeled by a subformula of the initial formula $\alpha$.  The node labeled by $\alpha$ is the \textit{root} of the tree; that is its indegree equals 0. Suppose a node is labeled by  $\beta\in\SubF{\alpha}$. If $\beta\in\Var\cup\Cons$, we call this node a \textit{leaf}; its outdegree equals 0. Now assume that $\beta=F_{i}\beta_1\ldots\beta_n$, that is $\#(F_i)=n$. Then outdegree of this node is $n$. Since $\alpha$ may have more than one occurrence of a subformula $\beta$, all nodes labeled by the last formula will have outdegree equal to $\#(F_i)$. Then, $n$ directed edges will connect this node to $n$ nodes labeled by the subformulas $\beta_1,\ldots,\beta_n$, wherein the weight $i$ is assigned to the edge connecting $\beta$ with $\beta_i$. These $n$ nodes are distinct from one another even if some of $\beta_1,\ldots,\beta_n$ may be equal. We note that each occurrence of every subformula of $\alpha$ is a label of a formula tree for $\alpha$. Also, we observe that the degree of each $\beta_i$ is less than the degree of $\beta$. A formula tree for a given formula $\alpha$ is completed if each complete path from the root ends with a leaf. This necessarily happens, for the degree of a formula is a nonnegative integer.

We note that we use weighted edges only in order to identify a (unique) path from the root of a tree to a particular node.

It should be obvious that all formula trees for a given formula are isomorphic as graphs. So one can speak of \textit{the} formula tree for a given formula.

Let us illustrate the last definition for the formula $F_{i}F_{j}app$ with $\#(F_i)=\#(F_j)=2$.
Below we identify the nodes with their labels.

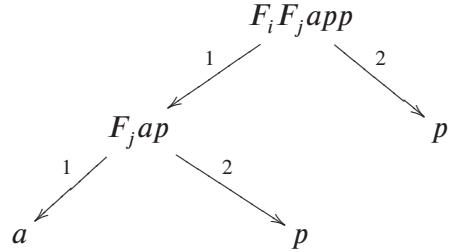
\begin{figure}[h!]
	\[
	\xymatrix{
		&&F_{i}F_{j}app \ar[dl]_1 \ar[dr]^2\\
		&F_{j}ap \ar[dl]_1 \ar[dr]^2 &&p\\
		a &&p
	}
	\]
	\caption{The formula tree for $F_{i}F_{j}apq$}
	\label{Fig-1}
\end{figure}

It should be clear that the formula tree for an initial formula allows us to work with each occurrence of every subformula of this formula separately and also to follow through the subformula relation even within a designated occurrence of a subformula of the initial formula. This observation will be used in the operation of replacement below.

Let us consider any leaf of the formula tree of Figure~\ref{Fig-1}, say the one labeled by $a$. We observe that the nodes of the path from the root to the selected leaf consist of all subformulas of the initial formula which contain $a$ as their subformula. Now let us color the first occurrence of $p$ in red and the second in blue so that we obtain
\[
F_{i}F_{j}a\textcolor{red}{p}\textcolor{blue}{p}. \tag{initial formula}
\]

We use color only for the sake of clarity to identify the path we are working with at the moment (see below); we could use weights of edges instead.

Next we consider the two paths: the first,
\[
F_{i}F_{j}a\textcolor{red}{p}\textcolor{blue}{p}\stackrel{1}{\longrightarrow} F_{j}a\textcolor{red}{p}
\stackrel{2}{\longrightarrow}\textcolor{red}{p},\tag{red}
\]
corresponds to the red occurrence of $p$ in the sense that  its nodes constitute all subformulas of the initial formula, which contain this red occurrence of $p$; the second path,
\[
F_{i}F_{j}a\textcolor{red}{p}\textcolor{blue}{p}\stackrel{1}{\longrightarrow}
\textcolor{blue}{p}, \tag{blue}
\]
corresponds to the blue occurrence of $p$ and satisfies the same property regarding this occurrence.
Now we substitute for both occurrences of $p$ in the initial formula any formula $\alpha$ so that we obtain the formula
\[
F_{i}F_{j}a\textcolor{red}{\alpha}\textcolor{blue}{\alpha}. \tag{red-blue $\alpha$}
\]
Instead of the tree of Figure 1 we will get the new tree:
\begin{figure}[h!]
	\[
	\xymatrix{
		&&F_{i}F_{j}a\textcolor{red}{\alpha}\textcolor{blue}{\alpha}  \ar[dl]_1 \ar[dr]^2\\
		&F_{j}a\textcolor{red}{\alpha} \ar[dl]_1 \ar[dr]^2 &&\textcolor{blue}{\alpha}\\
		a &&\textcolor{red}{\alpha}
	}
	\]
	\caption{The transformed tree after the substitution $p\mapsto\alpha$}
	\label{Fig-2}
\end{figure}
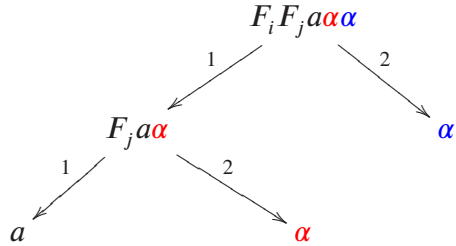

If $\alpha\notin\Var\cup\Cons$, the tree of Figure~\ref{Fig-2} is not a formula tree. However, the last tree shows that the $\Lan$-word
$F_{i}F_{j}a\alpha\alpha$ obtained by this substitution is an $\Lan$-formula. Instead of the paths (red) and (blue) we observe the following paths:
\[
F_{i}F_{j}a\textcolor{red}{\alpha}\textcolor{blue}{\alpha}\stackrel{1}{\longrightarrow} F_{j}a\textcolor{red}{\alpha}
\stackrel{2}{\longrightarrow}\textcolor{red}{\alpha}\tag{red $\alpha$}
\]
and 
\[
F_{i}F_{j}a\textcolor{red}{\alpha}\textcolor{blue}{\alpha}\stackrel{2}{\longrightarrow}
\textcolor{blue}{\alpha}. \tag{blue $\alpha$}
\]
And again, we notice the same property: the nodes of (red $\alpha$) form the set of subsets of the red occurrence of $\alpha$ in the formula (red-blue $\alpha$); the same property regarding the blue occurrence of $\alpha$ is true for the path (blue $\alpha$).

This observation will help us discuss a syntactic transformation which is vital for the theory of sentential
logic. Essentially, this transformation makes speak of logic as a system of forms of judgments rather than that of individual judgments about facts.

First we define a transformation related to a formula $\alpha$ and a sentential variable $p$, denoting this operation by $\sigma^{\alpha}_{p}$. Before proceeding, we note that the initial formula $F_iF_japp$ was chosen rather arbitrarily; it might be that it does not contain $p$.
Having this in mind, given any $\Lan$-formula $\alpha$ and any sentential variable $p$, we define a
\textit{\textbf{unary substitution}}\index{substitution!unary} $\sigma^{\alpha}_{p}:\Forms_{\mathcal{L}}\longrightarrow\Forms_{\mathcal{L}}$ as follows:
\[
\begin{array}{cl}
\bullet & \text{if $\beta\in\Var\cup\Cons$ and $\beta\neq p$, then $\sigma^{\alpha}_{p}(\beta)=\beta$};\\
\bullet & \text{if $\beta = p$, then $\sigma^{\alpha}_{p}(\beta)=\alpha$};\\
\bullet &\text{if $\beta=F_{i}\beta_1\ldots\beta_n$, then $\sigma^{\alpha}_{p}(\beta)=F_{i}\sigma^{\alpha}_{p}(\beta_1)\ldots\sigma^{\alpha}_{p}(\beta_n)$}.
\end{array}
\]

At first sight, it may seem that unary substitution is not powerful enough. For example, let us take a formula $F_{i}F_{j}apq$, where  connectives $F_i$ and $F_j$ have the arity equal to 2. Assume that we want to produce a formula $F_{i}F_{j}a\alpha\beta$. Suppose we first obtain:
\[
\sigma^{\alpha}_{p}(F_{i}F_{j}apq)=F_{i}F_{j}a\alpha q.
\]
And we face a problem if $\alpha$ contains $q$, for then we get
\[
\sigma^{\beta}_{q}(F_{i}F_{j}a\alpha q)=F_{i}F_{j}a\sigma^{\beta}_{q}(\alpha)\beta,
\]
which is not necessarily equal to $F_{i}F_{j}a\alpha\beta$.
To overcome this difficulty, we select a variable which occurs neither in $\alpha$ nor in the initial formula $F_{i}F_{j}apq$. Suppose this variable is $r$. Then, we first produce:
\[
\sigma^{r}_{q}(F_{i}F_{j}apq)=F_{i}F_{j}apr.
\]
Now the desirable formula can be obtained in two steps:
\[
\begin{array}{l}
\sigma^{\alpha}_{p}(F_{i}F_{j}apr)=F_{i}F_{j}a\alpha r,\\
\sigma^{\beta}_{r}(F_{i}F_{j}a\alpha r)=F_{i}F_{j}a\alpha\beta.
\end{array}
\]
Thus we observe:
\[
\sigma^{\beta}_{r}\circ\sigma^{\alpha}_{p}\circ\sigma^{r}_{q}(F_{i}F_{j}apq)=F_{i}F_{j}a\alpha\beta.
\]

This leads to the conclusion that the composition of two and more unary substitutions may result in substitution of a different kind; namely, we can obtain the formula $F_{i}F_{j}a\alpha\beta$ from the formula $F_{i}F_{j}apq$ by  substituting simultaneously 
the formulas $\alpha$ and $\beta$ for the variables $p$ and $q$, respectively. Thus we arrive at the following definition.
\begin{defn}[substitution]\label{D:substitution} \index{substitution}
	Given a formal language $\Lan$, a syntactic transformation {\em$\sigma:\Forms_{\mathcal{L}}\longrightarrow\Forms_{\mathcal{L}}$} is called a \textit{\textbf{uniform}}\index{substitution!uniform} $($or \textit{\textbf{simultaneous}}$)$ \textit{\textbf{substitution}} $($or simply a \textit{\textbf{substitution}}\index{substitution!simultaneous} or $\Lan$-\textit{\textbf{substitution}}$)$ if it is an extension of a map {\em$\sigma_{0}:\Var\longrightarrow
		\Forms_{\mathcal{L}}$}, satisfying the following conditions$\,:$
	{\em\[
		\begin{array}{cl}
		(\text{a}) &\text{if $\alpha\in\Cons$, then $\sigma(\alpha)=\alpha$};\\
		(\text{b}) &\text{if $\alpha=F_{i}\alpha_1\ldots\alpha_n$, then $\sigma(\alpha)=F_{i}\sigma(\alpha_1)\ldots\sigma(\alpha_n)$}.
		\end{array}
		\]}
	Given a substitution $\sigma$ and a formula $\alpha$, $\sigma(\alpha)$ is called a \textbf{substitution instance}\index{substitution instance} of $\alpha$. The set of $\Lan$-substitutions is denoted by $\SbsL$.\index{$\SbsL$}
\end{defn}

For a set $X$ of $\Lan$-formulas and an $\Lan$-substitution $\sigma$, we denote:
\[
\sigma(X):=\set{\sigma(\alpha)}{\alpha\in X}~\text{and}~\sigma^{-1}(X):=\set{\alpha}{\sigma(\alpha)\in X}.
\]
Thus, for any substitution $\sigma$,
\begin{equation}\label{E:substitution-for-X}
\sigma(\varnothing)=\varnothing~\text{and}~\sigma^{-1}(\varnothing)=\varnothing.
\end{equation}
Also, we observe:
\begin{equation}\label{E:substitution-inequalities}
i)~\sigma(\sigma^{-1}(X))\subseteq X\quad\text{and}\quad ii)~ X\subseteq\sigma^{-1}(\sigma(X)).
\end{equation}
(Exercise~\ref{section:languages}.\ref{EX:substitution-inequalities}) 

We leave for the reader to prove that, given a fixed formula $\alpha$, any simultaneous substitution applied to $\alpha$ results in the composition of singular substitutions applied to $\alpha$. (See Exercise~\ref{section:languages}.\ref{EX:substitution-2}.)
\begin{prop}\label{P:substitution}
	For any map {\em$\sigma_{0}:\Var\longrightarrow\Forms_{\mathcal{L}}$}, there is a unique substitution {\em$\sigma:\Forms_{\mathcal{L}}\longrightarrow\Forms_{\mathcal{L}}$} extending $\sigma_0$.
\end{prop}
\noindent\textit{Proof}~is left to the reader. (Exercise~\ref{section:languages}.\ref{EX:substitution}.)\\

\textit{The last proposition allows us to think of a substitution simply as a map} $\sigma:\Var\longrightarrow\Forms_{\mathcal{L}}$, which will be convenient when we want to define a particular substitution. On the other hand, the property (b) of Definition~\ref{D:substitution} allows us to look at any substitution as a mapping satisfying a special algebraic property, which will also be convenient when we want to use this property.

For now, we want to note that composition of finitely many substitutions is a substitution. Moreover, for any substitutions $\sigma_1$, $\sigma_2$ and $\sigma_3$,
\begin{equation}\label{E:substitution-associative}
\sigma_1\circ(\sigma_2\circ\sigma_3)=(\sigma_1\circ\sigma_2)\circ\sigma_3
\end{equation}
(We leave this property for the reader; see Exercise~\ref{section:languages}.\ref{EX:substitution-associative}.) This implies that the set of all substitutions, related to a language $\Lan$, along with composition $\circ$ constitute a monoid with the identity element
\[
\iota(p)=p,~\text{for any $p\in\Var$}
\]
which will be calling the \textit{\textbf{identity substitution}}.\index{substitution!identity} It must be clear that for any formula $\alpha$,
\[
\iota(\alpha)=\alpha.
\]
Because of \eqref{E:substitution-associative}, any grouping in
\[
\sigma_{1}\circ\ldots\circ\sigma_n
\]
leads to the same result. For this reason, parentheses can be omitted, if we do want to be specific. 

We will call a set $F\subseteq\Forms_{\mathcal{L}}$ \textit{\textbf{closed under substitution}}\index{closed under substitution}
if for any $\Lan$-formula $\alpha$ and substitution $\sigma$,
\[
\alpha\in F\Longrightarrow\sigma(\alpha)\in F.
\]

We have already noted that the tree of Figure~\ref{Fig-2} is not a formula tree if $\alpha\notin\Var\cup\Cons$. However, we can build the formula tree for $\alpha$. It does not matter how it looks; so we depict it schematically as follows:

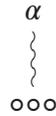
\begin{figure}[h!]
	\[
	\xymatrix{
		\alpha \ar@{~}[d] \\
		\circ\circ\circ 
	}
	\]	
	\caption{The formula tree for $\alpha$}
	\label{Fig-3}
\end{figure}

Attaching two copies of the last tree to the tree of Figure~\ref{Fig-2}, we obtain the tree for $F_{i}F_{j}a\alpha\alpha$.

\begin{figure}[h!]
	\[
	\xymatrix{
		&&F_{i}F_{j}a\textcolor{red}{\alpha}\textcolor{blue}{\alpha}  \ar[dl]_1 \ar[dr]^2\\
		&F_{j}a\textcolor{red}{\alpha} \ar[dl]_1 \ar[dr]^2 &&\textcolor{blue}{\alpha} \ar@{~}[d]\\
		a &&\textcolor{red}{\alpha} \ar@{~}[d] &\circ\circ\circ\\
		&&\circ\circ\circ
	}
	\]
	\caption{The formula tree for $F_{i}F_{j}a\alpha\alpha$}
	\label{Fig-4}
\end{figure}
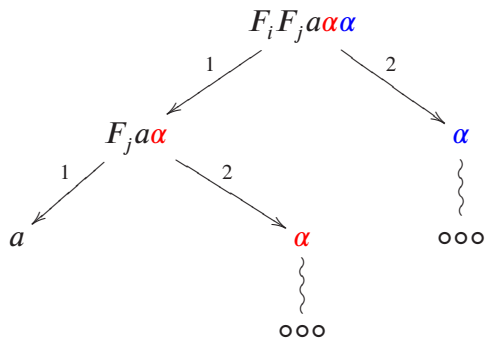
\pagebreak

We will call one copy of the formula tree for $\alpha$ \textit{the tree for red} $\alpha$ and the other copy \textit{the tree for the blue} $\alpha$.

Now let us take any formula $\beta$  and construct the formula tree for it depicted as follows:
\begin{figure}[h!]
	\[
	\xymatrix{
		\beta \ar@{~}[d] \\
		\bullet\bullet\bullet 
	}
	\]	
	\caption{The formula tree for $\beta$}
	\label{Fig-5}
\end{figure}
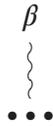

Now in the path
\[
F_{i}F_{j}a\textcolor{red}{\alpha}\textcolor{blue}{\alpha}\stackrel{1}{\longrightarrow}
F_{j}a\textcolor{red}{\alpha}\stackrel{2}{\longrightarrow} \textcolor{red}{\alpha}
\]
\pagebreak
of the tree of Figure~\ref{Fig-4}, we replace all red occurrences of $\alpha$ by $\beta$ and replace the tree for red alpha by the tree for $\beta$ (Figure~\ref{Fig-5}), arriving thus at the following tree:

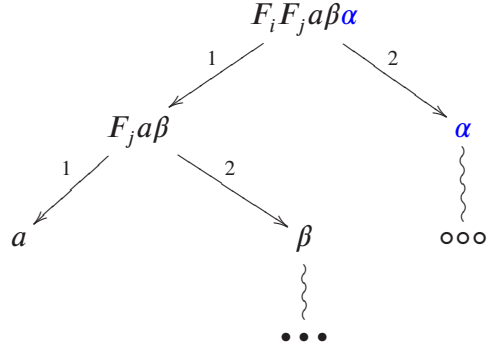
\begin{figure}[h!]
	\[
	\xymatrix{
		&&F_{i}F_{j}a\beta\textcolor{blue}{\alpha}  \ar[dl]_1 \ar[dr]^2\\
		&F_{j}a\beta \ar[dl]_1 \ar[dr]^2 &&\textcolor{blue}{\alpha} \ar@{~}[d]\\
		a &&\beta \ar@{~}[d] &\circ\circ\circ\\
		&&\bullet\bullet\bullet
	}
	\]
	\caption{The formula tree for $F_{i}F_{j}a\beta\alpha$}
	\label{Fig-6}
\end{figure}

This is the formula tree for $F_{i}F_{j}a\beta\alpha$. The syntactic transformation just described shows that replacing a designated subformula of a formula, we obtain a formula. Such a transformation is called \textit{\textbf{replacement}}\index{replacement}. Designating a subformula $\beta\in\SubF{\alpha}$, we write $\alpha[\beta]$. The result of the replacement of designated occurrences (one or more) of a subformula $\beta$ with a formula $\gamma$, we denote by $\alpha[\gamma]$. We note that the last notation does not specify, which occurrences of $\beta$ were replaced. It is clear that the procedure of replacement can be extended. Namely when we want to replace designated occurrences of subformulas $\beta_1,\ldots,\beta_n$ of a formula $\alpha$, symbolically $\alpha[\beta_1,\ldots,\beta_n]$, with formulas $\gamma_1,\ldots,\gamma_n$, respectively, the resulting formula will be denoted by $\alpha[\gamma_1,\ldots,\gamma_n]$. Since the procedure of replacement depends on the selection of occurrences of subformulas of a given formula, it does not matter, whether a resulting formula is obtained by simultaneous replacement or by a finite sequence of replacements, in each of which only one subformula is replaced.\\

The following operator will be useful in the sequel. We define an operator
$\Sb:\mathcal{P}(\Forms_{\mathcal{L}})\longrightarrow\mathcal{P}(\Forms_{\mathcal{L}})$ as follows:
\[
\Sb X=\set{\sigma(\alpha)}{\sigma\in\SbsL~\text{and}~\alpha\in X}.
\]

\begin{prop}\label{P:Sb-properties}
	The following properties hold:
	{\em\[
		\begin{array}{cl}
		(\text{a}) &X\subseteq \Sb X;\\
		(\text{b}) &\Sb(X\cup Y)=\Sb X\cup\Sb Y;\\
		(\text{c}) &\Sb(X\cap Y)=\Sb X\cap\Sb Y;\\
		(\text{d}) &X\subseteq Y\Longrightarrow\Sb X\subseteq\Sb Y;\\
		(\text{e}) &\textit{for any $\sigma\in\SbsL$},~\sigma(\Sb X)\subseteq\Sb X;\\
		(\text{f}) & \Sb(\Sb X)=\Sb X.
		\end{array}
		\]}
\end{prop}
\begin{proof}
	(a) follows from the observation that $\iota(X)=X$. (b), (c) and (d) are obvious. (e) follows from the fact that a compositions of two (and more) substitutions is a substitution. (f) follows from (a) and (e).
\end{proof}

Before going over examples of some sentential languages, we discuss how any metaformula can be seen in the light of the notion of substitution. It must be clear that the intended interpretation of any metaformula $\phi$ is the set of formulas of a particular shape.

Similar to Definition~\ref{D:substitution}, any map $\bm{\sigma}_0:\Mvar\longrightarrow
\Forms_{\mathcal{L}}$ can uniquely be extended to $\bm{\sigma}:\Mform\longrightarrow\Forms_{\mathcal{L}}$.  We will call both $\bm{\sigma}_0$ and its extension $\bm{\sigma}$ a \textit{\textbf{realization}}\index{realization}, or \textit{\textbf{instantiation}}\index{instantiation}, of metaformulas in $\Lan$. The following observation will be useful in the sequel. An instantiation $\bm{\sigma}$ is called \textit{\textbf{simple}}\index{instantiation!simple} if it is the extension of $\bm{\sigma}_0:\Mvar\longrightarrow
\Var_{\mathcal{L}}$.
\begin{prop}\label{P:metaformula-instantiations}
	Let $\phi$ be an $\Lan$-metaformula and $\Sigma_\phi$ be the class of all instantiations of $\phi$. Then for any $\Lan$-substitution $\sigma$ and any $\Lan$-formula $\alpha$,
	\[
	\alpha\in\Sigma_\phi\Longrightarrow\sigma(\alpha)\in\Sigma_\phi.
	\]
\end{prop}
\noindent\textit{Proof}~is leaft to the reader. (Exercise~\ref{section:languages}.\ref{EX:metaformula-instantiations})\\

Below we give examples of languages that are specifications of the schematic language $\Lan$.
All these languages include a denumerable set $\Var$ as the set of sentential variables. 
\begin{description}
	\item[$\Lan_{A}$:]\index{$\Lan_{A}$} In this language, there is no sentential constants; there are four (assertoric) logical connectives: $\wedge$ (\textit{conjunction}), $\vee$ (\textit{disjunction}), $\rightarrow$ (\textit{implication} or \textit{conditional}), and $\neg$ (\textit{negation}), having the arities $\#(\wedge)=\#(\vee)=\#(\rightarrow)=2$ and
	$\#(\neg)=1$. It is customary to use an infix notation for this language; however, the use of
	infix notation for this language requires two additional symbols: `(' and `)' called the left and right parenthesis, respectively. Accordingly, the definition of formula will be different: $(\alpha\wedge\beta)$, $(\alpha\vee\beta)$, $(\alpha\rightarrow\beta)$ and $\neg\alpha$ are formulas, whenever $\alpha$ and $\beta$ are. \textit{The parentheses do not apply to atomic formulas and negations}. Sometimes, we will not employ all the connectives of $\Lan_{A}$ but some of them, continuing to call this restricted language $\Lan_{A}$. Sometimes, $\Lan_{A}$ is considered in its extended version, when a connective $\leftrightarrow$ is added; if not, then the last connective is used as an abbreviation:
	\[
	\alpha\leftrightarrow\beta:=(\alpha\rightarrow\beta)\wedge(\beta\rightarrow\alpha).
	\]
	\item[$\Lan_{B}$:]\index{$\Lan_{B}$} This language extends $\Lan_{A}$ by adding to it two constants: $\top$ (\textit{truth} or \textit{verum}) and $\bot$ (\textit{falsehood} or \textit{falsity} or \textit{falsum} or \textit{absurdity}). Sometimes, only one of these constants will be employed.
	\textit{The parentheses do not apply to atomic formulas, negations and to the constants $\top$ and $\bot$.} We will call all these variants by $\Lan_{B}$, indicating every time the list of constants. 
	\item[$\Lan_{C}$:]\index{$\Lan_{C}$} This language extends either $\Lan_{A}$ or $\Lan_{B}$ with two sentential connectives (\textit{modalities}): $\Box$ (\textit{necessity}) and $\Diamond$ (\textit{possibility}) both of the arity 1, that is $\#(\Box)=\#(\Diamond)=1$. In addition to the corresponding requirements above, \textit{the parentheses do not apply to the formulas that start with the modalities $\Box$ and $\Diamond$.} Sometimes, when only one of these modalities is employed, we continue to call this language $\Lan_{C}$.
\end{description}

In the sequel we will employ more formal languages, which do not fit the proposed framework but closely related to it.\\

\pagebreak
\noindent{\textbf{Exercises~\ref{section:languages}}}
\begin{enumerate}
	\item\label{EX:formula-tree} Here we give another definition of a formula tree. This new definition is convenient when we deal with an algorithmic problem, in which case the sets $\Var$, $\Cons$ and $\Func$ are regarded countable and effectively decidable. One of such problems will be raised below as an exercise. We note that in this definition a designated path can be identified exclusively by the nodes it comprises.
	
	The \textit{\textbf{formula tree}}\index{formula tree} for a formula $\alpha$ is a finite directed labeled tree where each node is labeled with a pair $\langle n,\beta\rangle$, where $n$ is a natural number and $\beta$ is a subformula of $\alpha$.  
	
	There is a unique node, called the \textit{root}, labeled by $\langle\p{0},\alpha\rangle$, the ingdegree of which is $0$.\footnote{We remind the reader that $\p{0},\p{1},\p{2},\ldots$ is the sequence of prime numbers.} The ingdegree of all other nodes of the tree, if any, equals $1$.
	Let a node $v$ be labeled with $\langle n,\beta\rangle$, where $\beta\in\Var\cup\Cons$. Then the outdegree of $v$ is always $0$. Such nodes are called \textit{leaves}. (As the reader will see below, leaves are always present in a formula tree.)
	Thus if $\alpha\in\Var\cup\Cons$, the formula tree of $\alpha$ consists of only one node, the root, which is also a leaf.  Assume that $\beta=F\beta_1\ldots\beta_k$ with $F\in\mathcal{F}$ is a proper subformula of $\alpha$ and a node $v$ is labeled with $\langle n,\beta\rangle$. Suppose that $\p{j}$ is the greatest prime number that divides $n$.  Then the outdegree of $v$ is $\#(F)$, that is in our case, equals $k$. And then $v$ will be connected via $k$ edges to $k$ different nodes, $v_1,\ldots,v_k$. (The direction of each such edge is from $v$ to $v_i$.) We label $v_i$ with $\langle n\cdot\p{j+i},\beta_i\rangle$. Thus, 
	if $\beta_i=\beta_j$, then the two distinct nodes, $v_i$ and $v_j$, will be labeled with labels distinct in their first place and equal in the second. Also, we note that the degree of each $\beta_i$ is less than the degree of $\beta$. 
	
	Use the new definition of formula tree to prove the following.
	\begin{quote}
		{\em 	There is a mechanical procedure $($an algorithm$)$ which, given a list of pairs $\langle n_1,\alpha_1\rangle,\ldots,\langle n_k,\alpha_k\rangle$, where each $n_i$ is a positive integer, can assemble a formula tree, where the pairs are regarded as the nodes of the tree, or states that such a tree does not exist.}	
	\end{quote}
	\item \label{EX:substitution-inequalities} Prove the inclusions \eqref{E:substitution-inequalities}.
	\item\label{EX:substitution-2} Prove that $(a)$ for any fixed formula $\alpha$ and any fixed substitution $\sigma$, there are unary substations $\sigma_1,\ldots,\sigma_n$ such that $\sigma(\alpha)=\sigma_1\circ\ldots\circ\sigma_n(\alpha)$; however, $(b)$ for any fixed substitution $\sigma$, there is no finitely many fixed unary substitutions $\sigma_1,\ldots,\sigma_n$ such that for any formula $\alpha$, $\sigma(\alpha)=\sigma_1\circ\ldots\circ\sigma_n(\alpha)$.	
	\item\label{EX:substitution} Prove Proposition~\ref{P:substitution}.
	\item \label{EX:substitution-associative} Prove \eqref{E:substitution-associative}.
	\item \label{EX:metaformula-instantiations}Prove Proposition~\ref{P:metaformula-instantiations}.
	
\end{enumerate}

\section{Semantics of sentential formal language}\label{section:semantics}
One of our concerns in this book is to show how truth can be transmitted through formulas of a formal language by formally presented deduction. Calling a formula true, we will mean a kind of truth which does not depend on any factual content. Moreover, sharing the view of the multiplicity of the nature of truth, we will be dealing rather with degrees of truth, calling them \textit{truth values}. Then, the idea of truth will find its realization in the idea of \textit{validity}.

Concerning the independence of any factual, contingent content in relation to truth values, let us consider the following two sentences --- `Antarctica is large' and `Antarctica is large or Antarctica is not large'. If our formal language contains a symbol for the grammatical disjunction ``$\ldots_{1}~\text{or}~\ldots_{2}$'', say ``$\ldots_{1}\vee\ldots_{2}$''
(as in the languages $\Lan_{A}$--$\Lan_{C}$ of Section~\ref{section:languages}), and a symbol for creating the denial of sentence, say ``$\neg\ldots$'' (also as in $\Lan_{A}$--$\Lan_{C}$), we can write in symbolic form the first sentence as $p$ and the second as $p\vee\neg p$. From the viewpoint of Section~\ref{section:languages}, both formal expressions are formulas of the language $\Lan_{A}$. Now, depending on familiarity with geography, as well as of personal experience, one can count $p$ true or false. The second formula $p\vee\neg p$, however, in some semantic systems takes steadily the truth value \textit{true}, regardless of any personal relation towards $p$. The \textit{universality} of this kind will be stressed in the idea of validity.

From the last example, it should be clear the distinction between the statements, the truth value of which (understood as the  degree of trust or that of informativeness) is dependent on the contingency of facts, the level of knowledge or personal beliefs, and the formal judgments, which obtain their logical status or meaning determined by \textit{semantic rules}. These rules will allow us to interpret the formal judgments of a language in abstract structures, reflecting the structural character of these judgments, as well as the universal character of their validity, thereby understanding their invalidity as the failure to satisfy this universality.

In this book, semantics will play pragmatic rather than explanatory role;\footnote{However, since we state validity conditions for formal judgments, these conditions generate the meanings of the judgments, according to the point of view that has been emphasized by Wittgenstein in \textit{Tractatus}; cf.~\cite{wittgenstein2001}, 4.024. Also, compare with Curry's statement  that ``the meaning of a concept is determined by the conditions under which it is introduced into discourse;'' cf.~\cite{curry1966}, section II.2.} namely we will be stressing the use of semantic structures as \textit{separating tools}. This will be explained in Section~\ref{section:separating-means}. Moreover, we want to emphasize that in our approach to semantics we adhere to the relativistic position, according to which a single formal judgment can be valid according to one semantic system and invalid with respect to another.\\

Now we turn to the definition of semantics of the schematic language $\Lan$ exemplified by the languages $\Lan_{A}$--$\Lan_{C}$ in Section~\ref{section:languages}. Depending on one's purposes, or depending on the lines, along which the semantics in question can be employed, the semantics we propose can use the items of the following two categories --- the \textit{algebras of type} $\Lan$ and expansions of the latter called \textit{logical matrices}. 
Later on, we will generalize the notion of logical matrix, when we discuss consequence relation in Chapter~\ref{chapter:consequence}.

The following semantic concepts will be specified later on as needed.
\begin{defn}[algebra of type $\Lan$]\index{algebra!of type $\Lan$}
	An \textbf{algebra of type} $\Lan$ $($or an $\Lan$-\textbf{algebra} or simply an \textbf{algebra} when a type is specified ahead and confusion is unlikely$)$ is a structure {\em$\alg{A}=\langle\textsf{A};\Func,\Cons\rangle$}, where {\em\textsf{A}} is a nonempty set called the \textbf{carrier} $($or \textbf{universe}$)$ \textbf{of} {\em\alg{A}}$;$ $\Func$, which in the notion of $\Lan$-formulas represents logical connectives, now represents the set of finitary operations on {\em\textsf{A}} of the corresponding arities$;$ and $\Cons$, which is a set of logical constants, now represents a set of nullary operations on {\em\textsf{A}}. As usual, we use the notation {\em$|\alg{A}|=\textsf{A}$}. The part $\langle\Func,\Cons\rangle$ of the structure {\em\alg{A}} is called its \textbf{signature}. Thus, referring to an algebra {\em\alg{A}}, we can say that it is an algebra of signature $\langle\Func,\Cons\rangle$ or of signature $\Func$, if $\Cons=\varnothing$.
\end{defn}

Thus any $\Lan$-formula in relation to any algebra of type $\Lan$ becomes a \textit{\textbf{term}}\index{term}. Naturally, the sentential variables of a formula understood as term become the \textit{\textbf{individual variables}}\index{variable!individual} of this term. Since each algebraic term of, say, $n$ variables can be read as an $n$-ary derivative operation in an algebra of type $\Lan$, so can be understood any $\Lan$-formula containing $n$ sentential variables. 
\begin{defn}[inessential expansion of an algebra, equalized\index{algebra!inessential expansion} algebras]\label{D:equalized-algebras}
	Given an algebra {\em$\alg{A}=\langle\textsf{A};\Func,\Cons\rangle$}, we obtain an \textbf{inessential expansion of} {\em$\alg{A}$} if we add to $\Func$ one or more terms regarding them as new signature operations, or if we add to $\Cons$ one or more symbols with their intended interpretation as elements of {\em\textsf{A}}. Two algebras $($or their signatures$)$ are called \textbf{equalized} if they have inessential expansions of one and the same type.
\end{defn}

The central notion of Lindenbaum method is that of formula algebra associated with a particular language. The \textbf{\textit{formula algebra}}\index{formula algebra} of type $\Lan$ has $\textbf{Fm}_{\mathcal{L}}$ as its carrier and the signature operations as defined in Definition~\ref{D:L-formulas}. This algebra is denoted by $\mathfrak{F}_{\mathcal{L}}$. Thus we have:\index{$\mathfrak{F}_{\mathcal{L}}$}
\[
\mathfrak{F}_{\mathcal{L}}=\langle\Forms_{\mathcal{L}};\Func,\Cons\rangle.
\]

Analogously,  the set of $\Lan$-metaformulas along with the set $\Cons$ regarded as nullary operations and the set
$\Func$ regarded as non-nullary operations constitute a \textit{\textbf{metaformula algebra}}\index{metaformula algebra} $\MformAl$ of type $\Lan$.\\

A relation between the $\Lan$-formulas and the algebras of type $\Lan$ is established through the following notion.
\begin{defn}[valuation, value of a formula in an algebra]\label{D:valuation}\index{valuation}\index{formula!valuation}
	Let {\em\alg{A}} be an algebra of type $\Lan$. Then any map {\em$v:\Var\longrightarrow|\alg{A}|$}	is called a \textbf{valuation} $($or an \textbf{assignment}$)$ in the algebra {\em\alg{A}}. Then, given a valuation $v$, the value of an $\Lan$-formula
	$\alpha$ with respect to $v$, in symbols $v[\alpha]$, is defined inductively as follows$\,:$
	{\em\[
		\begin{array}{cl}
		(\text{a}) &v[p]=v(p),~\textit{for any $p\in\Var$};\\
		(\text{b}) &v[c]~\textit{is the value of the nullary operation of}~\alg{A} ~\textit{corresponding to}\\ 
		&\text{the constant $c$ of $\Lan$};\\
		(\text{c}) &v[\alpha]=Fv[\alpha_1]\ldots v[\alpha_n]~\textit{if $\alpha=F\alpha_1\ldots\alpha_n$}.
		\end{array}
		\]}
	If $\Var^{\prime}\subset\Var$, then a map {\em$v:\Var^{\prime}\longrightarrow|\alg{A}|$}	is called a $\Var^{\prime}$-\textbf{valuation}, or a \textbf{valuation restricted to} $\Var^{\prime}$.
\end{defn}

It makes sense to denote the mapping 
\[
\alpha\mapsto v[\alpha]
\]
by the same letter $v$. The following observation is a justification for this usage.
\begin{prop}\label{P:valuation}
	Let $v$ be a valuation in an algebra {\em\alg{A}}. Then the mapping $v:\alpha\mapsto v[\alpha]$ defines a homomorphism {\em$\mathfrak{F}_{\mathcal{L}}\longrightarrow\alg{A}$}. Conversely, each homomorphism {\em$\mathfrak{F}_{\mathcal{L}}\longrightarrow\alg{A}$} is a valuation.
\end{prop}
\noindent\textit{Proof}~of the first part is based on the clauses (a)--(c) of Definition~\ref{D:valuation}. The second part is obvious.\\

Thus, given a valuation $v$ in \alg{A} and a formula $\alpha(p_1,\ldots,p_n)$ where $p_1,\ldots,p_n$ are all variables occurring in $\alpha$, we observe that
\[
v[\alpha]=\alpha(v[p_1],\ldots,v[p_n]).
\]

Grounding on the last observation base in turn on Proposition~\ref{P:valuation}, very often a valuation is defined as a 
homomorphism $\mathfrak{F}_{\mathcal{L}}\longrightarrow\alg{A}$.\\

Comparing (b)--(c) of Definition~\ref{D:valuation} with (a)--(b) of Definition~\ref{D:substitution}, we obtain the following.
\begin{prop}\label{P:substitution-as-homomorphism}
	Any substitution is a valuation in $\mathfrak{F}_{\mathcal{L}}$ and vice versa. Hence any substitution is an endomorphism $\mathfrak{F}_{\mathcal{L}}\longrightarrow\mathfrak{F}_{\mathcal{L}}$.
\end{prop}

The converse of the last conclusion is also true: any endomorphism $\mathfrak{F}_{\mathcal{L}}\longrightarrow\mathfrak{F}_{\mathcal{L}}$ can be regarded as a substitution. (Exercise~\ref{section:semantics}.\ref{EX:endomorphism}.)\\

The next property, which a combination of Proposition~\ref{P:valuation} and Proposition~\ref{P:validity-homomorphism}, will be useful in the sequel.
\begin{prop}\label{P:valuation-and-substitution}
	Given an $\Lan$-algebra {\em$\alg{A}$} and a valuation $v$ in {\em$\alg{A}$}, for any substitution $\sigma$, there is a valuation $v_{\sigma}$ in {\em\alg{A}} such that for any formula $\alpha$, $v[\sigma(\alpha)]=v_{\sigma}[\alpha]$.
\end{prop}
\begin{proof}
	Indeed, we simply define
	\[
	v_{\sigma}:=v\circ\sigma.
	\]
	
	By Proposition~\ref{P:substitution-as-homomorphism}, $v_{\sigma}$ is a homomorphism, and, in virtue of Proposition~\ref{P:substitution}, $v_{\sigma}$ is a valuation in $\alg{A}$.
\end{proof}

\begin{prop}\label{P:valuation-in-algebra-generated-by-X}
	Let $v_0$ be a valuation in an algebra {\em\alg{A}} of type $\Lan$. Assume that
	{\em\alg{A}} is generated by the set $\set{v_0[p]}{p\in\Var_{\mathcal{L}}}$.
	Then for any valuation $v$ in {\em\alg{A}}, there is a substitution $\sigma$ such that $v=v_{0}\circ\sigma$. 
\end{prop}
\begin{proof}
	Let $v:p\mapsto a_p$ be a valuation in \alg{A}. In virtue of Proposition~\ref{P:subalgebra-generated-by-X}, for each $p\in\Var_{\mathcal{L}}$, there is an $\Lan$-formula, $\alpha_{p}$, such that
	$a_p=\alpha_{p}[p_i\backslash v_0[p_j],\ldots,q_i\backslash v_0[q_J]]$, where $p_i,\ldots,q_i$ are all variables occurring in $\alpha_{p}$. Without loss of generality , one can count that $p_i=p_j,\ldots,q_i=q_j$. Hence $a_p=v_0[\alpha_{p}]$, that is $v[p]=v_0[\alpha_{p}]$, for any $p\in\Var_{\mathcal{L}}$.
	Now we define a substitution as follows:  $\sigma(p)=\alpha_p$, where $p\in\Var_{\mathcal{L}}$. Thus, we have: $v[p]=v_0(\sigma(p))$, that is $v=v_{0}\circ\sigma$.
\end{proof}

\begin{prop}\label{P:valuation-in-F_L/theta}
	Let $\theta$ be a congruence on $\FormAl$. Then for valuation $v$ in ${\FormAl\slash\theta}$, there is an $\Lan$-substitution $\sigma$ such that for any $\Lan$-formula $\alpha$, $v[\alpha]=\sigma(\alpha)\slash\theta$.
	Conversely, for any $\Lan$-substitution $\sigma$, there is a unique valuation $v_\sigma$ in $\FormAl\slash\theta$ such that for any formula $\alpha$, $v_\sigma[\alpha]=\sigma(\alpha)\slash\theta$.
\end{prop}
\begin{proof}
	Suppose according to $v$, $v:p\mapsto \gamma\slash\theta$, where $p\in\Var_{\mathcal{L}}$. Then we define: $\sigma: p\mapsto\gamma$, where $p\in\Var_{\mathcal{L}}$.\footnote{Of course, the definition of such a substitution is not unique, and the axiom of choice is required here.} Thus, in virtue of Proposition~\ref{P:substitution-as-homomorphism}, we have: $v[\alpha]=\sigma(\alpha)\slash\theta$, for any $\Lan$-formula $\alpha$.
	
	The converse-part is obvious.
\end{proof}

Now we turn to a central notion of any semantics --- the notion of\textit{ validity}.\footnote{W.~and M.~Kneale note, ``[\dots] logic is not simply valid argument but the reflection upon principles of validity;''~\cite{kneales1962}, section I.1.} This notion in turn is defined in the framework of the notion of logical matrix.

\begin{defn}[logical matrix of type $\Lan$]\label{D:matrix-type-L}\index{logical matrix}\index{matrix!logical}
	A \textbf{logical matrix} $($or simply a \textbf{matrix}$)$ of type $\Lan$ is a structure {\em$\mat{M}=\langle\alg{A},D(x)\rangle$},  where {\em\alg{A}} is an algebra of type $\Lan$, which, if considered as part of {\em\mat{M}}, is called an {\em\mat{M}}-\textbf{algebra}, and $D(x)$ is a predicate on {\em$|\alg{A}|$}. The matrix {\em\mat{M}} is \textbf{nontrivial} if the algebra {\em\alg{A}} is nontrivial, and \textbf{trivial} otherwise. We will also denote a matrix as a structure {\em$\langle\textsf{A};\Func,\Cons,D(x)\rangle$} or as a structure
	{\em$\langle\textsf{A};\Func,D(x)\rangle$}, if $\Cons=\varnothing$. We will also denote a matrix as {\em$\langle\alg{A},D\rangle$} or as {\em$\langle\textsf{A};\Func,D\rangle$}
	or as {\em$\langle\textsf{A};\Func,D\rangle$}, if $\Cons=\varnothing$, where {\em$D\subseteq|\alg{A}|$}. In the last three notations $D$ is called  the set of \textbf{designated elements} or a \textbf{logical filter} $($or simply a \textbf{filter}$)$ in $($or of$)$ {\em\mat{M}}. 
\end{defn}

We should remark that using ``logical'' in the last definition was a bit hasty, for the first connection with logic emerges, at least in the case of logical matrices, in the next definition
when we separate semantically valid $\Lan$-formulas from the other formulas.
\begin{defn}\label{D:validity-1}\index{formula!valid}
	Let {\em$\alg{M}=\langle\alg{A},D\rangle$} be a logical matrix of type $\Lan$. Then a formula $\alpha$ is \textbf{satisfied} by a valuation $v$ in {\em\alg{A}} if $v[\alpha]\in D$, in which case we say that $v$ \textbf{satisfies} $\alpha$ in {\em\mat{M}}{\em;} further, $\alpha$ is \textbf{semantically valid} $($or simply \textbf{valid}$)$ in {\em\alg{M}}, symbolically {\em$\alg{M}\models\alpha$} or {\em$\models_{\textbf{M}}\alpha$}, if for any valuation $v$ in {\em\alg{A}},
	$v[\alpha]\in D$. Given a set $X$ of formulas, we write {\em$\alg{M}\models X$} if {\em$\alg{M}\models\alpha$}, for all $\alpha\in X$. If a formula is not valid in a matrix, it is called \textbf{invalid}\index{formula!invalid} in this matrix or is \textbf{rejected} by the matrix or one can say that this matrix \textbf{rejects} the formula. If $v[\alpha]\notin D$ in {\em\alg{M}}, $v$ is a \textbf{valuation rejecting} or \textbf{refuting}\index{valuation!refuting}\index{valuation!rejecting} $\alpha$. The set of all formulas valid in a matrix {\em\alg{M}} is denoted by {\em$L\alg{M}$} and is called the \textbf{logic} of this matrix.
\end{defn}

\begin{prop}\label{P:matrix-substitution-closedness}
	Let {\em$\alg{M}=\langle\alg{A},D\rangle$} be a matrix $($of type $\Lan$$)$. Then for any substitution $\sigma$, if {\em$\alg{M}\models\alpha$}, so is {\em$\alg{M}\models\sigma(\alpha)$}.
\end{prop}
\begin{proof}
	Let $v$ be any valuation in \alg{A}. According to Propositions~\ref{P:substitution-as-homomorphism} and~\ref{P:valuation}, $\sigma:\mathfrak{F}_{\mathcal{L}}\longrightarrow\mathfrak{F}_{\mathcal{L}}$ and
	$v:\mathfrak{F}_{\mathcal{L}}\longrightarrow\alg{A}$ are homomorphisms. Then $v\circ\sigma$ is also a homomorphism $\mathfrak{F}_{\mathcal{L}}\longrightarrow\alg{A}$ and hence, by Proposition~\ref{P:valuation}, is a valuation. Therefore,
	$\alg{M}\models\alpha$ implies that $v[\sigma(\alpha)]\in D$.
\end{proof}

The last proposition induces a way of how, given a matrix, to obtain a matrix on the formula algebra so that both matrices validate the same set of formulas.

Let $\alg{M}=\langle\alg{A},D\rangle$ be a matrix. We denote
\[
D_{\textsf{M}}=\set{\alpha}{\alg{M}\models\alpha}
\]
and then define a matrix
\[
\alg{M}_{\textit{Fm}}=\langle\mathfrak{F}_{\mathcal{L}},D_{\textsf{M}}\rangle.
\]
We observe that for any formula $\alpha$,
\[
\alg{M}_{\textit{Fm}}\models\alpha\Longleftrightarrow\alg{M}\models\alpha.
\tag{\textit{Lindenbaum equivalence}}
\]

Indeed, let $\alg{M}_{\textit{Fm}}\models\alpha$. Then for any substitution $\sigma$ (which, according to Proposition~\ref{P:substitution-as-homomorphism}, is a valuation in $\mathfrak{F}_{\mathcal{L}}$), $\sigma(\alpha)\in D_{\textsf{M}}$; in particular, $\iota(\alpha)=\alpha\in D_{\textsf{M}}$, that is $\alg{M}\models\alpha$.

Now suppose that $\alg{M}\models\alpha$, that is $\alpha\in D_{\textsf{M}}$. In virtue of Proposition~\ref{P:matrix-substitution-closedness}, $D_{\textsf{M}}$ is closed under substitution. This completes the proof of the equivalence above.

Next let us consider $\alg{M}=\langle\mathfrak{F}_{\mathcal{L}},D\rangle$. We observe that if $\alpha\in D_{\textsf{M}}$, that is $\alg{M}\models\alpha$, then, since, by virtue of Proposition~\ref{P:substitution-as-homomorphism}, each substitution is a valuation in $\mathfrak{F}_{\mathcal{L}}$ and vice versa, $\alpha\in D$.

Now let $D^{\prime}\subseteq D$ and $D^{\prime}$ be closed under substitution. Assume that $\alpha\in D^{\prime}$. Then, obviously, $\alg{M}\models\alpha$ and hence $\alpha\in  D_{\textsf{M}}$.

Thus we have obtained the following.
\begin{prop}
	For any matrix {\em\alg{M}}, there is a matrix on the formula algebra, namely {\em$\alg{M}_{\textit{Fm}}$}, such that both matrices validate the same set of $\Lan$-formulas. Moreover, given a matrix {\em$\alg{M}=\langle\mathfrak{F}_{\mathcal{L}},D\rangle$}, the filter of the matrix
	{\em$\alg{M}_{\textit{Fm}}=\langle\mathfrak{F}_{\mathcal{L}},D_{\textsf{M}}\rangle$}is the largest subset of $D$ closed under substitution. Thus, if $D$ itself is closed under substitution, then {\em$D=D_{\textsf{M}}$}.
\end{prop}

The importance of the last proposition for what follows is that it shows that the matrices of the form 
$\langle\mathfrak{F}_{\mathcal{L}},D\rangle$, where a filter $D$ is closed under substitution, suffice for obtaining all sets of valid $\Lan$-formulas, which can be obtained individually  by the matrices of type $\Lan$. This observation will be clarified below as \textit{Lindenbaum's Theorem}; see Proposition~\ref{P:lindenbaum-theorem}.\\

Now we turn to some special matrices which have played their important role in the development of sentential logic. In all these matrices, the filter of designated elements consists of one element. In this connection, we note the following equivalence, which will be used throughout:
If $D=\{c\}$, then
\begin{equation}\label{E+one-element-filter}
x\in D\Longleftrightarrow x=c.
\end{equation}

In the sequel, we find useful the following property.
\begin{prop}\label{P:validity-homomorphism}
	Let {\em$\alg{M}_{1}=\langle\alg{A}_{1},D_{1}\rangle$} and {\em$\alg{M}_{2}=\langle\alg{A}_{2},D_{2}\rangle$} be matrices of the same type.
	Also, let $h$ be an epimorphism of {\em$\alg{A}_{1}$} onto {\em$\alg{A}_{2}$} such that
	$h[D_1]\subseteq D_2$. Then {\em$L\alg{M}_{1}\subseteq L\alg{M}_{2}$}
\end{prop}
\begin{proof}
	Let $v$ be a valuation in $\alg{A}_2$. Since $h$ is an epimorphisms, for each variable $p_{\gamma}$, there is an element $x_{\gamma}\in|\alg{A}_1|$ such that $h(x_\gamma)=v[p_\gamma]$. We define a valuation $v^{\prime}$ in $\alg{A}_1$ as follows:~$v^{\prime}[p_\gamma]=x_\gamma$, for any variable $p_\gamma$. It is obvious that $h(v^{\prime}[p_\gamma])=v[p_\gamma]$. Assume that $\alpha\in L\alg{M}_{1}$; in particular, $v^{\prime}[\alpha]\in D_1$. This implies that $v[\alpha]\in D_2$.
\end{proof}

\subsection{The logic of a two-valued matrix}\label{S:two-valued}

We begin with a well-known matrix of the so-called ``two-valued logic.'' The significance of this matrix in the history of logic is difficult to overestimate. The two-valued matrix was assumed, first, to distinguish true statements from false ones; and, second, to validate only the former. As we will see below, the main logic laws such as, for example, the law of identity, the law of contradiction, the law of the excluded middle (see below) and others can be demonstrated as valid formulas in this matrix. For this reason, it is worth considering this matrix first. Later on, we will explain another significance of this matrix, purely algebraic in nature, which leads to a very important view on this matrix as a refuting and separating tools. \\

The two-valued matrix assumes two main variants of formulation --- as a matrix of type $\Lan_{A}$ and as that of type $\Lan_{B}$. As will be seen below, the difference between these variants concerns only semantic rules and does not affect essentially the set of valid formulas in the corresponding matrices.
For this reason, we denote both variants by $\booleTwo$.\index{$\booleTwo$} The elements of $\booleTwo$ are denoted by $\zero$ (\textit{false}) and $\one$ (\textit{true}), which constitute a Boolean algebra with $\zero$ as the least element and $\one$ as the greatest with respect to the partial ordering $\zero\leq \one$. According to the intended interpretation, $\one$ represents a truth value \textit{true}, $\zero$ represents \textit{false}. Being viewed as a 2-element Boolean algebra of type $\Lan_{A}$, the connectives $\wedge$, $\vee$ and $\neg$ are interpreted as meet, join and complementation, respectively, while
\begin{equation}\label{E:material-implication}
x\rightarrow y:=\neg x\vee y.\footnote{We remind the reader that, in terms where the complementation (or pseudo-complementation) occurs, it precedes both meet and join when we interpret such terms in the lattices with complementation or pseudo-complementation.}
\end{equation}
This leads to the familiar truth tables:
\begin{center}
	\begin{tabular}{c|cc} \hline
		&\multicolumn{2}{|c}{$x\wedge y$}\\
		$x\backslash y$ &$\zero$  &$\one$\\ \hline
		$\zero$ &$\zero$
		&$\zero$ \\
		$\one$ &$\zero$
		&$\one$\\ \hline
	\end{tabular}
	\quad\begin{tabular}{c|cc} \hline
		&\multicolumn{2}{|c}{$x\vee y$}\\
		$x\backslash y$ &$\zero$  &$\one$\\ \hline
		$\zero$ &$\zero$ &$\one$\\
		$\one$ &$\one$
		&$\one$ \\ \hline
	\end{tabular}
	\quad\begin{tabular}{c|c} \hline
		$x$ &$\neg x$\\ \hline
		$\zero$ &~~$\one$\\
		$\one$ &~~$\zero$\\ \hline
	\end{tabular}
	\quad\begin{tabular}{c|cc} \hline
		&\multicolumn{2}{|c}{$x\rightarrow y$}\\
		$x\backslash y$ &$\zero$  &$\one$\\ \hline
		$\zero$ &$\one$ &$\one$\\
		$\one$ &$\zero$
		&$\one$ \\ \hline
	\end{tabular}
\end{center}

We observe that for any $x,y\in|\booleTwo|$,
\begin{equation}\label{E:implication}
\begin{rcases*}
\begin{array}{c}
x\rightarrow y=\one\Longleftrightarrow x\leq y;\\
x\rightarrow\zero=\neg x;\\
\one\rightarrow x=x.
\end{array}
\end{rcases*}
\end{equation}

If one regards $\booleTwo$ (more exactly, its counterpart) as an algebra of type $\Lan_{B}$, two nullary operations, $\zero$ and $\one$, should be added to the signature, along with the intended interpretation of $\bot$ and $\top$ as $\zero$ and $\one$, respectively. Making of $\booleTwo$ a logical matrix, we add a logical filter $\{\one\}$. Thus a formula $\alpha$ of type $\Lan_{A}$ (or of type $\Lan_{B}$) is valid in $\booleTwo$ if for any valuation $v$,
\[
v[\alpha]=\one.
\] 
We denote both, the algebra and the matrix, grounded on this algebra, by $\booleTwo$. It must be clear that any $\Lan_{A}$-formula is valid in the matrix $\booleTwo$ of type $\Lan_{B}$ if and only if it is valid in $\booleTwo$ of type $\Lan_{A}$.

A formula $\alpha$ valid in $\booleTwo$ is called a \textit{\textbf{classical tautology}}\index{classical tautology}. It is not difficult to see that the following formulas are classical tautologies:
\[
\begin{array}{cl}
p\rightarrow p, &(\textit{law of identity})\\
\neg(p\wedge\neg p), &(\textit{law of contradiction})\\
p\vee\neg p; &(\textit{law of the excluded middle})
\end{array}
\]
as well as the following:
\[\label{E:paradoxes-implicaion}
\begin{array}{cl}
\begin{rcases*}
p\rightarrow(q\rightarrow p);\\
\neg p\rightarrow(p\rightarrow q).
\end{rcases*} &(\text{``paradoxes of material implication''})
\end{array}
\]
The above tautologies had been known from ancient times. In modern times, such tautologies as the following have been in focus:
\[
\begin{array}{cl}
((p\rightarrow q)\rightarrow p)\rightarrow p; &(\textit{Peirce's law})\\
\neg\neg p\leftrightarrow p; &(\textit{double negation law})\\
\begin{rcases}
\neg(p\wedge q)\leftrightarrow\neg p\vee\neg q;\\
\neg(p\vee q)\leftrightarrow\neg p\wedge\neg q.
\end{rcases} &(\textit{De Morgan's laws})
\end{array}
\]

By the \textit{\textbf{two-valued logic}}\index{logic!two-valued}, one usually means the set of all classical tautologies.

\subsection{The {\L}ukasiewicz logic of a three-valued matrix}\label{S:lukasiewicz}\index{logic!{\L}ukasiewicz}

Perhaps, the ``{\L}ukasiewicz 3-valued logic'' is the best-known of all multiple-valued systems alternative to the set of classical tautologies.
The matrix which defines this ``logic'' has three values, $\zero$ (\textit{false}), $\tb$ (\textit{indeterminate}) and
$\one$ (\textit{true}), which are arranged by a linear ordering as follows: $\zero\leq\tb\leq\one$. On this carrier, an algebra of type $\Lan_{A}$ is defined in such a way that $\wedge$ and $\vee$ are interpreted as meet and join in the 3-element distributive lattice; then, $\neg$ and $\rightarrow$ are defined according to the following tables:

\begin{center}
	\begin{tabular}{c|c} \hline
		$x$ &$\neg x$\\ \hline
		$\zero$ &~~$\one$\\
		$\tb$  &~~$\tb$\\
		$\one$ &~~$\zero$\\ \hline
	\end{tabular}
	\quad\begin{tabular}{c|ccc} \hline
		&\multicolumn{3}{|c}{$x\rightarrow y$}\\
		$x\backslash y$ &$\zero$  &$\tb$ &$\one$\\ \hline
		$\zero$ &$\one$ &$\one$ &$\one$\\
		$\tb$  &$\tb$ &$\one$ &$\one$\\
		$\one$ &$\zero$ &$\tb$	&$\one$ \\ \hline
	\end{tabular}
\end{center}

We denote this algebra by $\lukasThree$.
Actually, when we have defined $\neg$ in $\lukasThree$ as above, the operation $\rightarrow$ can be defined according to \eqref{E:implication}. Moreover, the following identities can be easily checked:
\begin{equation}\label{E:lukasiewicz}
\begin{array}{c}
\begin{rcases*}
x\vee y=(x\rightarrow y)\rightarrow y;\\
x\wedge y=\neg(\neg x\vee\neg y).
\end{rcases*}
\end{array}
\end{equation}
(Exercise~\ref{section:semantics}.\ref{EX:lukasiewicz}.)
Also, we note that, as algebras of type $\Lan_{A}$, $\booleTwo$ is a subalgebra $\lukasThree$. 
(See Exercise~\ref{section:semantics}.\ref{EX:subalgebra}.)

Now taking $\{\one\}$ as a logical filter, we obtain the {\L}ukaciewicz 3-valued matrix which is also denoted by $\lukasThree$.\index{$\lukasThree$}

We note that the law of the excluded middle and Peirce's law, which are tautologies of the two-valued logic, fail to be tautologous in $\lukasThree$. (Exercise~\ref{section:semantics}.\ref{EX:failure}.)
The law of contradiction is also invalid in $\lukasThree$, for the set $\{\tb\}\subseteq|\lukasThree|$ is closed under $\wedge$ and $\neg$.

By the {\L}ukasiewicz 3-valued logic one means the set of formulas which are validated by $\lukasThree$. These formulas are called $\lukasThree$-\textit{\textbf{tautologies}}\index{tautology}. It should be clear that the set 
$\lukasThree$-tautologies is a proper subset of the set of classical tautologies.
(Exercise~\ref{section:semantics}.\ref{EX:luk-three-included-boole-two}.) One time, it had been assumed that the matrix $\langle\lukasThree,\{\tb,\one\}\rangle$ validates the classical tautologies and only them. If the first part of the last statement is obviously true, the second is false. The following counterexample is due to A.~R. Turquette: the formula $\neg(p\rightarrow\neg p)\vee\neg(\neg p\rightarrow p)$ is a classical tautology but is refuted in the last matrix by the assignment $v[p]=\tb$.

\subsection{The {\L}ukasiewicz modal logic of a three-valued matrix}\label{S:lukasiewicz-modal}

Both the algebra and matrix $\lukasThree$ can be expanded to type $\Lan_{C}$ as follows:
\begin{center}
	\begin{tabular}{c|c|c} \hline
		$x$ &$\Box x$ &$\Diamond x$\\ \hline
		$\zero$ &~$\zero$ &~$\zero$\\
		$\tb$ &~$\zero$ &$\one$\\
		$\one$ &~$\one$ &$\one$\\ \hline
	\end{tabular}
\end{center}

The last definitions of the modalities $\Box$ and $\Diamond$ are not independent of the other assertoric connectives, for in the expanded algebra $\lukasThree$ the following identities hold:
\begin{equation}\label{E:lukasiewicz-modal}
\begin{array}{c}
\begin{rcases*}
\Diamond x=\neg x\rightarrow x;\\
\Box x=\neg\Diamond\neg x.
\end{rcases*}
\end{array}
\end{equation}
We leave for the reader to check this. (See Exercise~\ref{section:semantics}.\ref{EX:lukasiewicz-modal}.)

We observe that the formulas
\[
\Box p\rightarrow p,~p\rightarrow\Diamond p~\text{and}~\Box p\rightarrow\Diamond p
\]
are valid in the expanded matrix, while their converses,
\[
p\rightarrow\Box p,~\Diamond p\rightarrow p~\text{and}~\Diamond p\rightarrow\Box p,
\]
are rejected by it. (See Exercise~\ref{section:semantics}.\ref{EX:lukasiewicz-modal-2}.)

\subsection{The G\"{o}del $n$-valued logics}\label{section:goedel}
As in the case $\lukasThree$, let us consider the three distinct ``truth values,'' $\zero$, $\tb$ and
$\one$, arranged by a linear ordering as follows: $\zero\leq\tb\leq\one$. Regarding this chain as a Heyting algebra with $\wedge$, $\vee$, $\rightarrow$ and $\neg$ as meet, join, relative pseudo-complementation and pseudo-complementation, respectively, we get an algebra $\godelThree$ of type $\Lan_{A}$. Thus we arrive at the following truth tables:
\begin{center}
	\begin{tabular}{c|ccc} \hline
		&\multicolumn{3}{|c}{$x\wedge y$}\\
		$x\backslash y$ &$\zero$ &$\tb$  &$\one$\\ \hline
		$\zero$ &$\zero$
		&$\zero$ &$\zero$\\
		$\tb$ &$\zero$ &$\tb$ &$\tb$\\
		$\one$ &$\zero$ &$\tb$
		&$\one$\\ \hline
	\end{tabular}
	\quad\begin{tabular}{c|ccc} \hline
		&\multicolumn{3}{|c}{$x\vee y$}\\
		$x\backslash y$ &$\zero$ &$\tb$ &$\one$\\ \hline
		$\zero$ &$\zero$ &$\tb$ &$\one$\\
		$\tb$ &$\tb$ &$\tb$ &$\one$\\
		$\one$ &$\one$
		&$\one$ &$\one$\\ \hline
	\end{tabular}
	\quad\begin{tabular}{c|c} \hline
		$x$ &$\neg x$\\ \hline
		$\zero$ &~~$\one$\\
		$\tb$ &~~$\zero$\\
		$\one$ &~~$\zero$\\ \hline
	\end{tabular}
	\quad\begin{tabular}{c|ccc} \hline
		&\multicolumn{3}{|c}{$x\rightarrow y$}\\
		$x\backslash y$ &$\zero$ &$\tb$ &$\one$\\ \hline
		$\zero$ &$\one$ &$\one$ &$\one$\\
		$\tb$ &$\zero$ &$\one$ &$\one$\\
		$\one$ &$\zero$ &$\tb$
		&$\one$ \\ \hline
	\end{tabular}
\end{center}

We notice that all properties \eqref{E:implication} are satisfied in $\godelThree$.\index{$\godelThree$}

Taking into account the relation $\leq$, the above truth tables can be read through the following equalities.
\begin{equation}\label{E:LC-algebra-operations}
\begin{rcases}
\begin{array}{l}
x\wedge y=\min(x,y)\\
x\vee y=\max(x,y)\\
x\rightarrow y=\begin{cases}
\begin{array}{cl}
\one &\text{if $x\le y$}\\
y &\text{otherwise}
\end{array}
\end{cases}\\
\neg x=\begin{cases}
\begin{array}{cl}
\zero &\text{if $x\neq\zero$}\\
\one &\text{otherwise},
\end{array}
\end{cases}
\end{array}
\end{rcases}
\end{equation}

Taking $\{\one\}$ as a logical filter, we obtain a logical matrix $\godelThree$.
The last matrix can be generalized if instead of 3 elements we take 4 and more elements and use \eqref{E:LC-algebra-operations} to define the same signature operations as in $\godelThree$. Thus we obtain G\"{o}del's $n$-valued logics $\textbf{G}_n$,\index{$\textbf{G}_n$} where $n\geq 3$. We note that \eqref{E:LC-algebra-operations} are also valid in $\booleTwo$; therefore, it is sometimes convenient to rename  the algebra $\booleTwo$ by $\textbf{G}_2$.

The matrices $\textbf{G}_n$, $n\ge 2$, have been in particular important in investigation of the lattice of intermediate logics, that is the logics the valid formulas of which contain all \textit{intuitionistic tautologies} and are included into the set of classical tautologies. In particular, all the matrices $\textbf{G}_n$, $n\ge 3$, fail Peirce's law, the double negation law, as well as both De Morgan's laws. (Exercise~\ref{section:semantics}.\ref{EX:Peirce-double-negation}.)

\subsection{The Dummett denumerable matrix}\label{section:dummett}
An denumerable generalization of $\textbf{G}_n$, $n\ge 2$, can be obtained if we take take a denumerable set and arrange it by a linear ordering $\le$ according to type $1+\omega^{\ast}$. We denote the least element by $\zero$ and the greatest by $\one$. (We note: the ordinal notation of the former is $\omega$ and of the latter is $0$.) Then, introducing the operations $\wedge$, $\vee$, $\rightarrow$ and $\neg$ as meet, join, relative pseudo-complementation and pseudo-complementation, respectively, of a Heyting algebra and equipping the latter algebra with a logical filter $\{\one\}$, we receive the \textit{Dummett denumerable matrix} \mat{LC}. As for $\godelThree$, the algebra of the Dummett matrix satisfies the properties~\eqref{E:implication} and also~\eqref{E:LC-algebra-operations}. We note that \mat{LC} validates the formula
\[
(p\rightarrow q)\vee(q\rightarrow p).
\]

The set of all formulas valid in this matrix, $L\mat{LC}$, is called the \textit{\textbf{Dummett logic}}.\index{logic!Dummett} 

M. Dummett shows in~\cite{dummett1959}, theorem 2, that
\begin{equation}\label{E:LC-characterization}
L\mat{LC}=\bigcap_{n\ge 2}L\mat{G}_{n}.
\end{equation}

In the sequel, we will need another logical matrix whose logic coincides with the Dummett logic $\LC$. 

Let us consider a denumerable set arranged linearly according to type $\omega+1$.
We treat this structure as a Heyting algebra with a least element $\zero$, greatest element $\one$ and operations $\wedge$, $\vee$, $\rightarrow$ and $\neg$ that are defined according to~\eqref{E:LC-algebra-operations}. Announcing $\lbrace\one\rbrace$ a logical filter, we define the logical matrix $\mat{LC}^{\ast}$.
\begin{prop}\label{P:LC=LC-ast}
	{\em$L\mat{LC}=L\mat{LC}^{\ast}$}.
\end{prop}
\begin{proof}
	We show that
	\[
	L\mat{LC}^{\ast}=\bigcap_{n\ge 2}L\mat{G}_{n}.
	\]
	It is clear that for any $n\ge 2$, the algebra of $\mat{G}_n$ is a homomorphic image of the algebra of $\mat{LC}^{\ast}$; see Section~\ref{section:heyting-algebra}. In virtue of Proposition~\ref{P:validity-homomorphism}, $L\mat{LC}^{\ast}\subseteq\bigcap_{n\ge 2}L\mat{G}_{n}$.
	
	Now, assume that a valuation $v$ in $\mat{LC}^{\ast}$ refutes a formula $\alpha$, that is $\alpha\notin L\mat{LC}^{\ast}$. Let $\beta_1,\ldots,\beta_n$ be all subformulas of $\alpha$. Then, according to~\eqref{E:LC-algebra-operations}, the set $\lbrace v[\beta_1],\ldots,v[\beta_n],\zero,\one\rbrace$ is the carrier of a subalgebra of the $\mat{LC}^{\ast}$-algebra, which in turn is isomorphic to a $\mat{G}_{n}$-algebra, for some $n\ge 2$. This implies that $\alpha\notin\bigcap_{n\ge 2}L\mat{G}_{n}$. 
\end{proof}

\subsection{Two infinite generalizations of {\L}ukasiewicz's $\lukasThree$}

The {\L}ukasiewicz algebra $\lukasThree$ can be easily generalized if we take for a carrier the interval $[0,1]$ or the rational scale $\mathbb{Q}$ or on the real scale $\mathbb{R}$ arranged in both cases by the ordinary less-than-or-equal relation $\le$. The set of designated elements in both cases is $\{1\}$.
The operations $\neg$ and $\rightarrow$ are defined as follows:
\[
\begin{array}{l}
\neg x:=1-x,\\
x\rightarrow y:=\min(1,1-x+y);
\end{array}
\] 
then, the other operations are defined according to \eqref{E:lukasiewicz}. This leads to the identities
\[
\begin{array}{c}
x\wedge y=\min(x,y)~\text{and}~x\vee y=\max(x,y),
\end{array}
\]
which can be taken as definitions of $\wedge$ and $\vee$ as well.

The set of formulas valid in the first matrix, that is when $[0,1]$ is considered in the rational scale, is denoted by $\mat{\L}_{\aleph_{0}}$ and the formulas valid in the second matrix by $\mat{\L}_{\aleph_{1}}$.\\

So far, we have discussed the semantic rules for the language $\Lan$. As we introduce other schematic languages later on (in Chapters~\ref{chapter:equational-con}), we will have to adjust these rules for them. However, although semantic rules allow us to distinguish valid formulas from refuted ones, they, taken alone, do not suffice to
analyze reasoning. This echoes the remark due to Aristotle that ``it is absurd to discuss refutation (\textgreek{>'elegqoc}) without discussing reasoning (\textgreek{sullogism'os}).''\footnote{Cf.~\cite{kneales1962}, chapter I, section 4.} We discuss the relation of formal deduction between $\Lan$-formulas  in the form of \textit{consequence relation}. This is a topic of the next chapter.\\

\noindent{\textbf{Exercises~\ref{section:semantics}}}
\begin{enumerate}
	\item\label{EX:endomorphism}Show that any endomorphism $f:\mathfrak{F}_{\mathcal{L}}\longrightarrow\mathfrak{F}_{\mathcal{L}}$ can be regarded as a substitution; that is $f$ is a unique extension of $v:\Var\longrightarrow\mathfrak{F}_{\mathcal{L}}$ such that $v(p)=f(p)$, for any $p\in\Var$.
	\item Prove the second part of Proposition~\ref{P:valuation}.
	\item Prove that the laws of identity, of contradiction, of the excluded middle, as well as Peirce's and double negation laws are tautologies in $\booleTwo$.
	\item Prove the properties~\eqref{E:implication}.
	\item Show that for any two $\Lan_{A}$-formulas $\alpha$ and $\beta$, $\alpha\leftrightarrow\beta$ is a classical tautology if and only if for any $\booleTwo$-valuation $v$, $v[\alpha]=v[\beta]$.
	\item\label{EX:lukasiewicz}Prove the identities \eqref{E:lukasiewicz}.
	\item\label{EX:subalgebra} Prove that $\booleTwo$ is a subalgebra of $\lukasThree$.
	\item\label{EX:failure} Show that the law of the excluded middle and Peirce's law are not tautologous in $\lukasThree$.
	\item\label{EX:luk-three-included-boole-two}Prove that the set of $\lukasThree$-tautologies is properly included in the set of classical tautologies. 
	\item Prove that the ``paradoxes of material implications'' (p.~\pageref{E:implication}) are $\lukasThree$-tautologous.
	\item\label{EX:wajsberg-axioms} Show that the following formulas are $\lukasThree$-tautologies.
	\[
	\begin{array}{cl}
	(a) &p\rightarrow(q\rightarrow p);\\
	(b) &(p\rightarrow q)\rightarrow((q\rightarrow r)\rightarrow(p\rightarrow r));\\
	(c) &(\neg p\rightarrow\neg q)\rightarrow(q\rightarrow p);\\
	(d) &((p\rightarrow\neg p)\rightarrow p)\rightarrow p.
	\end{array}
	\]
	\item\label{EX:lukasiewicz-modal}Prove that the identities~\eqref{E:lukasiewicz-modal} hold the algebra $\lukasThree$ expanded by the modalities $\Box$ and $\Diamond$.
	\item\label{EX:lukasiewicz-modal-2}Show that the logical matrix of type $\Lan_{C}$ of Section~\ref{S:lukasiewicz-modal} validates the formulas $\Box p\rightarrow p,~p\rightarrow\Diamond p~\text{and}~\Box p\rightarrow\Diamond p$ and rejects their converses $p\rightarrow\Box p,~\Diamond p\rightarrow p~\text{and}~\Diamond p\rightarrow\Box p$.
	\item\label{EX:Peirce-double-negation}Show that Peirce's law, the double negation law and both De Morgan's laws can be refuted in $\godelThree$.
	\item Let us designate on the algebra $\godelThree$ the logical filter $\{\tb,\one\}$. We denote the new matrix by $\godelThree^{\prime}$. Show that $\godelThree^{\prime}$ validates a formula $\alpha$ if and only if $\alpha$ is a classical tautology.
\end{enumerate}

\section{Historical notes}\label{languages-historical-notes}
The idea that logic, as the science of reasoning, deals with the universal rather than the particular can be traced back to Plato and Aristotle. According to the Platonic doctrine of Forms, a single Form is associated with each group of things to which a common name can be applied. Plato believed that each Form differs from the particulars it represents but, at the same time, is a particular of its own kind.
Applying this doctrine to reasoning, he believed that correct reasoning is based on the connections of the involved Forms. This is how W. Kneale and M. Kneale summarize Plato's view on reasoning.
\begin{quote}
	``For Plato necessary connexions hold between Forms, and inference is presumably valid when we follow in thought the connexions between Forms as they are. How is this view related to the theory which assigns truth and falsity to sentences made up of nouns and verbs? Plato seems to hold that a sentence is true if the arrangement of its parts reflects or corresponds to connexion between Forms.''
	\cite{kneales1962}, section I.5.
\end{quote}

However, what Plato wrote about correct thinking can be regarded as philosophy of logic rather than logic itself. The first construction of formal logic we encounter in Aristotle. Although Aristotle was influenced by Plato's doctrine of Forms, he rejected it. Nevertheless, we dare to assert that Aristotle's doctrine of the syllogism is \textit{formal} in the Platonic sense. Namely Aristotle was the first to introduce variables into logic. In his influential book on Aristotle's syllogistic J. {\L}ukasiewicz makes the following comment about this fact.
\begin{quote}
	``It seems that Aristotle regarded his invention as entirely plain and requiring no explanation, for there is nowhere in his logical works any mention of variables. It was Alexander who first said explicitly that Aristotle presents his doctrine in letters, \textgreek{stoiqe{\^{i}}a}, in order to show that we get the conclusion not in consequence of the matter of the premisses, but in consequence of their form and combination; the letters are marks of universality and show that such a conclusion will follow always and for any term we may choose.'' \cite{lukasiewicz1951}, {\S} 4.
\end{quote}

Although in the course of the development of logic, when after Aristotle the generalized nature of formal judgments was reflected in the use of letters, and the latter were perceived, at least gradually, as sentential variables,\footnote{It is appropriate to quote J. {\L}ukasiewicz, who wrote, ``[\dots] the Stoics prepared the way for formalism, and they cannot be credited highly enough for that. They held strictly to \emph{words} and not to their \emph{meanings}, which is the principle requirement of formalization, and they even did so in conscious opposition to the Peripatetics.''~\cite{lukasiewicz1934}} and therefore it would be natural to consider the uniform substitution, as defined above, as a legitimate \textit{operation} that would apply to all formal judgments, however, the uniform substitution was first proposed as an \textit{inference rule}, that is, only in conjunction with (and for constructing) logical inferences. A brief history of how the rule of  uniform substitution has paved its way to existence will be discussed in more detail in Section~\ref{section:consequence-historical-notes}.

Nonetheless, the realization, though implicit and at least for propositional languages such as $\Lan$-languages above, of that the set of formulas is closed under substitution and that the notion of validity should be defined in such a way to maintain this closedness for the formulas identified as valid one can find in George Boole. As to the first part of this observation, Boole puts it as follows.
\begin{quote}
	``Let us imagine any known or existing language freed from idioms and divested of superfluity and let us express in that language any given proposition in a manner the most simple and literal [\ldots] The transition from such a language to the notation of analysis would consist of no more than the substitution of one set of signs for another, without essential change either in form or character.'' \cite{boole1854}, chapter 11, section 15.
\end{quote}
As to the second, he writes:
\begin{quote}
	``I shall in some instances modify the premises by [\ldots] substitution of some new proposition, and shall determine how by such change the ultimate conclusions are affected.'' \cite{boole1854}, chapter 13, section 1.
\end{quote}

It appears that this idea has become fully realized by the 1920s in the Lvov-Warsaw School, for without this, Lindenbaum's Theorem (Proposition~\ref{P:lindenbaum-theorem}) would be impossible. Although we can only guess, but it seems plausible to suggest that the confusion between the concepts of predicate and propositional term was an obstacle for such a long delay.\\

With merely some reservation, one can say that George Boole was the first who interpreted
symbolic expressions representing ``laws of thought'' as algebraic terms. However,
his algebraic treatment of logic was merely the beginning of a long way. This is how T. Hailperin estimates Boole's contribution.   
\begin{quote}
	``Boole wrote his \textit{Laws of Thought} before the notion of an abstract formal system, expressed within a precise language, was fully developed. And, at that time, still far in the future was the contemporary view which makes a clear distinction between a formal system and its realization or models. Boole's work contributed to bringing these ideas to fruition.'' \cite{hailperin1986}, chapter 1, {\S} 1.0.
\end{quote}

We find interpretation of formal judgments in an abstract logical matrix, namely $\booleTwo$, in~\textit{Tractatus}; cf.~\cite{wittgenstein2001}, 4.441 and further. Also, Wittgenstein used the term \textit{tautology}, applying it to any classical tautology. Later on, logical matrices turned to be a natural way to define numerous deductive systems. We owe to Alfred Tarski a systematic approach to deductive systems through logical matrices. (See more about this in Section~\ref{section:consequence-historical-notes}.)  

At first matrices were employed in the search for an independent axiomatic system for the set of classical tautologies in~\cite{bernays1926} and~\cite{lukasiewicz1929}, section 6; in particular, Bernays coined the term a ``designated element.''\footnote{Cf.~\cite{bernays1926}, p. 316.}
Also, J. C. C. McKinsey used in~\cite{mckinsey1939} matrices to prove independence of logical connectives in intuitionistic propositional logic.\\

The idea of the use of a 3-valued matrix was expressed by J. {\L}ukasiewicz in his  inaugural address as Chancellor of the University of Warsaw in 1922. Later on, a revised version of this speech was published as~\cite{lukasiewicz1961}, where one can read:
\begin{quote}
	``I maintain that there are propositions which are neither true nor false but \textit{indeterminate}. All sentences about future facts which are not yet decided belong to this category. Such sentences are neither true at the present moment, for they have no real correlate, nor are they false, for their denials too have no real correlate. If we make use of philosophical terminology which is not particularly clear, we could say that ontologically there corresponds to these sentences neither being nor non-being but possibility. Indeterminate sentences, which ontologically have possibility as their correlate, take the third truth-value.'' \cite{luk70}, p. 126.
\end{quote}

Technically, however, {\L}ukasiewicz proposed $\lukasThree$ already in 1920; cf.~\cite{lukasiewicz1920}. In~\cite{lukasiewicz1930} he expands $\lukasThree$ to include modalities and thus obtain the matrix of Section~\ref{S:lukasiewicz-modal}. In fact, he defines only $\Diamond$, but it is easy to check that the definition
\[
\Box x:=\neg\Diamond\neg x
\] 
gives operation $\Box$ as introduced in Section~\ref{S:lukasiewicz-modal} directly.

It is worth noting that Lukasiewicz's interest in "multivalued logic" (in the sense of logical matrices having more than two values) was motivated by his interest in modal contexts, and not vice versa. Indeed, he wrote:
\begin{quote}
	``I can assume without contradiction that my presence in Warsaw at a certain moment of next year, e.g. at noon on 21 December, is at the present time determined neither positively nor negatively, Hence it is \textit{possible}, but not \textit{necessary}, that I shall be present in Warsaw at the given time. On this assumption the proposition `I shall be in Warsaw at noon on 21 December of next year', can at the present time be neither true nor false. For if it were true now, my future presence in Warsaw would have to be necessary, which is contradictory to the assumption. If it were false now, on the other hand, my future presence in Warsaw would have to be impossible, which is also contradictory to the assumption. Therefore the proposition considered is at the moment \textit{neither true nor false} and must possess a third value, different from `0' or falsity and `1' or truth. This value we can designate by `$\frac{1}{2}$'. It represents `the possible', and joins `the true' and `the false' as a third value.'' \cite{lukasiewicz1930}, {\S} 6; cf. English translation in~\cite{mccall1967}, pp. 40--65.
\end{quote}

The logical matrix $\textbf{\L}_{\aleph_{1}}$ and its logic for the first time were discussed by {\L}ukasiewicz in his address delivered at a meeting of the Polish Philosophical Society in Lvov on October 14, 1922. The report of this lecture was published in {\em Ruch Filozoficzny}, vol. 7 (1923), no. 6,  pp. 92--93; see English translation in \cite{luk70}, pp. 129--130. According to this report, {\L}ukasiwicz's motivation for this matrix was as follows.
\begin{quote}
	``The application of this interpretation is twofold: 1) It can be demonstrated that if those verbal rules, or directives, which are accepted by the authors of {\em Principia Mathematica} (the rule of deduction and the rule of substitution) are adopted, then no set of the logical laws that have the numerical value 1 can yield any law that would have a lesser numerical value. 2) If 0 is interpreted as falsehood, 1 as truth, and other numbers in the interval 0--1 as the degrees of probability corresponding to various possibilities, a many-valued logic is obtained, which is an expansion of three-valued logic and differs from the latter in certain details.'' \cite{lukasiewicz1961}, p. 130.
\end{quote}

It is clear that $\textbf{\L}_{\aleph_{0}}$ is an obvious simplification of $\textbf{\L}_{\aleph_{1}}$.\\

Answering a question of Hahn, K. G\"{o}del considered~\cite{godel1932} (see English translation in~\cite{godel1986}) an infinite series of finite logical matrices and their logics, where the matrix $\booleTwo$ was the second and $\godelThree$ the third in the series. All these matrices were partially ordered chains with the maximal element designated. G\"{o}del's goal was twofold: to show that 1) the intuitionistic propositional calculus cannot be determined by a single matrix; and that 2) there are at least countably many logics between the \textit{intuitionistic} and \textit{classical logics}. (These logics will be descussed in more detail in Chapter~\ref{chapter:consequence}.)\\

The aforementioned series of G\"{o}del's matrices produces a descending chain of logics all containing the intuitionistic propositional logic. Michael Dummett showed~\cite{dummett1959} that \textbf{LC} equals the intersection of all G\"{o}del's logics. He also gave an axiomatization of \textbf{LC}, which will be discussed in Chapter~\ref{chapter:consequence}.

\chapter[Logical Consequence]{Logical Consequence}\label{chapter:consequence}	
The conception of logical consequence which we discuss in this chapter can be traced to
Aristotle's distinction between demonstrative and dialectical argument. The essence of their difference lies in the status of the premises of the argument in question. In a demonstrative argument, its premises are regarded as true. Aristotle called the argument with such premises ``didactic'' and compared them with what a teacher lays down as a starting point of the development of his subject. We would call such premises the axioms of a theory. In a dialectical argument, according to Aristotle, any assumptions can be employed for premises. More than that, in the \textit{reductio ad absurdum} argument one starts with a premise which is aimed to be proven false.

Thus, the focus in the dialectical reasoning is transferred on the argument itself. This involves considering such an argumentation as a binary relation between a set of premises and what can be obtained from them in the course of this argumentation. And the first question raised by this point of view is this: What characteristics of the dialectical argument could qualify it as sound?

\section{Consequence relations}\label{section:consequence-relation}
From the outset, we are faced with a choice between two paths each addressing one of the questions: Should consequence relation state determination between a set of premises (which can be empty) and a set of conclusions (which can also be empty), or should it do it between a set of premises (possibly the empty one) and a single formula? The first type is known as a \emph{multiple-conclusion}, or \emph{poly-conclusion}, consequence, the second --- as a \emph{single-conclusion} consequence, or simply a \emph{logical consequence}.  In this book, we discuss only the second option. As to the first, we refer the reader to~\cite{shoesmith-smiley2008}. \\

Talking about consequence relation in general, we denote it by ``$\vdash$'' and equip this sign with subscript when referring to a specific consequence relation.

In the context of consequence relation, we use the following notation:
\begin{equation}\label{E:agreement}
	\begin{rcases}
		\begin{array}{c}
			X,Y,\ldots,Z:= X\cup Y\cup\ldots\cup Z;\\
			X,\alpha,\ldots,\beta:=X\cup\{\alpha,\ldots,\beta\},
		\end{array}
	\end{rcases}
\end{equation}
for any sets $X, Y,\ldots, Z$ of $\Lan$-formulas and any $\Lan$-formulas $\alpha,\ldots,\beta$.
\begin{defn}[single-conclusion consequence relation]\label{D:consequnce-relation-single}\index{consequence relation!single conclusion}
	A relation {\em$\vdash\subseteq\mathcal{P}(\FormsL)\times\FormsL$} is an $\Lan$-\textbf{single-conclusion consequence relation} $($or simply an $\Lan$-\textbf{consequence relation} or just a \textbf{consequence relation} when  consideration occurs in a fixed, though unspecified, formal language$)$ if the following conditions are satisfied:
	{\em\[
		\begin{array}{cll}
			(\text{a}) &\alpha\in X~\textit{implies}~X\vdash\alpha; &(\textit{reflexivity})\\
			(\text{b}) &X\vdash\alpha~\textit{and}~X\subseteq Y~\textit{imply}~Y\vdash\alpha;
			&(\textit{monotonicity})\\
			(\text{c}) &\textit{if $X\vdash\beta$, for all $\beta\in Y$, and
				$Y,Z\vdash\alpha$, then $X,Z\vdash\alpha$}, &(\textit{transitivity})\\
		\end{array} 
		\]}
	for any sets $X$, $Y$  and $Z$ of $\Lan$-formulas.
\end{defn}

In case $Y=\{\beta\}$, the property  (c) reads:
\[
\begin{array}{cll}
	(\text{c}^{\ast}) &\textit{$X\vdash\beta$ and $Z,\beta\vdash\alpha$ imply $X,Z\vdash\alpha$}. &\quad\quad~~~~~~(\textit{cut})
\end{array}
\]

A consequence relation $\vdash$ is called \textit{\textbf{nontrivial}}\index{consequence relation!nontrivial} if
\[
\vdash~\subset\,\mathcal{P}(\FormsL)\times\FormsL;
\]
otherwise it is \textit{\textbf{trivial}}. \index{consequence relation!trivial}

A consequence relation $\vdash$ is called \textit{\textbf{finitary}}\index{consequence relation!finitary} if
\[
\textit{$X\vdash\alpha$ implies that there is $Y\Subset X$ such that $Y\vdash\alpha$}.
\tag{\textit{finitariness}}
\]

The last notion can be generalized as follows. Let $\kappa$ be an infinite cardinal. A consequence relation $\vdash$ is called $\kappa$-\textit{\textbf{compact}}\index{consequence relation!compact} if
\[
\textit{$X\vdash\alpha$ implies that there is $Y\subseteq X$ such that $Y\vdash\alpha$
	and $\card{Y}<\kappa$}.
\tag{$\kappa$-\textit{compactness}}
\]
Thus finitariness is simply $\aleph_{0}$-compactness.

A consequence relation $\vdash$ is \textit{\textbf{structural}}\index{consequence relation!structural} if
{\em\[
	X\vdash\alpha~\Longrightarrow~\sigma(X)\vdash\sigma(\alpha), \tag{\textit{structurality}}\\ 
	\]}
\text{for any set $X\cup\{\alpha\}\subseteq\FormsL$ and any $\Lan$-substitution $\sigma$}. $\vdash$ is \textit{\textbf{structural in a weak sense}}\index{consequence relation!structural in weak sense} if the previous implication holds for all nonempty sets $X$.\footnote{The concept of structurality in a weak sense will be used when we want to show that abstract logics obtained in a certain way have a very weak connection with structurality.}\\

We will be using the term \textit{\textbf{sequent}}\index{sequent} for a statement `$X\vdash\alpha$'.

\begin{rem}[Note on monotonicity of $\vdash$].
	{\em The property (b) of Definition~\ref{D:consequnce-relation-single} characterizes the category of \textit{monotone consequence relations} as opposed to that of \textit{non-monotone consequence relations} when the property (b) does not hold. In this book we consider only the first category of consequence relations.}
\end{rem}

Other, specific, properties in addition to (a)--(c) have been considered in literature. For instance, for the language $\Lan_{B}$ and its expansions the following property is typical for many consequence relations:
\[
X,\alpha\vdash\bot~\Longrightarrow~X\vdash\neg\alpha. \tag{\textit{reductio ad absurdum}}
\]

We will be using indices to distinguish in the same formal language either two different consequence relations or two equal consequence relation that are defined differently. The following three propositions will be used in the sequel.
\begin{prop}\label{P:con-relation-intersection}
	Let $\lbrace\vdash_i\rbrace_{i\in I}$ be a family of consequence relations in a language $\Lan$. Then the relation $\vdash$ defined by $\bigcap_{i\in I}\lbrace\vdash_i\rbrace$ is also a consequence relation. Moreover, if each $\vdash_i$ is structural, so is $\vdash$. In addition, if $I$ is finite and each $\vdash_i$ is finitary, then $\vdash$ is also finitary.
\end{prop}
\noindent\textit{Proof}~is left to the reader. (See Exercise~\ref{section:consequence-relation}.\ref{EX:con-relation-intersection}.)\\

Given a consequence relation $\vdash$ and a set $X_0$ of $\Lan$-formulas. We define a relation as follows:
\begin{equation}\label{E:relative-consequence}
	X\vdash_{X_0}\alpha\stackrel{\text{df}}{\Longleftrightarrow}X, X_0\vdash\alpha.
\end{equation}
\begin{prop}\label{P:relative-consequence}
	For any consequence relation $\vdash$ and arbitrary fixed set $X_0$, $\vdash_{X_0}$ is a consequence relation such that
	\[
	\varnothing\vdash_{X_0}\alpha~\Longleftrightarrow~ X_0\vdash\alpha.
	\]
\end{prop}
\noindent\textit{Proof}~is left to the reader. (Exercise~\ref{section:consequence-relation}.\ref{EX:relative-consequence})\\

For any consequence relation $\vdash$, we define:
\begin{equation}\label{E:nonempty-consequence}
	X\vdash^{\circ}\alpha~\stackrel{\text{df}}{\Longleftrightarrow}~X\vdash\alpha~~\text{and}~~X\neq\varnothing.
\end{equation}

\begin{prop}\label{P:nonempty-consequence}
	Let $\vdash$ be a consequence relation. Then $\vdash^{\circ}$ is also a consequence relation. In addition, if $\vdash$ is structural, $\vdash^{\circ}$ is also structural; further, if $\vdash$ is $\kappa$-compact, so is $\vdash^{\circ}$.
\end{prop}
\noindent\textit{Proof}~is left to the reader. (Exercise~\ref{section:consequence-relation}.\ref{EX:nonempty-consequence})\\

Now for any consequence relation $\vdash$, we define:
\begin{equation}\label{E:finitarized-consequence}
	X\vdash^{\bullet}\alpha~\define~\text{there is a set $X_0\Subset X$ such that}~X_0\vdash\alpha.
\end{equation}
\begin{prop}\label{P:finitarized-consequence}
	Given a consequence relation $\vdash$, the relation $\vdash^{\bullet}$ is the largest finitary consequence relations among those finitary consequence relations that are included in $\vdash$. Hence, if $\vdash$ is finitary then $\vdash=\,\vdash^{\bullet}$. In addition, if $\vdash$ is structural, then $\vdash^{\bullet}$ is also structural.
\end{prop}
\noindent\textit{Proof}~ is left to the reader. (Exercise~\ref{section:consequence-relation}.\ref{EX:finitarized-consequence})\\

As is seen from the Proposition~\ref{P:con-relation-intersection}, Definition~\ref{D:consequnce-relation-single} is essentially extensional, because it does not show how the transition from premises to conclusion can be carried out.
However, before addressing this issue, we will discuss how the concept of single-conclusion consequence relation can be represented purely in terms of sets of $\Lan$-formulas.

\paragraph{Exercises~\ref{section:consequence-relation}}
\begin{enumerate}
	\item Consider the properties (a), (b) and (c) of Definition~\ref{D:consequnce-relation-single}. Show that (a) and (c) imply (b).
	\item\label{EX:con-relation-intersection} Prove Proposition~\ref{P:con-relation-intersection}.
	\item \label{EX:relative-consequence}Prove Proposition~\ref{P:relative-consequence}.
	\item \label{EX:nonempty-consequence}Prove Proposition~\ref{P:nonempty-consequence}.
	\item\label{EX:finitarized-consequence}Prove Proposition~\ref{P:finitarized-consequence}.
\end{enumerate}

\section{Consequence operators}\label{section:consequence-operator}

Let $\vdash$ be a consequence relation for $\Lan$-formulas. Then, we define:
\begin{equation}\label{E:consequence-definition}
	\Con{X}:=\set{\alpha}{X\vdash\alpha},
\end{equation}\index{$\Con$}
for any set $X$ of $\Lan$-formulas. We observe the following connection: \begin{equation}\label{E:consequence-interconnection}
	X\vdash\alpha~\Longleftrightarrow~\alpha\in\Con{X}.
\end{equation}

The last equivalence induces the following definition.
\begin{defn}\label{D:consequence-operator}
	A map {\em$\textbf{Cn}:\mathcal{P}(\FormsL)\longrightarrow\mathcal{P}(\FormsL)$} is called a \textbf{consequence operator}\index{consequence operator} if it satisfies the following three conditions:
	{\em\[
		\begin{array}{cll}
			(\text{a}^{\dag}) &X\subseteq\Con{X}; &(\textit{reflexivity})\\
			(\text{b}^{\dag}) &X\subseteq Y~\textit{implies}~\Con{X}\subseteq\Con{Y}; &(\textit{monotonicity})\\
			(\text{c}^{\dag}) &\Con{\Con{X}}\subseteq\Con{X}, &(\textit{closedness})\\
		\end{array}
		\]}
	for any sets $X$ and $Y$ of $\Lan$-formulas. 
\end{defn}

A consequence operator \textbf{Cn} is called \textit{\textbf{finitary}}\index{consequence operator!finitary} (or \textit{\textbf{compact}}) \index{consequence operator!compact}if
\[
\begin{array}{cl}
	(\text{d}^{\dag}) &\Con{X}\subseteq\bigcup\set{\Con{Y}}{Y\Subset X}. \tag{\textit{finitariness}~\text{or}~\textit{compactness}}
\end{array}
\]

As in case of consequence relation, we generalize the last notion for any infinite cardinal $\kappa$: an consequence operator \textbf{Cn} is $\kappa$-\textit{\textbf{compact}} if
\[
\Con{X}\subseteq\bigcup\set{\Con{Y}}{Y\subseteq X~\textit{and}~\card{Y}<\kappa}. \tag{$\kappa$-\textit{compactness}}
\]

And a consequence operator \textbf{Cn} is called \textit{\textbf{structural}}\index{consequence operator!structural}, if for any set $X$ and  any substitution $\sigma$,
\[
\sigma(\Con{X})\subseteq\Con{\sigma(X)}. \tag{\textit{structurality}}
\]
We say that $\textbf{Cn}$ is \textit{\textbf{structural in a weak sense}}\index{consequence operator!structural in weak sense} if the previous inclusion is true for all nonempty $X$.

\begin{rem}[Note on monotonicity of $\textbf{Cn}$]
	{\em A consequence operator as defined in Definition~\ref{D:consequence-operator} is called \textit{monotone} as opposed to \textit{non-monotone} consequence operators.}
\end{rem}

A set $X$ of formulas is called \textit{\textbf{closed}}, or is a \textit{\textbf{theory}}\index{theory}, (with respect to $\textbf{Cn}$) if $X=\Con{X}$. This, in view of $(\text{a}^{\dag})$ and $(\text{c}^{\dag})$, any set $\Con{X}$ is a theory; we call it the \textit{\textbf{theory generated by}}\index{theory!generated by} $X$ (with respect to \textbf{Cn}).

A theory $X$ (with respect to  \textbf{Cn}) is called \textit{\textbf{finitely axiomatizable}}\index{theory!finitely axiomatizable} if there is a set $X_{0}\Subset X$ such that $X=\textbf{Cn}(X_0)$.\\

For any consequence operator \textbf{Cn}, the theory $\textbf{Cn}(\varnothing)$ is finitely axiomatizable. Indeed, it is obvious when $\textbf{Cn}(\varnothing)=\varnothing$. If $\varnothing\subset\textbf{Cn}(\varnothing)$, then for any $X$ with $\varnothing\subset X\Subset\textbf{Cn}(\varnothing)$, we have: $\textbf{Cn}(\varnothing)\subseteq\textbf{Cn}(X)\subseteq\textbf{Cn}(\varnothing)$.

On the other hand, even if every theory of a consequence operator is finitely axiomatizable, it does not guarantee that the consequence operator is finitary.

Indeed, let $\card{\Lan}=\aleph_{0}$; thus
\[
\FormsL=\lbrace\alpha_{0},\alpha_{1},\ldots\rbrace.
\]
We define:
\[
\textbf{Cn}(X):=\begin{cases}
	\begin{array}{cl}
		\FormsL &\text{if $\alpha_{0}\in X$};\\
		\FormsL &\text{if $\alpha_{0}\notin X$ and $\card{X}=\aleph_{0}$};\\
		X &\text{if $\alpha_{0}\notin X$ and $\card{X}<\aleph_{0}$}.
	\end{array}
\end{cases}
\]

It must be clear that $\textbf{Cn}(X)=X$ if, and only if, either $X=\FormsL$ or both $\alpha_{0}\notin X$ and $\card{X}<\aleph_{0}$. This implies that every theory with respect to \textbf{Cn} is finitely axiomatizable. For if $X=\FormsL$, then
$X=\textbf{Cn}(\alpha_{0})$. If $X$ is a theory with $\card{X}<\aleph_{0}$, then $X$ is obviously finitely axiomatizable.

On the other hand, the operator \textbf{Cn} is not finitary. Indeed, let $X$ be any infinite set with $\alpha_{0}\notin X$. Then $\alpha_{0}\in\textbf{Cn}(X)$, but $\alpha_{0}\notin\bigcup\set{\textbf{Cn}(Y)}{Y\Subset X}$.\\

There is a very simple and effective criterion for finitely axiomatizable theories.
\begin{prop}[Tarski's criterion]\label{P:Tarski-criterion}
	Let $X$ be a theory with respect to a consequence operator {\em\textbf{Cn}} in a language $\Lan$ with $\card{\Lan}=\aleph_{0}$. Then $X$ is finitely axiomatizable if, and only if, there is no countable-infinite family 
	$\lbrace X_i\rbrace_{i<\omega}$ of theories $($all with respect to {\em\textbf{Cn}}$)$ such that
	\[
	X_0\subseteq X_1\subseteq\ldots\subset X
	\]
	and $\bigcup_{i<\omega}X_i=X$.
\end{prop}
\begin{proof}
	The `if' part is true, even if we drop the restriction on the cardinality of $\Lan$. We leave for the reader to prove this part of the criterion.
	(Exercise~\ref{section:consequence-operator}.\ref{EX:Tarski-criterion})
	
	To prove the `only-if' part, we assume that a theory $X$ is not finitely axiomatizable. Since $\card{X}=\aleph_{0}$, we can arrange all formulas of $X$ as an infinite sequence:
	\[
	\alpha_{0}, \alpha_{1}, \ldots
	\]
	Then, we define:
	\[
	X_0:=\textbf{Cn}(\lbrace\alpha_{0}\rbrace),~X_1:=\textbf{Cn}(\lbrace\alpha_{0},\alpha_{1}\rbrace)~\text{and so on}.
	\]
	It is clear that the so obtained theories $X_i$ satisfy the conditions of the proposition.
\end{proof}

We observe the following.
\begin{prop}\label{P:intersection-con}
	Given a consequence operator {\em\textbf{Cn}}, any set {\em$\bigcap_{i\in I}\lbrace\Con{X_i}\rbrace$} is closed.
\end{prop}
\noindent\textit{Proof}~is left to the reader. (See Exercise~\ref{section:consequence-operator}.\ref{EX:intersection-con}.) Compare this proposition with Proposition~\ref{P:con-relation-intersection}.\\

The next proposition establishes connection between the notion of consequence relation and that of consequence operator.
\begin{prop}\label{P:con-operation=con-relation}
	Let $\vdash$ be a consequence relation. The operator {\em\textbf{Cn}} defined by~\eqref{E:consequence-definition} is a consequence operator. Conversely, for a given consequence operator {\em\textbf{Cn}}, the corresponding relation defined by~\eqref{E:consequence-interconnection} is a consequence relation. Moreover, {\em\textbf{Cn}} is finitary $($$\kappa$-compact$)$ if, and only if, so is $\vdash$.
	In addition, {\em\textbf{Cn}} is structural if, and only if, so is $\vdash$.
\end{prop}
\noindent\textit{Proof}~can be obtained by routine check and is left to the reader. (See Exercise~\ref{section:consequence-operator}.\ref{EX:con-operation=con-relation}.)
\begin{quote}
	{\em Proposition~\ref{P:con-operation=con-relation} allows us to use the notion of consequence relation and that of consequence operator interchangeably.}	
\end{quote}

The next observation confirms the last remark; it just mirrors Proposition~\ref{P:con-relation-intersection} 
\begin{prop}\label{P:con-operator-intersection}
	Let {\em$\lbrace\textbf{Cn}_i\rbrace_{i\in I}$} be a family of consequence operators in a language $\Lan$. Then the operator {\em$\textbf{Cn}$} defined by 
	the equality
	{\em\[
		\textbf{Cn}(X):=\bigcap_{i\in I}\textbf{Cn}_{i}(X)
		\]}is also a consequence operator. Moreover, if each {\em$\textbf{Cn}_i$} is structural, so is {\em $\textbf{Cn}$}. In addition, if $I$ is finite and each {\em$\textbf{Cn}_i$} is finitary, then {\em$\textbf{Cn}$} is also finitary.
\end{prop}
\noindent\textit{Proof}~is left to the reader. (Exercise~\ref{section:consequence-operator}.\ref{EX:con-operator-intersection})\\

Besides the properties $(\text{a}^{\dag})$--$(\text{c}^{\dag})$, finitariness and $\kappa$-compactness, other properties of a map $\textbf{Cn}:\mathcal{P}(\FormsL)\longrightarrow\mathcal{P}(\FormsL)$ related to logical consequence have been considered in literature. These are some of them.
\[
\begin{array}{cl}
	(\text{e}^{\dag}) &X\subseteq \Con{Y}~\textit{implies}~\Con{X}\subseteq\Con{Y}; \quad(\textit{cumulative transitivity})\\
	(\text{f}^{\dag}) &X\subseteq \Con{Y}~\textit{implies}~\Con{X\cup Y}\subseteq\Con{Y}; \quad(\textit{strong cumulative transitivity})\\
	(\text{g}^{\dag}) &X\subseteq Y\subseteq\Con{X}~\textit{implies}~\Con{Y}\subseteq\Con{X}; 
	\quad(\textit{weak cumulative transitivity})\\
	(\text{h}^{\dag}) &\bigcup\set{\Con{Y}}{Y\Subset X}\subseteq\Con{X};
	\quad(\textit{finitary inclusion})\\
	(\text{i}^{\dag}) &\textit{if}~X\neq\varnothing,~\textit{then}~\Con{X}\subseteq
	\bigcup\set{\Con{Y}}{Y\neq\varnothing~\mbox{and}~Y\Subset X};
	\quad(\textit{strong finitariness})\\
	(\text{j}^{\dag}) &\textit{if {\em$\alpha\notin\Con{X}$}, then there is
		a maximal $X^{\ast}$ such that $X\subseteq X^{\ast}$}\\
	&\textit{and}~\alpha\notin\Con{X^{\ast}}.
	\quad(\textit{maximalizability})
\end{array}
\]
\begin{prop}\label{P:con-connections}
	Let {\em\textbf{Cn}} be a map from {\em$\mathcal{P}(\FormsL)$} into {\em$\mathcal{P}(\FormsL)$}, Then the following implications hold$\,:$
	{\em\[
		\begin{array}{rl}
			i) &(\text{a}^{\dag})~\textit{and}~(\text{g}^{\dag})~\textit{imply}~(\text{c}^{\dag});
			\\
			ii) &(\text{b}^{\dag})~\textit{and}~(\text{c}^{\dag})~\textit{imply}~(\text{e}^{\dag});
			\\
			iii) &(\text{e}^{\dag})~\textit{implies}~(\text{g}^{\dag});
			\\
			iv) &(\text{f}^{\dag})~\textit{implies}~(\text{g}^{\dag});
			\\
			v) &(\text{a}^{\dag})~\textit{and}~(\text{g}^{\dag})~\textit{imply}~(\text{f}^{\dag});
			\\
			vi) &(\text{a}^{\dag})~\textit{and}~(\text{f}^{\dag})~\textit{imply}~(\text{c}^{\dag});
			\\
			vii) &(\text{a}^{\dag})~\textit{and}~(\text{b}^{\dag})~\textit{and}~(\text{c}^{\dag})~
			\textit{imply}~(\text{f}^{\dag});
			\\
			viii) &(\text{i}^{\dag})~\textit{implies}~(\text{d}^{\dag});
			\\
			ix) &(\text{b}^{\dag})~\textit{and}~(\text{d}^{\dag})~\textit{imply}~(\text{i}^{\dag});
			\\
			x) &(\text{b}^{\dag})~\textit{implies}~(\text{h}^{\dag});
			\\
			xi) &(\text{b}^{\dag})~\textit{and}~(\text{d}^{\dag})~\textit{imply}~(\text{j}^{\dag});
			\\
			xii) &(\text{d}^{\dag})~\textit{and}~(\text{h}^{\dag})~\textit{imply}~(\text{b}^{\dag}).
			\\
		\end{array}	
		\]}
\end{prop}
\begin{proof}
	We prove $(vi)$, $(xi)$ and $(xii)$, in this order, and leave the rest to the reader. (See Exercise~\ref{section:consequence-operator}.\ref{EX:con-connections}.)
	
	Suppose $X\subseteq\Con{Y}$, Then, by virtue of $(\text{a}^\dag)$, $X\cup Y\subseteq\Con{Y}$. According to $(\text{b}^\dag)$ and $(\text{c}^\dag)$, we have: 
	$\Con{X\cup Y}\subseteq\Con{\Con{Y}}\subseteq\Con{Y}$.
	
	To prove $(xi)$, we assume that $\alpha\notin\Con{X}$, for some set $X$.
	
	Let us denote
	\[
	C_{\alpha}:=\set{Y}{X\subseteq Y~\text{and}~\alpha\notin\Con{Y}}.
	\]
	We aim to show that $C_\alpha$ contains a maximal set.
	
	First, we observe that $C_\alpha$ is nonempty, for $X\in C_\alpha$. Regarding $C_\alpha$ as a partially ordered set with respect to $\subseteq$, we consider a chain $\lbrace Y_i\rbrace_{i\in I}$ in $C_\alpha$. That is, either $Y_i\subseteq Y_j$ or $Y_i\subseteq Y_j$, and $\alpha\notin Y_i$, for any $i\in I$.  Next we denote
	\[
	Y_0:=\bigcup_{i\in I}\lbrace Y_i\rbrace.
	\]
	We aim to show that $\alpha\notin\Con{Y_0}$. For contradiction, assume that $\alpha\in\Con{Y_0}$. Then, by $(\text{d}^{\dag})$, there is a set $Y^{\prime}\Subset Y_0$ such that $\alpha\in\Con{Y^{\prime}}$. We note that $Y^{\prime}\neq\varnothing$, for otherwise, by $(\text{b}^{\dag})$, we would have that $\alpha\in\Con{X}$. Since $Y^\prime$ is finite, there is ${i_0}\in I$ such that $Y^\prime\subseteq Y_{i_0}$. Then, by virtue of $(\text{b}^{\dag})$, $\alpha\in\Con{Y_{i_0}}$. We have obtained a contradiction. Thus $C_\alpha$ satisfies the conditions of Zorn's lemma, according to which $C_\alpha$ contains a maximal set.
	
	Finally, to prove ($xii$), we assume that $X\subseteq Y$ and $\alpha\in\Con{X}$. In virtue of $(\text{d}^\dag)$, there is a finite $Y^{\prime}\subseteq X$ such that $\alpha\in\Con{Y^\prime}$. We note that, by premise, $Y^{\prime}\Subset Y$. Then, by $(\text{h}^\dag)$, $\alpha\in\Con{Y}$.
\end{proof}

In the sequel, we employ a letter  `$\mathcal{S}$' to denote both $\vdash$ and \textbf{Cn} interrelated in the sense of~\eqref{E:consequence-interconnection}.
More specifically, we call a consequence relation $\vdash_{\mathcal{S}}$ and the corresponding consequence operator $\textbf{Cn}_{\mathcal{S}}$ an \textit{\textbf{abstract logic}} $\mathcal{S}$.\index{logic!abstract}
Thus, given an abstract logic $\mathcal{S}$, for any sets $X$ of $\Lan$-formulas and any $\Lan$-formula $\alpha$, we have:
\begin{equation}\label{E:Cn-S-connection}
	X\vdash_{\mathcal{S}}\alpha\Longleftrightarrow\alpha\in\textbf{Cn}_{\mathcal{S}}(X).
\end{equation}

Further, we denote 
\[
\bm{T}_{\mathcal{S}}:=\textbf{Cn}_{\mathcal{S}}(\varnothing)
\]
and call the elements of the latter set \textbf{\textit{theorems}}\index{theorem}, or \textit{\textbf{theses}}\index{thesese}, \textit{\textbf{of}} $\mathcal{S}$, or simply $\mathcal{S}$-\textit{\textbf{theorems}}, or $\mathcal{S}$-\textit{\textbf{theses}}.
Any set $\textbf{Cn}_{\mathcal{S}}(X)$ is call the $\mathcal{S}$-\textit{\textbf{theory generated by a set}} $X$.  Thus $\bm{T}_{\mathcal{S}}$ is the least $\mathcal{S}$-theory with respect to $\subseteq$. 

We note that definition~\eqref{E:relative-consequence} in conjunction with Proposition~\ref{P:relative-consequence} shows a simple way to define an abstract logic without theorems. Indeed, let $\aLog$ be any abstract logic.
Then the abstract logic corresponding to the consequence relation $\vdashS^{\circ}$ has no theorems.\\

A set $X$, as well as the theory $\textbf{Cn}_{\mathcal{S}}(X)$, is called \textit{\textbf{inconsistent}}\index{theory!inconsistent} (relative to $\mathcal{S}$) if
$\textbf{Cn}_{\mathcal{S}}(X)=\FormsL$; otherwise both are \textit{\textbf{consistent}}\index{theory!consistent}.
An abstract logic $\mathcal{S}$ is \textit{\textbf{trivial}}\index{theory!trivial} if $\vdash_{\mathcal{S}}$ is trivial, and otherwise is \textit{\textbf{nontrivial}}\index{theory!nontrivial}. We observe that any structural $\mathcal{S}$ is trivial if, and only if, for any arbitrary variable $p$,
$p\in\ConS{\varnothing}$. (Exercise~\ref{section:realizations-abstract-logic}.\ref{EX:trivial-1})
We note that for any nontrivial abstract logic $\mathcal{S}$, $\bm{T}_{\mathcal{S}}\subset\FormsL$. (Exercise~\ref{section:realizations-abstract-logic}.\ref{EX:trivial-2})
) We also note that a theory is not required to be closed under substitution. A bit more we can say about substitution in the structural abstract logics: every theory of a structural abstract logic is closed under substitution if, and only if, this theory is generated by a set itself closed under substitution. (Exercise~\ref{section:realizations-abstract-logic}.\ref{EX:theory})

We denote
\[
\theory:=\set{X}{\ConS{X}=X}\tag{\textit{the set of $\mathcal{S}$-theories}}
\]
and
\[
\theoryC:=\set{X}{\ConS{X}=X\subset\FormsL}\tag{\textit{the set of consistent $\mathcal{S}$-theories}}
\]

Let $\aLog$ be an abstract logic and $X_0\subseteq\FormsL$. For future reference, we denote:
\[
\textbf{Cn}_{\mathcal{S}[X_0]}(X):=\ConS{X_0\cup X},
\]
where $X\subseteq\FormsL$.

It must be clear that $\textbf{Cn}_{\mathcal{S}[X_0]}$ is a consequence operator. The corresponding abstract logic we denote by $\mathcal{S}[X_0]$.
(See Exercise~\ref{section:consequence-operator}.\ref{EX:extension}.)

\paragraph{Exercises~\ref{section:consequence-operator}}
\begin{enumerate}
	\item \label{EX:Tarski-criterion}Prove that if $X$ is a theory with respect to to \textbf{Cn} and there is a countable-infinite family $\lbrace X_i\rbrace_{i<\omega}$ of theories (all with respect to \textbf{Cn}) such that
	\[
	X_0\subseteq X_1\subseteq\ldots\subset X
	\]
	and $\bigcup_{i<\omega}=X_i=X$, then $X$ is not finitely axiomatizable.
	\item \label{EX:intersection-con}Prove Proposition~\ref{P:intersection-con}.
	\item \label{EX:con-operation=con-relation}Prove Proposition~\ref{P:con-operation=con-relation}.
	\item Prove that a consequence operator \textbf{Cn} is structural if, and only if, $\Con{\sigma(\Con{X})}=\Con{\sigma(X)}$, for any set $X$ and any substitution $\sigma$.
	\item  Let \textbf{Cn} be a structural operator. Prove that a set $\sigma(\Con{X})$ is closed (or is a theory) if, and only if, 
	$\sigma(\Con{X})=\Con{\sigma(X)}$, for any set $X$ and any substitution $\sigma$. (Cf.~\cite{zygmunt2012}, section 2.2.)
	\item \label{EX:con-operator-intersection}Prove Proposition~\ref{P:con-operator-intersection}.
	\item \label{EX:con-connections}Complete the proof of Proposition~\ref{P:con-connections}.
	\item \label{EX:extension}Given an abstract logic $\aLog$ and $X_0\subseteq\FormsL$, show that
	$\aLog[X_0]$ is an abstract logic.
\end{enumerate}

\section{Realizations of abstract logic}\label{section:realizations-abstract-logic}
Abstract logic, understood either as a consequence relation or as a consequence operator, as we will see below, can be implemented in different forms. The implementation of abstract logic by the closure system (see Section~\ref{section:Con-via-closure-system}) is quite technical and does not reveal the meaning of the term `consequence'. Secondly, we note that perhaps historically arose first in the discourse of logic, namely, a logical consequence that can be determined by a deductive system, that is, through the concepts of a rule of inference and formal deduction. It is no coincidence that A. Tarski~\cite{tarski1936b} called this method of defining the logical consequence ``the old concept of consequence'' in the following passage:
\begin{quote}
	``It is perhaps not superfluous to point out in advance that --- in comparison with the new --- the old concept of consequence as commonly used by mathematical logicians in no way loses its importance. This concept will probably always have a decisive significance for the practical construction of deductive theories, as an instrument which allows us to prove or disprove particular sentences of these
	theories.'' (quoted from~\cite{tarski1956}, p. 413) 
\end{quote}

We present this way of defining abstract logic in Section~\ref{section:inference-rules}, and then use it throughout. What Tarski called above the ``new'' concept of consequence is connected with semantics. This will be our main way of representing abstract logic. Therefore, we dwell on this in more detail. A general exposition of semantic approach begins in Section~\ref{section:defining-con-semantically-general}, a detailed implementation of it continues in Section~\ref{section:con-via-matrices}. This brings us to the conception of separating tools, which we will discuss in Section~\ref{section:separating-means}.

The presence of two forms of representation of abstract logic as a consequence relation and as consequence operator, gives us a certain flexibility. Thus, we can choose one form or another that is more convenient in the given circumstances. We start with a fairly general way of defining a new abstract logic from a given one using the concept of substitution. 

\subsection{Defining new abstract logics from given ones by substitution}

Let $\aLog$ be an abstract logic in a language $\Lan$. We define a binary relation $\vdash_{\aLogSub}$ on $\mathcal{P}(\FormsL)\times\FormsL$ as follows:
\begin{equation}\label{E:def-of-aLogSub}
	X\vdash_{\aLogSub}\alpha~\define~\Sb X\vdash_{\aLog}\alpha.
\end{equation}

First, we notice the following.
\begin{prop}\label{P:S_sub-closed under-substitution}
	Let $\aLog$ be a structural abstract logic.
	Then for any set {\em$X\cup\lbrace\alpha\rbrace\subseteq\FormsL$} and any $\Lan$-substitution $\sigma$,
	{\em\[
		X\vdash_{\aLogSub}\alpha~\Longrightarrow~X\vdash_{\aLogSub}\sigma(\alpha).
		\]}
\end{prop}
\begin{proof}
	Assume that $X\vdash_{\aLogSub}\alpha$, that is $\Sb X\vdash_{\aLog}\alpha$.
	Since $\aLog$ is structural, we  have that 
	$\sigma(\Sb X)\vdash_{\aLog}\sigma(\alpha)$. And, in virtue of Proposition~\ref{P:Sb-properties}--e, we derive that $\Sb X\vdash_{\aLog}\sigma(\alpha)$, that is $X\vdash_{\aLogSub}\sigma(\alpha)$.
\end{proof}

\begin{prop}
	If $\aLog$ is a structural abstract logic, then {\em$\vdash_{\aLogSub}$} is
	a consequence relation. 
\end{prop}
\begin{proof}
	To check the properties (a)--(c) of Definition~\ref{D:consequnce-relation-single} for $\vdash_{\aLogSub}$, we will use Proposition~\ref{P:Sb-properties}.
	
	To prove the condition (a), we notice that if $\alpha\in X$, then $\alpha\in\Sb X$ (Proposition~\ref{P:Sb-properties}--a) and hence, be definition~\eqref{E:def-of-aLogSub}, $X\vdash_{\aLogSub}\alpha$.
	
	Next, suppose that $X\vdash_{\aLogSub}\alpha$ and $X\subseteq Y$. In virtue of
	Proposition~\ref{P:Sb-properties}--d, $\Sb X\subseteq\Sb Y$. Using~\eqref{E:def-of-aLogSub} twice, we conclude that  $Y\vdash_{\aLogSub}\alpha$. Using the first premise and Proposition~\ref{P:S_sub-closed under-substitution}
	
	Now, assume that $X\vdash_{\aLogSub}\beta$, for all $\beta\in Y$, and $Y\cup Z\vdash_{\aLogSub}\alpha$. Applying the property Proposition~\ref{P:Sb-properties}--b to the second premise, we get:
	$\Sb Y\cup\Sb Z\vdash_{\aLogSub}\alpha$. Using the first premise and Proposition~\ref{P:S_sub-closed under-substitution}, we obtain that
	$\Sb X\vdash_{\aLog}\gamma$, for all $\gamma\in\Sb Y$. Therefore,
	$\Sb X\cup\Sb Z\vdash_{\aLog}\alpha$, that is $X\cup Z\vdash_{\aLogSub}\alpha$.
\end{proof}

We denote by $\aLogSub$ the abstract logic corresponding to the relation $\vdash_{\aLogSub}$, providing that $\aLog$ is a structural abstract logic.

\begin{prop}\label{P:aLog-implies-aLogSub}
	Let $\aLog$ be a structural abstract logic. Then
	{\em\[
		X\vdash_{\aLog}\alpha~\Longrightarrow~X\vdash_{\aLogSub}\alpha.
		\]}
\end{prop}
\begin{proof}
	Assume that $X\vdash_{\aLog}\alpha$. Then, obviously, $X\cup\Sb X\vdash_{\aLog}\alpha$. According to the property Proposition~\ref{P:Sb-properties}--a, $\Sb X\vdash_{\aLog}\alpha$, that is
	$X\vdash_{\aLogSub}\alpha$.
\end{proof}

The converse of the implication in Proposition~\ref{P:aLog-implies-aLogSub} in general does not hold. See Exercise~\ref{section:realizations-abstract-logic}.\ref{EX:demonstration} below.

In order to define $\mathcal{S}_{\text{sub}}$, we need a structural logic $\mathcal{S}$. However,  $\mathcal{S}_{\text{sub}}$ itself is almost never structural, although for any substitution $\sigma$, $\sigma(\textbf{Cn}_{\mathcal{S}_{\text{sub}}}(\varnothing))\subseteq\textbf{Cn}_{\mathcal{S}_{\text{sub}}}(\varnothing)$. More precisely, $\mathcal{S}_{\text{sub}}$ is structural in a weak sense if, and only if, $X\vdash_{\mathcal{S}_{\text{sub}}}\alpha$, for any nonempty set $X$ and any formula $\alpha$.

Indeed, in view of Proposition~\ref{P:S_sub-closed under-substitution}, we have: $p\vdash_{\aLogSub}q$. Let $X$ be a nonempty set with $\beta\in X$ and $\alpha$ be an arbitrary formula. Then, if $\mathcal{S}_{\text{sub}}$ is structural in a weak sense, we have:  $\beta\vdash_{\aLogSub}\alpha$ and hence $X\vdash_{\aLogSub}\alpha$.

The converse is obvious.
\begin{prop}
	Let $\mathcal{S}$ be a structural abstract logic. If $\mathcal{S}$ is finitary, then {\em$\mathcal{S}_{\text{sub}}$} is also finitary.
\end{prop}
\begin{proof}
	Suppose $X\vdash_{\aLogSub}\alpha$, that is $\Sb X\vdash_{\mathcal{S}}\alpha$. Since $\mathcal{S}$ is finitary, there is a set $Y_0\Subset\Sb X$ such that $Y_0\vdash_{\mathcal{S}}\alpha$. The latter implies that $\Sb Y_0\vdash_{\mathcal{S}}\alpha$. Further, in virtue of Proposition~\ref{P:Sb-properties}--e, $\Sb Y_0\subseteq\Sb X$.
	
	Now, we define: $X_0:=Y_0\cap X$. Using the property (c) of Proposition~\ref{P:Sb-properties} and the above inclusion, we conclude that $\Sb X_0=\Sb Y_0$. This implies that $X_0\vdash_{\mathcal{S}}\alpha$ and, hence, $X_0\vdash_{\aLogSub}\alpha$.
\end{proof}

\subsection{Consequence operator via a closure system}\label{section:Con-via-closure-system}
A set $\pmb{\mathcal{A}}\subseteq\mathcal{P}(\FormsL)$ which contains $\FormsL$ and is closed under any nonempty intersections in $\mathcal{P}(\FormsL)$ is called a \textit{\textbf{closure system over}}\index{closure system} $\FormsL$. That is, for any closure system $\mathcal{A}$, $\FormsL\in\pmb{\mathcal{A}}$ and for any $I\neq\varnothing$,
\[
\lbrace X_i\rbrace_{i\in I}\subseteq\pmb{\mathcal{A}}\Longrightarrow\bigcap_{i\in I}
\lbrace X_i\rbrace\in\pmb{\mathcal{A}}.
\]

All operators in this subsection are considered over a fixed language $\Lan$.
\begin{prop}\label{P:closure-system}
	Let $\pmb{\mathcal{A}}$ be a closure system. Then the operator
	{\em\[
		\textbf{Cn}_{\pmb{\mathcal{A}}}(X):=\bigcap\set{Y}{X\subseteq Y~\textit{and}\,~Y\in\pmb{\mathcal{A}}}
		\]}is a consequence operator and {\em$\pmb{\mathcal{A}}=\set{X}{X=\textbf{Cn}_{\pmb{\mathcal{A}}}(X)}$}. Conversely, if {\em\textbf{Cn}} is a consequence operator and {\em$\pmb{\mathcal{A}}_{\textbf{Cn}}$} is the family of its closed sets, then
	{\em$\textbf{Cn}=\textbf{Cn}_{\pmb{\mathcal{A}}_{\textbf{Cn}}}$}.
\end{prop}\index{$\textbf{Cn}_{\pmb{\mathcal{A}}}(X)$}
\begin{proof}
	Assume that $\pmb{\mathcal{A}}$ is a closure system. It is easy to check that the operator $\textbf{Cn}_{\mathcal{A}}$ satisfies the properties $(\text{a}^{\dag})$ and $(\text{b}^{\dag})$ of Definition~\ref{D:consequence-operator}. We prove that it also satisfies $(\text{c}^{\dag})$. For this, given a set $X$, we denote
	\[
	X_0:=\textbf{Cn}_{\pmb{\mathcal{A}}}(X)
	\]
	and suppose that $\alpha\in \textbf{Cn}_{\pmb{\mathcal{A}}}(X_0)$. Aiming to show that $\alpha\in X_0$, we note that
	\[
	\alpha\in X_0\Longleftrightarrow\forall Y\in\pmb{\mathcal{A}}.~X\subseteq Y\Rightarrow
	\alpha\in Y.
	\]
	Thus, we assume that $X\subseteq Y$, for an arbitrary $Y\in\pmb{\mathcal{A}}$. In virtue of 
	$(\text{b}^\dag)$, $X_0\subseteq\textbf{Cn}_{\pmb{\mathcal{A}}}(Y)$ and hence $\Con{X_0}\subseteq Y$, for $Y\in\pmb{\mathcal{A}}$. This implies that $\alpha\in Y$ and, using the equivalence above, we obtain that $\alpha\in X_0$.
	
	Finally, one can easily observe that
	\[
	X\in\pmb{\mathcal{A}}\Longleftrightarrow\textbf{Cn}_{\mathcal{A}}(X)=X.
	\]
	
	Conversely, suppose \textbf{Cn} is a consequence operator and
	\[
	\pmb{\mathcal{A}}_{\textbf{Cn}}:=\set{X}{\textbf{Cn}(X)=X}.
	\]\index{$\pmb{\mathcal{A}}_{\textbf{Cn}}$}
	Assume that $\lbrace X_i\rbrace_{i\in I}\subseteq\pmb{\mathcal{A}}_{\textbf{Cn}}$ and denote
	\begin{equation}\label{E:conjunction}
		X_0:=\bigcap_{i\in I}\lbrace X_i\rbrace.
	\end{equation}
	Since $X_0\subseteq X_i$, for each $i\in I$, $\Con{X_0}\subseteq\Con{X_i}=X_i$, also for each $i\in I$. Hence $\Con{X_0}\subseteq X_0$, that is $X_0\in\pmb{\mathcal{A}}_{\textbf{Cn}}$.
	
	Finally, we show that for any formula $\alpha$ and any set $X$ of formulas,
	\[
	\alpha\in\textbf{Cn}(X)\Longleftrightarrow\alpha\in
	\textbf{Cn}_{\pmb{\mathcal{A}}_{\textbf{Cn}}}(X).
	\]
	
	First, we notice that
	\[
	\alpha\in
	\textbf{Cn}_{\pmb{\mathcal{A}}_{\textbf{Cn}}}(X)\Longleftrightarrow
	\forall Y.~(X\subseteq Y~\text{and}~\Con{Y}=Y)\Rightarrow\alpha\in Y. \tag{$\ast$}
	\]
	Then, in virtue of $(\ast)$, we have:
	\[
	\begin{array}{rl}
		\alpha\in
		\textbf{Cn}_{\pmb{\mathcal{A}}_{\textbf{Cn}}}(X)\!\!\!\! &\Longrightarrow
		[(X\subseteq\Con{X}~\text{and}~\Con{\Con{X}}=\Con{X})\Rightarrow\alpha\in\Con{X}]\\
		&\Longrightarrow\alpha\in\Con{X}.
	\end{array}
	\]
	
	Now we suppose that $\alpha\in\Con{X}$, $X\subseteq Y$ and $\Con{Y}=Y$. From the second premise, we derive that $\Con{X}\subseteq\Con{Y}$ and hence $\Con{X}\subseteq
	Y$. Thus $\alpha\in Y$. Since $Y$ is an arbitrary set satisfying the aforementioned conditions, using $(\ast)$, we conclude that $\alpha\in\textbf{Cn}_{\pmb{\mathcal{A}}_{\textbf{Cn}}}(X)$.
\end{proof}

Given an abstract logic $\mathcal{S}$, we denote:
\[
\pmb{\mathcal{A}}_{\mathcal{S}}:=\set{X}{X=\ConS{X}}.
\]

We conclude with the following observation.
\begin{prop}[Brown-Suszko theorem]\label{P:brown-suszko-theorem}
	An abstract logic $\mathcal{S}$ in a language $\Lan$ is structural if, and only if, for any set {\em$X\subseteq\FormsL$} and any substitution $\sigma$, the following holds:
	\[
	X\in\pmb{\mathcal{A}}_{\mathcal{S}}\Longrightarrow\sigma^{-1}(X)\in\pmb{\mathcal{A}}_{\mathcal{S}}.
	\]
\end{prop}
\begin{proof}
	Suppose $\mathcal{S}$ is structural. Let $\sigma$ be a substitution and $X\in\pmb{\mathcal{A}}_{\mathcal{S}}$. Assume that $\sigma^{-1}(X)\vdash_{\mathcal{S}}\alpha$.
	Since $\mathcal{S}$ is structural, we have: $\sigma(\sigma^{-1}(X))\vdash_{\mathcal{S}}\sigma(\alpha)$. In virtue of (\ref{E:substitution-inequalities}--$i$), $X\vdash_{\mathcal{S}}\sigma(\alpha)$, that is $\sigma(\alpha)\in X$. This means that $\alpha\in\sigma^{-1}(X)$.
	
	To prove the converse, we assume that the conditional above holds. Then we, with help of~(\ref{E:substitution-inequalities}--$ii$), obtain that $X\subseteq\sigma^{-1}(\sigma(X))\subseteq\sigma^{-1}(\ConS{\sigma(X)})$. This implies that $\ConS{X}\subseteq\ConS{\sigma^{-1}(\ConS{\sigma(X)})}$. Applying the assumption, we get that $\ConS{X}\subseteq\sigma^{-1}(\ConS{\sigma(X)})$. Thus, if $X\vdash_{\mathcal{S}}\alpha$, that is $\alpha\in\ConS{X}$, then $\alpha\in\sigma^{-1}(\ConS{\sigma(X)})$, that is $\sigma(X)\vdash_{\mathcal{S}}\sigma(\alpha)$.
\end{proof}

\subsection{Defining abstract logic semantically: a general approach}\label{section:defining-con-semantically-general}
Consider an abstract object $\fM$ that does not belong to the language $\Lan$, which can \emph{verify} $\Lan$-formulas. The term `verify' will remain undefined during consideration until the end of this subsection, when we propose two explications for this term and for the object $\fM$. We call $\fM$ a \textit{model}. If the model $\fM$, verifying an $\Lan$-formula $\alpha$, \emph{accepts} it, we write $\fM\models\alpha$, if it does not, we write $\fM\not\models\alpha$ and say that $\fM$ \textit{rejects} $\alpha$. Thus a metalogic here is two-valued, which is useful to remember. The relation of verification of an $\Lan$-formula in a model $\fM$ is extended to any set $X$ of $\Lan$-formulas in a usual way:
\[
\fM\models X~\define~(\fM\models\alpha,~\text{for any $\alpha\in X$}).
\]\index{$\models$}

The relation of \textit{semantic consequence} is defined as follows: Given a model $\fM$, for any set $X\cup\lbrace\alpha\rbrace\subseteq\Forms_{\Lan}$,
\begin{equation}\label{E:semantic-consequence}
	X\models_{\fM}\alpha~\define~(\fM\models X~\Longrightarrow~\fM\models\alpha).
\end{equation}

As an auxiliary notation, we also define:
\[
X\models_{\fM}Y~\define~(\fM\models X~\Longrightarrow~\fM\models Y),
\]
where $X\cup Y\subseteq\Forms_{\Lan}$.\index{$\models_{\fM}$}

\begin{prop}\label{P:semantic-consequence}
	The semantic consequence {\em\eqref{E:semantic-consequence}} is a consequence relation.
\end{prop}
\begin{proof}
	Let $\fM$ be a model that can verify or reject every $\Lan$-formula. We have to show that the properties (a)--(c) of Definition~\ref{D:consequnce-relation-single} are valid for the relation $\models_{\fM}$. In terms of this relation, these properties can be written as follows.
	\[
	\begin{array}{cl}
		(\text{a}^\ast) &X\models_{\fM}\alpha,~\text{whenever $\alpha\in X$};\\
		(\text{b}^\ast) &(X\models_{\fM}\alpha~\text{and}~X\subseteq Y)~\Longrightarrow~Y\models_{\fM}\alpha;\\
		(\text{c}^\ast) &(X\models_{\fM}Y~\text{and}~Y\cup Z\models_{\fM}\alpha)~\Longrightarrow X\cup Z\models_{\fM}\alpha.
	\end{array}
	\]
	
	We prove $(\text{c}^\ast)$ and leave the properties $(\text{a}^\ast)$ and $(\text{b}^\ast)$ to the reader. (Exercise~\ref{section:realizations-abstract-logic}.\ref{EX:semantic-consequence})
	
	Assume that $X\models_{\fM}Y,~Y\cup Z\models_{\fM}\alpha$ and $\fM\models X\cup Z$. Then the first and third premises imply that $\fM\models Y$ and hence $\fM\models Y\cup Z$. This, in conjunction with the second premise, implies that $\fM\models\alpha$.
\end{proof}

Let $\bm{T}_{\fM}$ be the set of all $\Lan$-formulas that are verified by a model $\fM$, that is,\index{$\bm{T}_{\fM}$}
\[
\bm{T}_{\fM}:=\set{\alpha\in\Forms_{\Lan}}{\fM\models\alpha}.
\]

It is obvious that, given a model $\fM$, the consequence relation $\models_{\fM}$ can be formulated in terms of $\bm{T}_{\fM}$ as follows:
\begin{equation}\label{E:semantic-consequence-equivalent}
	X\models_{\fM}\alpha~\Longleftrightarrow~X\not\subseteq\bm{T}_{\fM}~\text{or}~
	\alpha\in\bm{T}_{\fM}.
\end{equation}

The last observation gives rise to the following definitions:
\[
\bm{T}_{\fM}^{\ast}:=\set{(X,\alpha)}{X\cup\lbrace\alpha\rbrace\subseteq\Forms_{\Lan}~\text{and}~\alpha\in\bm{T}_{\fM}}
\]
and
\[
\bm{T}_{\fM}^{\ast\ast}:=\set{(X,\alpha)}{X\cup\lbrace\alpha\rbrace\subseteq\Forms_{\Lan}~\text{and}~X\not\subseteq\bm{T}_{\fM}}.
\]

Using~\eqref{E:semantic-consequence-equivalent}, one can easily conclude that
\[
\bm{T}_{\fM}^{\ast}\cup\bm{T}_{\fM}^{\ast\ast}\subseteq\,\models_{\fM}.
\]
On the other hand, if $\vdash$ is a consequence relation with $\bm{T}_{\fM}^{\ast}\cup\bm{T}_{\fM}^{\ast\ast}\subseteq\,\vdash$, then
$\models_{\fM}\subseteq\,\vdash$. (See Exercise~\ref{section:realizations-abstract-logic}.\ref{EX:semantic-consequence-2}.)

Thus, we have proved the following.
\begin{prop}\label{P:semantic-consequence-2}
	Given a model $\fM$, the consequence relation $\models_{\fM}$ is the least consequence relation $\vdash$ with $\bm{T}_{\fM}^{\ast}\cup\bm{T}_{\fM}^{\ast\ast}\subseteq\,\vdash$.
\end{prop}

The discussion above can be easily extended to the case when instead of one model we have a class of models, each of which can verify or reject each $\Lan$-formula.

Let $\cM$ be such a class of models. Then we define:
\begin{equation}\label{E:semantic-consequence-generalized}
	X\models_{\cM}\alpha~\define~(X\models_{\fM}\alpha,~\text{for any $\fM\in\cM$}).
\end{equation}
\begin{prop}\label{P:semantic-consequence-generalized}
	Let $\cM$ be a class of models that can verify or reject $\Lan$-formulas. Then the relation $\models_{\cM}$ defined in {\em\eqref{E:semantic-consequence-generalized}} is a consequence relation.
\end{prop}
\begin{proof}
	In virtue of Proposition~\ref{P:semantic-consequence} and Proposition~\ref{P:con-relation-intersection}, the relation $\models_{\cM}$ is a consequence relation.
\end{proof}

We note that, in virtue of Proposition~\ref{P:semantic-consequence-2}, we have:
\begin{equation}\label{E:semantic-consequence-2-equivalent}
	\models_{\cM}=\bigcap_{\fM\in\cM}\set{\vdash}{\vdash~\text{is a consequence relation with}~\bm{T}_{\fM}^{\ast}\cup\bm{T}_{\fM}^{\ast\ast}\subseteq\,\vdash}.
\end{equation}
(Exercise~\ref{section:realizations-abstract-logic}.\ref{EX:semantic-consequence-2-equivalent})\\

Now, let us fix a set $X\cup\lbrace\alpha\rbrace\subseteq\Forms_{\Lan}$ and a class $\cM$ of models. We denote:\index{$\textbf{Mod}_{\cM}$}
\[
\textbf{Mod}_{\cM}(X,\alpha):=\set{\fM\in\cM}{X\models_{\fM}\alpha},
\]
\[
\textbf{Mod}_{\cM}(X):=\set{\fM\in\cM}{\fM\models X},
\]
and as a particular case of the last definition
\[
\textbf{Mod}_{\cM}(\alpha):=\set{\fM\in\cM}{\fM\models\alpha}.
\]
\begin{prop}\label{P:semantic-consiquence-equality-models}
	Given a set {\em $X\cup\lbrace\alpha\rbrace\subseteq\Forms_{\Lan}$} and a class
	$\cM$ of models,
	{\em\[
		\textbf{Mod}_{\cM}(X,\alpha)=\textbf{Mod}_{\cM}(\alpha)\cup (\cM\setminus\textbf{Mod}_{\cM}(X)).
		\]}
\end{prop}
\begin{proof}
	Indeed, we have:
	\[
	\begin{array}{rl}
		\fM\in\textbf{Mod}_{\cM}(X,\alpha)\!\! &\Longleftrightarrow~\fM\models\alpha~\text{or}~X\not\subseteq\bm{T}_{\fM},
		~\text{and}~\fM\in\cM
		\quad[\text{in virtue of~\eqref{E:semantic-consequence-equivalent}}]\\
		&\Longleftrightarrow~\fM\in\textbf{Mod}_{\cM}(\alpha)\cup(\cM\setminus
		\textbf{Mod}_{\cM}(X)).
	\end{array}
	\]
\end{proof}

Now, if we do not limit ourselves to the class of $\cM$ models, but take into account all models that have the ability to accept or reject $\Lan$-formulas, we can offer similar definitions.
\[
\textbf{Mod}(X,\alpha):=\set{\fM}{X\models_{\fM}\alpha},
\]
\[
\textbf{Mod}(X):=\set{\fM\in}{\fM\models X},
\]
and 
\[
\textbf{Mod}(\alpha):=\set{\fM}{\fM\models\alpha}.
\]

Also, as in Proposition~\ref{P:semantic-consiquence-equality-models}, we have the following equality:
\begin{equation}\label{E:semantic-consiquence-equality-models}
	\textbf{Mod}(X,\alpha)=\textbf{Mod}(\alpha)\cup\set{\fM}{\fM\notin\textbf{Mod}(X)}.
\end{equation}

Let $R\subset\mathcal{P}(\Forms_{\Lan})\times\Forms_{\Lan}$. We define:
\[
\textbf{Mod}(R):=\bigcap_{(X,\alpha)\in R}\textbf{Mod}(X,\alpha).
\]

In particular, we have:
\[
\textbf{Mod}(\models_{\cM})=\bigcap\set{\textbf{Mod}(X,\alpha)}{X\models_{\cM}\alpha}.
\]

For brevity, we denote:
\[
\cM^{\ast}:=\textbf{Mod}(\models_{\cM}).
\]

The, we observe that 
\begin{equation}\label{E:M_subseteq_M^ast}
	\cM\subseteq\cM^{\ast}.
\end{equation}
and
\begin{equation}\label{E:vdash_M=vdash_M^ast}
	\models_{\cM}~=~\,\models_{{\cM}^{\ast}}
\end{equation}
(We leave for the reader to prove~\eqref{E:M_subseteq_M^ast} and~\eqref{E:vdash_M=vdash_M^ast}; see Exercise~\ref{section:realizations-abstract-logic}.\ref{EX:M_subseteq_M^ast}  and Exercise~\ref{section:realizations-abstract-logic}.\ref{EX:vdash_M=vdash_M^ast}, respectively.)

This takes us to the conception of separating tools which will be discussed in Section~\ref{section:separating-means}. Indeed, if we aim to show that $X\not\models_{\cM}\alpha$, it may be easier to find a model that \emph{separates} $X$ from $\alpha$ in $\cM^\ast$, not $\cM$.\\

In this book, the term \emph{model} will be refined in two different ways. These two refinements will lead to two different types of consequence relation given semantically.

Let $\mat{M}=\langle\alg{A},D\rangle$ be a logical matrix of type $\Lan$ and let $\Phi_\textbf{M}$ be the set of all valuations in the algebra $\alg{A}$.
We call a\emph{ model of type} (A) a pair $\langle\mat{M},v\rangle$, where $v\in\Phi_\textbf{M}$. The model $\langle\mat{M},v\rangle$ \emph{accepts} a formula $\alpha$ if $v[\alpha]\in D$; otherwise it \emph{rejects} $\alpha$. The models of type (A) are used in the definition of matrix consequence for the formulas of language $\Lan$ (see Section~\ref{section:con-via-matrices}), in the definition of $\E$-consequence for $\Lan$-equalities (see Section~\ref{section:E-consequence}), and in the definition of $\mathcal{Q}$-consequence for $\Lan_{\mathcal{Q}}$-formulas (see Section~\ref{section:Q-consequence}).  

A \emph{model of type} (B) is a pair $\langle\mat{M},\Phi_\textbf{M}\rangle$.
The model $\langle\mat{M},\Phi_\textbf{M}\rangle$ \emph{accepts} a formula $\alpha$ if, and only if, $v[\alpha]\in D$, for every $v\in\Phi_\textbf{M}$. This type is used in the definition of $\E_\text{L}$-consequence for $\Lan$-equalities (see Section~\ref{section:equational-L-consequence}).

\subsection{Consequence relation via logical matrices}\label{section:con-via-matrices}
The importance of this section is that the consequence operator of any abstract logic $\mathcal{S}$ is determined by $\mathcal{S}$-theories.

Given a (logical) matrix $\mat{M}=\langle\alg{A},D\rangle$, we define a relation
$\models_{\textbf{M}}$ on $\mathcal{P}(\FormsL)\times\FormsL$ as follows:
\[
X\models_{\textbf{M}}\alpha\stackrel{\text{df}}{\Longleftrightarrow}(v[X]\subseteq D
\Rightarrow v[\alpha]\in D,~\text{for any valuation $v$ in \alg{A}}),
\]
where
\[
v[X]:=\set{v[\beta]}{\beta\in X}.
\]

Let $\mathcal{M}$ be a nonempty class of matrices. We define the relation of \textit{\textbf{matrix consequence}}\index{consequence!matrix} (with respect to a class $\mathcal{M}$) as follows:
\begin{equation}\label{E:matrix-consequence}
	X\models_{\mathcal{M}}\alpha\stackrel{\text{df}}{\Longleftrightarrow}
	(X\models_{\textbf{M}}\alpha,~\text{for all $\mat{M}\in\mathcal{M}$}).
\end{equation}
Accordingly, the relation $\models_{\textbf{M}}$ is called a \textit{\textbf{single-matrix consequence}}\index{consequence!single matrix}, or $\mat{M}$-\textit{\textbf{consequence}} for short.
If a class $\mathcal{M}$ is known, we call a matrix consequence associated with a class $\mathcal{M}$ an $\mathcal{M}$-\textit{\textbf{consequence}}. Given a nonempty set of matrices $\mathcal{M}$, we will be denoting the $\mathcal{M}$-consequence by $\mathcal{S}_{\mathcal{M}}$ and will be saying that an abstract logic $\mathcal{S}$ is \textit{\textbf{determined}}\index{consequence!determined by} by $\mathcal{M}$ if $\mathcal{S}=\mathcal{S}_{\mathcal{M}}$; and we use the notation $\aLog_{\mat{M}}$ if $\mathcal{M}=\{\mat{M}\}$.

The justification for the terminology `$\mathcal{M}$-consequence' and `\mat{M}-consequence' is the following. 

\begin{prop}\label{P:matrix-con-is-con-relation}
	Any matrix consequence is a structural consequence relation.
\end{prop}
\begin{proof}
	Let $\mathcal{M}=\lbrace\textsf{\textbf{M}}_i\rbrace_{i\in I}$ be a nonempty set of matrices. We have to check that the conditions (a)--(c) of Definition~\ref{D:consequnce-relation-single} are satisfied for the relation $\models_{\mathcal{M}}$. In addition, we have to show that this consequence relation is structural. We leave this routine check to the reader. (Exercise~\ref{section:realizations-abstract-logic}.\ref{EX:matrix-con-is-con-relation})
\end{proof}

Single-matrix consequence gives us examples of both finitary and non-finitary consequence relations. The first type is exemplified by $\models_{\textbf{B}_{2}}$, a single-matrix consequence with respect to the matrix $\textbf{B}_{2}$ (Section~\ref{S:two-valued}). This will be explained in Section~\ref{section:inference-rules}. Now, we show that the single-matrix consequence $\models_{\mat{LC}^{\ast}}$ (with respect to the matrix $\mat{LC}^{\ast}$ (Section~\ref{section:dummett})) is not finitary.

Indeed, let us denote:
\[
X:=\set{(p_{i}\leftrightarrow p_{j})\rightarrow p_0}{0<i<j<\omega}.
\]

Let $v$ be a valuation in the $\mat{LC}^{\ast}$-algebra. Assume that
$v[X]\subseteq\lbrace\one\rbrace$. If for some $i$ and $j$ with $i<j$, $v[p_i]=v[p_j]$, then $v[p_0]=\one$. If for any $i$ and $j$ with $i<j$,
$v[p_i]\neq v[p_j]$, then each $v[p_i\leftrightarrow p_j]=v[p_i]\land v[p_j]$. Hence, there is an ascending chain, all element of which are less than or equal to $v[p_0]$. This implies that $v[p_0]=\one$.  From this we conclude that $X\models_{\mat{LC}^{\ast}}p_0$.

On the other hand, let $X_0$ be any finite subset of $X$. If $X_0=\varnothing$, then any valuation $v^{\prime}$ with $v^{\prime}[p_0]\neq\one$ ``separates'' $X$ from $p_0$.  Suppose that $X_0\neq\varnothing$. Let us take a valuation $v^{\ast}$ with the property: $v^{\ast}[p_i]=a_i$, where
\[
\zero=a_1<a_2<\ldots<\one
\]
are the elements of the matrix $\mat{LC}^{\ast}$. It must be clear that for any $i$, $v[p_i\leftrightarrow p_j]=a_i$. Thus $v^{\ast}[X_0]=\lbrace a_{i_1},\ldots,a_{i_k}\rbrace$. If we define $v^{\ast}[p_0]$ to satisfy the conditions $\max\lbrace a_{i_1},\ldots,a_{i_k}\rbrace<v^{\ast}[p_0]<\one$,
then $v^{\ast}[X_0]\subseteq\lbrace\one\rbrace$, but $v^{\ast}[p_0]\neq\one$.\\

Let $\mathcal{S}$ be an abstract logic and $\mat{M}=\langle\alg{A},D\rangle$ be a matrix. Accordingly, we have the consequence relation $\vdash_{\mathcal{S}}$ and a single-matrix consequence $\models_{\textbf{M}}$. If $\vdash_{\mathcal{S}}\subseteq\models_{\textbf{M}}$, we call the filter $D$ an $\mathcal{S}$-\textit{\textbf{filter}}\index{matrix!filter}, in which case the matrix \mat{M} is called an $\mathcal{S}$-\textit{\textbf{matrix}} (or an $\mathcal{S}$-\textit{\textbf{model}}).

\begin{defn}[Lindenbaum matrix]\label{D:lindebaum-matrix}\index{matrix!Lindenbaum}
	Let $\mathcal{S}$ be a structural abstract logic and $D_{\mathcal{S}}$ be an $\mathcal{S}$-theory. The matrix $\langle\FormAl,D_{\mathcal{S}}\rangle$ is called a Lindenbaum matrix $($relative to $\mathcal{S}$$)$. For an arbitrary {\em$X\subseteq\FormsL$}, we denote:\index{$\Lin_{\,\mathcal{S}}$}
	{\em\[
		\Lin_{\,\mathcal{S}}[X]:=\langle\FormAl,\ConS{X}\rangle
		\]}
	and
	{\em\[
		\Lin_{\,\mathcal{S}}:=\langle\FormAl,\ConS{\varnothing}\rangle.
		\]}
\end{defn}

We remind the reader that 
\begin{center}
	\textit{any $\Lan$-substitution is a valuation in $\mathfrak{F}_{\mathcal{L}}$ and vice versa.}	
\end{center}
(Cf.~Proposition~\ref{P:substitution-as-homomorphism}.)\\

It is obvious that, given a set $X\subseteq\FormsL$, for any formula $\alpha$,
\[
X\vdashS\alpha
~\Longleftrightarrow~\LinS[X]\models_{\iota}\alpha,
\]
where $\iota$ is the identity substitution.

This shows the prospect that we need to explore if we want to find a method for determining whether $X\vdashS\alpha$ or $X\not\vdashS\alpha$.

We begin with the following.
\begin{prop}[Lindenbaum's theorem]\label{P:lindenbaum-theorem}
	Let $\mathcal{S}$ be a structural abstract logic. Then {\em$\Lin_{\mathcal{S}}$} is an $\mathcal{S}$-model. Also, if a set {\em$X\subseteq\FormsL$} is closed under arbitrary $\Lan$-substitution, then for any $\Lan$-formula $\alpha$,
	{\em\[
		X\vdash_{\mathcal{S}}\alpha~\Longleftrightarrow~\LinS[X]\models\alpha;
		\]}
	in particular,
	{\em\[
		\alpha\in\bm{T}_{\mathcal{S}}~\Longleftrightarrow~\LinS\models\alpha.
		\]}
\end{prop}
\begin{proof}
	Let us take any set $X\cup\lbrace\alpha\rbrace\subseteq\FormsL$ such that $X\vdash_{\mathcal{S}}\alpha$. Then, for an arbitrary $\Lan$-substitution $\sigma$, $\sigma(X)\vdash_{\mathcal{S}}\sigma(\alpha)$. Therefore, if $\sigma(X)\subseteq\ConS{\varnothing}$, then $\sigma(\alpha)\in\ConS{\varnothing}$. Thus
	$\Lin_{\,\mathcal{S}}$ is an $\mathcal{S}$-model. 
	
	Now, if, in addition, $\sigma(X)\subseteq X$, then $\sigma(\alpha)\in\ConS{X}$. Since $\sigma$ is an arbitrary $\Lan$-substitution, we have that
	${\textsf{\textbf{Lin}}_{\,\mathcal{S}}[X]}\models\alpha$.
	
	Conversely, if ${\textsf{\textbf{Lin}}_{\,\mathcal{S}}[X]}\models\alpha$, then, in particular, $\iota(\alpha)\in\ConS{X}$, that is $X\vdash_{\mathcal{S}}\alpha$.
	
	The second equivalence is a specification of the first.
\end{proof}

We say that a matrix is \textit{\textbf{weakly adequate}}\index{matrix!weakly adequate} for an abstract logic $\mathcal{S}$ if this matrix validates all $\mathcal{S}$-theorems and only them. Thus, by virtue of the last equivalence of Proposition~\ref{P:lindenbaum-theorem},\index{Theorem!Lindebaum}
we obtain the following.
\begin{cor}\label{C:lindenbaum-matrix}
	Let $\mathcal{S}$ be structural abstract logic. Then any Lindenbaum matrix relative to $\mathcal{S}$
	is an $\mathcal{S}$-model.  Moreover, {\em$\Lin_{\,\mathcal{S}}$} is weakly adequate for $\mathcal{S}$.
\end{cor}
\noindent{\em Proof}~is left to the reader. (Exercise~\ref{section:realizations-abstract-logic}.\ref{EX:lindenbaum-matrix})\\

Proposition~\ref{P:lindenbaum-theorem} can be generalized as follows.
\begin{prop}\label{P:Lindenbaum-theorem-generalized}
	Given a structural abstract logic $\mathcal{S}$, let $\mathcal{M}_{\mathcal{S}}$ be the set of all Lindenbaum matrices relative to $\mathcal{S}$. Then $\vdash_{\mathcal{S}}\,=\,\models_{\mathcal{M}_{\mathcal{S}}}$; that is for any set $X$ of formulas and any formula $\alpha$,
	\[
	X\vdash_{\mathcal{S}}\alpha\Longleftrightarrow X\models_{\mathcal{M}_{\mathcal{S}}}
	\alpha.
	\]
\end{prop}
\begin{proof}
	First we suppose that $X\vdash_{\mathcal{S}}\alpha$. Let us take any $\langle\FormAl,T\rangle$, where $T\in\Sigma_{\mathcal{S}}$. Assume that for some substitution $\sigma$, $\sigma(X)\subseteq T$.
	Since $\mathcal{S}$ is structural, $\sigma(X)\vdash_{\mathcal{S}}\sigma(\alpha)$.
	This implies that $\sigma(\alpha)\in\textbf{Cn}_{\mathcal{S}}(\sigma(X))$ and hence
	$\sigma(\alpha)\in T$. 
	
	Conversely, assume that $X\not\vdash_{\mathcal{S}}\alpha$. This means (see~\eqref{E:Cn-S-connection}) that $\alpha\notin\textbf{Cn}_{\mathcal{S}}(X)$.
	Let us denote $\textbf{M}:=\langle\FormAl,\textbf{Cn}_{\mathcal{S}}(X)\rangle$ and take the identity substitution $\iota$. It is obvious that $\iota(X)\subseteq
	\textbf{Cn}_{\mathcal{S}}(X)$ but $\iota(\alpha)\notin\textbf{Cn}_{\mathcal{S}}(X)$; that is $X\not\models_{\textbf{M}}\alpha$.
\end{proof}

The last proposition inspires us for the following (complex) definition.
\begin{defn}\label{D:atlas-bundle}\index{matrix!bundle}\index{matrix!atlas}
	A $($nonempty$)$ family {\em$\mathcal{B}=\lbrace\langle\alg{A}_{i},D_i\rangle\rbrace_{i\in I}$} of matrices of type $\Lan$ is called a \textbf{bundle} if for any $i,j\in I$, the algebras {\em$\alg{A}_i$} and {\em$\alg{A}_j$} are isomorphic. A pair {\em$\langle\alg{A},\lbrace D_i\rbrace_{i\in I}\rangle$}, where each $D_i$ is a logical filter in $($or of$)$ {\em$\alg{A}$}, is called an \textbf{atlas}. By definition, the \textbf{matrix consequence relative to an atlas}
	{\em$\langle\alg{A},\lbrace D_i\rbrace_{i\in I}\rangle$} is the same as the matrix consequence related the bundle {\em$\lbrace\langle\alg{A},D_i\rangle\rbrace_{i\in I}$}, which in turn is defined according to~{\em \eqref{E:matrix-consequence}}.
	Given an abstract logic $\mathcal{S}$, the bundle  $\lbrace\langle\FormAl,T\rangle\rbrace_{T\in\Sigma_{\mathcal{S}}}$, where $\Sigma_{\mathcal{S}}$ is the set of all $\mathcal{S}$-theories, is called a \textbf{Lindenbaum bundle}\index{matrix!Lindenbaum bundle} relative to $\mathcal{S}$ and the atlas {\em$\mat{Lin}[\Sigma_{\mathcal{S}}]:=\langle\FormAl,\Sigma_{\mathcal{S}}\rangle$}, is called a \textbf{Lindenbaum atlas}\index{matrix!Lindenbaum atlas} relative to $\mathcal{S}$.
\end{defn}

We say that a given matrix (bundle or atlas) is \textit{\textbf{adequate}}\index{matrix!adequate} for an abstract logic $\mathcal{S}$ if $\vdash_{\mathcal{S}}$ and the corresponding matrix relation coincide.\label{adequate-matrix} 

We define:
\[
\begin{array}{rl}
	X\models_{\textsf{\textbf{Lin}}[\Sigma_{\mathcal{S}}]}\alpha \stackrel{\text{df}}{\Longleftrightarrow}
	&\!\!\!\text{for any $\Lan$-substitution $\sigma$ and any $T\in\Sigma_{\mathcal{S}}$},\\
	&\!\!\!\text{$\sigma(X)\subseteq T$ implies $\sigma(\alpha)\in T$},
\end{array}
\]
for any $X\cup\lbrace\alpha\rbrace\subseteq\FormsL$.

Thus we conclude with the following.
\begin{cor}\label{C:Lindenbaum-completeness}\index{Theorem!Lindenbaum copleteness}
	The Lindenbaum atlas relative to a structural abstract logic is adequate for this logic.
	That is, for any {\em$X\cup\lbrace\alpha\rbrace\subseteq\FormsL$},
	{\em\[
		X\vdash_{\mathcal{S}}\alpha~\Longleftrightarrow~X\models_{\textsf{\textbf{Lin}}[\Sigma_{\mathcal{S}}]}\alpha.
		\]}
\end{cor}
\noindent\textit{Proof}~is left to the reader. (Exercise~\ref{section:realizations-abstract-logic}.\ref{EX:Lindenbaum-completeness})

\begin{cor}\label{C:S-matrix-completeness}
	Let $\mathcal{S}$ be a structural abstract logic and $\mathcal{M}$ be the set of all $\mathcal{S}$-models. Then $\mathcal{S}=\mathcal{S}_{\mathcal{M}}$.
\end{cor}
\noindent\textit{Proof}~is left to the reader. (Exercise~\ref{section:realizations-abstract-logic}.\ref{EX:S-matrix-completeness})\\

Any $\Lan$-atlas (or $\Lan$-bundle) defines a structural consequence relation. The converse is also true as the next proposition shows.
\begin{prop}\label{P:matrix-con=atlas-con}
	Let {\em$\mathcal{M}=\lbrace\langle\alg{A}_i,D_i\rangle\rbrace_{i\in I}$} be a nonempty set of $\Lan$-matrices. There is an atlas $\mathcal{M}^{\ast}$ such that for any set {\em$X\cup\lbrace\alpha\rbrace\subseteq\FormsL$},
	\[
	X\vdash_{\mathcal{M}}\alpha~\Longleftrightarrow~X\vdash_{\mathcal{M}^{\ast}}\alpha.
	\]
\end{prop}
\begin{proof}
	We define an atlas  $\mathcal{M}^{\ast}=\langle\alg{A},\lbrace D_{i}^{\ast}\rbrace_{i\in I}\rangle$ as follows.
	\[
	\alg{A}:=\prod_{i\in I}\alg{A}_i
	\]
	and for every $i\in I$,
	\[
	D_{i}^{\ast}:=\prod_{j\in I} H_j,
	\]
	where, given $i\in I$,
	\[
	H_j:=\begin{cases}
		\begin{array}{cl}
			D_i &\text{if $j=i$}\\
			|\alg{A}_j| &\text{if $j\neq i$}.
		\end{array}
	\end{cases}
	\]
	
	We aim to prove that
	\[
	X\not\vdash_{\mathcal{M}}\alpha~\Longleftrightarrow~X\not\vdash_{\mathcal{M}^{\ast}}\alpha.
	\]
	
	Assume first that $X\not\vdash_{\mathcal{M}}\alpha$, that is, for some valuation $v$ in $\alg{A}_{i_0}$, although $v[X]\subseteq D_{i_0}$, $v[\alpha]\notin D_{i_0}$. Since each set $|\alg{A}_i|$ is not empty, we select elements $a_i\in|\alg{A}_i|$, foe each $i\in I$. Next, we define a valuation $v^{\ast}$ in \alg{A} as follows:
	\[
	v^{\ast}[p]=(a_{i}(p))_{i\in I},~\text{for every $p\in\VarL$},
	\]
	where
	\[
	a_{i}(p):=\begin{cases}
		\begin{array}{cl}
			a_i &\text{if $i\neq i_0$}\\
			v[p] &\text{if $i=i_0$}.
		\end{array}
	\end{cases}
	\]
	
	It must be clear that
	\[
	v^{\ast}[X]\subseteq D^{\ast}_{i}~\text{and}~v^{\ast}[\alpha]\notin D^{\ast}_{i}.\tag{Pr~\ref{P:matrix-con=atlas-con}--$\ast$}
	\]
	
	We leave for the reader to prove (Pr~\ref{P:matrix-con=atlas-con}--$\ast$). (See Exercise~\ref{section:realizations-abstract-logic}.\ref{EX:matrix-con=atlas-con}.)\\
	
	Now, we suppose that at some valuation $v^\circ$ in \alg{A}, $v^{\circ}[X]\subseteq D^{\ast}_{j_0}$ and $v^{\circ}[\alpha]\notin D^{\ast}_{j_0}$, for some $j_0\in I$. This means that for any $\beta\in X$, the projection of $v^{\circ}[\beta]$ on the $j_{0}$th component of $\alg{A}$ belongs to $D_{j_0}$ but the projection of $v^{\circ}[\alpha]$ does not; in symbols
	\[
	(v^{\circ}[\beta])_{j_0}\in D_{j_0},~\text{for every $\beta\in X$, 
		but $(v^{\circ}[\alpha])_{j_0}\notin D_{j_0}$}.\tag{Pr~\ref{P:matrix-con=atlas-con}--$\ast\ast$}
	\]
	We leave for the reader to prove (Pr~\ref{P:matrix-con=atlas-con}--$\ast\ast$). (See Exercise~\ref{section:realizations-abstract-logic}.\ref{EX:matrix-con=atlas-con}.)
	
	Thus we obtain that $X\not\vdash_{\mathcal{M}}\alpha$.
\end{proof}

\subsection{Consequence relation via inference rules}\label{section:inference-rules}
Since we perceive the sentential formulas as forms of judgments (Section 2.1), our understanding of a rule of inference, at least in one  important case, will be a relation in the class of forms of formulas. As in the first event, where substitution plays a key role in seeing any sentential formula as an abstract representation of infinitely many judgments, it will play a similar role in the concept of inference rule, at least in this important class of inference rules. It is a way, by which we create an abstraction of infinitely many instances of one element that falls under the name of \textit{structural inference rule}\footnote{This notion (see definition below) should not be confused with the notion of a structural rule in structural proof theory; for the latter see, e.g., in~\cite{negri-von-plato2001}, section 1.3.}. 

\subsubsection{Rules and hyperrules}\label{section:rules-and-hyperrules}
An \textit{\textbf{inference rule}}\index{inference rule} is a set $R\subseteq\mathcal{P}(\FormsL)\times\FormsL$. An inference rule receives its significance, when it is a subset of some consequence relation, say $\vdash$. We say that a rule $R$ is \textit{\textbf{sound with respect to}}\index{inference rule!sound} $\vdash$ if for any $\langle\Gamma,\alpha\rangle\in R$, $\Gamma\vdash\alpha$.

Let $\mathcal{R}$ be a set of inference rules. Taking into account that $\mathcal{P}(\FormsL)\times\FormsL$ is the trivial consequence relation (relative to a language $\Lan$), we define
\begin{equation}\label{E:consequence-by-rules}
	\vdash_{[\mathcal{R}]}:=\bigcap\set{\vdash_{\mathcal{S}}}
	{\forall R\in\mathcal{R}.~\text{$R$ is sound with respect to $\mathcal{S}$}}.
\end{equation}
(We note that the right-hand set is nonempty.)

It is not difficult to see that $\vdash_{[\mathcal{R}]}$ is a consequence relation. (Exercise~\ref{section:realizations-abstract-logic}.\ref{EX:S_R-consequence}) We say that an abstract logic $\mathcal{S}$ is \textit{\textbf{determined by}}\index{logic!determined by} a set $\mathcal{R}$ of inference rules if $\vdash_{\mathcal{S}}\,=\,\vdash_{[\mathcal{R}]}$. Given an abstract logic $\mathcal{S}$, one can treat the consequence relation $\vdash_{\mathcal{S}}$ as an inference rule that determines $\mathcal{S}$.\\

Now we turn to structural rules. 

Let us consider a pair $\langle\Gamma,\phi\rangle\in\mathcal{P}(\Mform)\times\Mform$. We call the set
\[
R_{\langle\Gamma,\phi\rangle}:=\set{\langle\bm{\sigma}(\Gamma),\bm{\sigma}(\phi)\rangle}{\bm{\sigma}~\text{is an instantiation of}~\Mvar}
\]
a \textit{\textbf{structural rule}} (in $\Lan$) with the premises $\Gamma$ and conclusion $\phi$. Each pair $\langle\bm{\sigma}(\Gamma),\bm{\sigma}(\phi)\rangle$ is an \textit{\textbf{instance of the rule}}\index{inference rule!instance} $R_{\langle\Gamma,\phi\rangle}$ (with the premises $\bm{\sigma}(\Gamma)$ and conclusion $\bm{\sigma}(\phi)$). It is convenient (and we follow this practice) to identify the rule $R_{\langle\Gamma,\phi\rangle}$ with its label $\langle\Gamma,\phi\rangle$, using them interchangeably. It is customary to write a structural rule $\langle\Gamma,\phi\rangle$ in the forms $\Gamma\slash\phi$ or
$\frac{\Gamma}{\phi}$. (It is customary to omit the set-builder notation ``$\{\ldots\}$'', when a set of premises is finite, and leave the space blank if $\Gamma=\varnothing$.)

Given an abstract logic $\mathcal{S}$, specifically, a structural rule $\langle\Gamma,\phi\rangle$ is \textit{\textbf{sound with respect to}} $\mathcal{S}$ if for any instantiation $\bm{\sigma}$, $\bm{\sigma}(\Gamma)\vdash_{\mathcal{S}}\bm{\sigma}(\phi)$.

A structural rule is called \textit{\textbf{finitary}}\index{inference rule!finitary} (or a \textit{\textbf{modus}}\index{inference rule!modus} or \textit{\textbf{modus rule}}) if $\Gamma\Subset\Mform$. If an abstract logic $\mathcal{S}$ is determined by a finite set of modus rules, one says that set of modus rules is a \textit{\textbf{calculus}}\index{calculus} for $\mathcal{S}$. An abstract logic may be defined by more than one calculus, if any. An example of a modus rule is \textit{modus ponens} which is the rule
$\bm{\alpha},\bm{\alpha}\rightarrow\bm{\beta}\slash\bm{\beta}$.

Not all useful inference rules are structural.

Let $\sigma$ be a fixed $\Lan$-substitution. Then we define
\[
R_{\sigma}:=\set{\langle\lbrace\alpha\rbrace,\sigma(\alpha)\rangle}{\alpha\in\FormsL}.
\]
Then the \textit{\textbf{rule of substitution}}\index{rul!of substitution} is the set
\[
R_{\text{sub}}:=\bigcup\set{R_{\sigma}}{\sigma~\text{is an $\Lan$-substitution}}.
\]
Each $\langle\alpha,\sigma(\alpha)\rangle$ is called an \textit{\textbf{application of the substitution}}\index{substitution!application} $\sigma$ \textit{\textbf{to}} $\alpha$. A customary utterance is: $\sigma(\alpha)$ \textit{is obtained from $\alpha$ by substitution} $\sigma$. 

Although substitution is an important component of any structural rule (see Proposition~\ref{P:metaformula-instantiations}), the rule of substitution itself is not structural, providing that $\Lan$ contains sentential connectives.

Indeed, for contradiction, we assume that 
\[
R_{\text{sub}}=R_{\langle\Gamma,\phi\rangle},
\]
for some $\Gamma\cup\lbrace\phi\rbrace\subseteq\Mform$.
Let us fix an arbitrary $p\in\Var_{\Lan}$. By definition, it is clear that for any $\alpha\in\FormsL$, $\langle\lbrace p\rbrace,\alpha\rangle\in R_{\text{sub}}$.
This implies that for any $\alpha\in\FormsL$, there is an instantiation $\bm{\xi}_{\alpha}$ such that $\bm{\xi}_{\alpha}(\Gamma)=\lbrace p\rbrace$ and $\bm{\xi}_{\alpha}(\phi)=\alpha$. Thus, we conclude that $\Gamma\cup\lbrace\phi\rbrace\subseteq\Mvar$ and $\Gamma\cap\lbrace\phi\rbrace=\varnothing$; moreover, it must be clear that $\Gamma$ cannot have more than one metavariables, because if it were otherwise, the claim 
$\langle\lbrace p\rbrace,\alpha\rangle\in R_{\langle\Gamma,\phi\rangle}$ would not be true. Thus, $\langle\Gamma,\phi\rangle=\langle\lbrace\gamma\rbrace,\beta\rangle$, for two distinct metavariables $\gamma$ and $\beta$. Then, providing that $\Lan$ contains sentential connectives, there is an instantiation $\bm{\xi}$ such that the degree of $\bm{\xi}(\gamma)$ is greater than the degree of $\bm{\xi}(\beta)$. A contradiction.\\

In the sequel (Chapter~\ref{chapter:L-equational-logic}), we will discuss other non-structural rules of inference.

\begin{prop}\label{P:structural-rules-for-structural-logic}
	Let $\mathcal{S}$ be a structural abstract logic in $\Lan$. Then there are structural rules that determine $\mathcal{S}$.
\end{prop}
\begin{proof}
	Since $\card{\Var_{\Lan}}\le\card{\Mvar}$, there is a one-one map $\Var_{\Lan}\longrightarrow\Mvar$ which can, obviously, be extended to $f:\FormAl\longrightarrow\MformAl$. (We employ the same notation for both maps.)
	
	Next, we define:
	\[
	\mathcal{R}:=\set{\langle f(X),f(\alpha)\rangle}{X\vdash_{\mathcal{S}}\alpha}.
	\]
	
	It must be clear that any rule $R\in\mathcal{R}$ is structural. And if $R=\langle f(X),f(\alpha)\rangle$, then there is an instantiation $\bm{\xi_0}$ such that $\bm{\xi_0}(f(X))=X$
	and $\bm{\xi_0}(f(\alpha))=\alpha$. This implies that any rule $R\in\mathcal{R}$ is sound with respect to $\mathcal{S}$. Hence $\vdash_{[\mathcal{R}]}\subseteq\vdash_{\mathcal{S}}$.
	
	Now we assume that any rule $R\in\mathcal{R}$ is sound with respect to some $\mathcal{S}^{\prime}$. Suppose $X\vdash_{\mathcal{S}}\alpha$. Then 
	$\langle f(X),f(\alpha)\rangle\in\mathcal{R}$. Therefore, for any instantiation $\bm{\eta}$, $\bm{\eta}(f(X))\vdash_{\mathcal{S}^{\prime}}\bm{\eta}(f(\alpha))$.
	In particular, $\bm{\xi_0}(f(X))\vdash_{\mathcal{S}^{\prime}}\bm{\xi_0}(f(\alpha))$.
	Hence $X\vdash_{\mathcal{S}^{\prime}}\alpha$.
\end{proof}

For defining consequence relations, especially structural consequence relations,  sometimes even more powerful machinery is convenient. We mean the use of hyperrules.

An (n+1)-tuple
\[
\Theta:=\langle\langle\Gamma_1,\phi_1\rangle,\ldots,\langle\Gamma_n,\phi_n\rangle,
\langle\Delta,\psi\rangle\rangle,
\]
where $\Gamma_{1}\cup\ldots\cup\Gamma_{n}\cup\Delta\cup\lbrace\phi_1,\ldots\phi_n,\psi
\rbrace\subseteq\Mform$, is called a (\textit{\textbf{structural}})\index{hyperrule!structural} \textit{\textbf{hyperrule}}. 

A hyperrule
$\Theta$ is \textit{\textbf{sound with respect to}}\index{hyperrule!sound} an abstract logic $\mathcal{S}$ if for any $\Lan$-instantiation $\bm{\xi}$ and any set $X\subseteq\FormsL$,
\[
[X\cup\,\bm{\xi}(\Gamma_{1})\vdash_{\mathcal{S}}\bm{\xi}(\phi_1),\ldots,X\cup\,\bm{\xi}(\Gamma_{n})\vdash_{\mathcal{S}}\bm{\xi}(\phi_n)]\Longrightarrow
X\cup\,\bm{\xi}(\Delta)\vdash_{\mathcal{S}}\bm{\xi}(\psi).
\]

If a hyperrule $\Theta$ is employed in relation to an abstract logic $\mathcal{S}$, it is customary to write this hyperrule in the following form:
\[
\dfrac{X,\Gamma_{1}\vdash_{\mathcal{S}}\phi_1;\ldots; X,\Gamma_{n}\vdash_{\mathcal{S}}\phi_n}{X,\Delta\vdash_{\mathcal{S}}\psi}.
\] 

These are examples of rules and hyperrules. (Later in this section we will introduce other inference rules.) We will use them to define consequence relations.\\

\textsf{Inference rules:}
\[
\begin{array}{cllll}
	(\text{a}) &~~~~~i)~\dfrac{\bm{\alpha},\bm{\beta}}{\bm{\alpha}\wedge\bm{\beta}} &~~~~~ii)~\dfrac{\bm{\alpha}}{\bm{\alpha}\vee\bm{\beta}}
	&~~~~~iii)~\dfrac{\bm{\beta}}{\bm{\alpha}\vee\bm{\beta}} 
	&~~~~~iv)~\dfrac{\bm{\alpha}}{\neg\neg\bm{\alpha}}\\\\
	(\text{b}) &~~~~~i)~\dfrac{\bm{\alpha}\wedge\bm{\beta}}{\bm{\alpha}}
	&~~~~~ii)~\dfrac{\bm{\alpha}\wedge\bm{\beta}}{\bm{\beta}}
	&~~~~~iii)~\dfrac{\bm{\alpha},\bm{\alpha}\rightarrow\bm{\beta}}{\bm{\beta}}
	&~~~~~iv)~\dfrac{\neg\neg\bm{\alpha}}{\bm{\alpha}}.
\end{array}
\]

\textsf{Hyperrules:}
\[
\begin{array}{cccc}
	(\text{c}) &i)~\dfrac{X,\bm{\alpha}\vdash\bm{\beta}}{X\vdash\bm{\alpha}\rightarrow\bm{\beta}}
	&ii)~\dfrac{X,\bm{\alpha}\vdash\bm{\beta};~X,\bm{\alpha}\vdash\neg\bm{\beta}}{X\vdash\neg\bm{\alpha}}  &iii)~\dfrac{X,\bm{\alpha}\vdash\bm{\gamma};~X,\bm{\beta}\vdash\bm{\gamma}}
	{X,\bm{\alpha}\vee\bm{\beta}\vdash\bm{\gamma}},
\end{array}
\]
where $\vdash$ in the hyperrules above denotes a fixed binary relation on
$\mathcal{P}(\FormsL)\times\FormsL$.

\subsubsection{Consequence relations $\vdash_1$, $\vdash_2$, $\vdash_{\textsf{Cl}}$, $\vdash_{\textsf{Cl}^{\ast}}$, and $\vdash_3$}
First we define a relation $\vdash_1$ on $\mathcal{P}(\FormsL)\times\FormsL$ and then prove that it is a consequence relation. Then, we do the same for relations $\vdash_2$ and $\vdash_3$. In all three definitions, we use the concept of (\textit{\textbf{formal}}) \textit{\textbf{derivation}}\index{formal derivation} (or \textit{\textbf{formal proof}})\index{formal proof}. The relations being defined are written as $X\vdash_i\alpha$, where $i=1,2,3$. (For convenience we call $X\vdash_i\alpha$ a \textit{sequent}, even when we may not know yet that $\vdash_i$ is a consequent relation.) If $X\vdash_i\alpha$ is justified according to the clauses below, we say that \textit{$X$ derives $\alpha$}, or \textit{$\alpha$ is derivable from $X$} (all with respect to $\vdash_i$). It is for such a justification we use the notion of derivation.

\paragraph{Relation $\vdash_1$} Given a set $X$ of $\Lan_{A}$-formulas and an $\Lan_{A}$-formula $\alpha$, derivation confirming (or establishing or justifying) a sequent $X\vdash_1\alpha$ is a finite (nonempty) list of $\Lan_{A}$-formulas
\begin{equation}\label{E:list}
	\alpha_1,\ldots,\alpha_n
\end{equation}
such that
\[
\begin{array}{cl}
	(\text{d}_1) &\alpha_n=\alpha;~\text{and}\\
	(\text{d}_2) &\text{each $\alpha_i$ is either in $X$ or  obtained from one or more formulas}\\ 
	&\text{preceeding it on the list \eqref{E:list} by one of the rules (a)--(b)}.
\end{array}
\]

It is easy to see that $\vdash_1$ satisfies the properties (a)--(c) of Definition~\ref{D:consequnce-relation-single} and thus is a consequence relation. Moreover, this relation is obviously finitary and structural. One peculiarity of it is that the set of the theses of this relation is empty. Another property of $\vdash_1$ is that
\begin{equation}\label{E:classical-1}
	X\vdash_1\alpha\Longrightarrow X\models_{\textbf{B}_{2}}\alpha,
\end{equation}
where $\models_{\textbf{B}_{2}}$ is a single-matrix consequence with respect to matrix $\textbf{B}_{2}$ (Section~\ref{S:two-valued}).  To prove~\eqref{E:classical-1}, one needs to show that each rule of (a)--(b) preserves validity in $\textbf{B}_2$.
We leave this task to the reader. (Exercise~\ref{section:realizations-abstract-logic}.\ref{EX:classical-1})

Clearly, the converse of \eqref{E:classical-1} is not true, for if $\alpha$ is a classical tautology, then 
$\varnothing\models_{\textbf{B}_{2}}\alpha$ but $\varnothing\not\vdash_1\alpha$.

\paragraph{Relation $\vdash_{2}$} Our next example, relation $\vdash_2$, is a modification of $\vdash_1$. Namely we define:
\[
X\vdash_2\alpha\stackrel{\text{df}}{\Longleftrightarrow} X,L\textbf{B}_{2}\vdash_1\alpha,
\]
where $L\textbf{B}_{2}$ is the logic of matrix $\textbf{B}_{2}$ (Definition~\ref{D:validity-1}), that is the set of classical tautologies. According to Proposition~\ref{P:relative-consequence}, $\vdash_2$ is a consequence relation whose set of theses is $L\textbf{B}_{2}$. 

In addition, we claim that
\begin{equation}\label{E:classical-2}
	X\vdash_2\alpha\Longrightarrow X\models_{\textbf{B}_{2}}\alpha.
\end{equation}

The proof of \eqref{E:classical-2} follows from the proof of \eqref{E:classical-1}, where we have to add one more (obvious) case --- when $\alpha_i$ in \eqref{E:list} is a classical tautology, that is $\alpha_i\in L\textbf{B}_2$. 

The converse of \eqref{E:classical-2} is also true, that is
\begin{equation}\label{E:classical-2-converse}
	X\models_{\textbf{B}_{2}}\alpha\Longrightarrow X\vdash_2\alpha.
\end{equation}

To prove \eqref{E:classical-2-converse}, we first prove that the consequence relation $\models_{\textbf{B}_{2}}$ is finitary. Indeed, the premise $X\models_{\textbf{B}_{2}}\alpha$ is equivalent to that the set $X\cup\lbrace\neg\alpha\rbrace$ is not \textit{satisfiable} in $\textbf{B}_2$, that is for any valuation $v$ in $\textbf{B}_2$, $v[X\cup\lbrace\neg\alpha\rbrace]\not\subseteq\lbrace\one\rbrace$. 
Then, by compactness theorem,\footnote{Cf.~\cite{barwise1977}, theorem 4.2. See also~\cite{chang-keisler1990}, corollary 1.2.12.} $X\cup\lbrace\neg\alpha\rbrace$ is not \textit{finitely satisfiable} in $\textbf{B}_2$, that is, there is $X_{0}\Subset X$ such that $X_0\cup\lbrace\neg\alpha\rbrace$ is not satisfiable in $\textbf{B}_2$.
The latter is equivalent to $X_{0}\models_{\textbf{B}_{2}}$.\footnote{We will consider this argument in a more general setting in Section~\ref{section:finitary-matrix-consequence}.}

In the next step, having $X_{0}\models_{\textbf{B}_{2}}\alpha$ with $X_{0}\Subset X$, we show that $X_{0}\vdash_2\alpha$.

If $X_0=\varnothing$, then $\alpha$ is a classical tautology and the conclusion is obvious. Assume that
\[
X_0:=\lbrace\beta_1,\ldots,\beta_n\rbrace.
\]
We denote
\[
\wedge X_0:=(\ldots(\beta_1\wedge\beta_2)\wedge\ldots\wedge\beta_n)
\]
and observe that $\wedge X_0\rightarrow\alpha$ is a classical tautology, that is
$\wedge X_0\rightarrow\alpha\in L\textbf{B}_2$. Further, it is easy to see that
$X_0\vdash_2 \wedge X_0$. The two last observations imply that $X_0\vdash_2\alpha$ and hence $X\vdash_2\alpha$.

According to Proposition~\ref{P:matrix-con-is-con-relation} and the consideration above, $\models_{\textbf{B}_{2}}$, as well as $\vdash_2$, is a finitary structural consequence relation.\\

\paragraph{Relation $\vdash_{\textsf{Cl}}$} We begin with the remark that  $\vdash_2$ is not a calculus, for the set $L\textbf{B}_2$ of classical tautologies infinite. However, $\vdash_2$ can be reformulated in such a way that in a new formulation, we denote it by $\vdash_{\textsf{Cl}}$, the latter is a calculus. First of all, we remove from the definition of $\vdash_2$ all rules except $\alpha,\alpha\rightarrow\beta\slash\beta$ (known as \textit{modus ponens} or \textit{detachment}). The removed rules are actually replaced with rules without premises, that is rules of type $\varnothing\slash\varphi$, where $\varphi$ is one of the following forms of formulas (which are called  \textit{axiom schemata}):
\[
\begin{array}{rl}
	\text{ax1} &\bm{\alpha}\rightarrow(\bm{\beta}\rightarrow\bm{\alpha}),\\
	\text{ax2} &(\bm{\alpha}\rightarrow\bm{\beta})
	\rightarrow((\bm{\alpha}\rightarrow(\bm{\beta}\rightarrow\bm{\gamma}))\rightarrow
	(\bm{\alpha}\rightarrow\bm{\gamma})),\\
	\text{ax3} &\bm{\alpha}\rightarrow(\bm{\beta}\rightarrow(\bm{\alpha}\wedge\bm{\beta})),\\
	\text{ax4} &(\bm{\alpha}\wedge\bm{\beta})\rightarrow\bm{\alpha},\\
	\text{ax5} &(\bm{\alpha}\wedge\bm{\beta})\rightarrow\bm{\beta},\\
	\text{ax6} &\bm{\alpha}\rightarrow(\bm{\alpha}\vee\bm{\beta}),\\
	\text{ax7} &\bm{\beta}\rightarrow(\bm{\alpha}\vee\bm{\beta}),\\
	\text{ax8} &(\bm{\alpha}\rightarrow\bm{\gamma})\rightarrow((\bm{\beta}\rightarrow\bm{\gamma})
	\rightarrow((\bm{\alpha}\vee\bm{\beta})\rightarrow\bm{\gamma})),\\
	\text{ax9} &(\bm{\alpha}\rightarrow\bm{\beta})
	\rightarrow((\bm{\alpha}\rightarrow\neg\bm{\beta})\rightarrow\neg\bm{\alpha}), \\
	\text{ax10} &\neg\neg\bm{\alpha}\rightarrow\bm{\alpha}.
\end{array}
\]

The abstract logic \textsf{Cl}, corresponding to $\vdash_{\textsf{Cl}}$ is called \textit{\textbf{classical propositional logic}}\index{logic!classical propositional}. It is obvious that {\textsf{Cl} is a calculus. 
	\begin{prop}\label{P:Cl-is-finitary}
		The abstract logic {\em$\Cl$} is finitary.
	\end{prop}
	\noindent\textit{Proof}~is left to the reader. (Exercise~\ref{section:realizations-abstract-logic}.\ref{EX:Cl-is-finitary})\\	
	
	To prove that $\vdash_{2}~=~\vdash_{\textsf{Cl}}$, it suffices to show that
	\begin{equation}\label{E:B_2=Cl}
		L\booleTwo=\set{\alpha}{\varnothing\vdash_{\textsf{Cl}}\alpha}.
	\end{equation}
	
	For the last equality, we refer to~\cite{kleene1952}, {\S} 28, theorem 9, and {\S} 29, theorem $10^\circ$. We leave the remaining part of the proof to the reader. (Exercise~\ref{section:realizations-abstract-logic}.\ref{EX:B_2=Cl})
	
	We note the aforementioned removal of rules is not necessary; it would suffice to add to the rules (a)--(b) nine new rules of type $\varnothing\slash\varphi$, where $\varphi$ is one of the schemata ax1--ax10; in this case, the notion of derivation should be modified from $\vdash_1$ (namely the clause ($\text{d}_{2}$)) by including the new rules. We treat the forms ax1--ax10 as \textit{\textbf{premiseless modus rules}}\index{rule!premissless modus}. 
	
	\paragraph{Relation $\vdash_3$} For the next example, we employ the hyperrules (c). We state that $X\vdash_3\alpha$ when the latter can be confirmed by one of the following types of \textit{confirmation} or their combination.  A formal definition is the following.
	\begin{itemize}
		\item[$1^{\circ}$] The first type is a derivation of $\alpha$ from $X$ by using steps determined by ($\text{d}_1$)--($\text{d}_2$), that is a derivation confirming $X\vdash_1\alpha$.
		\item[$2^{\circ}$] The second type is when confirmed premises of one of the hyperrules (c) are satisfied for $\vdash_3$, we can confirm the conclusion of this hyperrule for $\vdash_3$. (This means an application of the hyperrule in question.)
		\item[$3^{\circ}$] The third type consists in the following.
		When we have confirmed $X\vdash_3\alpha_1,\ldots,X\vdash_3\alpha_n$, we confirm
		$X\vdash_3\alpha$ if $X,\alpha_1,\ldots,\alpha_n\vdash_1\alpha$.
	\end{itemize}
	\begin{rem}
		{\em
			In practice, it is assumed that for each confirmed premise of a hyperrule there is a deduction (confirming it), which is called in~\cite{kleene1952}, {\S} 21,  
			a \emph{subsidiary deduction}, and the deduction confirming the conclusion of the rule a \emph{resulting deduction}.}
	\end{rem}
	
	Thus, for any set $X$ of $\Lan_{A}$-formulas and any $\Lan_{A}$-formula $\alpha$, $X\vdash_3\alpha$ is true if, and only if, there is a consecutive finite number of confirmations of the types of $1^{\circ}$--$3^{\circ}$. In particular,
	\begin{equation}\label{E:classic-1-implies-classic-3}
		X\vdash_1\alpha\Longrightarrow X\vdash_3\alpha.
	\end{equation}
	
	We illustrate the last definition on the following example. Using $X=\lbrace\alpha,\beta\rbrace$ (for concrete formulas $\alpha$ and $\beta$), we confirm $X\vdash_3\alpha\wedge\beta$ by the derivation
	\[
	\alpha,\beta,\alpha\wedge\beta.
	\]
	Then, we rewrite the last sequent as $\lbrace\alpha\rbrace,\beta\vdash_3\alpha\wedge\beta$ and apply (c-$i$) to obtain $\lbrace\alpha\rbrace\vdash_3\beta\rightarrow(\alpha\wedge\beta)$ (confirmation of the second type). Now, we rewrite the last sequent as $\varnothing,\alpha\vdash_3\beta\rightarrow(\alpha\wedge\beta)$ and apply this hyperrule one more time. The resulting sequent is $\varnothing\vdash_3\alpha\rightarrow(\beta\rightarrow(\alpha\wedge\beta))$, or
	simply $\vdash_3\alpha\rightarrow(\beta\rightarrow(\alpha\wedge\beta))$.\\
	
	It should be clear that the relation $\vdash_3$ is monotone; that is
	\begin{equation}\label{E:classical-3-is-monotone}
		(X\vdash_3\alpha~\textit{and}~X\subseteq Y)~\Longrightarrow~Y\vdash_3\alpha.
	\end{equation}
	At the point, we do not know yet if $\vdash_3$ is a consequence relation. Our aim is to prove this. More than that, we shall show that $\vdash_3$ and $\vdash_2$ are equal as relations, though defined differently. We discuss advantage and disadvantage of $\vdash_3$ in relation to $\vdash_2$ at the end of this subsection.
	
	Aiming to prove $\vdash_3\,=\,\vdash_2$, we notice that the relation $\vdash_3$ enjoys the following variant of cut (the property $(\text{c}^{\ast})$ from Definition~\ref{D:consequnce-relation-single}).
	\begin{equation}\label{E:calculus-3-has-cut}
		(X,\beta\vdash_3\alpha\textit{ and }\vdash_3\beta)~\Longrightarrow~X\vdash_3\alpha.
	\end{equation}
	
	Indeed, applying the hyperrule of (c-$i$), we confirm that $X\vdash_3\beta\rightarrow\alpha$. Using \eqref{E:classical-3-is-monotone}, we obtain that $X\vdash_3\beta$. Finally, we use a confirmation of type $3^{\circ}$ to conclude that $X\vdash_3\alpha$.
	
	On the next step, we prove that
	\begin{equation}\label{E:cal-2-implies-cal-3}
		X\vdash_2\alpha\Longrightarrow X\vdash_3\alpha.
	\end{equation}
	
	Our proof of~\eqref{E:cal-2-implies-cal-3} is based on the following observation.
	\begin{equation}\label{E:classical-tautology-in-cal-3}
		\alpha\in L\textbf{B}_2~\Longrightarrow~\vdash_3\alpha.
	\end{equation}
	
	To prove~\eqref{E:classical-tautology-in-cal-3}, first, we show that for any $\Lan_{A}$-formula $\alpha$,
	\begin{equation}\label{E:excluded-middle}
		\vdash_3\alpha\vee\neg\alpha.
	\end{equation}
	
	Indeed,\footnote{The following steps 1--6 are the adaptation of the proof of $51^{\ast}$ from~\cite{kleene1952}, {\S} 27.}
	\[
	\begin{array}{cl}
		1. &\lbrace\neg(\alpha\vee\neg\alpha)\rbrace,\alpha\vdash_3 \alpha\vee\neg\alpha\quad[\text{for $\neg(\alpha\vee\neg\alpha),\alpha\vdash_1 \alpha\vee\neg\alpha$ and~\eqref{E:classic-1-implies-classic-3}}]\\
		2. &\lbrace\neg(\alpha\vee\neg\alpha)\rbrace,\alpha\vdash_3 \neg(\alpha\vee\neg\alpha)\quad[\text{for $\neg(\alpha\vee\neg\alpha),\alpha\vdash_1 \neg(\alpha\vee\neg\alpha)$ and~\eqref{E:classic-1-implies-classic-3}}]\\
		3. &\neg(\alpha\vee\neg\alpha)\vdash_3\neg\alpha\quad[\text{by (c-$ii$) with 1 and 2 as premises}]\\
		4. &\neg(\alpha\vee\neg\alpha)\vdash_3\neg\neg\alpha\quad[\text{analogously to 1--3}]\\
		5. &\vdash_3\neg\neg(\alpha\vee\neg\alpha)\quad[\text{by (c-$ii$) with 3 and 4 as premises}]\\
		6. &\vdash_3\alpha\vee\neg\alpha\quad[\text{by applying \eqref{E:calculus-3-has-cut} to 5 and $\neg\neg(\alpha\vee\neg\alpha)\vdash_1\alpha\vee\neg\alpha$}]
	\end{array}
	\]
	
	Also, we will use lemma 13 of~\cite{kleene1952}, {\S} 29, which we formulate in a paraphrased form as follows.
	\begin{quote}
		Let $\alpha$ be an $\Lan_{A}$-formula in the distinct propositional variables $p_1,\ldots,p_n$; let $v$ be a valuation in $\textbf{B}_2$; and let a sequence of $\Lan_{A}$-formulas is produced according to the rule: if $v(p_j)=\one$, then
		$\beta_j=p_j$, and if $v(p_j)=\zero$, then $\beta_j=\neg p_j$. Then 
		$\beta_1,\ldots,\beta_n\vdash_3\alpha$ or $\beta_1,\ldots,\beta_n\vdash_3\neg\alpha$, according as for the given valuation $v$ $\alpha$ takes the value $\one$ or the value $\zero$.
	\end{quote}
	
	Now, to prove ~\eqref{E:classical-tautology-in-cal-3}, we assume that $\alpha\in L\textbf{B}_2$ and $p_1,\ldots,p_n$ are all variables that occur in $\alpha$. Then, according to the last property, $p_1,\ldots,p_n\vdash_3\alpha$ and
	$p_1,\ldots,\neg p_n\vdash_3\alpha$. This, in virtue of (c-$iii$), implies that
	$p_1,\ldots,p_{n-1},p_{n}\vee\neg p_n\vdash_3\alpha$. Since, by~\eqref{E:calculus-3-has-cut}, $\vdash_3 p_{n}\vee\neg p_n$, we obtain that
	$p_1,\ldots,p_{n-1}\vdash_3\alpha$. Continuing the elimination of the variables $p_j$, we come to the conclusion that $\vdash_3\alpha$.  
	
	Now we obtain the proof of~\eqref{E:cal-2-implies-cal-3} as follows.
	
	Assume that $X\vdash_2\alpha$; that is for some $\alpha_1,\ldots,\alpha_n\in L\textbf{B}_2$, $X,\alpha_1,\ldots,\alpha_n\vdash_1\alpha$. According to \eqref{E:classic-1-implies-classic-3}, $X,\alpha_1,\ldots,\alpha_n\vdash_3\alpha$,
	where for each $\alpha_j$, in virtue of~\eqref{E:classical-tautology-in-cal-3},
	$\vdash_3\alpha_j$. Using \eqref{E:calculus-3-has-cut} $n$ times, we obtain that
	$X\vdash_3\alpha$.
	
	Finally, we show that
	\begin{equation}\label{E:cal-3-implies-cal-2}
		X\vdash_3\alpha\Longrightarrow X\vdash_2\alpha.
	\end{equation}
	The last implication will be proven if we show that
	\begin{equation}\label{E:cal-3-implies-B-2}
		X\vdash_3\alpha\Longrightarrow X\models_{\textbf{B}_{2}}\alpha
	\end{equation}
	and then apply~\eqref{E:classical-2-converse}. 
	
	We note  that~\eqref{E:cal-3-implies-B-2} will be proven if we establish that every rule of (a)--(b) and every hyperrule of (c) preserves the truth value $\one$. Let us check this property for the hyperrules.
	
	We begin with (c-$i$). Assume that for any valuation $v^{\prime}$ in $\textbf{B}_2$, 
	\[
	v^{\prime}[X,\alpha]\subseteq\lbrace\one\rbrace\Longrightarrow v^{\prime}[\beta]=\one.
	\]
	Suppose for some valuation $v$, $v[X]\subseteq\lbrace\one\rbrace$. We observe that if $v[\alpha]=\zero$, then $v[\alpha\rightarrow\beta]=\one$, regardless of the value of $v[\beta]$. However, if $v[\alpha]=\one$, then $v[\alpha\rightarrow\beta]=\one$, by premise. 
	
	Now we turn to (c-$ii$). Suppose that for any valuation $v^{\prime}$ in $\textbf{B}_2$, 
	\[
	v^{\prime}[X,\alpha]\subseteq\lbrace\one\rbrace\Longrightarrow
	\text{both $v^{\prime}[\beta]=\one$ and $v^{\prime}[\neg\beta]=\one$}.
	\]
	It is clear that this premise is equivalent to that for any valuation $v^{\prime}$ in $\textbf{B}_2$, 
	\[
	v^{\prime}[X,\alpha]\not\subseteq\lbrace\one\rbrace.
	\]
	Thus, if for some valuation $v$, $v[X]\subseteq\lbrace\one\rbrace$, then obviously
	$v[\alpha]=\zero$, that is $v[\neg\alpha]=\one$.
	
	To prove this preservation property for (c-$iii$), we suppose that for any valuation $v^{\prime}$ in $\textbf{B}_2$, both
	\[
	\begin{array}{l}
		v^{\prime}[X,\alpha]\subseteq\lbrace\one\rbrace\Longrightarrow v^{\prime}[\gamma]=\one~\text{and}\\
		v^{\prime}[X,\beta]\subseteq\lbrace\one\rbrace\Longrightarrow v^{\prime}[\gamma]=\one
	\end{array}
	\]
	hold. Now assume that for some valuation $v$, $v[X,\alpha\vee\beta]\subseteq\lbrace\one\rbrace$. It is clear that either $v[\alpha]=\one$ or $v[\beta]=\one$ (or both). In each case, we obtain that $v[\gamma]=\one$.
	
	Thus we have obtained the equality $\vdash_2\,=\,\vdash_3$. As a byproduct of this equality, we conclude that $\vdash_3$ is a structural consequence relation.\\
	
	\paragraph{Comparison of $\vdash_2$ and $\vdash_3$} 
	Although $\vdash_3$ is not a calculus, it is very convenient in confirming sequents that are true with respect to $\vdash_{\textsf{Cl}}$. For instance, we show that the schema ax8 above is a thesis of $\vdash_3$. Indeed, we have:
	\[
	\begin{array}{cl}
		1. &\lbrace\alpha\rightarrow\gamma,\beta\rightarrow\gamma\rbrace,\alpha
		\vdash_3\gamma;~\lbrace\alpha\rightarrow\gamma,\beta\rightarrow\gamma\rbrace,\beta
		\vdash_3\gamma\quad[\text{two $\vdash_1$-derivations}]\\
		2. &\lbrace\alpha\rightarrow\gamma,\beta\rightarrow\gamma\rbrace,
		\alpha\vee\beta\vdash_3\gamma\quad[\text{from 1 by (c-$iii$)}]\\
		3. &\lbrace\alpha\rightarrow\gamma\rbrace,\beta\rightarrow\gamma
		\vdash_3(\alpha\vee\beta)\rightarrow\gamma\quad[\text{from 2 by (c-$i$)}]\\
		4. &\alpha\rightarrow\gamma\vdash_3(\beta\rightarrow\gamma)\rightarrow((\alpha\vee\beta)\rightarrow\gamma)\quad[\text{from 3 by (c-$i$)}]\\
		5. &\vdash_3(\alpha\rightarrow\gamma)\rightarrow((\beta\rightarrow\gamma)\rightarrow((\alpha\vee\beta)\rightarrow\gamma)).\quad[\text{from 4 by (c-$i$)}]
	\end{array}
	\]
	However, there is a price for such a convenience. While in proving $X\vdash_{\textsf{Cl}}\alpha$ (or in proving $X\vdash_{2}\alpha$) a set $X$ stays the same in each application of a rule; in proving $X\vdash_3\alpha$, when we need to apply a hyperrule, the set $X$ in the premise(s) of a hyperrule may differ from the original $X$ in the sequent we are proving. In addition, in each application of a hyperrule, the $X$ in its premise(s) may vary from hyperrule to hyperrule or from one application to another of the same hyperrule.
	
	\paragraph{The idea of semantic completeness}\label{paragraph:idea-completeness}
	The equalities like $\vdash_2\,=\,\models_{\textbf{B}_{2}}$ or $\vdash_{\textsf{Cl}}\,=\,\models_{\textbf{B}_{2}}$ or $\vdash_3\,=\,\models_{\textbf{B}_{2}}$, but not $\vdash_2\,=\,\vdash_{\textsf{Cl}}$ or $\vdash_2\,=\,\vdash_3$, (see Section~\ref{section:rules-and-hyperrules}) are typical examples of semantic completeness. That is, a completeness result is claimed when one and the same consequence relation admits two different definitions --- one purely syntactic and the other with the help of semantic concepts, such as, e.g., logical matrices. The former is regarded as a convenient way to define an abstract logic, say $\mathcal{S}$, the latter as an effective way to show that for a particular set $X$ and formula $\alpha$, the sequent $X\vdash_{\mathcal{S}}\alpha$ does not hold. This issue will be further discussed in the next chapter (Section~\ref{section:separating-means}).
	
	To demonstrate, e.g., the validity of $\vdash_2\,=\,\models_{\textbf{B}_{2}}$, we proved the properties~\eqref{E:classical-2} and~\eqref{E:classical-2-converse}. The former is known in literature as \textit{soundness}, the latter as (properly) \textit{semantic completeness}. Thus we can state that the abstract logic $\vdash_2$ is \textit{complete with respect to} the logical matrix $\mathbf{B}_2$, or that
	$\vdash_2$ is \textit{determined by} $\mathbf{B}_2$.
	
	Although the main concept of this book is consequence relation, the equalities
	like 
	\[
	\begin{array}{c}
		\bm{T}_{\textsf{LC}}=L\mat{LC}~~\text{and}~~\bm{T}_{\textsf{LC}}=L\mat{LC}^{\ast}\
	\end{array}
	\]
	(compare with the $\LC$ completeness above) will attract our attention too.
	We call them also completeness results but in a weaker sense of the term ``completeness''.
	
	\paragraph{Exercises~\ref{section:realizations-abstract-logic}}
	\begin{enumerate}
		\item \label{EX:trivial-1}Prove that a structural logic $\mathcal{S}$ is trivial if, and only if, for any arbitrary variable $p$, $p\in\ConS{\varnothing}$.
		\item\label{EX:trivial-2} Prove that for any abstract logic $\mathcal{S}$, $\mathcal{S}$ is trivial if, and only if, $\bm{T}_{\mathcal{S}}$ is inconsistent.
		\item \label{EX:theory}Let $\mathcal{S}$ be a structural abstract logic. Prove that an $\mathcal{S}$-theory is closed under substitution if, and only if, this theory is generated by a set closed under substitution. Prove that $\bm{T}_{\mathcal{S}}$ is the least $\mathcal{S}$-theory closed under any substitution.
		\item \label{EX:semantic-consequence}Prove the properties $(\text{a}^\ast)$ and $(\text{b}^\ast)$ of Proposition~\ref{P:semantic-consequence}.
		\item\label{EX:semantic-consequence-2}Prove that, given a model $\fM$ and a consequence relation $\vdash$, if $\bm{T}_{\fM}^{\ast}\cup\bm{T}_{\fM}^{\ast\ast}\subseteq\,\vdash$, then
		$\models_{\fM}\subseteq\,\vdash$. (\emph{Hint}: use~\eqref{E:semantic-consequence-equivalent}.)
		\item\label{EX:semantic-consequence-2-equivalent}Prove the equality~\eqref{E:semantic-consequence-2-equivalent}. (\textit{Hint}: use the equality 
		\[
		\models_{\fM}=\bigcap\set{\vdash}{\vdash~\text{is a consequence relation with}~\bm{T}_{\fM}^{\ast}\cup\bm{T}_{\fM}^{\ast\ast}\subseteq\,\vdash};
		\] 
		see Proposition~\ref{P:semantic-consequence-2}.)
		\item\label{EX:M_subseteq_M^ast}Prove the inclusion~\eqref{E:M_subseteq_M^ast}.
		\item\label{EX:vdash_M=vdash_M^ast}Prove the equality~\eqref{E:vdash_M=vdash_M^ast}.
		\item \label{EX:demonstration} Let $\mathcal{S}$ be the logic given by the rule (b--$iii$). It is clear that $p,p\rightarrow q\vdash_{\mathcal{S}}q$ and $p,p\rightarrow q\vdash_{\aLogSub}r$. Prove that $p,p\rightarrow q\not\vdash_{\mathcal{S}}r$. This would demonstrate that the converse of the implication in Proposition~\ref{P:aLog-implies-aLogSub} in general does not hold.
		\item\label{EX:demonstration-2} Let $\mathcal{S}$ be the same as in the previous exercise. Show that $\mathcal{S}_{\text{sub}}$ is not trivial.
		\item\label{EX:matrix-con-is-con-relation} Prove Proposition~\ref{P:matrix-con-is-con-relation}.
		\item\label{EX:lindenbaum-matrix} Prove Corollary~\ref{C:lindenbaum-matrix}.
		\item \label{EX:Lindenbaum-completeness}
		Prove Corollary~\ref{C:Lindenbaum-completeness}
		\item \label{EX:S-matrix-completeness}
		Prove Corollary~\ref{C:S-matrix-completeness}
		\item \label{EX:matrix-con=atlas-con}Complete the proof of Proposition~\ref{P:matrix-con=atlas-con} by proving (Pr~\ref{P:matrix-con=atlas-con}--$\ast$) and (Pr~\ref{P:matrix-con=atlas-con}--$\ast\ast$).
		\item \label{EX:S_R-consequence} Show that $\vdash_{[\mathcal{R}]}$ is a consequence relation.
		\item \label{EX:S-as-a-rule}Prove that for any abstract logic $\mathcal{S}$, the consequence relation $\vdash_{\mathcal{S}}$ is determined by $\lbrace\vdash_{\mathcal{S}}\rbrace$.
		\item \label{EX:classical-1}Prove \eqref{E:classical-1}.
		\item \label{EX:Cl-is-finitary}Prove that the abstract logic $\Cl$ is finitary.
		\item\label{EX:B_2=Cl} Using~\eqref{E:B_2=Cl}, prove the equality $\vdash_{2}~=~\vdash_{\textsf{Cl}}$.
		\item Show that the schemata ax1--ax7 and ax9--ax10 are theses of $\vdash_3$.
		\item Although for any consequence operator $\textbf{Cn}$, $\Con{X\cap Y}\subseteq\Con{X}\cap\Con{Y}$ for any sets $X$ and $Y$, the converse inclusion is not necessarily true. Define such an abstract logic $\mathcal{S}$ that $\ConS{X}\cap\ConS{Y}\not\subseteq\ConS{X\cap Y}$ for some sets $X$ and $Y$.
	\end{enumerate}
	
	\section{Abstract logics defined by modus rules}\label{section:modus-rules}
	When we use or investigate calculi, we inevitably deal with modus rules.
	As we will see below (Corollary~\ref{C:logic-by-modus-rules-criterion}), for a structural logic to be finitary  and defined as calculus, are equivalent properties. On the other hand, the modus rules are naturally modeled in first order logic.
	
	\subsection{General characterization}\label{section:modus-rules-general}
	If $\mathcal{R}$ is a set of modus rules, there is an effective way to determine the consequence relation $\vdash_{[\mathcal{R}]}$ defined in~\eqref{E:consequence-by-rules} (Section~\ref{section:inference-rules}). For this, we define a relation $\vdash_{\mathcal{R}}$ on $\mathcal{P}(\FormsL)\times\FormsL$, by generalizing (for $\Lan$-formulas) the notion of derivation which earlier was defined for $\Lan_{A}$-formulas via ($\text{d}_1$)--($\text{d}_2$).
	
	The relation $X\vdash_{\mathcal{R}}\alpha$ holds, by definition, if there is a list
	\begin{equation}\label{E:list-2}
		\alpha_1,\ldots,\alpha_n,\alpha
	\end{equation}
	(the list $\alpha_1,\ldots,\alpha_n$ can be empty), where  each formula on \eqref{E:list-2} is 
	\[
	\begin{array}{cl}
		(\text{d}^{\ast}_{1}) &\text{either in $X$ or}\\
		(\text{d}^{\ast}_{2}) &\text{obtained from the formulas preceding it on \eqref{E:list-2} by one of the rules $\mathcal{R}$.}
	\end{array}
	\]
	We note that so defined $\vdash_{\mathcal{R}}$ is finitary. About~\eqref{E:list-2}, we say that it is an $\mathcal{R}$-\textit{derivation} of $\alpha$ from $X$.
	\begin{prop}\label{P:modus-rules}
		Let $\mathcal{R}$ be a set of modus rules. Then $\vdash_{\mathcal{R}}\,=\,\vdash_{[\mathcal{R}]}$. Hence $\vdash_{\mathcal{R}}$ is a consequence relation.
	\end{prop}
	\begin{proof}
		We shall prove successively the following three claims.
		\[
		\begin{array}{cl}
			1^{\circ} &\vdash_{\mathcal{R}}~\text{is a structural consequence relation};\\
			2^{\circ} &\text{for any $R\in\mathcal{R}$, $R$ is sound with respect to $\vdash_{\mathcal{R}}$};\\
			3^{\circ} &\text{if an abstract logic $\mathcal{S}$ is such that if for any $R\in\mathcal{R}$,}\\
			&\text{$R$ is sound with respect to $\mathcal{S}$, then $\vdash_{\mathcal{R}}\subseteq\vdash_{\mathcal{S}}$.}
		\end{array}
		\]
		
		We notice that the properties (a)--(b) of Definition~\ref{D:consequnce-relation-single} are obviously true for $\vdash_{\mathcal{R}}$; and we leave for the reader to check the property (c) and  structurality property. (Exercise~\ref{section:modus-rules}.\ref{EX:modus-rules})
		
		Next, suppose that $\langle\Gamma,\phi\rangle\in\mathcal{R}$. Let $\bm{\sigma}$ be an arbitrary instantiation in $\Lan$. Since $\Gamma$ is finite, we denote
		\[
		\bm{\sigma}(\Gamma):=\lbrace\alpha_1,\ldots,\alpha_n\rbrace.
		\]
		According to the clauses $(\text{d}^{\ast}_{1})$--$(\text{d}^{\ast}_{2})$, the list
		\[
		\alpha_1,\ldots,\alpha_n,\bm{\sigma}(\phi)
		\]
		is an $\mathcal{R}$-derivation of $\bm{\sigma}(\phi)$ from $\bm{\sigma}(\Gamma)$. Hence $\bm{\sigma}(\Gamma)\vdash_{\mathcal{R}}\bm{\sigma}(\phi)$.
		
		Finally, let $\mathcal{S}$ be an abstract logic and for any $R\in\mathcal{R}$, $R$ is sound with respect to $\mathcal{S}$. We show that then
		\[
		X\vdash_{\mathcal{R}}\alpha\Longrightarrow X\vdash_{\mathcal{S}}\alpha.
		\]
		
		Suppose~\eqref{E:list-2} is a derivation with respect to $\mathcal{R}$. We prove that for each formula $\alpha_i$ of~\eqref{E:list-2}, $X\vdash_{\mathcal{S}}\alpha_i$.
		We employ induction on the number of the formulas preceding $\alpha$ in~\eqref{E:list-2}.
		Indeed, we see that, if $\alpha=\alpha_{1}$, then $X\vdash_{\mathcal{S}}\alpha_1$, for either $\alpha\in X$ or $\alpha$ is a conclusion of a premiseless rule.
		
		Now, suppose $\alpha$ is obtained by application of a rule  $R_{\langle\Gamma,\bm{\beta}\rangle}$. That is, for some formulas $\alpha_{i_1},\ldots,
		\alpha_{i_k}$ preceding $\alpha$ and some instantiation $\bm{\sigma}$,
		$\bm{\sigma}(\Gamma)=\lbrace\alpha_{i_1},\ldots,
		\alpha_{i_k}\rbrace$ and $\bm{\sigma}(\bm{\beta})=\alpha$. Since, by induction, $X\vdash_{\mathcal{S}}\alpha_{i_j}$, $1\le j\le k$, and, by soundness of $R_{\langle\Gamma,\bm{\beta}\rangle}$ with respect to $\mathcal{S}$, $\bm{\sigma}(\Gamma)\vdash_{\mathcal{S}}\alpha$, we conclude that
		$X\vdash_{\mathcal{S}}\alpha$.
	\end{proof}
	
	Given a set $\mathcal{R}$ of modus rules, we denote by $\mathcal{S}_{\mathcal{R}}$ the abstract logic corresponding to $\vdash_{[\mathcal{R}]}$ (see the definition~\eqref{E:consequence-by-rules}) or, equivalently (due to Proposition~\ref{P:modus-rules}), to $\vdash_{\mathcal{R}}$.
	\begin{cor}\label{C:modus-rules-finitary}
		If $\mathcal{R}$ is a set of modus rules, then the abstract logic $\mathcal{S}_{\mathcal{R}}$ is structural and finitary.
	\end{cor}
	\noindent\textit{Proof}~is left to the reader. (Exercise~\ref{section:modus-rules}.\ref{EX:modus-rules-finitary})
	
	\begin{cor}\label{C:logic-by-modus-rules-criterion}
		Let $\mathcal{S}$ be a structural abstract logic. Then $\mathcal{S}$ is finitary if, and only if, there is a set $\mathcal{R}$ of modus rules such that $\mathcal{S}=\mathcal{S}_{\mathcal{R}}$.
	\end{cor}
	\begin{proof}
		The if-part follows from Corollary~\ref{C:modus-rules-finitary}.
		
		To prove the `only-if' part, we assume that $\mathcal{S}$ is finitary. Next we define:
		\[
		\mathcal{R}:=\set{\langle\Gamma,\phi\rangle}{\Gamma\cup\lbrace\phi\rbrace\Subset
			\Mform~\text{and $\bm{\xi}(\Gamma)\vdash_{\mathcal{S}}\bm{\xi}(\phi)$, for any instantiation $\bm{\xi}$}}.
		\]
		
		We aim to prove that $\mathcal{S}=\mathcal{S}_{\mathcal{R}}$.
		
		First, assume that $X\vdash_{\mathcal{S}}\beta$, for a particular set $X\cup\lbrace\beta\rbrace\subseteq\FormsL$.
		Since $\mathcal{S}$ is finitary, there is a set
		$Y\Subset X$ such that $Y\vdash_{\mathcal{S}}\beta$. Now, replacing the variables $\Var(Y\cup\lbrace\beta\rbrace)$ with metavariables according to arbitrary one-one correspondence, we obtain a finite set $\Gamma\cup\lbrace\phi\rbrace$ of metaformulas such that for any instantiation $\bm{\xi}$,
		$\bm{\xi}(\Gamma)\vdash_{\mathcal{S}}\bm{\xi}(\phi)$ and for a particular instantiation $\bm{\xi}_{0}$, $\bm{\xi}_{0}(\Gamma)=Y$ and 
		$\bm{\xi}_{0}(\phi)=\beta$. This implies that
		$\langle\Gamma,\phi\rangle\in\mathcal{R}$. The latter in turn allows us to conclude that $Y\vdash_{\mathcal{R}}\beta$ and hence $X\vdash_{\mathcal{R}}\beta$. 
		
		Now, we suppose that $X\vdash_{\mathcal{R}}\beta$. That is, there is an $\mathcal{R}$--derivation, determined by $(\text{d}^{\ast}_{1})$--$(\text{d}^{\ast}_{2})$:
		\[
		\alpha_1,\ldots,\alpha_n,\beta. \tag{$\ast$}
		\]
		
		It is obvious that the list $(\ast)$ is not empty. We prove (by induction) that for each formula $\gamma$ of $(\ast)$, $X\vdash_{\mathcal{S}}\gamma$.
		
		It is true for the first formula $\gamma$ in $(\ast)$, in which case $\gamma=\alpha_1=\beta$, for, then, $\gamma\in X$ and hence $X\vdash_{\mathcal{S}}\gamma$.
		
		Suppose we proved that for each $\alpha_{i}$, $X\vdash_{\mathcal{S}}\alpha_{i}$. If $\beta\in X$, we are done. If $\beta\notin X$, then $\beta$ is obtained from $\alpha_{i_1},\ldots,\alpha_{i_k}$ by a rule
		$\langle\Gamma,\phi\rangle\in\mathcal{R}$. By induction, this implies that $\lbrace\alpha_{i_1},\ldots,\alpha_{i_k}\rbrace
		\vdash_{\mathcal{S}}\beta$. By the transitivity property, we conclude that $X\vdash_{\mathcal{S}}\beta$.
	\end{proof}
	
	\subsection{Abstract logics $\mathcal{S}^{\star}$}\label{section:S^star}
	
	Let an abstract logic $\mathcal{S}$ be equal $\mathcal{S}_{\mathcal{R}}$, for some set $\mathcal{R}$ of modus rules. We define a relation $\vdash_{\mathcal{S}^{\star}}$ on the set $\mathcal{P}(\FormsL)\times\FormsL$ by adding to $(\text{d}_{1}^{\ast})$--$(\text{d}_{2}^{\ast})$ one more clause:
	any formula of`\eqref{E:list-2} can be obtained from one of the preceding formulas by (uniform) substitution.
	
	\begin{prop}\label{P:S=S^star}
		{\em$\vdash_{\mathcal{S}^{\star}}~=~\vdash_{\mathcal{S}_{\text{sub}}}$}. Hence {\em$\mathcal{S}^{\star}$} is a finitary abstract logic.
	\end{prop}
	\begin{proof}
		Suppose for some set $X\cup\lbrace\alpha\rbrace\subseteq\FormsL$, 
		$X\vdash_{\mathcal{S}^{\star}}\alpha$. This means that there is a derivation
		\[
		\alpha_{1},\ldots,\alpha_n,
		\]
		where $\alpha_n=\alpha$ and each $\alpha_{i}$ is either obtained by one of a given set of modus rules (which includes instantiations of axioms, if any) belong to $X$ or obtained by substitution from one of the preceding formulas in the derivation.
		
		It is clear that only a finite set $X_0$ of the formulas of $X$ may occur in such a derivation. According to a technique developed by B. Soboci{\'n}ski~\cite{sobocinski1974} (see also~\cite{lambros1979}), the above derivation can be reconstructed in such a way that a new derivation will have the last formula identical with $\alpha$ and, in addition, all substitution applications can only occur with the formulas  of $ X_0 $. This implies
		that $\Sb X_0 \vdash_{\mathcal{S}}\alpha$, that is $X_0 \vdash_{\mathcal{S}_{\text{sub}}}\alpha$; hence $X \vdash_{\mathcal{S}_{\text{sub}}}\alpha$.
		
		Conversely, assume that $X \vdash_{\mathcal{S}_{\text{sub}}}\alpha$, that is $\Sb X\vdash_{\mathcal{S}}\alpha$. Since $\mathcal{S}$ is finitary (Corollary~\ref{C:modus-rules-finitary}), there is a set $Y_0\Subset\Sb X$ such that $Y_0 \vdash_{\mathcal{S}}\alpha$. This implies that $X\vdash_{\mathcal{S}^{\star}}\alpha$.
		
		The equality $\vdash_{\mathcal{S}^{\star}}~=~\vdash_{\mathcal{S}_{\text{sub}}}$ in turn implies that $\mathcal{S}^{\star}$ is an abstract logic (because $\mathcal{S}_{\text{sub}}$ is); and from the definition of $\mathcal{S}^{\star}$ we derive that the latter is finitary.
	\end{proof}
	
	Although a logic $\mathcal{S}$ defined by modus rules is always structural (Corollary~\ref{C:modus-rules-finitary}), $\mathcal{S}^{\star}$ is not necessarily so. We conclude this subsection with the following example.
	
	Since the abstract logic $\Cl$  is defined by modus rules, we define an abstract logic $\Cl^{\star}$. The latter logic is not structural. (See Exercise~\ref{section:modus-rules}.\ref{EX:Cl^star-is-not-structural}.) It is noteworthy that $\Cl^{\star}$ has only two theories. Namely the following proposition holds.
	\begin{prop}\label{P:Cl^star-theories}
		For any set {\em$X\subseteq\Forms_{\Lan_{A}}$},
		{\em\[
			\textbf{Cn}_{\Cl^{\star}}(X)=\begin{cases}
				\begin{array}{cl}
					\bm{T}_{\Cl} &\text{if $X\subseteq	\bm{T}_{\Cl}$}\\
					\Forms_{\Lan_{A}} &\text{otherwise}.
				\end{array}
			\end{cases}	
			\]}
	\end{prop}
	\noindent\textit{Proof}~is left to the reader. (Exercise~\ref{section:modus-rules}.\ref{EX:CL^star-theories})
	
	\subsection{The class of $\mathcal{S}$-models}\label{section:S-models}
	The class of $\mathcal{S}$-models, when $\aLog$ is given semantically, was discussed in~Section~\ref{section:realizations-abstract-logic}.
	In this subsection we focus on the class of $\mathcal{S}$-models, when an abstract logic $\mathcal{S}$ can be determined by a set of modus rules. Given a structural abstract logic
	$\mathcal{S}$ in a language $\Lan$, we employ for our characterization of $\mathcal{S}$-models a special language of first order. This language depends on $\Lan$ and we denote it by $\textbf{FO}_{\Lan}$.
	
	The denumerable set of the individual variables of $\textbf{FO}_{\Lan}$ consists of the
	symbols $\bm{x_{\alpha}}$, for each metavariable $\bm{\alpha}$. The functional
	connectives coincide with the connectives of $\Lan$. Also, there is one unary predicate
	symbol, $D$, two logical connectives, $\Land$ (conjunction) and $\Rarrow$
	(implication), universal quantifier $\forall$, as well the parentheses, `(' and `)'. The rules of formula formation are usual, as well as is usual the notion of an $\textbf{FO}_{\Lan}$-structure (including the notions of \textit{satisfactory} and \textit{validity}). If $\Theta$ is an $\textbf{FO}_{\Lan}$-formula containing only free variables, its universal closure is denoted by $\forall\ldots\forall\Theta$.
	Main reference here is~\cite{chang-keisler1990}.\\
	
	It must be clear that if $\phi$ is a metaformula relevant to $\Lan$, then, replacing each metavariable $\bm{\alpha}$ in $\phi$ with the individual variable $\bm{x_\alpha}$, we obtain a term
	of $\textbf{FO}_{\Lan}$. We denote this term (obtained from $\phi$) 
	by $\phi^{\ast}$.
	
	Further, any modus rule $R:=\dfrac{\psi_1,\ldots,\psi_n}{\phi}$ is translated to the formula
	\begin{equation}\label{E:horn-formula-1}
		R^{\ast}:=\Forall\ldots\Forall((D(\psi_{1}^{\ast})\Land\ldots\Land D(\psi_{n}^{\ast}))\Rarrow D(\phi^{\ast}))
	\end{equation}
	and any premiseless rule $T:=\varnothing\slash\phi$ to the formula
	\begin{equation}\label{E:horn-formula-2}
		T^{\ast}:=\Forall\ldots\Forall D(\phi^{\ast}).
	\end{equation}
	
	The formulas of types~\eqref{E:horn-formula-1} and~\eqref{E:horn-formula-2} belong to the class of \textit{basic Horn sentences}; cf.~\cite{chang-keisler1990}, section 6.2.\\
	
	It is clear that any matrix $\mat{M}=\langle\alg{A},D\rangle$ of type $\Lan$ can be regarded as a model for $\textbf{FO}_{\Lan}$-formulas.  We denote the \textit{validity} of a closed formula $\Theta$ in $\mat{M}$ by $\mat{M}\Vdash\Theta$.
	
	Let $\psi$ be a metaformula of $\Mform$ containing only metavariables $\alpha_1,\ldots,\alpha_n$. We note that there is a strong correlation between the standard valuation of the term $\psi^{\ast}$ (containing only individual variables $x_{\alpha_1}, \ldots,x_{\alpha_n}$) in a first order structure $\mat{M}$ and a valuation of a formula $\bm{\xi}(\psi)$, obtained from $\psi$ by a simple instantiation $\bm{\xi}$, in a logical matrix $\mat{M}$. For instance, let us denote the value of $\psi^{\ast}$ on a sequence $a_1,\ldots,a_n$ in $\mat{M}$ (as a first order structure) by $\psi^{\ast}[a_1,\ldots,a_n]$. Now let $\bm{\xi}$ be a simple instantiation with
	$\bm{\xi}(\alpha_i)=p_i$. Then, if $v$ is a valuation in $\alg{A}$ with $v[p_i]=a_i$,
	then $\psi^{\ast}[a_1,\ldots,a_n]=v[\bm{\xi}(\psi)]$. This can be proven by induction on the degree of $\psi$.
	
	Grounding on the last conclusion, we obtain The first order formula $D(\psi^{\ast})$ is satisfied on the sequence $\alpha_1,\ldots,\alpha_n$ in the first order structure $\mat{M}$ if, and only if, $v[\bm{\xi}(\psi)]\in D$ in the logical matrix $\mat{M}$.
	
	This argument leads to the following lemma.
	
	\begin{lem}\label{L:R-and-R*-in-matrix}
		Let {\em$\mat{M}=\langle\alg{A},D\rangle$} be a matrix of type $\Lan$. For any modus rule $R$, $R$ is sound with respect to {\em$\models_{\textbf{M}}$} if, and only if, {\em$\mat{M}\Vdash R^{\ast}$}.
	\end{lem}
	\noindent\textit{Proof}~is left to the reader. (Exercise~\ref{section:modus-rules}.\ref{EX:R-and-R*-in-matrix}) 
	
	\begin{prop}\label{P:S-models-modus-rules}
		Let $\mathcal{R}$ be a set of modus rules in $\Lan$. Then for any $\Lan$-matrix {\em$\mat{M}$}, {\em$\mat{M}$} is an $\mathcal{S}_{\mathcal{R}}$-model if, and only if, 
		{\em$\mat{M}\Vdash R^{\ast}$}, for any $R\in\mathcal{R}$.
	\end{prop}
	\begin{proof}
		We fix an $\Lan$-matrix $\mat{M}=\langle\alg{A},D\rangle$. In view of Proposition~\ref{P:modus-rules}, we can associate with the abstract logic $\mathcal{S}_{\mathcal{R}}$ the consequence relation $\vdash_{\mathcal{R}}$.
		
		We begin with the `only-if' part. Thus we assume that the inclusion $\vdash_{\mathcal{R}}\,\subseteq\,\models_{\textbf{M}}$ holds. This implies that any $R\in\mathcal{R}$ is sound with respect to $\models_{\textbf{M}}$. By Lemma~\ref{L:R-and-R*-in-matrix}, this in turn implies that $\mat{M}\Vdash R^{\ast}$, for any $R\in\mathcal{R}$.
		
		To prove the `if' part, we assume that the last conclusion holds. By Lemma~\ref{L:R-and-R*-in-matrix}, we derive that any $R\in\mathcal{R}$ is sound with respect to
		$\models_{\textbf{M}}$.
		
		Now, suppose that $X\vdash_{\mathcal{R}}\alpha$. That is, there is an $\mathcal{R}$-derivation~\eqref{E:list-2}. By induction on the length of~\eqref{E:list-2},
		one can prove that for any valuation $v$ in $\alg{A}$, if $v[X]\subseteq D$, then $v[\alpha]\in D$. We leave this task to the reader. (Exercise~\ref{section:modus-rules}.\ref{EX:S-models-modus-rules})
	\end{proof}
	\begin{cor}\label{C:S-models-modus-rules}
		Let $\mathcal{R}$ be a set of modus rules in $\Lan$. Then the class of all 
		$\mathcal{S}_{\mathcal{R}}$-models is closed under reduced products, in particular under ultraproducts.
	\end{cor}
	\begin{proof}
		Let $\mathcal{M}$ be the class of all $\mathcal{S}_{\mathcal{R}}$-models. According to Proposition~\ref{P:S-models-modus-rules}, an $\Lan$-matrix $\mat{M}\in\mathcal{M}$ if, and only if, for any $R\in\mathcal{R}$, $R\in\mathcal{R}$. This means that $\mathcal{M}$ is determined by a set of Horn formulas. It is well known that every Horn sentence is preserved under reduced products, in particular under ultraproducts; cf.~\cite{chang-keisler1990}, proposition 6.2.2.
	\end{proof}
	
	\subsection{Modus rules vs. non-finitary rules}\label{section:modus-rules-vs-rules}
	
	Proposition~\ref{P:modus-rules} and Corollary~\ref{C:modus-rules-finitary} show limitation of modus rules in defining consequent relations. It should not be surprising that there are consequence relations that cannot be defined by any set of modus rules. Below we discuss an example of the consequence relation of that kind; moreover, the theorems of that abstract logic, as we will show, coincide with the theorems of a calculus, \textsf{LC}. (See \textsf{LC} completeness below.)
	
	In the language $\Lan_{A}$, we formulate two more premiseless structural rules: 
	\[
	\begin{array}{cl}
		\text{ax11} &\bm{\beta}\rightarrow(\neg\bm{\beta}\rightarrow\bm{\alpha}),\\
		\text{ax12} &(\bm{\alpha}\rightarrow\bm{\beta})\vee(\bm{\beta}\rightarrow\bm{\alpha}).
	\end{array}
	\]
	
	A calculus \textsf{LC} is defined through derivations according to $(\text{d}_1)$--$(\text{d}_2)$ where in $(\text{d}_2)$, however, we use only the rules 
	ax1--ax9, ax11-ax12 and (\text{b}-$iii$) (modus ponens).
	
	We note two properties of \textsf{LC}. The first is that for any set $X_{0}\cup\lbrace\alpha\rbrace\Subset\Forms_{\Lan_{A}}$,
	\begin{equation}\label{E:LC-property-1}
		X_0\vdash_{\textsf{LC}}\alpha~\Longleftrightarrow~\wedge X_0\vdash_{\textsf{LC}}\alpha.
	\end{equation}
	(See Exercise~\ref{section:modus-rules}.\ref{EX:LC-property-1}.)
	
	The second property is known as the ``deduction theorem'' (for \textsf{LC}), which in its $\Rightarrow$-part, in essence, asserts that the hyperrule (c-$i$) is sound with respect to $\vdash_{\textsf{LC}}$. Namely, for any set $X\cup\lbrace\alpha,\beta\rbrace\subseteq\Forms_{\Lan_{A}}$,
	\begin{equation}\label{E:LC-deduction-theorem}
		X,\alpha\vdash_{\textsf{LC}}\beta~\Longleftrightarrow~X\vdash_{\textsf{LC}}\alpha\rightarrow\beta.
	\end{equation}
	(See Exercise~\ref{section:modus-rules}.\ref{EX:LC-deduction-theorem}.)
	
	In our consideration, we employ the single-matrix relation $\models_{\textbf{LC}}$  which is grounded on the logical matrix \textbf{LC}; see Section~\ref{section:dummett}. 
	
	We will need the following property.
	
	For any $\Lan_{A}$-formulas $\alpha$ and $\beta$,
	\[
	\vdash_{\textsf{LC}}\alpha~\Longleftrightarrow~\models_{\textbf{LC}}\alpha.
	\tag{$\LC$ completeness}
	\]
	
	The last equivalence is due to M. Dummett~\cite{dummett1959}, theorem 1.\\
	
	Next, we define a consequence relation, $\vdash_{\textsf{LC}}$ and, accordingly, an abstract logic $\LC$), analogously to the relation $\vdash_{\textsf{Cl}}$ (Section~\ref{section:consequence-defining}), but by using the rules ax1--ax9, ax11--ax12 and modus ponens (the rule (b-$iii$)) for the notion of \textsf{LC}-derivation defined by the clauses ($\text{d}_{1}^{\ast}$)--($\text{d}_{2}^{\ast}$) above. 
	
	In addition, we consider the non-modus structural rule $R^{\ast}:=R_{\langle\Gamma^{\ast},\bm{\alpha}_0\rangle}$ with
	\[
	\Gamma^{\ast}:=\set{(\bm{\alpha}_{i}\leftrightarrow \bm{\alpha}_{j})\rightarrow \bm{\alpha}_0}{0<i<j<\omega},
	\]
	where $\bm{\alpha}_i$, $0\leq i<\omega$, are pairwise distinct metavariables. Also, we define an infinite set of modus rules
	\[
	\mathcal{R}:=\set{\Gamma\slash\bm{\alpha}_0}{\varnothing\subset\Gamma\Subset\Gamma^{\ast}}.
	\]
	
	Then, two consequence relations, $\vdash^{\ast}$ and $ \vdash^{\star} $, are defined in $\mathcal{P}(\Forms_{\Lan_{A}})\times\Forms_{\Lan_{A}}$
	similarly to the notion of ${\textbf{LC}}$-derivation but also for $\vdash^{\ast}$ with $R^{\ast}$ as an additional inference rule and for $\vdash^{\star}$ with the rules of $\mathcal{R}$ as additional inference rules. 
	
	We aim to show that, although for any set $X_{0}\cup\lbrace\alpha\rbrace\Subset\Forms_{\Lan_{A}}$,
	\begin{equation}\label{E:finitary-equivalence}
		X_0\vdash^{\star}\alpha~\Longleftrightarrow~X_0\vdash^{\ast}\alpha,
	\end{equation}
	the consequence relation $\vdash^{\ast}$ cannot be defined by any set of structural modus rules, for, as we will show, it is not finitary; see Corollary~\ref{C:modus-rules-finitary}.
	
	As concerns~\eqref{E:finitary-equivalence}, it is obvious, for in any $\vdash^{\ast}$-derivation with a finite set $X_0$ as a set of premises,
	the rule $R^{\ast}$ can be applied only under such an instantiation $\bm{\xi}$, for which
	there is an instantiation $\bm{\zeta}$ such that \mbox{$\bm{\xi}(\Gamma^\ast)\slash\bm{\xi}(\bm{\alpha_{0}})=\bm{\zeta}(\Gamma)
		\slash\bm{\zeta}(\bm{\alpha})$}, for some $\Gamma\Subset\Gamma^{\ast}$; and vice versa.

	Finally, we demonstrate that the consequence relation $\vdash^{\ast}$ differs from $\vdash_{\textsf{LC}}$, for the latter is finitary (Corollary~\ref{C:modus-rules-finitary}), while the former cannot be defined by any set of modus rules. 
	
	For contradiction, assume that $\vdash^{\ast}$ is finitary. Let $\bm{\sigma}$ be an instantiation such that $\bm{\sigma}(\bm{\alpha}_{i})=p_i$, $0\le i<\omega$. We note that $\bm{\sigma}(\Gamma^{\ast})\vdash^{\ast}p_0$. By assumption, for some finite $X_0\Subset\bm{\sigma}(\Gamma^{\ast})$, $X_{0}\vdash^{\ast}p_0$. In virtue of~\eqref{E:finitary-equivalence}, $X_0\vdash_{\textsf{LC}}p_0$. This in turn,
	according to~\eqref{E:LC-property-1} and~\eqref{E:LC-deduction-theorem}, implies
	that $\vdash_{\textsf{LC}}\wedge X_{0}\rightarrow p_0$. By \textsf{LC} completeness,
	we obtain that $\models_{\textbf{LC}}\wedge X_{0}\rightarrow p_0$ and hence
	$X_{0}\models_{\textbf{LC}}p_0$. Thus, by monotonicity, we have that
	$\bm{\sigma}(\Gamma^{\ast})\models_{\textbf{LC}} p_0$. However, let $v$ be such a valuation that $v[p_0]$ is a pre-last element $a$ of $\mat{LC}$ with respect to $\leq$ and $v[p_i]\neq v[p_j]$, whenever $i<j$. Then, since
	$v[p_i\leftrightarrow p_j]\leq a$, we have: $v[\bm{\sigma}(\Gamma^{\ast})]\subseteq\lbrace\one\rbrace$, but $v[p_0]=a\neq\one$.
	
	\paragraph{Exercises~\ref{section:modus-rules}}
	\begin{enumerate}
		\item \label{EX:modus-rules-finitary}Prove Corollary~\ref{C:modus-rules-finitary}.
		\item \label{EX:modus-rules}Show that if $\mathcal{R}$ is a set of structural modus rules, then the relation $\vdash_{[\mathcal{R}]}$ satisfies the property (c) and structurality property of Definition~\ref{D:consequnce-relation-single}.
		\item \label{EX:Cl^star-is-not-structural} The logic $\Cl^{\star}$ is obtained from $\Cl$ as explained in Section~\ref{section:S^star}.
		Show that $\Cl^{\star}$ is not structural.
		\item \label{EX:CL^star-theories} Prove Proposition~\ref{P:Cl^star-theories}.
		\item \label{EX:LC-property-1}Prove the equivalence~\eqref{E:LC-property-1}.
		\item\label{EX:LC-deduction-theorem} Prove the equivalence~\eqref{E:LC-deduction-theorem} (the deduction theorem for \textsf{LC}).
		\item \label{EX:R-and-R*-in-matrix}Prove Lemma~\ref{L:R-and-R*-in-matrix}.
		\item \label{EX:S-models-modus-rules}Finish the proof of Proposition~\ref{P:S-models-modus-rules}.
	\end{enumerate}
	
	\section{Extensions of abstract logics}\label{section:extensions}
	In this section we consider interaction of abstract logics defined in different formal languages. We begin with the following definition.
	
	\begin{defn}[extension, conservative extension]\label{D:conservative-extension}\index{language!conservative extension}
		Let a language $\Lan^{+}$ be an extension of a language $\Lan$. Also, let $\mathcal{S}^{+}$ and $\mathcal{S}$ be abstract logics in the languages $\Lan^{+}$ and
		$\Lan$, respectively. The abstract logic $\mathcal{S}^{+}$ is an \textbf{extension of}  $($or \textbf{over}$)$ 
		$\mathcal{S}$ if for any set $X$ of  $\Lan$\!-formulas,
		{\em\[
			\textbf{Cn}_{\mathcal{S}}(X)\subseteq\textbf{Cn}_{\mathcal{S}^{+}}(X).
			\]}
		We call $\mathcal{S}^{+}$ a \textbf{conservative extension of} $($or \textbf{over}$)$ $\mathcal{S}$ if for any set $X$ of $\Lan$\!-formulas,
		{\em\begin{equation}\label{E:conservative-extension}
				\textbf{Cn}_{\mathcal{S}^{+}}(X)\cap\FormsL=\textbf{Cn}_{\mathcal{S}}(X).
		\end{equation}}
	\end{defn}
	
	The following two observation are almost obvious.
	\begin{prop}\label{P:trivial-extension}
		Let a language $\Lan^{+}$ be an extension of a language $\Lan$. Assume that $\mathcal{S}$ is a structural logic in $\Lan$ and $\mathcal{S}^{+}$ is its structural extension in $\Lan^{+}$. Then if $\mathcal{S}$ is trivial, then so is $\mathcal{S}^{+}$.
	\end{prop}
	\begin{proof}
		If $\mathcal{S}$ is trivial, then $p\in\ConS{\varnothing}$, for an arbitrary $p\in\Var_{\Lan}$. Since $\mathcal{S}^{+}$ is an extension of $\mathcal{S}$, $p\in\textbf{Cn}_{\mathcal{S}^{+}}(\varnothing)$. And since $\mathcal{S}^{+}$ is structural, the latter implies that $\mathcal{S}^{+}$ is trivial.
	\end{proof}
	\begin{prop}\label{P:intersection-conservative-extensions}
		Let a language $\Lan^{+}$ be an extension of a language $\Lan$ and let $\mathcal{S}$ be an abstract logic in $\Lan$.
		The intersection of any nonempty set of conservative extensions of $\mathcal{S}$ in $\Lan^{+}$ is a conservative extensions of $\mathcal{S}$ in $\Lan^{+}$.
	\end{prop}
	\noindent\textit{Proof}~is left to the reader. (Exercise~\ref{section:extensions}.\ref{EX:intersection-conservative-extensions})\\

	In the sequel, we will need to keep in mind that, given a language $\Lan$ and its extension $\Lan^{+}$,  for any structural logic $\aLog$ in $\Lan$, there is
	a conservative structural extension $\aLog^{+}$.
	
	Indeed, assume that $\mathcal{S}$ is a structural logic in $\Lan$ and $\Lan^{+}$ is an extension of $\Lan$. According to Corollary~\ref{C:Lindenbaum-completeness}, the atlas $\Lin[\Sigma_{\mathcal{S}}]=\left\langle \FormAl,\Sigma_{\mathcal{S}}\right\rangle $ is adequate for $\mathcal{S}$.
	Let us select an arbitrary formula $\alpha\in\FormsL$ and expand the signature of $\FormAl$ to $\Lan^{+}$, denoting a new algebra by $\FormAl[\alpha]^{+}$, as follows. We interpret in $\FormAl[\alpha]^{+}$ any constant $c\in\Cons_{\Lan^{+}}\setminus\Cons_{\Lan}$ as the element $\alpha\in |\FormAl|$. Similarly, we interpret any term $F(\beta_1,\ldots,\beta_n)$, where $F\in\Func_{\Lan^{+}}\setminus\Func_{\Lan}$, as $\alpha$. Then, we define an abstract logic $\mathcal{S}^{+}$ as the one determined by the atlas $\left\langle \FormAl[\alpha]^{+},\Sigma_{\mathcal{S}}\right\rangle $.
	
	According to Proposition~\ref{P:matrix-con-is-con-relation}, $\mathcal{S}^{+}$ is structural. Further, $\mathcal{S}^{+}$ is a conservative extension of $\mathcal{S}$, because for any set $X\cup\lbrace\beta\rbrace\subseteq\FormsL$,
	\begin{equation}\label{E:S^+-conservativity}
		X\vdash_{\mathcal{S}^{+}}\beta~\Longleftrightarrow~X\models_{\textsf{\textbf{Lin}}[\Sigma_{\mathcal{S}}]}\beta.
	\end{equation}
	
	We leave the reader to prove~\eqref{E:S^+-conservativity}; see Exercise~\ref{section:extensions}.\ref{EX:S^+-conservativity}.
	
	In virtue of Proposition~\ref{P:trivial-extension}, $\mathcal{S}^{+}$ is trivial if $\mathcal{S}$ is trivial. \\
	
	We summarize as follows.
	\begin{prop}\label{P:S^+-conservativity}
		Let $\mathcal{S}$ be an structural abstract logic in $\Lan$ and $\Lan^{+}$ be an extension of $\Lan$. Then there is a conservative structural extension $\mathcal{S}^{+}$ of $\mathcal{S}$. Moreover, if $\mathcal{S}$ is trivial, $\mathcal{S}^{+}$ is also trivial.
	\end{prop}
	
	\begin{cor}
		Let a language $\Lan^{+}$ be an extension of a language $\Lan$. Any structural abstract logic in $\Lan$ has a least structural conservative extension in $\Lan^{+}$.
	\end{cor}
	\begin{proof}
		We use Proposition~\ref{P:S^+-conservativity} and  Proposition~\ref{P:intersection-conservative-extensions}. 	
	\end{proof}
	
	Our next definition will play a key part in the sequel.
	\begin{defn}[relation $\vdash_{(\mathcal{S})^{+}}$]\label{D:(S)^+}\index{$\vdash_{(\mathcal{S})^{+}}$}
		Let a language $\Lan^{+}$ be an extension of a language $\Lan$. Also, let $\mathcal{S}$ be an abstract logic in $\Lan$. We define a relation $\vdash_{(\mathcal{S})^{+}}$ in {\em$\mathcal{P}(\Forms_{\Lan^{+}})\times\Forms_{\Lan^{+}}$} as follows:
		{\em\[
			\begin{array}{rl}
				X\vdash_{(\mathcal{S})^{+}}\alpha\stackrel{\text{df}}{\Longleftrightarrow}
				\!\!\!&\textit{there is a set {\em$Y\cup\lbrace\beta\rbrace\subseteq\Forms_{\Lan}$} and an $\Lan^{+}$\!-substitution}\\
				&\textit{$\sigma$ such that $\sigma(Y)\subseteq X$, $\alpha=\sigma(\beta)$ and $Y\vdash_{\mathcal{S}}\beta$}.
			\end{array}
			\]}
	\end{defn}
	
	We find a justification for the last definition in the following proposition.
	\begin{prop}\label{P:(S)^+}
		Let a language $\Lan^{+}$ be an extension of a language $\Lan$ and $\mathcal{S}$ be a structural abstract logic in $\Lan$. If the relation $\vdash_{(\mathcal{S})^{+}}$ is a consequence relation, then $(\mathcal{S})^{+}$ is the least structural conservative extension of $\mathcal{S}$ in $\Lan^{+}$. In addition, $(\mathcal{S})^{+}$ is nontrivial if, and only if, $\mathcal{S}$ is nontrivial.
	\end{prop}
	\begin{proof}
		Assume that $\vdash_{(\mathcal{S})^{+}}$ is a consequence relation.
		Then it is obvious that $(\mathcal{S})^{+}$ is a conservative extension of $\mathcal{S}$.
		(The structurality of $\mathcal{S}$ should be used.)
		
		Next, we prove the structurality of $(\mathcal{S})^{+}$.
		For this, assume that $X\vdash_{(\mathcal{S})^{+}}\alpha$. By definition, there are a set $Y\cup\lbrace\beta\rbrace\subseteq\Forms_{\Lan^{+}}$ and an $\Lan^{+}$-substitution $\sigma$ such that $Y\vdash_{\mathcal{S}}\beta$, $\sigma(Y)\subseteq X$ and $\sigma(\beta)=\alpha$.
		Let $\xi$ be any $\Lan^{+}$-substitution. We observe: $\xi\circ\sigma(Y)\subseteq\xi(X)$ and $\xi\circ\sigma(\beta)=\xi(\alpha)$. Then, by definition, $\xi(X)\vdash_{(\mathcal{S})^{+}}\xi(\alpha)$.
		
		Now, let $\mathcal{S}^{+}$ be an arbitrary structural conservative extension of $\mathcal{S}$ in $\Lan^{+}$. Suppose that $X\vdash_{(\mathcal{S})^{+}}\alpha$. Then there is a set
		$Y\cup\lbrace\beta\rbrace\subseteq\FormsL$ and an $\Lan^{+}$-substitution $\sigma$ such that $Y\vdash_{\mathcal{S}}\beta$, $\sigma(Y)\subseteq X$ and $\sigma(\beta)=\alpha$. By premise, $Y\vdash_{\mathcal{S}^{+}}\beta$. Since 
		$\mathcal{S}^{+}$ is structural, we also have that $\sigma(Y)\vdash_{\mathcal{S}^{+}}\sigma(\beta)$. This implies that $X\vdash_{\mathcal{S}^{+}}\alpha$.
		
		Finally, the last claim of the proposition follows from the conservativeness of $(\mathcal{S})^{+}$ over $\mathcal{S}$.
	\end{proof}
	
	In the next proposition, we give a sufficient condition for $(\mathcal{S})^{+}$ to be a consequence relation.
	
	\begin{prop}\label{P:(S)^{+}-existence}
		Let $\Lan$ be a language with $\kappa=\card{\Var_{\Lan}}\ge\aleph_{0}$ and $\Lan^{+}$ be a primitive extension of $\Lan$. Assume that $\mathcal{S}$ is a structural abstract logic in $\Lan$. Then the relation $\vdash_{(\mathcal{S})^{+}}$ is a consequence relation in $\Lan^{+}$. Moreover, the abstract logic $(\mathcal{S})^{+}$ is $\kappa$-compact. Moreover, if $\mathcal{S}$ is finitary, then  $(\mathcal{S})^{+}$ is also finitary. 
	\end{prop}
	\begin{proof}
		First, we have to show that $\vdash_{(\mathcal{S})^{+}}$ satisfies the properties of (a)--(c) of Definition~\ref{D:consequnce-relation-single}.
		Leaving the check of the properties (a)--(b) to the reader, we turn to the property (c). (Exercise~\ref{section:extensions}.\ref{EX:(S)_plus})
		
		Assume that for a set $X\cup Y\cup Z\cup\lbrace\alpha\rbrace\subseteq\Forms_{\Lan^{+}}$,
		\[
		\begin{array}{cl}
			(\text{P1}) &X\vdash_{(\mathcal{S})^{+}}\gamma,~\text{for all $\gamma\in Y$};\\
			(\text{P2}) &Y,Z\vdash_{(\mathcal{S})^{+}}\alpha.
		\end{array}
		\]
		
		In virtue of (P2) and  Definition~\ref{D:(S)^+} (further in this proof, simply definition),  there are sets $Y_0$ and $Z_0$ with $Y_{0}\cup Z_{0}\cup\lbrace\beta_{0}\rbrace\subseteq\FormsL$ and an $\Lan^{+}$-substitution $\sigma_0$ such that $\sigma_0(Y_0)\subseteq Y$, $\sigma_0(Z_0)\subseteq Z$, $\sigma_0(\beta_{0})=\alpha$ and $Y_0, Z_0\vdash_{\mathcal{S}}\beta_{0}$. 
		
		We denote:
		\[
		\lbrace\gamma_i\rbrace_{i\in I}:=\sigma_0(Y_0).
		\]
		
		In virtue of (P1) and definition, for each $i\in I$, there is a set $X_{i}\cup\lbrace\beta_i\rbrace\subseteq\FormsL$ and an $\Lan^{+}$-substitution $\sigma_i$ such that $\sigma_i(X_i)\subseteq X$, $\sigma_i(\beta_i)=\gamma_i$ and
		$X_i\vdash_{\mathcal{S}}\beta_i$.
		
		Now the assumptions (P1)--(P2) can be refined as follows.
		\[
		\begin{array}{cl}
			(\text{P1}^{\ast}) &\bigcup_{i\in I}\sigma_i(X_i)\vdash_{(\mathcal{S})^{+}}\gamma_i,
			~\text{for all $i\in I$};\\
			(\text{P2}^{\ast}) &\lbrace\gamma_i\rbrace_{i\in I},\sigma_0(Z_0)\vdash_{(\mathcal{S})^{+}}\alpha.
		\end{array}
		\]
		
		We denote:
		\[
		X^{\ast}:=\bigcup_{i\in I}\sigma_i(X_i)~\text{and}~Z^{\ast}:=\sigma_0(Z_0).
		\]
		
		We aim to show that $X^{\ast}, Z^{\ast}\vdash_{(\mathcal{S})^{+}}\alpha$. Since $X^{\ast}\subseteq X$ and $Z^{\ast}\subseteq Z$, we will thus obtain that $X,Z\vdash_{(\mathcal{S})^{+}}\alpha$, for we have already proved that $\vdash_{(\mathcal{S})^{+}}$ is monotone.
		
		Let us make notes about the cardinalities of the $\Lan^+$-variables of some formula sets. (This is where we use that the extension $\Lan^{+}$ is primitive.)
		\begin{itemize}
			\item $\card{\Var(\lbrace\gamma_i\rbrace_{i\in I})}\le\card{\Var_{\Lan}}$;
			\item $\card{\Var(Z^\ast)}=\card{\Var(\sigma_0(Z_0))}\le\card{\Var_{\Lan}}$;
			\item $\card{\Var(X^\ast)}=\card{\Var(\bigcup_{i\in I}\sigma_i(X_i))}\le\card{\Var_{\Lan}}$.
		\end{itemize}
		
		Thus we have: $\card{\Var(\lbrace\gamma_i\rbrace_{i\in I}\cup X^{\ast}\cup Z^{\ast}\cup\lbrace\alpha\rbrace)}\le\card{\Var_{\Lan}}$. This allows us to claim that there is a one-one map
		\[
		f:\Var(\lbrace\gamma_i\rbrace_{i\in I}\cup X^{\ast}\cup Z^{\ast}\cup\lbrace\alpha\rbrace)
		\stackrel{\text{1--1}}{\longrightarrow}\Var_{\Lan}.
		\]
		
		For simplicity, we denote:
		\[
		\Var^{\ast}:=\Var(\lbrace\gamma_i\rbrace_{i\in I}\cup X^{\ast}\cup Z^{\ast}\cup\lbrace\alpha\rbrace).
		\]
		
		Next we define two $\Lan^+$-substitutions. Let $q$ be an arbitrary fixed $\Lan$-variable.
		Then, we define:
		\[
		\mu(p):= \begin{cases}
			\begin{array}{cl}
				f(p) &\text{if $p\in\Var^\ast$}\\
				q &\text{otherwise};
			\end{array}
		\end{cases}
		\]
		\[
		\nu(p):=\begin{cases}
			\begin{array}{cl}
				f^{-1}(p) &\text{if $p\in f(\Var^\ast)$}\\
				p &\text{otherwise}.
			\end{array}
		\end{cases}
		\]
		
		We observe the following:
		\begin{itemize}
			\item $\mu\circ\sigma_i$, for all $i\in I$, and $\mu\circ\sigma_0$ are $\Lan$-substitutions;
			\item $\nu\circ\mu\circ\sigma_i(X_i)=\sigma_i(X_i)$ and $\nu\circ\mu(\gamma_i)=\gamma_i$, for all $i\in I$, $\nu\circ\mu\circ\sigma_0(Z_0)=\sigma_0(Z_0)$ and 
			$\nu\circ\mu\circ\sigma_0(\beta_0)=\sigma_0(\beta_0)$;
			\item thus $\nu(\bigcup_{i\in I}\mu\circ\sigma_i(X_i))=\nu(\mu(X^\ast))=X^\ast$, $\nu(\mu\circ\sigma_0(Z_0))\nu(\mu(Z^\ast))=Z^\ast$ and $\nu(\mu\circ\sigma_0(\beta_0))=\nu(\mu(\alpha))=\alpha$. And, by definition, we conclude that $X^{\ast},Z^{\ast}\vdash_{(\mathcal{S})^{+}}\alpha$.
		\end{itemize}
		
		From this, since $\mathcal{S}$ is structural, we obtain:
		\begin{itemize}
			\item $\mu\circ\sigma_i(X_i)\vdash_{\mathcal{S}}\mu\circ\sigma_i(\beta_i)$, for each $i\in I$, which implies $\mu(\bigcup_{i\in I}\sigma_i(X_i))\vdash_{\mathcal{S}}\mu(\gamma_i)$, for each $i\in I$;
			\item $\lbrace\mu(\gamma_i)\rbrace_{i\in I},\mu\circ\sigma_0(Z_0)\vdash_{\mathcal{S}}\mu\circ\sigma_0(\beta_{0})$.
		\end{itemize}
		Then, from the two last assertions and the transitivity of $\mathcal{S}$, we derive:
		$\mu(X^{\ast}),\mu(Z^{\ast})\vdash_{\mathcal{S}}\mu(\alpha)$. 
	\end{proof}

	\paragraph{Exercises~\ref{section:extensions}}
	\begin{enumerate}
		\item \label{EX:intersection-conservative-extensions} Prove Proposition~\ref{P:intersection-conservative-extensions}
		\item \label{EX:S^+-conservativity}Prove the equivalence~\eqref{E:S^+-conservativity}, that is for any set $X\cup\lbrace\beta\rbrace\subseteq\Forms_{\Lan^{+}}$,
		\[
		X\models_{\langle\FormAl[\alpha]^{+},\Sigma_{\mathcal{S}}\rangle}\beta
		~\Longleftrightarrow~X\models_{\langle\FormAl,\Sigma_{\mathcal{S}}\rangle}\beta.
		\]
		
		(\textit{Hint}: Show that for any valuation $v$ in $\FormAl$, there is a valuation $v^{\ast}$ in $\FormAl[\alpha]^{+}$ such that $v^{\ast}\!\!\upharpoonright\!\!\FormsL=v$.)
		\item\label{EX:(S)_plus} Let $\Lan^{+}$ be a primitive extensions of a language $\Lan$ and let $\mathcal{S}$ be an abstract logic in $\Lan$. Show that the relation $\vdash_{(\mathcal{S})^{+}}$ of Definition~\ref{D:(S)^+} satisfies the properties (a)--(b) of Definition~\ref{D:consequnce-relation-single}.
		\item \label{EX:kappa-compactness_(S)-plus} Prove that $(\mathcal{S})^{+}$ in Proposition~\ref{P:(S)^{+}-existence} is $\kappa$-compact.
	\end{enumerate}
	
	\section{Historical notes}\label{section:consequence-historical-notes}
	
	As was mentioned in Section~\ref{section:key-concepts} and in the beginning of this chapter, the distinction in the argumentation of reasoning with predetermined premises from arguments with arbitrary premises goes back to Aristotle, who called the former a \emph{demonstrative} argument, and the latter a \emph{dialectical} one. It is very likely that the notion of demonstration has gradually developed in connection with geometry from ``empirical study'' to ``a demonstrative \emph{a priory} science;''~\cite{kneales1962}, chapter I, section 2. Dialectical arguments, on the other hand, were characteristic of metaphysics. In the early twentieth century, this distinction was emphasized by C. I. Lewis:
	\begin{quote}
		``The consequences of this difference between the ``implies'' of the algebra and the ``implies'' of valid inference are most serious.'' \cite{lewis1913a}, p. 242.~\footnote{It should be noted that Lewis did not distinguish between the “algebra of implication” and  ``calculus of propositions.” See, for example, the first line of~\cite{lewis1913a} or of~\cite{lewis1914}.}
	\end{quote}
	
	The formalization of logical consequence (a dialectical argument in the sense of Aristotle) in the framework of modern logic was undertaken by A. Tarski in the 1930s, by introducing the notion of consequence operator, not that of consequence relation; see~\cite{tarski1930a,tarski1930b}\footnote{The English translations of these important papers were later on published in~\cite{tarski1956,tarski1983}}. It seems that an earlier attempt by P. Hertz~\cite{hertz1922} was made from an essentially different perspective, since
	he interpreted `$(a_1,\ldots,a_n)\rightarrow b$' (where $a_1,,\ldots,a_n,b$ are sentences) as
	\begin{quote}
		``If $(a_1,\ldots,a_n)$ altogether holds, so does $b$,''\\ 
		(quoted from~\cite{antology2012}, p. 11)
	\end{quote}
	while Tarski would interpret `$\beta\in\Con{\alpha_{1},\ldots,\alpha_n}$' that $\beta$ is obtained from $\alpha_{1},\ldots,\alpha_n$ with respect to \textbf{Cn}, regardless whether any $\alpha_{i}$ (or they altogether) holds or not. The following passage shows that Hertz regarded `$\rightarrow$' as a relation of demonstration, and not as a dialectical argument when he wrote:
	\begin{quote}
		``Whenever a system of sentences is recognized to be valid, it is often not necessary to convey each or every sentence to memory; it is sufficient to choose some of them from which the rest can follow. Such sentences, as it is generally known, are called axioms.'' Quoted from~\cite{antology2012}, p. 11.
	\end{quote}
	
	Also, it seems that D. Scott was the first to start using consequence as a relation; cf.~\cite{scott1974}. He also stated explicitly therein that the two notions can be used interchangeably; cf.~\cite{shoesmith-smiley2008}, Historical note, pp. 19--20.\footnote{D. J. Shoesmith and T. J. Smiley also used consequence relation but did not relate it to Tarski's consequence operator; cf.~\cite{shoesmith-smiley1971}.} 
	
	It is worthwhile to dwell on the works of Tarski in the 1930s in a little more detail. As can be assumed from the title of~\cite{tarski1930b}, Tarski had in mind a far-reaching program. So he writes in~\cite{tarski1930b}, {\S} 1,
	\begin{quote}
		``\textit{The two concepts --- of sentence and of consequence --- are the only primitive concepts which appear in this discussion.}''\\ 
		(italics highlighted by Tarski)
	\end{quote}
	
	The first concept in Tarski's passage is a nonempty set of cardinality less than or equal to $\aleph_{0}$, and the second is a map, defined on the power set of all sentences and satisfying $(\text{a}^{\dagger})$--$(\text{c}^{\dagger})$ of Definition~\ref{D:consequence-operator}. Interestingly, Tarski admits, at least implicitly, that the set of sentences may be finite, a condition that was changed to `greater than or equal to $\aleph_{0}$' in later expositions; cf., e.g.,~\cite{czelakowski2001,fjp09,font2016}.
	
	It may seem surprising that there is no structure in Tarski's ``meaningful sentences'' (from the first concept). Perhaps he wanted to propose a concept of logical consequence independent of the grammatical structure of the objective language. Indeed, criticizing Carnap's definitions of logical consequence, he wrote,
	\begin{quote}
		``These definitions, in so far as they set up on the basis of `general syntax', seem to me to be materially inadequate, just because the defined concepts depend essentially, in their extension, on the richness of the language investigated.'' \cite{tarski1936b}
	\end{quote}
	However, the deprivation of sentences from a grammatical structure immediately blocks the Lindenbaum method, not allowing to define the formula algebra. 
	
	In~\cite{tarski1930a} and~\cite{tarski1930b}, the property $(\text{d}^{\dagger})$ (finitariness) is a mandatory axiom, while the structurality property is not mentioned at all, even as an option.
	Structurality is defined through the notion of substitution. However, substitution plays a dual role in logic.
	
	Although Aristotle was familiar with the idea of a variable, but he ``never argued by means of substitutions, as we have been doing, except in proofs of invalidity,''  that is, only using instantiation; ``for example, no thought of substituting variables for variables'' \cite{bochenski1970}, {\S} 13, sections B and C, respectively. Alexander of Aphrodisias, who lived at the turn of the fourth century AD, already knew the operation of identification of sentential variables and Boethius (c. 477--524 AD) knew substitution of sentential formulas for variables. He would probably be the first philosopher to agree, if he had a chance, that any map $\sigma:\FormsL\longrightarrow\FormsL$ satisfying the conditions of Definition~\ref{D:substitution} conveys the idea that logic deals with forms rather than with their instances.\footnote{I. M. Boche\'{n}ski calls this ``the clear distinction between propositional functions and propositions themselves;'' cf.~\cite{bochenski1970}, {\S} 24, section I.} Boethius also showed his ``aspiration to the rule of substitution'' but never used it as a rule of inference; cf.~\cite{bochenski1970}, {\S} 24, section E.
	
	G. Frege was perhaps the first to start using substitution as a rule of inference, though implicitly. This is what van Heijenoort wrote in his introductory note to Frege's \textit{Begriffsschrift} (1879), 
	\begin{quote}
		``His axioms for the propositional calculus (they are not independent) are formulas [\dots].  His rules of inference are the rule of detachment and an unstated rule of substitution.'' Quoted from~\cite{heijenoort2002}, p. 2.
	\end{quote}
	
	Although A. Church credits L. Couturat for being the first to explicitly formulate the substitution rule ``specifically for the propositional calculus,'' this is hardly justified. Indeed in~\cite{couturat1905},  the first edition of which was published in 1905, Couturat wrote,
	\begin{quote}
		``That is not all: to that principle one must add the substitution principle which is described thusly: ``In a general formula, to a general or indeterminate term one may substitute a particular or individual term.'' That is evident since a general formula has no value nor even meaning except in the way it can be applied to particular terms'' 
		\cite{couturat1905}, chapter 1.~\footnote{We are indebted to Dr. L. Wolffe for translating this passage from the French.}
	\end{quote}
	Fair enough, right after his remark, Church adds,``but his [Couturat's] statement is perhaps insufficient as failing to make clear that the expression substituted for a (propositional) variable may itself contain variables.''\footnote{Cf.~\cite{church1996}, {\S} 29.}
	
	B. Russell did use the substitution rule in~\cite{russell1906} in formal derivations of sentential formulas, but there is still confusion for the author in using the term ``substitution'' in the sense of Definition~\ref{D:substitution} and as universal instantiation in higher-order logic. After Russell, Lewis illustrates in~\cite{lewis1913b,lewis1914} the use of the substitution rule in the sense of Section~\ref{section:inference-rules}, but without such generality; and after~\cite{lewis1960} (the first publication appeared in 1918) substitution as an operation and as a rule of inference became standard.
	
	Yet, the structurality principle understood as implication
	\[
	X\vdash\alpha~\Longrightarrow~\sigma(X)\vdash\sigma(\alpha),
	\]
	for any substitution $\sigma$, does not come down to the rule of substitution since it is a characteristic of an entire deductive system with this property.
	The structurality as a important property of a consequence operator was introduced by J. {\L}o\'{s} and R. Suszko in~\cite{los-suszko1958}, although was anticipated by  Lewis in~\cite{lewis1913b}, when he wrote,
	\begin{quote}
		``Proof takes place through the collusion of two factors; first, postulates or propositions of the particular mathematical system in hand; secondly, postulates of the logical or mathematical type and the desired conclusion. A mathematical operation is ideally no more than this: the substitution of the variables or functions of variables of the particular system --- say, of cardinal number --- for the logical variables in some proposition about implication. This proposition is more than a rule for inference; when the substitution is made, it \emph{states} the implication.''
	\end{quote}
	
	We conclude the theme of substitution with the following passage by B. Russell in~\cite{russell1993} (the second edition appeared in 1920), where he wrote,
	\begin{quotation}
		``The primitive propositions, whatever they may be, are to be regarded as asserted for all possible values of the variable propositions $p,q,r$ which occur in them. We may therefore substitute for (say) $p$ any expression whose value is always a proposition, \emph{e.g.} not-$p$, ``$s$ implies $t$,'' and so on. By means of such substitutions we really obtain sets of special cases of our original proposition, but from a practical point of view we obtain what are virtually new propositions. The legitimacy of substitutions of this kind has to be insured by means of a non-formal principle of inference. [Footnote: No such principle is enunciated in \textit{Principia Mathematica} [\dots] But this would seem to be an omission.]'' (ibid, chapter XIV)
	\end{quotation}
	
	From this passage, we learn that although the authors of \textit{Principia Mathematica} were unable to appreciate the importance of the concept of substitution, by 1920 Russell already realized (or was close to this) that propositions themselves can be used as (truth) values for propositions, which leads directly to Lindenbaum method.\\
	
	The properties $(\text{e}^{\dagger})$ and $(\text{j}^{\dagger})$ of \textbf{Cn} in Section~\ref{section:consequence-operator} were considered earlier by D. Makinson in~\cite{makinson2005}, section 1.2.\\
	
	Given a consequence operator \textbf{Cn}, the notion of a theory (under the name ``deductive system'') with respect to \textbf{Cn} is due to~\cite{tarski1930a}, as well as Proposition~\ref{P:Tarski-criterion}, ibid, theorem 20.\footnote{A revised text of this article appeared in English translation in~\cite{tarski1983}, pp. 30--37, under the title `On some fundamental concepts of metamathematics'.}
	
	A. Tarski introduced a matrix method for studying ``deductive systems,'' that is, theories of abstract logics.\footnote{We read in~\cite{lukasiewicz-tarski1930}, footnote 2, ``The view of matrix formation as a general method of constructing systems is due to Tarski.''} In print form, the concept of a logical matrix was defined in~\cite{lukasiewicz-tarski1930}, definition 3.
	
	The theorem, which was originally known as the \emph{Lindenbaum theorem} first appeared in print without proof in~\cite{lukasiewicz-tarski1930}, theorem 3. This theorem reads that, given a consequence operator \textbf{Cn}, for any ``deductive system'' $X_0$, there exists ``a normal matrix,'' say \mat{M}, such that
	\begin{equation}\label{E:historic-lindenbaum-equivalence}
		\alpha\in\textbf{Cn}(X_0)~\Longleftrightarrow~X_0\models_{\textbf{M}}\alpha.
	\end{equation}
	Since $X_0$ in~\eqref{E:historic-lindenbaum-equivalence} was tacitly assumed to be fixed, defining
	an abstract logic $\mathcal{S}$ by the equivalence
	\[
	X\vdashS\alpha~\define~\alpha\in\Con{X_{0}\cup X},
	\]
	for any set $X$ and formula $\alpha$, we can reduce~\eqref{E:historic-lindenbaum-equivalence} to the second part of Proposition~\ref{P:lindenbaum-theorem}.
	
	The condition that $X_0$ is an $\mathcal{S}$-theory, that is $X_0=\ConS{\varnothing}$, is not important, but it is important to add that $X_0$ is closed under substitution, which must have been tacitly assumed. This observation implies two conclusions. First, participants in the seminar for mathematical logic at the University of Warsaw, where, according to~\cite{lukasiewicz-tarski1930}, footnote $\ddagger$, the result was reported, perhaps focused on the set of theorems of a deductive system, and not on the consequence relation associated with this system; and secondly, the set of theorems was assumed to be closed under substitution. 
	
	Further, in~\cite{lukasiewicz-tarski1930}, there was no indication of what kind of matrix \mat{M} is, in particular, whether $\mat{M}=\LinS$. 
	Moreover, although the nature of \mat{M} was very likely known (to Lindenbaum at least), the statement of theorem 3 did not contribute to the take-off of the Lindenbaum method.  The birth of the concept of a matrix consequence was also delayed until the publication of~\cite{los-suszko1958}; see definition (8.1) therein. The proof of Lindenbaum's theorem was published in~\cite{los1949}, proposition 10, and, independently, in~\cite{hermes1951}. (To prove~\eqref{E:historic-lindenbaum-equivalence}, {\L}o\'{s} used $\LinS$ as \mat{M}.) Further, since the authors of~\cite{lukasiewicz-tarski1930} did not require explicitly the closedness under substitution for an $\mathcal{S}$-theory $X$ (as {\L}o\'{s} did~\footnote{Cf.~\cite{los1949}, definition 3.}), they missed to notice that any Lindenbaum matrix $\LinS[X]$ is an $\mathcal{S}$-model, which was noted in~\cite{los1949}, proposition 7.\footnote{Actually, $\LinS[X]$ is an $\mathcal{S}$-model, even if $X$ is not closed under substitution, but $\mathcal{S}$ is structural; see Corollary~\ref{C:lindenbaum-matrix}.}
	
	The last omission and the omission of structurality in the 1930s however delayed the result, due to~\cite{wojcicki1969}, theorem 4, that any structural abstract logic is determined by its Lindenbaum matrices (Corollary~\ref{C:Lindenbaum-completeness}). After this, the role of substitution was clearly understood, since any valuation in the formula algebra is a substitution.
	
	The W\'{o}jcicki's observation opened up a way to define consequence relations  using bundles or atlases of logical matrices. Instead of the term `atlas' W\'{o}jcicki used the term `generalized matrix'. But even before him, T. Smiley~\cite{smiley1962} used a similar concept, known as a `Smiley matrix'\index{matrix!Smiley}. The term `atlas' was coined by M. J. Dunn and G. M. Hardegree in~\cite{dunn-hardegree2001}. The concept and term of the Lindenbaum bundle were used by G. Malinowski in~\cite{malinowski2007}.\\
	
	Almost certainly, introducing the concept of a logical consequence, A. Tarski had in mind the usual practice of obtaining new hypothetical truths in mathematics and other deductive sciences. So in his very first article~\cite{tarski1930a} on logical consequence, he wrote,
	\begin{quote}
		``From the sentences of any set $X$ cernain other sentences can be obtained by means of certain operations called \emph{rules of inference}. These sentences are called the \emph{consequences of the set} $X$. The set of all consequences is denoted by the symbol `\textit{Cn}($X$).'' Quoted from~\cite{tarski1983}, p. 30. 
	\end{quote}
	
	In antiquity, according to Plato, in the framework of his Theory of Forms, ``inference was [is] presumably valid when we follow in thought the connexions between Forms as they are.''\footnote{Cf.~\cite{kneales1962}, chapter I, section 5.} And in later antiquity, Galen and Boethius presented Aristotle's view on inference ``as a theory of inference schemata, i.e. as a theory in which rules of inference are thought by means of skeleton arguments.''\footnote{Cf.~\cite{kneales1962}, chapter II, section 6.} And later, in the development of logic, ``the Stoics generally concentrated attention on valid inference schemata (\textgreek{sunaktik\`{a} sq\'{h}mata}).''\footnote{Cf.~\cite{kneales1962}, chapter III, section 5.}
	The modern name for ``inference schemata'' is \emph{structural} (\emph{inference}) \emph{rules}. The last concept and term was introduced in~\cite{los-suszko1958}, section 3.
	
	However, any system of inference rules, structural (such as modus ponens) or non-structural (such as substitution), does not suffice to define a consequence relation. Perhaps G. Leibniz was the first logician who understood the importance of the notion of formal derivation, or formal proof, or formal deduction.\footnote{See~\cite{kneales1962}, chapter V, section 2, where we read: ``The most fruitful idea that Leibniz derived from his study of Aristotelian logic was the notion of formal proof.''} Also, Bolzano is recognized as ``a noteworthy precursor of modern proof-theory.''\footnote{Cf.~\cite{bochenski1970}, {\S} 38, section B1.}
	G. Boole used derivations in abstract manner to obtain consequences without referring to interpretation. G. Frege was close to formalizing the concept of deduction when he wrote:
	\begin{quote}
		``Inference is conducted in my symbolic system (\textit{Begriffsschrift}) according to a kind of calculation. I do not mean this in the narrow sense, as though an algorithm was in control, the same as or similar to that of ordinary addition and multiplication, but in the sense that the whole is algorithmic, with a complex of rules which so regulate the passage from one proposition or from two such to another, that nothing takes place but what is in accordance with these rules.'' Quoted from~\cite{bochenski1970}, {\S} 38, section B2.
	\end{quote}
	
	But only D. Hilbert and P. Bernays (1934) gave an exact definition of formal derivation, or formal proof, as a sequence of formulas satisfying certain conditions based on the rules of detachment (modus ponens) and substitution.~\footnote{Cf.~\cite{hilbert-bernays1968}, chapter III, {\S} 2.}
	Moreover, they formulated the first rule in terms of metavariables, that is, as a structural inference rule. The importance of such a definition is that it makes it possible to study formal derivations as words in a formal alphabet, using mathematical induction on the length of the derivation in question.\\

	The inference rules (a)--(b) were introduced by G. Gentzen; cf.~\cite{gentzen1964}.
	The hyperrules (c) were used by S. C. Kleene in~\cite{kleene1952}, {\S} 23, although as derived rules, not as postulated ones.~\footnote{Consult~\cite{scott1974b} about the difference between (postulated) rules and derived rules; the idea of a `derived rule of inference' is also explained on the example of uniform substitution in~\cite{church1996}, {\S} 19. It seems that the idea of derived inference rules was familiar to the Stoics; cf.~\cite{lukasiewicz1934}.} The inference rules (a)--(b) and the hyperrules (c) are similar to the proper rules and, respectively, to deduction rules of~\cite{prawitz1965}, chapter 1, {\S} 1.\\
	
	Proposition~\ref{P:brown-suszko-theorem} first appeared in~\cite{brown-suszko1973} as theorem 1 of section XII.

\chapter[Matrix Consequence]{Matrix Consequence}\label{chapter:matrix-consequence}	

\section{Single-matrix consequence}\label{section:single-matrix-consequence}
Two circumstances draw our attention to a single-matrix consequence. The first one is Lindenbaum's theorem (Proposition~\ref{P:lindenbaum-theorem}) which provides a universal way of definition of a weakly adequate logical matrix, namely the Lindenbaum matrix, which being defined for a given consequence relation, validates all of the theses of this consequence relation and only them. The other circumstance is related to the consequence relation $\vdash_2$ (Section~\ref{section:rules-and-hyperrules}) which turned out to be equal to the single-matrix relation $\models_{\textbf{B}_{2}}$. These circumstances lead to the the following concept.

We recall (Section~\ref{section:con-via-matrices}) that
a logical matrix $\mat{M}$ is called adequate for an abstract logic $\mathcal{S}$ if $\vdash_{\mathcal{S}}\,=\,\models_{\textbf{M}}$; that is, for any set $X$ of $\Lan$-formulas and any $\Lan$-formula $\alpha$,
\[
X\vdash_{\mathcal{S}}\alpha\Longleftrightarrow X\models_{\textbf{M}}\alpha.
\]

In connection with the notion of an adequate matrix, two questions arise. Is it true that any abstract logic has an adequate matrix? And if the answer to the first question is negative, what are the properties of an abstract logic that make it have an adequate matrix?

The first question can be answered in the negative as follows. Suppose $\mathcal{S}$ is a nontrivial abstract logic in $\Lan$ with $\card{\Var_{\Lan}}\ge\aleph_{0}$, which has a nonempty consistent set, say $X_0$, and such that $\bm{T}_{\mathcal{S}}\neq\varnothing$. Now, let
$\vdash_{\mathcal{S}}^{\circ}$ be a relation obtained from $\vdash_{\mathcal{S}}$ according to~\eqref{E:nonempty-consequence}. In virtue of Proposition~\ref{P:nonempty-consequence}, $\vdash_{\mathcal{S}}^{\circ}$ is a consequence relation. We claim that there is no adequate matrix for $\vdash_{\mathcal{S}}^{\circ}$. (Such a relation $\vdash_{\mathcal{S}}^{\circ}$ can be exemplified, e.g., by $\models_{\textbf{B}_2}^{\circ}$ with $X_0=L\booleTwo$, where $\models_{\textbf{B}_2}$ is a single-matrix consequence of Section~\ref{section:con-via-matrices}.)

Indeed, for contradiction, assume that a matrix $\mat{M}=\langle\alg{A}, D\rangle$ is adequate for $\vdash_{\mathcal{S}}^{\circ}$. Since $X_0$ is consistent with respect to $\mathcal{S}$, it is also consistent with respect to $\vdash_{\mathcal{S}}^{\circ}$.
This and the assumption that $X_0$ is nonempty imply that $D\neq\varnothing$.
Assume that $a\in D$.
By premise, there is a thesis $\alpha\in\bm{T}_{\mathcal{S}}$. Let us take a variable $p\in\Var_{\Lan}\setminus\Var(\alpha)$. By definition, $\varnothing\not\vdash_{\mathcal{S}}^{\circ}\alpha$ and $p\vdash_{\mathcal{S}}^{\circ}\alpha$. The former implies that there is a valuation $v$ in $\alg{A}$ such that $v[\alpha]\notin D$. Now we define a valuation:
\[
v^{\prime}[q]:=\begin{cases}
	\begin{array}{cl}
		v[q] &\text{if $q\neq p$}\\
		a &\text{if $q=p$}.
	\end{array}
\end{cases}
\]
It is clear that $v^{\prime}[\alpha]=v[\alpha]$ and hence $p\not\models_{\textbf{M}}\alpha$.\\

Trying to find a characteristic for abstract logic to have an adequate matrix, we consider any logical matrix $\mat{M}:=\langle\alg{A},D\rangle$ of type $\Lan$.

Now, let $\alpha, \beta$ and $\gamma$ be $\Lan$-formulas with $(\Var(\alpha)\cup\Var(\beta))\cap\Var(\gamma)=\varnothing$. Assume that for some valuation $v_0$, $v_0(\gamma)\in D$. (This assumption is equivalent to the condition that the set $\lbrace\gamma\rbrace$ is consistent with respect to the consequence relation $\models_{\textbf{M}}$.) We prove that
\[
\gamma,\alpha\models_{\textbf{M}}\beta~~\textit{implies}~~\alpha\models_{\textbf{M}}\beta.
\]

Indeed, for contrapositive, suppose that $\alpha\not\models_{\textbf{M}}\beta$;
that is, for some valuation $v$, $v(\alpha)\in D$, but $v(\beta)\notin D$.

Let us define a valuation
\[
v_1(p)=\begin{cases}
	\begin{array}{cl}
		v_0[p] &\text{if $p\in\Var(\gamma)$}\\
		v[p] &\text{otherwise}.
	\end{array}
\end{cases}
\]

We note that $v_1[\gamma]=v_0[\gamma]$, $v_1[\alpha]=v[\alpha]$ and  $v_1[\beta]=v[\beta]$. This implies that $\gamma,\alpha\not\models_{\textbf{M}}\beta$.\\

The above consideration leads to the following definition.
\begin{defn}[uniform abstract logic]\label{D:uniform-logic}\index{logic!uniform abstract}
	An abstract logic $\mathcal{S}$ is said to be \textbf{uniform} if 	for any set
	$X\cup Y\cup\lbrace\alpha\rbrace$ of $\Lan$-formulas with $\Var(X\cup\lbrace\alpha\rbrace)\cap\Var(Y)=\varnothing$ and {\em$\textbf{Cn}_{\mathcal{S}}(Y)\neq\Forms_{\mathcal{L}}$}, if $X,Y\vdash_{\mathcal{S}}\alpha$, then $X\vdash_{\mathcal{S}}\alpha$.
\end{defn}

We note that the consequence relation $\models_{\textbf{B}_{2}}^{\circ}$ mentioned above is finitary but not uniform. (Exercise~\ref{section:single-matrix-consequence}.\ref{EX:b-two-circle})
The importance of the last definition becomes clear from Proposition~\ref{P:los-suszko} below. We prove this proposition by means of two lemmas.

\begin{lem}\label{L:los-suszko-1}
	If an abstract logic $\mathcal{S}$	has an adequate matrix, then $\mathcal{S}$ is uniform.
\end{lem}
\begin{proof}
	Let a matrix {\em$\mat{M}=\langle\alg{A},D\rangle$} be adequate for $\mathcal{S}$.
	Suppose $X,Y\vdash_{\mathcal{S}}\alpha$, where $\Var(X,\alpha)\cap\Var(Y)=\varnothing$
	and $\textbf{Cn}_{\mathcal{S}}(Y)\neq\Forms_{\mathcal{L}}$. We will show that $X\models_{\textbf{M}}\alpha$, which, in virtue of the adequacy of $\mat{M}$ for $\mathcal{S}$, will in turn imply that $X\vdash_{\mathcal{S}}\alpha$.
	
	For contradiction, assume that for some valuation $v$ in $\alg{A}$, $v[X]\subseteq D$ but $v[\alpha]\notin D$. Since $Y$ is consistent with respect to $\mathcal{S}$, there is a valuation $v^{\prime}$ such that $v^{\prime}[Y]\subseteq D$. Then, we define a valuation $v^{\prime\prime}$ as follows:
	\[
	v^{\prime\prime}[p]:=\begin{cases}
		\begin{array}{cl}
			v[p] &\text{if $p\in\Var(X,\alpha)$}\\
			v^{\prime}[p] &\text{if $p\in\Var(Y)$}.
		\end{array}
	\end{cases}
	\]
	
	It must be clear that $v^{\prime\prime}[\beta]=v[\beta]$, whenever $\beta\in X\cup\lbrace\alpha\rbrace$, and $v^{\prime\prime}[\beta]=v^{\prime}[\beta]$, if $\beta\in Y$. Therefore, $X,Y\not\models_{\textbf{M}}\alpha$; that is $X,Y\not\vdash_{\mathcal{S}}\alpha$. A contradiction.
\end{proof}

As preparation to the next proposition we observe the following.
\begin{lem}\label{L:(S)_plus-uniformity}
	Let $\mathcal{S}$ be a uniform structural abstract logic in a language $\Lan$ with $\card{\Var_{\Lan}}\ge\aleph_{0}$ and $\Lan^{+}$ be a primitive extension of $\Lan$. Then the abstract logic $(\mathcal{S})^{+}$ is also uniform.
\end{lem}
\begin{proof}
	Suppose for a set $X\cup Y\cup\lbrace\alpha\rbrace\subseteq\Forms_{\mathcal{L}^{+}}$,
	$X,Y\vdash_{(\mathcal{S})^{+}}\alpha$, providing that $\Var(X\cup\lbrace\alpha
	\rbrace)\cap\Var(Y)=\varnothing$ and $Y$ is consistent with respect to $(\mathcal{S})^{+}$.
	Then, by definition, there is $X_0\cup Y_0\cup\lbrace\beta\rbrace\subseteq\Forms_{\mathcal{L}}$ and an $\Lan^{+}$-substitution $\sigma$ such that $\sigma(X_0)\subseteq X$, $\sigma(Y_0)\subseteq Y$, $\sigma(\beta)=\alpha$ and $X_0,Y_0\vdash_{\mathcal{S}}\beta$. It is clear that
	$\Var(X_0\cup\lbrace\beta\rbrace)\cap\Var(Y_0)=\varnothing$. Also, $Y_0$ is consistent with respect to $\mathcal{S}$, for if it were otherwise, in virtue of conservativity of $(\mathcal{S})^{+}$ over $\mathcal{S}$, we would have $\sigma(Y_0)\vdash_{(\mathcal{S})^{+}}\gamma$, for any arbitrary $\Lan^{+}$-formula $\gamma$. This would imply the inconsistency of $Y$ in $(\mathcal{S})^{+}$.
	
	Since $\mathcal{S}$ is uniform, $X_0\vdash_{\mathcal{S}}\beta$ and hence $X_0\vdash_{(\mathcal{S})^{+}}\beta$, for $(\mathcal{S})^{+}$ is a conservative extension of $\mathcal{S}$. And since $(\mathcal{S})^{+}$ is structural (Proposition~\ref{P:(S)^+}), $\sigma(X_0)\vdash_{(\mathcal{S})^{+}}\alpha$ and hence
	$X\vdash_{(\mathcal{S})^{+}}\alpha$.
\end{proof}

\begin{prop}[{\L}o\'{s}-Suszko-W\'{o}jcicki theorem]\label{P:los-suszko}\index{Theorem!{\L}o\'{s}-Suszko-W\'{o}jcicki }
	Let $\mathcal{S}$ be a finitary structural abstract logic in a language $\Lan$ with $\card{\Var_{\Lan}}\ge\aleph_{0}$. Then $\mathcal{S}$ has an adequate matrix if, and only if, $\mathcal{S}$ is uniform.
\end{prop}
\begin{proof}
	The `only-if' part follows from Lemma~\ref{L:los-suszko-1}. So we continue with the proof of the `if' part.
	
	We distinguish two cases: $\mathcal{S}$ is trivial and $\mathcal{S}$ is nontrivial. If the former is the case, then the matrix $\langle\FormAl,\FormsL\rangle$ is adequate for $\mathcal{S}$.
	
	Now, we assume that $\mathcal{S}$ is nontrivial; that is, the set $\theoryC$ of consistent sets with respect to $\mathcal{S}$ is nonempty. 
	
	We introduce for each $X\in\theoryC$ and each $p\in\Var_{\Lan}$, a new sentential variable $p_X$ so that for any $p,q\in\Var_{\Lan}$ and any $X,Y\in\theoryC$,
	\[
	p_X=q_Y\Longleftrightarrow p=q~\text{and}~X=Y.
	\]
	Accordingly, for an arbitrary $X\in\theoryC$, we denote:
	\[
	\Var_X=\set{p_X}{p\in\Var_{\Lan}}.
	\]
	Thus for any $X,Y\in\theoryC$,
	\begin{equation}\label{E:implication-1}
		X\neq Y\Longrightarrow\Var_X\cap\Var_Y=\varnothing.
	\end{equation}
	
	A language $\Lan^{+}$ is defined as a primitive extension of $\Lan$, by adding the new sentential variables to $\Var_{\Lan}$. Then, we introduce $(\mathcal{S})^{+}$ according to Definition~\ref{D:(S)^+}.
	
	Next, for each $X\in\theoryC$, we define two $\Lan^{+}$-substitutions:
	\[
	\begin{array}{l}
		\sigma_X:p\mapsto p_{X},~\text{for any $p\in\Var_{\Lan}$},~\text{and}~\sigma_X:p_Y\mapsto p_Y,
		~\text{for any $p_Y\in\Var_{\Lan^{+}}\setminus\Var_{\Lan}$};\\
		\xi_X: p_X\mapsto p,~\text{for any $p_X\in\Var_X$ and $\xi_X:q\mapsto q$, for any $q\in\Var_{\Lan^{+}}\setminus\Var_X$}.
	\end{array}
	\]
	
	We observe: for any $X\in\theoryC$,
	\begin{equation}\label{E:zig-zag}
		\xi_{X}\circ\sigma_X(\beta)=\beta,~\text{for all $\beta\in\Forms_{\mathcal{L}}$}.
	\end{equation}
	
	Then, we make two observations.
	First, since for any $X\in\theoryC$, $\Var(\sigma_X(X))\subseteq\Var_X$, \eqref{E:implication-1} yields that for any $X,Y\in\theoryC$,
	\begin{equation}\label{E:implication-2}
		X\neq Y\Longrightarrow\Var(\sigma_X(X))\cap\Var(\sigma_Y(Y))=\varnothing.
	\end{equation}
	
	Second, 
	\begin{equation}\label{E:second-observation}
		\textit{for each $X\in\theoryC$, the set $\sigma_X(X)$ is consistent with respect to $(\mathcal{S})^{+}$.}
	\end{equation}
	
	Indeed, for contradiction, assume that $\sigma_{X}(X)$ is inconsistent. Let $\beta$ be any formula of $\Forms_{\mathcal{L}}$. By assumption, $\sigma_{X}(X)\vdash_{(\mathcal{S})^{+}}\sigma_{X}(\beta)$. Since $(\mathcal{S})^{+}$ is structural (Proposition~\ref{P:(S)^+}), we have that $\xi_{X}(\sigma_{X}(X))\vdash_{(\mathcal{S})^{+}}\xi_{X}(\sigma_{X}(\beta))$ and hence, in virtue of \eqref{E:zig-zag}, $X\vdash_{(\mathcal{S})^{+}}\beta$. Since $(\mathcal{S})^{+}$ is conservative over $\mathcal{S}$, we obtain that $X\vdash_{\mathcal{S}}\beta$. A contradiction.
	
	Now, we define:
	\[
	\mat{M}:=\langle\mathfrak{F}_{\mathcal{L}^{+}},\textbf{Cn}_{(\mathcal{S})^{+}}(\bigcup
	\set{\sigma_X(X)}{X\in\theoryC})\rangle.
	\]
	
	We aim to show that $\mat{M}$ is adequate for $\mathcal{S}$; that is, for any set $Y\cup\lbrace\alpha\rbrace\subseteq\Forms_{\mathcal{L}}$,
	\[
	Y\vdash_{\mathcal{S}}\alpha\Longleftrightarrow Y\models_{\textbf{M}}\alpha.
	\]
	
	Suppose $Y\vdash_{\mathcal{S}}\alpha$. Since $(\mathcal{S})^{+}$ is conservative over $\mathcal{S}$, this implies that $Y\vdash_{(\mathcal{S})^{+}}\alpha$. Now, let $\sigma$ be any $\Lan^{+}$-substitution.
	Assume that $\sigma(Y)\subseteq\textbf{Cn}_{(\mathcal{S})^{+}}(\bigcup
	\set{\sigma_X(X)}{X\in\theoryC})$. Since $(\mathcal{S})^{+}$ is structural, we have: $\sigma(\alpha)\in\textbf{Cn}_{(\mathcal{S})^{+}}(\sigma(Y))$ and hence
	$\sigma(\alpha)\in\textbf{Cn}_{(\mathcal{S})^{+}}(\bigcup
	\set{\sigma_X(X)}{X\in\theoryC})$.
	
	Conversely, suppose that $Y\not\vdash_{\mathcal{S}}\alpha$. We denote:
	\[
	X_0:=\textbf{Cn}_{\mathcal{S}}(Y).
	\]
	Thus $\alpha\notin X_0$ and hence $X_0\in\theoryC$. The latter implies that
	$\sigma_{X_0}(X_0)\subseteq\textbf{Cn}_{(\mathcal{S})^{+}}(\bigcup
	\set{\sigma_X(X)}{X\in\theoryC})$.
	
	For contradiction, we assume that
	\[
	\sigma_{X_0}(\alpha)\in\textbf{Cn}_{(\mathcal{S})^{+}}(\bigcup
	\set{\sigma_X(X)}{X\in\theoryC}). \tag{$\ast$}
	\]
	In virtue of finitariness of $(\mathcal{S})^{+}$ (Proposition~\ref{P:(S)^{+}-existence}), there is a finite number of sets $X_1,\ldots, X_n\in\theoryC$ such that 
	\[
	\sigma_{X_0}(\alpha)\in\textbf{Cn}_{(\mathcal{S})^{+}}(\sigma_{X_0}(X_0)\cup\sigma_{X_1}
	(X_1)\cup\ldots\cup\sigma_{X_n}(X_n)).
	\]
	Since, according to~\eqref{E:second-observation}, $\sigma_{X_n}(X_n)$ is consistent and
	\[
	\Var(\sigma_{X_0}(X_0)\cup\sigma_{X_1}
	(X_1)\cup\ldots\cup\sigma_{X_{n-1}}(X_{n-1}))\cap\Var(\sigma_{X_{n}}(X_{n}))=\varnothing,
	\]
	in virtue of the uniformity of $(\mathcal{S})^{+}$ (Lemma~\ref{L:(S)_plus-uniformity}),
	\[
	\sigma_{X_0}(X_0),\sigma_{X_1}
	(X_1),\ldots\sigma_{X_{n-1}}(X_{n-1})\vdash_{(\mathcal{S})^{+}}\sigma_{X_0}(\alpha)
	\]
	Repeating this argument several times, we arrive at the conclusion that 
	$\sigma_{X_0}(X_0)\vdash_{(\mathcal{S})^{+}}\sigma_{X_0}(\alpha)$. Then, applying the substitution $\xi_{X_0}$, we get $X_0\vdash_{(\mathcal{S})^{+}}\alpha$. Finally, the conservativity of $(\mathcal{S})^{+}$ over $\mathcal{S}$ yields $X_0\vdash_{\mathcal{S}}\alpha$. A contradiction.
\end{proof}

The case of nonfinitary abstract logics needs an additional investigation.
\begin{defn}[couniform abstract logic]\label{D:couniform-logic}
	An abstract logic $\mathcal{S}$ in a language $\Lan$ is called couniform if for any collection $\lbrace X_i\rbrace_{i\in I}$ of formula sets with $\Var(X_i)\cap\Var(X_j)=\varnothing$, providing that $i\neq j$, and $\Var(\bigcup_{i\in I}\lbrace X_i\rbrace) \neq\Var_{\Lan}$, there is at least one $X_i$ which is inconsistent, whenever the union $\bigcup_{i\in I}\lbrace X_i\rbrace$ is inconsistent.	
\end{defn}

The next lemma shows that couniformity is necessary for an abstract logic to have an adequate matrix.
\begin{lem}\label{L:couniform}
	If an abstract logic $\mathcal{S}$ has an adequate matrix, then $\mathcal{S}$ is couniform.
\end{lem}
\begin{proof}
	Let a matrix $\mat{M}=\langle\alg{A},D\rangle$ be adequate for $\mathcal{S}$. We assume that a family $\lbrace X_i\rbrace_{i\in I}$ of formula sets satisfies all the premises of couniformity. Suppose each $X_i$ is consistent with respect to $\mathcal{S}$. Then for each $i\in I$, there is a formula $\alpha_i$ such that $X_i\not\vdash_{\mathcal{S}}\alpha_i$. This implies that there is a valuation $v_i$ in $\alg{A}$ such that $v_i[X_i]\subseteq D$ but
	$v_i[\alpha_i]\notin D$. Let us select either one of $v_i[\alpha_i]$ and denote
	\[
	a:=v_i[\alpha_i].
	\]
	Then, we select any variable $p\in\Var_{\Lan}\setminus\Var(\bigcup_{i\in I}\lbrace X_i\rbrace)$ and for each $v_i$, define:
	\[
	v_{i}^{\ast}[q]:=\begin{cases}
		\begin{array}{cl}
			v_i[q] &\text{if $q\in\Var(\bigcup_{i\in I}\lbrace X_i\rbrace)$}\\
			a &\text{if $q\in\Var_{\Lan}\setminus\Var(\bigcup_{i\in I}\lbrace X_i\rbrace)$}.
		\end{array}
	\end{cases}
	\]
	
	It must be clear that for each $i\in I$, $v_{i}^{\ast}[X_i]\subseteq D$. And since $\Var(X_i)\cap\Var(X_j)=\varnothing$, providing $i\neq j$, we also have that $v_{i}^{\ast}[\bigcup_{i\in I}\lbrace X_i\rbrace]
	\subseteq D$. Thus $\bigcup_{i\in I}\lbrace X_i\rbrace\not\vdash_{\mathcal{S}} p$, that is the set $\bigcup_{i\in I}\lbrace X_i\rbrace$ is consistent.
\end{proof}

The next lemma is rather technical and we will need it for the proof of Proposition~\ref{P:wojcicki} below.
\begin{lem}\label{L:(S)^{+}-couniform}
	Let $\mathcal{S}$ be a couniform structural abstract logic in a language $\Lan$ with $\card{\Var_{\Lan}}\ge\aleph_{0}$ and $\Lan^{+}$ be a primitive extension of $\Lan$. Then the abstract logic $(\mathcal{S})^{+}$ is also couniform.
\end{lem}
\begin{proof}
	Let $\lbrace X_i\rbrace_{i\in I}$ be a nonempty family of nonempty sets of $\Lan^{+}$-formulas, about which we assume that the union $\bigcup_{i\in I}\lbrace X_i\rbrace$ is inconsistent and $\Var(\bigcup_{i\in I}\lbrace X_i\rbrace)\neq\Var_{\Lan^{+}}$. Let us select any $p\in\Var_{\Lan^{+}}\setminus\Var(\bigcup_{i\in I}\lbrace X_i\rbrace)$. Since $\bigcup_{i\in I}\lbrace X_i\rbrace$ is inconsistent, we have that $\bigcup_{i\in I}\lbrace X_i\rbrace\vdash_{(\mathcal{S})^{+}}p$. By definition of $(\mathcal{S})^{+}$ (Definition~\ref{D:(S)^+}), there are a set $Y$ of $\Lan$-formulas, an $\Lan$-variable $q$ and an $\Lan$-substitution $\sigma$ such that
	$\sigma(Y)\subseteq\bigcup_{i\in I}\lbrace X_i\rbrace$, $\sigma(q)=p$ and $Y\vdash_{\mathcal{S}}q$. 
	
	Next we define: for each $i\in I$,
	\[
	Y_i:=\set{\gamma\in Y}{\sigma(\gamma)\in X_i}.
	\]
	
	Firstly, we notice that for each $i\in I$, $\sigma(Y_i)\subseteq X_i$, and $\bigcup_{i\in I}\lbrace Y_i\rbrace\vdash_{\mathcal{S}}q$.
	(The latter is true because $\bigcup_{i\in I}\lbrace Y_i\rbrace=Y$.)
	Secondly, $q\in\Var_{\Lan}\setminus\Var(\bigcup_{i\in I}\lbrace Y_i\rbrace)$. Indeed, if it were the case that $q\in\Var(\bigcup_{i\in I}\lbrace Y_i\rbrace)$, then we would have that $p\in
	\Var(\bigcup_{i\in I}\lbrace X_i\rbrace)$.
	
	Thus we obtain that $\bigcup_{i\in I}\lbrace Y_i\rbrace$ is inconsistent with respect to $\mathcal{S}$. Since, by definition, $\mathcal{S}$ is couniform, there is $i_0\in I$ such that $Y_{i_0}$ is inconsistent with respect to $\mathcal{S}$. This implies that  Since $(\mathcal{S})^{+}$ is a conservative extension of $\mathcal{S}$, we receive that $Y_{i_0}\vdash_{(\mathcal{S})^{+}}q$. And since 
	$(\mathcal{S})^{+}$ is structural (Proposition~\ref{P:(S)^+}), $\sigma(Y_{i_0})\vdash_{(\mathcal{S})^{+}}p$ and hence $X_{i_0}\vdash_{(\mathcal{S})^{+}}p$. We recall that $p\notin\Var(X_{i_0})$. Therefore, the structurality of $(\mathcal{S})^{+}$ implies that $X_{i_0}$ is inconsistent.
\end{proof}

\begin{prop}[W\'{o}jcicki theorem]\label{P:wojcicki}\index{Theorem!W\'{o}jcicki}
	Let $\mathcal{S}$ be a structural abstract logic in a language $\Lan$ with $\card{\Var_{\Lan}}\ge\aleph_{0}$. Then $\mathcal{S}$ has an adequate matrix if, and only if, $\mathcal{S}$ is both uniform and couniform.
\end{prop}
\begin{proof}
	The `only-if' part follows from Lemmas~\ref{L:los-suszko-1} and~\ref{L:couniform}. 
	
	The proof of the `if' part repeats the steps of the ``if'' part of the proof of the Lo\'{s}-Suszko theorem (Proposition~\ref{P:los-suszko}) up to the assumption $(\ast)$.
	Then, we continue as follows.
	
	First, we notice that, in virtue of~\eqref{E:second-observation} and of the couniformity 
	of $(\mathcal{S})^{+}$ (Lemma~\ref{L:(S)^{+}-couniform}), the set $\bigcup\set{\sigma_{X}(X)}{X\in\theoryC}$ is consistent. Thus, we have:
	\[
	\sigma_{X_0}(X_0),\bigcup\set{\sigma_{X}(X)}{X\in\theoryC~\text{and}~X\neq X_0}
	\vdash_{(\mathcal{S})^{+}}\sigma_{X_0}(\alpha),
	\]
	where the set $\bigcup\set{\sigma_{X}(X)}{X\in\theoryC~\text{and}~X\neq X_0}$ is consistent with respect to $(\mathcal{S})^{+}$.
	
	Let us denote
	\[
	Y_0:=\bigcup\set{\sigma_{X}(X)}{X\in\theoryC~\text{and}~X\neq X_0}.
	\]
	With the help of~\eqref{E:implication-1} and~\eqref{E:implication-2}, we conclude that $\Var(\sigma_{X_0}(X_0)\cup\lbrace\sigma_{X_0}(\alpha)\rbrace)\cal\Var(Y_0)=\varnothing$.
	Then, using the uniformity of $(\mathcal{S})^{+}$, we derive that $\sigma_{X_0}(\alpha)
	\vdash_{(\mathcal{S})^{+}}\sigma_{X_0}(\alpha)$. This, in virtue of structurality of $(\mathcal{S})^{+}$ (Proposition~\ref{P:(S)^+}), implies that
	$\xi_{X_0}(\sigma_{X_0}(\alpha))
	\vdash_{(\mathcal{S})^{+}}\xi_{X_0}(\sigma_{X_0}(\alpha))$ and hence, according to~\eqref{E:zig-zag}, $X_0\vdash_{(\mathcal{S})^{+}}\alpha$. Because of the conservativeness of $(\mathcal{S})^{+}$ over $\mathcal{S}$, the latter leads to $X_0\vdash_{\mathcal{S}}\alpha$. A contradiction.
\end{proof}

\begin{defn}\label{D:cancellation-property}\index{logic!cancellation property}
	An abstract logic $\mathcal{S}$ has the \textbf{cancellation property} if for any system $\lbrace X_i\rbrace_{i\in I}$ of consistent sets with $\Var(X_i)\cap \Var(X_j)=\varnothing$, whenever $i\neq j$, and a set $X\cup\lbrace\alpha\rbrace$ of $\Lan$-formulas such that for all $i\in I$,   $\Var(X_i)\cap\Var(X\cup\lbrace\alpha\rbrace)=\varnothing$, if
	$X\cup\bigcup_{i\in I}\lbrace X_i\rbrace\vdashS\alpha$, then $X\vdashS\alpha$.
\end{defn}

\begin{lem}\label{L:shoesmith-smiley}
	Let an abstract logic $\mathcal{S}$ be structural. Then $\mathcal{S}$ has the cancellation property if, and only if, $\mathcal{S}$ is uniform and couniform.
\end{lem}
\begin{proof}
	Assume first that $\mathcal{S}$ has the cancellation property. Then, obviously, $\mathcal{S}$ is also uniform. 
	
	Next, let a system $\lbrace X_i\rbrace_{i\in I}$ of consistent sets satisfy all the premises of couniformity. If $\card{I}=1$, that is $\lbrace X_i\rbrace_{i\in I}=\lbrace X\rbrace$, then, by premise $\ConS{\bigcup_{i\in I}\lbrace X_i\rbrace}=\ConS{X}\neq\FormsL$. Otherwise, let us fix any $X_{i_0}$, where $i_0\in I$. 
	
	Now, let us take $p\in\VarL\setminus\bigcup_{i\in I}\lbrace X_i\rbrace$.
	By premise, $X_{i_0}$ is consistent. For contradiction, assume that $\bigcup_{i\in I}\lbrace X_i\rbrace$ is inconsistent, that is, $\ConS{X_{i_0}\cup\bigcup_{j\in I\setminus\lbrace i_0\rbrace}\lbrace X_j\rbrace}=\FormsL$. By the cancellation property,
	$X_{i_0}\vdashS p$. Since $\mathcal{S}$ is structural, we derive that $X_{i_0}$ is inconsistent. A contradiction.
	
	Conversely, suppose $\mathcal{S}$ is uniform and couniform. Let $\lbrace X_i\rbrace_{i\in I}$ and $X\cup\lbrace\alpha\rbrace$ be a system of consistent sets and, respectively, a set of formulas satisfying all the premises of the cancellation property. 
	
	We denote
	\[
	Y:=\bigcup_{i\in I}\lbrace X_i\rbrace.
	\]
	
	Since $\mathcal{S}$ is couniform, $Y$ is consistent. Now assume that $X,Y\vdashS\alpha$. Since $\mathcal{S}$ is uniform, $X\vdashS\alpha$.
\end{proof}

\begin{prop}[Shoesmith-Smiley theorem]\label{P:shoesmith-smiley-theorem}\index{Shoesmith-Smiley}
	Let $\mathcal{S}$ be a structural abstract logic in a language $\Lan$ with $\card{\VarL}\ge\aleph_{0}$. Then $\mathcal{S}$ has an adequate matrix if, and only if, $\mathcal{S}$ has the cancellation property.	
\end{prop}
\begin{proof}
	We apply Lemma~\ref{L:shoesmith-smiley} to Proposition~\ref{P:wojcicki}.
\end{proof}

The following two criteria for single-matrix consequence are formulated in terms of matrix consequence.

Suppose we have an $\Lan$-matrix $\mat{M}=\langle\alg{A},D\rangle$. We define:
\[
X\in\mathfrak{S}(\mat{M})~\define~\text{there is a valuation $v$ in \alg{A} such that $X=\set{\alpha}{v[\alpha]\in D}$}.
\]

Then, we define:
\begin{equation}\label{E:R(M)-definition}
	X\in\mathfrak{R}(\mat{M})~\define~\text{there is $Y\in\mathfrak{S}(\mat{M})$ such that $X\subseteq Y$}.
\end{equation}

According to~\eqref{E:R(M)-definition}, we have:
\begin{equation*}\label{E:S(M)-included-R(M)}
	\mathfrak{S}(\mat{M})\subseteq\mathfrak{R}(\mat{M}).
\end{equation*}

Next, given $\Lan$-matrices $\mat{M}$ and $\mat{N}$, we observe:
\begin{equation*}
	\mathfrak{S}(\mat{M})\subseteq\mathfrak{S}(\mat{N})~\Longrightarrow~\mathfrak{R}(\mat{M})\subseteq\mathfrak{R}(\mat{N}).
\end{equation*}

Indeed, let $X\in\mathfrak{R}(\mat{M})$. In virtue of~\eqref{E:R(M)-definition}, $X\subseteq Y$, for some $Y\in\mathfrak{S}(\mat{M})$. By premise, $Y\in\mathfrak{S}(\mat{N})$ and hence, by~\eqref{E:R(M)-definition},
$X\in\mathfrak{R}(\mat{N})$.

The last implication obviously implies that
\begin{equation}\label{E:R(M)-and-S(M)-2}
	\mathfrak{S}(\mat{M})=\mathfrak{S}(\mat{N})~\Longrightarrow~\mathfrak{R}(\mat{M})=\mathfrak{R}(\mat{N}).
\end{equation}

First, we use the $\mathfrak{R}$-operator to characterize single-matrix consequence.

\begin{prop}\label{P:R-operator-criterion}
	A structural consequence relation in a language $\Lan$ is a single-matrix consequence if, and only if, it can be determined by a class {\em$\mathcal{M}=\lbrace\mat{M}_i\rbrace_{i\in I}$} of $\Lan$-matrices such that for any $i,j\in I$, {\em$\mathfrak{R}(\mat{M}_i)=\mathfrak{R}(\mat{M}_j)$}.
\end{prop}
\begin{proof}
	The `only-if' implication is obvious.
	
	So, assume that the given consequence relation is determined by a class $\mathcal{M}=\lbrace\mat{M}_i\rbrace_{i\in I}$ with $\mathfrak{R}(\mat{M}_i)=\mathfrak{R}(\mat{M}_j)$, for any $i,j\in I$.
	
	Using notation $\mat{M}_{i}=\langle\alg{A}_i,D_i\rangle$, we define the matrix
	\[
	\mat{M}^{\ast}:=\langle\prod_{i\in I}\alg{A}_i,\prod_{i\in I} D_i\rangle,
	\]
	where the first component is the direct product of the algebras $\alg{A}_i$ and the second is the Cartesian product of the sets $D_i$. 
	
	We intend to prove that the $\mathcal{M}$-consequence equals $\mat{M}^{\ast}$-consequence.
	
	Suppose $X\vdash_{\mathcal{M}}\alpha$ and let $v$ be a valuation with $v[X]\subseteq\prod_{i\in I} D_i$. We denote by $v_i$ the projection of $v$ on $\alg{A}_i$. Since  each $v_i$ is a valuation, respectively, in $\alg{A}_i$ with
	$v_{i}[X]\subseteq D_i$, then, by premise, $v_{i}[\alpha]\in D_i$, for any $i\in I$. Hence $v[\alpha]\in\prod_{i\in I} D_i$. 
	
	Next, assume that $X\vdash_{\mat{M}^{\ast}}\alpha$. Let $v_{i_0}$ be a valuation in $\mat{M}_{i_0}$ with $v_{i_0}[X]\subseteq D_{i_0}$, for arbitrary $i_0\in I$. This means that $X\in\mathfrak{R}(\mat{M}_{i_0})$. By premise, for all other $j\in I$, $X\in\mathfrak{R}(\mat{M}_j)$ with respect to corresponding valuations $v_j$ in $\mat{M}_j$. This allows us to define the valuation $v$ in $\mat{M}^\ast$ so that for each $i\in I$, the projection of $v$ on $\alg{A}_i$ equals $v_i$.
	Then, by premise, we have that $v[\alpha]\in\prod_{i\in I} D_i$ which in turn implies that $v_i[\alpha]\in D_i$, for each $i\in D_i$, in particular, $v_{i_0}[\alpha]\in D_{i_0}$.
\end{proof}

An advantage of the $\mathfrak{S}$-operator over the $\mathfrak{R}$-operator is that in terms of the former, matrix consequence can be easily formulated.

Indeed, let $\mathcal{M}:=\lbrace\mat{M}_i\rbrace_{i\in I}$ be a set of $\Lan$-matrices. It is obvious that
\begin{equation}\label{E:equivalence-matrix-con}
	X\vdash_{\mathcal{M}}\alpha~\Longleftrightarrow~\text{for every $i\in I$, for every $Y\in\mathfrak{S}(\mat{M}_i)$, if $X\subseteq Y$, then $\alpha\in Y$}.
\end{equation}
(Exercise~\ref{section:single-matrix-consequence}.\ref{EX:equivalence-matrix-con})
\begin{prop}\label{P:S-operator-criterion}
	A structural consequence relation in a language $\Lan$ is a single-matrix consequence if, and only if, it can be determined by a class {\em$\mathcal{M}=\lbrace\mat{M}_i\rbrace_{i\in I}$} of $\Lan$-matrices such that for any $i,j\in I$, {\em$\mathfrak{S}(\mat{M}_i)=\mathfrak{S}(\mat{M}_j)$}.
\end{prop}
\begin{proof}
	Indeed, the `only-if' part is obvious.
	
	To prove the `if' part, we notice that, in virtue of~\eqref{E:R(M)-and-S(M)-2},
	for any $i,j\in I$, $\mathfrak{R}(\mat{M}_i)=\mathfrak{R}(\mat{M}_j)$. It remains to apply Proposition~\ref{P:R-operator-criterion}.
\end{proof}

We note that the restriction on the language `$\card{\VarL}\ge\aleph_{0}$' imposed in Proposition~\ref{P:los-suszko}, Proposition~\ref{P:wojcicki} and Proposition~\ref{P:shoesmith-smiley-theorem} in Proposition~\ref{P:R-operator-criterion} and Proposition~\ref{P:S-operator-criterion} can be dropped.

\paragraph{Exercises~\ref{section:single-matrix-consequence}}
\begin{enumerate}
	\item\label{EX:b-two-circle} Prove that the relation $\models_{\textbf{B}_{2}}^{\circ}$ obtained from $\models_{\textbf{B}_2}$ according to \eqref{E:nonempty-consequence} with $X_0=L\booleTwo$ is finitary but not uniform.
	\item Prove that a consequence relation defined by premiseless rules and modus ponens is uniform.
	\item\label{EX:equivalence-matrix-con}Prove the equivalence~\eqref{E:equivalence-matrix-con}.
\end{enumerate}

\section{Finitary matrix consequence}\label{section:finitary-matrix-consequence}
Given a (nonempty) class $\mathcal{M}$ of $\Lan$-matrices, an abstract logic $\mathcal{S}_{\mathcal{M}}$, according to Proposition~\ref{P:matrix-con-is-con-relation},
is structural and, according to Corollary~\ref{C:S-matrix-completeness}, is determined by the class of all $\mathcal{S}_{\mathcal{M}}$-models. The last class definitely contains $\mathcal{M}$. For the purpose of this section, it will not be an exaggeration to assume that $\mathcal{M}$ coincides with the class of all $\mathcal{S}_{\mathcal{M}}$-models.
Under this assumption, the question arises: What characteristics should $\mathcal{M}$ have in order for $\mathcal{S}_{\mathcal{M}}$ to be finitary? Proposition~\ref{P:ultraclosedness} below proposes an answer to this question.

However, first we discuss the property that the matrix consequence of any finite matrix is finitary. To prove this, we use topological methods discussed in Section~\ref{section:topology}. 

Let $\mat{M}=\langle\alg{A},D\rangle$ be a finite logical matrix of type $\Lan$ and let $\valA:=\alg{A}^{\Var_{\Lan}}$ (the set of all valuations in algebra $\alg{A}$). We regard $|\alg{A}|$ as a finite topological space with discrete topology and $\valA$ as the cartesian power of the space $|\alg{A}|$ relative to the product topology. According to~Corollary~\ref{C:product-finite-discrete-spaces}, the space $\valA$ is compact.

For any formula $\alpha\in\Forms_{\mathcal{L}}$, we denote:
\[
\valA(\alpha):=\set{v\in\valA}{v[\alpha]\in D}.
\]

We aim to show that each $\valA(\alpha)$ is open and closed in $\valA$. The idea is that, since an arbitrary formula $\alpha$ contains finitely many sentential variables and the set $|\alg{A}|$ is finite, there are only finitely many restricted assignments that validate $\alpha$ and there are finitely many restricted assignments that refute it in $\mat{M}$.

We show that each $\valA(\alpha)$ is an open set in the product topology. The closedness of it, that is the openness of its complement, is similar.

Given an arbitrary formula $\alpha$, we define a binary relation on $\valA(\alpha)$ as follows:
\[
v\equiv w~\stackrel{\text{df}}{\Longleftrightarrow}~v(p)=w(p),~\text{for any $p\in\Var(\alpha)$}.
\] 

It is clear that the relation $\equiv$ is an equivalence. For any $v\in\valA(\alpha)$, we denote:
\[
|v|:=\set{w\in\valA(\alpha)}{w\equiv v}. \tag{the equivalence class generated by $v$}
\]

It is well know that $\set{|v|}{v\in\valA(\alpha)}$ is a partition of $\valA(\alpha)$. If we show that each $|v|$ is open in the product topology, so will be $\valA(\alpha)$.

We observe that
\[
|v|=\!\!\!\prod_{~p\in\Var_{\Lan}}Z_p,
\]
where
\[
Z_p=\begin{cases}
	\begin{array}{cl}
		\lbrace v(p)\rbrace &\text{if $p\in\valA(\alpha)$}\\
		|\alg{A}| &\text{otherwise}.
	\end{array}
\end{cases}
\]

According to the description~\eqref{E:representation} of the product topology in Section~\ref{section:topology}, this demonstrates that each $|v|$ is open. In a similar manner, we prove that each
\[
\valA\setminus\valA(\alpha)=\set{v\in\valA}{v[\alpha]\notin D}
\]
is also open.\\

Now, assume that for a set $X\cup\lbrace\alpha\rbrace\subseteq\Forms_{\mathcal{L}}$,
$X\models_{\textbf{M}}\alpha$. This implies (is even equivalent to) that $\bigcap_{\beta\in X}\valA(\beta)\subseteq\valA(\alpha)$. The last inclusion in turn is equivalent to the following:
\begin{equation}\label{E:open-cover}
	\valA(\alpha)\cup\bigcup_{\beta\in X}(\valA\setminus\valA(\beta))=
	\valA.
\end{equation}

This means that $\lbrace\valA(\alpha\rbrace\cup\lbrace(\valA\setminus\valA(\beta))
\rbrace_{\beta\in X}$ is an open cover of $\valA$. Since the space 
$\valA$ is compact, it has a finite subcover 
$\lbrace\valA(\alpha\rbrace\cup\lbrace(\valA\setminus\valA(\beta))
\rbrace_{\beta\in Y}$, for some $Y\Subset X$. The latter is equivalent to $Y\models_{\textbf{M}}\alpha$. 

Thus we have proved the following.
\begin{prop}\label{P:con-finite matrix-is-finitary}
	A single-matrix consequence relative to a finite matrix is finitary.
\end{prop}

\begin{cor}\label{C:finitariness-M-consequence}
	Let $\mathcal{M}$ be a nonempty finite family of finite matrices of type $\Lan$. Then the $\mathcal{M}$-consequence is finitary.
\end{cor}
\noindent\textit{Proof}~is left to the reader. (Exercise~\ref{section:single-matrix-consequence}.\ref{EX:finitariness-M-consequence})\\

The next example shows that the converse of the statement of Proposition~\ref{P:con-finite matrix-is-finitary} in general is not true.
\begin{example}\label{Example:logic-without-finite-matrix}
	{\em Let $\Lan$ be a propositional language with any set (perhaps finite) of propositional variables and a single binary connective $\circledast$. An abstract logic $\aLog^{\circledast}$ is defined by the premiseless rule $\bm{\alpha}\circledast\bm{\alpha}$ and modus ponens, $\bm{\alpha},\,\bm{\alpha}\circledast\bm{\beta}\slash\bm{\beta}$. We claim that $\aLog^{\circledast}$ has no finite adequate matrix.
	}
\end{example}

We will prove that $\aLog^{\circledast}$ has no finite weakly adequate matrix; this will imply that it has no adequate matrix as well.

First of all, we note that the logic $\aLog^{\circledast}$ is structural and finitary. Secondly, we notice the set of $\aLog^{\circledast}$-theses consists of the formulas of the form $\alpha\circledast\alpha$, for any $\alpha\in\FormsL$.

Now, for contradiction, assume that a matrix $\mat{M}:=\langle\textsf{A};\circledast,D\rangle$ is weakly adequate for $\aLog^{\circledast}$ and has $n$ elements. The observations above (about $\aLog^{\circledast}$-theses) implies that for any $a\in\textsf{A}$, $a\circledast a\in D$.

Let us definite an infinite sequence of terms as follows:
\[
\begin{array}{c}
	t^{1}(x)=x\circledast x\\
	t^{n+1}(x)=t(t^{n}(x)).
\end{array}
\]

It must be clear that for an arbitrary $a\in\textsf{A}$, there is a natural number $k\ge 1$ such that $t^{k}(a)=a$. Assume that $a=v[\alpha]$, for some formula $\alpha$ at some valuation $v$ in \alg{A}. We note that $t^{k}(\alpha)$ is a formula unequal to  $\alpha$ and hence $t^{k}(\alpha)\circledast\alpha$ is not an $\aLog^{\circledast}$-theorem, though $v[t^{k}(\alpha)\circledast\alpha]\in D$. A contradiction.\\

Now we turn to a more general case.

\begin{prop}[Bloom-Zygmunt criterion]\label{P:ultraclosedness}\index{Theorem!Bloom-Zygmunt}
	Let $\mathcal{M}$ be the class of all models of an abstract logic $\mathcal{S}$ in $\Lan$. Then $\mathcal{S}$ is finitary if, and only if, the class $\mathcal{M}$ is closed with respect to ultraproducts.\footnote{Ultraproducts are understood in the sense of~\eqref{E:ultraproduct-df}.} 
\end{prop}
\begin{proof}
	Suppose $\mathcal{S}$ is finitary. Then, according to Corollary~\ref{C:logic-by-modus-rules-criterion}, there is a set $\mathcal{R}$ of modus rules such that $\mathcal{S}=\mathcal{S}_{\mathcal{R}}$. Thus, $\mathcal{M}$ is the class of all $\mathcal{S}_{\mathcal{R}}$-models. Then, in virtue of Corollary~\ref{C:S-models-modus-rules}, $\mathcal{M}$ is closed under ultraproducts.
	
	Next, let the class $\mathcal{M}$ be closed with respect to ultraproducts. For contradiction, assume that $\mathcal{S}$ is not finitary. That is, for some set $X\cup\lbrace\alpha\rbrace\subseteq\Forms_{\mathcal{L}}$, $\alpha\in\ConS{X}$ but for any $Y\Subset X$, $\alpha\notin Y$. 
	
	We denote:
	\[
	I:=\set{Y}{Y\Subset X}
	\]
	and for any $Y\in I$,
	\[
	\langle Y\rangle:=\set{Z\in I}{Y\subseteq Z}.
	\]
	
	We observe that for any $Y_1,\ldots,Y_n\in I$,
	\[
	\langle Y_1\rangle\cap\ldots\cap\langle Y_n\rangle=\langle Y_1\cap\ldots\cap Y_n\rangle.
	\]
	That is, the set $\set{\langle Y\rangle}{Y\in I}$ has the finite intersection property. According to the Ultrafilter Theorem (\cite{chang-keisler1990}, proposition 4.1.3), there is an ultrafilter $\nabla$ over $I$ such that
	\begin{equation}\label{E:ultrafilter-1}
		\set{\langle Y\rangle}{Y\in I}\subseteq\nabla.
	\end{equation}
	
	Given an $\Lan$-formula $\beta\in X$, we note that
	\[
	\langle\lbrace\beta\rbrace\rangle=\set{Y\in I}{\beta\in Y}\subseteq\set{Y\in I}{\beta\in\ConS{Y}}.
	\]
	This, in view of~\eqref{E:ultrafilter-1}, implies that for arbitrary formula $\beta\in X$,
	\begin{equation}\label{E:ultrafilter-2}
		\set{Y\in I}{\beta\in\ConS{Y}}\in\nabla.
	\end{equation}
	
	Using the notation 	$\Lin_{\,\mathcal{S}}[Y]:=\langle\FormAl,\ConS{Y}\rangle$ (Section~\ref{section:con-via-matrices}), we define:
	\[
	\mat{M}:=\prod_{F}\Lin_{\,\mathcal{S}}[Y].
	\]
	Thus
	\[
	\mat{M}=\langle(\FormAl)^{I}\slash\nabla,D\rangle,
	\]
	where a set $D\subseteq |(\FormAl)^{I}\slash\nabla|$ and such that for any $\theta\in|(\FormAl)^{I}|$,
	\begin{equation}\label{E:ultrafilter-3}
		\theta\slash\nabla\in D~\Longleftrightarrow~\set{Y\in I}{\theta_{Y}\in\ConS{Y}}\in\nabla.
	\end{equation}
	
	For any $\Lan$-formula $\beta$, we define:
	\[
	\overline{\beta}:I\longrightarrow\lbrace\beta\rbrace
	\]
	and then, for any $\beta\in\Forms_{\mathcal{L}}$,
	\[
	v:\beta\mapsto\overline{\beta}\slash\nabla.
	\]
	
	It must be clear that $v:\FormAl\longrightarrow(\FormAl)^{I}\slash\nabla$ is a homomorphism and, hence, can be regarded as an assignment of the $\Lan$-formulas in $(\FormAl)^{I}\slash\nabla$.
	
	We aim to show that $v[X]\subseteq D$ but $v[\alpha]\notin D$. This will imply that $\mat{M}$ is not an $\mathcal{S}$-model and hence does not belong to $\mathcal{M}$. The latter in turn will implies that $\mathcal{M}$ is not closed under ultraproducts, which is contrary to the assumption.
	
	Indeed, it follows from \eqref{E:ultrafilter-3} and \eqref{E:ultrafilter-2} that for any $\beta\in X$,
	\[
	v[\beta]=\overline{\beta}\slash\nabla\in D.
	\]
	On the other hand, 
	\[
	\set{Y\in I}{\alpha\in\ConS{Y}}=\varnothing.
	\]
	Therefore, in view of~\eqref{E:ultrafilter-3}, $v[\alpha]=\overline{\alpha}\slash\nabla\notin D$. 
\end{proof}

\paragraph{Exercises~\ref{section:finitary-matrix-consequence}}
\begin{enumerate}
	\item \label{EX:finitariness-M-consequence} Prove Corollary~\ref{C:finitariness-M-consequence}.
	\item Show that the class of all models of $\Cl$ (Section~\ref{section:inference-rules}) is closed with respect to ultraproducts.
	\item Let $\mathcal{R}$ be a set of modus rules. Prove that the abstract logic $\mathcal{S}_{\mathcal{R}}$ (Section~\ref{section:modus-rules-general}) is finitary.
\end{enumerate}

\section{The conception of separating tools}\label{section:separating-means}
Considering an abstract logic $\mathcal{S}$ as a consequence relation $\vdash_{\mathcal{S}}$, there are two main tasks: to confirm $X\vdash_{\mathcal{S}}\alpha$ or to refute it. Given a finitary abstract logic $\mathcal{S}$, these two tasks can be regarded as one twofold problem (of extreme importance): Is  there an algorithm that for any set $X\cup\lbrace\alpha\rbrace\Subset\Forms_{\mathcal{L}}$, decides whether $X\vdash_{\mathcal{S}}\alpha$ is true or false.

The importance of the second part of this dual-task was emphasized by J. {\L}ukasiewicz; namely, he wrote:
\begin{quote}
	Of two intellectual acts, to assert a proposition and to reject it,\footnote{[{\L}ukasiewicz's footnote]: ``I owe this distinction to Franz Brentano, who describes the act of believing as \emph{anerkennen} and \emph{verwerfen}.''} only the first has been taken into account in modern formal logic. Gottlob Frege introduced into logic the idea of assertion, and the sign of assertion ($\bm{\vdash}$), accepted afterwards by the authors of \emph{Principia Mathematica}. The idea of rejection, however, so far as I know, has been neglected up to the present day. \cite{luk57}, {\S} 27
\end{quote}

Trying to confirm $X\vdash_{\mathcal{S}}\alpha$, one can use rules of inference (sound with respect to $\mathcal{S}$) or a suitable adequate matrix or a set of matrices which determine $\mathcal{S}$. On the other hand, to refute $X\vdash_{\mathcal{S}}\alpha$, we need ``means'' that would separate $X$ from $\alpha$, for any $X\cup\lbrace\alpha\rbrace\subseteq\Forms_{\mathcal{L}}$ or at least for any $X\cup\lbrace\alpha\rbrace\Subset\Forms_{\mathcal{L}}$. According to Corollary~\ref{C:Lindenbaum-completeness}, the Lindenbaum atlas for $\mathcal{S}$ is a universal tool, for it works, at least in principal, for both tasks. However, it is too complicated and can be used mainly to prove very general theorems. For the refutation task, the Lindenbaum atlas is very inefficient. Nevertheless, in many cases, as we will see, it can be simplified to an extent to be useful.

A general view on the conception of \textit{separating tools} was formulated by A. Kuznetsov~\cite{kuznetsov1979}. According to him, if all formulas of a set $X$ stand in a relation $R$ to an object $\mathfrak{M}$, but a formula $\alpha$ does not stand in $R$ to $\mathfrak{M}$, we say that $\mathfrak{M}$ \textit{separates} $X$ from $\alpha$ with respect to $R$.
In this case, the character of the relation $R$ must be such that from the separation of $X$ from $\alpha$ it follows that $X\vdash_{\mathcal{S}}\alpha$ does not hold.\footnote{The conception of separating tools, proposed by Kuznetsov, was intended to be applied not only to the issue of derivability, but also to that of expressibility understood within the framework of a given calculus.} He also emphasized that the object $\mathfrak{M}$ which was employed as a separating tools in the above sense has often been of algebraic nature.\footnote{One important exception is Kleene's notion of realizability; see~\cite{kleene1952}, {\S} 82.}

While the Lindenbaum atlas relative to an abstract logic $\mathcal{S}$ can play a role of universal separating tools, but for a particular pair $X$ and $\alpha$, in case 
$X\vdash_{\mathcal{S}}\alpha$ does not hold, potentially any $\mathcal{S}$-matrix can be a separating tools. Sometimes a suitable separating model may not be the Lindenbaum matrix itself, but another matrix whose algebraic carrier is related to the algebraic carrier of the former. Here it is appropriate to recall the corelations~\eqref{E:M_subseteq_M^ast} and~\eqref{E:vdash_M=vdash_M^ast} and the final remark about them in Section~\ref{section:realizations-abstract-logic}.
Therefore, it would be interesting to find an algebraic connection between the Lindenbaum atlas or Lindenbaum matrix relative to $\mathcal{S}$ and any $\mathcal{S}$-matrix. This will be the topic of the next chapter.

\section{Historical notes}\label{section:consequence-historical-notes}

Although Proposition~\ref{P:los-suszko}, Proposition~\ref{P:wojcicki} and Proposition~\ref{P:shoesmith-smiley-theorem} were proven a bit longer than a half century ago, their history is somewhat complicated. Therefore, we break the statements of Proposition~\ref{P:los-suszko} and Proposition~\ref{P:wojcicki} into several part and add several new ones. This will allow us to see the logical structure of these statements, as well as a few more. Then, all statements will be analyzed from a historical point of view.
\[
\begin{array}{cl}
	(\text{A}) &\text{an abstract logic $\mathcal{S}$ has an adequate matrix};\\
	(\text{B}) &\text{an abstract logic $\mathcal{S}$ is structural};\\
	(\text{C}) &\text{an abstract logic $\mathcal{S}$ is finitary};\\
	(\text{D}) &\text{an abstract logic $\mathcal{S}$ is uniform (Definition~\ref{D:uniform-logic})};\\
	(\text{E}) &\text{an abstract logic $\mathcal{S}$ is couniform (Definition~\ref{D:couniform-logic})};\\
	(\text{F}) &\text{an abstract logic $\mathcal{S}$ has the cancellation property (Definition~\ref{D:cancellation-property})};\\
	(\text{G}) &\text{an abstract logic $\mathcal{S}$ is \emph{regular}
		(see the definition of a regular logic}\\
	&\text{in~\cite{wojcicki1969} or in~\cite{wojcicki1970}, section 1}.\\
	(\text{H}) &\text{an abstract logic $\mathcal{S}$ is \emph{absolutely separable}
		(see the definition}\\
	&\text{of an absolutely separable logic in~\cite{wojcicki1984}, section 34.6)}.\\
	
\end{array}
\]

The original theorem of Lindenbaum, as was stated in~\cite{lukasiewicz-tarski1930}, theorem 3, concerned, using the formulation of~\cite{zygmunt2012}, the ``(weak) adequacy problem.'' That is, as the aforementioned theorem 3 was clarified and proved in~\cite{los1949}, theorem 10, given an abstract logic $\mathcal{S}$, any  $\mathcal{S}$-theory $X$ is determined by a single logical matrix, namely $\LinS[X]$. And, again, according to~\cite{zygmunt2012}, \cite{los-suszko1958} turned this problem into the ``(strong) adequacy problem;'' namely, they raised the question of which abstract logic can be determined using a single logical matrix. 

Having introduced the notion of a uniform abstract logic, {\L}o\'{s} and Suszko then attempted to prove that a structural abstract logic defines a single-matrix consequence if, and only if, it is uniform.
The `only-if' part of this criterion was stated in~\cite{los-suszko1958}, section 8, p. 181, while the `if' part was presented as ``the main theorem.'' Thus, their combined statement was
\begin{equation}\label{E:los-suszko-original}
	(\text{A})~\Longleftrightarrow~(\text{B})\&(\text{D})
\end{equation}

However, R. W\'{o}jcicki noticed in~\cite{wojcicki1969} that, while
\[
(\text{A})~\Longrightarrow~(\text{B})\&(\text{D}) 
\]
is true, the converse,
\[
(\text{B})\&(\text{D})~\Longrightarrow~(\text{A}),
\]
in general is not. Since the implication
\[
(\text{A})~\Longrightarrow~(\text{B}) 
\]
(Lemma~\ref{L:los-suszko-1}) is true, the problem was with the implication
\[
(\text{D})~\Longrightarrow~(\text{A}),
\]
providing that the logic $\mathcal{S}$ is structural. Thus, in the sequel, we always assume (B) as a premise.

Fixing a lacuna in {\L}o\'{s} and Suszko's argumentation for~\eqref{E:los-suszko-original}, W\'{o}cicki
announced in~\cite{wojcicki1969}, theorem 2, and proved in~\cite{wojcicki1970}, theorem 2.2, that
\[
(\text{B})~\Longrightarrow~((\text{D})\&(\text{G})~\Longrightarrow~(\text{A})).
\]
Also, he proved in~\cite{wojcicki1970} that
\[
(\text{B})~\Longrightarrow~((\text{C})~\Longrightarrow~(\text{G}))
\]
(assertion 1.3) and
\[
(\text{B})~\Longrightarrow~((\text{A})~\Longrightarrow~(\text{D})\&(\text{G})).
\]
(assertion 2.1). Thus, it was proven therein that
\[
(\text{B})~\Longrightarrow~((\text{A})~\Longleftrightarrow~(\text{D})\&(\text{G})).
\]

Later on, in~\cite{wojcicki1984}, theorem 34.14, W\'{o}jcicki replaced regularity with absolute separability, i.e.
\[
(\text{B})~\Longrightarrow~((\text{A})~\Longleftrightarrow~(\text{D})\&(\text{H}))
\]
and in~\cite{wojcicki1988}, theorem 3.2.5, with finitariness, i.e.
\[
(\text{B})\&(\text{C})~\Longrightarrow~((\text{A})~\Longleftrightarrow~(\text{D})).
\]
(Proposition~\ref{P:los-suszko}) which we call the \emph{{\L}o\'{s}-Suszko-W\'{o}jcicki theorem}. 

In~\cite{wojcicki1988}, he introduced a couniform logic and proved (theorem 3.2.7) that
\[
(\text{B})~\Longrightarrow~((\text{A})~\Longleftrightarrow~(\text{D})\&(\text{E}))
\]
(Proposition~\ref{P:wojcicki}) which we call the \textit{W\'{o}jcicki theorem}.

In~\cite{shoesmith-smiley1971}, not referring to either~\cite{los-suszko1958} or~\cite{wojcicki1969,wojcicki1970}, D. J. Shoesmith and T. J. Smiley gave their version of the adequacy theorem: structural abstract logic has an adequate matrix if, and only if, it has the cancellation property, i.e.
\[
(\text{B})~\Longrightarrow~((\text{A})~\Longleftrightarrow~(\text{F}))
\] 
(Proposition~\ref{P:shoesmith-smiley-theorem}) which we call the \emph{Shoesmith-Smiley theorem}.

Proposition~\ref{P:con-finite matrix-is-finitary} is due to~\cite{los-suszko1958}, section 8, p.181. 

Example~\ref{Example:logic-without-finite-matrix} is due to~\cite{mckinsey-tarski1948}, section 2, Remark on p. 6.

The `only-if' part of Proposition~\ref{P:ultraclosedness} was proved in~\cite{zygmunt1974}, corollary; the full equivalence is due to~\cite{bloom1975}, theorem 2.6.

\chapter[Unital Abstract Logics]{Unital Abstract Logics}
\label{chapter:lindenbaum-tarski}

\section{Unital algebraic expansions}	\label{section:chapter-unital-preliminaries}
From the discussion in Section~\ref{section:separating-means}, the following questions arise. 
\begin{itemize}
	\item[(1)]~\,Is there a regular method of simplifying the Lindenbaum atlas (Lindenbaum matrix) relative to an abstract logic $\mathcal{S}$ to obtain an adequate atlas (or matrix), or at least weakly adequate matrix that would be a more working separating tool than $\Lin[\Sigma_{\mathcal{S}}]$ or $\LinS$?
	\item[(2)]~\,Is there a regular method of simplifying the Lindenbaum atlas (Lindenbaum matrix) relative to an abstract logic $\mathcal{S}$ to have a connection between these simplified forms of $\Lin[\Sigma_{\mathcal{S}}]$ or $\LinS$ and other $\mathcal{S}$-atlases and $\mathcal{S}$-matrices, respectively?
\end{itemize}

As to the question $(1)$, any finite matrix determines an abstract logic, the Lindenbaum matrix of which is always infinite, for the cardinality of any Lindenbaum matrix equals the cardinality of the set of formulas. Take, for instance, the matrix $\booleTwo$ of Section~\ref{S:two-valued}. In Section~\ref{section:rules-and-hyperrules} we have proved the equality $\models_{\textbf{B}_2}~=~\vdash_{\textsf{Cl}}$. The latter corresponds to the abstract logic $\Cl$ (classical propositional logic). Even a rough comparison between $\Lin_{\textsf{Cl}}$ and $\booleTwo$ demonstrates the difference between their complexities.

One thing is remarkable about the logical matrices discussed in Sections~\ref{S:two-valued}--\ref{section:dummett}. The logical filters of all of them consist of a single element. A logical matrix $\mat{M}=\langle\alg{A},D\rangle$ is called \textit{\textbf{unital}} if $D$ is a one-element set. We denote the designated element of a unital matrix by $\one$ (perhaps with subscript). Thus, if $\mat{M}=\langle\alg{A},\lbrace\one\rbrace\rangle$, we can consider  
the algebra $\langle\alg{A},\one\rangle$ as a semantics
for $\Lan$-formulas.  We call $\langle\alg{A},\one\rangle$ a \textit{\textbf{unital expansion}} of $\alg{A}$. Very often, such an expansion is not needed. Namely, when $D=\lbrace c\rbrace$, where $c$ is a constant from the signature of {\alg{A}}, we count that $\one=c$. In this event, we count $\langle\alg{A},\one\rangle$ the same as $\alg{A}$.

Thus, in case \mat{M} is unital, instead of saying that an $\Lan$-formula $\alpha$ is valid in $\mat{M}$, we say that $\alpha$ is \textit{\textbf{valid}} in $\langle\alg{A},\one\rangle$, meaning that for any valuation $v$ in $\alg{A}$, $v[\alpha]=\one$. Further, given an algebra $\langle\alg{A},\one\rangle$, for any set $X\cup\lbrace\alpha\rbrace\subseteq\Forms_{\mathcal{L}}$, we define: 
\[
X\models_{\langle\textbf{A},\one\rangle}\alpha\stackrel{\text{df}}{\Longleftrightarrow}v[X]\subseteq\lbrace\one\rbrace\Rightarrow v[\alpha]=\one,~\text{for any valuation $v$ in $\alg{A}$}.
\]

The last definition allows us to adapt the definition of $\mathcal{S}$-model to algebras
$\langle\alg{A},\one\rangle$. Namely such an algebra is an \textit{\textbf{algebraic}} $\mathcal{S}$-\textit{\textbf{model}} if for any set $X\cup\lbrace\alpha\rbrace\subseteq\Forms_{\mathcal{L}}$,
\[
X\vdash_{\mathcal{S}}\alpha~\Longrightarrow~X\models_{\langle\textbf{A},\one\rangle}\alpha.
\]

Further, we say that an algebra $\langle\alg{A},\one\rangle$ is \textit{\textbf{adequate}} for an abstract logic $\mathcal{S}$ if for any set $X\cup\lbrace\alpha\rbrace\subseteq\Forms_{\mathcal{L}}$,
\[
X\vdash_{\mathcal{S}}\alpha~\Longleftrightarrow~X\models_{\langle\textbf{A},\one\rangle}\alpha;
\]
and that $\alg{A}$ is \textit{\textbf{weakly adequate}} for $\mathcal{S}$ if for any $\alpha\in\Forms_{\mathcal{L}}$,
\[
\alpha\in\ThmS~\Longleftrightarrow\text{for any valuation $v$ in $\alg{A}$},
v[\alpha]=\one.
\]

The definition of $\mathcal{M}$-consequence (Section~\ref{section:con-via-matrices}) can be applied to the case when $\mathcal{M}=\lbrace\langle\alg{A}_i,\one_i\rangle\rbrace_{i\in I}$. Namely, for any (nonempty) such class $\mathcal{M}$, we define:
\[
X\models_{\mathcal{M}}\alpha~\stackrel{\text{df}}{\Longleftrightarrow}~
X\models_{\langle\textbf{A}_i,\one_i\rangle}\alpha,~\text{for each $\langle\alg{A}_i,\one_i\rangle$}.
\]
An abstract logic which can be defined as an $\mathcal{M}$-consequence of this type is called \textit{\textbf{assertional}}; cf~\cite{font2016}, definition 3.5. Thus the assertional abstract logics are those which are determined by classes of their algebraic models.\\

It will be convenient to treat the validity of an $\Lan$-formula $\alpha$ in 
$\langle\alg{A},\one\rangle$ as the validity of a first-order formula in a first-order model.

Namely, given a sentential language $\Lan$, we define a first-order language $\FOstar$ whose individual variables are the sentential variables of $\Lan$ and the connectives of $\Func_{\mathcal{L}}$ and the constants of $\Cons_{\mathcal{L}}$ are the functional symbols. (The constant symbols are treated as 0-ary functional symbols). In addition,  $\FOstar$ has a 0-ary functional symbol $\one$ which is interpreted in the $\FOstar$-models by a constant $\one$. Also, $\FOstar$ has symbol `$\approx$' which is interpreted in the $\FOstar$-models as equality, the first-order logical connectives
$\bm{\&}$ (conjunction), $\Rarrow$ (implication) and universal quantifier $\Forall$.
The parentheses `(' and `)' are used in the usual way.

Thus each $\Lan$-formula becomes in $\FOstar$ a term. In addition, $\FOstar$ has terms like
$\alpha(p,\ldots,\one,\ldots, q)$ obtained from $\alpha(p,\ldots,r,\ldots, q)$ by replacement of a variable $r$ with the constant $\one$. 

The atomic formulas of $\FOstar$, called also \textit{\textbf{equalities}}, are the expressions of the form:
\[
\mathfrak{r}\approx\mathfrak{t},
\]
where $\mathfrak{r}$ and $\mathfrak{t}$ are $\FOstar$-terms.

Thus $\FOstar$ as a logical system is a first-order logic with equality.\\

With each expansion $\langle\alg{A},\one\rangle$ we associate an $\FOstar$-model $\fA=\langle\alg{A},\one,\approx\rangle$ so that,
talking about validity of identities and quasiidentities in $\fA$, we, referring to~\eqref{E:equivalence-for-modeling}, will be writing
\[
\fA\models\alpha(p,\one,\ldots)\approx\beta(q,\one,\ldots)~~\text{and}~~
\fA\models(\phi_1\bm{\&}\ldots\bm{\&}\phi_n)\Rarrow\phi
\]
instead of, respectively, $\fA\models\Forall\ldots\Forall(\alpha(p,\one,\ldots)\approx\beta(q,\one,\ldots))$
and\\ $\fA\models\Forall\ldots\Forall((\phi_1\bm{\&}\ldots\bm{\&}\phi_n)\Rarrow\phi)$.
Thus, it must be clear that for any $\Lan$-formula $\alpha$ and any $\FOstar$-model
$\fA=\langle\alg{A},\one,\approx\rangle$,
\begin{equation}\label{E:model=algebra-validity-1}
	\models_{\langle\textbf{A},\one\rangle}\alpha~\Longleftrightarrow~\fA\models
	\Forall\ldots\Forall(\alpha\approx\one).
\end{equation}

Let $X:=\lbrace\alpha_1,\ldots,\alpha_n\rbrace$. We denote:
\[
X\approx\one\stackrel{\text{df}}{\Longleftrightarrow}\alpha_1\approx\one\&\ldots\&\alpha_n\approx\one.
\]
Thus, in this notation, for any set $X\cup\lbrace\alpha\rbrace\Subset\Forms_{\mathcal{L}}$, we can observe that
\begin{equation}\label{E:model=algebra-validity-2}
	X\models_{\langle\textbf{A},\one\rangle}\alpha~\Longleftrightarrow~\fA\models \Forall\ldots\Forall(X\approx\one\Rarrow\alpha\approx\one).
\end{equation}
(If $X=\varnothing$, we understand the last equivalence as~\eqref{E:model=algebra-validity-1}.)

\begin{prop}\label{P:S-algebraic-models=q-variety}
	Let $\mathcal{S}$ be a finitary assertional abstract logic. Then all algebraic $\mathcal{S}$-models form a quasi-variety.
\end{prop}
\begin{proof}
	Indeed, let $\left\langle \alg{A},\one\right\rangle $ be a unital algebraic expansion. Then, since $\mathcal{S}$ is finitary, $\left\langle \alg{A},\one\right\rangle$ is an $\mathcal{S}$-model if, and only if,
	$\left\langle \alg{A},\one,\approx\right\rangle\models X\approx\one\Rarrow \alpha\approx\one$, for any $X\cup\lbrace\alpha\rbrace\Subset\Forms_{\mathcal{L}}$ such that $X\vdash_{\mathcal{S}}\alpha$.	
\end{proof}

\paragraph{Exercises~\ref{section:chapter-unital-preliminaries}}
\begin{enumerate}
	\item Prove that if an abstract logic $\aLog$ has a finite unital adequate matrix, then the class of all algebraic $\aLog$-models is a quasi-variety.
	\item Show that the class of all algebraic $\Cl$-models is a quasi-variety.
\end{enumerate}

\section{Unital abstract logics}\label{section:unital-logics}
Given a set $X\subseteq\Forms_{\mathcal{L}}$, we denote by $\theta(X)$ the congruence on $\FormAl$, generated by the set $X\times X$. We will be interested in congruences $\theta(D)$, where $D$ is a theory of an abstract logic  $\mathcal{S}$, in particular in congruences $\theta(\ThmS)$.

The following observation will be useful in the sequel.
\begin{prop}\label{P:motivating}
	A structural abstract logic $\mathcal{S}$ is nontrivial if, and only if, it satisfies the following property: For any $p,q\in\Var_{\mathcal{L}}$,
	\begin{equation}\label{E:motivating}
		p\neq q~\Longrightarrow~p\slash\theta(\ThmS)\neq q\slash\theta(\ThmS).
	\end{equation}
\end{prop}
\begin{proof}
	To prove the `only-if' part, we assume that $\mathcal{S}$ is nontrivial. Then, since $\mathcal{S}$ is structural, $\ThmS$ does not contain sentential variables. Let $M$ be the congruence class with respect to $\theta(\ThmS)$, which contains the variable $p$. Now, we define: $M_1:=\lbrace p\rbrace$ and $M_2:=M\setminus\lbrace p\rbrace$. The new partition of $\Forms_{\mathcal{L}}$ induces the equivalence $\theta$.
	It is obvious that $\theta\subset\theta(\ThmS)$. We will show that $\theta$ is a congruence on $\FormAl$.
	
	Let $F$ be an arbitrary $n$-ary connective with $n\ge 1$. Assume that for formulas $\alpha_1,\ldots,\alpha_n,\beta_1,\ldots,\beta_n$, each pair $(\alpha_i,\beta_i)\in\theta$.
	Then each pair $(\alpha_i,\beta_i)\in\theta(\ThmS)$ and hence 
	$(F\alpha_1\ldots\alpha_n,F\beta_1\ldots\beta_n)\in\theta(\ThmS)$. This means that either both $F\alpha_1\ldots\alpha_n$ and $F\beta_1\ldots\beta_n$ belong to $M_2$ or they belong to a class of $\theta(\ThmS)$ that is different from $M$ and hence is a class of $\theta$.
	This implies that $\lbrace p\rbrace$ and $\lbrace q\rbrace$ are two distinct congruence classes with respect to $\theta(\ThmS)$.
	
	Conversely, suppose the implication \eqref{E:motivating} holds. Then, since one of the congruence classes with respect to $\theta(\ThmS)$ contains $\ThmS$ and there are more than two classes,
	$\ThmS\neq\Forms_{\mathcal{L}}$.
\end{proof}

\subsection{Definition and some properties}

\begin{defn}[unital abstract logic]\label{D:unital-logic}\index{logic!unital}
	A structural abstract logic $\mathcal{S}$ is called \textbf{unital} if any $\mathcal{S}$-theory $D$ is a congruence class with respect to the congruence $\theta(D)$.
\end{defn}

In the sequel (Chapter~\ref{chapter:Q-consequence}), we also use the following notation:
\[
\theta_{\aLog}(X)~\define~\theta(\ConS{X}).
\]
where $\aLog$ is  a given abstract logic.\\

We need the following property in the sequel.
\begin{prop}
	Given a unital logic $\aLog$ and a set {\em$X\subseteq\Forms_{\Lan}$}, the abstract logic $\aLog[X]$ {\em(Section~\ref{section:consequence-operator})} is also unital.
\end{prop}
\begin{proof}
	We observe that any $\aLog[X]$-theory is also an $\aLog$-theory.
\end{proof}

Since any congruence class cannot be empty, from Definition~\ref{D:unital-logic} and Proposition~\ref{P:motivating}, we immediately obtain the following.
\begin{prop}\label{P:lindenbaum-connection}
	For any unital abstract logic $\mathcal{S}$, $\ThmS\neq\varnothing$; in addition, $\mathcal{S}$ is nontrivial if, and only if,  \eqref{E:motivating} holds.
\end{prop}

As we will see below (Corollary~\ref{C:compatible-congruence-2}), the unitality of an abstract logic is connected to the following concept.

\begin{defn}\label{D:compatible-congruence}\index{compatible congruence}
	Let {\em$X\subseteq\Forms_{\mathcal{L}}$} and $X\neq\varnothing$. A congruence $\theta$ on $\FormAl$ is said to be \textbf{compatible} with $X$
	if for any $(\alpha,\beta)\in\theta$,
	\[  
	\alpha\in X~\Longleftrightarrow~\beta\in X;
	\]
	in other words, $X$ is the union of congruence classes with respect to $\theta$.
\end{defn}

It is obvious that 
\[
\theta_0:=\set{(\alpha,\alpha)}{\alpha\in\FormAl}
\]
is a congruence on $\FormAl$, which is compatible with any nonempty set $X\subseteq\Forms_{\Lan}$.

\begin{prop}\label{P:compatible-congruence}
	Let {\em$X\subseteq\Forms_{\mathcal{L}}$} and $X\neq\varnothing$.
	For any congruence $\theta$ on $\FormAl$, $X$ is a congruence class with respect to $\theta$ if, and only if, the following conditions are satisfied:
	{\em\[
		\begin{array}{cl}
			(\text{a}) &\theta(X)\subseteq\theta;\\
			(\text{b}) &\theta~\textit{is compatible with}~ X.
		\end{array}
		\]}
\end{prop}
\begin{proof}
	Suppose $X$ is a congruence class with respect to $\theta$. Then (a) and (b) are obviously true.
	
	Conversely, suppose both (a) and (b) are fulfilled. Then, according to (a), all formulas in $X$ are $\theta$-congruent to one another. In addition, in virtue of (b), each formula in $X$ cannot be $\theta$-congruent to a formula beyond $X$.
\end{proof}

The next two corollaries follow immediately from Proposition~\ref{P:compatible-congruence}.
\begin{cor}\label{C:compatible-congruence-1}
	Let $\mathcal{S}$ be an abstract logic with $\ThmS\neq\varnothing$ and $D$ be an $\mathcal{S}$-theory. Then $D$ is a congruence class with respect to $\theta(D)$ if, and only if, 
	$\theta(D)$ is compatible with $D$.
\end{cor}
\begin{cor}\label{C:compatible-congruence-2}
	Let $\mathcal{S}$ be a structural abstract logic with $\ThmS\neq\varnothing$. Then the following conditions are equivalent$\,:$
	{\em\[
		\begin{array}{cl}
			(\text{a}) & \textit{$\mathcal{S}$ is unital};\\
			(\text{b}) & \textit{$\theta(D)$ is compatible with $D$, for each $\mathcal{S}$-theory $D$};\\
			(\text{c}) & \textit{for any $\mathcal{S}$-theory $D$, there is a congruence $\theta$ on $\FormAl$}\\
			&\textit{such that $\theta(D)\subseteq\theta$ and $\theta$ is compatible with $D$}.
		\end{array}
		\]}
\end{cor}

\subsection{Some subclasses of unital logics}\label{section:subclasses-unital}
In this subsection, we introduce several classes of abstract logics, which were considered in the literature and which turned out to be included in the class of unital logics.

\begin{defn}[implicative logic]\label{D:weakly-implicative}\index{logic!implicative}
	A structural logic $\mathcal{S}$ in a language $\Lan$ is called \textbf{implicative} if $\Lan$ contains a binary connective $ \rightarrow $ such that for any $\mathcal{S}$-theory $D$, the following conditions are satisfied:
	{\em\[
		\begin{array}{cl}
			(\text{a}) &\text{$\alpha\rightarrow\alpha\in D$, for any formula $\alpha$};\\
			(\text{b}) &\text{$\alpha\rightarrow\beta\in D$ and $\beta\rightarrow\gamma\in D$ imply $\alpha\rightarrow\gamma\in D$, for any $\alpha,\beta, \gamma\in\Forms_{\mathcal{L}}$};\\
			(\text{c}) &\text{$\alpha\in D$ and $\alpha\rightarrow\beta\in D$ imply $\beta\in D$, for any $\alpha,\beta\in\Forms_{\mathcal{L}}$};\\
			(\text{d}) &\text{$\beta\in D$ implies $\alpha\rightarrow\beta\in D$, for any $\alpha,\beta\in\Forms_{\mathcal{L}}$};\\
			(\text{e}) &\text{if $\alpha_i\rightarrow\beta_i\in D$, 
				$\beta_i\rightarrow\alpha_i\in D$, $1\le i\le n$, and $F$ is any $n$-ary connective},\\
			&\text{then $F\alpha_1\ldots\alpha_n\rightarrow F\beta_1\ldots\beta_n\in D$, for any formulas $\alpha_1,\ldots,\alpha_n,\beta_1,\ldots,\beta_n$}.\\
		\end{array}
		\]}
\end{defn}

We note that if a logic $\mathcal{S}$ is implicative, then $\bm{T}_{\mathcal{S}}\neq\varnothing$.

Let $\mathcal{S}$ be (at least) an implicative logic and $D$ be an $\mathcal{S}$-theory. We define:
\[
(\alpha,\beta)\in\theta_{D}~\stackrel{\text{df}}{\Longleftrightarrow}~
\text{$\alpha\rightarrow\beta\in D$ and $\beta\rightarrow\alpha\in D$}.
\]
\begin{prop}\label{P:some-classes-unital-logics-1}
	Let $D$ be a theory of an implicative logic $\mathcal{S}$. Then $\theta_{D}$ is a congruence on $\FormAl$
	and $D$ is a congruence class with respect to $\theta_{D}$.
\end{prop}
\begin{proof}
	That $\theta_{D}$ is a congruence follows straightforwardly from (a), (b) and (e) of Definition~\ref{D:weakly-implicative}.
	
	Next, we first note that, by premise, $D\neq\varnothing$. Secondly, if both
	$\alpha$ and $\beta$ belong to $D$, then, according to (d) of Definition~\ref{D:weakly-implicative}, both $\alpha\rightarrow\beta$ and $\beta\rightarrow\alpha $ belong to $D$, that is $(\alpha,\beta)\in\theta_{D}$. And, thirdly, let $(\alpha,\beta)\in\theta_{D}$ and $\alpha\in D$. Then $\beta\in D$. In virtue of Proposition~\ref{P:compatible-congruence}, this implies that $D$ is a congruence class with respect to $\theta_{D}$. 
\end{proof}
\begin{cor}\label{C:some-classes-unital-logics-1}
	If an abstract logic $\mathcal{S}$ is implicative, then $\mathcal{S}$ is unital.
\end{cor}
\begin{proof}
	Let $D$ be any $\mathcal{S}$-theory.
	If $\alpha,\beta\in D$, then, in virtue of Proposition~\ref{P:some-classes-unital-logics-1}, $(\alpha,\beta)\in \theta_{D}$. Then means that $\theta(D)\subseteq\theta_{D}$. Further, by premise, $\theta_{D}$ is compatible with $D$. In virtue of Corollary~\ref{C:compatible-congruence-2},
	$\mathcal{S}$ is unital. 
\end{proof}

In Section~\ref{section:inference-rules}, we defined the calculus $\Cl$ (classical propositional logic) in the language $\Lan_A$, which is defined by the inference rules ax1--ax10 and (b-$iii$) (modus ponens) of Section~\ref{section:inference-rules}.
The calculus in $\Lan_A$, which is  defined by the rules ax1--ax9, ax11 (Section~\ref{section:modus-rules-vs-rules}) and modus ponens, is called the \textit{\textbf{intuitionistic propositional logic}}\index{logic!intuitionistic propositional} and denoted by $\Int$.
Also, it will be convenient to consider another calculus in the language $\Lan_A$, \textsf{P}, which is defined by the inference rules ax1--ax8 and (b-$iii$) (modus ponens) of Section~\ref{section:inference-rules}. We call $\Pos$ \textit{\textbf{positive logic}}.\index{logic!positive}

It must be clear that for any set $X\cup\lbrace\alpha\rbrace\subseteq\Forms_{\mathcal{L}_A}$
\begin{equation}\label{E:two-implications}
	X\vdash_{\textsf{P}}\alpha~\textit{implies both}~X\vdash_{\textsf{Int}}\alpha
	~\textit{and}~X\vdash_{\textsf{Cl}}\alpha.
\end{equation}

Further, it can be easily shown that {\Pos}, $\Cl$ and $\Int$ are implicative. (See Exercise~\ref{section:unital-logics}.\ref{EX:Cl-In-implicative}.)\\

Another class of abstract logics, the class of Fregean logics, stems from 
\textit{Frege's principle of compositionality} (FPC). 

An abstract logic $\mathcal{S}$ has the FPC if for for any set $X\cup\lbrace\gamma(p,\ldots),\alpha,\beta\rbrace\subseteq\Forms_{\mathcal{L}}$,
\[
\text{$X,\alpha\vdash_{\mathcal{S}}\beta$ and $X,\beta\vdash_{\mathcal{S}}\alpha$ imply $X,\gamma(\alpha,\ldots)\vdash_{\mathcal{S}}\gamma(\beta,\ldots)$}.
\]

This leads to the following (equivalent) definition.
\begin{defn}[Frege relation, Fregean logic]\label{D:fregean-relation}\index{logic!Fregean}
	Given a structural logic $\mathcal{S}$ and an $\mathcal{S}$-theory $D$, the relation on $\FormAl$ defined by the equivalence
	{\em\[
		(\alpha,\beta)\in\Lambda D~\stackrel{\text{df}}{\Longleftrightarrow}~
		\text{$D,\alpha\vdash_{\mathcal{S}}\beta$ and $D,\beta\vdash_{\mathcal{S}}\alpha$}
		\]}
	is called a \textbf{Frege relation relative to} $D$. A structural logic
	is \textbf{Fregean} if for every $\mathcal{S}$-theory $D$, $\Lambda D$ is a congruence on $\FormAl$.
\end{defn}

We leave the reader to prove that for any structural logic $\mathcal{S}$, the properties of satisfying the FPC and being Fregean are equivalent. (See Exercise~\ref{section:lindenbaum-algebra}.\ref{EX:Fregean-equivalence}.)
\begin{prop}
	If $\mathcal{S}$ is a Fregean logic with $\ThmS\neq\varnothing$, then $\mathcal{S}$ is unital.
\end{prop}
\begin{proof}
	Let $D$ be an arbitrary $\mathcal{S}$-theory. We observe that
	$\theta(D)\subseteq\Lambda D$. Also, $\Lambda D$ is compatible with $D$.
	Applying Corollary~\ref{C:compatible-congruence-2}, we obtain that $\mathcal{S}$ is unital.
\end{proof}

\paragraph{Exercises~\ref{section:unital-logics}}
\begin{enumerate}
	\item Prove that the $\booleTwo$- consequence and $\lukasThree$-consequence are unital.
	\item \label{EX:Cl-In-implicative}Prove that the abstract logics $\Pos$, $\Cl$ and $\Int$ are implicative.
	\item Prove that the abstract logics $\Pos$, $\Cl$ and $\Int$ are Fregean.
	\item\label{EX:Fregean-equivalence} Prove that a structural logic $\mathcal{S}$ is Fregean if, and only if, it satisfies the FPC.
\end{enumerate}

\section{Lindenbaum-Tarski algebras}\label{section:lindenbaum-algebra}

\subsection{Definition and properties}
The concept of unital abstract logic in turn induces the following important notion.
\begin{defn}[Lindenbaum-Tarski algebras]\label{D:LT-algebra}\index{algebra!Lindenbaum-Tarski}
	Let $\mathcal{S}$ be a unital logic and $D$ be an $\mathcal{S}$-theory. The algebra
	{\em$\LT[D]:=\langle\FormAl\slash\theta(D),\one_D\rangle$}, where $\one_{D}$ is $D$ understood as a congruence class with respect to $\theta(D)$, $($or very often only the first component $\FormAl\slash\theta(D)$$)$ is called a \textbf{Lindenbaum-Tarski algebra} $($of $\mathcal{S}$$)$ \textbf{relative to} $D$. The algebra {\em$\LT:=\langle\FormAl\slash\theta(\ThmS),\one_{\bm{T}_{\mathcal{S}}}\rangle$} $($or the first component $\FormAl\slash\theta(\bm{T}_{\mathcal{S}})$$)$ is called simply a \textbf{Lindenbaum-Tarski algebra} of $\mathcal{S}$.
\end{defn}

\begin{prop}\label{P:LT-algebras-1}
	Let $\mathcal{S}$ be a unital logic and and let $D$ and $D^{\prime}$ be $\mathcal{S}$-theories. If $D\subseteq D^{\prime}$, then the algebra {\em$\LT[D^{\prime}]$} is a homomorphic image of {\em$\LT[D]$}. Hence, for any $\mathcal{S}$-theory $D$, the algebra {\em$\LT[D]$} is a homomorphic image of {\em$\LT$}.
\end{prop}
\begin{proof}
	We observe that, given a unital logic $\mathcal{S}$ and any $\mathcal{S}$-theories $D$ and $D^{\prime}$, if $D^{\prime}\subseteq D$, then
	$\theta(D)\subseteq\theta(D^{\prime})$; see~\eqref{E:preliminaries-algebra-1}. Further, in virtue of Proposition~\ref{P:congruence-thm-auxiliary}, this implies that $\LT[D]$ is a homomorphic image of $\LT[D^{\prime}]$. Thus, for any $\mathcal{S}$-theory $D$, the algebra $\LT[D]$ is a homomorphic image of $\LT$.	
\end{proof}

\begin{lem}\label{L:congruence-substitution}
	Let $\theta$ be a congruence on $\FormAl$. Then for any valuation $v$ in \mbox{$\FormAl\slash\theta$},  there is an $\Lan$-substitution $\sigma$ such that for any formula $\alpha$,  \mbox{$v[\alpha]=\sigma(\alpha)\slash\theta$}.
\end{lem}
\begin{proof}
	Let $v$ be any valuation in $\FormAl\slash\theta$; that is, $v$ is a homomorphism of
	$\FormAl$ to $\FormAl\slash\theta$. Suppose $v[p]=\beta_p\slash\theta$, where $p\in\Var_{\mathcal{L}}$. Then we define a substitution as follows:
	\[
	\sigma:p\mapsto\beta_p,
	\]
	for all $p\in\Var_{\mathcal{L}}$, with any selection of $\beta_p$ from the class $\beta_p\slash\theta$. Thus, we have: $v[p]=\sigma(p)\slash\theta$, for all $p\in\Var_{\mathcal{L}}$.
	
	Now let us take any formula $\alpha$ and assume that 
	$p_1,\ldots,p_n$ are all variables  that occur in $\alpha$. Then we have:
	\[
	v[\alpha]=\alpha(v[p_1],\ldots,v[p_n])=\alpha(\sigma(p_1)\slash\theta,\ldots,\sigma(p_n)\slash\theta)=\sigma(\alpha)\slash\theta.
	\]
\end{proof}

\begin{prop}\label{P:S-algebraic-models}
	Let $\mathcal{S}$ be a unital logic. Then each {\em$\LT[D]$} is an algebraic $\mathcal{S}$-model.
\end{prop}
\begin{proof}
	Suppose that $X\vdash_{\mathcal{S}}\alpha$. Let us take any $\LT[D]$ and consider a valuation $v$ in it. In virtue of Lemma~\ref{L:congruence-substitution}, there is a substitution $\sigma$ such that for any formula $\beta$, $v[\beta]=\sigma(\beta)\slash\theta(D)$. Further, we note that, since $\mathcal{S}$ is structural, $\sigma(X)\vdash_{\mathcal{S}}\sigma(\alpha)$. 
	
	Now, suppose that for any $\gamma\in X$, $v[\gamma]=\one_{D}$. Then, since $\mathcal{S}$ is unital, for any $\gamma\in X$, $\sigma(\gamma)\in D$. This implies that $\sigma(\alpha)\in D$, that is $v[\alpha]=\one_{D}$.
\end{proof}
\begin{cor}\label{C:quasiidentity-sound}
	Let $\mathcal{S}$ be a unital logic. Then for any set \mbox{{\em$X\cup\lbrace\alpha\rbrace\Subset\Forms_{\mathcal{L}}$}},
	and for any $\mathcal{S}$-theory $D$, 
	{\em\[
		X\vdash_{\mathcal{S}}\alpha~\Longrightarrow~\LT[D]\models X\approx\one_{D}\Rarrow\alpha\approx\one_{D}.
		\]}
\end{cor}
\begin{proof}
	We apply successively Proposition~\ref{P:S-algebraic-models} and~\eqref{E:model=algebra-validity-2}.
\end{proof}

\begin{prop}\label{P:preLindenbaum-algebra}
	Let $\mathcal{S}$ be a unital logic and let {\em$D\!:=\ConS{X}$}. Then for any formula $\alpha$,
	{\em\[
		\models_{\LT[D]}\alpha\Longleftrightarrow~\models_{\LinS[X]}\alpha.
		\]}
	Hence, {\em$\LT$} is weakly adequate for $\mathcal{S}$.
\end{prop}
\begin{proof}
	First, assume that $\models_{\LinS[X]}\alpha$; that is for any substitution $\sigma$, $\sigma(\alpha)\in D$.
	Let $v$ be a valuation in $\LT[D]$. By Lemma~\ref{L:congruence-substitution}, there is a substitution $\sigma_v$ such that $v[\alpha]=\sigma_v(\alpha)\slash\theta(D)$. By premise, this implies that $v[\alpha]=\one_{D}$.
	
	Conversely, suppose for some substitution $\sigma$, $\sigma(\alpha)\notin D$. We define a valuation in $\LT[D]$ as follows: $v[p]:=\sigma(p)\slash\theta(D)$, for any $p\in\Var_{\mathcal{L}}$. Then $v[\alpha]=\sigma(\alpha)\slash\theta(D)\neq\one_D$.
\end{proof}

\begin{cor}
	Let $D$ be an $\mathcal{S}$-theory of a unital logic $\mathcal{S}$. Then for any
	$\mathcal{S}$-thesis $\alpha$, the algebra {\em$\LT[D]$} validates $\alpha$.
\end{cor}
\begin{proof}
	Having noticed that $\LT[D]$ is a homomorphic image of $\LT$, we apply the last part of Proposition~\ref{P:preLindenbaum-algebra}.
\end{proof}

We recall that $\Sigma_{\mathcal{S}}$ denotes the set of all $\mathcal{S}$-theories.
\begin{cor}\label{C:LT-determination}
	Any unital logic $\mathcal{S}$ is determined by the class \mbox{\em$\set{\LT[D]}{D\in\Sigma_{\mathcal{S}}}$}. Hence, any unital abstract logic is assertional.
\end{cor}
\begin{proof}
	Indeed, if $X\vdash_{\mathcal{S}}\alpha$, then, in virtue of Proposition~\ref{P:S-algebraic-models}, $X\models_{\LT[D]}\alpha$, for each $D\in\Sigma_{\mathcal{S}}$.
	
	On the other hand, if $X\not\vdash_{\mathcal{S}}\alpha$, then $\not\models_{\LinS[X]}\alpha$ and hence, according to Proposition~\ref{P:preLindenbaum-algebra}, $\not\models_{\LT[D]}\alpha$, where $D=\ConS{X}$.
\end{proof}

\begin{cor}\label{C:preLindenbaum-algebra}
	Let $\mathcal{S}$ be a unital logic.  If a set {\em$X\subseteq
		\Forms_{\mathcal{L}}$} is closed under arbitrary substitutions, then for any formula $\alpha$,
	{\em\[
		X\vdash_{\mathcal{S}}\alpha~\Longleftrightarrow~\models_{\LT[D]}
		\alpha,
		\]}
	where {\em$D=\ConS{X}$}.
\end{cor}
\begin{proof}
	We have:
	\[
	\begin{array}{rl}
		X\vdash_{\mathcal{S}}\alpha \Longleftrightarrow \!\!\!&\models_{\LinS[X]}\alpha~~[\text{Lindenbaum's theorem (Proposition~\ref{P:lindenbaum-theorem})}]\\
		\Longleftrightarrow \!\!\!&\models_{\LT[D]}\alpha.~~[\text{Proposition~\ref{P:preLindenbaum-algebra}}]\\
	\end{array}
	\]
\end{proof}

Now we define a special valuation in $\LT[D]$:
\[
\bm{v}_{D}: p\mapsto p\slash\theta(D),\tag{\textit{Lindenbaum valuation}}
\]
for any $p\in\Var_{\mathcal{L}}$. 
\begin{prop}\label{P:valuations-in-LT}
	Let $\mathcal{S}$ be a unital logic and {\em$D=\ConS{X}$}.
	Then for any formula $\alpha$, the following conditions are equivalent$\,:$
	{\em\[
		\begin{array}{cl}
			(\text{a}) & X\vdash_{\mathcal{S}}\alpha;\\
			(\text{b}) &X\models_{\LT[D]}\alpha;\\
			(\text{c}) &\bm{v}_D[\alpha]=\one_{D}.
		\end{array}
		\]}
	In particular, the following conditions are equivalent$\,:$
	{\em\[
		\begin{array}{cl}
			(\text{a}^{\ast}) & \alpha\in\ThmS;\\
			(\text{b}^{\ast}) &v[\alpha]=\one_{\bm{T}_{\mathcal{S}}},~\textit{for any valuation $v$ in}~ \LT;\\
			(\text{c}^{\ast}) &\bm{v}_{\bm{T}_{\mathcal{S}}}[\alpha]=\one_{\bm{T}_{\mathcal{S}}}.
		\end{array}
		\]}
\end{prop}
\begin{proof}
	The implication (\text{a})$\Rightarrow$(\text{b}) follows from Corollary~\ref{C:quasiidentity-sound}.
	
	Assume (b). Let us take any $\beta\in X$. We have: $\bm{v}_{D}[\beta]=\beta\slash\theta(D)=\one_{D}$. Hence (c). 
	
	Now assume (c), that is $\bm{v}_D[\alpha]=\one_{D}$. This means that $\alpha\slash\theta(D)=\one_{D}$, that is $\alpha\in D$. Hence (a).
\end{proof}

We aim to show that the Lindenbaum-Tarski algebras are good as separating tools (Section~\ref{section:separating-means}); or, better to say, they are useful for creating separating tools. With this in mind, given a unital, and thus assertional (Corollary~\ref{C:LT-determination}), logic $\mathcal{S}$ in a language $\Lan$, we define the following classes of unital expansions. These definitions are based on the conclusion of Corollary~\ref{C:LT-determination} that all unital logics are assertional.
\begin{itemize}
	\item $\PS$ is the class of all algebraic $\mathcal{S}$-models; according to Proposition~\ref{P:S-algebraic-models=q-variety}, $\PS$ is a quasi-variety.
	\item $\QS$ is the quasi-variety generated by the set $\set{\LT[D]}{D\in\Sigma_{\mathcal{S}}}$. According to Proposition~\ref{P:S-algebraic-models}, 
	\begin{equation}\label{E:LT-algebras-unital-1}
		\QS\subseteq\PS.
	\end{equation}
	Also, if $\mathcal{S}$ is finitary, then, in virtue of Corollary~\ref{C:quasiidentity-sound}, $\QS$ is determined by the quasiidentities $X\approx\one\Rarrow\alpha\approx\one$, for any $X\cup\lbrace\alpha\rbrace\Subset\Forms_{\mathcal{L}}$ such that $X\vdash_{\mathcal{S}}\alpha$.
	\item $\TS$ is the class of all unital expansions $\langle\alg{A},\one\rangle$ such that
	$\langle\alg{A},\one\rangle\models\alpha\approx\one$, for any $\alpha\in\ThmS$. In view of~\eqref{E:model=algebra-validity-1}, the class $\TS$ is a variety. This, the second part of Proposition~\ref{P:valuations-in-LT} and Proposition~\ref{P:LT-algebras-1} imply that
	\begin{equation}\label{E:LT-algebras-unital-2}
		\set{\LT[D]}{D\in\Sigma_{\mathcal{S}}}\subseteq\TS,
	\end{equation}
	which in turn implies that
	\begin{equation}\label{E:LT-algebras-unital-3}
		\QS\subseteq\TS.
	\end{equation} 	
\end{itemize}

We conclude this subsection with the following.
\begin{prop}\label{P:LT-unital-PS-QS-equivalence}
	Let $\mathcal{S}$ be a unital logic. Then the following conditions are equivalent$\,:$
	{\em\[
		\begin{array}{cl}
			(\text{a})	&X\vdash_{\mathcal{S}}\alpha;\\
			(\text{b})	&\text{$X\models_{\langle\textbf{A},\one\rangle}\alpha$, for any 
				$\langle\textbf{A},\one\rangle\in\PS$};\\
			(\text{c})	&\text{$X\models_{\langle\textbf{A},\one\rangle}\alpha$, for any 
				$\langle\textbf{A},\one\rangle\in\QS$};\\
			(\text{d})	&\text{$X\models_{\LT[D]}\alpha$, for any $D\in\Sigma_{\mathcal{S}}$}.
		\end{array}
		\]}
\end{prop}
\begin{proof}
	The implication (a)$\Rightarrow$(b) is due to definition of $\PS$; (b)$\Rightarrow$(c) is due to inclusion $\QS\subseteq\PS$; (c)$\Rightarrow$(d) is due to definition of $\QS$; and
	(d)$\Rightarrow$(a) is due to Corollary~\ref{C:LT-determination}.
\end{proof}

Thus, we see that if $X\not\vdash_{\mathcal{S}}\alpha$, an algebra $\langle\alg{A},\one\rangle$ separating $X$ from $\alpha$ 
can be found among the algebras of the class 
$\QS=\Su\Pred\lbrace\LT[D]\rbrace_{D\in\Sigma_{\mathcal{S}}}$; cf.~\cite{malcev1973}, {\S} 11, theorem 4. As we will see in the next subsection, under a restriction on $\mathcal{S}$, separating algebras (in the above sense) can be reached as homomorphic images of $\LT$. In the proof of this below, the class $\TS$ plays a key role. 
\begin{cor}\label{C:unital-logic-has-single-matrix}
	Let $\mathcal{S}$ be a unital logic. Then $\mathcal{S}$ has an adequate matrix.
\end{cor}
\begin{proof}
	First we assume that $\mathcal{S}$ is nontrivial.
	
	Basing on the equivalence (a)$\Leftrightarrow$(c) of Proposition~\ref{P:LT-unital-PS-QS-equivalence}, for any pair $(X,\alpha)$ with $X\not\vdash_{\mathcal{S}}\alpha$, we define: 
	\begin{center}
		$\langle\alg{A},\one\rangle_{(X,\alpha)}$ is such an algebra $\langle\alg{A},\one\rangle\in\QS$ that $X\not\models_{\langle\alg{A},\one\rangle}\alpha$.
	\end{center}
	
	Further, we define:
	\[
	\mathcal{U}_{\mathcal{S}}:=\set{\langle\alg{A},\one\rangle_{(X,\alpha)}}{(X,\alpha)\in\mathcal{P}(\FormsL)\times\FormsL~\text{and}~X\not\vdash_{\mathcal{S}}\alpha}.\footnote{Axiom of Choice is needed for this definition.}
	\]
	
	We observe that $\mathcal{U}_{\mathcal{S}}\neq\varnothing$ and $\mathcal{U}_{\mathcal{S}}\subseteq\QS$.
	Next, we denote:
	\[
	\langle\overline{\alg{A}},\overline{\one}\rangle:=\prod\mathcal{U}_{\mathcal{S}}.
	\]
	
	It is obvious that $\langle\overline{\alg{A}},\overline{\one}\rangle\in\QS$. Further, $\langle\overline{\alg{A}},\overline{\one}\rangle$ is an adequate matrix for $\mathcal{S}$. Indeed, if $X\vdash_{\mathcal{S}}\alpha$, then, since $\langle\overline{\alg{A}},\overline{\one}\rangle\in\QS$, in virtue of the equivalence (a)$\Leftrightarrow$(c) of Proposition~\ref{P:LT-unital-PS-QS-equivalence}, $X\models_{\langle\overline{\alg{A}},\overline{\one}\rangle}\alpha$. On the other hand, if $X\not\vdash_{\mathcal{S}}\alpha$, then for some $\langle\alg{A},\one\rangle\in\mathcal{U}_{\mathcal{S}}$, $X\not\models_{\langle\alg{A},\one\rangle}\alpha$ and, hence,
	$X\not\models_{\langle\overline{\alg{A}},\overline{\one}\rangle}\alpha$.
	
	If $\mathcal{S}$ is trivial, a degenerated algebra in $\QS$ is an adequate matrix.
\end{proof}

\begin{cor}\label{C:Cl-and-Int-uniform}
	Both {\em\Cl} and {\em\Int} have an adequate matrix. Therefore, both are uniform, couniform and have cancellation property.
\end{cor}

\subsection{Implicational unital logics}
We aim to show that, under above restriction on $\mathcal{S}$, $\LT$ is a free algebra over $\TS$.

\begin{defn}[implicational abstract logic]\label{D:implicational}\index{logic!implicational}
	Let $\mathcal{S}$ be a unital abstract logic  in a language $\Lan$. Assume that 
	there is an $\Lan$-formula of two sentential variables, $\im(p,q)$, such that for any $\Lan$-formulas $\alpha$ and $\beta$, any algebra {\em$\langle\alg{A},\one\rangle\in\TS$} and any valuation $v$ in $\alg{A}$,
	\[
	v[\alpha]=v[\beta]~\Longleftrightarrow~v[\im(\alpha,\beta)]=\one~\text{and}~
	v[\im(\beta,\alpha)]=\one.
	\]
	Then we call $\mathcal{S}$ \textbf{implicational} with respect to $\im(p,q)$.
\end{defn}

The abstract logics $\Pos$, $\Cl$ and $\Int$ are examples of implicational logics w.r. $p\rightarrow q$. (Exercise~\ref{section:lindenbaum-algebra}.\ref{EX:Cl-Int-implicational})
This means that if $\aLog$ is one of these abstract logics, then
\[
(\alpha,\beta)\in\theta(\ThmS)~\Longleftrightarrow~\vdash_{\mathcal{S}} \alpha\leftrightarrow\beta.
\]

The last observations leads to the following proposition.
\begin{prop}[replacement property]\label{P:replacement}
	Let $\aLog$ be one of the abstract logics {\em$\Pos$}, {\em$\Cl$} or {\em$\Int$}. Assume that a formula $\beta$ is a subformula of a formula $\alpha$ so that a designated occurrence of $\beta$ in $\alpha$ is denoted by $\alpha[\dots\beta\dots]$. Then
	\[
	\vdash_{\mathcal{S}}\beta\leftrightarrow\gamma~\Longrightarrow~\vdash_{\mathcal{S}}\alpha[\dots\beta\dots]\leftrightarrow\alpha[\dots\gamma\dots],
	\]
	where $\alpha[\dots\gamma\dots]$ is the formula obtained from $\alpha[\dots\beta\dots]$ by replacement of the designated occurrence of $\beta$ with $\gamma$.
\end{prop}
\begin{proof}
	To prove we notice that if, by premise, $\beta\slash\theta(\ThmS)=\gamma\slash\theta(\ThmS)$, then
	$\alpha[\dots\beta\slash\theta(\ThmS)\dots]=\alpha[\dots\gamma\slash\theta(\ThmS)\dots]$ and hence $\alpha[\dots\beta\dots]\slash\theta(\ThmS)=
	\alpha[\dots\gamma\dots]\slash\theta(\ThmS)$.	
\end{proof}

Below we will be considering implicational logics which are also unital. Since for such a logic $\mathcal{S}$, $\ThmS\neq\varnothing$ (Proposition~\ref{P:lindenbaum-connection}), we select an arbitrary $\Lan$-formula $\bigstar\in\ThmS$.

In virtue of Proposition~\ref{P:valuations-in-LT}, for any valuation $v$ in $\LT$,
\begin{equation}\label{E:T_S-star-thesis}
	v(\bigstar)=\one_{\bm{T}_{\mathcal{S}}}.
\end{equation}
(Exercise~\ref{section:lindenbaum-algebra}.\ref{EX:T_S-star-thesis})

We recall the signature of any unital expansion $\langle\alg{A},\one\rangle$ is of type $\Lan$ expanded by an additional constant symbol $\one$ which is interpreted by the 0-ary operation $\one$ in $\langle\alg{A},\one\rangle$ (denoted by $\one_{\bm{T}_{\mathcal{S}}}$ in $\LT$). We denote this extended type by $\Lan^{\star}$.

Let $\mathbf{t}$ be a term of type $\Lan^{\star}$. The \textbf{\textit{conversion}}
of $\mathbf{t}$ is the $\Lan$-formula (or term of type $\Lan$), denoted by $\mathbf{t}^{\star}$, which is obtained from $\mathbf{t}$ by simultaneous replacement of all occurrences of $\one$, if any, in $\mathbf{t}$ by $\bigstar$. 

Consequently, any valuation $v$ in an algebra $\alg{A}$ of type $\Lan$ is extended to
the valuation $v^{\star}$ in $\langle\alg{A},\one\rangle$ such that $v^{\star}(p)=v(p)$, for any $p\in\Var_{\mathcal{L}}$, and $v^{\star}(\one)=\one$, where the first occurrence of `$\one$' is a 0-ary constant of $\Lan^{\star}$, and the second is the 0-ary operation $\one$ of $\langle\alg{A},\one\rangle$. We refer to $v^{\star}$ as a $\star$-\textit{\textbf{extension}}\index{$\star$-extension} of $v$. We note that the valuation $v$ can be restored from $v^{\star}$ in the sense that if $w$ is any valuation in $\langle\alg{A},\one\rangle$, there is a unique valuation $v$ in $\alg{A}$ such that $v^{\star}=w$. (Exercise~\ref{section:lindenbaum-algebra}.\ref{EX:v-restored})

Using~\eqref{E:T_S-star-thesis}, we observe that for any $\Lan^{\star}$-terms $\mathbf{t}$ and $\mathbf{r}$ and any valuation $v^{\star}$ in $\LT$,
\begin{equation}\label{E:LT-term-valuation}
	v^{\star}(\mathbf{t})=v^{\star}(\mathbf{r})~\Longleftrightarrow~v(\mathbf{t}^{\star})=v(\mathbf{r}^{\star}).
\end{equation} 
(Exercise~\ref{section:lindenbaum-algebra}.\ref{EX:LT-term-valuation})

For any terms $\mathbf{t}$ and $\mathbf{r}$ of type $\Lan^{\star}$, we denote:
\begin{center}
	$\TS\models\mathbf{t}\approx\mathbf{r}~\stackrel{\text{df}}{\Longleftrightarrow}~$
	for any $\left\langle \alg{A},\one\right\rangle \in\TS$, $\models_{\langle\textbf{A},\one\rangle}\mathbf{t}\approx\mathbf{r}$.
\end{center}

The equivalence~\eqref{E:LT-term-valuation} can be generalized as follows.
\begin{lem}\label{L:LT=free-algebra}
	Let $\langle \alg{A},\one\rangle\in\TS$. Then for any term $\mathbf{t}$
	of type $\Lan^{\star}$ and any valuation $v^{\star}$ in $\langle \alg{A},\one\rangle$,
	$v^{\star}(\mathbf{t})=v(\mathbf{t}^{\star})$.
\end{lem}
\noindent\textit{Proof}~can be carried out by induction on the number of functional symbols of the language $\Lan^{\star}$. We leave it to the reader. (Exercise~\ref{section:lindenbaum-algebra}.\ref{EX:LT=free-algebra})\\

\begin{prop}\label{P:LT=free-algebra}
	Let $\mathcal{S}$ be nontrivial unital logic which is implicational with respect to
	$\im(p,q)$. Then for terms $\mathbf{t}$ and $\mathbf{r}$ of type $\Lan^{\star}$,
	\[
	\TS\models\mathbf{t}\approx\mathbf{r}~\Longleftrightarrow~\bm{v}_{0}^{\star}(\mathbf{t})=\bm{v}_{0}^{\star}(\mathbf{r}),
	\]
	where $\bm{v}_{0}^{\star}$ is the $\star$-extension of $\bm{v}_{0}$. Hence {\em$\LT$} is a free algebra of rank $\card{\Lan}$ over the variety $\TS$ and also over the quasi-variety $\QS$.
\end{prop}
\begin{proof}
	We successively obtain:
	\[
	\begin{array}{rcl}
		\TS\models\mathbf{t}\approx\mathbf{r}
		&\Longleftrightarrow &\TS\models\mathbf{t}^{\star}\approx\mathbf{r}^{\star}\quad[\text{Lemma~\ref{L:LT=free-algebra}}]\\
		&\Longleftrightarrow &\TS\models\im(\mathbf{t}^{\star},\mathbf{r}^{\star})\approx\one_{\bm{T}_{\mathcal{S}}}~\text{and}~\TS\models\im(\mathbf{r}^{\star},\mathbf{t}^{\star})\approx\one_{\bm{T}_{\mathcal{S}}}\\
		&&[\text{for $\mathcal{S}$ is implicational with respect to $\im(p,q)$}]\\
		&\Longleftrightarrow &\im(\mathbf{t}^{\star},\mathbf{r}^{\star})\in\ThmS~
		\text{and}~\im(\mathbf{r}^{\star},\mathbf{t}^{\star})\in\ThmS\\
		&&[\text{since $\LT\in\TS$ and Proposition~\ref{P:preLindenbaum-algebra}, part 2}]\\
		&\Longleftrightarrow &\bm{v}_{0}(\im(\mathbf{t}^{\star},\mathbf{r}^{\star}))=\one_{\bm{T}_{\mathcal{S}}}~
		\text{and}~\bm{v}_{0}(im(\mathbf{r}^{\star},\mathbf{t}^{\star}))=\one_{\bm{T}_{\mathcal{S}}}\\
		&&[\text{Proposition~\ref{P:preLindenbaum-algebra}, part 2}]\\
		&\Longleftrightarrow &\bm{v}_{0}(\mathbf{t}^{\star})=\bm{v}_{0}(\mathbf{r}^{\star})
		\quad[\text{for $\mathcal{S}$ is implicational with respect to $\im(p,q)$}]\\
		&\Longleftrightarrow &\bm{v}_{0}^{\star}(\mathbf{t}^{\star})=\bm{v}_{0}^{\star}(\mathbf{r}^{\star})
		\quad[\text{in virtue of~\eqref{E:LT-term-valuation}}].
	\end{array}
	\]
	
	According to~\cite{malcev1973}, {\S} 12.2, theorem 1, this implies that $\LT$ is a free algebra over the class $\TS$. In virtue of Proposition~\ref{P:motivating}, the rank of $\LT$ equals $\card{\Lan}$.
	
	Since the quasi-variety $\QS$ contains all algebras $\LT[D]$ and because of Corollary~\ref{C:LT-determination}, $\mathcal{S}$ is determined by the algebraic expansions of
	$\QS$. It is obvious that $\LT$ is a free algebra over $\QS$.
\end{proof}

\subsection{Examples of Lindenbaum-Tarski algebras}\label{section:some-lindenbaum-algebras}
Let $\mathcal{S}$ be a nontrivial unital logic in a language $\Lan$ with $\card{\Var_{\mathcal{L}}}=\kappa$. Then, according to~\eqref{E:motivating},
$\card{|\LT|}\ge\kappa$. In particular, $\card{|\LTCl|}=\card{|\LTInt|}=\aleph_0$.
It must be clear that $\set{p\slash\theta(\ThmS)}{p\in\Var_{\mathcal{L}}}$ is the set of free generators of $\LT$, if the latter can be considered as a free algebra of $\TS$.
For some questions like decision problem for the consequence relation or theoremhood of $\mathcal{S}$, it suffices to consider the formulas built of a limited set of variables. Thus, if $\Var\subseteq\Var_{\mathcal{L}}$ with $\card{\Var}=\kappa$ and our consideration of $\mathcal{S}$ is limited to the formulas built from $\Var$, we denote by $\LT(\kappa)$ a subalgebra of $\LT$, generated by the set $\set{p\slash\theta(\ThmS)}{p\in\Var}$. The notation is justified by the fact if $\Var^{\prime}\subseteq\Var_{\mathcal{L}}$ and
$\card{\Var^{\prime}}=\card{\Var}$, then the subalgebras of $\LT$, generated by
$\set{p\slash\theta(\ThmS)}{p\in\Var}$ and, respectively, by 
$\set{p\slash\theta(\ThmS)}{p\in\Var^{\prime}}$ are isomorphic; cf.~\cite{gra79}, {\S}24, corollary 2. Since each such $\LT(\kappa)$ is a free algebra in $\QS$ (or even in $\TS$), it is called the free algebra of rank $\kappa$. We note that, since $\LT$ is always has at least one 0-ary operation, $\one$, even it makes sense to talk about a free algebra of rank 0.

We consider the free algebras $\LTCl(0)$, $\LTInt(0)$, $\LTCl(1)$
$\LTCl(2)$ and $\LTInt(1)$, in order. The first two algebras should be understood as subalgebras of $\LTCl$ and $\LTInt$, generated by their unit, $\one$. The resulting algebras are isomorphic (see Section~\ref{section:boolean-algebra} and Section~\ref{section:heyting-algebra}) and can be depicted by the diagram:

\begin{figure}[!ht]	
	\[
	\ctdiagram{
		\ctnohead
		\ctinnermid
		\ctel 0,0,0,20:{}
		\ctv 0,0:{\bullet}
		\ctv 0,20:{\color{red}\bullet}
		\ctv 0,27:{\mathbf{1}}
	}
	\]
	\caption{Lindenbaum-Tarski algebra $\LTCl(0)$, or equivalently $\LTInt(0)$}
\end{figure}

Before discussing the structure of the other just mentioned Lindenbaum-Tarski algebras, we first consider the varieties
$\mathbb{T}_{\textsf{Cl}}$ and $\mathbb{T}_{\textsf{Int}}$. Namely, we aim to show that these varieties coincide with the varieties of Boolean and Heyting algebras, respectively.

Since the classes $\mathbb{T}_{\textsf{Cl}}$ and $\mathbb{T}_{\textsf{Int}}$ consist of the algebras validating the theses of $\Cl$ and $\Int$, respectively, (defined in Section~\ref{section:subclasses-unital}) we take a look at these calculi more closely. 

Technically, it is convenient to have a collective name for the calculi $\Pos, \Int$ and $\Cl$. Refereeing to any of these calculi, we sometimes use the latter $\C$.
\begin{lem}\label{L:rules-sound-wrt-P-Int-Cl}
	In addition to modus ponens, the inference rules {\em(\text{a}--$i$)--(\text{a}--$iii$)} and 
	{\em(\text{b}--$i$)--(\text{b}--$ii$)} of Section~\ref{section:inference-rules} are sound with respect to {\em\Pos}, {\em\Int} and {\em\Cl}.
\end{lem}
\begin{proof}
	To prove the soundness of the rules (a--$i$)--(a--$iii$), we have to show that for any $\Lan_A$-formulas $\alpha$ and $\beta$, the following hold:
	\[
	\alpha,\beta\vdash_{\textsf{C}}\alpha\wedge\beta,\quad\alpha\vdash_{\textsf{C}}
	\alpha\vee\beta,\quad\beta\vdash_{\textsf{C}}\alpha\vee\beta,\quad\alpha\wedge\beta
	\vdash_{\textsf{C}}\alpha~\text{and}~\alpha\wedge\beta
	\vdash_{\textsf{C}}\beta;
	\]
	and then, for the soundness of the rules (b--$i$)--(b--$ii$):
	\[
	\alpha\land\beta\vdash_{\textsf{C}}\alpha~~\text{and}~~\alpha\land\beta\vdash_{\textsf{C}}\beta.
	\]
	
	We leave this task to the reader. (Exercise~\ref{section:lindenbaum-algebra}.\ref{EX:lemma-rules-sound-wrt-P-Int-Cl})
\end{proof}

The next lemma is known as the deduction theorem.
\begin{lem}\label{L:deduction-theorem}
	The hyperrule {\em(\text{c}--$i$)} of Section~\ref{section:inference-rules} is sound with respect to {\em\Pos}, {\em\Int} and {\em\Cl}.
\end{lem}
\begin{proof}
	First of all, we notice that for any $\Lan_A$-formulas $\alpha$ and $\beta$,
	\[
	\begin{array}{cll}
		1. &\vdash_{\textsf{C}}\alpha\rightarrow(\beta\rightarrow\alpha) &[\text{an instance of ax1}]\\
		2. &\vdash_{\textsf{C}}\alpha\rightarrow((\beta\rightarrow\alpha)\rightarrow\alpha) &[\text{an instance of ax1}]\\
		3. &\vdash_{\textsf{C}}(\alpha\rightarrow(\beta\rightarrow\alpha))\rightarrow
		((\alpha\rightarrow((\beta\rightarrow\alpha)\rightarrow\alpha))
		\rightarrow(\alpha\rightarrow\alpha)) &[\text{an instance of ax2}]\\
		4. &\vdash_{\textsf{C}}(\alpha\rightarrow((\beta\rightarrow\alpha)\rightarrow\alpha))
		\rightarrow(\alpha\rightarrow\alpha) &[\text{from 1 and 3 by (b--$iii$)}]\\
		5. &\vdash_{\textsf{C}}\alpha\rightarrow\alpha. &[\text{from 2 and 4 by (b--$iii$)}]
	\end{array}
	\]
	
	Thus for any formula $\alpha$,
	\[
	\vdash_{\textsf{C}}\alpha\rightarrow\alpha. \tag{$\ast$}
	\]
	Now assume that $X,\alpha\vdash_{\textsf{C}}\beta$. That is, there is a $\C$-derivation of $\beta$ from $X,\alpha$:
	\[
	\alpha_1,\ldots,\beta.\tag{$\ast\ast$}
	\]
	Let us form a finite sequence:
	\[
	\alpha\rightarrow\alpha_1,\ldots,\alpha\rightarrow\beta
	\]
	and denote:
	$\gamma_i:=\alpha\rightarrow\alpha_i$, where $1\le i\le n$; so that $\gamma_n:=\alpha\rightarrow\beta$. By induction on the number $n$ of the formulas of $(\ast\ast)$, we prove that $\vdash_{\textsf{C}}\gamma_n$.
	
	If $n=1$, then either $\beta=\alpha$ or $\beta\in X$. If the former is the case, then we apply ($\ast$). If $\beta\in X$, then we successively obtain:
	$\vdash_{\textsf{C}}\beta\rightarrow(\alpha\rightarrow\beta)$ and $X\vdash_{\textsf{C}}\alpha\rightarrow\beta$.
	
	Next, assume that $\beta$ is obtained from $\alpha_i$ and $\alpha_j=\alpha_i\rightarrow\beta$. By induction hypothesis,
	$\vdash_{\textsf{C}}\alpha\rightarrow\alpha_i$ and 
	$\vdash_{\textsf{C}}\alpha\rightarrow(\alpha_i\rightarrow\beta)$. Using ax2 and
	(b--$iii$) twice, we receive that $\vdash_{\textsf{C}}\alpha\rightarrow\beta$.
\end{proof}

\begin{lem}\label{L:hyperrules-c_ii-c_iii}
	The hyperrules {\em(c--$ii$)} and {\em(c--$iii$)} of Section~\ref{section:inference-rules} are sound with respect to {\em\Int} and {\em\Cl}.
	Hence the rule {\em(a--$iv$)} is sound with respect to {\em\Int} and {\em\Cl}.
\end{lem}
\begin{proof}
	The proof of the first part is straightforward and is left to the reader. (Exercise~\ref{section:lindenbaum-algebra}.\ref{EX:lemma-hyperrules-c_ii-c_iii})
	
	The proof of the second part is as follows. Let $\mathcal{S}$ be either {\Int} or {\Cl}. Since $\alpha,\neg\alpha\vdash_{\mathcal{S}}\alpha$ and 
	$\alpha,\neg\alpha\vdash_{\mathcal{S}}\neg\alpha$, by (c-$ii$), we conclude that $\alpha\vdash_{\mathcal{S}}\neg\neg\alpha$.
\end{proof}
\begin{cor}
	The logic {\em\Cl} is an extension of {\em\Int}.
\end{cor}
\begin{proof}
	We show that $\vdash_{\textsf{Cl}}\text{ax11}$, that is, for any $\Lan_A$-formulas $\alpha$ and $\beta$,
	\begin{equation}\label{E:Int-axiom-valid-in-Cl}
		\vdash_{\textsf{Cl}}\beta\rightarrow(\neg\beta\rightarrow\alpha).
	\end{equation}
	
	Indeed, since $\beta,\neg\beta,\neg\alpha\vdash_{\textsf{Cl}}\beta$ and
	$\beta,\neg\beta,\neg\alpha\vdash_{\textsf{Cl}}\neg\beta$, by the hyperrrule (c--$ii$), we derive that $\beta,\neg\beta\vdash_{\textsf{Cl}}\neg\neg\alpha$.
	Also, since, in virtue of Lemma~\ref{L:hyperrules-c_ii-c_iii}, $\neg\neg\alpha\vdash_{\textsf{Cl}}\alpha$, by cut, we obtain that
	$\beta,\neg\vdash_{\textsf{Cl}}\alpha$. Then, we apply twice the hyperrule
	(c--$i$).
\end{proof}

\begin{lem}\label{L:distributive-lettice-sound}
	For any $\Lan_A$-formulas $\alpha$, $\beta$ and $\gamma$, the following formulas and their converses are theses of {\em\Pos}:
	{\em
		\[
		\begin{array}{cl}
			(\text{a}) &(\alpha\wedge\beta)\rightarrow(\beta\wedge\alpha),\\
			(\text{b}) &(\alpha\wedge(\beta\wedge\gamma))\rightarrow((\alpha\wedge\beta)\wedge\gamma),\\
			(\text{c}) &((\alpha\wedge\beta)\wedge\gamma)\rightarrow(\alpha\wedge(\beta\wedge\gamma)),\\
			(\text{d}) &((\alpha\wedge\beta)\vee\beta)\rightarrow\beta,\\
			(\text{e}) &(\alpha\vee\beta)\rightarrow(\beta\vee\alpha),\\
			(\text{f}) &(\alpha\vee(\beta\vee\gamma))\rightarrow((\alpha\vee\beta)\vee\gamma),\\
			(\text{g}) &((\alpha\wedge\beta)\vee\beta)\rightarrow\beta,\\
			(\text{h}) &\alpha\wedge(\alpha\vee\beta)\rightarrow\alpha,\\
			(\text{j}) &(\alpha\wedge(\beta\vee\gamma))\rightarrow((\alpha\wedge\beta)\vee(\alpha\wedge\gamma)),\\
			(\text{i}) &(\alpha\vee(\beta\wedge\gamma))\rightarrow((\alpha\vee\beta)\wedge(\alpha\vee\gamma)).
		\end{array}
		\]}
\end{lem}
\begin{proof}
	We show about one of these formulas, namely (j), that it and its converse are theses of \Pos. First, we prove that 
	\[
	\vdash_{\textsf{P}}(\alpha\wedge(\beta\vee\gamma))\rightarrow((\alpha\wedge\beta)\vee(\alpha\wedge\gamma)). 
	\]
	
	We obtain:
	\[
	\begin{array}{rll}
		1a. &\alpha\wedge(\beta\vee\gamma)\vdash_{\textsf{P}}\alpha
		&[\text{Lemma~\ref{L:rules-sound-wrt-P-Int-Cl} and rule (b--$i$)}]\\
		2a. &\alpha\wedge(\beta\vee\gamma)\vdash_{\textsf{P}}\beta\vee\gamma
		&[\text{Lemma~\ref{L:rules-sound-wrt-P-Int-Cl} and rule (b--$ii$)}]\\
		3a. &\alpha,\beta\vdash_{\textsf{P}}\alpha\wedge\beta 
		&[\text{Lemma~\ref{L:rules-sound-wrt-P-Int-Cl} and rule (a--$i$)}]\\
		4a. &\alpha\wedge(\beta\vee\gamma),\beta\vdash_{\textsf{P}}\alpha\wedge\beta
		&[\text{from 1$a$ and 3$a$ by cut (Section~\ref{section:consequence-relation})}]\\
		5a. &\alpha\wedge\beta\vdash_{\textsf{P}}(\alpha\wedge\beta)\vee(\alpha\wedge\gamma) &[\text{Lemma~\ref{L:rules-sound-wrt-P-Int-Cl} and rule (a--$ii$)}]\\
		6a. &\alpha\wedge(\beta\vee\gamma),\beta\vdash_{\textsf{P}}(\alpha\wedge\beta)
		\vee(\alpha\wedge\gamma) &[\text{from 4$a$ and 5$a$ by cut}]\\
		7a. &\alpha,\gamma\vdash_{\textsf{P}}\alpha\wedge\gamma
		&[\text{Lemma~\ref{L:rules-sound-wrt-P-Int-Cl} and rule (a--$i$)}]\\
		8a. &\alpha\wedge(\beta\vee\gamma),\gamma\vdash_{\textsf{P}}\alpha\wedge\gamma
		&[\text{from 1$a$ and 7$a$ by cut}]\\
		9a. &\alpha\wedge\gamma\vdash_{\textsf{P}}(\alpha\wedge\beta)\vee(\alpha\wedge\gamma) &[\text{Lemma~\ref{L:rules-sound-wrt-P-Int-Cl} and rule (a--$iii$)}]\\
		10a. &\alpha\wedge(\beta\vee\gamma),\gamma\vdash_{\textsf{P}}
		(\alpha\wedge\beta)\vee(\alpha\wedge\gamma) &[\text{from 8$a$ and 9$a$ by cut}]\\
		11a. &\alpha\wedge(\beta\vee\gamma),\beta\vee\gamma\vdash_{\textsf{P}}
		(\alpha\wedge\beta)\vee(\alpha\wedge\gamma) 
		&[\text{from 6$a$ and 10$a$ by hyperrule (c--$iii$)}]\\
		12a. &\alpha\wedge(\beta\vee\gamma)\vdash_{\textsf{P}}
		(\alpha\wedge\beta)\vee(\alpha\wedge\gamma) &[\text{from 2$a$ and 11$a$ by cut}]\\
		13a. &\vdash_{\textsf{P}}(\alpha\wedge(\beta\vee\gamma))\rightarrow((\alpha\wedge\beta)\vee(\alpha\wedge\gamma)). &[\text{from 12$a$ by hyperrule (c--$i$)}]
	\end{array}
	\]
	
	Further, we have:
	\[
	\begin{array}{rll}
		1b. &\alpha\wedge\beta\vdash_{\textsf{P}}\alpha &[\text{Lemma~\ref{L:rules-sound-wrt-P-Int-Cl} and rule (b--$i$)}]\\
		2b. &\alpha\wedge\beta\vdash_{\textsf{P}}\beta &[\text{Lemma~\ref{L:rules-sound-wrt-P-Int-Cl} and rule (b--$ii$)}]\\
		3b. &\beta\vdash_{\textsf{P}}\beta\vee\gamma &[\text{Lemma~\ref{L:rules-sound-wrt-P-Int-Cl} and rule (a--$ii$)}]\\
		4b. &\alpha\wedge\beta\vdash_{\textsf{P}}\beta\vee\gamma
		&[\text{from 2$b$ and 3$b$ by cut}]\\
		5b. &\alpha\wedge\gamma\vdash_{\textsf{P}}\alpha &[\text{Lemma~\ref{L:rules-sound-wrt-P-Int-Cl} and rule (b--$i$)}]\\
		6b. &\alpha\wedge\gamma\vdash_{\textsf{P}}\gamma &[\text{Lemma~\ref{L:rules-sound-wrt-P-Int-Cl} and rule (b--$ii$)}]\\
		7b. &\gamma\vdash_{\textsf{P}}\beta\vee\gamma &[\text{Lemma~\ref{L:rules-sound-wrt-P-Int-Cl} and rule (a--$iii$)}]\\
		8b. &\alpha\wedge\gamma\vdash_{\textsf{P}}\beta\vee\gamma
		&[\text{from 6$b$ and 7$b$ by cut}]\\
		9b. &(\alpha\wedge\beta)\vee(\alpha\wedge\gamma)\vdash_{\textsf{P}}\beta\vee\gamma
		&[\text{Lemma~\ref{L:hyperrules-c_ii-c_iii} and hyperrule (c--$iii$)}]\\
		10b. &(\alpha\wedge\beta)\vee(\alpha\wedge\gamma)\vdash_{\textsf{P}}\alpha
		&[\text{Lemma~\ref{L:hyperrules-c_ii-c_iii} and hyperrule (c--$iii$)}]\\
		11b. &\alpha,\beta\vee\gamma\vdash_{\textsf{P}}\alpha\wedge(\beta\vee\gamma) &[\text{Lemma~\ref{L:rules-sound-wrt-P-Int-Cl} and rule (a--$i$)}]\\
		12b. &(\alpha\wedge\beta)\vee(\alpha\wedge\gamma),\alpha\vdash_{\textsf{P}}
		\alpha\wedge(\beta\vee\gamma) &[\text{from 9$b$ and 11$b$ by cut}]\\
		13b. &(\alpha\wedge\beta)\vee(\alpha\wedge\gamma)\vdash_{\textsf{P}}
		\alpha\wedge(\beta\vee\gamma) &[\text{from 10$b$ and 12$b$ by cut}]\\
		14b. &\vdash_{\textsf{P}}((\alpha\wedge\beta)\vee(\alpha\wedge\gamma))\rightarrow
		(\alpha\wedge(\beta\vee\gamma)). &\text{from 13$b$ by hyperrule (c--$i$)}
	\end{array}
	\]
	
	We leave for the reader to complete the proof of this lemma. (Exercise~\ref{section:lindenbaum-algebra}.\ref{EX:distributive-lettice-sound})
\end{proof}

\begin{lem}\label{L:Cl-implications}
	For any $\Lan_A$-formulas $\alpha$ and $\beta$, the following formulas and their converses are theses of {\em\Cl}:
	{\em\[
		\begin{array}{cl}
			(\text{k}) &((\alpha\wedge\neg\alpha)\vee\beta)\rightarrow\beta,\\
			(\text{l}) &((\alpha\vee\neg\alpha)\wedge\beta)\rightarrow\beta.
		\end{array}
		\]}
\end{lem}
\begin{proof}
	We prove that (k) and its converse are theses of {\Cl} and leave the rest to the reader. (Exercise~\ref{section:lindenbaum-algebra}.\ref{EX:Cl-implications})
	
	First, we prove that
	\[
	\vdash_{\textsf{Cl}}((\alpha\wedge\neg\alpha)\vee\beta)\rightarrow\beta.
	\]
	
	Indeed, we have:
	\[
	\begin{array}{cll}
		1a. &\alpha\wedge\neg\alpha\vdash_{\textsf{Cl}}\beta &[\text{Lemma~\ref{L:rules-sound-wrt-P-Int-Cl}, rules (b-$i$), (b--$ii$), \eqref{E:Int-axiom-valid-in-Cl} and cut}]\\
		2a. &\beta\vdash_{\textsf{Cl}}\beta\\
		3a. &(\alpha\wedge\neg\alpha)\vee\beta\vdash_{\textsf{Cl}}\beta
		&[\text{Lemma~\ref{L:hyperrules-c_ii-c_iii} and hyperule (c--$iii$)}]\\
		4a. &\vdash_{\textsf{Cl}}((\alpha\wedge\neg\alpha)\vee\beta)\rightarrow\beta.
		&[\text{Lemma~\ref{L:hyperrules-c_ii-c_iii} and hyperule (c--$i$)}]
	\end{array}
	\]
	
	For the converse, we have:
	\[
	\begin{array}{cll}
		1b. &\vdash_{\textsf{Cl}}\beta\rightarrow((\alpha\wedge\neg\alpha)\vee\beta).
		&[\text{rule ax7}]
	\end{array}
	\]
\end{proof}

\begin{prop}
	{\em$\mathbb{T}_{\textsf{Cl}}$} coincides with the class of all Boolean algebras. 
\end{prop}
\begin{proof}
	Let $\alg{B}=\left\langle \textsf{B};\wedge,\vee,\neg,\one\right\rangle$ be a Boolean algebra. Expanding this algebra with an operation
	\[
	x\rightarrow y:=\neg x\vee y, \tag{$\ast$}
	\]
	according to Proposition~\ref{P:boolean-algebra-as-heyting}, we obtain a Heyting algebra $\alg{H}=\left\langle \textsf{B};\wedge,\vee,\rightarrow,\neg,\one\right\rangle$. In virtue of
	Proposition~\ref{P:Int-properties}, the rules ax1--ax9 are valid in $\alg{H}$ and hence in \alg{B}. In view of Corollary~\ref{C:Cl-property}, the rule ax10 is also valid in \alg{B}. And, because of the property $(\text{h}_2)$
	(Section~\ref{section:heyting-algebra}), modus ponens preserves this validity over any \Cl-derivation.
	
	Now, assume that an algebra $\alg{A}=\left\langle \textsf{A};\wedge,\vee,\neg,\one\right\rangle\in\mathbb{T}_{\textsf{Cl}}$. That is, for any $\Lan_A$-formula $\alpha\in\bm{T}_{\textsf{Cl}}$,
	$\alg{A}\models\alpha\approx\one$. This means that for any $\alpha\in\bm{T}_{\textsf{Cl}}$ and any $\Lan_A$-valuation $v$ in \alg{A}, $v(\alpha)=\one$. Therefore, in virtue of~Lemma~\ref{L:distributive-lettice-sound} and the fact that $\Pos$ is implicational w.r.t $\rightarrow$, the properties $(\text{l}_1)$--$(\text{l}_4)$ are valid in \alg{A}. Namely, we have:
	(a) and its converse prove the validity of ($\text{l}_1$--$i$); (b) and its converse prove the validity of ($\text{l}_1$--$ii$); (c) and its converse prove the validity of ($\text{l}_2$--$i$);
	(f) and its converse prove the validity of ($\text{l}_2$--$ii$);
	(g) and its converse prove the validity of ($\text{l}_3$--$i$);
	(h) and its converse prove the validity of ($\text{l}_3$--$ii$);
	(j) and its converse prove the validity of ($\text{l}_4$--$i$);
	(i) and its converse prove the validity of ($\text{l}_4$--$ii$).
	Thus $\langle \textsf{A};\wedge,\vee\rangle$ is a distributive lattice.
	
	Further, if $\bigstar$ is an arbitrary {\Cl}-thesis and $p$ is an arbitrary variable, then
	\begin{equation}\label{E:Cl-thesis}
		\vdash_{\textsf{Cl}}p\rightarrow\bigstar.
	\end{equation}
	This implies that for any element $x\in\textsf{A}$, $x\le\one$, that is
	the equalities $(\text{b}_{1})$ hold in $\langle \textsf{A};\wedge,\vee,\one\rangle$,
	where $\one$ turns to be a greatest element in the sense of $\le$ (Section~\ref{section:boolean-algebra}).
	
	Finally, we use Lemma~\ref{L:Cl-implications} and the fact that {\Cl} is and implicational abstract logic with respect to $\rightarrow$ to conclude that $\alg{A}$ is a Boolean algebra.
\end{proof}
\begin{cor}\label{C:LT_Cl(k)}
	{\em$\LTCl$} is a free Boolean algebra of rank $\card{\Var_{\mathcal{L}_A}}$.
	Hence, {\em$\LTCl(\kappa)$}, where $\kappa\le\card{\Var_{\mathcal{L}_A}}$,  is a free Boolean algebras of rank $\kappa$.
\end{cor}

According to Corollary~\ref{C:LT_Cl(k)} and Proposition~\ref{P:finitely-generated-boolean}, any $\LTCl(\kappa)$ of a finite rank $\kappa$ is finite. Let us take $\kappa=1$, assuming that $\Var=\lbrace p\rbrace$. For simplicity, we denote
\[
[p]:=p\slash\theta(\bm{T}_{\textsf{Cl}}).
\]

It is not difficult to show that the subalgebra of $\LTCl$ generated by $[p]$ consists of the four pairwise distinct elements: $[p]$, $[\neg p]$, $[p\vee\neg p]$ and $[p\wedge\neg p]$. (Exercise~\ref{section:some-lindenbaum-algebras}.\ref{EX:four-elements})This subalgebra is $\LTCl(1)$. The relation $\le$ between $[\alpha]$ and $[\beta]$ of this algebra is defined, according to \eqref{E:ordering-in-lattice}, as follows:
\[
[\alpha]\le[\beta]\stackrel{\text{df}}{\Longleftrightarrow}[\alpha]\wedge[\beta]=[\alpha].
\]

Thus, we arrive at the following diagram:

\begin{figure}[!ht]	
	\[
	\ctdiagram{
		\ctnohead
		\ctinnermid
		\ctel 0,0,20,20:{}
		\ctel 0,0,-20,20:{}
		\ctel 0,40,20,20:{}
		\ctel 0,40,-20,20:{}
		\ctv 0,0:{\bullet}
		\ctv 20,20:{\color{red}\bullet}
		\ctv 28,20:{[p]}
		\ctv -20,20:{\bullet}
		\ctv -32,20:{[\neg p]}
		\ctv 0,40:{\bullet}
		\ctv 0,47:{\mathbf{1}=[p\vee\neg p]}
		\ctv 0,-8:{\mathbf{0}=[p\wedge\neg p]}
	}
	\]
	\caption{Lindenbaum-Tarski algebra $\LTCl(1)$}
\end{figure}

If $\Var=\lbrace p,q\rbrace$, we obtain $\LTCl(2)$ as a subalgebra of $\LTCl$ generated by the classes $[p]$ and $[q]$. If we identify any congruence class of this subalgebra with its representative, we can depict 	$\LTCl(2)$ by the following diagram:
\pagebreak
\begin{figure}[!ht]
	
	\[
	\ctdiagram{
		\ctnohead
		\ctinnermid
		\ctel 0,0,30,30:{}
		\ctel 0,0,-30,30:{}
		\ctel 0,60,30,30:{}
		\ctel 0,60,-30,30:{}
		\ctel 0,30,30,60:{}
		\ctel 0,30,-30,60:{}
		\ctel 0,90,30,60:{}
		\ctel 0,90,-30,60:{}
		\ctel 0,0,0,30:{}
		\ctel -30,30,-30,60:{}
		\ctel 30,30,30,60:{}
		\ctel 0,60,0,90:{}
		\ctel 120,30,90,60:{}
		\ctel 120,30,150,60:{}
		\ctel 90,60,120,90:{}
		\ctel 150,60,120,90:{}
		\ctel 120,60,150,90:{}
		\ctel 120,60,90,90:{}
		\ctel 150,90,120,120:{}
		\ctel 90,90,120,120:{}
		\ctel 120,30,120,60:{}
		\ctel 90,60,90,90:{}
		\ctel 150,60,150,90:{}
		\ctel 120,90,120,120:{}
		\ctel 0,0,120,30:{}
		\ctel -30,30,90,60:{}	
		\ctel 30,30,150,60:{}	
		\ctel 0,60,120,90:{}	
		\ctel 0,30,120,60:{}		
		\ctel 30,60,150,90:{}
		\ctel -30,60,90,90:{}
		\ctel 0,90,120,120:{}
		\ctv 0,0:{\bullet}
		\ctv 30,30:{\bullet}
		\ctv -30,30:{\bullet}
		\ctv 0,60:{\color{red}\bullet}
		\ctv 0,53:{p}
		\ctv 0,30:{\bullet}
		\ctv -45,30:{p\land q}
		\ctv 30,60:{\bullet}
		\ctv -30,60:{\bullet}
		\ctv -80,60:{(p \land q) \lor (\neg p \land \neg q)}
		\ctv 0,90:{\bullet}
		\ctv -22,92:{\neg p \lor \neg q}
		\ctv 120,30:{\bullet}
		\ctv 136,31:{\neg p \land q}
		\ctv 150,60:{\bullet}
		\ctv 200,60:{ (p \land \neg q) \lor (\neg p \land q)}
		\ctv 90,60:{\color{red}\bullet}
		\ctv 120,90:{\bullet}
		\ctv 120,60:{\bullet}
		\ctv 128,60:{\neg p}
		\ctv 150,90:{\bullet}
		\ctv 170,90:{\neg p \lor \neg q}
		\ctv 90,90:{\bullet}
		\ctv 120,120:{\bullet}
		\ctv 46,27:{p \land \neg q}
		\ctv -4,23:{\neg p \land \neg q}
		\ctv 72,95:{\neg p \lor q}
		\ctv 122,96:{p \lor q}
		\ctv 38,58:{\neg q}
		\ctv 85,64:{q}
		\ctv 120,127:{\mathbf{1}}
		\ctv 0,-10:{\mathbf{0}}		
	}
	\]
	\caption{Lindenbaum-Tarski algebra $\LTCl(2)$}
\end{figure}	

Now we turn to Heyting algebras.
\begin{lem}\label{L:Int-implications}
	Let $\bigstar$ be an arbitrary {\em\Int}-thesis. For any $\Lan_A$-formulas $\alpha$, $\beta$ and $\gamma$, the following implications and their converses are {\em\Int}-theses:
	{\em\[
		\begin{array}{cl}
			(\text{a}) &(\alpha\wedge(\alpha\rightarrow\beta))\rightarrow(\alpha\wedge\beta),\\
			(\text{b}) &((\alpha\rightarrow\beta)\wedge\beta)\rightarrow\beta,\\
			(\text{c}) &((\alpha\rightarrow\beta)\wedge(\alpha\rightarrow\gamma)
			\rightarrow(\alpha\rightarrow(\beta\wedge\gamma))\\
			(\text{d}) &(\alpha\wedge(\beta\rightarrow\beta))\rightarrow\alpha,\\
			(\text{e}) &(\neg\bigstar\vee\alpha)\rightarrow\alpha,\\
			(\text{f}) &\neg\alpha\rightarrow(\alpha\rightarrow\neg\bigstar).
		\end{array}
		\]}
\end{lem}
\begin{proof}
	We prove that $\vdash_{\textsf{Int}}(\neg\bigstar\vee\alpha)\rightarrow\alpha$
	and $\vdash_{\textsf{Int}}\alpha\rightarrow(\neg\bigstar\vee\alpha)$ and leave for the reader to complete this proof. (Exercise~\ref{section:lindenbaum-algebra}.\ref{EX:Int-implications})
	
	Indeed, we obtain:
	\[
	\begin{array}{cll}
		1. &\vdash_{\textsf{Int}}\bigstar &[\text{by premise $\bigstar$ is an {\Int}-thesis}]\\
		2. &\vdash_{\textsf{Int}}\bigstar\rightarrow(\neg\bigstar\rightarrow\alpha)
		&[\text{rule ax11}]\\
		3. &\vdash_{\textsf{Int}}\neg\bigstar\rightarrow\alpha
		&[\text{from $1$ and $2$ by modus ponens]}\\
		4. &\vdash_{\textsf{Int}}\alpha\rightarrow\alpha
		&[\text{because of $(\ast)$ in the proof of Lemma~\ref{L:deduction-theorem}}]\\
		5. &\vdash_{\textsf{Int}}(\neg\bigstar\vee\alpha)\rightarrow\alpha.
		&[\text{from ax8, $3$ and $4$ by modus ponens}] 
	\end{array}
	\]
	
	Further, $\alpha\rightarrow(\neg\bigstar\vee\alpha)$ is an instantiation of ax7.
\end{proof}

\begin{prop}\label{P:T_Int=Heyting-algebras}
	{\em$\mathbb{T}_{\textsf{Int}}$} coincides with the class of all Heyting algebras.
\end{prop}
\begin{proof}
	Let $\alg{H}=\langle\textsf{H};\wedge,\vee,\rightarrow,\neg,\one\rangle$
	be a Heyting algebra. In virtue of
	Proposition~\ref{P:Int-properties}, the rules ax1--ax9 are valid in $\alg{H}$.
	In virtue of the properties \eqref{E:zero-in-heyting}, \eqref{E:pseudo-complementation} and~\eqref{E:less-than=implication}, we receive that ax11 is also valid in \alg{H}. Noticing, that modus ponens preserves validity in any Heyting algebra, we conclude that if $\alpha\in\bm{T}_{\textsf{Int}}$, then $\alg{H}\models\alpha\approx\one$.
	
	Next, assume that an algebra $\alg{H}=\langle\textsf{H};\wedge,\vee,\rightarrow,\neg,\one\rangle\in\mathbb{T}_{\textsf{Int}}$. That is, for any $\alpha\in\bm{T}_{\textsf{Int}}$ and any
	valuation $v$ in \alg{H}, $v(\alpha)=\one$. We show that the properties
	$(\text{l}_1)$--$(\text{l}_4)$, $(\text{b}_1)$, and $(\text{h}_1)$--$(\text{h}_6)$ are valid in \alg{H}.

	In virtue of Lemma~\ref{L:distributive-lettice-sound} and the fact that $\Pos$ is implicational w.r.t $\rightarrow$, the properties $(\text{l}_1)$--$(\text{l}_4)$ are valid in \alg{H}. Thus $\langle\textsf{H};\wedge,\vee\rangle$ is a distributive lattice.
	
	Further, similarly to~\eqref{E:Cl-thesis}, we derive that 
	\[
	\vdash_{\textsf{Int}}p\rightarrow\bigstar,
	\]
	for any variable $p$ and \textsf{Int}-thesis $\bigstar$. This implies that the equalities $(\text{b}_1)$ are true in \alg{H}.
	
	To prove the properties $(\text{h}_1)$--$(\text{h}_6)$, we apply Lemma~\ref{L:Int-implications}.
\end{proof}
\begin{cor}\label{C:LT_Int}
	{\em$\LTInt$} is a free Heyting algebra of rank $\card{\Var_{\mathcal{L}_A}}$.
	Hence, {\em$\LTInt(\kappa)$}, where $\kappa\le\card{\Var_{\mathcal{L}_A}}$,  is a free Heyting algebras of rank $\kappa$.
\end{cor}

Below we give a description of $\LTInt(1)$. For this, we consider the following formulas of an arbitrary variable $p$.

First, we define:
\[
p^0:=p\land\neg p, \qquad p^1:=\neg p, \qquad p^2:=p, \qquad p^{\omega}:=p\to p.
\]
Then for any $n\ge 0$, we define:
\begin{equation}\label{E:powers_p^n}
	\begin{array}{ll}
		i)\quad p^{2n+3}:=p^{2n+1}\to p^{2n} &ii)\quad p^{2n+4}:= p^{2n+1}\lor p^{2n+2}.
	\end{array}
\end{equation}
Thus `$m$' in `$p^m$' can take the value of any ordinal less than or equal to the first infinite ordinal, symbolically $m\le \omega$.

We conclude with the following observation, due to I. Nishimura~\cite{nishimura1960}; see also~\cite{rieger1957}.\\

A description of $\LTInt(1)$ can be outlined in terms of the following relations on the set of $\Forms_{\mathcal{L}_A}$-formulas.
\begin{equation}\label{E:RN-main-relations}
	\begin{array}{cl}
		\bullet &\alpha\preccurlyeq\beta~\define~\vdash_{\textsf{Int}}\alpha\rightarrow\beta;\\
		\bullet &\alpha\sim\beta~\define~\vdash_{\textsf{Int}}\alpha\rightarrow\beta~\text{and}~\vdash_{\textsf{Int}}\beta\rightarrow\alpha;\\
		\bullet &\alpha\prec\beta~\define~\alpha\preccurlyeq\beta~\text{and not}~\alpha\sim\beta;\\
	\end{array}
\end{equation}

It should be clear that
\[
\alpha\sim\beta~\Longleftrightarrow~\alpha\slash\theta(\bm{T}_{\Int})=\beta\slash\theta(\bm{T}_{\Int});
\]
and, since
\[
\alpha\sim\beta~\Longleftrightarrow~\alpha\preccurlyeq\beta~\text{and}~\beta\preccurlyeq\alpha,
\]
\[
\alpha\slash\theta(\bm{T}_{\Int})\le\beta\slash\theta(\bm{T}_{\Int})~\Longleftrightarrow~\alpha\preccurlyeq\beta,
\]
where $\le$ is the (ordinary) partial order in the Heyting algebra $\LTInt$ (Corollary~\ref{C:LT_Int}). Below we treat $\sim$ as a congruence on $\Forms_{\mathcal{L}_A}$.

The above equivalences will be used without reference. In particular, in the description of $\LTInt(1)$ as a diagram below, we depict each congruence class by its representative. Thus if a formula $\alpha$ can be connected by an upward path with a formula $\beta$, which means that $\alpha\slash\theta(\bm{T}_{\Int})\le\beta\slash\theta(\bm{T}_{\Int})$, then $\alpha\preccurlyeq\beta$. Also, $\alpha\slash\theta(\bm{T}_{\Int})<\beta\slash\theta(\bm{T}_{\Int})$ if, and only if, $\alpha\prec\beta$.

In accordance with this agreement, we have:
\[
\zero:=p^0~~\text{and}~~\one:= p^{\omega}.
\]

It should be clear that for any $\Forms_{\mathcal{L}_A}$-formula $\alpha$,
\[
\zero\preccurlyeq\alpha~\text{and}~\alpha\preccurlyeq\one.
\]

Thus the elements $\zero$ and $\one$ are a least and greatest elements, respectively, of any Lindenbaum-Tarski algebra of $\Int$ and hence are those of $\LTInt(1)$.
\begin{prop}\label{P:nishimura}
	Every $\Lan_A$-formula of one variable $p$ is equivalent in {\em\Int} to one and only one formula $p^i$, which are not equivalent to each other $($in the sense of $\sim$$)$ and related to each other $($in the sense of $\preccurlyeq$$)$ as shown in the diagram of Figure~\ref{fig-RN}.
\end{prop}
\pagebreak
\begin{figure}[ht!]
	\[
	\begin{array}{l}
		\ctdiagram{
			\ctnohead
			\ctinnermid
			\ctv 45,0:{\bullet}
			\ctv 45,10:{p^\omega}
		}\\	
		\quad\qquad .\quad .\quad .\\
		\quad\qquad .\quad .\quad .\\
		\quad\qquad .\quad .\quad .\\
		\ctdiagram{
			\ctnohead
			\ctinnermid
			\ctel -7,-7,40,40:{}
			\ctel 7,-7,-20,20:{}
			\ctel -20,20,40,80:{}
			\ctel 27,13,-20,60:{}
			\ctel -20,100,40,40:{}
			\ctel 30,110,-20,60:{}
			\ctel 10,110,40,80:{}
			\ctel -10,110,-20,100:{}
			\ctv 0,0:{\bullet}
			\ctv 14,0:{p^{2n}}
			\ctv 20,20:{\bullet}
			\ctv 37,20:{p^{2n+2}}
			\ctv -20,20:{\bullet}
			\ctv -35,20:{p^{2n+1}}
			\ctv 0,40:{\bullet}
			\ctv -20,40:{p^{2n+4}}
			\ctv 40,40:{\bullet}
			\ctv 55,40:{p^{2n+3}}
			\ctv 20,60:{\bullet}
			\ctv 37,60:{p^{2n+6}}
			\ctv -20,60:{\bullet}
			\ctv -35,60:{p^{2n+5}}
			\ctv 0,80:{\bullet}
			\ctv -20,80:{p^{2n + 8}}
			\ctv 40,80:{\bullet}
			\ctv 55,80:{p^{2n+7}}
			\ctv 20,100:{\bullet}
			\ctv 40,100:{p^{2n+10}}
			\ctv -20,100:{\bullet}
			\ctv -35,100:{p^{2n+9}}
		}\\
		\quad\qquad .\quad .\quad .\\
		\quad\qquad .\quad .\quad .\\
		\quad\qquad .\quad .\quad .\\
		\ctdiagram{
			\ctnohead
			\ctinnermid
			\ctel 0,0,40,40:{}
			\ctel 0,0,-20,20:{}
			\ctel -20,20,40,80:{}
			\ctel 20,20,-20,60:{}
			\ctel -10,90,40,40:{}
			\ctel 10,90,-20,60:{}
			\ctel 30,90,40,80:{}
			\ctv 0,0:{\bullet}
			\ctv 10,0:{p^{0}}
			\ctv 20,20:{\color{red}\bullet}
			\ctv 30,20:{p^{2}}
			\ctv -20,20:{\bullet}
			\ctv -30,20:{p^{1}}
			\ctv 0,40:{\bullet}
			\ctv -15,40:{p^{4}}
			\ctv 40,40:{\bullet}
			\ctv 50,40:{p^{3}}
			\ctv 20,60:{\bullet}
			\ctv 30,60:{p^{6}}
			\ctv -20,60:{\bullet}
			\ctv -30,60:{p^{5}}
			\ctv 0,80:{\bullet}
			\ctv -15,80:{p^{8}}
			\ctv 40,80:{\bullet}
			\ctv 50,80:{p^{7}}
		}
	\end{array}
	\]
	\caption{Lindenbaum-Tarski algebra $\LTInt(1)$}\label{figure:Rieger-Nishimura algebra} \label{fig-RN}
\end{figure}
\pagebreak

We call the above diagram \textit{\textbf{Rieger-Nishimura algebra}}; other names in the literature are \textit{Rieger-Nishimura ladder} and \textit{Rieger-Nishimura lattice}.\index{algebra!Rieger-Nishimura}\index{Rieger-Nishimura ladder}
The rest of this section is devoted to the proof of Proposition~\ref{P:nishimura}.

\subsection{Rieger-Nishimura algebra}\label{section:rieger-nichimura-algebra}
In the proofs of propositions and lemmas below in this subsection we use Proposition~\ref{P:replacement} and transitivity of $\sim$ without reference.
Also, we use rules of (a) and (b) and the hyperrules of (c) (Section~\ref{section:inference-rules}), as needed. The applications of those rules are guaranteed by Lemma~\ref{L:rules-sound-wrt-P-Int-Cl}, Lemma~\ref{L:deduction-theorem} and Lemma~\ref{L:hyperrules-c_ii-c_iii}.

At one point or another, we will need to use the $\Int$-theses that we have collected in the following lemma.
\begin{lem}\label{L:Int-theses}
	For arbitrary $\Lan_A$-formulas $\alpha$, $\beta$ and $\gamma$, the following properties hold:	
	{\em\[
		\begin{array}{cl}
			(\text{a}) &(\alpha\rightarrow\neg\beta)\sim\neg(\alpha\land\beta),\\
			(\text{b}) &(\neg\alpha\rightarrow\beta)\preccurlyeq\neg\neg(\alpha\lor\beta),\\
			(\text{c}) &(\alpha\rightarrow(\beta\rightarrow\gamma))\sim((\alpha\land\beta)\rightarrow\gamma),\\
			(\text{d}) &((\alpha\lor\beta)\rightarrow\gamma)\sim((\alpha\rightarrow\gamma)\land(\beta\rightarrow\gamma)),\\
			(\text{e}) &(\alpha\rightarrow(\beta\land\gamma))\sim((\alpha\rightarrow\beta)\land(\alpha\rightarrow\gamma)),\\
			(\text{f}) &(\alpha\rightarrow\neg\alpha)\sim\neg\alpha,\\
			(\text{g}) &(\neg\alpha\rightarrow\alpha)\sim\neg\neg\alpha,\\
			(\text{h}) &(\neg\alpha\rightarrow\neg\neg\alpha)\sim\neg\neg\alpha,\\
			(\text{i}) &(\neg\neg\alpha\rightarrow\neg\alpha)\sim\neg\alpha,\\
			(\text{j}) &\neg\alpha\sim\alpha\rightarrow(p\land\neg p),\\
			(\text{k}) &\neg\neg\neg\alpha\sim\neg\alpha;\\
			(\text{l}) &(\neg\alpha\land\alpha)\sim p^0;\\
			(\text{m}) &((\alpha\rightarrow\beta)\land(\beta\rightarrow\gamma))\rightarrow
			(\alpha\rightarrow\gamma);\\
			(\text{n}) &\alpha\land(\alpha\rightarrow\beta)\sim\alpha\land\beta.
		\end{array}
		\]}
\end{lem}
\begin{proof}
	We prove (a) and leave the rest to the reader as an exercise. (Exercise~\ref{section:lindenbaum-algebra}.\ref{EX:Int-theses})
	
	With the help of rules (b--$i$) (b--$iii$) and the definition of $\Int$, we obtain:
	\[
	\alpha\rightarrow\neg\beta,\alpha\land\beta\vdash_{\textsf{Int}}\alpha;~
	\alpha\rightarrow\neg\beta,\alpha\land\beta\vdash_{\textsf{Int}}\beta;~
	\text{and}~\alpha\rightarrow\neg\beta,\alpha\land\beta,\alpha\vdash_{\textsf{Int}}\neg\beta.
	\]
	Hence we have:
	\[
	\alpha\rightarrow\neg\beta,\alpha\land\beta\vdash_{\textsf{Int}}\beta~
	\text{and}~
	\alpha\rightarrow\neg\beta,\alpha\land\beta\vdash_{\textsf{Int}}\neg\beta.
	\]	
	By hyperrules (c--$ii$) and (c--$i$), we successively derive: 
	\[
	\alpha\rightarrow\neg\beta\vdash_{\textsf{Int}}\neg(\alpha\land\beta)~
	\text{and}~\vdash_{\textsf{Int}}(\alpha\rightarrow\neg\beta)\rightarrow\neg(\alpha\land\beta).
	\]
	
	Conversely, applying rule (a--$i$), we derive:
	\[
	\neg(\alpha\land\beta),\alpha,\beta\vdash_{\textsf{Int}}\alpha\land\beta.
	\]
	Hence, by hyperrule (c--$ii$),
	\[
	\neg(\alpha\land\beta),\alpha\vdash_{\textsf{Int}}\neg\beta.
	\]
	Applying hyperrule (c--$i$) twice, we conclude that
	$\neg(\alpha\land\beta)\rightarrow(\alpha\rightarrow\neg\beta)$ is an $\Int$-thesis. 
	
	Then, it remains to apply rule (a--$i$).
\end{proof}

Further, we add to the relations~\eqref{E:RN-main-relations} one more, namely we define:
\[
\alpha\ll\beta~\define~\vdash_{\textsf{Int}}(\beta\rightarrow\alpha)\rightarrow\beta.
\]

It should be clear that Proposition~\ref{P:nishimura} will be proven when we prove Corollary~\ref{C:nishimura-completeness} and Corollary~\ref{C:nishimura-power-difference} below. These corollaries follow from Proposition~\ref{P:nishimura-completeness} and, respectively, Proposition~\ref{P:nishimura-power-difference}. Thus, our argument consists of two almost separate paths. Both paths, however, pass through  Lemma~\ref{L:relations-specific-properties} which is about properties of the relations $\preccurlyeq$, $\sim$ and $\ll$ for concrete formulas, or say it better, for concrete elements of $\LTInt(1)$. But we start with general properties of these relations.

\begin{lem}\label{L:relations-general-properties}
	For arbitrary $\Lan_A$-formulas $\alpha$, $\beta$ and $\gamma$, the following conditionals and equivalences hold:
	{\em\[
		\begin{array}{cl}
			(\text{a}) &(\alpha\preccurlyeq\beta~\text{and}~\beta\preccurlyeq\gamma)\Longrightarrow\alpha\preccurlyeq\gamma;\\
			(\text{b}) &\alpha\ll\beta~\Longleftrightarrow~\vdash_{\textsf{Int}}(\alpha\rightarrow\beta)\land((\beta\rightarrow\alpha)\rightarrow\alpha);\\
			(\text{c}) &(\beta\rightarrow\alpha)\sim\alpha~\Longleftrightarrow~\vdash_{\textsf{Int}}(\beta\rightarrow\alpha)\rightarrow\alpha;\\
			(\text{d}) &\alpha\ll\beta~\Longleftrightarrow~\alpha\preccurlyeq\beta~\text{and}~\beta\rightarrow\alpha\sim\alpha;\\
			(\text{e}) &((\alpha\rightarrow\beta)\sim\beta~\text{and}~\beta\preccurlyeq\gamma)~\Longrightarrow~\alpha\rightarrow\gamma\sim\gamma;\\
			(\text{f}) &(\alpha\ll\beta~\text{and}~\beta\preccurlyeq\gamma)~\Longrightarrow~\alpha\ll\gamma;\\
			(\text{g}) &(\alpha\preccurlyeq\beta~\text{and}~\beta\ll\gamma)~\Longrightarrow~\alpha\ll\gamma;\\
			(\text{h}) &(\alpha\sim \beta~\text{and}~\beta\sim\gamma)~\Longrightarrow~\alpha\sim\gamma;\\
			(\text{i}) &\alpha\preccurlyeq\beta~\Longleftrightarrow~\alpha\land\beta\sim\alpha;\\
			(\text{j}) &\alpha\preccurlyeq\beta~\Longleftrightarrow~\alpha\lor\beta\sim\beta;\\
			(\text{k}) &\alpha\lor\beta\preccurlyeq\gamma~\Longleftrightarrow~\alpha\preccurlyeq\gamma~\text{and}~\beta\preccurlyeq\gamma;\\
			(\text{l}) &\alpha\land\beta\preccurlyeq\gamma~\Longleftrightarrow~\alpha\preccurlyeq\beta\rightarrow\gamma.
		\end{array}
		\]	}
\end{lem}
\begin{proof}
	We prove (d) and leave other properties for the reader to prove as exercises. (Exercise~\ref{section:lindenbaum-algebra}.\ref{EX:relations-general-properties})	
	
	Suppose that $\alpha\ll\beta$, that is, $\vdash_{\textsf{Int}}(\beta\rightarrow\alpha)\rightarrow\beta$. This means that there is an $\Int$-derivation without any premise:
	\[
	\ldots,~(\beta\rightarrow\alpha)\rightarrow\beta,\tag{d--$\ast$}
	\]
	which can be continued as an $\Int$-derivation with the premise $\beta\rightarrow\alpha$:
	\[
	\underbrace{\ldots,~(\beta\rightarrow\alpha)\rightarrow\beta}_{\text{derivation (d--$\ast$)}},~\underbrace{\beta\rightarrow\alpha}_{\text{premise}}, \beta, \alpha.
	\]
	Thus, we have: $\beta\rightarrow\alpha\vdash_{\textsf{Int}}\alpha$. Applying
	hyperrule (c--$i$), we obtain that $\beta\rightarrow\alpha\preccurlyeq\alpha$.
	On the other hand, it is obvious that $\alpha\preccurlyeq \beta\rightarrow\alpha$. Applying rule (a--$i$), we get the equivalence
	$\beta\rightarrow\alpha\sim\alpha$. The latter allows us to replace $\beta\rightarrow\alpha$ with $\alpha$ in $(\beta\rightarrow\alpha)\beta$ and obtain that $\vdash_{\textsf{Int}}\alpha\rightarrow\beta$, that is $\alpha\preccurlyeq\beta$.
	
	Conversely, assume that $\alpha\preccurlyeq\beta$ and $\beta\rightarrow\alpha\sim\alpha$. Applying rule (b--$ii$)  to the second assumption and then (b--$iii$)to the result, we get that $\beta\rightarrow\alpha\vdash_{\textsf{Int}}\alpha$. Hence, along with the first assumption, we have: $\beta\rightarrow\alpha\vdash_{\textsf{Int}}\beta$.
	Applying hyperrule (c--$i$), we get $\vdash_{\textsf{Int}}(\beta\rightarrow\alpha)\rightarrow\beta$, that is $\alpha\ll\beta$.
\end{proof}

Before moving on to the next lemma, recall the following relation, namely, 
\[
p^0\preccurlyeq p^n,
\] 
which will be used without reference.
\begin{lem}\label{L:relations-specific-properties}
	For any $n,i<\omega$, the following properties hold:
	{\em\[
		\begin{array}{cl}
			(\text{a}) &p^{2n}\preccurlyeq p^{2n+i};\\
			(\text{b}) &p^{2n+1}\preccurlyeq p^{2n+4+i};\\
			(\text{c}) &p^n\preccurlyeq p^{n+3+i};\\
			(\text{d}) &p^{2n+1}\lor p^{2n+3}\sim p^{2n+6};\\
			(\text{e}) &p^{2n+1}\land p^{2n+3}\sim p^{2n};\\
			(\text{f}) &p^{2n+1}\land p^{2n+2}\sim p^{2n};\\
			(\text{g}) &p^{2n+1}\rightarrow p^{2n+3}\sim p^{2n+3};\\
			(\text{h}) &p^{2n+3}\rightarrow p^{2n+1}\sim p^{2n+1};\\
			(\text{i}) &p^{2n+1}\rightarrow p^{2n+2}\sim p^{2n+3};\\
			(\text{j}) &p^{2n+2}\rightarrow p^{2n+1}\sim p^{2n+1};\\
			(\text{k}) &p^{2n+2}\rightarrow p^{2n}\sim p^{2n+1};\\
			(\text{l}) & p^k\ll p^n,~\text{providing that $k,n<\omega$ and $n-k\ge 4$};\\
			(\text{m}) &(p^{2n+1}\rightarrow p^n)\rightarrow p^n\sim p^{2n+1}.\\
		\end{array}
		\]	}
\end{lem}
\begin{proof}
	We prove successively each statement of (a)--(m).\\
	\noindent\textit{Proof of}~(a). We observe that for any $n\ge 0$,
	\[
	p^{2n}\preccurlyeq p^{2n+2} \tag{a--$\ast$}
	\]
	
	Indeed, it is obvious if $n=0$.
	Now, assuming that $n=l+1$, we have $p^{2(l+1)+2}=p^{2l+1}\lor p^{2l+2}$.
	
	Now we observe that for any $n\ge 0$,
	\[
	p^{2n+1}\preccurlyeq p^{2n+5} \tag{a--$\ast\ast$}
	\]
	Indeed, we successively obtain the following:
	\[
	\begin{array}{cll}
		(\text{a}) &p^{2n+1},~p^{2n+1}\rightarrow p^{2n}\vdash_{\textsf{Int}}p^{2n} &[\text{by rule (b-$iii$)}]\\
		(\text{b}) &p^{2n+1},~p^{2n+1}\rightarrow p^{2n}\vdash_{\textsf{Int}}p^{2n+2} &[\text{using (e--$\ast$)}]\\
		(\text{c}) &p^{2n+1},~p^{2n+3}\vdash_{\textsf{Int}}p^{2n+2} &[\text{by definition (\ref{E:powers_p^n}--$i$)}]\\
		(\text{d}) &p^{2n+1}\vdash_{\textsf{Int}}p^{2n+3}\rightarrow p^{2n+2}
		&[\text{by hyperrule (c--$i$)}]\\
		(\text{e}) &p^{2n+1}\vdash_{\textsf{Int}}p^{2n+5} &[\text{by definition (\ref{E:powers_p^n}--$i$)}]\\
		(\text{f}) &\vdash_{\textsf{Int}}p^{2n+1}\rightarrow p^{2n+5} &[\text{by hyperrule (c--$i$)}].
	\end{array}
	\]
	
	Next we observe that for any $n\ge 0$,
	\[
	p^{2n}\preccurlyeq p^{2n+1}. \tag{a-$\ast\ast\ast$}
	\]
	Omitting the obvious case when $n=0$, we notice that, when $n=1$, $\vdash_{\textsf{Int}}p\rightarrow p^{3}$, since $p^{3}\sim\neg\neg p$. Now, assuming that $n=l+2$, we, on the one hand, have: $p^{2n}=p^{2l+1}\lor p^{2l+2}$ and, on the other hand, $p^{2n+1}=p^{2(l+1)+1}\rightarrow p^{2l+2}$. It must be obvious that $p^{2l+2}\preccurlyeq p^{2l+3}\rightarrow p^{2(l+1)}$. On the other hand, in the light of the rule (b--$iii$), $p^{2l+1},~p^{2l+1}\rightarrow p^{2l}\vdash_{\textsf{Int}}p^{2l}$ and hence, in virtue of (a--$\ast\ast\ast$), $p^{2l+1},~p^{2l+1}\rightarrow p^{2l}\vdash_{\textsf{Int}}p^{2l+2}$. Applying the hyperrule (c--$i$), we obtain that $p^{2l+1}\preccurlyeq p^{2l+3}\rightarrow p^{2l+2}$. Then, we apply Lemma~\ref{L:relations-general-properties}--k.
	
	Now, we begin the proof of (a) by induction on $i$. The case when $i=0$ is obvious. The case when $i=1$ follows from (a--$\ast\ast\ast$) and the case when $i=2$ follows from (a--$\ast\ast$). 
	
	By induction hypothesis, assume that
	\[
	p^{2n}\preccurlyeq p^{2n+k},
	\]
	for any $k$ with $0\le k<i$, where $i\ge 2$. Thus $0\le k\le i-2$.
	
	We consider two cases: $2n+(i-1)$ is even and $2n+(i-1)$ is odd.
	
	If $2n+(i-1)$ is even, then, in virtue of (a--$\ast\ast\ast$), 
	$p^{2n+(i-1)}\preccurlyeq p^{2n+i}$ and hence, according Lemma~\ref{L:relations-general-properties}--k, $p^{2n}\preccurlyeq p^{2n+i}$.
	
	Now let $2n+(i-1)=2l+1$. This implies that $2n+i=2l+2$. By induction hypothesis, $p^{2n}\preccurlyeq p^{2n+(i-2)}=p^{2l}$. Applying (a--$\ast$), we receive that $p^{2l}\preccurlyeq p^{2l+2}= p^{2n+i}$.
	Then, we apply Lemma~\ref{L:relations-general-properties}--k.\\
	
	\noindent\textit{Proof of}~(b). Since $\LTInt(1)$ is a Heyting algebra (Corollary~\ref{C:LT_Int}) and by definition (\ref{E:powers_p^n}--$ii$), $p^{2n+1}\preccurlyeq p^{2n+1}\lor p^{2n+2}=p^{2n+4}$. Further, in virtue of (a),
	$p^{2n+4}=2^{2(n+2)}\preccurlyeq 2^{2(n+2)+i}$. Then, we apply Lemma~\ref{L:relations-general-properties}--a.\\
	
	\noindent\textit{Proof of}~(c). Considering two cases, whether $n$ is even or odd, we get: $p^{2l}\preccurlyeq p^{2l+3+i}$, by (a), and $p^{2l+1}\preccurlyeq
	p^{2l+1+3+i}$, by (b).\\
	
	\noindent\textit{Proof of}~(d). We successively obtain:
	\[
	\begin{array}{rl}
		p^{2n+6}=p^{2(n+1)+4}=p^{2(n+1)+1}\lor p^{2^{2(n+1)+2}} 
		&[\text{definition (\ref{E:powers_p^n}--$ii$)}]\\
		=p^{2n+3}\lor p^{2n+4}=p^{2n+3}\lor(p^{2n+1}\lor p^{2n+2})
		&[\text{definition (\ref{E:powers_p^n}--$ii$)}]\\
		\sim p^{2n+1}\lor(p^{2n+2}\lor p^{2n+3}) &[\text{associative law for $\lor$ in $\LTInt(1)$}]\\
		\sim p^{2n+1}\lor p^{2n+3} &[\text{by (a)}]. 
	\end{array}
	\]
	
	\noindent\textit{Proof of}~(e). Since $\LTInt(1)$ is a Heyting algebra (Corollary~\ref{C:LT_Int}) and by definition, we have: $p^{2n+1}\land p^{2n+3}\sim p^{2n+1}\land p^{2n}$. Then, we apply (a--$\ast\ast\ast$).\\
	
	\noindent\textit{Proof of}~(f). First, we notice that, in view of (a),
	\[
	p^{2n}\preccurlyeq p^{2n+1}\land p^{2n+2}.
	\]
	
	Then, in virtue of Lemma~\ref{L:Int-theses}--c, 
	\[
	(p^{2n+1}\land p^{2n+2})\rightarrow p^{2n}\sim p^{2n+2}\rightarrow(p^{2n+1}\rightarrow p^{2n})
	\sim p^{2n+2}\rightarrow p^{2n+3}.
	\]
	Applying (a), we conclude that $p^{2n+1}\land p^{2n+2}\preccurlyeq p^{2n}$.\\
	
	\noindent\textit{Proof of}~(g). Using definition (\ref{E:powers_p^n}--$i$) and Lemma~\ref{L:Int-theses}--c, we have:
	\[
	p^{2n+1}\rightarrow p^{2n+3}=p^{2n+1}\rightarrow(p^{2n+1}\rightarrow p^{2n})
	\sim (p^{2n+1}\land p^{2n+1})\rightarrow p^{2n}\sim (p^{2n+1})\rightarrow p^{2n}=p^{2n+3}.
	\]
	
	\noindent\textit{Proof of}~(h). If $n=0$, according to Lemma~\ref{L:Int-theses}--i, we have:
	\[
	p^{3}\rightarrow p^{1}\sim\neg\neg p\rightarrow\neg p\sim\neg p=p^1.
	\]
	
	And for any $n=l+1$, applying Lemma~\ref{L:Int-theses}--c, (a--$\ast\ast$) and definition (\ref{E:powers_p^n}--$i$),
	we have:
	\[
	\begin{array}{r}
		p^{2(l+1)+3}\rightarrow p^{2(l+1)+1}\sim p^{2l+5}\rightarrow(p^{2l+1}\rightarrow p^{2l})\sim 
		(p^{2l+5}\land p^{2l+1})\rightarrow p^{2l}\\
		\sim p^{2l+1}\rightarrow p^{2l}
		\sim p^{2l+3}.
	\end{array}
	\]
	
	\noindent\textit{Proof of}~(i). On the one hand, in virtue of (\ref{E:powers_p^n}--$i$), (a--$\ast$) and Lemma~\ref{L:Int-theses}--m, we have
	\[
	p^{2n+3}\sim (p^{2n+1}\rightarrow p^{2n})\land (p^{2n}\rightarrow p^{2n+2})
	\preccurlyeq p^{2n+1}\rightarrow  p^{2n+2}.
	\]
	
	On the other hand, in virtue of Lemma~\ref{L:Int-theses}--n and (f), we obtain:
	\[
	p^{2n+1}\land(p^{2n+1}\rightarrow p^{2n+2})\sim p^{2n+1}\land p^{2n+2}\sim p^{2n}.
	\]
	In particular,
	\[
	(p^{2n+1}\rightarrow p^{2n+2})\land p^{2n+1}\preccurlyeq p^{2n},
	\]
	which, according to Lemma~\ref{L:relations-general-properties}--l, is equivalent to
	\[
	p^{2n+1}\rightarrow p^{2n+2}\preccurlyeq\underbrace{p^{2n+1}\rightarrow p^{2n}}_{p^{2n+3}}.
	\]
	
	\noindent\textit{Proof of}~(j). If $n=0$, then, using Lemma~\ref{L:Int-theses}--f, we have:
	\[
	p^2\rightarrow p^1\sim\neg p=p^1.
	\]
	
	Now let $n=l+1$. Then we obtain:
	\[
	\begin{array}{rl}
		p^{2l+4}\rightarrow p^{2l+3}=(p^{2l+1}\lor p^{2l+2})\rightarrow p^{2l+3}
		&[\text{by definition (\ref{E:powers_p^n}--$ii$)}]\\
		\sim (p^{2l+1}\rightarrow p^{2l+3})\land(p^{2l+2}\rightarrow p^{2l+3})
		&[\text{Lemma~\ref{L:Int-theses}--d}]\\
		\sim p^{2l+1}\rightarrow p^{2l+3} &[\text{by (a--$\ast\ast\ast$)}]\\
		\sim p^{2l+1}\rightarrow (p^{2l+1}\rightarrow p^{2l}) 
		&[\text{by definition (\ref{E:powers_p^n}--$i$)}]\\
		\sim p^{2l+1}\rightarrow p^{2l} &[\text{Lemma~\ref{L:Int-theses}--c}]\\
		\sim p^{2l+3}=p^{2n+1} &[\text{definition (\ref{E:powers_p^n}--$i$)}].
	\end{array}
	\]
	
	\noindent\textit{Proof of}~(k). If $n=0$, according to Lemma~\ref{L:Int-theses}--j, we have:
	\[
	p^2\rightarrow p^0\sim\neg p=p^1.
	\]
	
	Now assume that $n=l+1$. Then we obtain:
	\[
	\begin{array}{rl}
		p^{2l+4}\rightarrow p^{2l+2}\sim (p^{2l+2}\lor p^{2l+1})\rightarrow p^{2l+2}
		&[\text{definition (\ref{E:powers_p^n}--$ii$)}]\\
		\sim (p^{2l+2}\rightarrow p^{2l+2})\land(p^{2l+1}\rightarrow p^{2l+2})
		&[\text{Lemma~\ref{L:Int-theses}--d}]\\
		\sim p^{2l+1}\rightarrow p^{2l+2} &[\text{($\ast$) in the proof of Lemma~\ref{L:deduction-theorem}}]\\
		\sim p^{2l+3} &[\text{by (i)}].
	\end{array}
	\]
	
	\noindent\textit{Proof of}~(l). We distinguish two cases: $k$ is even and $k$ is odd.
	
	Suppose $k=2l$; then $n=2l+4+i$. We notice that, on the one hand, according to (a), $p^{2l}\preccurlyeq p^{2l+2}$. On the other hand, we obtain:
	\[
	\begin{array}{rl}
		p^{2l+4}\rightarrow p^{2l}=(p^{2l+1}\lor p^{2l+2})\rightarrow p^{2l}&[\text{by definition (\ref{E:powers_p^n}--$ii$)}]\\
		\sim (p^{2l+1}\rightarrow p^{2l})\land(p^{2l+2}\rightarrow p^{2l})
		&[\text{Lemma~\ref{L:Int-theses}--d}]\\
		=p^{2l+3}\land (p^{2l+2}\rightarrow p^{2l}) &[\text{by definition (\ref{E:powers_p^n}--$ii$)}]\\
		\sim p^{2l+3}\land p^{2l+1} &[\text{by (k)}]\\
		\sim p^{2l} &[\text{by (e)}].
	\end{array}
	\]
	
	Applying Lemma~\ref{L:relations-general-properties}--d, we receive that $p^{2l}\ll p^{2l+4}$. And since, in virtue of (a), $p^{2l+4}\preccurlyeq p^{2l+4+i}$, we, applying Lemma~\ref{L:relations-general-properties}--f, derive that $p^{2l}\ll p^{2l+4+i}$.
	
	Now assume that $k=2l+1$; then $n=2l+5+i$. First, we prove that $p^{2l+1}\ll p^{2l+4}$.
	
	Indeed, on the one hand, in virtue of definition (\ref{E:powers_p^n}--$ii$), $p^{2l+1}\preccurlyeq p^{2l+1}\lor p^{2l+2}=p^{2l+4}$. On the other hand, we have:
	\[
	\begin{array}{rl}
		p^{2l+4}\rightarrow p^{2l+1}=(p^{2l+1}\lor p^{2l+2})\rightarrow p^{2l+1}&[\text{by definition (\ref{E:powers_p^n}--$ii$)}]\\
		\sim (p^{2l+1}\rightarrow p^{2l+1})\land(p^{2l+2}\rightarrow p^{2l+1})
		&[\text{Lemma~\ref{L:Int-theses}--d}]\\
		\sim p^{2l+2}\rightarrow p^{2l+1}\sim p^{2l+1} &[\text{by (j)}].
	\end{array}
	\]
	
	Then, applying Lemma~\ref{L:relations-general-properties}--d, we get that
	$p^{2l+1}\ll p^{2l+4}$. 
	
	Now, since, in virtue of (a), $p^{2l+4}\preccurlyeq p^{2l+4+i}$, according to
	Lemma~\ref{L:relations-general-properties}--f, we conclude that
	$p^{2l+1}\ll p^{2l+4+i}$, in particular, $p^{2l+1}\ll p^{2l+5+i}$.\\
	
	\noindent\textit{Proof of}~(m). If $n=0$, then we have:
	\[
	(p^1\rightarrow p^0)\rightarrow p^0\sim \neg\neg\neg p\sim \neg p
	\]
	(Lemma~\ref{L:Int-theses}, (g) twice and (h)).
	
	Next we take $n=1$. Then we have:
	\[
	(p^3\rightarrow p^2)\rightarrow p^2\sim (\neg\neg p\rightarrow p)\rightarrow p. 
	\]
	
	First, we prove that
	\[
	(\neg\neg p\rightarrow p)\rightarrow p\vdash_{\textsf{Int}}\neg\neg p\tag{m--$\ast$}
	\]
	and, then, it will remain to apply the hyperrule (c--$i$) (Lemma~\ref{L:deduction-theorem}).
	
	We have:
	\[
	\begin{array}{cl}
		1. &(\neg\neg p\rightarrow p)\rightarrow p\quad[\text{premise}]\\
		2. &\neg p\quad[\text{assumption}]\\
		3. &\neg p\rightarrow(\neg\neg p\rightarrow p)\quad[\text{an instance of ax11 (Section~\ref{section:modus-rules})}]\\
		4. &\neg\neg p\rightarrow p\quad[\text{from 2 and 3 by (b--$iii$)}]\\
		5. &p\quad[\text{from 1 and 4 by (b--$iii$)}].
	\end{array}
	\]
	Then, we apply (c--$ii$) to obtain (m--$\ast$) (Lemma~\ref{L:hyperrules-c_ii-c_iii}).
	
	On the other hand, applying the rule (b--$iii$) and hyperrule (c--$i$), we derive:
	\[
	\neg\neg p\vdash_{\textsf{Int}}(\neg\neg p\rightarrow p)\rightarrow p.
	\]
	
	Next we take $n=l+2$ and prove that
	\[
	(\underbrace{p^{2(l+2)+1}}_{p^{2(l+1)+3}}\rightarrow p^{2(l+2)})\rightarrow p^{2(l+2)}\sim
	p^{2(l+2)+1},
	\]
	that is,
	\[
	((p^{2(l+1)+1}\rightarrow p^{2(l+1)})\rightarrow p^{2(l+2)})\rightarrow p^{2(l+2)}\sim p^{2(l+1)+1}\rightarrow p^{2(l+1)}. 
	\]
	
	First, we prove that
	\[
	((p^{2(l+1)+1}\rightarrow p^{2(l+1)})\rightarrow p^{2(l+2)})\rightarrow p^{2(l+2)}, p^{2(l+1)+1}\vdash_{\textsf{Int}}p^{2(l+1)}\tag{m--$\ast\ast$}
	\]
	and then apply the hyperrule (c--$i$).
	
	To prove (m--$\ast\ast$) we first obtain:
	\[
	\begin{array}{cl}
		1. &((p^{2(l+1)+1}\rightarrow p^{2(l+1)})\rightarrow p^{2(l+2)})\rightarrow p^{2(l+2)}\quad[\text{premise}]\\
		2. & p^{2(l+1)+1}\quad[\text{premise}]\\
		3. & p^{2(l+1)+1}\rightarrow p^{2(l+1)}\quad[\text{assumption}]\\
		4. & p^{2(l+1)}\quad[\text{from 2 and 3 by (b--$iii$)}]\\
		5. & \underbrace{p^{2l+1}\lor p^{2(l+1)}}_{p^{2(l+2)}}\quad[\text{from 3 by (a--$iii$); see (Lemma~\ref{L:rules-sound-wrt-P-Int-Cl})}] .\\
	\end{array}
	\]
	
	Applying the hyperrule (c--$i$), we obtain:
	\[
	((p^{2(l+1)+1}\rightarrow p^{2(l+1)})\rightarrow p^{2(l+2)})\rightarrow p^{2(l+2)}, p^{2(l+1)+1}\vdash_{\textsf{Int}}(p^{2(l+1)+1}\rightarrow p^{2(l+1)})\rightarrow p^{2(l+2)}
	\]
	and hence
	\[
	((p^{2(l+1)+1}\rightarrow p^{2(l+1)})\rightarrow p^{2(l+2)})\rightarrow p^{2(l+2)}, \underbrace{p^{2(l+1)+1}}_{p^{2l+1}\rightarrow p^{2l}}\vdash_{\textsf{Int}}\underbrace{p^{2(l+2)}}_{p^{2l+1}\lor p^{2(l+1)}}. \tag{m--$\ast\ast\ast$}
	\]
	We notice that $p^{2l+1}, p^{2(l+1)+1}\vdash_{\textsf{Int}} p^{2l}$ and hence, in virtue of (a),
	\[
	p^{2l+1}, p^{2(l+1)+1}\vdash_{\textsf{Int}} p^{2(l+1)}.
	\]
	On the other hand, $p^{2(l+1)}\vdash_{\textsf{Int}}p^{2(l+1)}$. Applying the hyperrule (c--$iii$), we conclude that
	\[
	p^{2(l+1)+1}, p^{2(l+2)}\vdash_{\textsf{Int}}p^{2(l+1)}.
	\]
	Using this and (m--$\ast\ast$), we obtain (m--$\ast\ast\ast$).
	
	The converse, that is
	\[
	p^{2n+1}, p^{2n+1}\rightarrow p^{2n}\vdash_{\textsf{Int}}p^{2n},
	\]
	is obvious. Two applications of the hyperrule (--$i$) leads to the conclusion:
	$p^{2n+1}\rightarrow((p^{2n+1}\rightarrow p^{2n})\rightarrow p^{2n})$ is an $\Int$-thesis.
\end{proof}
\begin{prop}\label{P:nishimura-completeness} 
	For any $n,k\le\omega$, there is $m\le\omega$, individually its own for each of the statements {\em (a)--(d)} below, such that
	{\em\[
		\begin{array}{cl}
			(\text{a}) & \neg(p^n)\sim p^m;\\
			(\text{b}) & p^n\land p^k\sim p^m;\\
			(\text{c}) & p^n\lor p^k\sim p^m\\
			(\text{d}) & p^n\rightarrow p^k\sim p^m.\\
		\end{array}
		\]}
\end{prop}
\begin{proof}
	Case (a) follows from (d) because of Lemma~\ref{L:Int-theses}--j and Lemma~\ref{L:relations-general-properties}--h.\\
	
	Case (b). Suppose that $n=k$, then, obviously, $p^n\land p^k\sim p^n$. Assume that $n<k$. If $n$ is even, then, in virtue of Lemma~\ref{L:relations-specific-properties}--a, $p^n\preccurlyeq p^k$ and, hence (Lemma~\ref{L:relations-general-properties}--i), $p^n\land p^k\sim p^n$.
	
	Now assume that $n=2l+1$. Then $p^n\land p^{n+1}\sim p^{2l}$ (Lemma~\ref{L:relations-specific-properties}--f), $p^n\land p^{n+2}\sim p^{2l}$
	(Lemma~\ref{L:relations-specific-properties}--e) and $p^n\land p^{n+3+i}\sim p^n$ (Lemma~\ref{L:relations-specific-properties}--c).\\
	
	Case (c). If $n=k$, then, obviously, $p^n\lor p^k\sim p^n$. Thus we assume that $n<k$. If $n$ is even, then, in virtue of Lemma~\ref{L:relations-specific-properties}--a and Lemma~\ref{L:relations-general-properties}--j, $p^n\lor p^k\sim p^k$.
	
	Let $n=2l+1$. Then, by definition (\ref{E:powers_p^n}--$ii$), $p^n\lor p^{n+1}\sim p^{n+3}$. Also, $p^n\lor p^{n+2}\sim p^{n+5}$. Indeed, since $p^{2l+2}\preccurlyeq p^{2l+3}$ (Lemma~\ref{L:relations-specific-properties}--a), $p^{2l+1}\lor p^{2l+3}\sim
	p^{2l+1}\lor (p^{2l+2}\lor p^{2l+3})$. However, in virtue of lattice laws true for $\LTInt(1)$, $p^{2l+1}\lor (p^{2l+2}\lor p^{2l+3})\sim(p^{2l+1}\lor p^{2l+2})\lor p^{2l+3}$. Hence, by definition (\ref{E:powers_p^n}--$ii$),
	$p^{2l+1}\lor p^{2l+3}\sim p^{2l+4}\lor p^{2l+3}$. However, also by definition (\ref{E:powers_p^n}--$ii$),
	$p^{2l+3}\lor p^{2l+4}=p^{2l+6}$.
	
	Finally (in this case), in virtue of Lemma~\ref{L:relations-specific-properties}--c and Lemma~\ref{L:relations-general-properties}--j, $p^{n}\lor p^{k}\sim p^{k}$, for any $k\ge n+3$.\\
	
	Case (d). The cases when $n=k$ or $n=0$ or $n=\omega$ or $k=\omega$ are obvious. So, first we suppose that 
	$0<n<k<\omega$.
	
	If $n$ is even, then $p^n\rightarrow p^k\sim p^{\omega}$ (Lemma~\ref{L:relations-specific-properties}--a).
	
	Now let $n$ be odd, that is, $n=2l+1$, for some $l\ge 0$. Then we observe:
	
	\[
	\begin{array}{rl}
		p^{2l+1}\rightarrow p^{2l+2}\sim p^{2l} &[\text{Lemma~\ref{L:relations-specific-properties}--h}];\\
		p^{2l+1}\rightarrow p^{2l+3}\sim p^{2l+3} &[\text{Lemma~\ref{L:relations-specific-properties}--g}];\\
		p^{2l+1}\rightarrow p^{2l+4+i}\sim p^{\omega} &[\text{Lemma~\ref{L:relations-specific-properties}--b}].\\
	\end{array}
	\]
	
	Now we turn to the case when $0\le k<n<\omega$. We will employ course-of-value induction on $n$. As the base of induction, we observe the following specifications:
	\[
	\begin{array}{rl}
		p^1\rightarrow p^0=p^3 &[\text{by definition}];\\
		p^2\rightarrow p^0\sim p^1 &[\text{Lemma~\ref{L:Int-theses}--j}];\\
		p^2\rightarrow p^1\sim p^1&\text{[Lemma~\ref{L:relations-specific-properties}--k]};\\
		p^3\rightarrow p^0=(p^1\rightarrow p^0)\rightarrow p^0\sim p^1
		&[\text{Lemma~\ref{L:relations-specific-properties}--m}];\\
		p^3\rightarrow p^1\sim p^1 &[\text{Lemma~\ref{L:relations-specific-properties}--h}];\\
		p^3\rightarrow p^2=p^5 &\text{[by definition (\ref{E:powers_p^n}--$i$)]}.\\
	\end{array}
	\]
	
	Next, assume that $n=2l+4$. Then we have:
	\[
	\begin{array}{rl}
		p^{2l+4}\rightarrow p^k\sim (p^{2l+1}\lor p^{2l+2})\rightarrow p^k
		&[\text{by definition (\ref{E:powers_p^n}--$ii$)}]\\
		\sim(p^{2l+1}\rightarrow p^k)\land(p^{2l+2}\rightarrow p^k)
		&[\text{Lemma~\ref{L:Int-theses}--d}].\\
	\end{array}
	\]
	One of the following cases are possible: $k\ge 2l+2$, $2l+1\le k<2l+2$ and $k<2l+1$. In any event, either by the case which have been considered above when $n\le k$ or by induction, there are $m_1$ and $m_2$ such that
	\[
	p^{2l+4}\rightarrow p^k\sim p^{m_1}\land p^{m_2}.
	\]
	This, according to the case (b) above, implies that there is an $m$ such that
	\[
	p^{2l+4}\rightarrow p^k\sim p^{m}.
	\]
	
	Now suppose that $n=2l+5$ and $0\le k\le 2l+1$.
	
	If $k=2l+4$, then, by definition, $p^n\rightarrow p^k\sim p^{2(l+2)+3}$. If $k=2l+3$, then $p^{n}\rightarrow p^{k}=p^{2(l+1)+3}\rightarrow p^{2(l+1)+1}\sim
	p^{2(l+1)+1}$ (Lemma~\ref{L:relations-specific-properties}--h). And if $k=2l+2$, then $p^{n}\rightarrow p^{k}=p^{2(l+1)+3}\rightarrow p^{2(l+1)}=
	(p^{2(l+1)+1}\rightarrow p^{2(l+1)})\rightarrow p^{2(l+1)}\sim p^{2(l+1)+1}$
	(Lemma~\ref{L:relations-specific-properties}--m).
	
	Thus it remains to consider the case when $n=2l+5$ and $k\le 2l+1$. But then $n-k\ge 4$ and, in virtue Lemma~\ref{L:relations-specific-properties}--l, we conclude that $p^k\ll p^n$.  This, in virtue of Lemma~\ref{L:relations-general-properties}--d, implies that $p^n\rightarrow p^k\sim p^k$.
\end{proof}
\begin{cor}\label{C:nishimura-completeness}
	For any formula $\alpha(p)$ with a single variable $p$, there is $n\le\omega$ such that $\alpha(p)\sim p^n$.
\end{cor}
\begin{proof}
	Using Proposition~\ref{P:nishimura-completeness}, the proof can be carried out by induction on the number of logical connectives occurring in $\alpha(p)$. We leave this to the reader as an exercise. (Exercise~\ref{section:lindenbaum-algebra}.\ref{EX:nishimura-completeness})	
\end{proof}

Now we turn to the other path of our argument. Since to prove Proposition~\ref{P:nishimura-power-difference}, we have to separate the set of $\Int$-theses from the formulas $p^n$, for any $n<\omega$. If some initial powers of $p$ can be refuted in relatively simple Heyting algebras, in order to carry out proof by induction, we will need a o more transparent presentation of $\Int$-theses. $\Mapsto$\\

For our purposes, we borrow a sequential calculus G3 from~\cite{kleene1952}, chapter XV. The main formal objects are \textit{\textbf{sequents}} which are pairs of finite sets of $\Lan_A$-formulas, say $\Gamma$ and $\Theta$, divided by the symbol $\Mapsto$ --- $\Gamma\Mapsto\Theta$, where $\Theta$ is allowed to consist of not more than one formula. The system G3 is formulated in terms of metavariables for $\Lan_A$-formulas and variables $\Gamma$ and $\Theta$ for finite sets of formulas (with the aforementioned restriction on $\Theta$) as follows.\\

Axiom schema:~ $\Gamma,\bm{\alpha}\Mapsto\bm{\alpha}$.\\

Rules of inference:
\begin{equation*}
	\begin{array}{llrl}
		(\text{$\rightarrow$-1}) &\dfrac{\bm{\alpha},\Gamma\Mapsto\bm{\beta}}{\Gamma\Mapsto\bm{\alpha}\rightarrow\bm{\beta}} &\quad(\text{$\rightarrow$-2}) &\dfrac{\bm{\alpha}\rightarrow\bm{\beta},\Gamma\Mapsto\bm{\alpha}~\text{and}~\bm{\beta},\bm{\alpha}\rightarrow\bm{\beta},\Gamma\Mapsto\Theta}{\bm{\alpha}\rightarrow\bm{\beta},\Gamma\Mapsto\Theta}\\\\
		(\text{$\land$-1}) &\dfrac{\Gamma\Mapsto\bm{\alpha}~\text{and}~\Gamma\Mapsto\bm{\beta}}{\Gamma\Mapsto\bm{\alpha}\land\bm{\beta}} &\quad(\text{$\land$-2}) &\dfrac{\bm{\alpha},\bm{\alpha}\land\bm{\beta},\Gamma\Mapsto\Theta~\text{or}~\bm{\beta},\bm{\alpha}\land\bm{\beta},\Gamma\Mapsto\Theta}{\bm{\alpha}\land\bm{\beta},\Gamma\Mapsto\Theta}\\\\
		(\text{$\lor$-1}) &\dfrac{\Gamma\Mapsto\bm{\alpha}~\text{or}~\Gamma\Mapsto\bm{\beta}}{\Gamma\Mapsto\bm{\alpha}\lor\bm{\beta}} &\quad(\text{$\lor$-2}) &\dfrac{\bm{\alpha},\bm{\alpha}\lor\bm{\beta},\Gamma\Mapsto\Theta~\text{and}~\bm{\beta},\bm{\alpha}\lor\bm{\beta},\Gamma\Mapsto\Theta}{\bm{\alpha}\lor\bm{\beta},\Gamma\Mapsto\Theta}\\\\
		(\text{$\neg$-1}) &\dfrac{\bm{\alpha},\Gamma\Mapsto }{\Gamma\Mapsto\neg\bm{\alpha}} &\quad(\text{$\neg$-2}) &\dfrac{\neg\bm{\alpha},\Gamma\Mapsto\bm{\alpha}}{\neg\bm{\alpha},\Gamma\Mapsto\Theta}
	\end{array}
\end{equation*}

A derivation in G3 is defined in the form of a tree where each node is labeled with a sequent; the top sequents of the tree are instances of the axiom schema and the root is a sequent that is intended to be proved. (Cf.~\cite{kleene1952}, the end of {\S} 24.) The following proposition we take for granted.
\begin{prop}[\cite{kleene1952}, {\S} 80]\label{P:G3-completeness}
	For any $\Lan_A$-formula $\alpha$, {\em$\vdash_{\textsf{Int}}\alpha$} if, and only if, the sequent $\Mapsto\alpha$ is derived in {\em G3}.
\end{prop}

We use the last proposition to prove the following lemma.
\begin{lem}\label{L:G3-lemma}
	For any $\Lan_A$-formulas $\alpha$, $\beta$ and $\gamma$, if {\em$\vdash_{\textsf{Int}}(\alpha\rightarrow\beta)\rightarrow(\alpha\lor\gamma)$}, then either {\em$\vdash_{\textsf{Int}}(\alpha\rightarrow\beta)\rightarrow\alpha$} or
	{\em$\vdash_{\textsf{Int}}(\alpha\rightarrow\beta)\rightarrow\gamma$}.
\end{lem}
\begin{proof}
	Suppose that $\vdash_{\textsf{Int}}(\alpha\rightarrow\beta)\rightarrow(\alpha\lor\gamma)$. This, in virtue of Proposition~\ref{P:G3-completeness}, means that the sequent $\Mapsto(\alpha\rightarrow\beta)\rightarrow(\alpha\lor\gamma)$ is derivable in G3. According to the rules of G3, this sequent can only be obtained from the sequent $\alpha\rightarrow\beta\Mapsto\alpha\lor\gamma$ by the rule ($\rightarrow$-1). Regarding the last sequent, we have several options. Either it is obtained from the sequent $\alpha\rightarrow\beta\Mapsto\alpha$	or from the sequent $\alpha\rightarrow\beta\Mapsto\gamma$ by one of the two variants of the rule ($\lor$-1). In event of each of theses cases, we arrive at a desirable conclusion, applying the rule ($\rightarrow$-1) and then using Proposition~\ref{P:G3-completeness}.
	
	Now we consider the case: $\alpha\rightarrow\beta\Mapsto\alpha\lor\gamma$ is obtained from $\alpha\rightarrow\beta\Mapsto\alpha$ and $\beta,\alpha\rightarrow\beta\Mapsto\alpha\lor\gamma$ by the rule ($\rightarrow$-2). Thus we have that it is necessary that $\vdash_{\textsf{Int}}(\alpha\rightarrow\beta)\rightarrow\alpha$ (Proposition~\ref{P:G3-completeness}).
\end{proof}
\begin{cor}[disjunction property for $\Int$]\label{C:disjuction-property-Int}
	If ~{\em$\vdash_{\textsf{Int}}\alpha\lor\beta$}, then either {\em$\vdash_{\textsf{Int}}\alpha$} or {\em$\vdash_{\textsf{Int}}\beta$}.	
\end{cor}
\begin{proof}
	Let $\star$ be an arbitrary $\Int$-thesis, that is, $\vdash_{\textsf{Int}}\star$. Assume that $\vdash_{\textsf{Int}}\alpha\lor\beta$. Then both $\vdash_{\textsf{Int}}\alpha\rightarrow\star$ and $\vdash_{\textsf{Int}}(\alpha\rightarrow\star)\rightarrow(\alpha\lor\beta)$ hold. In virtue of Lemma~\ref{L:G3-lemma}, the latter implies that either 	$\vdash_{\textsf{Int}}(\alpha\rightarrow\star)\rightarrow\alpha$ or
	$(\alpha\rightarrow\star)\rightarrow\beta$, which in turn implies that either
	$\vdash_{\textsf{Int}}\alpha$ or $\vdash_{\textsf{Int}}\beta$.
\end{proof}

\begin{lem}\label{L:towards-p^n-unprovability}
	For every natural number $n$, if~~{\em$\vdash_{\textsf{Int}}p^{2n+7}$}, then either~~{\em$\vdash_{\textsf{Int}}p^{2n+3}$} or~~{\em$\vdash_{\textsf{Int}}p^{2n+1}$}.	
\end{lem}
\begin{proof}
	Suppose that $\vdash_{\textsf{Int}}p^{2n+7}$. According to Lemma~\ref{L:relations-specific-properties}, (i) and (d), we respectively have
	$\vdash_{\textsf{Int}}p^{2n+5}\rightarrow p^{2n+6}$ and $\vdash_{\textsf{Int}}p^{2n+5}\rightarrow(p^{2n+1}\lor p^{2n+3})$. Then, by definition, we have: $\vdash_{\textsf{Int}}(p^{2n+3}\rightarrow p^{2n+2})\rightarrow(p^{2n+1}\lor p^{2n+3})$. By Lemma~\ref{L:towards-p^n-unprovability}, this implies that either
	$\vdash_{\textsf{Int}}(p^{2n+3}\rightarrow p^{2n+2})\rightarrow p^{2n+3}$
	or $\vdash_{\textsf{Int}}(p^{2n+3}\rightarrow p^{2n+2})\rightarrow p^{2n+1}$; that is, by definition, $\vdash_{\textsf{Int}}p^{2n+5}\rightarrow p^{2n+3}$ or
	$\vdash_{\textsf{Int}}p^{2n+5}\rightarrow p^{2n+1}$. Applying Lemma~\ref{L:relations-specific-properties}--h, we obtain that either
	$\vdash_{\textsf{Int}}p^{2n+3}$ or $\vdash_{\textsf{Int}}p^{2n+1}$.
\end{proof}

\begin{prop}\label{P:nishimura-power-difference}
	For every $n<\omega$, $p^n$ is not an {\em$\Int$}-thesis; in other words, {\em$\not\vdash_{\textsf{Int}} p^n$}, for any natural number $n$.
\end{prop}
\begin{proof}
	We prove by induction on $n$. First, we notice that neither of the following formulas is an $\Int$-thesis:
	\[
	\underbrace{p\land\neg p}_{p^0},~\underbrace{\neg p}_{p^1},~\underbrace{p}_{p^2},~\underbrace{\neg\neg p}_{p^3},~\underbrace{p\lor\neg p}_{p^4},
	~\underbrace{\neg\neg p\rightarrow p}_{p^5},
	\]
	which are equal or equivalent (in the sense of $\sim$) respectively to the formulas $p^0,~~p^1,~~p^2,~~p^3,~~p^4,~~p^5$.
	
	Indeed, the first four are even not classical tautologies, since can be refuted in $\booleTwo$, and the last two, though are classical tautologies, are refuted in the Heyting algebra $\godelThree$.
	
	Further, for contradiction, we assume that $\vdash_{\textsf{Int}}p^n$, where $n\ge 6$. Next we consider the following two cases.
	
	Case: $n$ is even, that is, $n=2l+6$. By definition, we have that $\vdash_{\textsf{Int}}p^{2l+4}\lor p^{2l+1}$. This, in virtue of Corollary~\ref{C:disjuction-property-Int}, either $\vdash_{\textsf{Int}}p^{2l+4}$ or $p^{2l+1}$. Both options contradict the induction hypothesis.
	
	Case: $n$ is odd, that is, $n=p^{2l+7}$. According to Lemma~\ref{L:towards-p^n-unprovability}, the assumption that $\vdash_{\textsf{Int}}p^{2l+7}$ implies that either $\vdash_{\textsf{Int}}p^{2l+3}$ or $\vdash_{\textsf{Int}}p^{2l+1}$. Again, both options contradict the induction hypothesis.
	
	Therefore, $\not\vdash_{\textsf{Int}}p^n$.	
\end{proof}
\begin{rem}
	{\em In the light of Proposition~\ref{P:nishimura-power-difference}, the conditional in Lemma~\ref{L:towards-p^n-unprovability} should be understood as counterfactual rather than implication. That is, it would be more appropriate to formulas it as follows: If it were the case that $\vdash_{\textsf{Int}}p^{2n+7}$, then it would be that either $\vdash_{\textsf{Int}}p^{2n+3}$ or $\vdash_{\textsf{Int}}p^{2n+1}$.} 	
\end{rem}

\begin{cor}\label{C:nishimura-power-difference}
	For any $n,k\le\omega$,
	\[
	p^n\sim p^k~\Longleftrightarrow~n=k.
	\]
\end{cor}
\begin{proof}
	The $\Leftarrow$-implication is obvious. Now, for definiteness, suppose that $n<k$. If it were that $p^n\sim p^k$, then it would imply that $\vdash_{\textsf{Int}} p^k\rightarrow p^n$. This in turn would imply that there is a natural number $m$ such that $\vdash_{\textsf{Int}} p^m$. (Proposition~\ref{P:nishimura-completeness}--d). A contradiction.	
\end{proof}

\paragraph{Exercises~\ref{section:lindenbaum-algebra}}
\begin{enumerate}
	\item \label{EX:Cl-Int-implicational} Show that the abstract logics $\Pos$, $\Cl$ and $\Int$ are implicational with respect to $p\rightarrow q$.
	\item \label{EX:T_S-star-thesis}Prove~\eqref{E:T_S-star-thesis}.
	\item \label{EX:v-restored}Prove that for any valuation $w$ in $\langle\alg{A},\one\rangle$, there is a unique valuation $v$ in $\alg{A}$ such that $v^{\star}=w$.
	\item \label{EX:LT-term-valuation}Prove~\eqref{E:LT-term-valuation}.
	\item \label{EX:LT=free-algebra} Prove Lemma~\ref{L:LT=free-algebra}.
	\item \label{EX:lemma-rules-sound-wrt-P-Int-Cl}Complete the proof of Lemma~\ref{L:rules-sound-wrt-P-Int-Cl}.
	\item \label{EX:lemma-hyperrules-c_ii-c_iii}Show that the hyperrules (c--$ii$) and (c--$iii$) are sound with respect to {\Cl} and {\Int}.
	\item \label{EX:distributive-lettice-sound} Complete the proof of Lemma~\ref{L:distributive-lettice-sound}.
	\item\label{EX:Cl-implications} Complete the proof of Lemma~\ref{L:Cl-implications}.
	\item \label{EX:four-elements} Show that the elements $[p]$, $[\neg p]$, $[p\vee\neg p]$ and $[p\wedge\neg p]$ are pairwise distinct in $\LTCl$.
	(Hint: prove that the matrix $\booleTwo$ is a $\Cl$-model.)
	\item\label{EX:Int-implications} Complete the proof of Lemma~\ref{L:Int-implications}.
	\item \label{EX:Int-theses} Complete the proof of Lemma~\ref{L:Int-theses}. (Hint: consult~\cite{kleene1952}, {\S} 26 and {\S} 27 if necessary.)
	\item \label{EX:relations-general-properties} Complete the proof of Lemma~\ref{L:relations-general-properties}.
	\item\label{EX:nishimura-completeness} Prove Corollary~\ref{C:nishimura-completeness}.
\end{enumerate}

\section{Examples of application of Lindenbaum-Tarski algebras}
In this section, we discuss several examples where the notion of Lindenbaum-Tarski algebra plays a key role. 

\subsection{Glivenko's theorem}
The following theorem was established by Glivenko~\cite{glivenko1929}: For any $\Lan_A$-formula $\alpha$,
\[
\alpha\in\bm{T}_{\Cl}~\Longleftrightarrow~\neg\neg\alpha\in\bm{T}_{\Int}.\tag{\textit{Glivenko's theorem}}
\]
\noindent\textit{Sketch of proof}. Let us fix an $\Lan_A$-formula $\alpha$.

First, we notice that $\neg\neg\alpha\rightarrow\alpha\in\bm{T}_{\Cl}$ and $\bm{T}_{\Int}\subseteq\bm{T}_{\Cl}$. Taking into account that $\bm{T}_{\Cl}$ is closed under detachment rule, we obtain the implication:
\[
\neg\neg\alpha\in\bm{T}_{\Int}~\Longrightarrow~\alpha\in\bm{T}_{\Cl}.
\]

Conversely, assume that $\alpha\in\bm{T}_{\Cl}$. According to the second part of Proposition~\ref{P:valuations-in-LT}, $\alpha$ is valid in $\LTCl$ and hence (Corollary~\ref{C:LT_Cl(k)}) $\alpha$ is valid in any Boolean algebra.

Now, assume that $\neg\neg\alpha\notin\bm{T}_{\Int}$. This implies that $v[\neg\neg\alpha]$ is invalid in $\LTInt$ at $v$, where $v$ is the Lindenbaum valuation. This in turn implies that the element $v[\alpha]$ is not dense in the Heyting algebra $\LTInt$. Thus, $\alpha$ will be invalid in the homomorphic image with respect to the homomorphism determined by the filter of all dense elements of $\LTInt$, which, according Rasiowa and Sikorski~\cite{rs70}, chapter IV, theorem 5.8, is a Boolean algebra. Then, in virtue of Corollary~\ref{C:LT_Cl(k)}, $\alpha$ is invalid in $\LTCl$ and hence (Proposition~\ref{P:valuations-in-LT}) $\alpha\notin\bm{T}_{\Cl}$.

\subsection{Completeness}
In this subsection by `completeness' we understand the following. A unital logic $\aLog$ is complete with respect to that class $\TS$ (in this narrow sense) if the equivalence
\[
\alpha\in\bm{T}_{\mathcal{S}}~\Longleftrightarrow~\alpha~\text{is valid in all algebras of}~\TS, 
\]
holds for all formulas $\alpha$.

The above equivalence, in particular, is true if the conditions of Proposition~\ref{P:LT=free-algebra} are fulfilled. This allowed Rasiowa and Sikorski to prove in~\cite{rs70} that for any $\Lan_A$-formula $\alpha$,
\[
\alpha\in\bm{T}_{\Int}~\Longleftrightarrow~\text{$\alpha$ is valid in all Heyting algebras}.
\]

However, the first result of this kind was obtained in~\cite{mckinsey1941} in relation to the Lewis systems S2 and S4.
But it was A. Tarski who opened the door towards this direction, proving in~\cite{tarski1935}, theorem 4, that $\LTCl$ is a Boolean algebra.

\subsection{Admissible and derivable structural inference rules}
Let $\aLog$ be a unital logic. Since the logic $\aLog$ is structural, any modus rule sound with respect to $\vdash_{\mathcal{S}}$ can be written in terms of a pair $\langle X,\alpha\rangle$, where $X=\lbrace\alpha_1,\ldots,\alpha_n\rbrace$ is a finite (perhaps empty) set of formulas, not metaformulas as we formulated in Section~\ref{section:inference-rules} in general.  According to this definition, for any $\Lan$-substitution $\sigma$, $\sigma(X)\vdash_{\mathcal{S}}\sigma(\alpha)$; and, according to Proposition~\ref{P:preLindenbaum-algebra} and Proposition~\ref{P:valuations-in-LT}, $X\models_{\LT[X]}\alpha$, which, in virtue of~\eqref{E:model=algebra-validity-2}, implies that
\[
\langle\LT[X], \approx\rangle\models\Forall\ldots\Forall\,(\alpha_1\approx\one_{D}\&\ldots\&\alpha_n\approx\one_{D}\Rarrow\alpha\approx\one_{D}),
\]
where $D=\ConS{X}$. The latter is especially interesting when we investigate \textit{\textbf{admissible rules}}\index{rule!admissible} of a unital calculus, that is, those rules that, being added as postulated rules, do not change the set of theses of the calculus. Then, a modus rule $\langle X,\alpha\rangle$ is admissible in a unital calculus $\aLog$ if, and only if,
\[
\langle\LT,\approx\rangle\models\Forall\ldots\Forall\,(\alpha_1\approx\one_{\ThmS}\&\ldots\&\alpha_n\approx\one_{\ThmS}\Rarrow\alpha\approx\one_{\ThmS}),
\]

The first to pay attention to this equivalence with respect to $\Int$ was A. Citkin; cf.~\cite{citkin1977}. This idea was later picked up by other researches, which led to a correspondence between finitary structural abstract logics and quasi-varieties; for more information see~\cite{rybakov1997}, section 1.4, and~\cite{gorbunov1998}, section 1.5.2.\\

Given a consequence relation $\vdash$ and a set of rules $\mathcal{R}$, one can construct a new, stronger, consequence relation $\vdash_{\mathcal{R}}$: for any set of formulas $\Gamma$ and a formula $\alpha$, $\Gamma \vdash_{\mathcal{R}} \alpha$ if, and only if, there is a sequence of formulas $\alpha_0,\dots,\alpha_n$ such that $\alpha_n$ = $\alpha$ and for each  $i < n$, $\alpha_i \in \Gamma$, or $\alpha_0,\dots,\alpha_{i-1} \vdash \alpha_i$, or $\alpha_i$ can be derived from the preceding formulas by application of an instance of one of the rules from $R$. Such a sequence is called an $\vdash_{\mathcal{R}}$-\textit{inference}.  

It is clear that if each rule of $\mathcal{R}$ is admissible in a consequence relation $\vdash$, then consequence relations $\vdash$ and $\vdash_{\mathcal{R}}$ have the same sets of theorems.

\textbf{Example.} According to Lemma~\ref{L:towards-p^n-unprovability}, for any $n \ge 0$,
if
$\vdash_{\textsf{Int}} p^{2n+7}$, either $\vdash_{\textsf{Int}} p^{2n^3}$ or
$\vdash_{\textsf{Int}}  p^{2n+1}$. Hence, for each $n \ge 0$, the rule 
\begin{equation}
	R_n := p^{2n+7}/p^{2n^3} \lor p^{2n^1} \label{eq-rn}
\end{equation}
is admissible for $\vdash_{\textsf{Int}}$ and therefore, for any set $N$ of natural numbers, if $\mathcal{R} = \{R_n \ | \ n \in N\}$, the set of the theorems of the consequence relation $\vdash_{\mathcal{R}}$ coincides with $\bm{T}_{\Int}$. \\

A rule $R = \Gamma/\alpha$ is said to be \textit{derivable from a set of rules} $\mathcal{R}$ relative to a consequence relation $\vdash$ (in symbols $\mathcal{R} \vdash R$) if
\begin{equation}
	\Gamma \vdash_\mathcal{R} \alpha. \label{eq-der}
\end{equation} 
Sets of rules $\mathcal{R}_0$ and $\mathcal{R}_1$ are said to be \textit{equivalent relative to a consequence relation} $\vdash$ (in symbols $\mathcal{R}_0 \approx_\vdash \mathcal{R}_1$) if $\mathcal{R}_0 \vdash R$ for each $R \in \mathcal{R}_1$ and $\mathcal{R}_1 \vdash R$ for each $R \in \mathcal{R}_0$.

\begin{prop}\label{pr-equvr}
	Let $\vdash$ be a consequence relation and let $\mathcal{R}_0$ and $\mathcal{R}_1$ be equivalent sets of rules relative to $\vdash$. Then for each rule $R$,
	$\mathcal{R}_0 \vdash R$ entails $\mathcal{R}_1 \vdash R$, and conversely.
\end{prop}
\begin{proof}
	Suppose that $\mathcal{R}_0 \vdash R$, where $R = \Gamma/\alpha$. Then, $\Gamma \vdash_{\mathcal{R}_0} \alpha$, that is, there is a $\vdash_{\mathcal{R}_0}$-inference of $\alpha$ from $\Gamma$. The proof can be carried out by induction by induction on the minimal number of the applications of the rules from $R_0$ in the $\vdash_{\mathcal{R}_0}$-inferences of $\alpha$ from $\Gamma$. The basis is obvious, and the proof of the induction step fis based on the observation that in the $\vdash_{\mathcal{R}_0}$-inference of $\alpha$ from $\Gamma$ the last application of rule from $\mathcal{R}_0$, let say $\Delta/\beta$, can be replaced by an $\vdash_{\mathcal{R}_1}$-inference of $\beta$ from $\Delta$. A more detailed sketch of the proof of the induction step is the following. 
	
	Suppose that in a $\vdash_{\mathcal{R}_0}$-inference $\alpha_0,\dots,\alpha_n = \alpha$ of $\alpha$ from $\Gamma$, formula $\alpha_i$ is obtained by application of an instance $\sigma(\Delta)/\sigma(\beta)$ of the rule $\Delta/\beta \in \mathcal{R}_0$. By assumption, $\mathcal{R}_1 \vdash \Delta/\beta$, that is, by the definition, $\Delta \vdash_{\mathcal{R}_1} \beta$. Hence, there is a $\vdash_{\mathcal{R}_1}$-inference $\beta_0,\dots,\beta_m=\beta$ of $\beta$ from $\Delta$. Thus, $\sigma(\beta_0),\dots,\sigma(\beta_m)$ is a $\vdash_{\mathcal{R}_1}$-inference of $\sigma(\beta)$ from $\sigma(\Delta)$. Now we observe that we can replace the application of $\sigma(\Delta)/\sigma(\beta)$ with $\sigma(\beta_0),\dots,\sigma(\beta_m)$, and use the induction assumption to obtain a $\vdash_{\mathcal{R}_1}$-inference of $\alpha$ from $\Gamma$.  
\end{proof}
\begin{cor}\label{cor-eqvr}
	Let $\vdash$ be a consequence relation and let $\mathcal{R}_0$ and $\mathcal{R}_1$ be sets of rules. Then
	\begin{equation}
		\vdash_{\mathcal{R}_0} = ~\vdash_{\mathcal{R}_1}~\Longleftrightarrow~ \mathcal{R}_0 \approx_\vdash \mathcal{R}_1. \label{eq-ruleq}
	\end{equation}
\end{cor}
\begin{proof} Suppose that $\vdash_{\mathcal{R}_0} = \vdash_{\mathcal{R}_1}$. Let $r = \Gamma/\alpha$ and $R \in \mathcal{R}_0$. It is clear that $\Gamma \vdash_{\mathcal{R}_0} \alpha$. Hence, by the assumption, $\Gamma \vdash_{\mathcal{R}_1} \alpha$, which means that $\mathcal{R}_1 \vdash R$. Thus, for each $R \in \mathcal{R}_0$, $R_1 \vdash R$. By the same argument, $\mathcal{R}_0 \vdash R$ for each $R \in\mathcal{R}_1$, that is, $\mathcal{R}_0 \approx_\vdash \mathcal{R}_1$.
	
	Conversely, suppose that $\mathcal{R}_0 \approx_\vdash \mathcal{R}_1$. Let $\Gamma$ be a set of formulas and $\alpha$ be a formula. Let us consider the rule $R := \Gamma/\alpha$. If $\Gamma \vdash_{\mathcal{R}_0} \alpha$, then by \eqref{eq-der}, $\mathcal{R}_0 \vdash R$ and by Proposition \ref{pr-equvr}, $\mathcal{R}_1 \vdash R$, which means that $\Gamma \vdash_{\mathcal{R}_1} \alpha$. Thus, $\vdash_{\mathcal{R}_0} \ \subseteq \ \vdash_{\mathcal{R}_1}$. By the same argument, $\vdash_{\mathcal{R}_1} \ \subseteq \ \vdash_{\mathcal{R}_0}$.	
\end{proof}

Let $\mathcal{R} = \{R_n \ | \ n \ge 0\}$, where $R_n$ are the rules defined in \eqref{eq-rn}, and $\mathcal{R}_n := \mathcal{R} \setminus\set{R_n}{n \ge 0}$. It was observed in \cite{citkin1977} that $\mathcal{R}_n \not\vdash_{\textsf{Int}} R_n$.
Hence,  by Corollary~\ref{cor-eqvr}, for any $\mathcal{R}_0, \mathcal{R}_1 \subseteq \mathcal{R}$, if $\mathcal{R}_0 \not\approx_{\vdash_{\textsf{Int}}} \mathcal{R}_1$, then, $\vdash_{R_0} \ \neq \ \vdash_{R_1}$.
Consequently, the set of consequence relations each having the set of its theorems equal to $\bm{T}_{\Int}$, is not countable. 

\subsection{Classes of algebras corresponding and fully corresponding to a calculus}
\begin{defn}\label{D:normal-extension}\index{logic!unital!normal extension}
	Let $\aLog$ be a unital calculus. A unital logic $\aLog^{\prime}$ is said to be a \textbf{normal extension of} $\aLog$ if {\em$\ConS{\varnothing}\subseteq\textbf{Cn}_{\mathcal{S}^\prime}(\varnothing)$}
	and	{\em$\textbf{Cn}_{\mathcal{S}^\prime}(\varnothing)$} is closed under all rules of inference that are postulated in the calculus $\aLog$.
\end{defn}

For example, both $\Cl$ and $\LC$ are normal extension of $\Int$.

\begin{defn}\label{D:corresponds-and-fully-corresponds}
	Let $\aLog$ be a unital logic and $\mathcal{K}$ be a class of unital\index{logic!unital!corresponds} expansions. We say that $\mathcal{K}$ \textbf{corresponds to} $\aLog$ if for any formula $\alpha$,
	{\em\[
		\alpha\in\ConS{\varnothing}~\Longleftrightarrow~(\alg{A}\models\alpha,~\text{for any $\alg{A}\in\mathcal{K}$});
		\]}
	and $\mathcal{K}$ \textbf{fully corresponds to} $\aLog$ if for any
	{\em$X\cup\lbrace\alpha\rbrace\subseteq\FormsL$}\index{logic!unital!fully corresponds}
	{\em\[
		X\vdashS\alpha~\Longleftrightarrow~(\alg{A}\models X\Longrightarrow\alg{A}\models\alpha, \text{for any $\alg{A}\in\mathcal{K}$});
		\]}	
\end{defn}

\begin{prop}\label{P:corresponds}
	Given a unital logic $\aLog$, the class $\TS$ corresponds to $\aLog$. Moreover, $\TS$ is the largest such a class.	
\end{prop}
\begin{proof}
	In virtue of~\eqref{E:LT-algebras-unital-2}, $\LT\in\TS$. Then, we apply the second part of Proposition~\ref{P:valuations-in-LT}.
	
	The second part is obvious.
\end{proof}

It must be clear that
\begin{equation}\label{E:fully-corresponds-implies-corresponds}
	\textit{If a class $\mathcal{K}$ fully corresponds to $\mathcal{S}$, $\mathcal{K}$ also corresponds to $\mathcal{S}$.}
\end{equation}

The converse is not always true.\footnote{See examples in~\cite{muravitsky2014b}, section 7.3.}

Next we define: For any unital calculus $\aLog$,
\[
\PS^{\circ}:=\set{\alg{A}\in\PS}{\aLog[L\alg{A}]~\text{is a normal extension of $\aLog$}}.
\]
\begin{lem}\label{L:{T[D] | D-thoery}-subset-PS^circ}
	Let $\aLog$ be a unital logic and $\theory$ be the sett of all $\aLog$-theories. Then {\em$\set{\LT[D]}{D\in\theory}\subseteq\PS^\circ$}.	
\end{lem}
\begin{proof}
	It suffices to show that for any instance $\alpha_1,\ldots,\alpha_n\slash\beta$ of an inference rule of the calculus $\aLog$,
	\[
	\LT[D]\models\lbrace\alpha_1,\ldots,\alpha_n\rbrace~\Longrightarrow~
	\LT[D]\models\beta,
	\]
	for any $D\in\theory$. So, suppose $\LT[D]\models\lbrace\alpha_1,\ldots,\alpha_n\rbrace$. 	Since $\lbrace\alpha_1,\ldots,\alpha_n\rbrace\vdashS\beta$, in virtue of the first part of Proposition~\ref{P:valuations-in-LT}, $\LT[D]\models\beta$.
\end{proof}
\begin{prop}\label{P:fully-corresponds}
	Given a unital calculus $\aLog$, the class $\PS^{\circ}$ fully corresponds to $\aLog$. And if a class 	$\mathcal{K}\subseteq\PS$ fully corresponds to $\aLog$, then $\mathcal{K}\subseteq\PS^{\circ}$.
\end{prop}
\begin{proof}
	Let $X\vdashS\alpha$ and suppose that $\alg{A}\models X$, where $\alg{A}\in\PS^{\circ}$. Since $\alg{A}$ is an algebraic $\mathcal{S}$-model, 
	$\alg{A}\models\alpha$. 
	
	Next, assume that $X\not\vdash_{\mathcal{S}}\alpha$. We note that, in virtue of Proposition~\ref{P:valuations-in-LT}, $\LT[D]\models X$ and $\LT[D]\not\models\alpha$, where $D=\ConS{X}$. However, according to Proposition~\ref{P:LT-unital-PS-QS-equivalence}, $\LT[D]\in\PS$. It remains to show that $L\LT[D]$ is closed under any inference rule postulated in $\mathcal{S}$.
	
	Let $\alpha_1,\ldots,\alpha_n\slash\beta$ be an instance of such a rule and $\LT[D]\models\lbrace\alpha_1,\ldots,\alpha_n\rbrace$. In virtue of Proposition~\ref{P:valuations-in-LT}, $\lbrace\alpha_1,\ldots,\alpha_n\rbrace\subseteq\ThmS$. By premise, this implies that $\beta\in\ThmS$ and hence (Proposition~\ref{P:valuations-in-LT})
	$\LT[D]\models\beta$, that is, $L\LT[D]$ is closed under the given rule of inference. Thus $\LT[D]\in\PS^\circ$.
	
	This completes the proof that $\PS^\circ$ fully corresponds to $\mathcal{S}$.
	
	Now, suppose that a class $\mathcal{K}\subseteq\PS$ fully corresponds to $\aLog$ and $\alg{A}\in\mathcal{K}$. We need to show $L\alg{A}$ is closed under any inference rule postulated in $\aLog$, since $\textbf{Cn}_{\mathcal{S}[L\alg{A}]}(\varnothing)=L\alg{A}$.
	
	Indeed, let $\alpha_1,\ldots,\alpha_n\slash\beta$ be an instance of such a rule. This implies that $\lbrace\alpha_1,\ldots,\alpha_n\rbrace\vdashS\beta$.
	Thus, if $\alg{A}\models\lbrace\alpha_1,\ldots,\alpha_n\rbrace$, then, by premise, $\alg{A}\models\beta$.
\end{proof}

\begin{cor}\label{C:fully-corresponds}
	Given a unital calculus $\aLog$, $\PS^{\circ}\subseteq\TS$.
\end{cor}
\begin{proof}
	We prove the inclusion in question, by using Proposition~\ref{P:fully-corresponds},~\eqref{E:fully-corresponds-implies-corresponds} and  Proposition~\ref{P:corresponds}.	
\end{proof}

\begin{prop}
	Let $\mathcal{K}\subseteq\PS$ be an equational class defined by identities of the form $\gamma\approx\one$. If $\mathcal{K}$ fully corresponds to a unital calculus $\aLog$, then $\mathcal{K}=\PS^{\circ}=\TS$.
\end{prop}
\begin{proof}
	In virtue of Proposition~\ref{P:fully-corresponds} and Corollary~\ref{C:fully-corresponds}, we have: $\mathcal{K}\subseteq\PS^{\circ}\subseteq\TS$. For contradiction, we assume that $\mathcal{K}\subset\TS$. This implies that there is an identity $\alpha\approx\one$ that belongs the equational theory of $\mathcal{K}$ but $\alpha\notin\bm{T}_{\mathcal{S}}$. Since, according to~\eqref{E:fully-corresponds-implies-corresponds}, $\mathcal{K}$ corresponds to $\aLog$, there is an algebra  $\alg{A}\in\mathcal{K}$ such that $\alg{A}\not\models\alpha$; that is, $\alpha\approx\one$ does not belong to the equational theory of $\mathcal{K}$. A contradiction.
\end{proof}

\section{An alternative approach to the concept of Lindenbaum-Tarski algebra}
The approach to the concept of Lindenbaum-Tarski algebra, which we present concisely below, is not supposed to be used in the search of convenient separating tools in the sense of Section~\ref{section:separating-means}.

Let $\mat{M}=\langle\alg{A},D\rangle$ be a logical matrix. It is not difficult to see that a congruence $\theta$ on $\alg{A}$ is compatible with $D$ if, and only if
\[
D=\bigcup\set{x\slash\theta}{x\in D}.
\]
The set of all congruences compatible with $D$ has a least element, the equality, and, as was proved in~\cite{blok-pigozzi1989}, has also a largest congruence with respect to $\subseteq$. This congruence is called the \textit{Leibniz congruence} of $\mat{M}$; it is denoted by $\Omega_{\alg{A}}D$ and can be defined as follows:
\[
\begin{array}{rl}
	\Omega_{\mat{A}}D=
	&\lbrace(a,b)~|~\forall
	\alpha(p,p_{0},\ldots,p_{n})\in\FormsL\forall c_{0},\ldots, c_{n}\in|\alg{A}|.\\
	&\alpha(a,c_{0},\ldots,c_{n})\in D\Leftrightarrow \alpha(b,c_{0},\ldots,c_{n})\in D\rbrace.
\end{array}
\]
If the matrix in question is a Lindenbaum matrix $\langle\FormAl,D_{\mathcal{S}}\rangle$, where $D_{\mathcal{S}}$ is an $\aLog$-theory, then an example of a compatible equivalence on this matrix is the Frege relation relative to $D_{\mathcal{S}}$ (Definition~\ref{D:fregean-relation}).

Another example of a compatible relation on $\langle\alg{A},D_{\mathcal{S}}\rangle$ is the largest congruence of {\alg{A}} contained in $\Lambda D_{\mathcal{S}}$, which is referred to as a \textit{Suszko congruence}:
\[
\begin{array}{rr}
	(\alpha,\beta)\in\Tilde{\Omega}D_{\mathcal{S}}~\Longleftrightarrow~
	&\text{for every $\gamma(p)$, $D_{\mathcal{S}}, \gamma(\alpha/p)\vdashS \gamma(\beta/p)$}\\
	&\text{and $D_{\mathcal{S}},\gamma(\beta/p)\vdashS \gamma(\alpha/p)$}.
\end{array}
\tag{\textit{Suszko congruence relative to} $D_{\mathcal{S}}$}
\]

Obviously, a system $\aLog$ is Fregean if, and only if, $\Lambda D_{\mathcal{S}} =\Tilde{\Omega}D_{\mathcal{S}}$, for all $D_{\mathcal{S}}$.

The Leibniz congruence of a matrix $\langle\FormAl,D_{\mathcal{S}}\rangle$ is referred to as \textit{Leinbniz congruence relative to} $D_{\mathcal{S}}$. It  turns out that
\[
\Omega D_{\mathcal{S}}=\bigcap\set{\Tilde{\Omega}D^{\prime}_{\mathcal{S}}}{D_{\mathcal{S}}\subseteq D^{\prime}_{\mathcal{S}}}
\]
and, therefore, each Suszko congruence $\Tilde{\Omega}D_{\mathcal{S}}$ is compatible with $D_{\mathcal{S}}$. Also, given an abstract logic $\mathcal{S}$, one defines
\[
\Tilde{\Omega}_{\mathcal{S}}=\bigcap\set{\Tilde{\Omega}D_{\mathcal{S}}}{
	D_{\mathcal{S}}~\text{is an $\mathcal{S}$-theory}}.\tag{\emph{Tarski congruence}}
\]

Thus we have:
\[
\Tilde{\Omega}_{\mathcal{S}}\subseteq\Tilde{\Omega}	D_{\mathcal{S}}\subseteq\Lambda
D_{\mathcal{S}}\cap\Omega 	D_{\mathcal{S}}.
\]

Suszko, Leibniz and Tarski congruences give rise to the $\aLog$-matrices $\langle\FormAl/\Omega	D_{\mathcal{S}},	D_{\mathcal{S}}/\Omega D_{\mathcal{S}}\rangle$, $\langle\FormAl/\Tilde{\Omega}	D_{\mathcal{S}},D_{\mathcal{S}}/\Tilde{\Omega}	D_{\mathcal{S}}\rangle$, and the atlas $\langle\FormAl/\Tilde{\Omega}_{\mathcal{S}},\set{D_{\mathcal{S}}/\Tilde{\Omega}_{\mathcal{S}}}
{D_{\mathcal{S}}~\text{is an $\aLog$-theory}}\rangle$, whose first components, $ \FormAl/\Omega D_{\mathcal{S}} $,
$\FormAl/\Tilde{\Omega}	D_{\mathcal{S}} $ and $ \FormAl/\Tilde{\Omega}_{\mathcal{S}} $, in \textit{Algebraic abstract logic} are referred to as \textit{Lindenbaum-Tarski algebras}. (See~\cite{fjp03} and~\cite{fj09} with update in~\cite{fjp09} for comprehensive surveys.)

\section{Historical notes}
The notion of a unital matrix is ​​a natural bridge from (ordinary) logical matrices to algebras. Such a transition makes it possible to employ more methods of universal algebra than methods of model theory. This notion was formally introduced in~\cite{czelakowski2001} and extensively used therein, although many of the firs logical matrices in the history of many-valued logic were unital. Thus, by 2010, when Czelakowski's book was published, the concept, not the term, had been well known. 

The class of unital logics was defined in~\cite{mur2014a}, although the name for this class of logics was suggested by an anonymous referee of this paper. In fact, this notion is rooted in Definition 4 of~\cite{muravitsky2014b}, where ``a calculus $C$ [that] admits the Lindenbaum-Tarski algebra'' was defined,  which is somewhat weaker than the notion of a unital logic. 

It was established in~\cite{mur2014a} that the class of implicative logics (studied in~\cite{ras74}) and that of Freagean logics (studied in~\cite{czelakowski-pigozzi2004a}) are properly included in the class of unital logics.

In \textit{Algebraic Abstract Logic}, a different approach to the notion of a Lindenbaum-Tarski algebra has been developed. See~\cite{font2016} and references therein for detail; a brief outline of this approach can be found in~\cite{citkin-mur2013}. Unlike ours, this approach does not focus on the Lindenbaum-Tarski algebra as a source for constructing separating tools.

The second part of Corollary~\ref{C:Cl-and-Int-uniform} was proved in~\cite{shoesmith-smiley1971}, theorem 5.

Corollary~\ref{C:LT_Int} was established in~\cite{rieger1949}, theorem 4.4.

The diagram in Figure~\ref{fig-RN} was discovered by L. Rieger in~\cite{rieger1957} and I. Nishimura in~\cite{nishimura1960}. However, earlier in~\cite{rieger1949}, theorem 4.6, Rieger, as well as Nishimura in~\cite{nishimura1960}, {\S} 3, gave the definition of the elements of $\LTInt(1)$ similar to definition~\eqref{E:powers_p^n}. Both researchers merely sketched their proofs. Referring to theorem 4.6, Rieger wrote:
\begin{quote}
	``The proof is essentially not difficult but somewhat labourious.''\cite{rieger1949}, p. 33
\end{quote}

For his part, referring to the diagram, Nishimura wrote:
\begin{quote}
	``We can investigate inclusion, noninclusion and equivalence relationships between these logics [\dots]'' \cite{nishimura1960}, {\S} 3
\end{quote}

Both researchers only sketched the proof of Proposition~\ref{P:nishimura-completeness} in order to obtain with its help Corollary~\ref{C:nishimura-completeness}. However, they failed to notice that Corollary~\ref{C:nishimura-power-difference} is also important.\\

For the treatment of $\Int$-theses (and, in fact, for a broader framework), G. Gentzen introduced in~\cite{gentzen1964} a sequential calculus LJ.\footnote{An English translation of ~\cite{gentzen1964} can be found in~\cite{gentzen-Collected}, pp. 68--131.} Kleene calls in~\cite{kleene1952}, {\S} 77, this calculus G1 for the intuitionistic system.\footnote{The difference  between the classical and intuitionistic systems G1 consists in one restriction for two postulates of G1. See~\cite{kleene1952}, {\S} 77 for detail.} Kleene's sequential calculus G3 (see~\cite{kleene1952}, {\S} 80) was proved equivalent to G1 (in the sense of~\cite{kleene1952}, {\S} 80, Theorem 56). The treatment of $\Int$-theses in the framework of G3 is carried out through Proposition~\ref{P:G3-completeness}.

\chapter[Equational Consequence]{Equational Consequence}
\label{chapter:equational-con}
The type of consequence we are going to discuss in this chapter is perhaps one of the earliest kinds of reasoning found in mathematical practice. Indeed, suppose we want to solve an equation $x+a=b$. Assuming that $c$ is a solution to this equation, we, then, proceed as follows.
\begin{equation}\label{E:equation-con-1}
	\begin{array}{rl}
		c+a=b &\Longrightarrow c=b-a.
	\end{array}
\end{equation}

Thus we obtain a solution under the assumption that a solution exists. To see that a solution does exist, we observe that the step indicated by \eqref{E:equation-con-1}, in fact, is invertible, that is, actually, the equivalence
\begin{equation}\label{E:equation-con-2}
	\begin{array}{rl}
		c+a=b &\Longleftrightarrow c=b-a
	\end{array}
\end{equation}
holds. As we will see below, \eqref{E:equation-con-2} is the exception rather than the rule. Therefore, we focus on the type of consequence exemplified by \eqref{E:equation-con-1}.

Also, we would like to note that in the antecedent `$c+a=b$' of \eqref{E:equation-con-1} we only \emph{assume} that it is true. In other words, we will be treating equalities rather as \emph{formal} messages which can be true or false. This viewpoint requires a special formal language which we introduce in the next section.

\section{Language for equalities and its semantics}	\label{section:language-equalities}
Our language for equalities has the same schematic character as the language $\Lan$ of Chapter~\ref{chapter:languages}. Actually, the new language is based on $\Lan$ but is not its extension (in the sense of Section~\ref{section:languages}). What is called $\Lan$-formulas in Definition~\ref{D:L-formulas}, in this chapter we call $\Lan$-\textit{\textbf{terms}}; accordingly, sentential (or propositional) variables we call \textit{\textbf{term variables}} and subformulas now will be called \textit{\textbf{subterms}}. The set of all (term) variables is denoted by $\VarL$ (or simply by $\Var$), and the set of variables occurring in the formulas of a set $X$ by $\Var(X)$; in particular, given an equality $\epsilon$, $\Var(\epsilon)$ denotes the set of all variables occurring in $\epsilon$. If language $\Lan$ is specified as a language $\Lan^{\prime}$, $\Lan^{\prime\prime},\ldots$, accordingly, we have $\Lan^{\prime}$-terms, $\Lan^{\prime\prime}$-terms and so on.

For formulas, the language $\Lan$ is extended by one predicate symbol, $\approx$, which will be used in the infix notation; that is for any $\Lan$-terms $\alpha$ and $\beta$, $\alpha\approx\beta$ is a \textit{\textbf{formula}}. The formulas are only those words that can be obtained in this way. The formulas are also called $\Lan$-\textit{\textbf{equalities}} (or simply \textit{\textbf{equalities}}). The set of all equalities of $\Lan$-terms is denoted by $\EqL$. Slightly abusing notation, we denote the new formal language also by $\EqL$.

In accordance with this change of terminology, in this chapter we call the set $\FormsL$ of Section~\ref{section:languages} the \textit{\textbf{set of $\Lan$-terms}} and the algebra $\FormAl$ the \textit{\textbf{algebra of $\Lan$-terms}}.

The notion of substitution defined for $\Lan$-formulas in Definition~\ref{D:substitution} now is renamed as \textit{\textbf{substitution of}} $\Lan$-\textit{\textbf{terms}} {\textit{\textbf{for term variables}}}.
(We continue calling it $\Lan$-\textit{\textbf{substitution}}.) And the latter is extended to the equalities, so that if $\epsilon=\alpha\approx\beta$ and $\sigma$ is an $\Lan$-substitution, then
\[
\sigma(\epsilon):=\sigma(\alpha)\approx\sigma(\beta).
\]
As before, the \textit{\textbf{unit substitution}} is denoted by $\iota$; that is 
\[
\iota(\alpha\approx\beta):=\alpha\approx\beta.
\]

Similar to what we defined for $\Lan$-formulas in Section~\ref{section:languages}, the last definition is extended to any set $X$ of equalities as follows:
\[
\sigma(X):=\set{\sigma(\epsilon)}{\epsilon\in X}.
\]

In parallel, we treat $\Lan$-metaformulas of Definition~\ref{D:L-metaformula} as
$\Lan$-\textit{\textbf{metaterms}}. Then, given $\Lan$-metaterms $\phi$ and $\psi$, we call a \textit{\textbf{metaequality}} the following expression:
\[
\phi\approx\psi.
\]
If $\bm{\epsilon}$ is metaequality, say
\[
\bm{\epsilon}:=\phi\approx\psi,
\]
in which some occurrences of a metavariable $\bm{\alpha}$ are designated, symbolically $\bm{\epsilon}[\bm{\alpha}]:=\phi[\bm{\alpha}]\approx\psi[\bm{\alpha}]$, ($\bm{\alpha}$ may occur in $\phi$ or in $\psi$ or in both), then we denote by
\[
\bm{\epsilon}[\bm{\beta}]:=\phi[\bm{\beta}]\approx\psi[\bm{\beta}]
\]
the result of replacement of the designated occurrences of $\bm{\alpha}$ with a metavariable $\bm{\beta}$.
We will need metaequalities, when we discuss rules of inference for equational consequence.\\

Let $\alg{A}=\langle\textsf{A},\Func,\Cons\rangle$ be a (nontrivial) algebra of type $\Lan$ and $D(x,y)$ be a binary predicate on $|\alg{A}|$. We call the structure $\langle\alg{A}, D(a,y)\rangle$ a \textit{\textbf{binary logical matrix}}. Instead of predicate $D(x,y)$, we will also be using a binary relation
$D\subseteq\textsf{A}\times\textsf{A}$, in which case a binary logical matrix will be denoted in the form $\langle\alg{A}, D\rangle$. We continue to call $D$ a \textit{\textbf{logical filter}}.

A binary logical matrix $\langle\alg{A}, D\rangle$ is called  an \textit{\textbf{equational matrix}} (or an E-\textit{\textbf{matrix}} for short) if $D=\varDelta$, where $\varDelta$ is the diagonal  of $\textsf{A}\times\textsf{A}$.

We recall that, working with $\varDelta$ as a binary predicate, we have:
\[
(x,y)\in\varDelta~\Longleftrightarrow~x=y,
\]
for any $x, y\in|\alg{A}|$. We note that, given equational matrices $\mat{M}=\langle\alg{A},\varDelta\rangle$ and $\mat{N}=\langle\alg{B},\varDelta\rangle$, the first occurrence of $\varDelta$ is the diagonal of $|\alg{A}|\times|\alg{A}|$ and the second on $|\alg{B}|\times|\alg{B}|$. It should not cause any confusion, since we will be using $\varDelta$ always as the diagonal of the Cartesian square of the carrier of an indicated  algebra. However, sometimes, when an algebra $\alg{A}$ is introduced and we need to refer to the diagonal of $|\alg{A}|\times|\alg{A}|$, we will use $\varDelta$ with a subscript --- $\diagA$.

Given E- matrices $\mat{M}=\langle\alg{A},\varDelta\rangle$ and $\mat{N}=\langle\alg{B},\varDelta\rangle$, any homomorphism of $\alg{A}$ to $\alg{B}$ we treat as a homomorphism of $\mat{M}$ to $\mat{N}$; and we say $\mat{M}$ and $\mat{N}$ are \textit{\textbf{isomorphic}} if so are $\alg{A}$ and $\alg{B}$. 

Given an equality $\epsilon=\alpha\approx\beta$  and a valuation $v$ in algebra $\alg{A}$, by a \textit{\textbf{value of an equality}} $\epsilon$ in $\alg{A}$ \textit{\textbf{under}} $v$ we understand the pair $v[\epsilon]=(v[\alpha],v[\beta])$,
where $v[\alpha]$ and $v[\beta]$ are the corresponding values of $\alpha$ and $\beta$; see Definition~\ref{D:valuation}.

\section{Equational consequence}\label{section:E-consequence}
We begin with the following definition.
\begin{defn}[E-matrix semantics]\label{D:E-matrix-semantics}\index{E-matrix!semantics}
	Let {\em$\mat{M}=\langle\alg{A},D\rangle$} be a binary matrix. We say that an equality $\alpha\approx\beta$ is \textbf{satisfiable}\index{E-matrix!stisfiable equality} in {\em$\mat{M}$} if there is a valuation $v$ in {\em$\alg{A}$} such that $(v[\alpha],v[\beta])\in D$, in which case we say that $\alpha\approx\beta$ is \textbf{satisfied at $v$}, or $v$ \textbf{satisfies} $\alpha\approx\beta$ \textbf{in} {\em$\mat{M}$}. A set $X$ of equalities is \textbf{satisfiable} in {\em$\mat{M}$} if there is a valuation $v$ such that all equalities of $X$ are satisfied at $v$. An equality $\alpha\approx\beta$ is \textbf{valid} in {\em$\mat{M}$} if it is satisfied at any valuation $v$ in {\em$\mat{A}$}. Accordingly, a set $X$ of equalities is \textit{\textbf{valid}} in\index{E-matrix!valid equality} {\em$\mat{M}$} if each equality of $X$ is valid.
	We denote the satisfiability of $X$ in {\em$\mat{M}$} at $v$ by
	{\em\[
		\mat{M}\models_{v}X
		\]}
	and validity by
	{\em\[
		\mat{M}\models X.
		\]}
\end{defn}

This suggests to extend the definition (\ref{E:matrix-consequence}) of matrix consequence to comprise the new notion of formula (equality) and new notion of logical matrix (binary logical matrix). 

Thus for a (nonempty) class $\mathcal{M}$ of binary matrices, we define a relation $\models_{\mathcal{M}}\subseteq\mathcal{P}(\EqL)\times\EqL$ as follows:
\begin{eqnarray}\label{E:def-binary-matrix-con}
	\begin{array}{rcl}
		X\models_{\mathcal{M}}\alpha\approx\beta &\stackrel{\text{df}}{\Longleftrightarrow}
		&\text{for any matrix $\mat{M}\in\mathcal{M}$ and any valuation}\\ 
		&&\text{$v$ in $\alg{A}$, $\mat{M}\models_{v}X$ implies $\mat{M}\models_{v}\alpha\approx\beta$}.
	\end{array}
\end{eqnarray}

We write simply
\[
X\models_{\textsf{\textbf{M}}}\alpha\approx\beta
\]
if $\mathcal{M}=\lbrace\textsf{\textbf{M}}\rbrace$.
\begin{prop}\label{P:binary-matrix-structural}
	The relation $\models_{\mathcal{M}}$ as defined in~\eqref{E:def-binary-matrix-con} is a structural consequence relation.
\end{prop}
\begin{proof}
	That $\models_\mathcal{M}$ is a consequence relation we obtain with the help of Proposition~\ref{P:semantic-consequence-2}. We leave to the reader to show that this consequence is structural.  (See Exercise~\ref{section:E-consequence}.\ref{EX:binary-matrix-structural})	
\end{proof}

The definition~\eqref{E:def-binary-matrix-con} and Proposition~\ref{P:binary-matrix-structural} give rise to the following concept.  

Let $\E$ be a nonempty class of E-matrices. 
The matrices of $\E$ are called $\E$-\textit{\textbf{matrices}}.

We define a relation $\vdashE\subseteq\mathcal{P}(\EqL)\times\EqL$ as follows:
\begin{eqnarray}\label{E:def-equation-con}
	\begin{array}{rcl}
		X\vdashE\alpha\approx\beta &\stackrel{\text{df}}{\Longleftrightarrow}
		&\text{for any $\E$-matrix {\em$\mat{M}=\langle\alg{A},\varDelta\rangle$} and any}\\ 
		&&\text{valuation $v$ in $\alg{A}$, $\mat{M}\models_{v}X$ implies $\mat{M}\models_{v}\alpha\approx\beta$}.
	\end{array}
\end{eqnarray}

We call any such a relation $\vdashE$ an \textit{\textbf{equational consequence}}.\index{E-matrix!equational correspondence} If a set $\E$ of E-matrices is specified, we call the corresponding equational consequence $\E$-\textbf{\textit{consequence}}\index{E-matrix!E-correspondence}
or \textit{\textbf{equational}} (\textit{\textbf{abstract}}) \textit{\textbf{logic}} $\E$.\index{equational abstract logic}\index{E-matrix!logic} Thus will be using the notation $\E$ as a class of E-matrices and as the equational abstract logic this class generates.\\

Applying Proposition~\ref{P:binary-matrix-structural} to $\vdashE$, we obtain straightforwardly the following.
\begin{cor}\label{C:E-structural-con}
	Any $\E$-consequence is a structural consequence relation.
\end{cor}

\begin{rem}\label{R:los-suszko-wojcicki-theorems}
	{\em To the question of whether any equational consequence can be given by a single E-matrix, we can answer as follows.
		Since any equational consequence is structural, although in a new sense that differs from the old one in that instead of substitution instances of terms we use substitution instances of equalities, many, if not all, of the properties of structural abstract logic can be extended to equational consequence. In particular, this is true about {\L}o\'{s}-Suszko-W\'{o}jcicki and W\'{o}jcicki theorems $($Proposition~\ref{P:los-suszko} and Proposition~\ref{P:wojcicki}, respectively$)$ understood in the context of equational consequence.}
\end{rem}

We leave the task mentioned in Remark~\ref{R:los-suszko-wojcicki-theorems}  to the reader. (Exercise~\ref{section:E-consequence}.\ref{EX:remark}) \\

The following property is rather technical; it will be needed in the sequel.
\begin{prop}\label{P:variables-auxiliary}
	Let a language $\Lan^{\prime}$ be a primitive extension of $\Lan$. Also, let $f$ be a function from $\Var_{\mathcal{L}^{\prime}}\setminus\VarL$ into $\FormAl$ Assume that $X$ is a set of $\Lan$-equalities and a set $Y$ is obtained by adding to $X$ some or all equalities of the set $\set{q\approx f(q)}{q\in\Var_{\mathcal{L}^{\prime}}\setminus\VarL}$. Then for any equality $\epsilon$ with $\Var(\epsilon)\subseteq\VarL$, if $Y\vdashE\epsilon$, then $X\vdashE\epsilon$.
\end{prop}
\begin{proof}
	Let $\langle\alg{A},\varDelta\rangle$ be a matrix of $\E$. Assume that $v[X]\subseteq\varDelta$, where $v$ is a valuation of $\Lan$-terms in $\alg{A}$. Now we define:
	\[
	\begin{array}{cl}
		v^{\prime}[q]:=\begin{cases}
			v[q] &\text{if $q\in\VarL$}\\
			v[f(q)] &\text{if $q\in\Var_{\mathcal{L}^{\prime}}\setminus\VarL$}.
		\end{cases}
	\end{array}
	\]
	
	It must be clear that for any $\Lan$-term $\alpha$,
	\[
	v^{\prime}[\alpha]=v[\alpha]. \tag{\ref{P:variables-auxiliary}--$\ast$}
	\]
	This implies that $v^{\prime}[Y]\subseteq\varDelta$ and hence, by premise,
	$v^{\prime}[\epsilon]\in\varDelta$. Applying (\ref{P:variables-auxiliary}--$\ast$) again, we obtain that $v[\epsilon]\in\varDelta$.
\end{proof}
\begin{defn}[abstract logics $\E$ and $\hat{\E}$]\index{E-matrix!abstract logic}
	The abstract logic and consequence operator corresponding to $\vdashE$ we denote also by $\E$ and by {\em$\ConE$}, respectively. If $\E$ is the class of all E-matrices of type $\Lan$, we denote the corresponding abstract logic by $\hat{\E}$.
\end{defn}

We immediately observe that $\ConEE(\varnothing)\neq\varnothing$, for $\varnothing\vdashEE\alpha\!\approx\!\alpha$, that is $\alpha\!\approx\!\alpha\in\ConEE(\varnothing)$, for any $\Lan$-term $\alpha$. In fact, $\ConEE(\varnothing)$ does not contain any other formulas. Indeed, let $\alpha$ and $\beta$ be two distinct $\Lan$-terms. Consider the equational matrix $\langle\FormAl,\varDelta\rangle$. We note that all valuations in this matrix are substitutions. Let us take the identity substitution $\iota$. Since $\iota(\alpha)\neq\iota(\beta)$ in $\FormAl$, according to~\eqref{E:def-equation-con}, $\varnothing\not\vdashEE\alpha\approx\beta$, that is $\alpha\approx\beta\notin\ConEE(\varnothing)$.

Thus we have proved the following.
\begin{prop}\label{P:alpha=alpha}
	{\em$\ConEE(\varnothing)=\set{\alpha\approx\alpha}{\alpha\in\FormsL}$}.
\end{prop}

In the sequel, we will find the following remark useful.
\begin{cor}\label{C:Con(0)_E-closed-under-sub}
	The class {\em$\ConEE(\varnothing)$} is closed under any $\Lan$-substitution.
\end{cor}
\noindent\emph{Proof}~is left to the reader. (Exercise~\ref{section:E-consequence}.\ref{EX:Con(0)_E-closed-under-sub})\\

To prove the next proposition we introduce a quantifier-free firs-order language, $\FOV$, in which the terms are built from individual variables, the individual constants are taken from the set $\Cons_{\Lan}\cup\VarL$, and $\Func_{\mathcal{L}}$ is a set of functional constants. (Thus the set of individual constants of $\FOV$ is never empty.) The only elementary formulas are equalities of $\FOV$-terms. Thus, the closed elementary $\FOV$-formulas  are $\Lan$-equalities and only them. 

Let $\mat{M}=\langle\alg{A},\varDelta\rangle$ be an E-matrix. We extend the signature of $\alg{A}$ by the set $\VarL$ of $0$-ary operations and denote so obtained algebra by $\alg{A}[\VarL]$. The structure $\mat{M}[\VarL]:=\langle\alg{A}[\VarL],\varDelta\rangle$ is a $\FOV$-\textit{\textbf{structure}}.\index{$\FOV$-structure}\index{E-matrix!$\FOV$-structure}

Thus, given an $\E$-matrix $\mat{M}=\langle\alg{A},\varDelta\rangle$ of type $\Lan$ and a valuation $v$ in $\alg{A}$, we can extend the interpreting function associated with $\mat{M}$ over the $\FOV$-constants of $\VarL$ by means of the valuation $v$. This leads to the understanding of an $\FOV$-structure $\mat{M}[\VarL]$ as a pair $(\mat{M}, v)$, where $\mat{M}:=\langle\alg{A},\varDelta\rangle$ is an $\Lan$-matrix and $v$ is a valuation in $\alg{A}$. Employing the usual definition of \textit{satisfiability} in first-order models (see Section~\ref{section:preliminaries-model-theory}), we observe that, given $\FOV$-model  $(\mat{M}, v)$, for any $\Lan$-equality $\epsilon$, $\epsilon$ is true in $(\mat{M}, v)$ if, and only if, $\mat{M}\models_{v}\epsilon$, in symbols
\begin{equation}\label{E:matrix-FO-structure-correlation}
	(\mat{M},v)\models
	\epsilon\Longleftrightarrow\mat{M}\models_{v}\epsilon.
\end{equation}

A direct consequence of \eqref{E:matrix-FO-structure-correlation} is the following observation.
\begin{prop}\label{P:finitariness-E-consequence}
	Let $\E$ be a nonempty class of E-matrices such that the class $\mathcal{A}[\E]$ is closed under ultraproducts. Then the abstract logic $\E$ is finitary.
\end{prop}
\begin{proof}
	Let $X\cup\lbrace\epsilon\rbrace\subseteq\EqL$. Regarding  each E-matrices of $\E$ as a $\FOV$-model, we obtain:
	\[
	\begin{array}{rcl}
		X\vdashE\epsilon&\Longleftrightarrow 
		&\text{for any E-matrix $\mat{M}$ and any valuation $v$, $\mat{M}\models_{v}X$ implies $\mat{M}\models_v\epsilon$};\\ 
		&\Longleftrightarrow
		&\text{the set $X\cup\lbrace\Neg\epsilon\rbrace$ is not satisfiable in $\E$};\\
		&\Longleftrightarrow &\text{for some $X_0\Subset X$, the set $X_0\cup\lbrace\Neg\epsilon\rbrace$ is not satisfiable}\\
		&&\text{in $\E$ [by Proposition~\ref{P:compactness-thm-refined}]};\\
		&\Longleftrightarrow &\text{there is $X_0\Subset X$ such that for any $\mat{M}\in\E$ and any valuation $v$,}\\
		&&\text{$\mat{M}\models_v X$ implies $\mat{M}\models_v\epsilon$};\\
		&\Longleftrightarrow &\text{there is $X_0\Subset X$ such that
			$X_0\vdashE\epsilon$}.
	\end{array}
	\]
\end{proof}

The next corollary is straightforwardly derived from the last propostion.
\begin{cor}\label{C:EE-finitary}
	The abstract logic $\hat{\E}$ is finitary.
\end{cor}

In the sequel (Section~\ref{section:equational-based-on-implicational}),  we will see that Proposition~\ref{P:alpha=alpha} is not true for some abstract logics $\E$. Now we show that Proposition~\ref{P:finitariness-E-consequence} also is not true for some abstract logics $\E$.

Indeed, let $\mat{M}:=\langle\alg{A},\varDelta\rangle$ be an E-matrix, where $\alg{A}$ is the algebra of the matrix $\mat{LC}^{\ast}$ of Section~\ref{section:dummett}, extended with an additional 0-ary operation $\one$ (the greatest element of $|\alg{A}|$). And $\E$ be the abstract logic determined by this matrix.

We remind that \alg{A} is an algebra of type $\Lan_{A}$ extended with an additional constant $\top$ which is interpreted in $\alg{A}$ by $\one$.
Now we define:
\[
X:=\set{(p_i\leftrightarrow p_j)\rightarrow p_0\approx\top}{0<i<j<\omega}.
\]

We aim to show that
\begin{equation*}\label{E:relative-equivalence-counterexample-1}
	X\vdashE p_0\approx\top.
\end{equation*}

Let $v$ be an arbitrary valuation in $\alg{A}$. By premise, for any $i$ and $j$ with $0<i<j<\omega$,
\[
v[p_i]\leftrightarrow v[p_j]\le v[p_0].
\]

If for at least one pair $i$ and $j$, $v[p_i]=v[p_j]$, then obviously $v[p_0]=\one$, that is $\mat{M}\models_{v}p_0\approx\top$.

Now assume that for pairs $i$ and $j$ with  $0<i<j<\omega$, $v[p_i]\neq v[p_j]$.
Then the set $\set{v[p_i]}{0<i<\omega}$ is infinite and can be arranged as an increasing chain with respect to $<$. We observe that $v[p_0]$ is an upper bound of this set. Since there only one upper bound of any infinite set in $|\alg{A}|$, namely $\one$, $v[p_0]=\one$. 

On the other hand, it is obvious that
\[
\varnothing\not\vdashE p_0\approx\top.
\]

Further, let $X_0$ be any nonempty finite subset of $X$; that is, $X_0$ is included in the set
\[
X_1:=\set{(p_i\leftrightarrow p_j)\rightarrow p_0\approx\top}{0<i<j<n},
\]
for some positive integer $n$. We prove that
\[
X_1\not\vdashE p_0\approx\top,
\]
which, by monotonicity, implies that
\[
X_0\not\vdashE p_0\approx\top
\]
either.

Let us take any valuation $v$ in $\alg{A}$ such that for any $i$ with
$0<i<n$, $v[p_i]$ equals the $i$th element with respect to the linear order $\le$ in $\alg{A}$ and $v[p_0]$ equals the $n$th element. It is obvious that
$\mat{M}\models_{v}X_1$ and $\mat{M}\not\models_{v}p_0\approx\top$.\\

Given a nonempty class $\E$ of E-matrices, we denote by $\ThmE$ the set of $\E$-theorems; that is
\[
\ThmE=\ConE(\varnothing).
\]

Since, by definition, the logical filter of a binary matrix is a binary relation, we need the following adjuster:
\[
\textbf{t}(\alpha\approx\beta):=(\alpha,\beta);
\]
and for any $X\subseteq\EqL$,
\[
X^{\textbf{t}}:=\set{\textbf{t}(\epsilon)}{\epsilon\in X}.
\]
Sometimes it will be convenient to write $[X]^{\textbf{t}}$ instead of $X^{\textbf{t}}$.\\

Using this notation, we may seem somewhat pedantic when we note:
\[
\textbf{t}(\epsilon)\in X^{\textbf{t}}\Longleftrightarrow\epsilon\in X.
\]

By means of the above adjuster, we define \textit{\textbf{Lindenbaum matrices}} as follows:\index{E-matrix!Lindenbaum}
\begin{equation}\label{E:Lind-matrix-for-EE}
	\LinE[X]:=\langle\FormAl,[\ConE(X)]^{\textbf{t}}\rangle,
\end{equation}
for any $X\subseteq\EqL$; respectively, we denote:
\[
\LinE:=\langle\FormAl,[\ConE(\varnothing)]^{\textbf{t}}\rangle.
\]

We note that, in virtue of Proposition~\ref{P:alpha=alpha}, the matrix $\LinEE$ is equational. On the other hand, the matrix $\LinEE[X]$ is equational if, and only if, $X\subseteq\ConEE(\varnothing)$; thus, for $\hat{\E}$, the only Lindenbaum matrix which is equational (up to isomorphism) is $\LinEE$.\\

Further, we observe that any valuation in $\LinE[X]$ is an $\Lan$-substitution. This observation leads to an analogue of Lindenbaum's theorem (Proposition~\ref{P:lindenbaum-theorem}).

\begin{prop}\label{P:lindenbaum-theorem-for-E}
	For any set {\em$X\cup\lbrace\epsilon\rbrace\subseteq\EqL$}, the following conditions are equivalent.
	{\em\[
		\begin{array}{cl}
			(\text{a}) &X\vdashE\epsilon;\\
			(\text{b}) &X\models_{\LinE[X]}\epsilon;\\
			(\text{c}) &\LinE[X]\models_{\iota}\epsilon.
		\end{array}
		\]}
	
	Consequently,
	{\em\[
		\epsilon\in\ThmE \Longleftrightarrow \LinE\models\epsilon.
		\]}
\end{prop}
\begin{proof}
	We will prove the following implications: (a)$\Rightarrow$(b)$\Rightarrow$(c)$\Rightarrow$(a).
	
	Case:~(a)$\Rightarrow$(b). Assume that $X\vdashE\epsilon$ and let $\sigma$ be any $\Lan$-substitution. If $\sigma(X)\subseteq\ConE(X)$, then, first, 
	$\ConE(\sigma(X))\subseteq\ConE(X)$ and, second, since $\vdashE$ is structural
	(Corollary~\ref{C:E-structural-con}), $\sigma(\epsilon)\in\ConE(\sigma(X))$ and hence $\sigma(\epsilon)\in\ConE(X)$.
	
	Case:~(b)$\Rightarrow$(c) is obvious, for $\iota(X)=X\subseteq\ConE(X)$.
	
	Case:~(c)$\Rightarrow$(a). Assume that $\iota(\epsilon)\in\ConE(X)$, that is $\epsilon\in\ConE(X)$, which in turn means that $X\vdashE\epsilon$.
\end{proof}

A set $X\subseteq\EqL$ is called an $\E$-\textit{\textbf{theory}}\index{$\E$-tehory} if
$\ConE(X)=X$.
We denote:
\[
\theoryE:=\set{X}{X~\text{is an $\E$-theory}}
\]
and
\[
[\theoryE]^{\textbf{t}}:=\set{X^{\textbf{t}}}{X~\text{is an $\E$-theory}}
\]

The last notion initiates the following definition.
\begin{defn}[Lindenbaum $\E$-bundle/$\E$-atlas]\index{E-matrix!Lindenbaum $\E$-bundle}
	The bundle {\em$\lbrace\langle\FormAl,D^{\textbf{t}}\rangle\rbrace_{D\in\theoryE}$} is called a \textbf{Lindenbaum} $\E$-\textbf{bundle} and the atlas {\em$\Lin[\theoryE]:=
		\langle\FormAl,[\theoryE]^{\textbf{t}}\rangle$} is a \textbf{Lindenbaum} $\E$-\textbf{atlas}.	
\end{defn}

As earlier, we define for any set $X\cup\lbrace\epsilon\rbrace\subseteq\EqL$:
\begin{equation}\label{E:Lindenbaum-E-bundle}
	\begin{array}{rl}
		X\models_{\textsf{\textbf{Lin}}[\Sigma_{\mathcal{S}}]}\epsilon \stackrel{\text{df}}{\Longleftrightarrow}
		&\!\!\!\textit{for any $\Lan$-substitution $\sigma$ and any $D\in\Sigma_{\mathcal{E}}$},\\
		&\!\!\!\textit{$\sigma(X)\subseteq D$ implies $\sigma(\epsilon)\in D$}.
	\end{array}
\end{equation}
\begin{prop}\label{P:Lindenbaum-E-completeness}
	For any set {\em$X\cup\lbrace\epsilon\rbrace\subseteq\EqL$},
	{\em	\[
		X\vdashE\epsilon\Longleftrightarrow X\models_{\textsf{\textbf{Lin}}[\Sigma_{\mathcal{E}}]}\epsilon.
		\]}
\end{prop}
\begin{proof}
	Suppose $X\vdashE\epsilon$.	 Let $D$ be an $\E$-theory. Now assume that, given a substitution $\sigma$, $\sigma(X)\subseteq D$. This implies that $\ConE(\sigma(X))\subseteq D$ as well. Also, since $\E$ is structural,
	$\sigma(X)\vdashE\sigma(\epsilon)$. Hence, $\sigma(\epsilon)\in D$.
	
	Conversely, suppose that $X\models_{\textsf{\textbf{Lin}}[\Sigma_{\mathcal{E}}]}\epsilon$. This implies that for the unit substitution $\iota$ and the theory $D=\ConE(X)$, if $\iota(X)\subseteq D$, then $\iota(\epsilon)\in D$. Hence $X\vdashE\epsilon$.
\end{proof}

\paragraph{Exercises~\ref{section:E-consequence}}
\begin{enumerate}
	\item \label{EX:binary-matrix-structural} Complete the proof of Proposition~\ref{P:binary-matrix-structural}.
	\item \label{EX:remark}Prove Proposition~\ref{P:los-suszko} and Proposition~\ref{P:wojcicki} in the context of equational consequence.
	\item \label{EX:Con(0)_E-closed-under-sub}Prove Corollary~\ref{C:Con(0)_E-closed-under-sub}.
\end{enumerate}

\section{$\hat{\E}$ as a deductive system}\label{section:E-deductive-systems}
In this section we will show, how the abstract logic $\hat{\E}$ can be defined via inference rules. Below we define three deductive systems. It turns out that each of them defines a consequence relation equal to $\hat{\E}$. In our exposition we employ the procedure of replacement that was defined in Section~\ref{section:languages}.
\paragraph{System E1:}
\[
\begin{array}{rllll}
	\text{Axiom:} &\bm{\alpha}\approx\bm{\alpha}.\\\\
	\text{Rules of inference:} &a)~\dfrac{\bm{\alpha}\,\approx\,\bm{\beta}}{\bm{\beta}\,\approx\,\bm{\alpha}}
	&b)~\dfrac{\bm{\alpha}\,\approx\,\bm{\beta},~\bm{\beta}\,\approx\,\bm{\gamma}}{\bm{\alpha}\,\approx\,\bm{\gamma}}
	&c)~\dfrac{\bm{\alpha}_{1}\,\approx\,\bm{\beta}_1,\ldots,~\bm{\alpha}_{n}\,\approx\,\bm{\beta}_n}{F\bm{\alpha}_1\ldots\bm{\alpha}_{n}\,\approx\, F\bm{\beta}_1\ldots\bm{\beta}_{n}},\\\\
	&&&~~~~\text{for any $F\in\Func$ of arity $n>0$}.
\end{array}
\]
\paragraph{System E2:}
\[
\begin{array}{rll}
	\text{Axiom:} &\bm{\alpha}\approx\bm{\alpha}.\\\\
	\text{Rules of inference:} &a)~\dfrac{\bm{\alpha}\,\approx\,\bm{\beta}}{\bm{\beta}\,\approx\,\bm{\alpha}}
	&b)~\dfrac{\bm{\epsilon}[\bm{\alpha}_1,\ldots,\bm{\alpha}_{n}],~\bm{\alpha}_{1}\,\approx\,\bm{\beta}_1,\ldots,~\bm{\alpha}_{n}\,\approx\,\bm{\beta}_n}{\bm{\epsilon}[\bm{\beta}_1,\ldots,\bm{\beta}_{n}]},\\\\
	&&~~~~\text{for any metaequality $\bm{\epsilon}[\bm{\alpha}_1,\ldots,\bm{\alpha}_{n}]$}.
\end{array}
\]
\paragraph{System E3:}
\[
\begin{array}{rll}
	\text{Axiom:} &\bm{\alpha}\approx\bm{\alpha}.\\\\
	\text{Rules of inference:} &a)~\dfrac{\bm{\alpha}\,\approx\,\bm{\beta}}{\bm{\beta}\,\approx\,\bm{\alpha}}
	&b)~\dfrac{\bm{\epsilon}[\bm{\alpha}],~\bm{\alpha}\,\approx\,\bm{\beta}}{\bm{\epsilon}[\bm{\beta}]},\\\\
	&&~~~~\text{for any metaequality $\bm{\epsilon}[\bm{\alpha}]$}.
\end{array}
\]

The notion of  formal derivation is similar for all three E1--E3. Namely, let $\text{E}^{\ast}$ stand for any of E1--E3.  A finite list of equalities
\begin{equation}\label{E:derivation-in-E1-3}
	\epsilon_1,\ldots,\epsilon_n
\end{equation}
is a \textit{\textbf{formal derivation}} (or simply a \textbf{derivation})\index{E-logic!derivation} in $\text{E}^{\ast}$ from a set $X$ of equalities if each $\epsilon_i$ is either an instantiation of the axiom of $\text{E}^{\ast}$ or belongs to $X$ or can be obtained from any suitable group of the preceding equalities $\epsilon_1,\ldots,\epsilon_{i-1}$, by one of the inference rules of $\text{E}^{\ast}$. We say that an equality $\epsilon$ is (\textit{\textbf{formally}}) {\textit{\textbf{derivable from}}} $X$, in symbols $X\vdash_{\text{E}^{\ast}}\epsilon$,\index{$X\vdash_{\text{E}^{\ast}}\epsilon$} if there is a derivation~\eqref{E:derivation-in-E1-3}, where $\epsilon=\epsilon_n$.
Given a derivation, the number of the formulas in it is called the \textit{\textbf{length}} of this derivation.
\begin{rem}\label{R:remark-on-E1-E2-E3}
	{\em
		Prior to Proposition~\ref{P:EE-consequence-completeness}, we do not claim that the relations `$\vdash_{\text{E}1}$', `$\vdash_{\text{E2}}$' and`$\vdash_{\text{E3}}$' are examples of equational consequence, although Proposition~\ref{P:EE-consequence-completeness} claims that they are.
	}
\end{rem}
\begin{prop}[variables preservation property]\label{P:variable-preservation}
	Let {\em$X\vdash_{\text{E}^{\ast}}\epsilon$}. Then 
	$\Var(\epsilon)\subseteq\Var(X)$ if, and only if, there is a derivation of $\epsilon$ from $X$ in which any instantiation of the axiom, $\alpha\approx\alpha$, satisfies the condition $\Var(\alpha)\subseteq\Var(X)$.
\end{prop}
\noindent\textit{Proof}~is left to the reader. (Exercise~\ref{section:E-deductive-systems}.\ref{EX:variable-preservation})

\begin{prop}\label{P:formal-derivation-equivalence}
	For any set $X\cup\lbrace\epsilon\rbrace\subseteq\EqL$, the following conditions are equivalent.
	{\em\[
		\begin{array}{rl}
			(\text{a}) &X\vdash_{\text{E1}}\epsilon;\\
			(\text{b}) &X\vdash_{\text{E2}}\epsilon;\\
			(\text{c}) &X\vdash_{\text{E3}}\epsilon.
		\end{array}
		\]}
\end{prop}
\begin{proof}
	We prove in the following order: (a)$\Rightarrow$(b)$\Rightarrow$(c)$\Rightarrow$(a).
	
	Case:~(a)$\Rightarrow$(b). Suppose that the sequence~\eqref{E:derivation-in-E1-3} is a derivation E1. Assume that one of the equalities of this sequence, $\alpha\approx\gamma$, is obtained from equalities $\alpha\approx\beta$ and $\beta\approx\gamma$, preceding it in~\eqref{E:derivation-in-E1-3}, by the rule $(b)$. Then we treat $\alpha\approx\beta$ as $\epsilon[\beta]$ and obtain $\alpha\approx\gamma$ by
	replacement $\epsilon[\beta], \beta\approx\gamma\slash\epsilon[\gamma]$, which is an instance of the rule (E2--$b$).
	
	Now suppose~\eqref{E:derivation-in-E1-3} has an application of the rule (E1--$c$). And assume that 
	\[
	\dfrac{\alpha_{1}\,\approx\,\beta_1,\ldots,~\alpha_{n}\,\approx\,\beta_n}{F\alpha_1\ldots\alpha_{n}\,\approx\, F\beta_1\ldots\beta_{n}}
	\]
	is the first such application. Then we insert the equality $F\alpha_1\ldots\alpha_{n}\,\approx\, F\alpha_1\ldots\alpha_{n}$, which is an instance of the axiom, right before the equality $F\alpha_1\ldots\alpha_{n}\,\approx\, F\beta_1\ldots\beta_{n}$. Now we treat the former equality as $\epsilon[\alpha_{1},\ldots,\alpha_{n}]$, where the designated occurrences of terms appear after the equality sign $\approx$. Next we apply the rule (E2--$c$) to obtain $F\alpha_1\ldots\alpha_{n}\,\approx\, F\beta_1\ldots\beta_{n}$. Since we may only insert new equalities in~\eqref{E:derivation-in-E1-3}, the resulting sequence will be a derivation of $\epsilon_n$ from $X$ in E2.
	
	Case:~(b)$\Rightarrow$(c) is obvious, for (E2--$b$) is a generalization of (E3--$b$).
	
	Case:~(c)$\Rightarrow$(a). Now we assume that~\eqref{E:derivation-in-E1-3} is a derivation in E3. And suppose
	\[
	\dfrac{\epsilon[\alpha], \alpha\approx\beta}{\epsilon[\beta]}
	\]
	is the first application of this rule in~\eqref{E:derivation-in-E1-3}. Now, by induction on the degree of the complexity of $\epsilon$, we prove that this step can be replaced with the axiom and rules of E1. Indeed, assume that
	$\epsilon[\alpha]:=\alpha\approx\gamma$. Then, using (E1--$a$) and (E1--$b$), we obtain:
	\[
	\epsilon_1,\ldots,\alpha\approx\gamma,\ldots,\alpha\approx\beta,\beta\approx\alpha, \beta\approx\gamma,
	\] 
	which is a derivation in E1. The case, when $\epsilon[\alpha]:=\gamma\approx\alpha$, is considered analogously.
	
	Next, we assume that
	\[
	\epsilon[\alpha]:=F\zeta_1[\alpha]\ldots\zeta_m[\alpha]\approx
	G\eta_1[\alpha]\ldots\eta_k[\alpha]. 
	\]
	
	Then, before the first application of (E3-$b$) in~\eqref{E:derivation-in-E1-3}, we insert the sequence:
	\[
	\zeta_1[\alpha]\approx\zeta_1[\alpha],\zeta_1[\alpha]\approx\zeta_1[\beta],
	\ldots,\zeta_m[\alpha]\approx\zeta_m[\alpha],\zeta_m[\alpha]\approx\zeta_m[\beta]
	\]
	and 
	\[
	\eta_1[\alpha]\approx\eta_1[\alpha],\eta_1[\alpha]\approx\eta_1[\beta],
	\ldots,\eta_k[\alpha]\approx\eta_k[\alpha],\eta_k[\alpha]\approx\eta_k[\beta].
	\]
	
	Let us focus on the sequence:
	\[
	\epsilon_1,\ldots,\alpha\approx\beta,\ldots,\zeta_i[\alpha]\approx\zeta_i[\alpha],\zeta_i[\alpha]\approx\zeta_i[\beta],
	\]
	which is an initial subsequence of the sequence obtained from~\eqref{E:derivation-in-E1-3} by such an insertion.
	We see that it is an E3-derivation, where the last equality is obtained from
	$\alpha\approx\beta$ and $\zeta_i[\alpha]\approx\zeta_i[\alpha]$ by (E3--$b$).
	By induction hypothesis, this application of (E3--$b$) can be transformed to an E1-derivation. Applying this argument to all equalities $\zeta_i[\alpha]\approx\zeta_i[\beta]$ and all equalities $\eta_i[\alpha]\approx\eta_i[\beta]$, we, then, apply (E1-$c$) to obtain that
	\[
	X\vdash_{\text{E}1}F\zeta_1[\alpha]\ldots\zeta_m[\alpha]\approx F\zeta_1[\beta]\ldots\zeta_m[\beta]
	\]
	and
	\[
	X\vdash_{\text{E}1}G\eta_1[\alpha]\ldots\eta_k[\alpha]\approx G\eta_1[\beta]\ldots\eta_k[\beta].
	\]
	Finally, by means of application of (E1--$a$) and (E1--$b$), we obtain:
	\[
	X\vdash_{\text{E}1}F\zeta_1[\beta]\ldots\zeta_m[\beta]\approx
	G\eta_1[\beta]\ldots\eta_k[\beta].
	\]
\end{proof}

\begin{prop}\label{P:E-soundness}
	Let $\E$ be any nonempty class of E-matrices.
	Suppose $\dfrac{\epsilon_1,\ldots,\epsilon_n}{\epsilon}$, where $n\ge 1$, is one of the inference rules of any of {\em E1--E3}. Then $X\vdashE\epsilon$, whenever $X\vdashE\epsilon_{1},\ldots,X\vdashE\epsilon_{n}$. Consequently,
	{\em\[
		X\vdash_{\text{E}^{\ast}}\epsilon\Longrightarrow X\vdashE\epsilon.
		\]}
\end{prop}
\noindent\textit{Proof}~is left to the reader. (Exercise~\ref{section:E-deductive-systems}.\ref{EX:E[M]-soundness})\\

Given a set $X\cup\lbrace\alpha\approx\beta\rbrace\subseteq\EqL$, we define:
\begin{equation}\label{E:defn-congruence}
	(\alpha,\beta)\in\theta_{\hat{\mathcal{E}}}(X)~\stackrel{\text{df}}{\Longleftrightarrow}~X\vdash_{\text{E}^{\ast}}\alpha\approx\beta.
\end{equation}
\begin{prop}\label{P:theta_X=congruence}
	Given a set {\em$X\subseteq\EqL$}, the relation $\theta_{\hat{\mathcal{E}}}(X)$ is a congruence on  {\em$\FormsL$}.	
\end{prop}
\begin{proof}
	To prove this proposition we apply Proposition~\ref{P:formal-derivation-equivalence} and the definition of E1.
\end{proof}

We aim to prove the following.
\begin{prop}[soundness and completeness]\label{P:EE-consequence-completeness}
	For any set {\em$X\cup\lbrace\epsilon\rbrace\subseteq\EqL$},
	{\em\[
		X\vdash_{\text{E}^{\ast}}\epsilon~\Longleftrightarrow~X\vdashEE\epsilon.
		\]}
\end{prop}
\begin{proof}
	The $\Rightarrow$-implication (soundness) is the second part of Proposition~\ref{P:E-soundness}.
	
	Now we prove the $\Leftarrow$-implication (completeness). Let
	\[
	\epsilon:=\alpha\approx\beta
	\]
	and $X\not\vdash_{\text{E}^{\ast}}\alpha\approx\beta$. Then $(\alpha,\beta)\notin\theta_{\hat{\mathcal{E}}}(X)$ and, in virtue of Proposition~\ref{P:theta_X=congruence}, $\alpha\slash\theta_{\hat{\mathcal{E}}}(X)\neq\beta\slash\theta_{\hat{\mathcal{E}}}(X)$ in the algebra $\FormAl\slash\theta_{\hat{\mathcal{E}}}(X)$.
	
	For convenience, we denote
	\[
	\alg{A}_0:=\FormAl\slash\theta_{\hat{\mathcal{E}}}(X).
	\]
	
	Let us consider a valuation
	\[
	v:\FormAl\longrightarrow\FormAl\slash\theta_{\hat{\mathcal{E}}}(X):~p\mapsto p\slash\theta_{\hat{\mathcal{E}}}(X),
	\]
	for any $p\in\VarL$.
	
	It must be clear that for any $\epsilon^{\prime}\in X$, $v(\epsilon^{\prime})\in\varDelta_{\textsf{A}_0}$, but $(\alpha\slash\theta_{\hat{\mathcal{E}}}(X),\beta\slash\theta_{\hat{\mathcal{E}}}(X))\!\notin\!\varDelta_{\textsf{A}_0}$. Hence $X\not\vdashE\epsilon$.
\end{proof}

\begin{cor}\label{C:congrunece-via-conE}
	Given a set {\em$X\subseteq\EqL$}, for any equality $\alpha\approx\beta$,
	{\em	\[
		(\alpha,\beta)\in\theta_{\hat{\mathcal{E}}}(X)\Longleftrightarrow \alpha\approx\beta\in\ConEE(X).
		\]}
\end{cor}
\noindent\textit{Proof}~follows immediately from Proposition~\ref{P:EE-consequence-completeness}.\\

We note that Corollary~\ref{C:EE-finitary} also follows from Proposition~\ref{P:EE-consequence-completeness}.

\paragraph{Exercises~\ref{section:E-deductive-systems}}
\begin{enumerate}
	\item \label{EX:variable-preservation}Prove Proposition~\ref{P:variable-preservation} by induction on the length of derivation in $\text{E}^{\ast}$.
	\item \label{EX:E[M]-soundness}Prove Proposition~\ref{P:E-soundness}.	
\end{enumerate}

\section{Lindenbaum-Tarski matrices for logic $\E$}\label{section:Lindenbaum-Tarski-matrices}
We draw the reader's attention to the fact that Lindenbaum matrices of the definition~\eqref{E:Lind-matrix-for-EE} are not necessarily E-matrices. On the other hand, we remind the reader that the concept of Lindenbaum-Tarski algebra, which was introduced for the unital abstract logics in
Section~\ref{section:lindenbaum-algebra}, turned out to be a convenient tool for constructing separating tools. (The latter conception was discussed in Section~\ref{section:separating-means}.) In this section we implement the same idea in relation to equational consequence.\\

For an arbitrary equational logic $\E$, we define:
\begin{equation}\label{E:defn-theta_E[M]}
	(\alpha,\beta)\in\theta_{\mathcal{E}}(X)\stackrel{\text{df}}{\Longleftrightarrow} X\vdashE\alpha\approx\beta,
\end{equation}
for any set $X\subseteq\EqL$; in case $X=\varnothing$, we denote
\[
\theta_{\mathcal{E}}:=\theta_{\mathcal{E}}(\varnothing).
\]

We note that for any $X\cup Y\subseteq\EqL$,
\begin{equation}\label{E:theta(X)-included-theta(Y)}
	X\subseteq Y \Longrightarrow\theta_{\mathcal{E}}(X)\subseteq\theta_{\mathcal{E}}(Y);
\end{equation}
in particular, it is obvious that
\begin{equation}\label{E:theta_M-implies-thets_M[X]}
	(\alpha,\beta)\in\theta_{\mathcal{E}}\Longrightarrow(\alpha,\beta)\in\theta_{\mathcal{E}}(X).
\end{equation}

According to Proposition~\ref{P:E-soundness}, each $\theta_{\mathcal{E}}(X)$ is a congruence on $\FormAl$. 
\begin{defn}[Lindenbaum-Tarski matrices for $\E$]\index{E-matrix!Lindenbaum-Tarski matrix}
	Let $\E$ be an equational logic. For an arbitrary set {\em$X\subseteq\EqL$}, the E-matrix {\em$\LTE[X]:=\langle\FormAl\slash\theta_{\mathcal{E}}(X),\varDelta\rangle$} is called a \textbf{Lindenbaum-Tarski matrix for $\E$ relative to} $X$. We denote {\em$\LTE:=\LTE[\varnothing]$} and call the latter simply a \textbf{Lindenbaum-Tarski matrix for $\E$}.
\end{defn}

We note that, in virtue of Proposition~\ref{P:alpha=alpha}, the map $\alpha\mapsto\lbrace\alpha\rbrace$ is an isomorphism of $\FormAl$ onto $\LTE$-algebra.\label{pagenote-1} (See Exercise~\ref{section:Lindenbaum-Tarski-matrices}.\ref{EX:isomorphism}.)
\begin{prop}\label{P:homomorphism-LT-onto-LT[X]}
	Let $\E$ be an equational logic and {\em$X\subseteq\EqL$}. Then the map $p\slash\theta_{\mathcal{E}}\mapsto p\slash\theta_{\mathcal{E}}(X)$, where $p\in\VarL$, can be uniquely extended to an epimorphism of {\em$\LTE$} onto {\em$\LTE[X]$}.
\end{prop}
\begin{proof}
	We use~\eqref{E:theta_M-implies-thets_M[X]} and Proposition~\ref{P:congruence-thm-auxiliary}. (See Exercise~\ref{section:Lindenbaum-Tarski-matrices}.\ref{EX:homomorphism-LT-onto-LT[X]}.)
\end{proof}

The following lemma is an adaptation of Proposition~\ref{P:valuation-in-F_L/theta}.
\begin{lem}\label{L:equivalence-with-LT[X]}
	Let $\E$ be an equational logic and {\em$X\subseteq\EqL$}. For any valuation $v$ in $\FormAl\slash\theta_{\mathcal{E}}(X)$, there is an $\Lan$-substitution $\sigma$ such that for any term $\alpha$,
	$v[\alpha]=\sigma(\alpha)\slash\theta_{\mathcal{E}}(X)$. Conversely, for any $\Lan$-substitution $\sigma$, there is a unique valuation $v_\sigma$ in $\FormAl\slash\theta_{\mathcal{E}}(X)$ such that for any term $\alpha$, $v_\sigma[\alpha]=\sigma(\alpha)\slash\theta_{\mathcal{E}}(X)$.
\end{lem}

Now we define a valuation in $\FormAl\slash\theta_{\mathcal{E}}(X)$:
\[
v_{X}[p]:=p\slash\theta_{\mathcal{E}}(X),~\text{for any $p\in\VarL$}.
\]
Extending $v_X$ on all $\FormAl$, we get:
\[
v_{X}[\alpha]=\alpha\slash\theta_{\mathcal{E}}(X),~\text{for any $\alpha\in\FormAl$}.
\]
\begin{prop}\label{P:equivalence-with-LT[X]}
	Let $\E$ be an equational logic. For any set {\em$X\cup\lbrace\epsilon\rbrace\subseteq\EqL$}, the following conditions are equivalent.
	{\em\[
		\begin{array}{cl}
			(\text{a}) &X\vdashE\epsilon;\\
			(\text{b}) &X\models_{\LTE[X]}\epsilon;\\
			(\text{c}) &\LTE[X]\models_{v_{X}}\epsilon.
		\end{array}
		\]}
\end{prop}
\begin{proof} 
	We prove the desired equivalence in the following order: (a)$\Rightarrow$(b)$\Rightarrow$(c)$\Rightarrow$(a).
	
	Suppose that $\epsilon=\alpha\approx\beta$. 
	
	Case:~(a)$\Rightarrow$(b). Assume that $X\vdashE\alpha\approx\beta$. Let us take any valuation $v$ in $\FormAl\slash\theta_{\mathcal{E}}(X)$. Now we take any valuation $v$ in $\FormAl$ such that for any $\alpha^{\prime}\approx\beta^{\prime}\in X$, $v[\alpha^{\prime}]=v[\beta^{\prime}]$. In virtue of Lemma~\ref{L:equivalence-with-LT[X]}, there is a substitution $\sigma$ such that $\sigma_v(\alpha^{\prime})\slash\theta_{\mathcal{E}}(X)=\sigma(\beta^{\prime})\slash\theta_{\mathcal{E}}(X)$, that is $\sigma(\alpha^{\prime}\approx\beta^{\prime})\in\ConE(X)$. This implies that $\ConE(\sigma(X))\subseteq\ConE(X)$. Since $\E$ is structural
	(Corollary~\ref{C:E-structural-con}),
	$X\vdashE\sigma(\alpha)\approx\sigma(\beta)$. This in turn implies that
	$v[\alpha]=v[\beta]$.
	
	Case:~(b)$\Rightarrow$(c). It must be clear that for any term $\alpha^{\prime\prime}$, $v_{X}[\alpha^{\prime\prime}]=\alpha^{\prime\prime}\slash\theta_{\mathcal{E}}(X)$. Since for any $\alpha^{\prime}\approx\beta^{\prime}\in X$, in virtue of Proposition~\ref{P:EE-consequence-completeness} and~\eqref{E:defn-congruence},
	$v_{X}[\alpha^{\prime}]=v_{X}[\beta^{\prime}]$, by means of (b), we conclude that
	$v_{X}[\alpha]=v_{X}[\beta]$.
	
	Case:~(c)$\Rightarrow$(a). Suppose that $v_{X}[\alpha]=v_{X}[\beta]$, that is
	$\alpha\slash\theta_{\mathcal{E}}(X)=\beta\slash\theta_{\mathcal{E}}(X)$. According to~\eqref{E:defn-theta_E[M]},
	this implies that $X\vdashE\alpha\approx\beta$.
\end{proof}
\begin{prop}\label{P:extension-to-homomorphism-with-LT}
	Let $\E$ be an equational logic; also, let $X$ be a set of equalities and {\em$\mat{M}=\langle\alg{A},\varDelta\rangle$} be an E-matrix. Given a valuation $v$ in {\em$\alg{A}$}, {\em$\mat{M}\models_{v}X$} if, and only if, the map
	$f: p\slash\theta_{\mathcal{E}}(X)\mapsto v[p]$, for any $p\in\VarL$, can be extended to a homomorphism from $\FormAl\slash\theta_{\mathcal{E}}(X)$ to {\em$\alg{A}$}.
\end{prop}
\begin{proof}
	Suppose $\mat{M}\models_{v}X$. Let us define an extension of $f$ to all elements of $\FormAl\slash\theta_{\mathcal{E}}(X)$ as follows:
	\begin{equation}\label{E:Malcev's-homomorphism}
		f^{\ast}(\alpha\slash\theta_{\mathcal{E}}(X)):=v[\alpha],~\text{for any term $\alpha$}.
	\end{equation}
	
	This map is a homomorphism if the above definition is correct. To show correctness, we assume that $\alpha\slash\theta_{\mathcal{E}}(X)=\beta\slash\theta_{\mathcal{E}}(X)$.
	According to~\eqref{E:defn-theta_E[M]}, $X\vdashE\alpha\approx\beta$.
	This, by premise, implies that $v[\alpha]=v[\beta]$.
	
	Conversely, assume that $f: p\slash\theta_{\mathcal{E}}(X)\mapsto v[p]$ can be extended to a homomorphism. This means that \eqref{E:Malcev's-homomorphism} holds. Now let us take any equality $\alpha\approx\beta\in X$. Since $X\vdashE\alpha\approx\beta$, applying again the definition~\eqref{E:defn-theta_E[M]}, we derive $\alpha\slash\theta_{\mathcal{E}}(X)=\beta\slash\theta_{\mathcal{E}}(X)$. Then, by means of~\eqref{E:Malcev's-homomorphism}, we conclude that $v[\alpha]=v[\beta]$.
\end{proof}

So far, we have dealt with individual matrix algebras. Now we want to consider them as a class. 

Let $\E$ be a nonempty class of E-matrices of type $\Lan$. We denote:
\begin{center}
	\emph{$\A[\E]$ is the class of all matrix algebras of the class $\E$.}
\end{center}

Given a nonempty class of E-matrices $\E$, for every $X\subseteq\EqL$, we define:
\[
\Phi_{\mathcal{E}}(X):=\set{\theta\in\Congruence\,\FormAl}{X^{\textbf{t}}
	\subseteq\theta~\text{and}~\FormAl\slash\theta\in\Su\A[\E]}.
\]
\begin{lem}\label{L:LT-matrices-first}
	Let $\E$ be an equational logic and {\em$X\subseteq\EqL$}. If $\E$ contains a matrix that satisfies all equalities of $X$, then $\Phi_{\mathcal{E}}(X)\neq\varnothing$. 
\end{lem}
\begin{proof}
	Let $\E$ be an equational logic and a matrix $\mat{M}=\langle\alg{A},\varDelta\rangle$ of $\E$ satisfies all equalities of $X$. According to  Proposition~\ref{P:extension-to-homomorphism-with-LT}, there is a homomorphism of $\FormAl\slash\theta_{\mathcal{E}}(X)$ onto a subalgebra $\alg{B}$ of $\alg{A}$ (Proposition~\ref{P:homomorphic-image-is-subalgebra}). In virtue of Proposition~\ref{P:congruence-extension}, there is a congruence $\theta\in\Congruence\,\FormAl$ such that $\theta_{\mathcal{E}}(X)\subseteq\theta$ and $\FormAl\slash\theta$ is isomorphic to $\alg{B}$; that is $\theta\in\Phi_{\mathcal{E}}(X)$.
\end{proof}

Next we define:
\[
\theta_{\mathcal{A}[\mathcal{E}]}(X):=\bigcap\Phi_{\mathcal{E}}(X).
\]
\begin{lem}\label{L:LT-matrices-second}
	If $\E$ is an equational logic and $\E$ contains a matrix that satisfies all equalities of $X$, then	$\theta_{\mathcal{A}[\mathcal{E}]}(X)=\theta_{\mathcal{E}}(X)$.	
\end{lem}
\begin{proof}
	We need the condition of the lemma to make sure that $\Phi_{\mathcal{E}}$ is not empty.
	
	First we assume that $(\alpha,\beta)\in\theta_{\mathcal{E}}(X)$. Let us take any $\theta\in\Phi_{\mathcal{E}}(X)$ and define a valuation $v$ in $\FormAl\slash\theta$ as follows:
	\[
	v[p]:=p\slash\theta,~\text{for any $p\in\VarL$}.
	\]
	Now we take any equality $\gamma\approx\gamma^{\prime}\in X$. Since $(\gamma,\gamma^{\prime})\in\theta$, we conclude that $\langle\FormAl\slash\theta,\varDelta\rangle\models_{v}X$. By definition, there is an algebra $\alg{A}\in\mathcal{A}[\mathcal{E}]$ such that $\FormAl\slash\theta$ is a subalgebra of $\alg{A}$ (Proposition~\ref{P:homomorphic-image-is-subalgebra}). It must be clear that
	$\langle\alg{A},\varDelta\rangle\models_{v}X$ and hence, by premise,
	$\langle\alg{A},\varDelta\rangle\models_{v}\alpha\approx\beta$. Since $v[\alpha]$ and $v[\beta]$ belong to $|\FormAl\slash\theta|$, $(\alpha,\beta)\in\theta$. This allows us to conclude that $\theta_{\mathcal{E}}(X)\subseteq\theta_{\mathcal{A}[\mathcal{E}]}(X)$.
	
	Conversely, suppose that $(\alpha,\beta)\notin\theta_{\mathcal{E}}(X)$.
	This means that for some matrix $\mat{M}=\langle\alg{A},\varDelta\rangle$ of $\E$ and a valuation $v$ in $\alg{A}$, $\mat{M}\models_{v}X$ and $\mat{M}\not\models_{v}\alpha\approx\beta$, that is $v[\alpha]\neq v[\beta]$ in $\alg{A}$. By Proposition~\ref{P:extension-to-homomorphism-with-LT}, there is a homomorphism  $f^{\ast}$ of $\FormAl\slash\theta_{\mathcal{E}}(X)$ onto a subalgebra $\alg{B}$ of $\alg{A}$ such that $f^{\ast}(\gamma\slash\theta_{\mathcal{E}}(X))=v[\gamma]$, for any $\Lan$-term $\gamma$. This implies that there is a congruence $\theta$ such that $\theta_{\mathcal{E}}(X)\subseteq\theta$ and $\FormAl\slash\theta$ is isomorphic $\alg{B}$ (Proposition~\ref{P:congruence-extension}). It is obvious that $\alpha\slash\theta\neq\beta\slash\theta$. That is, $\theta\in\Phi_{\mathcal{E}}$ and $(\alpha,\beta)\notin\theta$. This implies that $(\alpha,\beta)\notin\theta_{\mathcal{A}[\mathcal{E}]}$.
\end{proof}

\begin{prop}[Birkhoff's theorem generalized]\label{P:Birkhoff's-generalized}
	Let {\em$X\subseteq\EqL$} and a class $\E$ contain a matrix which satisfies all equalities of $X$. Then {\em$\LTE[X]$}-algebra belongs to $\Su\Pro\mathcal{A}[\mathcal{E}]$. Thus if $\E$ contains a trivial matrix
	and $\Is\A[\E]\subseteq\A[\E]$, $\Su\A[\E]\subseteq\A[\E]$, and $\Pro\A[\E]\subseteq\A[\E]$, in particular if $\mathcal{A}[\mathcal{E}]$ is a quasi-variety, then {\em$\LTE[X]$}-algebra belongs to $\A[\E]$.
\end{prop}
\begin{proof}
	Taking into account Lemma~\ref{L:LT-matrices-second}, we then employ Proposition~\ref{P:Birkhpff's-for-algebras} directly.
\end{proof}

\paragraph{Exercises~\ref{section:Lindenbaum-Tarski-matrices}}
\begin{enumerate}
	\item \label{EX:isomorphism}Prove that the map $\alpha\mapsto\lbrace\alpha\rbrace$ is an isomorphism of $\FormAl$ onto $\LTE$-algebra.
	\item \label{EX:homomorphism-LT-onto-LT[X]} Give a detailed proof of Proposition~\ref{P:homomorphism-LT-onto-LT[X]}.
\end{enumerate}

\section{Mal'cev matrices}\label{section:malcev-matrices}
According to Proposition~\ref{P:equivalence-with-LT[X]}, given a set $X\cup\lbrace\epsilon\rbrace$ of equalities, the matrix $\LTE[X]$ is an adequate tool to determine whether $X\vdashE\epsilon$ or not. In case $\Var(\epsilon)\subseteq\Var(X)$, there is a more effective tool, especially when the set $\Var(X)$ is finite. This tool is a refinement of $\LTE[X]$ in the form of \emph{Mal'cev matrix}, $\ME[X]$ (see Definition~\ref{D:Malcev-matrix} below). Moreover, as Proposition~\ref{P:Birkhoff's-generalized} and Proposition~\ref{P:Malcev-second} show, given a class $\E$ of E-matrices,
if $\mathcal{A}[\mathcal{E}]$ satisfies certain closure conditions, both 
$\LTE[X]$ and $\ME[X]$ belong to $\E$.\\

Given a set $X\subseteq\EqL$, we denote by $\FormAl[\Var(X)]\slash\theta_{\mathcal{E}}(X)$ the subalgebra of $\FormAl\slash\theta_{\mathcal{E}[\mathcal{M}]}(X)$ which is generated by the elements of the set $\set{p\slash\theta_{\mathcal{E}}(X)}{p\in\Var(X)}$. 

\begin{defn}[Mal'cev matrices]\label{D:Malcev-matrix}\index{E-matrix!Mal'cev}
	Let $\E$ be a nonempty class of E-matrices.
	For any set {\em$X\subseteq\EqL$}, we define {\em$\ME[X]:=\langle\FormAl[\Var(X)]\slash\theta_{\mathcal{E}}(X),\varDelta\rangle$}, and call this matrix a \textbf{Mal'cev} \textbf{matrix relative to} $X$. We denote {\em$\ME:=\ME[\varnothing]$}.
\end{defn}

We note that $\ME=\LTE$. Hence, $\MEE$ is  is isomorphic to $\FormAl$; see Exercise~\ref{section:Lindenbaum-Tarski-matrices}.\ref{EX:isomorphism}.\\

Taking into account  Proposition~\ref{P:equivalence-with-LT[X]}, we derive the following.
\begin{prop}\label{P:equivalence-with-Mal[X]}
	For any set {\em$X\cup\lbrace\epsilon\rbrace\subseteq\EqL$} with $\Var(\epsilon)\subseteq\Var(X)$, the following conditions are equivalent.
	{\em\[
		\begin{array}{cl}
			(\text{a}) &X\vdashE\epsilon;\\
			(\text{b}) &X\models_{\ME[X]}\epsilon;\\
			(\text{c}) &\ME[X]\models_{v^{\prime}}\epsilon,~\text{providing that $v^{\prime}[p]=v_{X}(p)$, for any $p\in\Var(X)$}.
		\end{array}
		\]}
\end{prop}
\noindent\textit{Proof}~is left to the reader. (Exercise~\ref{section:malcev-matrices}.\ref{EX:equivalence-with-Mal[X]})\\

\begin{defn}\label{D:X-is-defining-for-M}\index{E-matrix!defining set}
	A set {\em$X\subseteq\EqL$} is called a \textbf{defining set for} a matrix {\em$\mat{M}=\langle\alg{A},\varDelta\rangle$} $($or $X$ \textbf{defines} $\mat{M}$$)$ \textbf{relative to} $\E$ if there is a valuation $v$ in {\em$\alg{A}$} such that
	{\em\[
		\begin{array}{cl}
			(\text{a}) &\text{$\alg{A}$ is generated by the elements of the set $\set{v[p]}{p\in\Var(X)}$};\\
			(\text{b}) &\mat{M}\models_{v}X;\\
			(\text{c}) &\text{for any equality $\epsilon$ with $\Var(\epsilon)\subseteq\Var(X)$, if $\mat{M}\models_{v}\epsilon$,
				then $X\vdashE\epsilon$}.
		\end{array}
		\]}
	We also say that $X$ is a \textbf{defining set for} {\em$\mat{M}$} $($or $X$ \textbf{defines} {\em$\mat{M}$}$)$ \textbf{with respect to} $v$.
\end{defn}
\begin{prop}\label{P:defining-set-1}
	For any set {\em$X\subseteq\EqL$}, $X$ defines {\em$\ME[X]$} relative to $\E$.
\end{prop}
\begin{proof}
	Let $v^{\prime}$ be any valuation in $\FormAl[\Var(X)]\slash\theta_{\mathcal{E}}(X)$ with
	\[
	v^{\prime}[p]=v_{X}[p],~\text{for any $p\in\Var(X)$}.
	\]
	Then, according to Proposition~\ref{P:equivalence-with-Mal[X]}, $X$ is a defining set for $\ME[X]$ with respect to $v^{\prime}$.
\end{proof}
\begin{prop}\label{P:defining-set-2}
	Let $X$ be a defining set for {\em$\mat{M}=\langle\alg{A},\varDelta\rangle$} relative to $\E$ with respect to a valuation $v$. Then the maps $f:\alpha\slash\theta_{\mathcal{E}}(X)\mapsto v[\alpha]$ and $g:v[\alpha]\mapsto:\alpha\slash\theta_{\mathcal{E}}(X)$, where $\Var(\alpha)\subseteq\Var(X)$, are mutually invertible  isomorphisms from $\FormAl[\Var(X)]\slash\theta_{\mathcal{E}}(X)$ onto {\em$\alg{A}$} and from {\em$\alg{A}$} onto $\FormAl[\Var(X)]\slash\theta_{\mathcal{E}}(X)$, respectively. Hence, {\em$\mat{M}$} and {\em$\ME[X]$} are isomorphic.	 
\end{prop}
\begin{proof}
	According to Proposition~\ref{P:extension-to-homomorphism-with-LT}, the map defined in~\eqref{E:Malcev's-homomorphism}, being  restricted to the terms built from $\Var(X)$, is an epimorphism of $\FormAl[\Var(X)]\slash\theta_{\mathcal{E}}(X)$ onto $\alg{A}$ and is coincident with $f$.
	
	On the other hand, if $v[\alpha]=v[\beta]$, then $\mat{M}\models_{v}\alpha\approx\beta$. This, by premise, implies that
	$X\vdashE\alpha\approx\beta$. Hence, by the definition~\eqref{E:defn-theta_E[M]}, $\alpha\slash\theta_{\mathcal{E}}(X)=\beta\slash\theta_{\mathcal{E}}(X)$. This shows that $g$ is defined correctly. It is a routine check that $g$ is an epimorphism. It is also obvious that $f$ and $g$ are inverses of one another.
\end{proof}

\begin{cor}
	All the E-matrices, for which $X$ is a defining set relative to $\E$, are isomorphic.
\end{cor}
\begin{prop}[Mal'cev's first theorem]\label{P:Malcev's-first}
	Let {\em$X\subseteq\EqL$} and {\em$\mat{M}=\langle\alg{A},\varDelta\rangle$} be an {\em E}-matrix satisfying the conditions {\em(a)} and {\em(b)} of Definition~\ref{D:X-is-defining-for-M} with respect to a valuation $v$ in {\em$\alg{A}$}.
	Then $X$ is a defining set for {\em$\mat{M}$} relative to $\E$ with respect to $v$ if, and only if, for any E-matrix {\em$\mat{N}=\langle\alg{B},\varDelta\rangle$}, if {\em$\mat{N}\models_{w}X$}, for some valuation $w$ in {\em$\alg{B}$}, then the map
	$h:v[p]\mapsto w[p]$, for any $p\in\Var(X)$, can be extended to a homomorphism
	of {\em$\alg{A}$} to {\em$\alg{B}$}.
\end{prop}
\begin{proof}
	First, we suppose that $X$ is a defining set for $\mat{M}$ with respect to $v$. Because of Proposition~\ref{P:defining-set-2}, without loss of generality, we can identify $\mat{M}$ with $\ME[X]$ and $v$ with any valuation $v^{\prime}_{X}$ that is coincident with $v_X$ for any $p\in\Var(X)$. Accordingly, we will have that $h: p\slash\theta_{\mathcal{E}}(X)\mapsto w[p]$. Now, we extend $h$ according to the rule
	\[
	h^{\ast}(\alpha)=w[\alpha],~\text{for any $\alpha$ with $\Var(\alpha)\subseteq\Var(X)$}. 
	\]
	
	The correctness of $h^{\ast}$ is proved as follows. Suppose $\alpha\slash\theta_{\mathcal{E}}(X)=\beta\slash\theta_{\mathcal{E}}(X)$. This means that $X\vdashE\alpha\approx\beta$. Then, by premise, $w[\alpha]=w[\beta]$. 
	
	That $h^{\ast}$ is a homomorphism is obvious. Since the conditions (a) and (b) of Definition~\ref{D:X-is-defining-for-M} are fulfilled, we have to show that the condition (c) is also true.
	
	Now we prove the `if'-part.
	For contradiction, assume that $\mat{M}\models_{v}\alpha\approx\beta$, but $X\not\vdashE\alpha\approx\beta$. The latter means that there is an E-matrix 
	$\mat{N}=\langle\alg{B},\varDelta\rangle$ and a valuation $w$ in $\alg{B}$ such that $\mat{N}\models_{w}X$ and $w[\alpha]\neq w[\beta]$. Using the homomorphism $h$ which exists by premise, we obtain: $h(v[\alpha])=w[\beta]$
	and $h(v[\beta])=w[\beta]$. This implies that $v[\alpha]\neq v[\beta]$. A contradiction.
\end{proof}
\begin{prop}\label{P:defining-set-3}
	Given an E-matrix {\em$\mat{M}=\langle\alg{A},\varDelta\rangle$} and an epimorphism
	{\em$g:\FormAl[\Var(X)]\slash\theta_{\mathcal{E}}(X)\longrightarrow\alg{A}$}, there is a valuation $v$ in {\em$\alg{A}$} such that the conditions {\em(a)} and {\em(b)} of Definition~\ref{D:X-is-defining-for-M} are satisfied.
\end{prop}
\begin{proof}
	Let $a$ be an element of $|\alg{A}|$. Then, we define:
	\[
	v[p]:=\begin{cases}
		\begin{array}{cl}
			g(p\slash\theta_{\mathcal{E}}(X)) &\text{if $p\in\Var(X)$}\\
			a &\text{if otherwise},
		\end{array}
	\end{cases}
	\]
	for any $p\in\VarL$.
	
	It must be clear that the condition (a) is fulfilled. According to Proposition~\ref{P:equivalence-with-Mal[X]}, $\ME[X]\models_{v^{\prime}_X}X$,
	where $v^{\prime}_X[p]=v_{X}[p]$, for any $p\in\Var(X)$. Since $\alg{A}$ is a homomorphic image of $\FormAl[\Var(X)]\slash\theta_{\mathcal{E}}(X)$ with respect to $g$, $\mat{M}\models_{v}X$. (See Exercise~\ref{section:malcev-matrices}.\ref{EX:defining-set-3}.)
\end{proof}

\begin{prop}[Dyck's theorem]\label{P:dyck-theorem}
	Let {\em$X\cup Y\subseteq\EqL$} with $\Var(X)\subseteq\Var(Y)$. If for every $\epsilon\in X$, $Y\vdashE\epsilon$, then the map $f:p\slash\theta_{\mathcal{E}}(X)\mapsto p\slash\theta_{\mathcal{E}}(Y)$ can be extended to an epimorphism of {\em$\ME[X]$}-algebra onto {\em$\ME[Y]$}-algebra. Moreover, if $\Var(X)=\Var(Y)$ and for any $\epsilon$ with $\Var(\epsilon)\subseteq\Var(X)$, the equivalence
	\[
	X\vdashE\epsilon\Longleftrightarrow Y\vdashE\epsilon
	\]
	holds, then the above epimorphism is an isomorphism.
\end{prop}
\begin{proof}
	We extend the map $f$ according to the rule:
	\[
	f^{\ast}(\alpha\slash\theta_{\mathcal{E}}(X)):=\alpha\slash\theta_{\mathcal{E}}(Y),
	\]
	for any term $\alpha$ with $\Var(\alpha)\subseteq\Var(X)$. The proof of the correctness of this definition is standard. Namely, assume that for terms $\alpha$ and $\beta$, where $\Var(\alpha)\cup\Var(\beta)\subseteq\Var(X)$,
	$\alpha\slash\theta_{\mathcal{E}}(X)=\beta\slash\theta_{\mathcal{E}}(X)$. This means that $X\vdashE\alpha\approx\beta$. By transitivity of $\vdashE$, $Y\vdashE\alpha\approx\beta$; that is $\alpha\slash\theta_{\mathcal{E}}(Y)=\beta\slash\theta_{\mathcal{E}}(Y)$.
	
	It is easy to see that $f^{\ast}$ is a homomorphism. (Exercise~\ref{section:malcev-matrices}.\ref{EX:dyck-theorem-1})
	
	The proof of the second part of the proposition is also left to the reader. (See Exercise~\ref{section:malcev-matrices}.\ref{EX:dyck-theorem-2}.)
\end{proof}

\begin{prop}[Tietze's theorem]\label{P:tietze's-theorem}
	Let $\E$ be a nonempty class of E-matrices.
	Also, let a language $\Lan^{\prime}$ of terms be a primitive extension of a language $\Lan$ and let $X$ be a set of equalities with $\Var(X)\subseteq\VarL$. Now we obtain a set $Y$ by adding to $X$ the set of equalities in the language $\Lan^{\prime}$ that are formed in the following way: for any term variable $p\in\Var_{\mathcal{L}^{\prime}}\setminus\VarL$, we select at random an $\Lan$-term $\alpha_p$ and form an equality $p\approx\alpha_p$. Then {\em$\ME[X]$} and {\em$\ME[Y]$} are isomorphic.
\end{prop}
\begin{proof}
	We note that the signatures of the languages $\Lan$ and $\Lan^{\prime}$ are the same. 
	
	We define a map $f^{\ast}$ from $\ME[Y]$-algebra to $\ME[X]$-algebra as follows. First, we define:
	\[
	f(p\slash\theta_{\mathcal{E}}(Y)):=\begin{cases}
		\begin{array}{cl}
			p\slash\theta_{\mathcal{E}}(X) &\text{if $p\in\VarL$}\\
			\alpha_p\slash\theta_{\mathcal{E}}(X) &\text{if $p\in\Var_{\mathcal{L}^{\prime}}\setminus\VarL$}.
		\end{array}
	\end{cases}
	\]
	
	Then, we extend $f$ in the following way:
	\[
	f^{\ast}(\alpha(p_1,\ldots,p_n)\slash\theta_{\mathcal{E}}(Y)):=\alpha(f(p_1),\ldots,f(p_n)),
	\]
	where $\alpha(p_1,\ldots,p_n)$ is any $\Lan^{\prime}$-term and 
	$\lbrace p_1,\ldots,p_n\rbrace=\Var(\alpha)$.
	
	First, we show that $f^{\ast}$ is defined correctly. For this, we take any two 
	$\Lan^{\prime}$-terms $\alpha(p_1,\ldots,p_k,q_1,\ldots,q_l)$ and
	$\beta(p_1,\ldots,p_k,q_1,\ldots,q_l)$, where $\lbrace p_1,\ldots,p_k\rbrace\subseteq\Var(X)$, $\lbrace q_1,\ldots,q_l\rbrace\subseteq\Var_{\mathcal{L}^{\prime}}\setminus\VarL$
	and $\lbrace p_1,\ldots,p_k,q_1,\ldots,q_l\rbrace=\Var(\alpha)\cup\Var(\beta)$.
	
	Now assume that $\alpha\slash\theta_{\mathcal{E}}(Y)=\beta\slash\theta_{\mathcal{E}}(Y)$, that is $Y\vdashE\alpha\approx\beta$. Since, by definition,
	\[
	Y\vdashE q_{1}\approx\alpha_{q_1},\ldots, Y\vdashE q_{l}\approx\alpha_{q_l},
	\]
	in virtue of (the first part of) Proposition~\ref{P:E-soundness}, we obtain that
	\[
	Y\vdashE\alpha(p_1,\ldots,p_k, \alpha_{q_1},\ldots,\alpha_{q_l})\approx\beta(p_1,\ldots,p_k,\alpha_{q_1},\ldots,\alpha_{q_l}).
	\]
	This, in virtue of Proposition~\ref{P:variables-auxiliary}, implies that
	\[
	X\vdashE\alpha(p_1,\ldots,p_k, \alpha_{q_1},\ldots,\alpha_{q_l})\approx\beta(p_1,\ldots,p_k,\alpha_{q_1},\ldots,\alpha_{q_l}).
	\]
	Hence,
	\[
	\alpha(p_1,\ldots,p_k, \alpha_{q_1},\ldots,\alpha_{q_l})\slash\theta_{\mathcal{E}}(X)=
	\beta(p_1,\ldots,p_k, \alpha_{q_1},\ldots,\alpha_{q_l})\slash\theta_{\mathcal{E}}(X).
	\]
	This completes the proof of the correctness of the definition of $f^{\ast}$. 
	
	It is easily seen that $f^\ast$ is a homomorphism  of
	$\ME[Y]$ onto  $\ME[X]$. We leave for the reader to show that $f^{\ast}$, in fact, is an isomorphism. (Exercise~\ref{section:malcev-matrices}.\ref{EX:tietze's-theorem})
\end{proof}

The following proposition is useful, if one wants to move from one set of equalities to another one that is expressed in terms of different term variables so that both sets would define the same E-matrix relative to a class $\E$.
\begin{prop}\label{P:mal'cev's-procedure}
	Let {\em$X\subseteq\EqL$} and $\Var^{\prime}\subseteq\VarL$ with $\Var^{\prime}\cap\Var(X)=\varnothing$.
	Also, let $X$ define an E-matrix {\em$\mat{M}=\langle\alg{A},\varDelta\rangle$} with respect to a valuation $v$.
	Assume that there is a valuation $w$ in {\em$\alg{A}$} such that the set $\set{w(q)}{q\in\Var^{\prime}}$ generates {\em$\alg{A}$}, that is, there is a map {\em$f:\Var(X)\longrightarrow
		\Forms_{\mathcal{L}_{\mathcal{V}^{\prime}}}$} with
	$\Var^{\prime}=\bigcup\set{\Var(f(p))}{p\in\Var(X)}$ such that for any $p\in\Var(X)$,
	$v[p]=w[f(p)]$. Now let a set $Y$ be obtained from $X$ by the simultaneous replacement of every occurrence of $p\in\Var(X)$ in each equality of $X$ by the term $f(p)$. Then $Y$ defines {\em$\mat{M}$} with respect to $w$.
\end{prop}
\begin{proof}
	We define:
	\[
	Y^{\prime}:=Y\cup\set{p\approx f(p)}{p\in\Var(X)}.
	\]
	
	In virtue of Tietze's theorem (Proposition~\ref{P:tietze's-theorem}), the matrices $\ME[Y]$ and $\ME[Y^{\prime}]$ are isomorphic.
	
	Now we define a valuation in $\alg{A}$:
	\[
	v^{\ast}:=\begin{cases}
		\begin{array}{cl}
			v[q] &\text{if $q\in\Var(X)$}\\
			w[q] &\text{if $q\notin\Var(X)$}.
		\end{array}
	\end{cases}
	\]
	
	We note that for any equality $\epsilon$ with $\Var(\epsilon)\subseteq\Var(X)$,
	\[
	\mat{M}\models_{v^{\ast}}\epsilon\Longleftrightarrow\mat{M}\models_{v}\epsilon.
	\]
	This implies that $X$ defines $\mat{M}$ with respect to $v^{\ast}$. In particular, 
	for any equality $\epsilon$ with $\Var(\epsilon)\subseteq\Var(X)$,
	\[
	\mat{M}\models_{v^{\ast}}\epsilon\Longrightarrow X\vdashE\epsilon. \tag{\ref{P:mal'cev's-procedure}--$\ast$}
	\]
	
	Next we define:
	\[
	X^{\prime}:=X\cup\set{p\approx f(p)}{p\in\Var(X)}.
	\]
	
	We notice that $\Var(X^{\prime})=\Var(Y^{\prime})=\Var(X)\cup\Var^{\prime}$.
	Also, it should be clear that, using Proposition~\ref{P:variables-auxiliary}, for any $\epsilon$ with $\Var(\epsilon)\subseteq\Var(X)\cup\Var^{\prime}$,
	\[
	Y^{\prime}\vdashE\epsilon\Longleftrightarrow X^{\prime}\vdashE\epsilon.\tag{\ref{P:mal'cev's-procedure}--$\ast\ast$}
	\]
	(See Exercise~\ref{section:malcev-matrices}.\ref{EX:mal'cev's-procedure-1}.) In virtue of Dyck's theorem (Proposition~\ref{P:dyck-theorem}), the matrices
	$\ME[Y^{\prime}]$ and $\ME[X^{\prime}]$ are isomorphic and, hence, the matrices
	$\ME[Y]$ and $\ME[X^{\prime}]$ are also isomorphic. It remains to show that the matrices $\ME[X^{\prime}]$ and $\ME[X]$ are isomorphic.
	
	First, we note that $\Var(X)\subseteq\Var(X^{\prime})$. Second, we observe that for any term $\epsilon$ with $\Var(\epsilon)\subseteq\Var(X)$, $X\vdashE\epsilon$ implies $X^{\prime}\vdashE\epsilon$. According to the proof of Dyck's theorem
	(Proposition~\ref{P:dyck-theorem}), this means that the map
	\[
	g:\alpha\slash\theta_{\mathcal{E}}(X)\mapsto\alpha\slash\theta_{\mathcal{E}}(X^{\prime}),~\text{for any term $\alpha$ with $\Var(\alpha)\subseteq\Var(X)$},
	\] 
	is an epimorphism of $\ME[X]$-algebra onto $\ME[X^{\prime}]$-algebra.
	We will show that $g$ is, in fact, an isomorphism. 
	
	Indeed, assume that for arbitrary terms $\alpha$ and $\beta$ with $\Var(\alpha)\cup\Var(\beta)\subseteq\Var(X)$, $\alpha\slash\theta_{\mathcal{E}}(X^{\prime})=
	\beta\slash\theta_{\mathcal{E}}(X^{\prime})$, that is $X^{\prime}\vdashE\alpha\approx\beta$.
	We observe that $\mat{M}\models_{v^{\ast}}X^{\prime}$ and hence $\mat{M}\models_{v^{\ast}}\alpha\approx\beta$. Since $X$ defines $\mat{M}$ with respect to $v^{\ast}$, $X\vdashE\alpha\approx\beta$, that is $\alpha\slash\theta_{\mathcal{E}}(X)=
	\beta\slash\theta_{\mathcal{E}}(X)$. Thus, we conclude that $\ME[x]$ and $\ME[X^{\prime}]$ are isomorphic. Hence $\ME[X]$, $\ME[Y]$ and $\mat{M}$ are also isomorphic; 
	that is $Y$ is a defining set for $\mat{M}$.

	What remains to finish the proof is to show that $Y$ defines $\mat{M}$ with respect to $w$.
	We leave for the reader to complete the proof. (See Exercise~\ref{section:malcev-matrices}.\ref{EX:mal'cev's-procedure-2}.)
\end{proof}

\begin{rem}\label{R:mal'cev's-procedure}
	{\em 	We note that Proposition~\ref{P:mal'cev's-procedure} is true, even when $\bigcup\set{\Var(f(p))}{p\in\Var(X)}\subset\Var^{\prime}$.  In this case, we, having selected a term variable $p_0\in\Var(X)$, add to $Y$ the set $\set{p_0\approx q}{q\in\Var^{\prime}\setminus\bigcup\set{\Var(f(p))}{p\in\Var(X)}}$. (Exercise~\ref{section:malcev-matrices}.\ref{EX:mal'cev's-procedure-3})}
\end{rem}

\begin{prop}[Mal'cev's second theorem]\label{P:Malcev-second}
	Let $\E$ be a class of E-matrices which contains a trivial matrix. Then if $\Is\A[\E]\subseteq\A[\E]$, $\Su\A[\E]\subseteq\A[\E]$, and $\Pro\A[\E]\subseteq\A[\E]$, in particular if $\A[\E]$ is a quasi-variety, then each {\em$\ME[X]$} belongs to $\A[\E]$.
\end{prop}
\begin{proof}
	We recall that $\ME[X]$-algebra is a subalgebra of $\LTE[X]$-algebra. Then, we apply Proposition~\ref{P:Birkhoff's-generalized}.
\end{proof}

\paragraph{Exercises~\ref{section:malcev-matrices}}
\begin{enumerate}
	\item \label{EX:equivalence-with-Mal[X]}Prove Proposition~\ref{P:equivalence-with-Mal[X]}.
	\item \label{EX:defining-set-3}Give details of the concluding step of the proof of Proposition~\ref{P:defining-set-3}.
	\item \label{EX:dyck-theorem-1} Prove that $f^{\ast}$ in the proof of Proposition~\ref{P:dyck-theorem} is a homomorphism.
	\item \label{EX:dyck-theorem-2} Prove that, given $X\cup Y\subseteq\EqL$, if  $\Var(X)=\Var(Y)$ and for any $\epsilon$ with $\Var(\epsilon)\subseteq\Var(X)$, the equivalence
	\[
	X\vdashE\epsilon\Longleftrightarrow Y\vdashE\epsilon
	\]
	holds, then $f:p\slash\theta_{\mathcal{E}}(X)\mapsto p\slash\theta_{\mathcal{E}}(Y)$ can be extended to an isomorphism of $\ME[X]$-algebra onto $\ME[Y]$-algebra.
	\item \label{EX:tietze's-theorem} Complete the proof of Proposition~\ref{P:tietze's-theorem}, by showing that the homomorphism $f^\ast$ is an isomorphism.
	\item \label{EX:mal'cev's-procedure-1}Prove a step in the proof of Proposition~\ref{P:mal'cev's-procedure}: for any $\epsilon\in X$, $Y\vdash_{\text{E}2}\epsilon$.	
	\item \label{EX:mal'cev's-procedure-2}Complete the proof of Proposition~\ref{P:mal'cev's-procedure}, by showing that the set $Y$, as defined in the proof, defines $\mat{M}$ with respect to the valuation $w$. \textit{Hint}: Proof can go as follows.
	\begin{enumerate}
		\item Similarly to (\ref{P:mal'cev's-procedure}--$\ast$), prove that for any equality $\epsilon$ with $\Var(\epsilon)\subseteq\Var(Y)$,
		\[
		\mat{M}\models_{v^{\ast}}\epsilon\Longleftrightarrow\mat{M}\models_{w}\epsilon.
		\]
		\item Show that $X^{\prime}$ defines $\mat{M}$ with respect to $v^{\ast}$. (This is a generalization of the statement in the proof that $X$ defines $\mat{M}$ with respect to $v^{\ast}$.)
		\item Use (\ref{P:mal'cev's-procedure}--$\ast\ast$) to prove that $Y^{\prime}$ defines $\mat{M}$ with respect to $v^{\ast}$.
		\item Prove that for any $\epsilon$ with $\Var(\epsilon)\subseteq\Var(Y)$,
		\[
		Y^{\prime}\vdashE\epsilon\Longrightarrow Y\vdashE\epsilon.
		\]
		(This implication is analogous to the one regarding $X^{\prime}$ and $X$ in the proof.)
		\item Use the last implication to derive that $Y$ defines $\mat{M}$ with respect to $w$. 
	\end{enumerate}
	\item \label{EX:mal'cev's-procedure-3}Prove Proposition~\ref{P:mal'cev's-procedure} when $\bigcup\set{\Var(f(p))}{p\in\Var(X)}\subset\Var^{\prime}$. See Remark~\ref{R:mal'cev's-procedure}.
\end{enumerate}

\section{Equational consequence based on\\ implicational logic}\label{section:equational-based-on-implicational}
Let $\mathcal{S}$ be an implicational logic with respect to $\im(p,q)$ in a language $\Lan$. 

In this section, we discuss equational consequence based on implicational assertional consequence. For any equation $\alpha\approx\beta$ in $\Lan$, we define the following transforming function:
\[
\textbf{i}(\alpha\approx\beta):=\lbrace\im(\alpha,\beta),\im(\beta,\alpha)\rbrace;
\] 
and for any $X\subseteq\EqL$,
\[
X^{\textbf{i}}:=\bigcup\set{\textbf{i}(\epsilon)}{\epsilon\in X}.
\]

To discuss equational consequence based on abstract logics that are both implicative and implicational with respect to $p\rightarrow q$, we define a transforming function:
\[
\textbf{s}(\alpha\approx\beta):\alpha\leftrightarrow\beta,
\]
where
\[
\alpha\leftrightarrow\beta:=(\alpha\rightarrow\beta)\land(\beta\rightarrow\alpha).
\]
The function \textbf{s} is extended to an arbitrary set $X\subseteq\EqL$ as follows:
\[
X^{\textbf{s}}:=\set{\textbf{s}(\epsilon)}{\epsilon\in X}.
\]

We aim to prove the following proposition.
\begin{prop}\label{P:equational-on-implicational-1}
	Let $\mathcal{S}$ be an implicational logic with respect to $\im(p,q)$ in $\Lan$. Also, let $\E$ be an equational logic with $\mathcal{A}[\E]$ such that {\em$\set{\langle\alg{A},\one\rangle}{\alg{A}\in\mathcal{A}[\E]}=\QS$}. Then for any {\em$X\cup\lbrace\epsilon\rbrace\subseteq\EqL$},
	{\em\[
		X\vdashE\epsilon~\Longleftrightarrow~\textbf{i}(\epsilon)\subseteq\ConS{X^{\textbf{i}}}.
		\]}
	Further, if $\mathcal{S}$ is implicative and implicational with respect to $p\rightarrow q$, and, in addition, if $\mathcal{S}$ is sound with respect to the inferences rules {\em(\text{a}--$i$)} and {\em(\text{b}--$i$)} $($Section~\ref{section:inference-rules} $)$, then for
	any {\em$X\cup\lbrace\epsilon\rbrace\subseteq\EqL$},
	{\em\begin{equation}\label{E:equational-on-implicational-1}
			X\vdashE\epsilon~\Longleftrightarrow~X^{\textbf{s}}\vdash_{\mathcal{S}}\textbf{s}(\epsilon).
	\end{equation}}
\end{prop}
\begin{proof}
	Let $\epsilon:=\alpha\approx\beta$. Also, we denote by $\alg{A}^{\ast}:=\langle\alg{A},\one\rangle$, an algebraic expansions of algebra $\alg{A}$ of type $\Lan$. We use this notation to refer to algebras of $\PS$.\footnote{We remind the reader that, if the constant $\one$ is already present in the signature of $\alg{A}$, we count $\langle\alg{A},\one\rangle$ the same as $\alg{A}$; see Section~\ref{section:chapter-unital-preliminaries}.}
	
	By premise, the logic $\mathcal{S}$ is assertional. We note that, according to Proposition~\ref{P:LT-unital-PS-QS-equivalence}, $\mathcal{S}$ is determined by $\QS$. 
	
	Assume that $X\vdashE\epsilon$. Let $v$ be a valuation in $\alg{A}$ such that $\alg{A}^{\ast}\in\QS$. Suppose that $v[X^{\textbf{i}}]\subseteq\lbrace\one\rbrace$. This implies that for each equality $\alpha^{\prime}\approx\beta^{\prime}\in X$, $\lbrace v[\im(\alpha^{\prime},\beta^{\prime})], v[\im(\beta^{\prime},\alpha^{\prime})]
	\rbrace\subseteq\lbrace\one\rbrace$. Taking into account~\eqref{E:LT-algebras-unital-3}, by Definition~\ref{D:implicational}, $v[\alpha^{\prime}]=v[\beta^{\prime}]$. This means that $v[X]\subseteq\varDelta$. By premise, we obtain that $v[\alpha]=v[\beta]$. Using Definition~\ref{D:implicational} again, we derive that $v[\textbf{i}(\epsilon)]\subseteq\lbrace\one\rbrace$.
	
	Conversely, suppose that $\textbf{i}(\epsilon)\subseteq\ConS{X^{\textbf{i}}}$
	that $v[X]\subseteq\lbrace\one\rbrace$, for any valuation $v$ in $\alg{A}$ with $\alg{A}^{\ast}\in\QS$. According to Definition~\ref{D:implicational},
	$v[X^{\textbf{i}}]\subseteq\lbrace\one\rbrace$. Since $\mathcal{S}$ is assertional, $v[\im(\alpha,\beta)], v[\im(\beta^{\prime},\alpha^{\prime})]
	\rbrace\subseteq\lbrace\one\rbrace$. Applying Definition~\ref{D:implicational} one more time, we get that $v[\alpha]=v[\beta]$. This completes the proof of the first part.
	
	We leave for the reader to prove the second claim of the proposition. (Exercise~\ref{section:equational-based-on-implicational}.\ref{EX:equational-on-implicational})
\end{proof}

The next proposition is a refinement of the last one.
\begin{prop}\label{P:equational-on-implicational-2}
	Let $\mathcal{S}$ be a finitary logic whose language contains binary connectives $\land$ and $\rightarrow$ such that $\mathcal{S}$ is implicative and implicational with respect to $p\rightarrow q$ and sound with respect to the hyperrule {\em(\text{c}--$i$)} and the inferences rules {\em(\text{a}--$i$)} and {\em(\text{b}--$i$)} $($Section~\ref{section:inference-rules}$)$.  Also, we assume that the restriction of any algebraic $\mathcal{S}$-model to the signature $\langle\land,\rightarrow,\one\rangle$ is a Brouwerian semilattice, or at least satisfies the properties {\em$(\text{d}_4)$},~\eqref{E:diego-algebra-1} and~\eqref{E:brouwerian-conjunction} $($Section~\ref{section:prelimenaries-algebra}$)$. Further, let $\E$ be an equational logic determined by a class $\mathcal{A}[\E]$ such that {\em$\set{\langle\alg{A},\one\rangle}{\alg{A}\in\mathcal{A}[\E]}=\QS$} and $\E^{\prime}$ be an equational logic determined by a class $\mathcal{A}[\E^{\prime}]$ with {\em$\QS\subseteq\set{\langle\alg{A},\one\rangle}{\alg{A}\in\mathcal{A}[\E^{\prime}]}\subseteq\TS$}. Then for any set {\em$X\cup\lbrace\epsilon\rbrace\subseteq\EqL$}, the following conditions are equivalent$\,:$
	{\em\[
		\begin{array}{cl}
			(\text{a}) & X\vdash_{\mathcal{E}^{\prime}}\epsilon;\\
			(\text{b}) & X\vdash_{\mathcal{E}}\epsilon;\\
			(\text{c}) &X^{\textbf{s}}\vdash_{\mathcal{S}}\textbf{s}(\epsilon).
		\end{array}
		\]}
\end{prop}
\begin{proof}
	The implication (a)$\Rightarrow$(b) is supported by the inclusion~\eqref{E:LT-algebras-unital-3}. The implication (b)$\Rightarrow$(c)
	is the $\Rightarrow$-part of the equivalence~\eqref{E:equational-on-implicational-1}. Thus, it remains to prove the implication (c)$\Rightarrow$(a).
	
	Assume that (c) holds. Since $\mathcal{S}$ is finitary, for some set $X_0\Subset X$, $X_{0}^{\textbf{s}}\vdash_{\mathcal{S}}\textbf{s}(\epsilon)$.
	For the logic $\mathcal{S}$ is sound with respect to the rules (\text{a}--$i$) and (\text{b}--$i$) and the hyperrule (\text{c}--$i$), we derive that
	$\vdash_{\mathcal{E}}\bigwedge X_{0}^{\textbf{s}}\rightarrow\textbf{s}(\epsilon)$, that is, $X_{0}^{\textbf{s}}\rightarrow\textbf{s}(\epsilon)\in\ThmS$. Since the logic $\mathcal{S}$ is unital, in virtue of the second part of Proposition~\ref{P:valuations-in-LT}, at the valuation $\bm{v}_{\bm{T}_{\mathcal{S}}}$ in the algebra $\LT$, $\bm{v}_{\bm{T}_{\mathcal{S}}}[\bigwedge X_{0}^{\textbf{s}}\rightarrow\textbf{s}(\epsilon)]=\one_{\bm{T}_{\mathcal{S}}}$, that is, 
	\[
	\bm{v}_{\bm{T}_{\mathcal{S}}}[\bigwedge X_{0}^{\textbf{s}}]\rightarrow\bm{v}_{\bm{T}_{\mathcal{S}}}[\textbf{s}(\epsilon)]=\one_{\bm{T}_{\mathcal{S}}}. \tag{\ref{P:equational-on-implicational-2}--$\ast$}
	\] 
	
	Now let $\alg{A}\in\mathcal{A}[\E^{\prime}]$ and $v$ be a valuation in $\alg{A}$ such that $v[X]\subseteq\varDelta$. The latter implies that
	$v[X_0]\subseteq\varDelta$. Using successively  ($\text{d}_4$) and~\eqref{E:brouwerian-conjunction}, this implies that $v[\bigwedge X_{0}^{\textbf{s}}]=\one$.
	
	Next, we define a map $f:\bm{v}_{\bm{T}_{\mathcal{S}}}[\alpha]\mapsto v[\alpha]$, for any $\Lan$-term $\alpha$. Since $\LT$ is a free algebra of rank $\VarL$ over $\TS$ (Proposition~\ref{P:LT=free-algebra}), the map $f$ is a homomorphism of $\LT$ into $\alg{A}$. Using (\ref{P:equational-on-implicational-2}--$\ast$), we derive that
	\[
	v[\bigwedge X_{0}^{\textbf{s}}]\rightarrow v[\textbf{s}(\epsilon)]=\one.
	\]
	Then, with the help of~\eqref{E:diego-algebra-1}, we conclude that $v[\textbf{s}(\epsilon)]=\one$, that is, in virtue of~\eqref{E:brouwerian-conjunction} and ($\text{d}_4$), $v[\epsilon]\in\varDelta$.
\end{proof} 

\begin{rem}
	{\em Besides applications of Proposition~\ref{P:equational-on-implicational-2} which are demonstrated below, the usefulness of this proposition is seen in the following. In order to prove $X^{\textbf{s}}\vdash_{\mathcal{S}}\textbf{s}(\epsilon)$ we have to examine the class $\QS$, but for refuting the deducibility assertion, namely for $X^{\textbf{s}}\not\vdash_{\mathcal{S}}\textbf{s}(\epsilon)$, we have at our disposal a larger class of algebraic separating tools.}
\end{rem}

Below we consider two applications of Proposition~\ref{P:equational-on-implicational-2}.

\subsection{Abstract logic $\EB$}\label{section:Boolean-equational}
Let $\EB$ be the class of all E-matrices whose algebras are Boolean algebras.
Since the class $\mathcal{A}[\EB]$ is a variety and, therefore, is closed under ultraproducts, in virtue of Proposition~\ref{P:finitariness-E-consequence}, the abstract logic $\EB$ is finitary. A more complete characterization of it is obtained below.
For this characterization, we employ the matrix consequence 
$\models_{\textbf{B}_{2}}$, that is the matrix consequence with respect to the logical matrix $\booleTwo$ (Section~\ref{S:two-valued}), and the classical propositional logic
\textsf{Cl} (Section~\ref{section:inference-rules}).
\begin{prop}\label{P:logic_EB}
	For any set {\em$X\cup\lbrace\epsilon\rbrace\subseteq\EqL$}, the following conditions are equivalent:
	{\em\[
		\begin{array}{cl}
			(\text{a}) & X\vdashEB\epsilon;\\
			(\text{b}) & X^{\textbf{s}}\models_{\textbf{B}_{2}}\textbf{s}(\epsilon);\\
			(\text{c}) & X^{\textbf{s}}\vdash_{\textsf{Cl}}\textbf{s}(\epsilon).\\
		\end{array}
		\]}
\end{prop}
\begin{proof}
	The implication (a)$\Rightarrow$(b) follows from~\eqref{E:equational-on-implicational-1} and the fact that $\booleTwo$-algebra is a simplest nontrivial Boolean algebra.
	
	To prove the implication (b)$\Rightarrow$(c)  we recall that, in virtue of Proposition~\ref{P:con-finite matrix-is-finitary}, the consequence $\models_{\textbf{B}_{2}}$ is finitary. Then, we apply the equality~\eqref{E:B_2=Cl}.
	
	To prove the implication (c)$\Rightarrow$(a) we apply Proposition~\ref{P:equational-on-implicational-2}.
\end{proof}
\begin{cor}\label{C:logic_EB}
	{\em$\set{\alpha\approx\alpha}{\alpha\in\FormsL}\subset\ConEB(\varnothing)$}.
\end{cor}
\begin{proof}
	Indeed, we notice that $\varnothing\vdash_{\textsf{Cl}}p\land(q\lor r)\leftrightarrow (p\land q)\lor(p\land r)$. Then, we apply Proposition~\ref{P:logic_EB}.
\end{proof}

\subsection{Abstract logic $\EH$}\label{section:Heyting-equational}
Let $\EH$ be the class of all E-matrices whose algebras are Heyting algebras.
Since the class of all Heyting algebras is a variety, the abstract logic $\EH$ is finitary. (See Exercise~\ref{section:equational-based-on-implicational}.\ref{EX:EH-is-finitary}.) We learn more about it in the following proposition.

\begin{prop}\label{P:logic_EH}
	For any set {\em$X\cup\lbrace\epsilon\rbrace\subseteq\EqL$}, the following conditions are equivalent:
	{\em\[
		\begin{array}{cl}
			(\text{a}) & X\vdashEH\epsilon;\\
			(\text{b}) & X^{\textbf{s}}\vdash_{\textsf{Int}}\textbf{s}(\epsilon).\\
		\end{array}
		\]}
\end{prop}
\begin{proof}
	Using Proposition~\ref{P:T_Int=Heyting-algebras} and Proposition~\ref{P:equational-on-implicational-2}, we obtain the equivalence in question.	
\end{proof}

Similarly to Corollary~\ref{C:logic_EB}, we obtain:
\begin{cor}
	{\em$\set{\alpha\approx\alpha}{\alpha\in\FormsL}\subset\ConEH(\varnothing)$}.
\end{cor}

\paragraph{Exercises~\ref{section:equational-based-on-implicational}}
\begin{enumerate}
	\item \label{EX:equational-on-implicational}Prove that if $\mathcal{S}$ is implicative and implicational with respect to $p\rightarrow q$, and, in addition,  $\mathcal{S}$ is sound with respect to the inferences rules (\text{a}--$i$) and (\text{b}--$i$) (Section~\ref{section:inference-rules}), then for
	any $X\cup\lbrace\epsilon\rbrace\subseteq\EqL$,
	\[
	X\vdashE\epsilon~\Longleftrightarrow~X^{\textbf{s}}\vdash_{\mathcal{S}}\textbf{s}(\epsilon).
	\]
	\item \label{EX:EH-is-finitary} Show that the equational logic $\EH$ is finitary.
\end{enumerate}

\section{Philosophical and historical notes}
The question about the nature of equality had been challenging in philosophy of logic, at least until the beginning of the twentieth century. To illustrate the core of the problem, we quote A. N. Whitehead.
\begin{quote}
	``Identity may be conceived as a special limiting case of equivalence. For instance in arithmetic we write, $2+3=3+2$. This means that, in so far as the total number of objects mentioned, $2+3$ and $3+2$ come to the same number, namely $5$. But $2+3$ and $3+2$ are not identical: the order of the symbols is different in the two combinations, and this difference of order directs different processes of thought. The importance of the equation arises from its assertion that these different processes of thought are identical as far as the total number of things thought of is concerned.'' \cite{whitehead2009}, section I.3
\end{quote}

In the above excerpt, it is clear that Whitehead, on the one hand, points at the difference between the terms `$2+3$' and `$3+2$' and on the other hand, indicates the sameness of what the two thought processes associated with these terms, 
lead, namely, the number $5$. Thus, Whitehead makes it clear, though there may be confusion in the way he uses the terminology, that there is the equality between the terms `$2+3$' and `$3+2$' (he uses the term \emph{equivalence}), at the same time pointing  at their non-identity.

The distinction between the concepts of \emph{equality} and \emph{identity}, as we will see below, was fully realized (e.g., by Frege) only by the end of the nineteenth century.
The inability to see this distinction was, as Quine put it, a ``confusion between the sign and the object.'' No wonder that Wittgenstein mockingly, in our opinion, stated:
\begin{quote}
	``[$\dots$] to say of \textit{two} things that they are identical is nonsense, and to say of \textit{one} thing that it is identical with itself is to say nothing at all.'' \cite{wittgenstein2001}, 5.5303
\end{quote}

For Aristotle, the distinction indicated in the last paragraph has not yet been recognized. He determines rather than defines things that should be considered one and the same, if they are indistinguishable by all the attributes that can be applied to them:
\begin{quote}
	``Only to things that are indistinguishable and one in being it generally agreed that all the same attributes belong.'' \textit{Topics}, vii.I ($152^{\text{a}}30$); quoted from~\cite{kneales1962}, section II.3
\end{quote}

As regards logic, in full accordance with Aristotle’s view, signs, and not things, are a working tool of argumentation. Then, applying Aristotle's observation, not to things, but to their names and replacing identity with equality, we could express Aristotle's observation in symbolic form as follows:
\begin{equation}\label{E:Leibniz-principle}
	\forall  F(F(\alpha)\Leftrightarrow F(\beta))\Rightarrow\alpha\approx\beta).
\end{equation}

In philosophical literature ~\eqref{E:Leibniz-principle} is known as \emph{Leibniz's principle of identity of indiscernibles}; cf., e.g.,~\cite{ishiguro1990}, section II.1. And again in accordance with Aristotle's view, the quantification in~\eqref{E:Leibniz-principle} is supposed to run over all applicable attributes, or, better to say (in accordance with a new linguistic tune), over all applicable contexts.

However, there is a problem with attributing~\eqref{E:Leibniz-principle} to Leibniz, since, as H. Ishiguro points out, ``many discussions on extensionality and intensionality or referential transparency and referential opacity are made by using the \emph{salva veritate} principle in the sense of `indiscernibility of identicals', and calling it Leibniz's law. But this indiscernibility principle is not one that Leibniz himself expressed;'' cf.~\cite{ishiguro1990}, section II.1. Leibniz formulated his ``salva veritate'' principle as follows:
\begin{quote}
	``Those terms of which one can be substituted for the other without affecting truth are identical. (\emph{Eadem sunt, quorum unum alteri substitui potest salva veritate})'' Quoted from~\cite{ishiguro1990}, section II.1; see also~\cite{kneales1962}, section V.3
\end{quote}

Quoting Leibniz's \emph{salva veritate} principle, Ishiguro adds the comment:
\begin{quote}
	``The \emph{salva veritate} principle is \ldots an attempt to determine the identity of concepts expressed by words, not by reference to mental images invoked, or any consideration that we can apply the words taken in isolation, but by the role they play in determining the truth of the proposition in which they occur. Intersubstitutivity discovers the identity of the role.'' \cite{ishiguro1990}, section II.1
\end{quote}

And a bit further, she continues:
\begin{quote}
	``The Leibnizian \emph{salva veritate} principle which occurs mainly in his logical works and in his treatises on conceptual analysis, is not, I suggest, a general definition of identity. It is simply a principle to determine the identity of concepts.'' \cite{ishiguro1990}, section II.1
\end{quote}

Note that, although Aristotle and Leibniz do not make a clear distinction between the \emph{equal} and the \emph{identical}, they employ the \emph{substitution} operation, that is our \emph{replacement} in the sense of Section~\ref{section:languages}, which we use in the formulations of the deductive systems E1--E3.

We believe that G. Frege was the first who not only made it clear that there is a distinction between the equal and the identical, but also pointed at two features of the former: 1) equality is a relation between terms of a formal language, and 2) the determination of equal terms is of semantic character.

He starts his famous paper \emph{\"{U}ber Sinn und Bedeutung} with the following discussion.
\begin{quote}
	``Equality$^{\ast}$ [Frege's footnote: $^{\ast}$I use this word in the sense of identity and understand `$a=b$' to have the sense of `$a$ is the same as $b$' or `$a$ and $b$ coincide.'] gives rise to challenging questions which are not altogether easy to answer. Is it a relation? A relation between objects, or between names or signs of objects? In my \emph{Begriffsschrift} I assume the latter. The reasons which seem to favour this are the following: $a=a$ and $a=b$ are obviously statements of differing cognitive value; $a=a$ holds \emph{a priori} and, according to Kant, is to be labeled analytic, while statements of the form $a=b$ often contain very valuable extensions of our knowledge and cannot always be established \emph{a priori}.'' Quoted from~\cite{heijenoort2002}, p. 56
\end{quote}

Frege further adds:
\begin{quote}
	``Now if we were to regard equality as a relation between that which the names `$a$' and `$b$' designate, it would deem that $a=b$ could not differ from $a=a$ (i.e. provided $a=b$ is true). A relation would thereby be expressed of a thing itself, and indeed one in which each thing stands to itself but to no other thing. What is intended to be said by $a=b$ seems to be that signs or names `$a$' and `$b$' designate the same thing, so that those signs themselves would be under discussion; a relation between them would be asserted. But this relation would hold between the names or signs only in so far as they named or designated something.'' Quoted from \textit{ibid.}, p. 56
\end{quote}

Following Frege's approach,  we use the sign `$\approx$' to denote a relation between equal terms and the sign `$=$' for denoting the fact of identity. The necessity for the employment of the two different signs comes from a ``delay'' when, in order to assert an equality $\alpha\approx\beta$, we have to leave the linguistic framework of a formal language and enter the semantic framework of an algebra with a valuation $v$ in it and, then, check whether $v(\alpha)$ and $v(\beta)$ are identical.
The delay occurs when the two evaluations, $\alpha\mapsto v[\alpha]$ and
$\beta\mapsto v[\beta]$, (``different processes of thought,'' as Whitehead calls them) are running to arrive at one element of algebra, which fact is stated by writing `$v[\alpha]=v[\beta]$'.

On the other hand, if we stay in the framework of a formal language and define equality axiomatically, we have to use inference rules that define $\approx$ as an equivalence relation (as was noted by Whitehead; see above),
which are exemplified by the axiom and rules E1--$a$ and E1--$b$, while the rules  E1--$c$, E2--$b$ and E3--$b$, which together represent the converse of~\eqref{E:Leibniz-principle} in a weakened form, exemplify the principle of \emph{substitutivity}, or that of \emph{indiscernibility of identicals}. This is how W. V. Quine explains this principle.
\begin{quote}
	``It provides that, \emph{given a true statement of identity, one of its two terms may be substituted for the other in any true statement and the result will be true}.'' \cite{quine1994}, VIII Reference and Modality~\footnote{Although Quine points out that this principle may fail if ``the occurrence to be supplanted is not \emph{purely referential},'' (see \emph{ibid.}) this remark is not related to our usage of the aforementioned rules.}
\end{quote}

The definition of the relation $\vdashE$ is a natural implementation of Ferege's view on equality in the framework of matrix consequence. We also encounter in~\cite{malcev1973}, section 11.2,  a similar definition of what it means that a first-order formula is a ``consequence in the class $\mathfrak{R}$ of the collection of formulas $\mathfrak{S}$.''
This similarity prompted us to adapt Mal'cev's approach, having implemented it within the framework of equational consequence. The equational consequence in a language without logical connectives and quantifiers has been studied in~\cite{blok-pigozzi1989}, chapter 2, however, for a different purpose.\\

The conception a defining set of equalities and generating symbols was introduced by W. von Dyck in~\cite{dyck1882}, introduction and part 1, {\S} 4 and {\S} 5, where $\mathcal{A}[\mathcal{E}]$ was the class of all groups. This idea was picked up by other authors, including H. Tietze~\cite{tietze1908}, part 4, {\S} 11. A. I. Mal'cev extended this approach to abstract algebras in~\cite{malcev1956} and to models of first-order language in~\cite{malcev1958} and especially in~\cite{malcev1973}. (The latter book is a posthumous publication of Mal'cev's unfinished manuscript.)

The conception of a defining set~(of equalities) for an E-matrix (Definition~\ref{D:X-is-defining-for-M}) is our adaptation of an analogous definition for first-order logic in~\cite{malcev1973}, section 11.2; see also~\cite{henkin-monk-tarski1985}, remark 0.4.65. Furthermore, 
Proposition~\ref{P:Malcev's-first} and Proposition~\ref{P:Malcev-second} correspond to theorem 1 and theorem 4 of~\cite{malcev1973}, section 11.2, respectively; and Proposition~\ref{P:dyck-theorem} and Proposition~\ref{P:tietze's-theorem} correspond, respectively, to theorem 5 and theorem 6 there. The two last theorems originally appeared in~\cite{dyck1882}  and~\cite{tietze1908}, respectively, where they were proved for groups. \\

Taken by themselves, Proposition~\ref{P:logic_EB} and Proposition~\ref{P:logic_EH} are part of folklore.

\chapter[Equational L-Consequence]{Equational L-Consequence}
\label{chapter:L-equational-logic}
Substitution does not play a significant role in equational consequence, at least not more than it plays in ordinary matrix consequence. But we get to a different world when we let substitution play a key role. Because of this, the type of consequence which we will discuss in this chapter, as it can be seen after examination, falls into the category of Carnap's L-concepts; see~\cite{carnap1942}, {\S} 14. Roughly speaking, while in equational consequence the equality $\alpha\approx\beta$ is interpreted as an equation, in equational L-consequence it is interpreted as an identity. 

\section{Equational L-consequence}\label{section:equational-L-consequence}
\begin{defn}
	Given a non-empty class $\E$ of {\em E}-matrices, we define: For any set {\em$X\cup\lbrace\epsilon\rbrace\subseteq\EqL$},
	{\em\[
		X\vdashELog\epsilon\stackrel{\text{df}}{\Longleftrightarrow}\text{$\mat{M}\models X$ \textit{implies} $\mat{M}\models\epsilon$, \emph{for any} $\mat{M}\in\E$}.
		\]}
\end{defn}

It is clear that
\begin{equation}\label{E:equational-logic-implication-1}
X\vdashE\epsilon\Longrightarrow X\vdashELog\epsilon,
\end{equation}
but not vice versa. The first part we leave to the reader. (See Exercise~\ref{section:equational-L-consequence}.\ref{EX:equational-logic-implication-1}.)

To show that the converse of~\eqref{E:equational-logic-implication-1} in general does not hold, we take the class $\E$ consisting only of the E-matrix $\mat{M}=\langle\booleTwo,\varDelta\rangle$. Let us consider the two $\Lan_B$-equalities --- $p\approx\top$ and $p\approx\neg\top$. It is obvious that
$\mat{M}\nvDash p\approx\top$ and $\mat{M}\nvDash p\approx\neg\top$. Hence $p\approx\top\vdashELog p\approx\neg\top$. On the other hand, for any valuation $v$ in $\booleTwo$ with $v[p]=\one$, we obtain that $v$ satisfies $p\approx\top$ and refutes $p\approx\neg\top$; that is $p\approx\top\not\vdash_{\mathcal{E}} p\approx\neg\top$.

The equivalence in~\eqref{E:equational-logic-implication-1} however occurs when $X=\varnothing$, namely
\begin{equation}
\varnothing\vdashE\epsilon\Longleftrightarrow\varnothing\vdashELog\epsilon.
\end{equation}

\begin{prop}\label{P:equational-L-consequence}
	For any nonempty class $\E$ of E-matrices, the relation {\em$\vdashELog$} is a consequence relation.
\end{prop}
\begin{proof}
	It suffices to apply Proposition~\ref{P:semantic-consequence-2} for $\Lan$-equalities and the concept of model of type (B); see Section~\ref{section:realizations-abstract-logic}.
\end{proof}

Given a nonempty class $\E$ of E-matrices, we denote by $\ELog$ the abstract logic corresponding to the consequence relation $\vdashELog$ and call it the \textit{\textbf{equational {\em L}-consequence relative to}} $\E$, or \index{equational {\em L}-consequence} L-\textit{\textbf{equational logic}}\index{logic!equational} (\textit{\textbf{relative to}} $\E$) for short. In particular, we employ the class $\EE$ of all E-matrices of type $\Lan$, we have the abstract logic $\EELog$.\\

We note that structurality is a very rare property of equational L-consequence.
\begin{prop}
	Suppose the set $\VarL$ of term variables contains at least three variables.
	Then, given a class $\E$ of E-matrices, the consequence relation {\em$\vdashELog$} is structural if, and only if, $\E$ contains only a trivial E-matrix.
\end{prop}
\begin{proof}
	Let an E-matrix \mat{M} be trivial. Then $\mat{M}\models\epsilon$, for any equality $\epsilon$.
	
	Conversely, assume that \mat{M} is an arbitrary nontrivial E-matrix from $\E$ and $p,q$ and $r$ are three distinct term variables. It is obvious
	that $\mat{M}\not\models p\approx q$ and $\mat{M}\not\models p\approx r$. Hence
	$p\approx q\vdashELog p\approx r$. Now let $\sigma$ be the substitution that maps $q$ to $p$ and is identical on the other variables. Thus, we have: $\mat{M}\models\sigma(p)\approx\sigma(q)$
	and $\mat{M}\not\models\sigma(p)\approx\sigma(r)$. That is, $\sigma(p)\approx \sigma(q)\not\vdashELog \sigma(p)\approx\sigma(r)$.
\end{proof}

Even if only the trivial equational L-consequence is structural, every abstract logic $\ELog$ always has the following property:
\begin{equation}\label{E:implication-1}
X\vdashELog\epsilon \Longrightarrow X\vdashELog\sigma(\epsilon),
\end{equation}
for any $\Lan$-substitution $\sigma$.

Indeed, suppose $X\vdashELog\epsilon$ and $\sigma$ is an arbitrary $\Lan$-substitution. Let $\mat{M}=\langle\alg{A},\varDelta\rangle$ be any matrix of $\E$. Assume that $\mat{M}\vDash X$ and consider any valuation $v$ in $\alg{A}$. In virtue of Proposition~\ref{P:valuation-and-substitution}, there is a valuation $v_{\sigma}$ in $\alg{A}$ such that for any equality $\eta$,
\[
v[\sigma(\eta)]\in\varDelta\Longleftrightarrow v_{\sigma}[\eta]\in\varDelta.
\]
Since, by premise, $\mat{M}\models\epsilon$, then, in particular, $\mat{M}\models_{v_\sigma}\epsilon$, that is $v_{\sigma}[\epsilon]\in\varDelta$. The latter, through the equivalence above, implies that $v[\sigma(\epsilon)]\in\varDelta$.

Thus we have proved the following proposition.
\begin{prop}\label{P:E_L-theory-close-under-substotution}
	Every {\em$\ELog$}-theory is closed under substitution.
\end{prop}

As it is the case with equational consequence, there is a correlation between equational L-consequence and validity in model theory. 

We relate an instantiation of term language $\Lan$ expanded with a binary predicate symbol $\approx$ to form a first-order language $\FOL$, where the term variables play a role of individual variables and elementary formulas like $\alpha\approx\beta$ are interpreted by the equality of term functions $\beta^{\textsf{A}}$ and $\alpha^{\textsf{A}}$  in the corresponding first-order model $\mat{M}=\langle\textsf{A},\Func_{\mathcal{L}},\Cons_{\Lan},\varDelta\rangle$. For any elementary formula $\epsilon$, we denote by $\Forall\epsilon$ the universal closure of this formula; accordingly for a set $X$ of equalities, we denote:
\[
\Forall X:=\set{\Forall\epsilon}{\epsilon\in X}.
\]

The following equivalence must be obvious:
\begin{equation}\label{E:L-con-equivalence-1}
\begin{array}{rl}
\mat{M}\models\Forall\epsilon &\Longleftrightarrow(\text{$\mat{M}\models_{v}\epsilon$, for any valuation $v$ in $\alg{A}$})\\
&\Longleftrightarrow\mat{M}\models\epsilon,
\end{array}
\end{equation}
where the first occurrence of `$\models$' denotes the validity in the first-order model $\mat{M}$ (Section~\ref{section:preliminaries-model-theory}), and the last one the validity in the sense of Definition~\ref{D:E-matrix-semantics}.

Naturally, \eqref{E:L-con-equivalence-1} implies that
\[
\mat{M}\models\Forall X\Longleftrightarrow\mat{M}\models X,
\]
with the same distinction of the usage of $\models$ in the two last occurrences. 

Using~\eqref{E:L-con-equivalence-1}, we observe: For any nonempty class $\E$ of E-matrices,
\begin{equation}\label{E:L-con-equivalence-2}
X\vdashELog\epsilon\Longleftrightarrow (\mat{M}\models\Forall X\Rightarrow
\mat{M}\models\Forall\epsilon,~\text{for any $\mat{M}\in\E$}).
\end{equation}

Now we obtain the following.
\begin{prop}
	Let $\E$ be a nonempty class of E-matrices such that the class $\mathcal{A}[\E]$ is closed under ultraproducts. Then the abstract logic {\em$\ELog$} is finitary.
\end{prop}
\begin{proof}
	Indeed, we begin, using~\eqref{E:L-con-equivalence-2}:
	\[
	\begin{array}{rcl}
	X\vdashELog\epsilon &\Longleftrightarrow &\text{the set $\Forall X\cup\lbrace\Neg\Forall\epsilon\rbrace$ is not satisfiable in $\E$};\\
	&\Longleftrightarrow &\text{for some $X_0\Subset X$, the set
		$\Forall X_0\cup\lbrace\Neg\Forall\epsilon\rbrace$ is not satisfiable}\\
	&&\text{in $\E$ [by Proposition~\ref{P:compactness-thm-refined}]};\\
	&\Longleftrightarrow &\text{for some $X_0\Subset X$, $X_0\vdashELog\epsilon$}~\text{[again in virtue of~\eqref{E:L-con-equivalence-2}]}.
	\end{array}
	\]
\end{proof}

The next corollary directly follows from the last proposition.
\begin{cor}\label{C:EELog-is-finitary}
	The abstract logic {\em$\EELog$} is finitary.
\end{cor}

\paragraph{Exercises~\ref{section:equational-L-consequence}}
\begin{enumerate}
	\item \label{EX:equational-logic-implication-1}Prove the implication~\eqref{E:equational-logic-implication-1}.
	\item \label{EX:equational-logic-consequence}Complete the proof of Proposition~\ref{P:equational-L-consequence}.
\end{enumerate}

\section{$\EELog$ as a deductive system}\label{section:E_L-deductive-systems}
In this section we describe $\EELog$ deductively, that is, by means of axioms and inference rules. Our treatment of $\EELog$ will similar to that of $\EE$ in Section~\ref{section:E-deductive-systems}.

Indeed, first, we define systems E1s--E3s which obtained, respectively, from the systems E1--E3 the rule of (uniform) substitution as an additional inference rule. We use notation E$^{\ast}$\!s to refer to any of the systems E1s--E3s.
\[
\begin{array}{rc}
\text{Substitution:} &\dfrac{\bm{\alpha}\approx\bm{\beta}}{\sigma(\bm{\alpha})\approx\sigma(\bm{\beta})}, 
\end{array}
\]
where $\sigma$ is an $\Lan$-substitution, and $\alpha$ and $\beta$ are $\Lan$-terms.

Each relation $\vdash_{\text{E$^{\ast}$\!s}}$ is defined through the notion of formal derivation similar to the one defined in Section~\ref{section:E-deductive-systems} with the only difference --- the substitution rule is allowed in E$^{\ast}$\!s-derivations.

Next we prove a proposition analogous to Proposition~\ref{P:formal-derivation-equivalence}.

\begin{prop}\label{P:formal-derivation-L-equivalence}
	For any set {\em$X\cup\lbrace\epsilon\rbrace\subseteq\EqL$}, the following conditions are equivalent.
	{\em\[
		\begin{array}{rl}
		(\text{a}) &X\vdash_{\text{E1s}}\epsilon;\\
		(\text{b}) &X\vdash_{\text{E2s}}\epsilon;\\
		(\text{c}) &X\vdash_{\text{E3s}}\epsilon.
		\end{array}
		\]}
\end{prop}\index{$X\vdash_{\text{E1s}}\epsilon$}\index{$X\vdash_{\text{E2s}}\epsilon$}\index{$X\vdash_{\text{E3s}}\epsilon$}
\begin{proof}
	Argumentation is literally the same as in the proof of Proposition~\ref{P:formal-derivation-equivalence}, with one additional step. We demonstrate it, proving the case (a)$\Rightarrow$(b).  Suppose we have that
	$\vdash_{\text{E1s}}\epsilon$ and want to show that $X\vdash_{\text{E2s}}\epsilon$. Considering any E1s-derivation, we show, by induction on the length of this E1s-derivation, all moves along this derivation can be carried out in E2s. Sine the substitution rule is common in both systems, we have to be concerned about other rules of inference. But this was done in the case `(a)$\Rightarrow$(b)' of the proof of Proposition~\ref{P:formal-derivation-equivalence}. 
	
	The other cases are considered analogously.
\end{proof}

\begin{prop}\label{P:E_L-soundness}	
	Let $\E$ be any nonempty class of E-matrices.
	Suppose $\dfrac{\epsilon_1,\ldots,\epsilon_n}{\epsilon}$, where $n\ge 1$, is one of the inference rules of any of {\em E1s--E3s}. Then for any matrix {\em$\mat{M}\in\E$}, if all {\em$\mat{M}\models\epsilon_1,\ldots, \mat{M}\models\epsilon_n$} hold, then {\em$\mat{M}\models\epsilon$} also holds. Consequently, if all {\em$X\vdashELog\epsilon_{1},\ldots,X\vdashELog\epsilon_{n}$}
	hold, then {\em$X\vdashELog\epsilon$} also holds. Hence,
	{\em\[
		X\vdash_{\text{E}^{\ast}\!\text{s}}\epsilon\Longrightarrow X\vdashELog\epsilon.
		\]}
\end{prop}
\begin{proof}
	We prove the first part only for the substitution rule and leave  the rest of the proof to the reader. (Exercise~\ref{section:E_L-deductive-systems}.\ref{EX:E_L-soundness})
	
	Assume that the term variables of an equality $\alpha\approx\beta$ constitute a set $\lbrace p_1,\ldots, p_n\rbrace$. Let $\sigma$ be an $\Lan$-substitution according to which $\sigma[p_i]=\gamma_i$. Let $\mat{M}\in\E$, for which we assume that $\mat{M}\models\alpha\approx\beta$. Now let $v$ be any valuation in \mat{M}. Next, we define a valuation $v^{\prime}$ in \mat{M} such that
	$v^{\prime}[p_i]=v[\gamma_i]$.  It must be cleat that $v^{\prime}[\alpha]=v[\sigma(\alpha)]$ and $v^{\prime}[\beta]=v[\sigma(\beta)]$. By premise, we obtain that
	$v[\sigma(\beta)]=v[\sigma(\beta)]$.
\end{proof}

Given a set $X\cup\lbrace\alpha\approx\beta\rbrace\subseteq\EqL$, we define:
\begin{equation}\label{E:defn-L-congruence}
(\alpha,\beta)\in\theta_{\hat{\mathcal{E}}_{\text{L}}}(X)~\stackrel{\text{df}}{\Longleftrightarrow}~X\vdash_{\text{E}^{\ast}\!\text{s}}\alpha\approx\beta.
\end{equation}

\begin{prop}\label{P:theta_L(X)=congruence}
	Given a set {\em$X\subseteq\EqL$}, the relation {\em$\theta_{\hat{\mathcal{E}}_{\text{L}}}(X)$} is a congruence on  {\em$\FormsL$}.	
\end{prop}
\begin{proof}
	To prove this proposition we apply Proposition~\ref{P:formal-derivation-L-equivalence} and the definition of E1s.
\end{proof}

Our next aim is to prove soundness and completeness. Our proof will almost literally repeat the proof of Proposition~\ref{P:EE-consequence-completeness}, though with $\EELog$ instead of $\EE$ and with $\vdash_{\text{E}^{\ast}\!\text{s}}$ instead of $\vdash_{\text{E}^{\ast}}$.
\begin{prop}[soundness and completeness]\label{P:EE-L-consequence-completeness}
	For any set {\em$X\cup\lbrace\epsilon\rbrace\subseteq\EqL$},
	{\em\[
		X\vdash_{\text{E}^{\ast}\!\text{s}}\epsilon~\Longleftrightarrow~X\vdashEELog\epsilon.
		\]}
\end{prop}
\begin{proof}
	The $\Rightarrow$-implication (soundness) is the third part of Proposition~\ref{P:E_L-soundness}.
	
	Now we prove the $\Leftarrow$-implication (completeness). Let
	\[
	\epsilon:=\alpha\approx\beta
	\]
	and $X\not\vdash_{\text{E}^{\ast}\!\text{s}}\alpha\approx\beta$. Then $(\alpha,\beta)\notin\theta_{\hat{\mathcal{E}}_{\text{L}}}(X)$ and, in virtue of Proposition~\ref{P:theta_L(X)=congruence}, $\alpha\slash\theta_{\hat{\mathcal{E}}_{\text{L}}}(X)\neq\beta\slash\theta_{\hat{\mathcal{E}}_{\text{L}}}(X)$ in the algebra $\FormAl\slash\theta_{\hat{\mathcal{E}}_{\text{L}}}(X)$.
	
	For convenience, we denote
	\[
	\alg{A}_1:=\FormAl\slash\theta_{\hat{\mathcal{E}_{\text{L}}}}(X).
	\]
	
	Now, let $\gamma\approx\delta\in X$ and $\sigma$ be any $\Lan$-substitution. By definition, $X\vdash_{\text{E}^{\ast}\!\text{s}}\sigma(\gamma)\approx\sigma(\delta)$;
	therefore, in virtue of~\eqref{E:defn-L-congruence}, $(\sigma(\gamma),\sigma(\delta))\in\theta_{\hat{\mathcal{E}}_{\text{L}}}(X)$,
	that is, $\sigma(\gamma)\slash\theta_{\hat{\mathcal{E}}_{\text{L}}}=
	\sigma(\delta)\slash\theta_{\hat{\mathcal{E}}_{\text{L}}}$. 
	
	Since any valuation in $\alg{A}_1$ is a composition an $\Lan$-substitution and the canonical map $\eta\mapsto\eta\slash\theta_{\hat{\mathcal{E}}_{\text{L}}}$, this explains that
	$\langle\alg{A}_1,\Delta\rangle\models X$. 
	
	Let us consider a valuation
	\[
	v:\FormAl\longrightarrow\FormAl\slash\theta_{\hat{\mathcal{E}}_{\text{L}}}(X):~p\mapsto p\slash\theta_{\hat{\mathcal{E}}_{\text{L}}}(X),
	\]
	for any $p\in\VarL$.
	
	It must be clear that that $v[(\alpha\approx\beta)]=(\alpha\slash\theta_{\hat{\mathcal{E}}_{\text{L}}},\beta\slash\theta_{\hat{\mathcal{E}}_{\text{L}}})$. Therefore, $\langle\alg{A}_1,\Delta\rangle\not\models\alpha\approx\beta$.
\end{proof}
\begin{rem}
	{\em
		A congruence $\theta$ on an algebra \alg{A} is called \emph{fully invariant} (or \emph{completely invariant}) if for any substitution $\sigma$,
		\[
		(a,b)\in\theta~\Longrightarrow~(\sigma(a),\sigma(b)\in\theta.
		\]
		
		Thus, we see that the fact that $\theta_{\hat{\mathcal{E}}_{\text{L}}}$ is fully invariant was essential in the proof of the $\Leftarrow$-implication (completeness) of Proposition~\ref{P:EE-L-consequence-completeness}.}
\end{rem}

The following two corollaries follow immediately from Proposition~\ref{P:EE-L-consequence-completeness}.

\begin{cor}\label{C:congrunece-via-conE_L}
	Given a set {\em$X\subseteq\EqL$}, for any equality $\alpha\approx\beta$,
	{\em	\[
		(\alpha,\beta)\in\theta_{\hat{\mathcal{E}}_{\text{L}}}(X)\Longleftrightarrow \alpha\approx\beta\in\ConEEL(X).
		\]}
\end{cor}

We have earlier proved that the logic $\EELog$ is finitary; see Corollary~\ref{C:EELog-is-finitary}. Another way to obtain this result is by
Proposition~\ref{P:EE-L-consequence-completeness}.

\paragraph{Exercises~\ref{section:E_L-deductive-systems}}
\begin{enumerate}
	\item \label{EX:E_L-soundness}Complete the proof of Proposition~\ref{P:E_L-soundness}.
	\item Prove Corollary~\ref{C:congrunece-via-conE_L}.
\end{enumerate}

\section{Lindenbaum-Tarski matrices for logic $\ELog$}\label{section:LT-algebras-for-ELog}
Given a nonempty class $\E$ of E-matrices, we define:
\begin{equation}\label{E:definition-there_EL(X)}
(\alpha,\beta)\in\theta_{\mathcal{E}_{\text{L}}}(X)~\define~X\vdashEL\alpha\approx\beta;
\end{equation}
and we denote
\[
\theta_{\mathcal{E}_{\text{L}}}:=\theta_{\mathcal{E}_{\text{L}}}(\varnothing).
\]

It is clear that
\[
(\alpha,\beta)\in\theta_{\mathcal{E}_{\text{L}}}~\Longleftrightarrow~\text{for every $\mat{M}\in\E$, }\mat{M}\models\alpha\approx\beta.
\]

In virtue of Proposition~\ref{P:E_L-soundness}, each $\theta_{\mathcal{E}_{\text{L}}}(X)$ is a congruence on $\FormAl$. 

We note that
\begin{equation}\label{E:theta_EL(X)-inclusion-theta_EL(Y)}
X\subseteqq Y~\Longrightarrow~\theta_{\mathcal{E}_{\text{L}}}(X)\subseteq\theta_{\mathcal{E}_{\text{L}}}(Y);
\end{equation}
in particular,
\begin{equation}\label{E}
(\alpha,\beta)\in\theta_{\mathcal{E}_{\text{L}}}~\Longrightarrow~(\alpha,\beta)\in\theta_{\mathcal{E}_{\text{L}}}(X).
\end{equation}

\begin{defn}[Lindenbaum-Tarski matrices for $\E_{\text{L}}$]\label{D:LT-for-ELog}
	Let {\em$\E_{\text{L}}$} be an {\em L}-equational logic. For an arbitrary set {\em$X\subseteq\EqL$}, the E-matrix \index{equational logic!Lindenbaum-Tarski matrix} {\em$\LTEL[X]:=\langle\FormAl\slash\theta_{\mathcal{E}_{\text{L}}}(X),\varDelta\rangle$} is called a \textbf{Lindenbaum-Tarski matrix for {\em$\ELog$} relative to} $X$. We denote {\em$\LTEL:=\LTEL[\varnothing]$} and call the latter simply a \textbf{Lindenbaum-Tarski matrix for} {\em$\ELog$}.
\end{defn}

The following lemma is similar to Lemma~\ref{L:equivalence-with-LT[X]} and also is an adaptation of Proposition~\ref{P:valuation-in-F_L/theta}.
\begin{lem}\label{L:L-equivalence-with-LTL[X]}
	Let {\em$\ELog$} be an {\em\text{L}}-equational logic and {\em$X\subseteq\EqL$}. For any valuation $v$ in {\em$\LTEL[X]$}, there is an $\Lan$-substitution $\sigma$ such that for any term $\alpha$,
	{\em$v[\alpha]=\sigma(\alpha)\slash\theta_{\mathcal{E}_{\text{L}}}(X)$}. Conversely, for any $\Lan$-substitution $\sigma$, there is a unique valuation $v_\sigma$ in {\em$\FormAl\slash\theta_{\mathcal{E}_{\text{L}}}(X)$} such that for any term $\alpha$, {\em$v_\sigma[\alpha]=\sigma(\alpha)\slash\theta_{\mathcal{E}_{\text{L}}}(X)$}.
\end{lem}
\noindent\emph{Proof}~is left to the reader. (Exercise~\ref{section:LT-algebras-for-ELog}.\ref{EX:L-equivalence-with-LTL[X]})\\

Analogous to what we did in Section~\ref{section:Lindenbaum-Tarski-matrices}, we define a special valuation in $\FormAl\slash\theta_{\mathcal{E}_{\text{L}}}(X)$:
\[
\overline{v}_{X}[p]:=p\slash\theta_{\mathcal{E}_{\text{L}}}(X),~\text{for any $p\in\VarL$}.
\]
Extending $\overline{v}_X$ on all $\FormAl$, we get:
\[
\overline{v}_{X}[\alpha]=\alpha\slash\theta_{\mathcal{E}_{\text{L}}}(X),~\text{for any $\alpha\in\FormAl$}.
\]
\begin{prop}\label{P:equivalence-with-LTEL[X]}
	Let {\em$\ELog$} be an {\em L}-equational logic. For any set $X\cup\lbrace\epsilon\rbrace\subseteq\EqL$, the following conditions are equivalent.
	{\em\[
		\begin{array}{cl}
		(\text{a}) &X\vdashEL\epsilon;\\
		(\text{b}) &\LTEL[X]\models\epsilon;\\
		(\text{c}) &\LTEL[X]\models_{\overline{v}_{X}}\epsilon.
		\end{array}
		\]}
	Hence, {\em$\LTEL[X]\models X$} and, certainly, {\em$\LTEL[X]\models_{\overline{v}_{X}} X$}.
\end{prop}
\begin{proof}
	Let us denote: $\epsilon:=\alpha\approx\beta$. We will consecutively prove
	(a)$\Rightarrow$(b)$\Rightarrow$(c)$\Rightarrow$(a).
	
	Case:~(a)$\Rightarrow$(b). 	Let $X\vdashELog\alpha\approx\beta$ and $v$ be any valuation in $\LTEL[X]$.
	In virtue of Lemma~\ref{L:L-equivalence-with-LTL[X]}, there is a substitution $\sigma$ such that $v[\alpha]=\sigma(\alpha)\slash\theta_{\mathcal{E}_{\text{L}}}(X)$ and $v[\beta]=\sigma(\beta)\slash\theta_{\mathcal{E}_{\text{L}}}(X)$. On the other hand, according to~\eqref{E:implication-1} and \eqref{E:definition-there_EL(X)}, $\sigma(\alpha)\slash\theta_{\mathcal{E}_{\text{L}}}(X)=\sigma(\beta)\slash\theta_{\mathcal{E}_{\text{L}}}(X)$, that is $v[\alpha]=v[\beta]$.
	
	Case:~(b)$\Rightarrow$(c). It is obvious.
	
	Case:~(c)$\Rightarrow$(a). If $\LTEL[X]\models_{\overline{v}_{X}}\alpha\approx\beta$, then, by definition, 
	$\alpha\slash\theta_{\mathcal{E}_{\text{L}}}(X)=\beta\slash\theta_{\mathcal{E}_{\text{L}}}(X)$, that is, according to~\eqref{E:definition-there_EL(X)}, 
	$X\vdashEL\alpha\approx\beta$.
\end{proof}
\begin{rem}
	{\em Comparing Proposition~\ref{P:equivalence-with-LTEL[X]} with Proposition~\ref{P:equivalence-with-LT[X]}, we notice the difference in the conditions (b) in both propositions. This difference becomes clear when we recall that $\LTEL[X]\models X$ and thus $\LTEL[X]\models\epsilon$ is equivalent to the implication `$\LTEL[X]\models X\Rightarrow\LTEL[X]\models\epsilon$'.}	
\end{rem}

Let $\E$ be a nonempty class of E-matrices and $X\subseteq\EqL$. We denote:
\[
\E[X]:=\set{\mat{M}\in\E}{\mat{M}\models X}.
\]

It is obvious that $\E[X]\subseteq\E$, and, in virtue of Proposition~\ref{P:equivalence-with-LTEL[X]},
\[
X\vdashELog\epsilon~\Longleftrightarrow~\varnothing\vdash_{\mathcal{E}[X]_{\text{L}}}\epsilon.
\]

Now we denote:
\[
\E^{\ast}[X]:=\text{Mod}_{\EQ}(\E[X]).
\]
According to Proposition~\ref{P:equational-theory}, the classes $\E[X]$ and $\E^{\ast}[X]$ have the same equational theory; in our case, it implies that for any $\epsilon\in\EqL$,
\begin{equation}\label{E:equivalence-EL-and-E^star}
X\vdashELog\epsilon~\Longleftrightarrow~\varnothing\vdash_{\mathcal{E}^{\ast}[X]_{\text{L}}}\epsilon.
\end{equation}
(Exercise~\ref{section:LT-algebras-for-ELog}.\ref{EX:equivalence-EL-and-E^star})

The importance of the last equivalence is that $\LTEL[X]\in\E^{\ast}[X]$, though it is not necessary that $\LTEL[X]\in\E$. Thus, to prove that $X\not\vdashELog\epsilon$, one can use $\LTEL[X]$ as a separating matrix
(in the sense of Section~\ref{section:separating-means}) to show that $\LTEL[X]\not\models\epsilon$ or even $\LTEL[X]\not\models_{\overline{v}_{X}}\epsilon$. Since the class $\mathcal{A}[\E^{\ast}[X]]$ is a variety (Section~\ref{section:preliminaries-model-theory}), one of the last two rejections (and hence both) can be achieved if we notice the following.

\begin{prop}\label{P:ETLog-is-free-algebra}
	Let, given a set {\em$X\subseteq\EqL$}, $\E^{\ast}[X]$ be a class of E-matrices as defined above. Then the algebra of {\em$\LTEL[X]$}, that if $\FormAl\slash\theta_{\mathcal{E}_{\text{L}}}(X)$, is a free algebra with generators $\set{p\slash\theta_{\mathcal{E}_{\text{L}}}(X)}{p\in\VarL}$ in the variety $\mathcal{A}[\E^{\ast}[X]]$.
\end{prop}
\begin{proof}
	First, it is obvious that $\FormAl\slash\theta_{\mathcal{E}_{\text{L}}}(X)$ is generated by the set $\set{p\slash\theta_{\mathcal{E}_{\text{L}}}(X)}{p\in\VarL}$.
	
	Now let $\alg{A}$ be an arbitrary algebra from $\mathcal{A}[\E^{\ast}[X]]$, that is the matrix $\mat{M}:=\langle\alg{A},\Delta\rangle$ belongs to $\E^{\ast}[X]$. Let us consider a map
	\[
	f:\set{p\slash\theta_{\mathcal{E}_{\text{L}}}(X)}{p\in\VarL}\longrightarrow|\alg{A}|.
	\]
	
	We define a valuation $v$ in \alg{A} as follows:
	\[
	v:p\mapsto f(p\slash\theta_{\mathcal{E}_{\text{L}}}(X)),
	\]
	for each $p\in\VarL$.
	
	According to Proposition~\ref{P:valuation}, $v$ can be extended to a homomorphism (denoted also by $v$) of $\FormAl$ into \alg{A}. It must be obvious that
	\[
	v[\alpha]=f(\alpha\slash\theta_{\mathcal{E}_{\text{L}}}(X)),
	\]
	for any $\alpha\in|\FormAl|$.
	
	We notice that 
	\[
	\theta_{\mathcal{E}_{\text{L}}}(X)\subseteq\text{ker}(v).
	\]
	
	Indeed, suppose $(\alpha,\beta)\in\theta_{\mathcal{E}_{\text{L}}}(X)$. Then
	$\alpha\slash\theta_{\mathcal{E}_{\text{L}}}(X)=\beta\slash\theta_{\mathcal{E}_{\text{L}}}(X)$ and, hence, $v[\alpha]=v[\beta]$.
	
	Adapting the definition~\eqref{E:definition-of-fraction-congruence} to our situation, we obtain:
	\[
	\text{ker}(v)\slash\theta_{\mathcal{E}_{\text{L}}}(X)=
	\set{(\alpha\slash\theta_{\mathcal{E}_{\text{L}}}(X),\beta\slash\theta_{\mathcal{E}_{\text{L}}}(X))}{\alpha,\beta\in\FormsL~\text{and}~v[\alpha]=v[\beta]}.
	\]
	
	From Section~\ref{section:prelimenaries-algebra}, we know that $\text{ker}(v)\slash\theta_{\mathcal{E}_{\text{L}}}(X)$ is a congruence on
	$\FormAl\slash\theta_{\mathcal{E}_{\text{L}}}(X)$.
	
	Next, we define the following maps:
	\[
	\begin{array}{l}
	g:\alpha\slash\theta_{\mathcal{E}_{\text{L}}}(X)\mapsto(\alpha\slash\theta_{\mathcal{E}_{\text{L}}}(X))\slash(\text{ker}(v)\slash\theta_{\mathcal{E}_{\text{L}}}(X));\\\\
	h:(\alpha\slash\theta_{\mathcal{E}_{\text{L}}}(X))\slash(\text{ker}(v)\slash\theta_{\mathcal{E}_{\text{L}}}(X))\mapsto\alpha\slash\text{ker}(v);\\\\
	j:\alpha\slash\text{ker}(v)\mapsto v[\alpha].
	\end{array}
	\]
	
	The map $g$ is obviously a homomorphism; in virtue of Proposition~\ref{P:isomorphism-second}, $h$ is an isomorphism; according to Proposition~\ref{P:isomorphism-first}, $j$ is a homomorphism. Then, we define:
	\[
	\overline{f}:=j\circ h\circ g.
	\] 
	
	It is obvious that $\overline{f}$ is a homomorphism of $\FormAl\slash\theta_{\mathcal{E}_{\text{L}}}(X)$ into \alg{A}.
	
	Finally, we observe:
	\[
	\begin{array}{rl}
	\overline{f}(p\slash\theta_{\mathcal{E}_{\text{L}}}(X))\!\!\!&=j(h(g(p\slash\theta_{\mathcal{E}_{\text{L}}}(X))))\\
	&=j(h((p\slash\theta_{\mathcal{E}_{\text{L}}}(X))\slash(\text{ker}(v)\slash\theta_{\mathcal{E}_{\text{L}}}(X))))\\
	&=j(p\slash\text{ker}(v))\\
	&=v[p]=f(p\slash\theta_{\mathcal{E}_{\text{L}}}(X)).
	\end{array}
	\] 
	
	Thus $\overline{f}$ is an extension of $f$.
\end{proof}
\paragraph{Exercises~\ref{section:LT-algebras-for-ELog}}
\begin{enumerate}
	\item\label{EX:L-equivalence-with-LTL[X]} Prove Lemma~\ref{L:L-equivalence-with-LTL[X]}.
	\item\label{EX:equivalence-EL-and-E^star}Prove the equivalence~\eqref{E:equivalence-EL-and-E^star}.
\end{enumerate}

\section{Interconnection between $\E$ and $\ELog$}
In this section we show that there is an interconnection between abstract logics $\E$ and $\EL$, where substitution play a key role.\\

We begin with expanding the operator $\Sb$ (Section~\ref{section:languages}) from the set $\FormsL$ to $\EqL$. Namely, for any $X\subseteq\EqL$, we define:
\[
\Sb X:=\set{\sigma(\alpha)\approx\sigma(\beta)}{\sigma\in\SbsL~\text{and}~\alpha\approx\beta\in X}.
\]
\begin{prop}\label{P:interconnection-E-and-EL}
	Let $\E$ be a nonempty class of E-matrices with $\Su(\mathcal{A}[\E])\subseteq\mathcal{A}[\E]$. Then for any set {\em$X\cup\lbrace\epsilon\rbrace\subseteq\EqL$},
	{\em\[
		X\vdashELog\epsilon~\Longleftrightarrow~\Sb X\vdashE\epsilon.
		\]}
\end{prop}
\begin{proof}
	Suppose that $X\vdashELog\epsilon$. That is, for any $\mat{M}\in\E$,
	\[
	\mat{M}\models X~\Longrightarrow~\mat{M}\models\epsilon.\tag{\ref{P:interconnection-E-and-EL}--$\ast$}
	\]
	
	Now let us take an arbitrary $\mat{M}\in\E$ and valuation $v_0$ in \alg{A}, the algebra of \mat{M}. Assume that $\mat{M}\models_{v_0}\Sb X$. Next we denote:
	\[
	\alg{A}_0:=\lbrack\set{v_0[p]}{p\in\VarL}\rbrack_{\alg{A}}~\text{and}~\mat{M}_0:=\langle\alg{A}_0,\Delta\rangle.
	\]
	
	Since $\alg{A}_0$ is a subalgebra of \alg{A}, by premise, $\alg{A}_0\in\mathcal{A}[\E]$ and, hence, $\mat{M}_0\in\E$. Also, it is clear that $\mat{M}_0\models_{v_0}\Sb X$. 
	
	Now, we aim to show that $\mat{M}_0\models X$. Let $v$ be any valuation in $\alg{A}_0$. According to Proposition~\ref{P:valuation-in-algebra-generated-by-X}, there is a substitution $\sigma$ such that $v=v_{0}\circ\sigma$. This implies that for any $\eta\in X$, $v[\eta]=v_0[\sigma(\eta)]$. Since $\sigma(\eta)\in\Sb X$, then, by premise, $v[\eta]\in\Delta$.
	
	Using (\ref{P:interconnection-E-and-EL}--$\ast$), we derive that $\mat{M}_0\models\epsilon$; in particular, $\mat{M}_0\models_{v_0}\epsilon$. The latter implies that $\mat{M}\models_{v_0}\epsilon$.
	
	Conversely, assume that $\Sb X\vdashE\epsilon$. Let us take an arbitrary $\mat{M}\in\E$ and suppose (as a premise) that $\mat{M}\models X$. Then, let us take any valuation $v_1$ in \mat{M} and any equality $\eta\in\Sb X$. By definition, there is a substitution $\sigma$ such that $\eta=\sigma(\zeta)$, for some $\zeta\in X$. Now we define: $v_2:=v_1\circ\sigma$. Since $v_2$ is a valuation in $\mat{M}$ and, by premise, $v_2[\zeta]\in\Delta$, $v_1[\eta]\in\Delta$. The latter means that $\mat{M}\models_{v_1}\Sb X$. By premise, we obtain that $\mat{M}\models_{v_1}\epsilon$.
\end{proof}

\section{$\EL$-models}
Let $\E$ be a class of E-matrices of type $\Lan$ and let $\EL$ be an abstract logic of equational L-consequence which is determined by $\E$. An E-matrix $\mat{M}$ of type $\Lan$ is an $\EL$-\textit{\textbf{model}} \index{$\EL$-model} if for any set $X\cup\lbrace\epsilon\rbrace\subseteq\EqL$, $X\vdashEL\epsilon$ implies
\[
\mat{M}\models X\Longrightarrow\mat{M}\models\epsilon.
\]

It is important to note that the class of all $\EL$-models can be larger than the class $\E$. And it is the class of all $\EL$-models should be taken into account if, in the light of the conception of separating tools (Section~\ref{section:separating-means}), we aim to refute $X\vdashEL\epsilon$. For if $\mat{M}$ is an E-matrix such that $\mat{M}\models X$ and $\mat{M}\not\models\epsilon$, then $X\not\vdashEL\epsilon$.\\

With this in mind, we define:
\[
\textbf{Mod}_{\hat{\E}}(X,\epsilon):=\set{\mat{M}\in\hat{\E}}{\mat{M}\models X\Rightarrow\mat{M}\models\epsilon},
\]
where, we remind, $\hat{\E}$ is the class of all E-matrices of type $\Lan$;
\[
\textbf{Mod}_{\hat{\E}}(X):=\set{\mat{M}\in\hat{\E}}{\mat{M}\models X};
\]
and as a particular case of the last definition
\[
\textbf{Mod}_{\hat{\E}}(\epsilon):=\set{\mat{M}\in\hat{\E}}{\mat{M}\models\epsilon}.
\]
(These definitions imitate definitions in Section~\ref{section:realizations-abstract-logic}.)\\

Applying Proposition~\ref{P:semantic-consiquence-equality-models} to the new definitions above, we obtain the equality:
\begin{equation}\label{E:El-models}
\textbf{Mod}_{\hat{\E}}(X,\epsilon)=\textbf{Mod}_{\hat{\E}}(\epsilon)\cup(\hat{\E}\setminus\textbf{Mod}_{\hat{\E}}(X)).
\end{equation}

The equality~\eqref{E:El-models} is especially useful if we intend to separated $X$ from $\epsilon$ by an E-matrix. With this in mind, we have:
\[
\begin{array}{rl}
\mat{M}\notin\textbf{Mod}_{\hat{\E}}(X,\epsilon)\!\!\! &\Longleftrightarrow
\mat{M}\in\textbf{Mod}_{\hat{\E}}(X)\setminus\textbf{Mod}_{\hat{\E}}(\epsilon)\\
&\Longleftrightarrow \mat{M}\models X~\text{and}~\mat{M}\not\models\epsilon;\\
&\text{applying Proposition~\ref{P:EQ[M]=EQ[{M}_si]}, we obtain:}\\
&\Longleftrightarrow \Si\lbrace\mat{M}\rbrace\models X~\text{and}
~\Si\lbrace\mat{M}\rbrace\not\models\epsilon.\\
\end{array}
\]

Thus an E-matrix $\mat{M}$ separates $X$ from $\epsilon$ if, and only if, every equality of $X$ is valid in each matrix of $\Si\lbrace\mat{M}\rbrace$ and there is a matrix in $\Si\lbrace\mat{M}\rbrace$ which rejects $\epsilon$. We remind the the class $\Si\lbrace\mat{M}\rbrace$ consists of those E-matrices whose carriers are subdirectly irreducible homomorphic images of the carrier of $\mat{M}$. (Proposition~\ref{P:characterisation-of-Si-operator})

\section{Examples of logics $\ELog$}
We consider two examples --- the L-equational consequence determined by the class $\EB$ and by the class $\EH$ of Sections~\ref{section:Boolean-equational} and~\ref{section:Heyting-equational}, respectively. We denote the corresponding L-equational logics by $\EBLog$ and $\EHLog$, respectively.\\

Let $X\subseteq\EqL$. It is clear that $\mathcal{A}[\EBLog[X]]$ is a variety of Boolean algebras. According to Proposition~\ref{P:varieties-Boolean-algebras}, either all algebras of $\mathcal{A}[\EBLog[X]]$ are degenerate or $\mathcal{A}[\EBLog[X]]=\mathcal{A}[\EB]$. Thus, in the second case, $X\vdash_{\EBLog}\epsilon$ if, and only if, $\EB\models\epsilon$. And the latter is equivalent to that a 2-element Boolean E-matrix validates $\epsilon$.
The problem is that if $X$ is infinite, we may not know what our case is. Of course, if $X$ is finite, then it is easy to check whether $\LTEL[X]$ is degenerate or not, since  $\LTEL[X]$ is not degenerate if, and only, if it contains a 2-element Boolean algebra as its subalgebra, and the last condition is equivalent to whether all equalities of $X$ are validated by the 2-element Boolean algebra.\\

Our next example is the abstract logic $\EHLog$. Given a set $X\subseteq\EqL$, we obtain a variety $\mathcal{A}[\EHLog[X]]$ of Heyting algebras. First of all, we note that there is a continuum of such varieties.\footnote{Cf.~\cite{blok1977} and references therein.}\label{footnote:varieties-of-Heyting-algebras} Even if we restrict to the finite sets $X$, the problem whether $X\vdash_{\EHLog}\epsilon$ is algorithmically undecidable;  cf.~\cite{shehtman1978},~\cite{chagrov-zakhar1997}, section 16.5, with reference to the previous footnote. And here it is appropriate to recall Proposition~\ref{P:ETLog-is-free-algebra} and Proposition~\ref{P:preservation}. For in the hope to prove that
$X\not\vdash_{\EHLog}\epsilon$, in the light of Proposition~\ref{P:ETLog-is-free-algebra}, we can try to find such an algebra $\alg{A}\in\mathcal{A}[\EHLog[X]]$ that the matrix
$\mat{M}=\langle\alg{A},\Delta\rangle$ refutes $\epsilon$. Such an algebra \alg{A} must be a homomorphic image of $\LTEL$ (Proposition~\ref{P:ETLog-is-free-algebra}), so all equalities of $X$ are validated by \mat{M} (Proposition~\ref{P:preservation}).

\section{L-Equational consequence based on\\ implicational logic}\label{section:L-equational-based-on-implicational}

In this section we use some notions and their notation of Section~\ref{section:equational-based-on-implicational}. As in that section, let $\mathcal{S}$ be an implicational logic with respect to $\im(p,q)$ in a language $\Lan$. \\

Let $\sigma$ be an arbitrary $\Lan$-substitution. Then the following equalities are obvious:
\[
\begin{array}{rl}
\sigma(\im(\alpha,\beta))\!\!\! &=\im(\sigma(\alpha),\sigma(\beta)),\\
\textbf{i}(\sigma(\alpha\approx\beta))=\textbf{i}(\sigma(\alpha)\approx\sigma(\beta))\!\!\! &=\lbrace\im(\sigma(\alpha),\sigma(\beta)),\im(\sigma(\beta),\sigma(\alpha))\rbrace,\\
&= \sigma(\textbf{i}(\alpha\approx\beta)).
\end{array}
\]

This implies that for any set $X\subseteq\EqL$,
\[
(\sigma(X))^{\textbf{i}}=(\sigma(X^{\textbf{i}}))
\]
and hence
\begin{equation}\label{E:EL-based-on-implicational}
(\Sb X)^{\textbf{i}}=\Sb X^{\textbf{i}}.
\end{equation}

Similarly, for the function \textbf{s}, we have:
\[
\textbf{s}(\sigma(\alpha\approx\beta))=\textbf{s}(\sigma(\alpha)\approx\sigma(\beta))=\sigma(\alpha)\leftrightarrow\sigma(\beta).
\]

This in turn implies that for any set $X\subseteq\EqL$,
\[
(\sigma(X))^{\textbf{s}}=(\sigma(X^{\textbf{s}}))
\]
and hence
\begin{equation}\label{E:EL-based-on-implicational-2}
(\Sb X)^{\textbf{s}}=\Sb X^{\textbf{s}}.
\end{equation}

\begin{prop}\label{P:L-equational-on-implicational-1}
	Let $\mathcal{S}$ be an implicational logic with respect to $\im(p,q)$ in $\Lan$. Also, let $\EL$ be an L-equational logic with $\mathcal{A}[\E]$ such that $\set{\langle\alg{A},\one\rangle}{\alg{A}\in\mathcal{A}[\E]}=\QS$. Then for any {\em$X\cup\lbrace\epsilon\rbrace\subseteq\EqL$},
	{\em\[
		X\vdashEL\epsilon~\Longleftrightarrow~\textbf{i}(\epsilon)\subseteq\textbf{Cn}_{\mathcal{S}_{\text{sub}}}(X^{\textbf{i}}).
		\]}
	Further, if $\mathcal{S}$ is implicative and implicational with respect to $p\rightarrow q$, and, in addition, if $\mathcal{S}$ is sound with respect to the inferences rules {\em(\text{a}--$i$)} and {\em(\text{b}--$i$)} $($Section~\ref{section:inference-rules}$)$, then for
	any {\em$X\cup\lbrace\epsilon\rbrace\subseteq\EqL$},
	{\em\begin{equation}\label{E:equational-on-implicational-1}
		X\vdashEL\epsilon~\Longleftrightarrow~X^{\textbf{s}}\vdash_{\mathcal{S}_{\text{sub}}}\textbf{s}(\epsilon).
		\end{equation}}
\end{prop}
\begin{proof}
	The first part is proved as follows:
	\[
	\begin{array}{rl}
	X\vdashEL\epsilon\!\!\! &\Longleftrightarrow \Sb X\vdashE\epsilon\quad[\text{by Proposition~\ref{P:interconnection-E-and-EL}}]\\\\
	&\Longleftrightarrow \textbf{i}(\epsilon)\subseteq\ConS{(\Sb X)^{\textbf{i}}}\quad[\text{by Proposition~\ref{P:equational-on-implicational-1}}]\\\\
	&\Longleftrightarrow \textbf{i}(\epsilon)\subseteq\ConS{\Sb X^{\textbf{i}}}\quad[\text{by~\eqref{E:EL-based-on-implicational}}]\\\\
	&\Longleftrightarrow \textbf{i}(\epsilon)\subseteq\textbf{Cn}_{\mathcal{S}_{\text{sub}}}(X^{\textbf{i}})\quad[\text{by definition~\eqref{E:def-of-aLogSub}}].
	\end{array}
	\]
	
	Next, we turn to the second part. Indeed, we obtain:
	\[
	\begin{array}{rl}
	X\vdashEL\epsilon\!\!\! &\Longleftrightarrow \Sb X\vdashE\epsilon\quad[\text{by Proposition~\ref{P:interconnection-E-and-EL}}]\\\\
	&\Longleftrightarrow (\Sb X)^{\textbf{s}}\vdash_{\mathcal{S}}\textbf{s}(\epsilon)\quad[\text{by Proposition~\ref{P:equational-on-implicational-1}}]\\\\
	&\Longleftrightarrow \Sb X^{\textbf{s}}\vdash_{\mathcal{S}}\textbf{s}(\epsilon)\quad[\text{by~\eqref{E:EL-based-on-implicational-2}}]\\\\
	&\Longleftrightarrow X^{\textbf{s}}\vdash_{\mathcal{S}_{\text{sub}}}\textbf{s}(\epsilon)\quad[\text{by definition~\eqref{E:def-of-aLogSub}}].
	\end{array}
	\]
\end{proof}

Below we consider two examples.\\

First, let us consider the two abstract logics --- $\Cl$ and $\Cl^{\star}$ (Section~\ref{section:inference-rules}). The former logic, being formulated in the language $\Lan_A$, is implicative and implicational with respect to $p\rightarrow q$. (See Exercise\ref{section:unital-logics}.\ref{EX:Cl-In-implicative} and Exercise~\ref{section:lindenbaum-algebra}.\ref{EX:Cl-Int-implicational}.) 
On the other hand, the abstract logic $\Cl^{\star}$ is not even structural (Exercise~\ref{section:modus-rules}.\ref{EX:Cl^star-is-not-structural}) However, in view of the second part of Proposition~\ref{P:L-equational-on-implicational-1} and Proposition~\ref{P:S=S^star}, we obtain the following.
\begin{prop}\label{P:logic-EBLog}
	For any set {\em$X\cup\lbrace\epsilon\rbrace\subseteq\EqL$}, the following conditions are equivalent:
	{\em\[
		\begin{array}{cl}
		(\text{a}) & X\vdash_{\EBLog}\epsilon;\\
		(\text{b}) & X^{\textbf{s}}\vdash_{\textsf{Cl}_{\text{sub}}}\textbf{s}(\epsilon);\\
		(\text{c}) & X^{\textbf{s}}\vdash_{\textsf{Cl}^{\star}}\textbf{s}(\epsilon).
		\end{array}
		\]}
\end{prop} 

In view if Proposition~\ref{P:Cl^star-theories} and the equivalence of (a) and (c) of Proposition~\ref{P:logic-EBLog}, we conclude that there are only two $\EBLog$-theories --- the set of equalities $\epsilon$ such that $\textbf{s}(\epsilon)$ is a classical tautology and the set $\EqL$. Therefore, for finite sets $X$, the problem whether $X\vdash_{\EBLog}\epsilon$ can be easily decided; this problem even is algorithmically decidable. (See Chapter~\ref{chapter:decidability} for other applications  of the last concept.) \\

Next, we define an abstract logic $\Int^{\star}$ that is based on the abstract logic $\Int$ (Section~\ref{section:subclasses-unital}). The $\Int^{\star}$-theories are known as \emph{superintuitionistic logics}. Using the same argumentation as in the proof of Proposition~\ref{P:logic-EBLog}, we get the corresponding proposition in this case.

\begin{prop}\label{P:logic_EHLog}
	For any set {\em$X\cup\lbrace\epsilon\rbrace\subseteq\EqL$}, the following conditions are equivalent:
	{\em\[
		\begin{array}{cl}
		(\text{a}) & X\vdash_{\EHLog}\epsilon;\\
		(\text{b}) & X^{\textbf{s}}\vdash_{\textsf{Int}_{\text{sub}}}\textbf{s}(\epsilon);\\
		(\text{c}) & X^{\textbf{s}}\vdash_{\textsf{Int}^{\star}}\textbf{s}(\epsilon).
		\end{array}
		\]}
\end{prop}

The last proposition shows that the abstract logic $\EBLog$ appears to be complicated. One of the reasons of its analysis is that, similarly to the logic $\Int^{\star}$, in view of the equivalence of (a) and (c) of Proposition~\ref{P:logic_EHLog}, the cardinality of $\EHLog$-theories equals $2^{\card{\textbf{Fm}_{\Lan_A}}}$, providing that $\card{\mathcal{V}_{\Lan_A}}\ge\aleph_0$.\footnote{See the footnote on p.~\pageref{footnote:varieties-of-Heyting-algebras}.} On the other hand, the problem of characterizing some of the
$\EHLog$-theories seem to be quite difficult. One of the examples of such a theory is the $\EHLog$-theory associated with \emph{Medvedev's logic} which is a finitely approximable but not finitely axiomatizable $\Int^{\star}$-theory.\footnote{About Medvedev's logic, see~\cite{chagrov-zakhar1997}, section 2.9.}

\section{Historical notes}
An L-equational consequence first appeared in~\cite{birkhoff1935}.
However, Birkhoff’s approach is not as direct as ours, as well as the ones shown~in~\cite{pigozzi1975}, section 1.3,  in~\cite{burris-sankappan81}, chapter II, {\S}14, and in~\cite{burris1998}, chapter 3. Therefore, it would be interesting to look at the Birkhoff approach from our point of view on the equational L-consequence.

First of all, Birkhoff does not distinguish E-matrices and their algebras. To facilitate our analysis, in this section, we do not do it too. This allows as to use the phrases such as `a homomorphic image of an E-matrix', `a subalgebra of an E-matrix', and `the direct product of a family of E-matrices'. Also, the classes $\mathcal{A}[\E]$ and $\E$ will not be distinguished.

One of the difficulties of reading~\cite{birkhoff1935} is that he does not make a clear distinction between substitution as an operation on terms and a valuation in algebra. In this regard, we must remember that any valuation in the formula algebra $\FormAl$ is a substitution.

The third remark, before we begin our analysis, is that Birkhoff deals mainly with $\EELog$, that is, when the class of algebras (of E-matrices) contains all algebras of type $\Lan$; and he uses $\E\subseteq\EE$ only to form the variety $\EQ[\E]$ or even the lattice of all such varieties. It was this trend that influenced the further development of Birkhoff's ``equational logic'';  cf.~\cite{tarski1966,pigozzi1975,taylor1979}.

The equivalence between formally inferred equalities from a given set $X$ and those which follow from $X$ in the logic $\EELog$ is established by Birkhoff in the first statement of Theorem 9. As a formal system, Birkhoff uses E1s. For convenience, we denote:
\[
\textbf{Cn}_{\text{E1s}}(X):=\set{\epsilon\in\EqL}{X\vdash_{\text{E1s}}\epsilon},
\]
where $X\subseteq\EqL$. We note that we do not claim at this moment that $\textbf{Cn}_{\text{E1s}}$ is a consequence operator, although afterwards it is; see~Remark~\ref{R:remark-on-E1-E2-E3}. Thus, the first statement of Theorem 9 reads:
\begin{equation}\label{E:birkhoff-1}
\epsilon\in\textbf{Cn}_{\text{E1s}}(X)~\Longleftrightarrow~
\EE[X]\models\epsilon.
\end{equation}

The $\Rightarrow$-implication of~\eqref{E:birkhoff-1} is essentially the $\Rightarrow$-implication of Proposition~\ref{P:EE-L-consequence-completeness}, that is the soundness of $\vdash_{\text{E1s}}$.

Turning to the $\Leftarrow$-implication of~\eqref{E:birkhoff-1}, we first recall that for any substitution $\sigma$,
\begin{equation}\label{E:birkhoff-2}
\epsilon\in\textbf{Cn}_{\text{E1s}}(X)~\Longrightarrow~\sigma(\epsilon)\in\textbf{Cn}_{\text{E1s}}(X).
\end{equation}

Further, since $\EE[X]$ is an equational class (or a variety), it contains the algebra $\textbf{LT}_{\hat{\mathcal{E}}_{\text{L}}}[X]$ (see~\eqref{E:defn-L-congruence} and Definition~\ref{D:LT-for-ELog}) which is a free algebra in this variety (Proposition~\ref{P:ETLog-is-free-algebra}).\footnote{Namely~\eqref{E:birkhoff-2} is used to show that $\textbf{LT}_{\hat{\mathcal{E}}_{\text{L}}}[X]$  validates all equalities of $X$.} This implies that if
$\epsilon\notin\textbf{Cn}_{\text{E1s}}(X)$, then, by definition, $\textbf{LT}_{\hat{\mathcal{E}}_{\text{L}}}[X]$ refutes $\epsilon$.\\

The proof of completeness theorem in the aforementioned works~\cite{pigozzi1975},  \cite{burris-sankapp1981}, and~\cite{burris1998} is similar to ours, where, E2s is used  as a formal system. In all these cases, it is demonstrated that the rule E1--$c$ is derivable in E2.\\

Perhaps, under the influence of Tarski for many years, studies of equational L-consequence were limited to considering equational classes, or varieties, of algebras of the same similarity type. Tarski wrote in~\cite{tarski1966}, section 1:
\begin{quote}
	``Two system of equational logic differ only in operational symbols.''
\end{quote}

The focus, then, was on equational theories, where
\begin{quotation}
	``[by] the \emph{equational theory} of an algebra $\mathfrak{A}$, or a class \textsf{K} of algebras, we understand the set of all equations, in the system of logic determined by [a similarity type] $\sigma$, which are identically satisfied in the algebra $\mathfrak{A}$, or in all algebras of \textsf{K}.'' \textit{Ibid}
\end{quotation}

Thus, in our notation, given a nonempty class $\E$ of E-matrices, the algebras of which form a variety, all questions are about the set $\textbf{Cn}_{\ELog}(\varnothing)$; for instance: whether it is decidable?
whether it consists of the equalities derivable from its finite subset (\emph{finitely based theories})? and some others. The reader can find a detailed discussion of these issues in~\cite{pigozzi1975} and~\cite{taylor1979}.\\

A. Selman in~\cite{selman1972} goes even further from the L-equational consequence.
He uses the format of Definition~\ref{D:LT-for-ELog} to prove a Birkhoff-style completeness result for any set $X$ of quasi-equalities as premises and a quasi-equality $\eta_1\land\ldots\land\eta_n\rightarrow\epsilon$ as a conclusion.

\chapter{$\mathcal{Q}$-Consequence}\label{chapter:Q-consequence}

\section{$\Q$-Languages}\label{section:Q-languages}
$\Q$-languages, as the reader will see below, admit logical constants of a new type. This, however, requires a restructuring of the whole grammar of the formation of judgment. Therefore, we will start with a new formal language and its grammar. As in Section~\ref{section:languages}, we deal most of the time with a schematic
$\Q$-language which we denote by $\LanQ$. This language includes the (\textit{\textbf{schematic}}) \textit{\textbf{language of terms}}, $\LanT$, which is extended to the \textbf{\textit{language of formulas}}, which is denoted by $\LanQ$. In this language, we designate a sublanguage, $\LanF$,  which includes only logical connectives and no quantifiers.\\

All grammatical categories are listed below as follows.
\begin{itemize}
	\item a denumerable (or infinitely countable) set $\VarT$ of \textit{\textbf{bound}} (individual) \textit{\textbf{variables}} that are ordered by type of $\omega$:
	\begin{equation}\label{E:bound-variables}
	x_1,~x_2,~\ldots;
	\end{equation}
	unspecified bound variables are denoted by the letters $x, y$ with or without subscripts; sometimes, we will need to use finite lists of bound variables, they are denoted by over-lined letters like $\overline{x},\ldots$; 
	\item a denumerable (or infinitely countable) set $\ParT$ of \textit{\textbf{free}} (individual) \textit{\textbf{variables}} (called also \textit{\textbf{parameters}}) that are ordered by type of $\omega$:
	\begin{equation}\label{E:free-variables}
	u_1,~u_2,~\ldots;
	\end{equation}
	unspecified free variables are denoted by $u$ with or without superscript or with an apostrophe or asterisk; finite (maybe void) lists of bound variables are denoted by over-lined letters like $\overline{u},\overline{u^\ast},\ldots$; 
	\item a set (maybe empty) of \textit{\textbf{function symbols}}\index{function symbol} is denoted by  $\FuncT$ and the elements of $\FuncT$ by letters $f, g$ (with or without subscripts); a nonnegative integer number is assigned to each element $f\in\FuncT$; this number is called the \textit{arity} of $f$; we denote the arity of $f$ by $\#(f)$; the function symbols of arity $0$, if any, are called  (individual) \textit{\textbf{constants}};\index{constant} the set of all constants is denoted by $\ConsT$ (maybe empty); thus $\ConsT\subseteq\FuncT$;
	\item a nonempty set $\PredF$ of \textit{\textbf{predicate symbols}\index{predicate symbol}} which will be denoted by letters $p,q,\ldots$ (with or without subscripts); each predicate symbol $p$ also has an arity which is greater than or equal to $1$ and is denoted by $\#(p)$;  the equality symbol `$\approx$' may or may not be in $\LanF$; but if it is present, its arity equals $2$;
	\item a set (maybe empty) $\FuncF$ of \textit{\textbf{logical connectives}};\index{logical connectives} unspecified logical connectives will be denoted by letters $F, G, H$; each logical connective has an arity whose value is a nonnegative integer; a logical connective of arity $0$ is called a \textit{\textbf{logical constant}}; the set of logical constants (maybe empty) is denoted by $\ConsF$; thus $\ConsF\subseteq\FuncF$;
	\item a nonempty set $\Q$ of (\textbf{\textit{generalized}}) \textit{\textbf{quantifiers}};\index{quantifier} each quantifier $Q\in\Q$ has an arity greater than $0$.
	\item as usual, the parentheses, `(' and `)', are used as punctuation marks.
\end{itemize}

All these categories, except quantifiers, are the categories of language $\LanF$.\\

The next definition addresses three concepts; the first one, the concept of a $\LanT$-\emph{term expression}, is an auxiliary one; it will, however, help us define two other concepts.
\begin{defn}[$\LanT$-prototerm and $\LanT$-term]\label{D:LQ-terms}	\index{$\LanT$-prototerm}
$\LanT$-prototerms are defined inductively by the following clauses:
\begin{itemize}	
\item[{\em(a)}]	each element of the set $\VarT\cup\ParT\cup\ConsT$ is an $\LanT$-prototerm;
\item[{\em(b)}]	if $t_1,\ldots,t_n$ are $\LanT$-prototerms and $f\in\FuncT$ has arity $n\ge 1$, then $ft_1\ldots t_n$ is an $\LanT$-prototerm;
\item[{\em(c)}]	the $\LanT$-prototerms are only those expressions that can be constructed by {\em (a)} and {\em(b)}.
\end{itemize}

An $\LanT$-\textbf{term}\index{$\LanT$-term} is an $\LanT$-prototerm which does not contain bound variables.
\end{defn}

The set of all terms is denoted by $\Tm$.\index{$\Tm$}

Similar to the definition of the formula algebra of type $\Lan$ of Section~\ref{section:semantics},
we define a \textit{\textbf{term algebra}}\index{algebra!term} of type $\LanT$, whose carrier is $\Tm$ and the fundamental operations are the function symbols of $\FuncT$ (including the constants of $\ConsT$). We denote this algebra by $\TmAl$; that is,\index{$\TmAl$}
\[
\TmAl:=\lr{\Tm;\FuncF}.
\]

It is allowed to substitute (in the sense of Section~\ref{section:languages}) only terms in those prototerms that contain free variables and only for free variables. 
Thus, although an $\LanT$-\textit{\textbf{substitution}} can be regarded as a map $\sigma:\ParT\longrightarrow\Tm$ which can, in the usual way, be extended to a homomorphism (denoted by the same symbol) $\sigma:\TmAl\longrightarrow\TmAl$, we also allow the situations when substitutions apply to prototerms, setting by definition for any $\LanT$-substitution $\sigma$,
\[
\sigma(t(x_{i_1},\ldots,x_{i_n},u_{i_1},\ldots,u_{i_m})):=t(x_{i_1},\ldots,x_{i_n},\sigma(u_{i_1}),\ldots,\sigma(u_{i_m})),
\]
where $x_{i_1},\ldots,x_{i_n}$ are all bound variables and $u_{i_1},\ldots,u_{i_m}$ are all parameters that occur in a prototerm $t$.

 The \textit{\textbf{identity substitution}} is denoted by $\iota$; that is, 
\[
\iota(t)=t,
\]
for any prototerm $t$.

Given an $\LanT$-substitution $\sigma$, we will often use another substitution which differs from $\sigma$ at finitely many elements of its domain.

Let $\overline{u}=(u_{i_1},\ldots,u_{i_n})$ be a nonvoid list of pairwise distinct free variables, $\overline{t}=(t_1,\ldots,t_n)$ be a list of terms, not necessarily pairwise distinct, and $\sigma$ be an $\LanT$-substitution. Then we define:
\[
\delta_{(\sigma,\overline{u}\backslash\overline{t})}(u^{\prime}):=
\begin{cases}
\begin{array}{cl}
t_k &\text{if $u^{\prime}=u_{i_k}$}\\
\sigma(u^{\prime}) &\text{otherwise}.
\end{array}
\end{cases}
\]

Also, just for given $\overline{u}$ and $\overline{t}$ as above, we define an $\LanT$-substitution:
\[
\sigma_{\overline{u}\backslash\overline{t}}(u^\prime):=\begin{cases}
t_k &\text{if $u^{\prime}=u_{i_k}$}\\
u^{\prime} &\text{otherwise}.
\end{cases}
\]

If $\overline{u}=(u)$ and $\overline{t}=(t)$, we will simply write $\sigma_{u\backslash t}$.

In the sequel, we will often use the notation:
\[
\sigma(\overline{t}):=(\sigma(t_1),\ldots,\sigma(t_n));
\]
in particular,
\[
\sigma(\overline{u})=(\sigma(u_{i_1}),\ldots,\sigma(u_{i_n})).
\]

Our next concepts are the (auxiliary) concept of $\LanQ$-\emph{formula expression}, $\LanQ$-protoformula and $\LanQ$-formula. These concepts are based on the signature that includes logical connectives and quantifiers of $\LanQ$. However, it will be convenient to refer to the logical connectives of $\LanQ$ as the \textit{\textbf{signature of type}} $\LanF$.

\begin{defn}[$\LanQ$-protoformula, free occurrences of bound variables, and $\LanQ$-formula]\label{D:LQ-formula}\index{$\LanQ$-protoformula}\index{$\LanQ$-formula}
The $\LanQ$-formula expressions can be built up according to the following clauses:
\begin{itemize}
	\item[{\em(a)}] $p(t_1,\ldots,t_n)$ is an $($atomic$)$ $\LanQ$-\textbf{protoformula} if $p$ is a predicate symbol of arity $n$ and $t_1,\ldots,t_n$ are $\LanT$-prototerms; all occurrences of bound variables, if any, that occur in $p(t_1,\ldots,t_n)$ \textbf{occur freely} in it, or are \textbf{free occurrences};
	\item[{\em(b)}] if $F$ is a logical connective of arity $n$ and $A_1,\ldots, A_n$ are $\LanQ$-protoformulas, then $FA_1\ldots A_n$ is an $\LanQ$-\textbf{protoformula}; all bound variables that occur freely in at least one of the protoformulas $A_i$, \textbf{occur freely} also in $FA_1\ldots A_n$; \index{protoformula!free occurence}
	\item[{\em(c)}] if $Q$ is a quantifier of arity $n$ and
	$A$ is a protoformula that contains free occurrences of bound variables $x_{i_1},\ldots,x_{i_n}$, then $Q x_{i_1},\ldots,x_{i_n} A$ is an $\LanQ$-\textbf{protoformula}, in which case $A$ is called the \textbf{scope}
	of $Q x_{i_1},\ldots,x_{i_n}$; the free occurrences of bound variables that occur freely in $A$, but are not listed among $x_{i_1},\ldots,x_{i_n}$, are \textbf{free occurrences} of the protoformula $Q x_{i_1},\ldots,x_{i_n} A$;
	\item[{\em(d)}] the $\LanQ$-\textbf{protoformulas} are only those expressions which can be built up by clauses {\em(a)--(c)}; \textbf{free occurrences} of bound variables in a protoformula are only those that can be determined by clauses {\em(a)--(c)}.
\end{itemize}

An $\LanQ$-formula expression $A$ is an $\LanQ$-\textbf{formula}\index{$\LanQ$-formula} if each bound variable occurring in $A$ is in the scope of some quantifier; in other words, a protoformula $A$ if a formula if it has no free occurrences of bound variables or, possibly, no bound variables at all.
\end{defn}

The set of all $\LanQ$-formulas is denoted by $\FmQ$.\\

Let $A(x_{i_1},\ldots,x_{i_n},u_{j_1},\ldots,u_{j_m})$ be a protoformula whose free occurrences of bound variables are $x_{i_1},\ldots,x_{i_n}$ and whose free variables are $u_{j_1},\ldots,u_{j_m}$. Then, denoting
$\overline{x}:=(x_{i_1},\ldots,x_{i_n})$ and $\overline{u}:=(u_{j_1},\ldots,u_{j_m})$, we occasionally write $A(\overline{x},\overline{u})$ or $A(x_{i_1},\ldots,x_{i_n},\overline{u})$ instead of $A(x_{i_1},\ldots,x_{i_n},u_{j_1},\ldots,u_{j_m})$.\\

Let $A(u_{i_1},\ldots,u_{i_n})$ be a formula containing free variables $u_{i_1},\ldots,u_{i_n}$ (pairwise distinct), written also as $A(\overline{u})$. Also, let $x_{i_1},\ldots,x_{i_n}$ be $n$ bound variables which do not occur in
$A(\overline{u})$. The expression $A(u_{i_1}\backslash x_{i_1},\ldots,u_{i_n}\backslash x_{i_n})$, written also as $A(\overline{u}\backslash\overline{x})$ or as
$A(x_{i_1},\ldots,x_{i_n})$ or as $A(\overline{x})$, is a protoformula which is obtained by the simultaneous replacement of all occurrences of each $u_{i_k}$
with $x_{i_k}$, respectively. On the other hand, if $\overline{t}=(t_{1},\ldots,t_n)$ are $\LanT$-terms (not necessarily pairwise distinct), we denote by $A(x_{i_1}\backslash t_1,\ldots,x_{i_n}\backslash t_n)$, or by $A(t_1,\ldots,t_n)$ or by $A(\overline{t})$, the formula obtained by the simultaneous replacement of each occurrence of $x_{i_k}$ in $A(x_{i_1},\ldots,x_{i_n})$ with $t_k$. The same resulting formula,
$A(u_{i_1}\backslash t_1,\ldots,u_{i_n}\backslash t_n)$, written also as $A(\overline{u}\backslash\overline{t})$, is obtained by substitution of each term $t_k$ for all occurrences of the free variable $u_{i_k}$ in the formula $A(u_{i_1},\ldots,u_{i_n})$. 

Thus, if $\sigma$ is an $\LanT$-substitution such that $\sigma(u_{i_k})=t_k$, then we have:
\[
A(\sigma(u_{i_1}),\ldots,\sigma(u_{i_n}))=A(t_1,\ldots,t_n)).
\]

This induces the following definition:
\begin{equation}\label{E:formula-substitution}
\sigma(A(u_{i_1},\ldots,u_{i_n})):=A(\sigma(u_{i_1}),\ldots,\sigma(u_{i_n})).
\end{equation}

It is obvious that for the identity substitution $\iota$, we have:
\[
\iota(A)=A.
\]

If $X$ is a set of $\LanQ$-formulas and $\sigma$ is an $\LanT$-substitution,
we denote:
\[
\sigma(X):=\set{\sigma(A)}{A\in X}.
\]

Similar to the definition of the term algebra $\TmAl$ (or the formula algebra $\FormAl$ from Section~\ref{section:semantics}), we define the $\LanQ$-\textit{\textbf{formula algebra}}\index{algebra!formula} $\FmAlQ$ as follows:\index{$\FmAlQ$}
\[
\FmAlQ:=\langle\FmQ;\FuncF\rangle.
\]

In the sequel, we will expand $\FmAlQ$ by adding infinite operations corresponding to quantifiers. Also, in the sequel, along with $\FmAlQ$, we will use its unital expansion.

\section{Semantics of $\LanQ$}\label{section:Q-semantics}

An algebra $\langle\textsf{A};\FuncT\rangle$ is called a \textit{\textbf{term structure}},\index{term structure} or $t$-\textit{\textbf{structure}}\index{t-structure} for short. In particular, $\TmAl=\langle\Tm,\FuncT\rangle$ is a $t$-structure.

Given a $t$-structure $\fA=\langle\textsf{A};\FuncT\rangle$, a map
$v:\ParT\longrightarrow\textsf{A}$ which we extend to a homomorphism (denoted by the same letter) $v:\TmAl\longrightarrow\langle\textsf{A};\FuncT\rangle$ is called a \textit{\textbf{term valuation}}\index{term valuation}
(or $t$-\textit{\textbf{valuation}} for short) in the algebra $\langle\textsf{A};\FuncT\rangle$. For each term $t\in\Tm$,
we denote by  $v[t]$ the value of $t$ with respect to $v$. The set of the values of all $\LanT$-terms is denoted by $v[\Tm]$. 

We note that any $t$-valuation in $\TmAl$ is an $\LanT$-substitution, and vice versa. Also,
given an $\LanT$-substitution $\sigma$ and a $t$-valuation in a term structure $\fA$, $v\circ\sigma$ is also a $t$-valuation in $\fA$.

In the sequel, we will occasionally use the following notation:
\[
v[\overline{t}]:=(v[t_1],\ldots,v[t_n]),
\]
where $\overline{t}=(t_1,\ldots,t_n)$ is a list of $\LanT$-terms.

Let $v$ be a $t$-valuation in a $t$-structure $\fA=\langle\textsf{A};\FuncT\rangle$ and let $\overline{u}=(u_{i_1},\ldots,u_{i_n})$, where all parameters of $\overline{u}$ are pairwise distinct. Also, let $\overline{a}=(a_1,\ldots,a_n)\in\textsf{A}^{n}$. We will often use the following $t$-valuation related to $v$, $\overline{u}$ and $\overline{a}$:
\[
w_{(v,\overline{u}\backslash\overline{a})}[u^{\prime}]:=
\begin{cases}
\begin{array}{cl}
a_k &\text{if $u^{\prime}=u_{i_k}$}\\
v[u^{\prime}] &\text{otherwise}.
\end{array}
\end{cases}
\]

In addition, given $\overline{u}=(u_{i_1},\ldots,u_{i_n})$, where the parameters $u_{i_k}$ are pairwise distinct, and $\overline{t}=(t_1,\ldots,t_n)$ is a list of $\LanT$-terms, we will sometimes use a $t$-valuation:
\[
w_{(v,\overline{u}\backslash\overline{t})}[u^\prime]:=\begin{cases}
	\begin{array}{cl}
		v[t_k] &\text{if $u^{\prime}=u_{i_k}$}\\
		v[u^\prime] &\text{otherwise}.
	\end{array}
\end{cases}
\]

If $\overline{u}=(u)$ and $\overline{t}=(t)$, we simply write $w_{(v,u\backslash t)}$.

Also, we define a \textit{\textbf{restricted $\bm{t}$-valuation}}\index{restricted $\bm{t}$-valuation}} as follows.

Assume that all parameters of a term $t$ are among $u_{i_1},\ldots,u_{i_m}$ and $\overline{a}=(a_{i_1},\ldots,a_{i_m})$ where each $a_{i_k}\in\textsf{A}$. Then we define $w_{\overline{u}\backslash\overline{a}}[t]$ by induction on the complexity of $t$:
\begin{itemize}
	\item $w_{\overline{u}\backslash\overline{a}}[t]=a_{i_k}$ if $t=u_{i_k}$;
	\item $w_{\overline{u}\backslash\overline{a}}[t]=f(w_{\overline{u}\backslash\overline{a}}[t_1],\ldots,w_{\overline{u}\backslash\overline{a}}[t_n])$
	if $t=ft_{1}\ldots t_n$.
\end{itemize}

\begin{defn}[$\mathcal{Q}$-structure]\label{D:Q-structure}\index{$\mathcal{Q}$-structure}
	A $\mathcal{Q}$-\textbf{structure} of type $\LanQ$ is a system 
	{\em$\fA=\langle\langle\textsf{A};\FuncT\rangle,\langle\textsf{B};\FuncF\rangle,\PredF,\Phi_{\Q}\rangle$}, where {\em$\langle\textsf{A};\FuncT\rangle$} is a $t$-structure, {\em$\langle\textsf{B};\FuncF\rangle$} is an algebra of type $\LanF$; further,
	each $p\in\PredF$ with $\#(p)=n$ {\em($n\ge\! 1$)} is associated with a same-name, but bold, map
	{\em $\bm{p}:\textsf{A}^{n}\longrightarrow\textsf{B}$}, and  each quantifier $Q\in\mathcal{Q}$ with $\#(Q)=n$ {\em($n\ge\! 1$)} is associated with a map $\varphi_{Q}\in\Phi_{\Q}$ that maps each set
	{\em\[
		\set{\bm{q}_{(A(\overline{x}),v)}(a_1,\ldots,a_n)}{(a_1,\ldots,a_n)\in \textsf{A}^{n}} 
		\]}to an element of {\em\textsf{B}}, where $\bm{q}_{(A(\overline{x}),v)}$ is a map from {\em$\textsf{A}^n$} into {\em\textsf{B}}, for each protoformula $A(x_{i_1},\dots,x_{i_n})$ with free occurrences of exactly $n$ bound variables $\overline{x}=(x_{i_1},\dots,x_{i_n})$. {\em (See definition below.)}
\end{defn}

\begin{rem}
	{\em
		Definition~\ref{D:Q-structure} does not specify what $\varphi_Q$ is, but it indicates that every $\varphi_Q$ is an infinite operation which is defined on a subset of $\textsf{B}$ which depends on the values of $\bm{q}_{(A(\overline{x}),v)}$; the latter in turn depends on a protoformula $A(\overline{x})$ and a $t$-valuation $v$. Also, as we show below, the mappings $\bm{q}_{(A(\overline{x}),v)}$ depend not only on protoformulas $A(\overline{x})$, but also on the mapping $\bm{p}$.}
\end{rem}

Given a $\mathcal{Q}$-structure $\fA=\langle\langle\textsf{A};\FuncT\rangle,\langle\textsf{B};\FuncF\rangle,\PredF,\Phi_{\Q}\rangle$ and a (restricted or full) term valuation $v$ in $\langle\textsf{A};\FuncT\rangle$, we define a \textit{\textbf{formula valuation}}\index{formula valuation} (or $f$-\textit{\textbf{valuation}}\index{$f$-valuation} for short) of all $\LanQ$-formulas in the $\mathcal{Q}$-structure $\fA$ as a 
map $\hat{v}:\FmQ\longrightarrow\textsf{B}$. The definition of a formula valuation $\hat{v}$ is grounded on the complexity of each $\LanQ$-formula. 

The following notation will be used below. Let $A(\overline{x},\overline{u^\ast})$ be a protoformula whose free occurrences of bound variables are $\overline{x}=(x_{i_1},\ldots,x_{i_n})$ and whose free variables, if any, are $\overline{u^\ast}=(u_{j_1},\ldots,u_{j_m})$. We assume that the protoformula $A(\overline{x},\overline{u^\ast})$ is obtained from a formula $A(\overline{u},\overline{u^\ast})$, where $\overline{u}=(u_{i_1},\ldots,u_{i_n})$, by replacement of each $u_{i_k}$ with $x_{i_k}$. Thus all parameters of the last formula partitioned in two lists, $\overline{u}$ and $\overline{u^\ast}$, where the last list can be empty. 

Now we turn to the definition of a formula valuation. In the course of this definition, we introduce two more functions, 
$\bm{p}_{A(\overline{u},\overline{u^{\ast}})}:\textsf{A}^{n+m}\longrightarrow\textsf{B}$, where $\overline{u^{\ast}}$ can be void and thus $m=0$, 
and $\bm{q}_{(A(\overline{x},\overline{u^\ast}),v)}:\textsf{A}^{n}\longrightarrow\textsf{B}$.\\

If $A(\overline{u})=p(t_1,\ldots,t_k)$, where $p\in\mathcal{P}_{F}$ and $t_1,\ldots,t_k$ are terms, whose parameters are among $u_{i_1},\ldots,u_{i_n}$; then
\begin{itemize}
	\item $\hat{v}[A(\overline{u})]:=\bm{p}(v[t_1],\ldots,v[t_k])$;
	\item $\bm{p}_{A(\overline{u})}(a_1,\ldots,a_n):=\bm{p}(	w_{\overline{u}\backslash\overline{a}}[t_1],\ldots,	w_{\overline{u}\backslash\overline{a}}[t_k])$,
	where $w$ is a restricted valuation whose domain is $\{u_{i_1},\ldots,u_{i_n}\}$, $\overline{u}=(u_{i_1},\ldots,u_{i_n})$ and $\overline{a}=(a_1,\ldots,a_n)$.
\end{itemize}

Next, let $A(\overline{u})=FA_1\ldots A_k$, where $F\in \mathcal{F}_F$ and $A_1,\ldots,A_k$ are $\LanQ$-formulas. Then we define:
\begin{itemize}
	\item $\hat{v}[FA_1\ldots A_k]:=F(\hat{v}[A_1],\ldots,\hat{v}[A_k])$;
	\item $\bm{p}_{A(\overline{u})}(a_1,\ldots,a_n):=F(\bm{p}_{A_1(\overline{u})}(\overline{a}),\ldots,\bm{p}_{A_k(\overline{u})}(\overline{a}))$,
	where $\overline{a}=(a_1,\ldots,a_n)$.
	\end{itemize}

Further, let $A(\overline{u^\ast})=Qx_{i_1}\ldots x_{i_n}A(x_{i_1},\ldots,x_{i_n},\overline{u^\ast})$, where $Q$ is a quantifier of arity $n$ which precedes a protoformula $A(x_{i_1},\ldots,x_{i_n},\overline{u^\ast})$. We first obtain a formula $A(\overline{u},\overline{u^\ast})$, where $\overline{u}=(u_{i_1},\ldots,u_{i_n})$ and $\overline{u^\ast}=(u_{j_1},\ldots,u_{j_m})$, by a simultaneous replacement of each occurrence of $x_{i_k}$ with $u_{i_k}$, providing that the parameters of $\overline{u}$ do not occur in $A(\overline{u^\ast})$. Then, we successively define:
\begin{itemize}
	\item $\bm{q}_{(A(\overline{x},\overline{u^\ast}),v)}(a_1,\ldots,a_n):=\bm{p}_{A(\overline{u},\overline{u^\ast})}(a_1,\ldots,a_n,v[\overline{u^{\ast}}])$;
	in particular,
	\begin{equation}\label{E:q-and-p-at-v}
	\bm{q}_{(A(\overline{x},\overline{u^\ast}),v)}(v[t_1],\ldots,v[t_n])=\bm{p}_{A(\overline{u},\overline{u^\ast})}(v[t_1],\ldots,v[t_n],v[\overline{u^{\ast}}]),
	\end{equation}
	where $(t_1,\ldots,t_n)$ is a list of $\LanT$-terms which are not necessarily pairwise distinct;
	\item $\hat{v}[Q x_{i_1}\ldots x_{i_n} A(x_{i_1},\ldots,x_{i_n},\overline{u^\ast})]:=\varphi_{Q}(\set{\bm{q}_{(A(\overline{x},\overline{u^\ast}),v)}(\overline{a})}{\overline{a}\in\textsf{A}^n})$, where  each $\varphi_{Q}$ is an infinite (perhaps partial) operation which is defined on subsets of $\textsf{B}$ and depends on the quantifier $Q$;
	\item $\bm{p}_{Q x_{i_1}\ldots x_{i_n} A(x_{i_1}\ldots x_{i_n},\overline{u^\ast})}(b_1,\ldots,b_m):=\widehat{w_{\overline{u^\ast}\backslash\overline{b}}}[Q x_{i_1}\ldots x_{i_n} A(x_{i_1}\ldots x_{i_n},\overline{u^\ast})]$, where  $\overline{b}=(b_1,\ldots,b_m)\in\textsf{A}^m$.
\end{itemize}

Perhaps the last definition, consisting of several concepts, needs clarification.

Given a $\Q$ structure $\fA$, it should be clear that an $f$-valuation  $\hat{v}$ in $\fA$ depends on a $t$-valuation $v$ in $\lr{\textsf{A};\FuncT}$, which is converted to $\hat{v}$ through the maps $\bm{p}$ and $\bm{q}$. The former, indexed with an $\LanQ$-formula $A(\overline{u})$ which contains exactly free variables $u_{i_1},\ldots,u_{i_n}$ arranged as $\overline{u}=(u_{i_1},\ldots,u_{i_n})$, is a map $\bm{p}_{A(\overline{u})}:\textsf{A}^{n}\longrightarrow\textsf{B}$, which makes $A(\overline{u})$ a multiple-valued predicate of arity $n$ defined on \textsf{A} with its values in \textsf{B}.\footnote{See Lemma~\ref{L:valuation-of-formula} below.} If a formula $A$ has no free variables, $\bm{p}_A$ is an element of \textsf{B}. We have to note that $\bm{p}$ and $\bm{p}_{A}$, where a formula $A$ is with or without free variables, are not necessarily the same thing, since $\bm{p}$ being associated with a predicate symbol $p$ has arity $n\ge 1$, while $\bm{p}_{A}$ may have the arity equal to $0$. 

A map $\bm{q}_{(A(x,\overline{u^\ast}),v)}$, where $A(x,\overline{u^\ast})$ is a protoformula with free occurrences of exactly one bound variable $x$ and $v$ is a $t$-valuation, is needed to generate the set $\set{\bm{q}_{(A(x,\overline{u^\ast}),v)}(a)}{a\in\textsf{A}}\subseteq\textsf{B}$, to which each of infinite unary operations, $\varphi_{\forall}$ and $\varphi_{\exists}$, applies to obtain  in \textsf{B} the values $\hat{v}[\forall xA(x,\overline{u^\ast})]$ or $\hat{v}[\exists xA(x,\overline{u^\ast})]$, respectively.
\begin{rem}
{\em
Although $\bm{p}$ in general depends on the carriers \textsf{A} and \textsf{B}, we use the same notation, $\bm{p}$, for the except of the cases when $\textsf{A}=\Tm$ and \textsf{B} is $\FmQ$ or a quotient of the latter; in the first of the last cases, we use the notation $\bm{p}^\mathbf{L}$ (defined in Example~\ref{example:first} below), and in the second $\bm{p}^\mathbf{LT}$ (defined in Section~\ref{section:Lindenbaum-Tarski-Q-structures}). Accordingly,
we employ the notations: $\bm{p}^\mathbf{L}_{A(\overline{u})}$ and $\bm{q}^\mathbf{L}_{(A(\overline{x},\overline{u^\ast}),v)}$ in the first case and $\bm{p}^\mathbf{LT}_{A(\overline{u})}$, and $\bm{q}^\mathbf{LT}_{(A(\overline{x},\overline{u^\ast}),v)}$ in the second case.
}
\end{rem}

Now we consider some properties pertaining to formula valuations.
The next lemma is a refinement of a definition above.
\begin{lem}\label{L:lemma-on-p-function}
Let $(A(\overline{u},\overline{u^{\ast}})$ be a formula whose parameters are partitioned in a list $\overline{u}=(u_{i_1},\ldots,u_{i_n})$ and a list $($maybe empty$)$ $\overline{u^{\ast}}=(u_{j_1},\ldots,u_{j_m})$. Let $\overline{a}=(a_1,\ldots, a_n)$ be an element of {\em$\textsf{A}^n$}. Then
\[
\bm{p}_{A(\overline{u},\overline{u^{\ast}})}(a_1,\ldots,a_n,v[\overline{u^{\ast}}])=\bm{p}_{(A(\overline{u},\overline{u^{\ast}})}(w_{(v,\overline{u}\backslash\overline{a})}[\overline{u}],w_{(v,\overline{u}\backslash\overline{a})}[\overline{u^{\ast}}]).
\]
\end{lem}
\noindent\emph{Proof}~is left to the reader. (Exercise~\ref{section:Q-semantics}.\ref{EX:lemma-on-p-function})

\begin{lem}\label{L:valuation-of-formula}
Let $A(\overline{u})$ be a formula, whose free variables are $\overline{u}=(u_{i_1},\ldots,u_{i_n})$. Then for any $t$-valuation $v$ in an arbitrary $\Q$-structure,
$$\hat{v}[A(\overline{u})]=\bm{p}_{A(\overline{u})}(v[\overline{u}]).$$
\end{lem}
\begin{proof}
We use induction on the complexity of $A(\overline{u})$. We will consider the basis case, leaving the rest of the proof to the reader. (Exercise~\ref{section:Q-semantics}.\ref{EX:valuation-of-formula})

Assume that $A(\overline{u})=p(t_{1}(\overline{u}),\ldots,t_{k}(\overline{u}))$, where $t_1,\ldots,t_k$ are $\LanT$-terms whose (free) variables are included in $\overline{u}$. By definition,
\[
\hat{v}[A(\overline{u})]=\bm{p}(v[t_1],\ldots,v[t_k])=\bm{p}_{(A(\overline{u}),v)}(v[t_1],\ldots,v[t_k]).
\]
\end{proof}

\begin{lem}\label{L:formula-as-predicate}
	Let $A(\overline{u})$ be a formula, whose free variables are $\overline{u}=(u_{i_1},\ldots,u_{i_n})$, and let $\sigma$ be an $\LanT$-substitution. Then for any $t$-valuation $v$ in any $\Q$-structure, 
	\begin{align*}
	\widehat{v\circ\sigma}[A(\overline{u})]=
	\hat{v}[\sigma(A(\overline{u}))]=\hat{v}[A(\sigma(\overline{u})].
	\end{align*}
\end{lem}
\begin{proof}
	The proof of the statement can be carried out by induction on the complexity of $A(u_{i_1},\ldots,u_{i_n})$. We show the basis and leave the rest of the proof to the reader. (Exercise~\ref{section:Q-semantics}.\ref{EX:formula-as-predicate})
	
	Indeed, we successively obtain:
	\[
	\begin{array}{rl}
	\widehat{v\circ\sigma}[p(t_1,\ldots,t_n)] & =
	~\bm{p}(v\circ\sigma[t_1],\ldots,v\circ\sigma[t_n])\\
	&=~ \bm{p}(v[\sigma(t_1)],\ldots,v[\sigma(t_n)])\\
	&= ~\hat{v}[p(\sigma(t_1),\ldots,\sigma(t_n))]\\
	&= ~\hat{v}[\sigma(p(t_1,\ldots,t_n))].
	\end{array}
	\]
\end{proof}

\begin{defn}[logical $\mathcal{Q}$-matrix]\index{$\mathcal{Q}$-matrix}
	A \textbf{logical} $\mathcal{Q}$-\textbf{matrix} $($or simply $\mathcal{Q}$-\textbf{matrix}$)$ of type $\LanQ$ is a system $\langle\fA,D\rangle$, where {\em$\fA=\langle\langle\textsf{A};\FuncT\rangle,\langle\textsf{B};\FuncF\rangle,\PredF,\Phi_{\Q},\rangle$} is a $\mathcal{Q}$-structure of type $\LanQ$ and $D$, called a \textbf{logical filter}, is a nonempty subset of {\em\textsf{B}}. If $D$ consists of a single element, we call such a $\Q$-matrix \textbf{unital}. Denoting the element of $D$ by $\one$ $($perhaps with a subscript$)$ in a unital $\Q$-matrices $\langle\fA,D\rangle$, we then replace {\em$\langle\textsf{B};\FuncF\rangle$} in $\fA$ with its unital expansion {\em$\langle\textsf{B};\FuncF,\one\rangle$} and treat the $\Q$-structure 
	$\langle\langle\textsf{A};\FuncT\rangle,\langle\textsf{B};\FuncF,\one\rangle,\PredF,\Phi_{\Q},\rangle$ as a $\Q$-matrix.
\end{defn}

Now we will consider three examples. The first example is the most important because it will later be developed into the concept of the Lindenbaum $\Q$-matrix (Definition~\ref{D:Lindenbaum-Q-matrix}) and then into the concept of the Lindenbaum-Tarski $\Q$-structure (Definition~\ref{D:LT-Q-structure}). The importance of this example is also highlighted in Proposition~\ref{P:semantics-for-formula-algebra}.

\begin{example}\label{example:first}
{\em
As our first example, we consider a $\Q$-structure $\langle\TmAl,\FmAlQ,\mathcal{P}_F,\Phi_{\Q}\rangle$, where  by definition, for any $p\in\mathcal{P}$ with $\#(p)=n$, $\bm{p}^{\mathbf{L}}:(t_1,\ldots,t_n)\mapsto p(t_1,\ldots,t_n)$. Since any $t$-valuation in this structure is an $\LanT$-substitution, we fix an $\LanT$-substitution $\sigma$. Now we want to illustrate what $\hat{\sigma}[A(\overline{u})]$ is in this $\Q$-structure under the valuation $\sigma$.

Assume that $A(\overline{u})=p(\overline{t})$ with $\#(p)=k$, where $\overline{u}=(u_{i_1},\ldots,u_{i_n})$ and $\overline{t}=(t_1,\ldots,t_k)$. Then we have:
\begin{itemize}
	\item $\hat{\sigma}[p(\overline{t})]=p(\sigma(\overline{t}))$;
	\item$\bm{p}^{\mathbf{L}}_{A(\overline{u})}(s_1,\ldots,s_n)=
	p(\delta_{\overline{u}\backslash\overline{s}}(t_1),\ldots,\delta_{\overline{u}\backslash\overline{s}}(t_k))=\delta_{\overline{u}\backslash\overline{s}}(p(t_1,\ldots,t_k))$,\\ where
	$\overline{s}=(s_1,\ldots,s_n)\in\Tm^n$.
\end{itemize}

Next, let $A(\overline{u})=FA_{1}\ldots A_k$. Then we have:
\begin{itemize}
	\item $\hat{\sigma}[A(\overline{u})]=F(\hat{\sigma}[A_{1}],\ldots,\hat{\sigma}[A_{k}])$;
	\item $\bm{p}^{\mathbf{L}}_{A(\overline{u})}(s_1,\ldots,s_n)=F(\bm{p}^{\mathbf{L}}_{A_1(\overline{u})}(\overline{s}),\ldots,\bm{p}^{\mathbf{L}}_{A_k(\overline{u})}(\overline{s}))$;
\end{itemize}

Now, let $A(\overline{u^\ast})=Qx_{i_1}\ldots x_{i_n}A(x_{i_1},\ldots,x_{i_n},\overline{u^\ast})$ and $A(\overline{u},\overline{u^\ast})$ be a formula, where $\overline{u}=(u_{i_1},\ldots,u_{i_n})$ and $\overline{u^\ast}=(u_{j_1},\ldots,u_{j_m})$, such that neither of $u_{i_k}$ occurs in the protoformula $A(x_{i_1},\ldots,x_{i_n},\overline{u^\ast})$ and all $u_{i_k}$ are pairwise distinct. Then, defining
\[
\bm{q}^{\mathbf{L}}_{(A(\overline{x},\overline{u^\ast}),\sigma)}(s_1,\ldots,s_n):=\bm{p}^{\mathbf{L}}_{A(\overline{u},\overline{u^\ast})}(s_1,\ldots,s_n,\sigma(\overline{u^{\ast}}))
\]
and 
\begin{align*}
	\varphi_{Q}(\set{\bm{q}^{\mathbf{L}}_{(A(\overline{x},\overline{u^\ast}),\sigma)}(\overline{s})}{\overline{s}\in\Tm^n}):=\sigma(Q x_{i_1}\ldots x_{i_n}A(x_{i_1},\ldots,x_{i_n},\overline{u^\ast})\\
	=Q x_{i_1}\ldots x_{i_n} A(x_{i_1},\ldots,x_{i_n},\sigma(\overline{u^{\ast}})),
\end{align*}
where $\varphi_{Q}\in\Phi_{\Q}$, we successively obtain:
\begin{itemize}
	\item $
	\hat{\sigma}[Qx_{i_1}\ldots x_{i_n}A(x_{i_1},\ldots,x_{i_n},\overline{u^\ast})]= \sigma(Q x_{i_1}\ldots x_{i_n}A(x_{i_1},\ldots,x_{i_n},\overline{u^\ast});
	$	
	\item $\bm{p}^{\mathbf{L}}_{Q\overline{x} A(\overline{x},\overline{u^\ast})}(t_1,\ldots,t_m):=\widehat{\delta_{\overline{u^\ast}\backslash\overline{t}}}[Q\overline{x} A(\overline{x},\overline{u^\ast})]=\delta_{\overline{u^\ast}\backslash\overline{t}}(Q\overline{x} A(\overline{x},\overline{u^\ast}))$.
\end{itemize}
}
\end{example}

Summarizing the definitions and their simplifications of the last example, we obtain the following.
\begin{prop}\label{P:semantics-for-formula-algebra}
	Let $\sigma$ be an $\LanT$-substitution, $A(\overline{u})$ and $B(\overline{u},\overline{u^\ast})$ be formulas such that $\overline{u}=(u_{i_1},\ldots,u_{i_n})$ is the list of all parameters of the former and the parameters of the latter are partitioned in two lists ---
	$\overline{u}$ and $\overline{u^\ast}=(u_{j_1},\ldots,u_{j_m})$ Then the following hold in {\em$\langle\TmAl,\FmAlQ,\mathcal{P}_F,\Phi_{\Q}\rangle$}:
	{\em\[
		\begin{array}{cl}
		(\text{a}) &\hat{\sigma}[A(\overline{u})] = \sigma(A(\overline{u}));\\
		(\text{b}) &\bm{p}^{\mathbf{L}}_{B(\overline{u},\overline{u^{\ast}})}(s_1,\ldots,s_n,\sigma(\overline{u^{\ast}}))=
		\delta_{(\sigma,\overline{u}\backslash\overline{s})}(B(\overline{u},\overline{u^{\ast}}));\\
		(\text{c}) &\bm{q}^{\mathbf{L}}_{(B(\overline{x},\overline{u^\ast}),\sigma)}(s_1,\ldots,s_n)=
		\delta_{(\sigma,\overline{u}\backslash\overline{s})}(B(\overline{u},\overline{u^\ast})).
		\end{array}
		\]}
\end{prop}
\begin{proof}
We leave for the reader to check the statement. (Exercise~\ref{section:Q-semantics}.\ref{EX:semantics-for-formula-algebra})
\end{proof} 

\begin{example}\label{example:ordinary-quantifiers}
	{\em
		Suppose an $\LanQ$-language contains a nonempty set $\FuncT$ of function symbols, a nonempty set $\PredF$ of predicate symbols, the binary logical connectives: $\land,\lor,\rightarrow,\neg$, and two unary quantifiers, $\forall$ and $\exists$.
		Consider a $\Q$-structure $\fA=\langle\langle\textsf{A};\FuncT\rangle,\langle\lbrace\zero,\one\rbrace;\land,\lor,\rightarrow,\neg\rangle,\PredF,\lbrace\varphi_{\forall},\varphi_{\exists}\rbrace\rangle$, where $\lr{\lbrace\zero,\one\rbrace;\land,\lor,\rightarrow}$ is a 2-element Boolean algebra. Further, given a $t$-valuation $v$ and a protoformula $A(x,\overline{u^{\ast}})$ containing free occurrences of only $x$, we define
		\[
		\varphi_{\forall}(\set{\bm{q}_{(A(x,\overline{u^{\ast}}),v)}}{a\in\textsf{A}}):=\begin{cases}
		\begin{array}{cl}
	\one &\text{if $\set{\bm{q}_{(A(x,\overline{u^{\ast}}),v)}}{a\in\textsf{A}}=\textsf{A}$}\\
	\zero &\text{otherwise},
		\end{array}
		\end{cases}
		\]
		and
			\[
		\varphi_{\exists}(\set{\bm{q}_{(A(x,\overline{u^{\ast}}),v)}(a)}{a\in\textsf{A}}):=\begin{cases}
		\begin{array}{cl}
		\one &\text{if $\set{\bm{q}_{(A(x,\overline{u^{\ast}}),v)}(a)}{a\in\textsf{A}}\neq\varnothing$}\\
		\zero &\text{otherwise}.
		\end{array}
		\end{cases}
		\]
		
		As a $\Q$-matrix, we use the unital expansion
		$\langle\lbrace\zero,\one\rbrace;\land,\lor,\rightarrow,\one\rangle$.
			
		It should be clear that the above maps $\varphi_{\forall}$ and $\varphi_{\exists}$ define the ordinary interpretation of the universal and existential quantifier, respectively, in an $\LanQ$-language exemplifying that of firs-order logic.
	}
\end{example}

The next example may seem artificial.\footnote{See a remark in~\cite{rosser-turquette1951}, p. 51: ``To the best of our knowledge, no one has so far found need for quantifiers with either $\beta_i>1$ or $\gamma_i>1$. However, such quantifiers are distinctly possible $\ldots$''}
\begin{example}
	{\em 
Suppose an $\LanQ$-language, besides the mandatory categories, contains only one predicate symbol $p$ of arity 3 and one binary quantifier $Q$. 
Let $\fA=\langle\langle\mathbb{Z}^+;\varnothing\rangle,\langle\mathbb{Z}^+;\varnothing\rangle,\lbrace p\rbrace,\lbrace\varphi_{Q}\rbrace\rangle$, where $\mathbb{Z}^+$ is the set of positive integers, is a $\mathcal{Q}$-structure for this language, where the predicate symbol $p$ is interpreted as the map	 $\bm{p}:(a,b,c)\mapsto~\gcd(a,b,c)$, where `$\gcd(a,b,c)$' stands for the greatest common divisor of $a,b,c$, and $\varphi_{Q}(X)$, if it is defined, is the maximal element of $X\subseteq\mathbb{Z}^+$. It should be clear that the operation $\varphi_{Q}$ is partial. However, since the set
\[
\set{\bm{q}_{(p(x,y,u),v)}(a,b)}{(a,b)\in\mathbb{Z}^{+}\times\mathbb{Z}^{+}}
=\set{\bm{p}(a,b,v[u])}{(a,b)\in\mathbb{Z}^{+}\times\mathbb{Z}^{+}}
\]
is always finite, for any $t$-valuation $v$,
\[
\hat{v}[Qxy\,p(x,y,u)]=\varphi_{Q}(\set{\bm{q}_{(p(x,y,u),v)}(a,b)}{(a,b)\in\mathbb{Z}^{+}\times\mathbb{Z}^{+}})=v[u].
\]	
}
\end{example}

More examples will be given in Section~\ref{section:Q-logics-three-examples}.\\

The vocabulary and contents of the next definition are similar to those of Definition~\ref{D:validity-1}.
\begin{defn}\label{D:Q-matrix}\index{Q-matrix!valuation}
	Let $\fM=\langle\fA,D\rangle$ be a $\mathcal{Q}$-matrix of type $\LanQ$. An $\LanQ$-formula $A$ is \textbf{satisfied} by a term valuation $v$ if $\hat{v}[A]\in D$, in which case we say that $v$ \textbf{satisfies} $A$ in $\fM$, symbolically $\fM\models_{v}A$. Further, a formula $A$ is \textbf{semantically valid} $($or simply \textbf{valid}$)$ in $\fM$, symbolically $\fM\models A$ or $\models_{\mathfrak{M}}A$, if $\hat{v}[A]\in D$, for any term valuation $v$ in $\fM$. If a formula $A$ is not valid in a $\mathcal{Q}$-matrix, it is called \textbf{invalid} in this matrix or is \textbf{rejected} by the matrix or one can say that this matrix \textbf{rejects}\index{Q-matrix!rejecting} the formula. If $\hat{v}[A]\notin D$, $v$ is a  \textbf{valuation rejecting} or \textbf{Q-matrix!refuting} $A$. \index{Q-matrix!valuation!refuting}
\end{defn}

Referring to a $t$-\textit{\textbf{valuation}} in a $\Q$-matrix $\fM=\langle\fA,D\rangle$, we mean a valuation in the term structure of $\fA$.

Given a valuation $v$ in a $\mathcal{Q}$-matrix $\fM$ and a set $X\subseteq\FmQ$, we denote
\[
\fM\models_{v}X
\]
if $\fM\models_{v}A$ for each $A\in X$. 

Also, we denote:
\[
X\models_{\mathfrak{M}}A~\define~(\fM\models_{v}X\Longrightarrow\fM\models_{v}A~\text{for every $t$-valuation $v$ in $\fM$}).
\]

\paragraph{Exercises~\ref{section:Q-semantics}}
\begin{enumerate}
	\item \label{EX:valuation-formula-without-free-variables} Prove that if a formula $A$ contains no free variables, then for any $t$-valuations $v_1$ and $v_2$, $\widehat{v_1}[A]=\widehat{v_2}[A]$.
	\item\label{EX:lemma-on-p-function}Prove Lemma~\ref{L:lemma-on-p-function}.
	\item\label{EX:valuation-of-formula}Complete the proof of Lemma~\ref{L:valuation-of-formula}.
	\item\label{EX:formula-as-predicate}Complete the proof of Lemma~\ref{L:formula-as-predicate}.
	\item Let $A(u_{i_1},\ldots,u_{i_n})$ be any formula whose all free variables are in the list $\overline{u}=(u_{i_1},\ldots,u_{i_n})$. And let $\overline{t}=(t_1,\ldots,t_n)$ and $\overline{t^\prime}=(t^{\prime}_1,\ldots,t^{\prime}_n)$ be lists of terms not all necessarily distinct. Prove that for any $t$-valuation $v$, if $v[t_k]=v[t^{\prime}_k]$, $1\le k\le n$, then $\widehat{v\circ\sigma_{\overline{u}\backslash\overline{t}}}[A]=\widehat{v\circ\sigma_{\overline{u}\backslash\overline{t^{\prime}}}}[A]$.
	\item\label{EX:semantics-for-formula-algebra}Check the statement of Proposition~\ref{P:semantics-for-formula-algebra}.
\end{enumerate}

\section{$\Q$-Consequence}\label{section:Q-consequence}
\begin{defn}[$\Q$-consequence]\label{D:Q-consequence}
A relation $\vdash_{\Q}$ on {\em$\mathcal{P}(\FmQ)\times\FmQ$} satisfying the conditions {\em(a)--(c)} of Definition~\ref{D:consequnce-relation-single} is called a $\Q$-consequence.\index{Q-consequence}
The consequence operator {\em$\textbf{Cn}_{\Q}:\mathcal{P}(\FmQ)\longrightarrow\FmQ$} corresponding to $\vdash_{\Q}$ is defined as follows:
{\em\[
A\in\textbf{Cn}_{\Q}(X)~\define~X\vdash_{\Q}A,
\]}
for any set {\em$X\cup\lbrace A\rbrace\subseteq\FmQ$}. A $\Q$-consequence is called \textbf{structural} if for any term substitution $\sigma$ and any set {\em$X\cup\lbrace A\rbrace\subseteq\FmQ$},
\[
X\vdash_{\Q}A~\Longrightarrow~\sigma(X)\vdash_{\Q}\sigma(A),
\]
or in terms of operator {\em$\textbf{Cn}_{\Q}$},
{\em\[
A\in\textbf{Cn}_{\Q}(X)~\Longrightarrow~\sigma(A)\in\textbf{Cn}_{\Q}(\sigma(X)).
\]}
\end{defn}

To refer to $\vdash_{\Q}$ or to $\textbf{Cn}_{\Q}$  we use the term  \textit{\textbf{abstract logic}} $\aLog_{\Q}$\index{$\aLog_{\Q}$}. Abstract logics $\aLog_{\Q}$ will also be called $\Q$-\textit{\textbf{logics}}.\index{Q-logic}

Given an abstract logic $\aLog_{\Q}$, a set $X$ of $\LanQ$-formulas satisfying the property $X=\ConSQ(X)$ is called an $\aLog_{\Q}$-\textit{\textbf{theory}}.\index{$\aLog_{\Q}$-theory}
The set of all $\aLog_{\Q}$-theories is denoted by $\Sigma_{\aLog_{\Q}}$.

The formulas of the theory
\[
\bm{T}_{\aLog_{\Q}}:=\ConSQ(\varnothing)
\]
are called $\aLog_{\Q}$-\textbf{\textit{theses}} or $\aLog_{\Q}$-\textbf{\textit{theorems}}.\index{$\aLog_{\Q}$-tehorem}\\

The concept of implicative logics (Section~\ref{section:subclasses-unital}) can be extended to abstract logics $\aLog_{\Q}$.
\begin{defn}[implicative $\Q$-logic]\label{D:implicative-Q-logics}\index{Q-logic!implicative}
	A structural logic $\aLog_{\Q}$ is called an \textbf{implicative $\Q$-logic} if $\LanQ$ contains a binary connective $ \rightarrow$ such that for any $\aLog_{\Q}$-theory $D$, the following conditions are satisfied:
\begin{itemize}
	\item[\emph{(a)}]  $A\rightarrow A\in D$;
	\item[\emph{(b)}] $A\rightarrow B\in D$ and $B\rightarrow C\in D$ imply $A\rightarrow C\in D$;
	\item[\emph{(c)}] $A\in D$ and $A\rightarrow B\in D$ imply $B\in D$;
	\item[\emph{(d)}] $B\in D$ implies $A\rightarrow B\in D$;
	\item[\emph{(e)}] if $A_i\rightarrow B_i\in D$, 
	$B_i\rightarrow A_i\in D$, where $1\le i\le n$, and $F$ is any $n$-ary connective from $\LanF$, then $FA_1\ldots A_n\rightarrow FB_1\ldots B_n\in D$.
\end{itemize}
\end{defn}

Given an abstract $\Q$-logic $\aLog_{\Q}$, for any $\aLog_{\Q}$-theory $D=\ConSQ(X)$, we define:
\[
(A,B)\in\theta_{\aLog_{\Q}}(X)~\define~A\rightarrow B~\text{and}~B\rightarrow A\in D.
\]

We denote:
\[
\theta_{\aLog_{\Q}}:=\theta_{\aLog_{\Q}}(\varnothing).
\]

Employing the argument used in the proof of Proposition~\ref{P:some-classes-unital-logics-1}, we obtain:
\begin{prop}\label{P:congruence-property-implicative-logics}
	Let {\em$D=\ConSQ(X)$} be a theory of an implicative $\Q$-logic $\aLog_{\Q}$. Then $\theta_{\aLog_{\Q}}(X)$ is a congruence on $\FmAlQ$
	and $D$ is a congruence class with respect to $\theta_{\aLog_{\Q}}(X)$.
\end{prop}

\begin{rem}\label{remark:unital-Q-logics}
{\em
From the point of view of Definition~\ref{D:unital-logic}, the implicative
$\Q$-logics are \textit{\textbf{unital}}.\index{unital $\Q$-logic} We will not develop this concept for all structural $\Q$-logics, but we will use the property indicated in Proposition~\ref{P:congruence-property-implicative-logics} for particular ones; see Section~\ref{section:Q-logics-three-examples}.}
\end{rem}

\subsection{$\Q$-consequence defined by $\Q$-matrices}
Let $\mathcal{M}^{\Q}$ be a (nonempty) set of $\mathcal{Q}$-matrices of type $\LanQ$. For any set $X\cup\lbrace A\rbrace\subseteq\FmQ$, we define:
\begin{equation}\label{E:vdash_M^Q}
\begin{array}{r}
X\models_{\mathcal{M}^{\Q}}A~\define~\text{for any $\fM\in\mathcal{M}^{\Q}$ and any valuation $v$ in $\fM$},\\
\text{$\fM\models_{v}X$ implies $\fM\models_{v}A$}
\end{array}
\end{equation}

\begin{prop}\label{P:MQ-relation-is-structural-Q-con}
	The relation $\models_{\mathcal{M}^{\Q}}$ defined in~\eqref{E:vdash_M^Q} is a structural $\Q$-consequence.
\end{prop}
\begin{proof}
That the relation  is a consequence relation can be obtained by applying Proposition~\ref{P:semantic-consequence-generalized} to $\LanQ$-formulas
and $\Q$-matrices. Or this can be done directly; this last task is left to the reader. (Exercise~\ref{section:Q-consequence}.\ref{EX:MQ-relation-is-structural-Q-con-1}) We prove below that $\models_{\mathcal{M}^{Q}}$ is structural.

Let $\fA=\langle\langle\textsf{A};\FuncT\rangle,\langle\textsf{B};\FuncF\rangle,\PredF,\mathcal{Q},\rangle$ be a $\mathcal{Q}$-structure and $\langle\fA,D\rangle$ be an arbitrary $\mathcal{Q}$-matrix of $\mathcal{M}^{\Q}$. Also, assume that $v$ is an arbitrary $t$-valuation and $\sigma$ is an arbitrary $\LanT$-substitution. Since the former is a homomorphism $v:\TmAl\longrightarrow
\langle\textsf{A},\FuncT\rangle$ and the latter a homomorphism $\sigma:\TmAl\longrightarrow\TmAl$, the composition $v\circ \sigma:\TmAl\longrightarrow\langle\textsf{A},\FuncT\rangle$ is also a homomorphism; that is, the map $v\circ\sigma$ is a term valuation in $\langle\textsf{A},\FuncT\rangle$.

Let us denote
\[
v_1:=v\circ\sigma.
\]

We recall that for any $\LanT$-term $t$,
\[
v_1[t]=v[\sigma(t)]. \tag{$\ast$}
\]
(This can be proved by induction on the construction of term $t$; see Exercise~~\ref{section:Q-consequence}.\ref{EX:Q-consequence-is-consequence-1}.)

In virtue of Lemma~\ref{L:formula-as-predicate}, the property ($\ast$) further implies that for any $\LanQ$-formula $A$,
\[
\widehat{v_1}[A]=\hat{v}[\sigma(A)] \tag{$\ast\ast$}
\]

Then the equality ($\ast\ast$) obviously implies that for any $\LanQ$-formula $A$,
\[
\widehat{v_1}[A]\in D~\Longleftrightarrow~\hat{v}[\sigma(A)]\in D.
\]

The last equivalence in turn implies that the relation $\models_{\mathcal{M}^{\Q}}$ is structural.
\end{proof}

We call the consequence relation defined in~\eqref{E:vdash_M^Q} $\mathcal{M}^{\Q}
$-\textit{\textbf{consequence}}, for a given class $\mathcal{M}^{\Q}$ of $\mathcal{Q}$-matrices. 

Given a $\Q$-matrix $\fM$ and an abstract logic $\aLog_{\Q}$, $\fM$ is an $\aLog_{\Q}$-\textit{\textbf{model}} \index{$\aLog_{\Q}$-model}if for any set $X\cup\lbrace A\rbrace\subseteq\FmQ$,
\[
X\vdash_{\aLog_{\Q}}A~\Longrightarrow~X\models_{\mathfrak{M}}A.
\]

\subsection{$\mathcal{Q}$-consequence defined via inference rules}\label{section:Q-con-via-inference-rules}
In this subsection we limit ourselves with a formal language  $\Lan_{\forall\exists}$ and its fragment $\Lan_{(\forall\exists)^+}$; the former has already been employed in Example~\ref{example:ordinary-quantifiers}.
We recall that the set of the logical connectives of $\Lan_{\forall\exists}$ consists of
 binary connectives $\land, \lor, \rightarrow$, and a unary connective $\neg$; the language $\Lan_{(\forall\exists)^+}$ contains only the first three binary connectives. Also, both languages have only two unary quantifiers:
$\forall$ and $\exists$. Thus, $\Lan_{\forall\exists}$ is an extension of $\Lan_{(\forall\exists)^+}$. Since we will consider the formulas in these languages together, it will be convenient to call the formulas of both languages $\Lan_{\forall\exists}$-\textit{\textbf{formulas}}. With the last remark in mind, we denote by $\FmAE$ the set of all $\Lan_{\forall\exists}$-formulas.
 Unspecified $\LanAE$-formulas will be denoted by letters $A,B,C$ (with indices or without).

We  will consider three $\mathcal{Q}$-consequence relations given by three different sets of axiom schemata. To do this, we reproduce the axiom schemata that first appeared in Section~\ref{section:inference-rules} and Section~\ref{section:modus-rules-vs-rules}, provided that the metavariables are replaced with arbitrary $\LanAE$-formulas $A, B$ and $C$. To emphasize this change, we add the letter `Q' to the name of each schema below.
\[
\begin{array}{rl}
\text{ax1Q} &A\rightarrow(B\rightarrow A),\\
\text{ax2Q} &(A\rightarrow B)
\rightarrow((A\rightarrow(B\rightarrow C))\rightarrow
(A\rightarrow C)),\\
\text{ax3Q} &A\rightarrow(B\rightarrow(A\wedge B)),\\
\text{ax4Q} &(A\wedge B)\rightarrow A,\\
\text{ax5Q} &(A\wedge B)\rightarrow B,\\
\text{ax6Q} &A\rightarrow(A\vee B),\\
\text{ax7Q} &B\rightarrow(A\vee B),\\
\text{ax8Q} &(A\rightarrow C)\rightarrow((B\rightarrow C)
\rightarrow((A\vee B)\rightarrow C)),\\
\text{ax9Q} &(A\rightarrow B)
\rightarrow((A\rightarrow\neg B)\rightarrow\neg A), \\
\text{ax10Q} &\neg\neg A\rightarrow A,\\
\text{ax11Q} &B\rightarrow(\neg B\rightarrow A).
\end{array}
\]
To the above list, we add four new axiom schemata:
\begin{itemize}
	\item[ax12Q] $\forall x A(x)\rightarrow A(x\backslash t)$,\\
	where, given an arbitrary term $t$, the formula $A(x\backslash t)$ is obtained from the protoformula $A(x)$ by replacing of all free occurrences
	of $x$ with $t$;
	\item[ax13Q] $A(t)\rightarrow\exists x A(t\backslash x)$,\\
	where a protoformula $A(t\backslash x)$ is obtained from the formula $A(t)$ by replacing all occurrences of $t$ with a bound variable $x$ such that all new occurrences of $x$, if any, are free; this happens in particular if $x$ does not occur in $A(t)$;
	\item[ax14Q] $\forall x(A\rightarrow B(x))\rightarrow(A\rightarrow\forall x B(x))$, \\
	where $A$ is a formula that does not contain free occurrences of $x$
	and $B$ is a protoformula that contains free occurrences of $x$;
	\item[ax15Q] $\forall x(A(x)\rightarrow B)\rightarrow(\exists x A(x)\rightarrow B)$,\\
	where $A$ is a protoformula containing free occurrences of $x$
	and $B$ is any formula.
\end{itemize}

Each axiom schema becomes a corresponding axiom when concrete formulas are employed for $A, B$ and $C$.  

In addition to the premiseless rules ax1Q--ax15Q, we formulate two other modus rules. First of all, we have
\[
\dfrac{A, A\rightarrow B}{B},\tag{\textit{modus ponens}}
\]
as well as 
\begin{equation}\label{E:Generalization-rule}
\dfrac{A(u)}{\forall x A(u\backslash x)},
\end{equation}
where a protoformula $A(u\backslash x)$ is obtained from a formula $A(u)$ by replacing all occurrences of a free variable $u$ with a bound variable $x$ which does not occur freely in $A(u)$ and all new occurrences of $x$ in $A(u\backslash x)$ are free.\\

The concept of \textit{\textbf{derivation}} in calculi based on the axiom schemata ax1Q--ax15Q and rules modus ponens and~\eqref{E:Generalization-rule} is slightly different from what was discussed in Section~\ref{section:inference-rules}. We will consider three types of derivation, all of which are based on one and the same scheme. They differ only in what axiom schemata from the above list can be employed. We use a notation `$\CQ$' to denote one of the three calculi which we are going to deal with in this chapter:\index{$\CQ$}
\begin{itemize}
	\item a calculus $\PosQ$ is formulated in $\Lan_{(\forall\exists)^+}$ and is based on the use of the schemata ax1Q--ax8Q and ax12Q--ax15Q;
	\item a calculus $\ClQ$ is in the language $\LanAE$ and is based on the use of the schemata ax1Q--ax10Q and ax12Q--ax15Q;
	\item a calculus $\IntQ$\index{$\IntQ$} is also in $\LanAE$ and uses the schemata ax1--ax9Q, ax11Q and ax12Q--ax15Q.
\end{itemize}

The fact of derivability, that is, the existence of a derivation of the corresponding type will be denoted by a binary relation $\Vdash_{\CQ}$\index{$\Vdash_{\CQ}$} between a finite set (possibly empty) of formulas and a formula. Thus, for any set $X\cup\{A\}\Subset\FmAE$, we define $X\Vdash_{\CQ}A$ to mean that there is a (nonempty) list of $\LanAE$-formulas 
\[
A_1,\ldots, A_n,
\]
where the last formula, $A_n$, is $A$ and each formula in the list is either from $X$ or one of the axQ-axioms of $\CQ$ or can be obtained from two preceding formulas in the list by modus ponens or, if $A(u)$ is among preceding formulas  and a free variable $u$ does not occur in the formulas of $X$, by the rule~\eqref{E:Generalization-rule}.\index{rule!generalization} If $X\Vdash_{\CQ}A$, we say that $A$ is \textit{\textbf{derivable}} from $X$ in $\CQ$. We call any such derivation
a $\CQ$-\textit{\textbf{derivation}}.\index{$\CQ$-derivation}

Let $X\Vdash_{\CQ}A$. Assume that there is a $\CQ$-derivation whose free variables related to every application of~\eqref{E:Generalization-rule} in this derivation constitute a finite list $u_{i_1},\ldots,u_{i_n}$ of pairwise distinct symbols of~\eqref{E:free-variables}. We denote an arbitrary such derivation by
\[
X\Vdash_{\CQ}A\upharpoonright\{u_{i_1},\ldots,u_{i_n}\}.
\]

\begin{lem}\label{L:lemma-one-in-Q-con}
Let 
{\em\[
X\Vdash_{\CQ}A\upharpoonright\{u_{i_1},\ldots,u_{i_n}\}.\tag{$\ast$}
\]}
Assume that $u_{j_1},\ldots,u_{j_n}$ are pairwise distinct parameters that do not occur in $X$ and $\{u_{i_1},\ldots,u_{i_n}\}\cap\{u_{j_1},\ldots,u_{j_n}\}=\varnothing$. Then there is a derivation
{\em\[
X\Vdash_{\CQ}A\upharpoonright\{u_{j_1},\ldots,u_{j_n}\}.\tag{$\ast\ast$}
\]} 
\end{lem}
\begin{rem}
{\em
We note that parameters $u_{i_k}$ may occur in $(\ast\ast)$, but in this case, they are not eliminated by any application of the rule~\eqref{E:Generalization-rule} in $(\ast\ast)$.
}
\end{rem}
\begin{proof}
We prove by induction on the number of applications of the rule~\eqref{E:Generalization-rule} in  $(\ast)$ related to parameters in the list $\overline{u}=(u_{i_1},\ldots,u_{i_n})$. 

The basis case, that is, when $n=0$, is obvious.

Now we assume that there is at least one application of the rule~\eqref{E:Generalization-rule} in the derivation $(\ast)$. Suppose the transition from $B(u_{i_k})$ to $\forall x B(x)$ is the first such an application of~\eqref{E:Generalization-rule}. Then, right after  each formula $C(u_{i_k})$ that contains $u_{i_k}$ and precedes $B(u_{i_k})$, including $B(u_{i_k})$ itself, we place a formula $C(u_{j_k})$ obtained from $C(u_{i_k})$ by replacement of all occurrences of $u_{i_k}$ with $u_{j_k}$. The list of formulas from the very beginning up to $B(u_{j_k})$ thus obtained is a $\CQ$-derivation where only modus ponens is used. Next, we deduce $\forall x B(x)$ from $B(u_{j_k})$.

Since all formulas of the original derivation $(\ast)$ have been preserved, we have a new formal derivation of $A$ from $X$, in which all parameters of the list $\overline{u}$  can still be used in applications of \eqref{E:Generalization-rule}, but less number of times, because one application of this rule, related to $u_{i_k}$, has been replaced with an application of it related to $u_{j_k}$.
\end{proof}
\begin{cor}\label{C:lemma-one-in-Q-con}
Let {\em$X\subseteq Y\Subset\FmAE$}. Then
{\em\[
X\Vdash_{\CQ} A~\Longrightarrow~Y\Vdash_{\CQ} A.
\]}
\end{cor}
\begin{proof}
Assume that
\[
X\Vdash_{\CQ}A\upharpoonright\{u_{i_1},\ldots,u_{i_n}\}.
\]
Let us take parameters $u_{j_1},\ldots,u_{j_n}$ that do not occur in $Y$. By Lemma~\ref{L:lemma-one-in-Q-con},
\[
X\Vdash_{\CQ}A\upharpoonright\{u_{j_1},\ldots,u_{j_n}\}.
\]
Since neither of $u_{j_k}$ occurs in $Y$, we have that $Y\Vdash_{\CQ} A$.
\end{proof}

As before, we write sometimes $X, A$ instead of $X\cup\{A\}$ and
$\bigwedge X$ instead of $A_1\land\ldots\land A_n$,
where $X=\{A_1,\ldots,A_n\}$, that is a nonempty set of $\LanQ$-formulas. 

\begin{lem}\label{L:lemma-two-in-Q-con}
	For any finite set {\em$X\cup\{A,B\}\subseteq\FmAE$},
{\em	\[
	(X\Vdash_{\CQ}A~\text{and}~X\Vdash_{\CQ}B)~\Longrightarrow~X\Vdash_{\CQ}A\land B.
	\]}
\end{lem}
\begin{proof}
Suppose we have:
\[
X\Vdash_{\CQ}A\upharpoonright\{u_{i_1},\ldots,u_{i_n}\}~\text{and}~X\Vdash_{\CQ}B\upharpoonright\{u_{j_1},\ldots,u_{j_m}\}.
\]
Applying Lemma~\ref{L:lemma-one-in-Q-con} to each of the two last $\CQ$-derivations, we obtain:
\[
X\Vdash_{\CQ}A\upharpoonright\{u^{\prime}_{i_1},\ldots,u_{i_n}^{\prime}\}~\text{and}~X\Vdash_{\CQ}B\upharpoonright\{u_{j_1}^{\prime},\ldots,u_{j_m}^{\prime}\}.
\]
Without loss of generality, we can count that $\{[u_{i_1^{\prime}},\ldots,u_{i_n}^{\prime}\}\cap\{u_{j_1}^{\prime},\ldots,u_{j_m}^{\prime}\}=\varnothing$. Then, if we take the concatenation of the last formal proofs and add to the resulting list the axiom ax3Q, $B\rightarrow(A\land B)$ and $A\land B$, we obtain a $\CQ$-derivation of $A\land B$ from $X$ in $\CQ$.
\end{proof}
\begin{lem}\label{L:lemma-three-in-Q-con}	
For any nonempty finite set {\em$X\subseteq\FmAE$}, {\em$X\Vdash_{\CQ}\bigwedge X$}.
\end{lem}
\begin{proof}
We prove by induction on the number of formulas in $X$. If $X$ consists of one formula, then the conclusion is obvious. 

Suppose $X=Y\cup\{A\}$, where $Y\neq\varnothing$. By induction hypothesis,
$Y\Vdash_{\CQ}\bigwedge Y$ and hence $Y,A\Vdash_{\CQ}\bigwedge Y$. On the other hand, it is obvious that $Y,A\Vdash_{\CQ} A$. Applying Lemma~\ref{L:lemma-two-in-Q-con}, we receive that $X\Vdash_{\CQ}\bigwedge Y\land A$.
\end{proof}

The next lemma is similar to Lemma~\ref{L:deduction-theorem}.
\begin{lem}[deduction theorem]\label{L:deduction-theorem-in-Q-con}
For any set {\em$X\cup\{A,B\}\subseteq\FmAE$},
{\em\[
X,A\Vdash_{\CQ}B~\Longrightarrow~X\Vdash_{\CQ}A\rightarrow B.
\]}
\end{lem}
\begin{proof}
In the proof of this lemma, we use the argument given in the proof of Lemma~\ref{L:deduction-theorem}. Indeed, the first steps of the both proofs are the same, until, in this proof, we arrive at the consideration when in a given derivation a formula $\forall xA_k(u\backslash x)$ is obtained from a preceding formula $A_k(u)$, where a parameter $u$ does not occur in $X,A$. Then, by induction hypothesis, $X\Vdash_{\CQ}A\rightarrow A_k(u)$. Applying~\eqref{E:Generalization-rule}, we get that $X\Vdash_{\CQ}\forall x(A\rightarrow A_k(u\backslash x))$. Applying modus ponens to the last formula and ax14Q, we derive that $X\Vdash_{\CQ}A\rightarrow\forall x A_k(u\backslash x))$.
\end{proof}\begin{lem}\label{L:lemma-four-for-CQ}
For any finite nonempty set {\em$X\subseteq\FmAE$} and a formula $A$,
{\em\[
	X\Vdash_{\CQ}A~\Longleftrightarrow~\Vdash_{\CQ}\bigwedge X\rightarrow A.
	\]}
\end{lem}
\noindent\textit{Proof}~can be obtained by using Lemma~\ref{L:deduction-theorem-in-Q-con}, Lemma~\ref{L:lemma-three-in-Q-con}, and the property:
\[
\Vdash_{\CQ}(B\land C)\rightarrow A~\Longleftrightarrow~\Vdash_{\CQ}B\rightarrow(C\rightarrow A),
\]
for any $\LanAE$-formulas $A,B$ and $C$. We leave for the reader to complete the proof of this lemma. (Exercise~\ref{section:Q-consequence}.\ref{EX:lemma-four-for-CQ})
\\

Using the concept of a $\CQ$-derivation, we define three relations $\vdash_{\CQ}$ on
$\mathcal{P}(\FmAE)\times\FmAE$, one for each of the realizations of $\CQ$, as follows:
\begin{equation}\label{E:definition-Q-con-via-derivation}
X\vdash_{\CQ}A~\define~\text{there is $X_0\Subset X$ such that $X_0\Vdash_{\CQ}A$}.
\end{equation}

Lemmas~\ref{L:lemma-one-in-Q-con}--\ref{L:deduction-theorem-in-Q-con} suffice to show that each $\vdash_{\CQ}$ is a $\Q$-consequence relation. However, we postpone this conclusion and first prove one more lemma.
\begin{lem}\label{L:lemma-on-structurality-for-Q-con}
For any finite set $X$ of $\LanAE$-formulas and any $\LanT$-substitution $\sigma$,
{\em\[
	X\Vdash_{\CQ}A~\Longrightarrow~	\sigma(X)\Vdash_{\CQ}\sigma(A).
	\]}
\end{lem}
\begin{proof}
Suppose that
\[
X\Vdash_{\CQ}A\upharpoonright\{u_{i_1},\ldots,u_{i_n}\}\tag{$\ast$}
\]
and denote by $\mathcal{U}$ the set of parameters that occur in $\sigma(X)\cup\{\sigma(A)\}$. First, we assume that $\mathcal{U}\neq\varnothing$, that is
\[
\mathcal{U}=\{u_{j_1},\ldots,u_{j_k}\}.
\]

Taking into account Lemma~\ref{L:lemma-one-in-Q-con}, without loss of generality, we can assume that
\[
\{u_{i_1},\ldots,u_{i_n}\}\cap\{u_{j_1},\ldots,u_{j_k}\}=\varnothing.
\]

Then, we define a $\LanT$-substitution:
\[
\sigma^{\prime}(u):=\begin{cases}
\begin{array}{cl}
u &\text{if $u\in\{u_{i_1},\ldots,u_{i_n}\}$}\\
\sigma(u) &\text{if $u\notin\{u_{i_1},\ldots,u_{i_n}\}$}.
\end{array}
\end{cases}
\]

We note that for every formula $B\in X\cup\{A\}$, 
\[
\sigma^{\prime}(B)=\sigma(B).\tag{$\ast\ast$}
\]

Further, suppose that $(\ast)$ is the following finite sequence:
\[
A_1,\ldots, A_i,\ldots, A.\tag{$\ast\ast\ast$}
\]
Let us consider the sequence:
\[
\sigma^{\prime}(A_1),\ldots, \sigma^{\prime}(A_i),\ldots, \sigma^{\prime}(A).\tag{$\ast\ast\ast\ast$}
\]
We aim to show that $(\ast\ast\ast\ast)$ is a $\CQ$-derivation of $\sigma(A)$ from $\sigma(X)$. In view of the remark $(\ast\ast)$ above, it suffices to show that $(\ast\ast\ast\ast)$ is a $\CQ$-derivation.

Indeed, let $\sigma^{\prime}(A_i)$ of $(\ast\ast\ast\ast)$ is obtained from
$A_i$ of $(\ast\ast\ast)$.  We consider here the two main cases: when $A_i$ is obtained from $A_k$ and $A_k\rightarrow A_i$ by modus ponens and when $A_i=\forall x B(u\backslash x)$ and $B(u)$ precedes $A_i$ in $(\ast\ast\ast)$ (in which case $u\in\{u_{i_1},\ldots,u_{i_n}\}$).

The first case is obvious and we turn to the second. We notice that the parameter $u$ occurs in $\sigma^{\prime}(B(u))$ and does not occur in $\sigma^{\prime}(X)$. Therefore, 
\[
\forall x\sigma^{\prime}(B(u))(u\backslash x)=\sigma^{\prime}(\forall xB(u\backslash x)).
\]
\end{proof}

Finally, we arrive at the main conclusion of this section.
\begin{prop}\label{P:QC-is-Q-con}
Each {\em$\vdash_{\CQ}$} is a finitary structural $\mathcal{Q}$-consequence relation. In addition, each {\em$\CQ$} is an implicative $\Q$-logic.
\end{prop}
\begin{proof}
The main task is obviously to prove that the property (c) of Definition~\ref{D:consequnce-relation-single} holds for $\vdash_{\CQ}$.

To do this, assume that for some sets $X,Y$ and $Z$ of $\LanQ$-formulas,
$X\vdash_{\CQ}B$, for all formulas $B\in Y$, where $Y\neq\varnothing$, and $Y,Z\vdash_{\CQ}A$. By definition, we first have that 
\[
Y_{0},Z_{0}\Vdash_{\CQ}A\upharpoonright\{u_{i_1},\ldots,u_{i_n}\}
\] 
for some $Y_{0}\cup Z_{0}\Subset\FmAE$, where $Y_0$ is assumed to be unequal to $\varnothing$.

We denote
\[
Y_0:=\{B_1,\ldots,B_k\}.
\]

By premise, we have that
\[
X_1\Vdash_{\CQ}B_1,\ldots,X_n\Vdash_{\CQ}B_k,
\]
where $X_1\cup\ldots\cup X_k\Subset X$. Apply Lemma~\ref{L:lemma-one-in-Q-con}, we first obtain that
\[
X_1\cup\ldots\cup X_n\Vdash_{\CQ}B_1,\ldots,X_1\cup\ldots\cup X_k\Vdash_{\CQ}B_k.
\]

In virtue of Lemma~\ref{L:lemma-one-in-Q-con}, without loss of generality, we can assume that
\[
X_1\cup\ldots\cup X_n\Vdash_{\CQ}B_1\upharpoonright\mathcal{U},\ldots,X_1\cup\ldots\cup X_n\Vdash_{\CQ}B_n\upharpoonright\mathcal{U},
\]
where $\mathcal{U}$ is a finite set of parameters with $\mathcal{U}\cap\{u_{i_1},\ldots,u_{i_n}\}=\varnothing$. This allows us to form a $\CQ$-derivation of $A$ from $X_0\cup Z_0$. The latter in turn implies that $X,Z\vdash_{\CQ} A$.

The finitariness and structurality of the consequence relation $\vdash_{\CQ}$ follows from definition~\eqref{E:definition-Q-con-via-derivation} and Lemma~\ref{L:lemma-on-structurality-for-Q-con}, respectively.

We leave for the reader to prove that each abstract logic $\CQ$ is implicative. (Exercise~\ref{section:Q-consequence}.\ref{EX:QC-is-Q-con})
\end{proof}

In the sequel we will need the following property.
\begin{prop}\label{P:finite-set-equivalence}
For any set {\em$X\Subset\FmAE$},
{\em\[
X\vdash_{\CQ}A~\Longleftrightarrow~X\Vdash_{\CQ}A.
\]}
\end{prop}
\begin{proof}
Using definition~\eqref{E:definition-Q-con-via-derivation}, we observe that
\[
\varnothing\vdash_{\CQ}A~\Longleftrightarrow~\varnothing\Vdash_{\CQ}A.
\]
Therefore we assume that $X\neq\varnothing$.

The $\Leftarrow$-implication is obvious.

To prove the $\Rightarrow$-implication, we assume that
$X\vdash_{\CQ}A$. This implies that there is a set $Y\subseteq X$ such that
$Y\Vdash_{\CQ}A$. If $Y=\varnothing$, then we also have $Y\vdash_{\CQ}A$ and hence $X\vdash_{\CQ}A$.

Now, assume that $Y\neq\varnothing$. In virtue of Lemma~\ref{L:lemma-four-for-CQ}, $\Vdash_{\CQ}\bigwedge Y\rightarrow A$, and, by Lemma~\ref{L:lemma-three-in-Q-con}, 
$X\Vdash_{\CQ}\bigwedge X$. 

Further, we claim that $\Vdash_{\CQ}\bigwedge X\rightarrow\bigwedge Y$. (We leave for the reader to prove the last property; see Exercise~\ref{section:Q-consequence}.\ref{EX:finite-set-equivalence}.) Using the last derivability, we conclude that $X\Vdash_{\CQ}A$.
\end{proof}

\paragraph{Exercises~\ref{section:Q-consequence}}
\begin{enumerate}{\tiny }
	\item\label{EX:Q-consequence-is-consequence-1}Prove that for any $t$-valuation $v$ and any $\LanT$-substitution $\sigma$, for any term $t$, $v\circ\sigma[t]=v[\sigma(t)]$.
	\item\label{EX:MQ-relation-is-structural-Q-con-1} Prove that $\models_{\mathcal{M}^{\Q}}$ is a consequence relations, by using definition~\eqref{E:vdash_M^Q} and Definition~\ref{D:Q-matrix}.
	\item\label{EX:MQ-relation-is-structural-Q-con-2}Complete the proof of Proposition~\ref{P:MQ-relation-is-structural-Q-con}, by proving that for any $\LanQ$-formula $A$, any $t$-valuation $v$ and any $\LanT$-substitution $\sigma$, $\widehat{v\circ\sigma}[A]=\hat{v}[\sigma(A)]$.
	\item\label{EX:lemma-four-for-CQ}Complete the proof of Lemma~\ref{L:lemma-four-for-CQ}
	\item\label{EX:finite-set-equivalence} Suppose $Y\subseteq X\Subset\FmAE$ and $Y\neq\varnothing$. Then $\Vdash_{\CQ}\bigwedge X\rightarrow\bigwedge Y$; moreover, there is a derivation which does not use the axioms ax12Q--ax15Q and the rule~\eqref{E:Generalization-rule}.
	\item\label{EX:QC-is-Q-con}Prove that each logic $\CQ$ is implicative.
\end{enumerate}

\section{Lindenbaum $\Q$-matrices}\label{section:LT_Q-matrices}\index{Q-matrix!Lindenbaum}
The $\Q$-matrices we discuss in this section are based on the $\Q$-structure of Example~\ref{example:first}.
\begin{defn}[Lindenbaum $\Q$-matrix]\label{D:Lindenbaum-Q-matrix}
Given an abstract logic $\aLog_Q$, for any set {\em$X\subseteq\FmQ$}, we define a $\Q$-matrix
{\em$\LinSQ[X]:=\langle\langle\TmAl,\FmAlQ,\PredF,
\Phi_{\Q}\rangle,\ConSQ(X)\rangle$}, where for every $p\in\PredF$ with $\#(p)=n$ and any
{\em$t_1,\ldots,t_n\in\Tm$, $\bm{p}^{\mathbf{L}}(t_1,\ldots,t_n)=p(t_1,\ldots,t_n)$}; and the mappings
$\bm{p}^{\mathbf{L}}_{A(\overline{u})}$, $\bm{q}^{\mathbf{L}}_{(A(\overline{x},\overline{u^\ast}))}$ and $\varphi_{Q}$ are defined as in Example~\ref{example:first}.
Each {\em$\LinSQ[X]$} is called a \textbf{Lindenbaum} $\Q$-\textbf{matrix} for {\em$\aLog_Q$} $($relative to $X)$. We denote the Lindenbaum $Q$-matrix {\em$\LinSQ[\varnothing]$} by {\em$\LinSQ$}.
\end{defn}

Reviewing Example~\ref{example:first}, we conclude the following.
\begin{prop}\label{P:LinSQ-is-Q-matrix}
Every {\em$\LinSQ[X]$} is a $\Q$-matrix.
\end{prop}
\noindent\emph{Proof} is obvious and is left to the reader. (Exercise~\ref{section:LT_Q-matrices}.\ref{EX:LinSQ-is-Q-matrix})\\

The next proposition is similar to Proposition~\ref{P:lindenbaum-theorem}.

\begin{prop}\label{P:Lindenbaum-theorem-for-Q-con-1}
Let $\aLog_Q$ be a structural abstract $\mathcal{Q}$-logic. Then $\LinSQ$ is an $\aLog_Q$-model. In addition, if $X$ is closed under arbitrary
term substitution, then for any $\LanQ$-formula $A$,
{\em\[
X\vdash_{\mathcal{S}_Q}A~\Longleftrightarrow~\LinSQ[X]\models A;
\]}
in particular,
{\em
\[
A\in\bm{T}_{\mathcal{S}_Q}~\Longleftrightarrow~\LinSQ\models A.
\]}
\end{prop}
\begin{proof}
Suppose $X\vdash_{\mathcal{S}_Q}A$ and let $\sigma$ be any $t$-substitution.
Since the logic $\aLog_Q$ is structural, $\sigma(X)\vdash_{\mathcal{S}_Q}\sigma(A)$, that is, $\sigma(A)\in\ConSQ(\sigma(X))$. Thus, if $\sigma(X)\subseteq\ConSQ(\varnothing)$, then, obviously, $\sigma(A)\in\ConSQ(\varnothing)$. This shows that $\LinSQ$ is an $\aLog_Q$-model. 

Next, if $X$ is closed under $\LanT$-substitutions, then, in particular, $\sigma(X)\subseteq X$. This implies that $\ConSQ(\sigma(X))\subseteq\ConSQ(X)$ and hence $\sigma(A)\in\ConSQ(X)$.

Conversely, assume that $\LinSQ[X]\models A$. Taking the identity $\LanT$-substitution $\iota$, we have: $A=\iota(A)\in\ConSQ(X)$, that is, $X\vdash_{\mathcal{S}_Q}A$.
\end{proof}

The next definition is similar to Definition~\ref{D:atlas-bundle}.
\begin{defn}\label{D:Q-bundle-atlas}\index{Q-matrix!Q-bundle}
A $($nonempty$)$ family {\em$\mathcal{B}=\lbrace\langle\fA,D_i\rangle\rbrace_{i\in I}$} of $\Q$-matrices, which may differ by logical filters $D_i$,  is called a $\Q$-\textbf{bundle}. A pair {\em$\langle\fA,\lbrace D_i\rbrace_{i\in I}\rangle$}, where each $D_i$ is a logical filter in $\fA$, is called an $\Q$-\textbf{atlas}.\index{Q-matrix!Q-atlas} By definition, the $\Q$-\textbf{consequence relative to an atlas}
{\em$\langle\fA,\lbrace D_i\rbrace_{i\in I}\rangle$} is the same as the  $\Q$-consequence related the bundle {\em$\lbrace\langle\fA,D_i\rangle\rbrace_{i\in I}$}, which in turn is defined according to~{\em \eqref{E:vdash_M^Q}}.
Given an abstract logic $\aLog_{\Q}$, the bundle  $\lbrace\langle\FmAlQ,D\rangle\rbrace_{D\in\Sigma_{\aLog_{\Q}}}$, where $\Sigma_{\aLog_{\Q}}$ is the set of all $\aLog_{\Q}$-theories, is called a \textbf{Lindenbaum $\Q$-bundle} relative to $\aLog_{\Q}$ \index{Q-matrix!Lindenbaum Q-bundle} and the atlas {\em$\mat{Lin}_{\aLog_{\Q}}[\Sigma_{\aLog_{\Q}}]:=\langle\FmAlQ,\Sigma_{\aLog_{\Q}}\rangle$}, is called a \textbf{Lindenbaum $\Q$-atlas} relative to $\aLog_{\Q}$.\index{Q-matrix!Lindenbaum Q-atlas}
\end{defn}

According to the last definition and the main definition~\eqref{E:vdash_M^Q}, for any set $X\cup\lbrace A\rbrace\subseteq\FmQ$, 
\begin{equation}
\begin{array}{r}
X\models_{\mat{Lin}_{\aLog_{\Q}}[\Sigma_{\aLog_{\Q}}]}A~\Longleftrightarrow~
\text{for any $D\in\Sigma_{\aLog_{\Q}}$ and any $\LanT$-substitution $\sigma$},\\
\text{if $\sigma(X)\subseteq D$, then $\sigma(A)\in D$}.
\end{array}
\end{equation}

\begin{prop}\label{P:Lindenbaum-theorem-for-Q-con-2}
Let $\aLog_{\Q}$ be a structural $\Q$-consequence. Then for any set {\em$X\cup\lbrace A\rbrace\subseteq\FmQ$},
{\em\[
X\vdash_{\mathcal{S}_Q}A~\Longleftrightarrow~X\models_{\mat{Lin}_{\aLog_{\Q}}[\Sigma_{\aLog_{\Q}}]} A.
\]}
\end{prop}
\noindent\emph{Proof}~is similar to the proof of Proposition~\ref{P:Lindenbaum-theorem-generalized} and is left to the reader. (Exercise~\ref{section:LT_Q-matrices}.\ref{EX:Lindenbaum-theorem-for-Q-con-2})
\\

\paragraph{Exercises~\ref{section:LT_Q-matrices}}
\begin{enumerate}
	\item\label{EX:LinSQ-is-Q-matrix}Prove Proposition~\ref{P:LinSQ-is-Q-matrix}.
	\item \label{EX:Lindenbaum-theorem-for-Q-con-2}Prove Proposition~\ref{P:Lindenbaum-theorem-for-Q-con-2}.
\end{enumerate}

\section{Lindenbaum-Tarski $\Q$-structures}\label{section:Lindenbaum-Tarski-Q-structures}

Let $\aLog_{\Q}$ be an arbitrary $\Q$-logic. For any set $X\subseteq\FmQ$, we denote:
\begin{equation}\label{E:Theta_SQ-definition}
(A,B)\in\theta_{\aLog_{\Q}}(X)~\define~(A,B)\in\theta(\ConSQ(X)).
\end{equation}

Thus $\theta_{\aLog_{\Q}}(X)$ is a congruence on $\FmAlQ$. We employ the above definition when $\aLog_{\Q}$ is a unital $\Q$-logic.
\begin{defn}[Lindenbaum-Tarski $\Q$-structure]\label{D:LT-Q-structure}\index{Lindenbaum-Tarski Q-structure}
Let $\aLog_{\Q}$ be a unital $\Q$-logic and  {\em$X\subseteq\FmQ$}.
A structure
{\em\[
\LTSQ[X]:=\langle\TmAl,\langle\FmAlQ/\theta_{\aLog_{\Q}}(X), \one_{D}\rangle,\mathcal{P}_{F},\Q\rangle,
\]}where {\em$D=\ConSQ(X)$} and $\one_D$ is $D$ understood as a congruence class with respect to the congruence $\theta_{\aLog_{\Q}}(D)$ on $\FmAlQ$, is called
\textit{\textbf{Lindenbaum-Tarski $\Q$-structure}} $($of $\aLog_{\Q}$$)$ \textit{\textbf{relative to}} $D$ $($or \textit{\textbf{relative to}} $X$$)$
if each $p\in\mathcal{P}_{F}$ with $\#(p)=n$ is associate with a map
$\bm{p}^\mathbf{LT}:(t_1,\ldots,t_n)\mapsto p(t_1,\ldots,t_n)/\theta_{\aLog_{\Q}}(X)$. Accordingly, for any $t$-substitution $\sigma$, the value of a formula $A(\overline{u})$ in {\em$\LTSQ[X]$} under $\sigma$ is denoted by $\sigma_{X}[A(\overline{u})]$ and is defined inductively on the complexity of this formula in the following way.
If $A(\overline{u})=p(\overline{t})$ with $\#(p)=k$, where $\overline{u}=(u_{i_1},\ldots,u_{i_n})$ and $\overline{t}=(t_1,\ldots,t_k)$, then we have:
\begin{itemize}
	\item $\sigma_{X}[p(\overline{t})]:=\bm{p}^{\mathbf{L}}(\sigma(\overline{t}))/\theta_{\aLog_{\Q}}(X)=\bm{p}^{\mathbf{LT}}(\sigma(\overline{t}))$;
	\item$\bm{p}^{\mathbf{LT}}_{A(\overline{u})}(s_1,\ldots,s_n):=
	\bm{p}^{\mathbf{L}}(\delta_{\overline{u}\backslash\overline{s}}(\overline{t}))/\theta_{\aLog_{\Q}}(X)=\delta_{\overline{u}\backslash\overline{s}}(p(\overline{t}))/\theta_{\aLog_{\Q}}(X)$,\\ where
	{\em$\overline{s}=(s_1,\ldots,s_n)\in\Tm^n$}.
\end{itemize}

If $A(\overline{u})=FA_{1}\ldots A_k$, then
\begin{itemize}
	\item $\sigma_{X}[A(\overline{u})]:=F(\hat{\sigma}[A_{1}],\ldots,\hat{\sigma}[A_{k}])=\sigma(A(\overline{u}))/\theta_{\aLog_{\Q}}(X)$;
	\item $\bm{p}^{\mathbf{LT}}_{A(\overline{u})}(\overline{s})=F(\bm{p}^{\mathbf{LT}}_{A_1(\overline{u})}(\overline{s}),\ldots,\bm{p}^{\mathbf{LT}}_{A_k(\overline{u})}(\overline{s}))$;
\end{itemize}
and for any protoformula $A(\overline{x},\overline{u^{\ast}})$, where $\overline{x}=(x_{i_1},\ldots,x_{i_n})$,
\begin{itemize}
	\item {\em$\bm{q}^{\mathbf{LT}}_{(A(\overline{x},\overline{u^\ast}),\sigma)}(\overline{s}):=\bm{p}^{\mathbf{LT}}_{A(\overline{u},\overline{u^\ast})}(\overline{s},\sigma(\overline{u^{\ast}}))$};
\end{itemize}
and for any $Q\in\Q$,
\begin{itemize}
	\item 
	{\em	$\varphi_{Q}(\set{\bm{q}^{\mathbf{LT}}_{(A(\overline{x},\overline{u^\ast}),\sigma)}(\overline{s})}{\overline{s}\in\Tm^n}):=\sigma(Q x_{i_1}\ldots x_{i_n}A(x_{i_1},\ldots,x_{i_n},\overline{u^\ast})/\theta_{\aLog_{\Q}}(X)$}\\
		$=Q x_{i_1}\ldots x_{i_n} A(x_{i_1},\ldots,x_{i_n},\sigma(\overline{u^{\ast}}))/\theta_{\aLog_{\Q}}(X)$.	
\end{itemize}

\begin{itemize}
	\item $
	\sigma_{X}[Qx_{i_1}\ldots x_{i_n}A(x_{i_1},\ldots,x_{i_n},\overline{u^\ast})]:=\hat{\sigma}[Q x_{i_1}\ldots x_{i_n}A(\overline{x},\overline{u^\ast})]/\theta_{\aLog_{\Q}}(X)= \sigma(Q x_{i_1}\ldots x_{i_n}A(\overline{x},\overline{u^\ast}))/\theta_{\aLog_{\Q}}(X);
	$	
	\item $\bm{p}^{\mathbf{LT}}_{Q\overline{x} A(\overline{x},\overline{u^\ast})}(t_1,\ldots,t_m):=\widehat{\delta_{\overline{u^\ast}\backslash\overline{t}}}[Q\overline{x} A(\overline{x},\overline{u^\ast})]/\theta_{\aLog_{\Q}}(X)=\delta_{\overline{u^\ast}\backslash\overline{t}}(Q\overline{x} A(\overline{x},\overline{u^\ast}))/\theta_{\aLog_{\Q}}(X)$.
\end{itemize}
\end{defn}

 We will also use the notations $\LTSQ[D]$ for $\LTSQ[X]$ and $\one_{\aLog_{\Q}}[X]$ for $\one_{D}$, where $D=\ConSQ(X)$, and $\LTSQ:=\LTSQ[\varnothing]$.  
 
In the above definition, $\langle\FmAlQ/\theta_{\aLog_{\Q}}(X), \one_{D}\rangle$ is a unital expansion of $\FmAlQ/\theta_{\aLog_{\Q}}(X)$. In the sequel, we treat $\LTSQ[X]$ as a $\Q$-matrix
\[
\langle\langle\TmAl,\FmAlQ/\theta_{\aLog_{\Q}}(D), \mathcal{P}_{F},\Q\rangle,\{\one_{D}\}\rangle.
\]

Below, instead of $\sigma_{X}$, we occasionally write $\sigma_{D}$, where $D=\ConSQ(X)$. It is not difficult to see that $\sigma_{X}=\sigma_{D}$.
(Exercise~\ref{section:Lindenbaum-Tarski-Q-structures}.\ref{EX:sigma_X=sigma_D})

\begin{lem}\label{L:i_D-valuates-formulas-of-X-generating-D}
Let $\aLog_{\Q}$ be a unital $\Q$-logic and $D$ be an $\aLog_{\Q}$-theory. Then $\iota_{D}$ satisfies a formula $A$ in {\em$\LTSQ[D]$} if, and only if, $A\in D$.
\end{lem}
\noindent\emph{Proof} is left to the reader. (Exercise~\ref{section:Lindenbaum-Tarski-Q-structures}.\ref{EX:LT-algebra-for-Q-con})\\

Similarly to Proposition~\ref{P:valuations-in-LT}, we prove the following.
\begin{prop}\label{P:LT-algebra-for-Q-con}
		Let $\aLog_{\Q}$ be a unital logic and {\em$D=\ConSQ(X)$}.
	Then for any $\LanQ$-formula $A$, the following conditions are equivalent$\,:$
	{\em\[
		\begin{array}{cl}
		(\text{a}) & X\vdash_{\aLog_{\Q}}A;\\
		(\text{b}) &X\models_{\LTSQ[D]}A;\\
		(\text{c}) &\iota_{X}[A]=\one_{D}.
		\end{array}
		\]}
	In particular, the following conditions are equivalent$\,:$
	{\em\[
		\begin{array}{cl}
		(\text{a}^{\ast}) & A\in\ThmS;\\
		(\text{b}^{\ast}) &\sigma_{X}[A]=\one_{\bm{T}_{\aLog_{\Q}}},~\textit{for any $t$-substitution $\sigma$};\\
		(\text{c}^{\ast}) &\iota_{\bm{T}_{\aLog_{\Q}}}[A]=\one_{\bm{T}_{\aLog_{\Q}}}.
		\end{array}
		\]}
\end{prop}
\begin{proof}
Suppose $ X\vdash_{\aLog_{\Q}}A$. Let $\sigma$ be an $t$-substitution.
Assume that $\sigma_X$ satisfies in $\LTSQ[D]$ every $B\in X$. This means that for every $B\in X$, $\sigma_{X}[B]=\one_{D}$, that is, $\sigma_{X}[B]\in\ConSQ(X)$. In virtue of Proposition~\ref{P:semantics-for-formula-algebra}-a, for every $B\in X$, $\sigma(B)\in\ConSQ(X)$, that is, $\sigma(X)\subseteq\ConSQ(X)$. Since $\aLog_{\Q}$ is structural (being unital), $\sigma(A)\in\ConSQ(X)$. Hence $\sigma_X$ satisfies $A$ in $\LTSQ[D]$.\\

To prove the implication (b)$\Rightarrow$(c), we use Lemma~\ref{L:i_D-valuates-formulas-of-X-generating-D} to observe that $\iota_{D}$ satisfies every $B\in X$.\\

Now, assume that $\iota_{D}[A]=\one_{D}$, that is,
$\iota_D[A]/\theta_{\aLog_{\Q}}(D)=\one_{D}$. In virtue of Lemma~\ref{L:i_D-valuates-formulas-of-X-generating-D}, we conclude that $A\in\ConSQ(X)$.
\end{proof}

\paragraph{Exercises~\ref{section:Lindenbaum-Tarski-Q-structures}}
\begin{enumerate}
	\item\label{EX:sigma_X=sigma_D} Let $D=\ConSQ(X)$. Show that $\sigma_{X}=\sigma_{D}$.
	\item \label{EX:LT-algebra-for-Q-con}Let $\aLog_{\Q}$ be a unital $\Q$-logic and $D$ be an $\aLog_{\Q}$-theory. Prove that $\iota_{D}$ valuates a formula $A$ in $\LTSQ[D]$ if, and only if, $A\in D$. (\emph{Hint}: Use Proposition~\ref{P:semantics-for-formula-algebra}-a.)
\end{enumerate}

\section{Three $\Q$-matrix consequences}\label{section:Q-logics-three-examples}
In this section, we return to the languages $\Lan_{\forall\exists}$ and $\Lan_{(\forall\exists)^+}$ of Section~\ref{section:Q-con-via-inference-rules}.
We use these languages to further analyze the consequence relations $\vdash_{\CQ}$, where $\textsf{C}$ stands for $\PosQ$, $\ClQ$ or $\IntQ$. Therefore, it is convenient to denote the set of formulas, corresponding to a particular relation $\vdash_{\CQ}$, by $\textbf{Fm}_{\CQ}$ and call the element of this set $\CQ$-\textit{\textbf{formulas}}.

Accordingly, we consider $\Q$-structures of the following two forms:
\begin{equation}\label{E:two-structures-for-CQ}
	\begin{array}{rl}
		i) &\langle\langle\textsf{A};\FuncT\rangle,\langle\textsf{B};\land,\lor,\rightarrow\rangle,\PredF,\{\varphi_{\forall},\varphi_{\exists}\}\rangle\\
		ii) &\langle\langle\textsf{A};\FuncT\rangle,\langle\textsf{B};\land,\lor,
		\rightarrow,\neg\rangle,\PredF,\{\varphi_{\forall},\varphi_{\exists}\}\rangle.
	\end{array}
\end{equation}

We call a structure of the first form a $\Q_p$-\textit{\textbf{structure}} (or \textit{\textbf{positive $\Q$-structure}}), and that of the second form a $\Q_a$-\textit{\textbf{structure}} (or \textit{\textbf{assertoric $\Q$-structure}}). Accordingly, we define the notions of $\Q_p$-\textit{\textbf{matrix}} and
$\Q_a$-\textit{\textbf{matrix}} and, respectively, the notions of \textit{\textbf{unital $\Q_p$-matrix}} and
\textit{\textbf{unital $\Q_a$-matrix}}.

We denote by $\Lin_{\PosQ}[X]$, $\Lin_{\ClQ}[X]$ and $\Lin_{\IntQ}[X]$ the Lindenbaum $\Q$-matrices for $\PosQ$, $\ClQ$, and $\IntQ$, respectively; also, we employ the notation $\Lin_{\CQ}[X]$ to refer to one of these three types of Lindenbaum matrices.
Accordingly, we denote the formula algebra associated with $\CQ$ by $\mathfrak{F}_{\CQ}$, that is,
\[
\mathfrak{F}_{\CQ}:=\langle\Forms_{\CQ};\mathcal{F}_{\CQ}\rangle,
\]
where $\mathcal{F}_{\PosQ}:=\{\land,\lor,\rightarrow\}$ and $\mathcal{F}_{\ClQ}=\mathcal{F}_{\IntQ}:=\{\land,\lor,\rightarrow,\neg\}$.
The set of terms is denoted by $\textbf{Tm}_{\CQ}$.\\

\begin{prop}\label{P:theta_CQ-is-a-congruence}
Each abstract logic {\em$\CQ$} is implicative and, therefore, is unital. 	
\end{prop}
\begin{proof}
When the first statement is proved, we apply Proposition~\ref{P:some-classes-unital-logics-1}.

 To establish the first statement, it suffices to prove that: For any finite set $Z$,
\begin{itemize}
	\item $Z\Vdash_{\CQ} A\rightarrow B$ and $Z\Vdash_{\CQ} C\rightarrow D$ imply
$Z\Vdash_{\CQ} (A\land C)\rightarrow (B\land D)$\\ 
and $Z\Vdash_{\CQ} (A\lor C)\rightarrow (B\lor D)$;
\item $Z\Vdash_{\CQ}A\rightarrow B$ and $Z\Vdash_{\CQ}B\rightarrow C$ imply
$Z\Vdash_{\CQ}A\rightarrow C$; and
	\item if $\CQ=\ClQ$ or $\CQ=\IntQ$, then $Z\Vdash_{\CQ} A\rightarrow B$ implies $Z\Vdash_{\CQ} \neg B\rightarrow\neg A$.
\end{itemize}

We leave the task of proving these implications to the reader. (Exercise~\ref{section:Q-logics-three-examples}.\ref{EX:theta_CQ-is-a-congruence})\\

Using Proposition~\ref{P:QC-is-Q-con} and Corollary~\ref{C:some-classes-unital-logics-1}, we conclude that each $\CQ$ is unital.
\end{proof}

Applying Proposition~\ref{P:theta_CQ-is-a-congruence}, Proposition~\ref{P:some-classes-unital-logics-1}
and Proposition~\ref{P:Lindenbaum-theorem-for-Q-con-1}, we obtain the following.
\begin{cor}\label{C:theta_CQ-congruence}
For the congruence {\em$\theta_{\CQ}(X)$} defined in~\eqref{E:Theta_SQ-definition}, the following equivalences hold:
{\em\[
\begin{array}{rl}
	(A,B)\in\theta_{\CQ}(X)~&\Longleftrightarrow~X\vdash_{\CQ}A\rightarrow B~\textit{and}~X\vdash_{\CQ}B\rightarrow A,\\
	&\Longleftrightarrow~\Lin_{\CQ}[X]\models A\rightarrow B~\textit{and}~\,\Lin_{\CQ}[X]\models B\rightarrow A.\\
\end{array}
\]}
\end{cor}

We note that
\begin{equation}
	X\subseteq Y~\Longrightarrow~\theta_{\CQ}(X)\subseteq\theta_{\CQ}(Y).
\end{equation}

We denote:
\[
\LTCQ[X]:=\langle\TmAl,\langle\mathfrak{F}_{\CQ}/\theta_{\CQ}(X), \one_{\CQ}[X]\rangle,\mathcal{P}_{F},\{\varphi_{\forall},\varphi_{\exists}\}\rangle.
\]

According to Definition~\ref{D:LT-Q-structure}, we have:
\[
\bm{p}^{\mathbf{LT}}:(t_1,\ldots, t_n)\mapsto p(t_1,\ldots,t_n)/\theta_{\CQ}(X),
\]
where $(t_1,\ldots,t_n)\in\Tm^n$ and $p\in\PredF$ with $\#(p)=n$. Also, specifically for $\LanAE$, according to this definition, for any $t$-substitution $\sigma$, we have:\\
First, if $A(\overline{u})=p(\overline{t})$, where $\overline{u}=(u_{i_1},\ldots,u_{i_n})$ and all $u_{i_k}$ are pairwise distinct, then
\begin{itemize}
	\item $\sigma_{X}[p(\overline{t})]:=\bm{p}^{\mathbf{LT}}(\sigma(\overline{t}))=\bm{p}^{\mathbf{L}}(\sigma(\overline{t}))/\theta_{\CQ}(X)=p(\sigma(\overline{t}))/\theta_{\CQ}(X)$;
	\item$\bm{p}^{\mathbf{LT}}_{A(\overline{u})}(s_1,\ldots,s_n):=
	\bm{p}^{\mathbf{L}}(\delta_{\overline{u}\backslash\overline{s}}(\overline{t}))/\theta_{\CQ}(X)=\delta_{\overline{u}\backslash\overline{s}}(p(\overline{t}))/\theta_{\CQ}(X)$,\\ where {\em$\overline{s}=(s_1,\ldots,s_n)\in\Tm^n$}.
\end{itemize}

If $A(\overline{u})=FA_{1}\ldots A_k$, then
\begin{itemize}
	\item $\sigma_{X}[A(\overline{u})]:=F(\sigma_{X}[A_{1}],\ldots,\sigma_{X}[A_{k}])=\hat{\sigma}[A(\overline{u})]/\theta_{\CQ}(X)$;
	\item $\bm{p}^{\mathbf{LT}}_{A(\overline{u})}(\overline{s})=F(\bm{p}^{\mathbf{LT}}_{A_1(\overline{u})}(\overline{s}),\ldots,\bm{p}^{\mathbf{LT}}_{A_k(\overline{u})}(\overline{s}))$;
\end{itemize}
and for any protoformula $A(x,\overline{u^{\ast}})$ and $s\in\Tm$, 
\begin{itemize}
	\item $\bm{q}^{\mathbf{LT}}_{(A(x,\overline{u^\ast}),\sigma)}(s):=\bm{p}^{\mathbf{LT}}_{A(u,\overline{u^\ast})}(s,\sigma(\overline{u^{\ast}}))$;
\end{itemize}
and for any $Q\in\{\forall,\exists\}$,
\begin{itemize}
	\item 
	$\varphi_{Q}(\set{\bm{q}^{\mathbf{LT}}_{(A(x,\overline{u^\ast}),\sigma)}(s)}{s\in\Tm}):=\varphi_{Q}(\set{\bm{q}^{\mathbf{L}}_{(A(x,\overline{u^\ast}),\sigma)}(s)}{s\in\Tm})/\theta_{\CQ}(X)=\sigma(Q x A(x,\overline{u^\ast}))/\theta_{\CQ}(X)
	=Q x A(x,\sigma(\overline{u^{\ast}}))/\theta_{\CQ}(X)$;	
\end{itemize}
\begin{itemize}
	\item $
	\sigma_{X}[Q x A(x,\overline{u^\ast})]:=\hat{\sigma}[Q x A(x,\overline{u^\ast})]/\theta_{\CQ}(X)= \sigma(Q x A(x,\overline{u^\ast}))/\theta_{\CQ}(X);
	$	
	\item $\bm{p}^{\mathbf{LT}}_{Q x A(x,\overline{u^\ast})}(t_1,\ldots,t_m):=\widehat{\delta_{\overline{u^\ast}\backslash\overline{t}}}[Q x A(x,\overline{u^\ast})]/\theta_{\CQ}(X)=\delta_{\overline{u^\ast}\backslash\overline{t}}(Q x A(x,\overline{u^\ast}))/\theta_{\CQ}(X)$.
\end{itemize}

Similarly to Proposition~\ref{P:semantics-for-formula-algebra}, we summarize the multiple definition above in the following proposition
\begin{prop}\label{P:quotients-characteristics}
	Let $\sigma$ be an $\FuncT$- substitution, $A(\overline{u})$ and $B(u,\overline{u^\ast})$ be formulas such that $\overline{u}=(u_{i_1},\ldots,u_{i_n})$ is the list of all parameters of the former and the parameters of the latter are partitioned in two lists ---
	$(u)$ and $\overline{u^\ast}=(u_{j_1},\ldots,u_{j_m})$ with $u\notin\overline{u^{\ast}}$. Then the following hold in {\em$\langle\TmAl,\mathfrak{F}_{\CQ}/\theta_{\CQ}(X),\mathcal{P}_F,\{\varphi_{\forall},\varphi_{\exists}\}\rangle$}:
	{\em\[
		\begin{array}{cl}
			(\text{a}) &\sigma_{X}[A(\overline{u})] = \sigma(A(\overline{u}))/\theta_{\CQ}(X);\\
			(\text{b}) &\bm{p}^{\mathbf{LT}}_{B(u,\overline{u^{\ast}})}(s,\sigma(\overline{u^{\ast}}))=
			\delta_{(\sigma,u\backslash s)}(B(u,\overline{u^{\ast}}))/\theta_{\CQ}(X);\\
			(\text{c}) &\bm{q}^{\mathbf{LT}}_{(B(x,\overline{u^\ast}),\sigma)}(s)=
			\delta_{(\sigma,u\backslash s)}(B(u,\overline{u^\ast}))/\theta_{\CQ}(X).
		\end{array}
		\]}
\end{prop}

\vskip 0.1in
For the next proposition, we also need the following relation on $\mathfrak{F}_{\CQ}/\theta_{\CQ}(X)$.
\[
A/\theta_{\CQ}(X)\le_{\CQ}B/\theta_{\CQ}(X)~\define~X\vdash_{\CQ}A\rightarrow B.
\]

The soundness of $\le_{\CQ}$ can be easily shown. For this, one needs to prove that if $A/\theta_{\CQ(X)}=C/\theta_{\CQ}(X)$ and $B/\theta_{\CQ}(X)=D/\theta_{\CQ}(X)$, then
\[
X\vdash_{\CQ}A\rightarrow B~\Longleftrightarrow~X\vdash_{\CQ}C\rightarrow D.
\]
(Exercise~\ref{section:Q-logics-three-examples}.\ref{EX:less-then-or-equal-soundness})
\begin{lem}
The relation {\em$\le_{\CQ}$} is a partial order on {\em$\langle\mathfrak{F}_{\CQ}/\theta_{\CQ}(X), \one_{\CQ}[X]\rangle$} so that
{\em$\one_{\CQ}[X]$} is a greatest element.
\end{lem}
\begin{proof}
We prove that the relation $\le_{\CQ}$ is transitive and that $\one_{\CQ}[X]$ is a greatest element.

Indeed, suppose $A/\theta_{\CQ}(X)\le_{\CQ}B/\theta_{\CQ}(X)$ and $B/\theta_{\CQ}(X)\le_{\CQ}C/\theta_{\CQ}(X)$. By definition, it means that
$X\vdash_{\CQ}A\rightarrow B$ and $X\vdash_{\CQ}B\rightarrow C$, that is,
for some $X_1,X_2\Subset X$, $X_1\Vdash_{\CQ}A\rightarrow B$ and $X_2\Vdash_{\CQ}B\rightarrow C$. In virtue of Corollary~\ref{C:lemma-one-in-Q-con}, $X_1\cup X_2\Vdash_{\CQ}A\rightarrow B$ and $X_1\cup X_2\Vdash_{\CQ}B\rightarrow C$. This implies that $X_1\cup X_2\Vdash_{\CQ}A\rightarrow C$. Hence $X\vdash_{\CQ}A\rightarrow C$.

Now, if $A$ is an $\CQ$-thesis, then for any formula $B$, $\varnothing\Vdash B\rightarrow A$. Hence $B/\theta_{\CQ}(X)\le_{\CQ}A/\theta_{\CQ}(X)$. It remains to note that $A/\theta_{\CQ}(X)=\one_{\CQ}[X]$.
\end{proof}

\begin{prop}\label{P:threre-Q-logics-main-theorem}
Let {\em$X\subseteq\textbf{Fm}_{\CQ}$}. Then the following are true:
{\em\[
\begin{array}{cl}
(\text{a}) &\fA_1=\langle\mathfrak{F}_{\PosQ}/\theta_{\PosQ}(X),\one_{\PosQ}[X]\rangle~\textit{is a relatively pseudo-complemented lattice};\\
(\text{b}) &\fA_2=\langle\mathfrak{F}_{\ClQ}/\theta_{\ClQ}(X),\one_{\ClQ}[X]\rangle~\textit{is a Boolean algebra, where $x\rightarrow y:=\neg x\lor y$};\\
(\text{c}) &\fA_3=\langle\mathfrak{F}_{\IntQ}/\theta_{\IntQ}(X),\one_{\IntQ}[X]\rangle~\textit{is a Heyting algebra}.
\end{array}
\]}
In addition, the relation $\le_i$ defined in each $\fA_i$ according to~\eqref{E:diego-algebra-ordering}, or equivalently according to~\eqref{E:ordering-in-lattice}, coincides with the corresponding relation {\em$\le_{\CQ}$}.
Moreover, if the set of parameters not occurring in the formulas of $X$ is infinite, then for any $t$-substitution $\sigma$ and  any protoformula $A(x)$ with free occurrences of only $x\in\VarT$, both
{\em\[
\bigwedge\set{\bm{q}^{\mathbf{LT}}_{(A(x),\sigma)}(t)}{t\in \Tm}~\textit{and}~\bigvee\set{\bm{q}^{\mathbf{LT}}_{(A(x),\sigma)}(t)}{t\in \Tm}
\]}exist, where $\bigwedge$ and $\bigvee$ are the operations of a greatest lower bound and a least upper bound, respectively, with respect to {\em$\le_{\CQ}$} in 
{\em$\mathfrak{F}_{\CQ}/\theta_{\CQ}(X)$}. Consequently, defining
$\varphi_{\forall}$ as $\bigwedge$ and $\varphi_{\exists}$ as $\bigvee$ in a 
{\em$\Q$}-matrix {\em$\langle\TmAl,\langle\mathfrak{F}_{\CQ}/\theta_{\CQ}(X),\one_{\CQ}[X]\rangle,
\mathcal{P}_{F},\{\varphi_{\forall},\varphi_{\exists}\}\rangle$}, we have:
{\em
\begin{equation}\label{E:three-Q-logics-main-theorem-1}
\sigma_{X}[\forall x A(x)]=\bigwedge\set{\bm{q}^{\mathbf{LT}}_{(A(x),\sigma)}(t)}{t\in \Tm}
\end{equation}}
and
{\em
\begin{equation}\label{E:three-Q-logics-main-theorem-2}
\sigma_{X}[\exists x A(x)]=\bigvee\set{\bm{q}^{\mathbf{LT}}_{(A(x),\sigma)}(t)}{t\in \Tm}. 
\end{equation}
}
\end{prop}
\begin{proof}
Since the restriction of any Heyting algebra to the signature $\langle\land,\lor,\rightarrow\rangle$ is a relatively pseudo-complemented lattice and in view of that the sentential fragments  of $\PosQ, \ClQ$ and $\IntQ$ are coincident, respectively, with $\Pos, \Cl$ and $\Int$, in virtue of Proposition~\ref{P:LT-algebra-for-Q-con}, we observe that (a)--(c) above hold.\\

To show that the relation $\le_{\CQ}$ is coincident with the corresponding $\le_i$, we consider the following equivalences:
\[
\begin{array}{rl}
	A/\theta_{\CQ}(X)\le_{\CQ}B/\theta_{\CQ}(X) &\Longleftrightarrow~X\vdash_{\CQ}A\rightarrow B\\
	&\Longleftrightarrow~A\rightarrow B\in\ConCQ(X)\\
	&\Longleftrightarrow~A\rightarrow B/\theta_{\CQ}(X)
	=\one_{\CQ}[X]\\
	&\Longleftrightarrow~A/\theta_{\CQ}(X)\rightarrow B/\theta_{\CQ}(X)
	=\one_{\CQ}[X].
\end{array}
\]

To prove~\eqref{E:three-Q-logics-main-theorem-1} and~\eqref{E:three-Q-logics-main-theorem-2}, we fix a protoformula $A(x,\overline{u^\ast})$ with free occurrences of a single bound variable $x$ and a finite list
$\overline{u^\ast}=(u_{j_1},\ldots,u_{j_m})$ of all parameters that occur in this protoformula. Further, we consider a formula $A(u,\overline{u^\ast})$ obtained from $A(x,\overline{u^\ast})$ by replacing each free occurrence of $x$ with a parameter $u$, providing that $u$ does not occur in this protoformula. 

In view of Proposition~\ref{P:quotients-characteristics}, we have:
\[
\bm{q}^{\mathbf{LT}}_{(A(x,\overline{u^{\ast}}),\sigma)}(t)=\widehat{\delta_{(\sigma,u\backslash t)}}[A(u,\overline{u^{\ast}})]/\theta_{\CQ}(X)=A(t,\sigma(\overline{u^{\ast}}))/\theta_{\CQ}(X).
\]
That is,
\[
\set{\bm{q}^{\mathbf{LT}}_{(A(x,\overline{u^{\ast}}),\sigma)}(t)}{t\in\Tm}=\set{A(t,\sigma(\overline{u^{\ast}}))/\theta_{\CQ}(X)}{t\in\Tm}.
\]

We recall that, in virtue of ax12Q,
\[
\Vdash_{\CQ}\forall x A(x,\sigma(\overline{u^{\ast}}))\rightarrow A(t,\sigma(\overline{u^{\ast}})),
\]
for an arbitrary $t\in\Tm$. Hence,
\[
\forall x A(x,\sigma(\overline{u^{\ast}}))/\theta_{\CQ}(X)\le_{\CQ}
A(t,\sigma(\overline{u^{\ast}}))/\theta_{\CQ}(X).
\]

On the other hand, assume that for some formula $B$,
\[
X\vdash_{\CQ}B\rightarrow A(t,\sigma(\overline{u^{\ast}})),
\]
for any $t\in\Tm$. In particular,
\[
X\vdash_{\CQ}B\rightarrow A(u^{\prime},\sigma(\overline{u^{\ast}})),
\]
where $u$ is a parameter which occurs neither in $X$ nor in $B$. This means that for some set $Y\Subset X$,
\[
Y\Vdash_{\CQ}B\rightarrow A(u,\sigma(\overline{u^{\ast}})).
\]
We observe that $u$ does not occur in the formulas of $Y$. Then, by the rule~\eqref{E:Generalization-rule},
\[
Y\Vdash_{\CQ}\forall x(B\rightarrow A(x,\sigma(\overline{u^{\ast}})).
\]
Applying modus ponens to the last formula and ax14Q, we derive that
\[
Y\Vdash_{\CQ}B\rightarrow\forall x A(x,\sigma(u^{\ast})),
\]
that is,
\[
X\vdash_{\CQ}B\rightarrow\forall x A(x,\sigma(u^{\ast})).
\]
Hence,
\[
B/\theta_{\CQ}(X)\le_{\CQ}\forall x A(x,\sigma(\overline{u^{\ast}}))/\theta_{\CQ}(X).
\]

This completes the proof of~\eqref{E:three-Q-logics-main-theorem-1}. We leave for the reader to prove~\eqref{E:three-Q-logics-main-theorem-2}. (Exercise~\ref{section:Q-logics-three-examples},\ref{EX:threre-Q-logics-main-theorem})
\end{proof}

The last proposition induces the following definition.
\begin{defn}[$\QP$-matrix, $\QB$-matrix, and $\QH$-matrix]\index{$\QP$-matrix}\index{$\QB$-matrix}\index{$\QH$-matrix}
	A unital $\Q_p$-matrix  {\em$\langle\langle\textsf{A};\FuncT\rangle,\langle\textsf{B};\land,\lor,\rightarrow,\one\rangle,\PredF,\{\varphi_{\forall},\varphi_{\exists}\}\rangle$} is called a {\em$\QP$}-\textbf{matrix} if {\em$\langle\textsf{B};\land,\lor,\rightarrow,\one\rangle$} is a relatively pseudo-complemented lattice, and for any $t$-valuation $v$ in {\em$\langle\textsf{A};\FuncT\rangle$} and any formulas $\forall x A(x,\overline{u})$ and $\exists x A(x,\overline{u})$, both
	{\em$\bigwedge\set{\bm{q}_{(A(x,\overline{u}),v)}(a,v[\overline{u}])}{a\in \textsf{A}}$} and {\em$\bigvee\set{\bm{q}_{(A(x,\overline{u}),v)}(a,v[\overline{u}])}{a\in \textsf{A}}$}, where $\bigwedge$ is a greatest lower bound and $\bigvee$ is a least upper bound, exist so that $\varphi_{\forall}$ and $\varphi_{\exists}$ are defined as $\bigwedge$ and $\bigvee$, respectively. 
	
	A unital $\Q_a$-matrix  {\em$\langle\langle\textsf{A};\FuncT\rangle,\langle\textsf{B};\land,\lor,\rightarrow,\neg,\one\rangle,\PredF,\{\varphi_{\forall},\varphi_{\exists}\}\rangle$} is called a {\em$\QB$}-\textbf{matrix} if {\em$\langle\textsf{B};\land,\lor,\rightarrow,\neg,\one\rangle$}, where $x\rightarrow y:=\neg x\lor y$, is a Boolean algebra;
	it is called a {\em$\QH$}-\textbf{matrix} if {\em$\langle\textsf{B};\land,\lor,\rightarrow,\neg,\one\rangle$} is a Heyting algebra, providing that the other conditions and definitions in the two last cases are the same as for {\em$\QP$}-matrices.
\end{defn}

Referring collectively to each of these $\Q$-matrices, we use the term
$\QC$-\textit{\textbf{matrix}}.

Thus, according to the definition of an $f$-valuation (Section~\ref{section:Q-semantics}), in all $\QC$-matrices, for any $t$-valuation $v$ and protoformula $A(x,\overline{u})$,
\begin{equation}\label{E:universal-quantifier}
	\hat{v}[\forall x A(x,\overline{u})]:=\bigwedge\set{\bm{q}_{(A(x,\overline{u}),v)}(a,v[\overline{u}])}{a\in \textsf{A}}
\end{equation}
and
\begin{equation}\label{E:existential-quantifier}
	\hat{v}[\exists x A(x,\overline{u})]:=\bigvee\set{\bm{q}_{(A(x,\overline{u}),v)}(a,v[\overline{u}])}{a\in \textsf{A}}.
\end{equation}
\begin{cor}\label{C:LT_CQ-as-QC-matrix}
Let {\em$X\subseteq\FmAE$} such that the set of parameters not occurring in the formulas of $X$ is infinite. Then each {\em$\LTCQ[X]$} is a {\em$\QC$}-matrix.
\end{cor}
\begin{defn}[logics $\QP$, $\QB$ and $\QH$]\index{logic!$\QP$}\index{logic!$\QB$}\index{logic!$\QH$}
Abstract logics {\em$\QP$}, {\em$\QB$} and {\em$\QH$} are defined as the consequence relations that are determined, respectively, by the classes of {\em$\QP$}-matrices, {\em$\QB$}-matrices, and {\em$\QH$}-matrices. The corresponding consequence relations are denoted by {\em$\vdash_{\QP}$}, {\em$\vdash_{\QB}$}, and {\em$\vdash_{\QH}$}.\index{$\vdash_{\QP}$}\index{$\vdash_{\QB}$}\index{$\vdash_{\QH}$} Referring to each of these consequence relations, we employ the notation {\em$\vdash_{\QC}$}.
\end{defn}

We aim to show that every $\QC$-matrix is a $\CQ$-model. We will do it by means of the following lemma which exemplifies the simplest case of this property.
\begin{lem}\label{L:CQ-vs-QC-theses}
Every {\em$\CQ$}-thesis is valid in any {\em$\QC$}-matrix; in other words,
{\em\[
\vdash_{\CQ}A~\Longrightarrow~\vdash_{\QC}A.
\]}
\end{lem}
\begin{proof}
In view of Proposition~\ref{P:finite-set-equivalence}, it suffices to prove the implication:
\[
\Vdash_{\CQ}A~\Longrightarrow~\vdash_{\QC}A;
\]
in other words, every formula $A$ derivable in $\CQ$ is valid in any $\QC$-matrix.\\

We prove the last implication by induction on the length of a $\CQ$-derivation of a formula $A$. The ``sentential axiom schemata'' ax1Q--ax11Q are easy and we leave the base case, concerning these schemata, to the reader (Exercise~\ref{section:Q-logics-three-examples}.\ref{EX:CQ-vs-QC-theses}) and move to the schemata ax12Q--ax15Q.\\

Let us fix a $\QC$-matrix:
\[
\fM=\lr{\lr{\textsf{A};\FuncT},\lr{\textsf{B};\land,\lor,\rightarrow,\neg,
\one},\PredF,\{\varphi_{\forall},\varphi_{\exists}\}}\footnote{If $\CQ=\PosQ$, instead of $\lr{\textsf{B};\land,\lor,\rightarrow,\neg,\one}$, we have $\lr{\textsf{B};\land,\lor,\rightarrow,\one}$.}
\]
and a $t$-valuation $v$ in $\fM$.

We aim to show that $\hat{v}[\forall xA(x)\rightarrow A(u\backslash t)]
=\one$, that is, $\hat{v}[\forall x A(x)]\le\hat{v}[A(u\backslash t)]$, where $\le$ is the partial order in $\lr{\textsf{B};\FuncF,\one}$ which is determined by $\land$ (or by $\lor$) from $\FuncF$. The formula $A(u)$ may have other free variables, different from $u$; we denote them by $\overline{u^\ast}=(u_{j_1},\ldots,u_{j_m})$.

According to~\eqref{E:universal-quantifier},
\[
\hat{v}[\forall x A(x,\overline{u^\ast})]=\bigwedge\set{\bm{q}_{(A(x,\overline{u^\ast}),v)}(a)}{a\in\textsf{A}}, \tag{$\ast$}
\]
and, in virtue of Lemma~\ref{L:formula-as-predicate} and Lemma~\ref{L:valuation-of-formula},
\begin{align*}
	\hat{v}[A(u\backslash t)]=\bm{p}_{A(u,\overline{u^\ast})}(v\circ\sigma_{u\backslash t}[u],v\circ\sigma_{u\backslash t}[\overline{u^{\ast}}])\\
	=\bm{p}_{A(u,\overline{u^\ast})}(v[t],v[\overline{u^{\ast}}]),
\end{align*}

In virtue of~\eqref{E:q-and-p-at-v}, the desired conclusion is obvious.\\

Since 
\[
\hat{v}[\exists x A(x,\overline{u^\ast})]=\bigvee\set{\bm{q}_{(A(x,\overline{u^\ast}),v)}(a)}{a\in\textsf{A}},
\]
in virtue of Lemma~\ref{L:valuation-of-formula}, $\hat{v}[A(u\backslash t)]\le \hat{v}[\exists x A(x)]$.\\

We leave for the reader to show that the axioms ax14Q and ax15Q are also valid in any $\QC$-matrix. (Exercise~\ref{section:Q-logics-three-examples}.\ref{EX:Q_C-consequence-relations})\\

Now suppose a formula $B$ is obtained by modus ponens from $A$ and $A\rightarrow B$. By induction hypothesis, for any $t$-valuation $v$ is any $\QC$-structure
$\hat{v}[A]=\hat{v}[A]\rightarrow\hat{v}[B]=\one$. This implies that $\hat{v}[B]=\one$.\\

Finally, suppose, by induction hypothesis, $w[A(u,\overline{u^{\ast}})]=\one$, where $u\notin\overline{u^{\ast}}$, for any valuation $w$ in $\fM$.

Assume, for contradiction, $\hat{v}[\forall x A(x,\overline{u^\ast})]\neq\one$.
Then, because of $(\ast)$, there is an element $a\in\textsf{A}$ such that $\bm{q}_{(A(x,\overline{u^\ast}),v)}(a)\neq\one$. With the purpose of deriving a contradiction, we define:
\[
w_{(v,u\backslash a)}[u^\prime]=\begin{cases}
	\begin{array}{cl}
		a &\text{if $u^\prime=u$}\\
		v[u^\prime] &\text{otherwise}.
	\end{array}
\end{cases}
\]

By definition,
\begin{align*}
	\widehat{w_{(v,u\backslash a)}}[A(u,\overline{u^\ast})]=\bm{p}_{A(u,\overline{u^\ast})}(a,v[\overline{u^{\ast}}])\\
	=\bm{q}_{(A(x,\overline{u^\ast}),v)}(a)\neq\one.
\end{align*}
A contradiction.
\end{proof}

\begin{prop}[soundness]\label{P:CQ-soundness}
Every {\em$\QC$}-matrix is a {\em$\CQ$}-model; in other words,
{\em\[
	X\vdash_{\CQ}A~\Longrightarrow~X\vdash_{\QC}A.
	\]}
\end{prop}
\begin{proof}
Suppose $X\vdash_{\CQ}A$, that is, for some set $Y\Subset X$, $Y\Vdash_{\CQ} A$. In virtue of Lemma~\ref{L:lemma-four-for-CQ}, $\Vdash_{\CQ}\bigwedge Y\rightarrow A$. By Lemma~\ref{L:CQ-vs-QC-theses}, $\vdash_{\QC}\bigwedge Y\rightarrow A$, which in particular implies that $Y\vdash_{\QC}A$. And by monotonicity, we conclude that $X\vdash_{\QC}A$.
\end{proof}

\begin{prop}[limited completeness]\label{P:limited-completeness}
Let {\em$X\subseteq\FmAE$} such that the set of parameters not occurring in the formulas of $X$ is infinite. Then
{\em\[
X\vdash_{\QC}A~\Longrightarrow~X\vdash_{\CQ}A.
\]}
\end{prop}
\begin{proof}
Assume that $X\not\vdash_{\CQ}A$. According to Proposition~\ref{P:LT-algebra-for-Q-con}, the algebra $\LTCQ[X]$ validates all formulas of $X$ and refutes $A$. On the other hand, in virtue of Corollary~\ref{C:LT_CQ-as-QC-matrix}, $\LTCQ[X]$ is a $\QC$-matrix. Thus $X\not\vdash_{\QC}A$.
\end{proof}

\begin{rem}
{\em
The reader should remember that each abstract logic $\QC$ is a $\Q$-matrix consequence. Let $\mathcal{M}^{\Q}$ be the abstract class of unital $\Q$-matrices of the type $\langle\langle\textsf{A};\FuncT\rangle,\booleTwo,\mathcal{P}_F,\{\varphi_{\forall},\varphi_{\exists}\}\rangle$. Suppose that the semantics of language $\Lan_{\forall\exists}$ is defined as in Example~\ref{example:ordinary-quantifiers}. Then the well-known completeness of first-order logic reads:
\[
X\vdash_{\ClQ}A~\Longleftrightarrow~X\models_{\mathcal{M}^{Q}}A.
\]

Comparing the last equivalence with Proposition~\ref{P:limited-completeness} when $\textsf{C}=\Cl$, we note that the completeness expressed by the last equivalence is reached at the cost of narrowing the class of determining $\Q$-matrices. However, due to the fact that the relation $\Vdash_{\ClQ}$ is known to be undecidable, to show that $X\not\vdash_{\ClQ}A$, it could be easier if we have a class of $\ClQ$-models as large as possible.
}
\end{rem}

\paragraph{Exercises~\ref{section:Q-logics-three-examples}}
\begin{enumerate}
	\item\label{EX:less-then-or-equal-soundness}Show that the definition of relation $\le_{\CQ}$ is sound. (\emph{Hint}: use Corollary~\ref{C:lemma-one-in-Q-con}.)
		\item \label{EX:theta_CQ-is-a-congruence}Complete the proof of Proposition~\ref{P:theta_CQ-is-a-congruence} (\textit{Hint}: use Corollary~\ref{C:lemma-one-in-Q-con} and Lemma~\ref{L:deduction-theorem-in-Q-con}).
	\item\label{EX:threre-Q-logics-main-theorem} Complete the proof of Proposition~\ref{P:threre-Q-logics-main-theorem} by proving~\eqref{E:three-Q-logics-main-theorem-2}.
	\item\label{EX:CQ-vs-QC-theses}Prove that every $\CQ$-axiom in the list ax1Q--ax11Q is valid in any $\QC$-matrix.
	\item\label{EX:Q_C-consequence-relations}Show that the axioms ax14Q and ax15Q are valid in any $\QC$-matrix. 
\end{enumerate}

\section{Historical notes}

Algebraic treatment of predicate language systems, which was discussed in this chapter,\footnote{The reader can find other algebraic approaches to first-order languages ​​in the theory of polyadic algebras and of cylindrical algebras; see~\cite{halmos1962} and~\cite{henkin-monk-tarski1985}, respectively.} began with A. Mostowski, who demonstrated in~\cite{mostowski1948} that if ordinary quantifiers, $\forall$ and $\exists$, are interpreted as a greatest lower bound and a least upper bound, respectively, in a complete lattice, then 
\begin{equation}\label{E:historical-Q-con-1}
\Vdash_{\IntQ}A~\Longrightarrow~\models_{\mathcal{M}^{\Q}_{1}}A
\end{equation}
holds, where $A$ is a sentence and $\mathcal{M}^{\Q}_{1}$ is the class of $\Q$-structures
$\langle\langle\textsf{A};\FuncT\rangle,\langle\textsf{B};\land,\lor,\rightarrow,\neg\rangle,\PredF,\{\varphi_{\forall},\varphi_{\exists}\},\rangle$, where $\langle\textsf{B};\land,\lor,\rightarrow,\neg\rangle$ is a complete Heyting algebra.

Mostowski states clearly the aim of his approach:
\begin{quote}
	``$\dots$ to outline a general method which permits us to establish the intuitionistic non-deducibility of many formulas $\dots$'' \cite{mostowski1948}, p. 204 
\end{quote}

As the reader can see, Mostowski's approach is consistent with the conception of separating tools (Section~\ref{section:separating-means}) and, therefore, the spirit of this book.

Using his method, Mostowski shows that the formula $\neg\neg(\forall x\,\neg\neg p(x)\rightarrow\neg\neg\forall x\, p(x))$ is not derivable in $\IntQ$, the result which had earlier been established by S. C. Kleene and D. Nelson in a different way.

Regarding~\eqref{E:historical-Q-con-1}, we note that the difference between it and the $\Rightarrow$-implication of Lemma~\ref{L:CQ-vs-QC-theses} is that, firstly, the formula $A$ in Lemma~\ref{L:CQ-vs-QC-theses} is not necessarily a sentence, and secondly, a $\QH$-matrix is not necessarily a complete Heyting algebra, for greatest lower bounds and least upper bounds must exist but for some subsets of \textsf{B}.\\

A decisive step towards strengthening Mostowski’s approach was made in~\cite{henkin1950}. At first glance, this does not look like a strengthening, because Henkin used the calculus $\ImQ$, which is based on the schemata ax1Q--ax2Q and ax12Q--ax15Q and, therefore, is weaker than $\IntQ$. Henkin showed however that the implication of type~\eqref{E:historical-Q-con-1} can be enhanced by equivalence:
\begin{equation}\label{E:historical-Q-con-2}
\Vdash_{\ImQ}A~\Longleftrightarrow~\models_{\mathcal{M}^{\Q}_{2}}A,
\end{equation}
for any sentence $A$, where $\mathcal{M}^{\Q}_{2}$ is the class of $\Q$-matrices, where $\langle\textsf{B};\rightarrow,\one\rangle$ is a complete positive implication algebra; see~\cite{ras74}, chapter II, section 2.

Lemma~\ref{L:CQ-vs-QC-theses} does not discuss the relationship between $\Vdash_{\ImQ}$ and a consequence relation based on corresponding $\Q$-matrices. However, such an equivalence could be easily proved.
In this event, the difference between~\eqref{E:historical-Q-con-2} and  such an enhancement of the equivalence in Lemma~\ref{L:CQ-vs-QC-theses} would be that, firstly, in the enhanced lemma, $A$ could be any formula, and secondly, the $\Q$-matrices involved would be not necessarily complete implication algebras, but implication algebras where greatest lower bounds and least upper bounds but partially exist.

In the final section of~\cite{henkin1950}, it is shown how the main results can be extended to $\IntQ$, $\ClQ$ and a modal predicate calculus $\textsf{S4}\text{Q}$ (based on sentential modal system \textsf{S4}). However, the remark concerning the comparison of the equivalence of Lemma~\ref{L:CQ-vs-QC-theses} and the equivalences regarding these systems remains in force.\\

Independently of Henkin, H. Rasiowa established in~\cite{rasiowa1951} the following equivalences:
\[
\begin{array}{c}
\Vdash_{\IntQ}A~\Longleftrightarrow~~\models_{\mathcal{M}^{\Q}_{1}}A,~\text{for any sentence $A$};\\
\Vdash_{\textsf{S4}\text{Q}}A~\Longleftrightarrow~~\models_{\mathcal{M}^{\Q}_{3}}A,~\text{for any sentence $A$}.
\end{array}
\]
Here $\mathcal{M}^{\Q}_{3}$ is the class of complete topological Boolean algebras; see~\cite{ras74}, supplement. Rasiowa makes an important remark:
\begin{quote}
	``The assumption that the lattices under consideration are complete  was made only in order to ensure that all the least upper bounds and greatest lower bounds appearing in an interpretation of formulas exist. However, it may happen that, by a certain interpretation of $m$-argument predicates as mappings from $J^m$ into a suitable lattice $\fA$, all the least upper bounds and all the greatest lower bounds appearing in the interpretation of formulas exist in spite of the fact that $\fA$ is incomplete.'' \cite{ras74}, supplement, introduction
\end{quote}

There have been obtained numerous results, mostly by H. Rasiowa and R. Sikorski, according to which completeness theorems for the theses of the above-mentioned abstract logics, as well as for the predicate calculus of minimal logic and of constructive logic with strong negation, are proved for specific classes of $\Q$-structures. (See references in~\cite{ras74}, supplement.) In this regard, we would like, however, to quote D. Scott, who wrote:
\begin{quote}
	``The `first-order disease' is most plainly seen in the book by Rasiowa and Sikorski (1963). They had for a number of years considered Boolean-valued models of first-order sentences (as had Tarski, Mostowski, Halmos, and many others working on algebraic logic). Unfortunately they spent most of the time considering \emph{logical validity} (truth in all models) rather than the construction of possibly interesting \emph{particular} models.'' (D. Scott, Foreword to~\cite{bell2011})
\end{quote}
The models mentioned by Scott are discussed, for example, in~\cite{bell2011}; see also references there.\\

The concept of unary generalized quantifiers is rooted in the numerical quantifiers of A. Tarski; cf.~\cite{tarski1994}, section 20. A further explication of unary generalized quantifiers was given in~\cite{mostowski1957}. While the quantifiers of both Tarski and Mostowski operate with a set of two truth values, N. Rescher proposed the use of multivalued logic for a set of truth values; see~\cite{rescher1969}, chapter 2, section 29, and references therein.\\

Our definition of $\mathcal{M}^{\Q}$-consequence (definition~\eqref{E:vdash_M^Q}) is similar to a more restricted definition of consequence that was introduced in~\cite{malcev1973}, section 11.2. Namely, in his definition of defining collection of first-order atomic formulas, Mal'cev employed the notion of $\mathcal{M}^{\Q}$-consequence, $X\models_{\mathcal{M}^{\Q}}A$, in the framework of Example~\ref{example:ordinary-quantifiers} and for the case when X is a set of atomic formulas.

\chapter[Decidability]{Decidability}\label{chapter:decidability}
	
From a very general point of view, one can see the following two ways of defining abstract logic: (a) using a set of axioms or axiom schemata and a set of inference rules, that is, by defining a deductive system, and (b) using a class of logical matrices. In each case, the following questions naturally arise.

\begin{itemize}
	\item[(d.1)] Is there an effective procedure (in the sense of Section~\ref{section:prelimanaries-computability}) that, for a given finite deductive system, decides whether the abstract logic $\aLog$ defined by this system is consistent, that is, whether $\ThmS \neq \FormsL$?
	
	\item[(d.2)] Is there an effective procedure that for two given finite deductive systems decides whether these deductive systems are weakly equivalent, that is, do the abstract logics defined by these systems have the same set of theorems?
	
	\item[(d.3)] Is there an effective procedure that, for two given finite deductive systems, decides whether these deductive systems are equivalent, that is, do they define the same abstract logic?	
\end{itemize} 

Similarly, if abstract logics are defined by finite matrices, the following questions arise. 

\begin{itemize}
	\item[(m.1)] Is there an effective procedure, which, for a given finite matrix $\mat{M}$, decides whether the abstract logic $\aLog_{\mat{M}}$ defined by this matrix is nontrivial, that is, whether $L\mat{M}\neq \FormsL$?
	
	\item[(m.2)] Is there an effective procedure that for two given finite matrices decides whether these matrices are weakly equivalent, that is, whether the abstract logics defined by these matrices have the same set of theorems?
	
	\item[(m.3)] Is there an effective procedure that given two finite sets of finite $\Lan$-matrices decides whether these sets define the same abstract logic? In view of Proposition~\ref{P:matrix-con=atlas-con}, this problem is equivalent to the following problem: Is there an effective procedure that, given two finite $\Lan$-atlases, decides whether these atlases are equivalent, that is whether they define the same abstract logic? 	
\end{itemize} 

We do not ask whether a finite deductive system has theorems because to have theorems, it has to contain axioms or rules without premises. On the other hand, we do not ask whether an abstract logic defined by a finite matrix is nontrivial, for if the matrix contains at least one non-designated element, formula $p$, where $p$ is a variable, is not a theorem. In addition, in questions m.1 and m.2 we take a finite matrix and not a finite set of finite matrices because any finite set of finite matrices has the same set of theorems as their direct product.

Problems (m.1) and (m.2) can be formulated for finite atlases. Let us observe that if $\mathcal{M} = \langle{\alg{A},\D}\rangle$ is an atlas, then
\[
\bm{T}_{\aLog_{\mathcal{M}}}\! = \bm{T}_{\aLog_{\mat{M}}}=L\mat{M}, 
\]		
where $\mat{M} = \langle\alg{A},\bigcap\D\rangle$; that is, $\textbf{Cn}_{\aLog_{\mathcal{M}}}(\varnothing)=\textbf{Cn}_{\aLog_{\mat{M}}}(\varnothing)$, although the identity $\textbf{Cn}_{\aLog_{\mathcal{M}}}(X)\!=\textbf{Cn}_{\aLog_{\mat{M}}}(X)$, for all $X\subseteq\FormsL$, is not necessarily true.
Thus, the problems (m.1) and (m.2) for finite atlases can be reduced to finite matrices. \\

It was observed in \cite{kuznetsov1963} that the problems (d.1)--(d.3) are undecidable for a broad class of naturally defined deductive systems. On the other hand, the problems (m.1)--(m.3) have positive solutions for finite matrices: to the questions (m.1) and (m.2) are due to J. {\L}o{\'s} (cf. \cite[Propositions 25 and 27*]{los1949}) and J. Kalicki (cf. \cite{kalicki1948,kalicki1952}); the positive solutions to (m.3) were found in~\cite{citkin1975},~\cite{zygmunt1983} and~\cite{devyatkin2013}.

In this chapter we will look at solutions to the problems (m.1)--(m.3).

\section{Abstract logics defined by finite atlases}\label{section:logics-defined-by-atlases}

In this chapter, we suppose that a language $\Lan$ contains only a finite number of connectives and countably infinite set of variables. For each $n >0$, we denote by $\FormsLn$ the set of all $\Lan$-formulas containing variables only from the set $\{p_1,\dots,p_n\} \subseteq \Lan$. If $X \subseteq \FormsL$, we denote:
\[
\rton{X} := X \cap \FormsLn.
\]

In what follows, for abstract logics defined by atlases and the notions of indistinguishably and restricted Lindenbaum algebra will play a key role. 

\subsection{Indistinguishable formulas} 

\begin{defn}\label{def-indist}\index{indistinguishable formulas}\index{$\meqvA$}
	Given an algebra {\em$\alg{A}$} of type $\Lan$, formulas $\alpha,\beta \in \FormsL$ are {\em$\alg{A}$}-\textbf{indistinguishable}, in symbols {\em$\alpha\meqvA\beta$}, if
	{\em\begin{equation}
			v[\alpha] = v[\beta], \textit{ for every valuation } v \text{ in } \alg{A}. \label{eq-def-indist}
	\end{equation}}
\end{defn}

It is not hard to see that \alg{A}-indistinguishability relation~\eqref{eq-def-indist} is an equivalence relation on $\FormsL$. (Exercise~\ref{section:logics-defined-by-atlases}.\ref{EX:eq-def-indist}) 

\begin{prop}\label{pr-indtheor}
	Let {\em$\M = \lr{\alg{A},\D}$} be an $\Lan$-atlas. Then, for any {\em$\{\alpha, \alpha^\prime,\beta, \beta^\prime \} \cup X \subseteq \FormsL$}, whenever {\em$\alpha \meqvA \alpha^\prime$} and {\em$\beta \meqvA \beta^\prime$},
	\[
	X,\alpha \models_{\M} \beta \iff X,\alpha^\prime \models_{\M} \beta^\prime.
	\]
\end{prop}
\noindent\emph{Proof} immediately follows from definition of $\meqvA$ and definition of the consequence relation defined by an atlas. (Exercise~\ref{section:logics-defined-by-atlases}.\ref{EX:pr-indtheor})\\

Thus, if  $\M = \lr{\alg{A},\D}$ is an $\Lan$-atlas, then, for each $\aLog_{\mathcal{M}}$-theory $D$ and any formulas $\alpha, \beta$ such that $\alpha \meqvA \beta$, 
\begin{equation}
	\alpha \in D \iff \beta \in D. \label{eq-equth}
\end{equation}
That is, for each $\alpha \in \FormsL$,
\[
\alpha/\mmeqvA \ \subseteq D \  \text{ or } \alpha/\mmeqvA  \cap \ D =\varnothing,
\]
i.e. the $\meqvA$-equivalence classes (aka $\meqvA$-\textit{\textbf{cosets}}) do not cross the boundaries of $D$ (cf. Fig. \ref{fig-meqv}).\index{coset}

\begin{figure}[ht]
	\[
	\begin{array}{lcr}
		\ctdiagram{
			\ctnohead
			\cten 0,0,0,90:{}
			\cten 0,90,90,90:{}
			\cten 90,90,90,0:{}
			\cten 90,0,0,0:{}	
			\cten 0,70,90,70:{}	
			\ctnohead\ctdash
			\def\zzctdrawdashedge{\drawdashedge{1pt}{1pt}11}
			\cten 0,70,10,90:{}	
			\cten 10,70,20,90:{}	
			\cten 20,70,30,90:{}
			\cten 30,70,40,90:{}
			\cten 40,70,50,90:{}
			\cten 50,70,60,90:{}
			\cten 60,70,70,90:{}
			\cten 70,70,80,90:{}	
			\cten 80,70,90,90:{}	
			\ctv -10,79:{D}
			\ctv 38,-12:{\FormsL}
		}
		& \qquad
		&
		\ctdiagram{
			\ctnohead
			\cten 0,0,0,90:{}
			\cten 0,90,90,90:{}
			\cten 90,90,90,0:{}
			\cten 90,0,0,0:{}	
			\cten 0,70,90,70:{}	
			\ctnohead\ctdash
			\def\zzctdrawdashedge{\drawdashedge{1pt}{1pt}11}
			\cten 0,70,10,90:{}	
			\cten 10,70,20,90:{}	
			\cten 20,70,30,90:{}
			\cten 30,70,40,90:{}
			\cten 40,70,50,90:{}
			\cten 50,70,60,90:{}
			\cten 60,70,70,90:{}
			\cten 70,70,80,90:{}	
			\cten 80,70,90,90:{}	
			\ctnohead\ctdash
			\def\zzctdrawdashedge{\drawdashedge{3pt}{2.5pt}11}
			\cten 0,10,90,10:{} 	 	
			\cten 0,20,90,20:{} 	 	
			\cten 0,30,90,30:{} 	 
			\cten 0,50,90,50:{} 		
			\cten 0,40,90,40:{} 	 	
			\cten 0,60,90,60:{} 	 	
			\cten 0,80,90,80:{} 	 	
			\cten 10,0,10,90:{} 	 	
			\cten 20,0,20,90:{} 	 	
			\cten 30,0,30,90:{} 	 	
			\cten 40,0,40,90:{} 	 	
			\cten 50,0,50,90:{} 	 	
			\cten 60,0,60,90:{} 	 	
			\cten 70,0,70,90:{} 	 	
			\cten 80,0,80,90:{} 	 	
			\ctv -20,79:{D/\mmeqvA}
			\ctv 38,-12:{\FormsL/\mmeqvA}
		}
	\end{array}
	\]
	\caption{$\meqv{_\alg{A}}$-Equivalence partitioning}\label{fig-meqv}
\end{figure}

The notion of a depth of formula is instrumental in describing the complete sets of representatives of $\meqvA$.

\begin{defn}\index{formula!depth}
	Let $\alpha$ be an $\Lan$-formula. The \textbf{depth} of $\alpha$ $($denoted by $d(\alpha)$$)$ is a number defined inductively:
	$d(p) = 0$ for each variable $p$ of $\Lan$, $d(c) = 0$ for each constant $c$ of $\Lan$; also, if $F$ is an $n$-ary connective of $\Lan$ and $\alpha=F\alpha_1,\dots,\alpha_n$, then
	\[
	d(\alpha) = 1 + \max(d(\alpha_1),\dots,d(\alpha_n)).
	\]
\end{defn}

It is not hard to see that $d(\alpha)$ is a maximal nestedness of connectives in $\alpha$, or that $d(\alpha)$ is a length of the longest branch in the formula-tree of $\alpha$. (See Section~\ref{section:languages}.)

\begin{example}
	{\em	Let $\alg{A}$ be an arbitrary nontrivial Boolean algebra. Then, every formula $\alpha$ is $\meqvA$-equivalent to its disjunctive normal form $\beta$. It is not hard to see that $d(\beta) \leq 3$: variables have depth equal to 0, the negations of variables have depth equal to 1, the conjunctions of variables and negations of variables have depth equal to at most 2, and the disjunctions of the conjunctions of the variables and their negations have depth equal to at most 3. We note that, even though the depth of disjunctive normal form does not exceed 3, the  formula may be very long.} 	
\end{example}

\begin{lem}\label{L:pr-complsr}
	Suppose that {\em$\alpha,\beta \in \FormsL$} such that {\em$\alpha \meqvA \beta $}. Also, let $\gamma[\alpha]$ be a formula having an occurrence of $\alpha$, while $\gamma[\beta]$ be a formula obtained from $\gamma[\alpha]$ by replacing this occurrence of $\alpha$ with $\beta$. Then {\em$\gamma[\alpha] \meqvA \gamma[\beta]$}.
\end{lem}
\noindent\emph{Proof}~is left to the reader. (Exercise~\ref{section:logics-defined-by-atlases}.\ref{EX:lemma-pr-complsr})\\

\begin{prop}\label{pr-complsr}
	Let {\em$\alg{A}$} be an $\Lan$-algebra and {\em$X \subseteq \FormsL$} be a set of formulas of depth $< n$. If for each formula {\em$\alpha \in \FormsL$} of depth $\leq n$, there is formula $\alpha^{\prime} \in X$ such that {\em$\alpha \meqvA \alpha^{\prime}$}, then for any formula {\em$\beta \in \FormsL$} there is a formula $\beta^{\prime} \in X$ such that {\em$\beta \meqvA \beta^{\prime}$}. 
\end{prop}
\begin{proof}
	Let $\alpha$ be a formula of depth $m$. If $m \leq n$ the statement is trivial. Let $m = n + k$ and we prove the statement by induction on $k$.
	
	\textbf{Basis}. The case $k = 0$ is trivial.
	
	\textbf{Assumption}. Assume that for each formula $\alpha$ of depth $ \leq n + s$, there is a formula $\alpha^{\prime} \in X$ such that $\alpha \meqvA \alpha^{\prime}$.
	
	\textbf{Step}. Let $\alpha$ be a formula of depth $n + s +1$. Then for some $m$-ary connective $F \in \Lan$, $\alpha = F\alpha_1,\dots,\alpha_m$. By the induction assumption, for each $i \in \{1,\ldots,m\}$, because $d(\alpha_i)  \leq n + s$, there is a formula $\alpha_{i}^{\prime} \in X$ such that $\alpha_i \meqvA \alpha_{i}^{\prime}$, that is for every valuation $v$ in $\alg{A}$,  $v[\alpha_i] = v[\alpha_i^\prime] $. Then, in virtue of Lemma~\ref{L:pr-complsr}, for each valuation $v$ in $\alg{A}$,
	\[
	v[F\alpha_1,\dots,\alpha_m] = F(v[\alpha_1],\dots,v[\alpha_m]) = F(v[\alpha_{1}^{\prime}],\dots,v[\alpha^{\prime}_m]) = v[F\alpha^{\prime}_1,\dots,\alpha^{\prime}_n].
	\]
	and hence, $F\alpha_1,\dots,\alpha_m \meqvA F\alpha_{1}^{\prime},\dots,\alpha_{m}^{\prime}$. Let us observe that because all $\alpha_i'$ are in $X$, we have $d(\alpha_i') < n$ and hence, $d(F\alpha_{1}^{\prime},\dots,\alpha_{m}^{\prime}) \leq n$. By the assumption of the theorem, there is a formula $\beta \in X$ such that $F\alpha_{1}^{\prime},\dots,\alpha_{m}^{\prime} \meqvA \beta$ and, hence, $F\alpha_1,\dots,\alpha_m \meqvA \beta$. 
\end{proof}

Let \alg{A} be an algebra. Any set $\Forms_{\meqvA}\subseteq\FormsL$ satisfying the properties
\begin{itemize}
	\item $\bigcup\set{\alpha/\meqvA}{\alpha\in \Forms_{\meqvA}}=\FormsL$, and
	\item for all $\alpha,\beta\in \Forms_{\meqvA}$, $\alpha/\meqvA =\beta/\meqvA$ if, and only if, $\alpha=\beta$,
\end{itemize}
is called a \textbf{\textit{complete set of representatives}}\index{complete set of representatives} of $\FormsL/\meqvA$.\\

We observe the following.

\begin{prop}\label{pr-repr}
	Let {\em$\M_1 = \lr{\alg{A},\D_1}$} and {\em$\M_2 = \lr{\alg{A},\D_2}$} be two $\Lan$-atlases. And let {\em$\Forms_{\meqvA}$} be a complete set of representatives of {\em$\FormsL/\meqvA$}. Then the following properties hold:
	{\em
		\[
		\begin{array}{cl}
			(\text{a}) &\Thm_{\aLog_{\mathcal{M}_1}} = \Thm_{\aLog_{\mathcal{M}_2}},\textit{ if, and only if, }\Thm_{\aLog_{\mathcal{M}_1}} \cap \Forms_{\meqvA} = \Thm_{\aLog_{\mathcal{M}_2}}  \cap \Forms_{\meqvA};\\
			(\text{b}) &\vdash_{\aLog_{\mathcal{M}_1}} \subseteq\, \vdash_{\aLog_{\mathcal{M}_2}}~\textit{if, and only if, for each}~\alpha \cup X \subseteq \Forms_{\meqvA},~X \models_{\mathcal{M}_1} \alpha~\textit{entails}~X \models_{\mathcal{M}_2} \alpha.
		\end{array}
		\]}		
\end{prop} 
\noindent\emph{Proof}~is left to the reader. (Exercise~\ref{section:logics-defined-by-atlases}.\ref{EX:pr-repr})\\

\subsection{Matrix equivalences.} 
We recall the definition of matrix equivalence introduced in~\cite{zygmunt1983}. 

\begin{defn}\index{matrix!congrquence}\index{matrix!equivalence\index{atlas!congruence}}\index{atlas!equivalence}
	Suppose that {\em$\matM = \lr{\alg{A},D}$} is an $\Lan$-matrix. An equivalence relation $\sim$ on {\em$\alg{A}$} is a \textbf{matrix equivalence} $($on $\mat{M}$$)$ if for every {\em$a,b \in |\alg{A}|$}, $a \in D$  and $a \sim b$ entails $b \in D$. Respectively, a congruence $\theta$ on {\em$\alg{A}$} is a \textbf{matrix congruence} if it is a matrix equivalence. 
	Accordingly, given an $\Lan$-atlas {\em$\M = \lr{\alg{A},\D}$}, an equivalence relation $\sim$ on {\em$\alg{A}$} is an \textbf{atlas equivalence} if for each $D \in \D$ and every {\em$a,b \in |\alg{A}|$}, if $a \in D$ and $a \sim b$, then $b \in D$, and a congruence $\theta$ on {\em$\alg{A}$} is an \textbf{atlas congruence} if $\theta$ is an atlas equivalence.  
\end{defn}

It should be clear that if $\M = \lr{\alg{A},\D}$ is an $\Lan$-atlas, a relation $\sim$ on $\alg{A}$ is an atlas equivalence (or matrix congruence) if, and only if, $\sim$ is a matrix equivalence (or matrix congruence) for each matrix $\lr{\alg{A},D}$, where $D \in \D$.

The identity relation is an example of matrix congruence.\\

If $\aLog$ is an abstract logic, then an equivalence $\sim$ on $\FormsLn$ is a matrix equivalence if the following holds: For any $\alpha,\beta \in \FormsL$,
\[
\alpha \sim \beta  \implies \alpha \vdashS \beta \text{ and } \beta \vdashS \alpha.
\]  
It is not hard to see that relation $\sim$ defined as 
\[
\alpha \sim \beta \define \alpha \vdash_S \beta \text{ and } \beta \vdash_S \alpha
\]
is a matrix equivalence on $\LinS$, that is, on the Lindenbaum matrix of $\aLog$.\\

Suppose that $\matM = \lr{\alg{A},D}$ is an $\Lan$-matrix. Consider the partition $\{D, |A| \setminus D\}$. This partition defines a matrix equivalence which obviously is the greatest matrix equivalence on $\matM$. If $\M = \lr{\alg{A},\D}$ is an $\Lan$-atlas and $\sim_D$ is a greatest matrix equivalence on $\lr{\alg{A},D}$, where $D \in \D$, then  $\bigcap\set{\sim_D}{D \in \D}$ is the greatest atlas equivalence on $\M$. (See Exercise~\ref{section:logics-defined-by-atlases}.\ref{EX:greates-atlas-equivalence}.)  

Thus, each $\Lan$-atlas ($\Lan$-matrix) has the greatest atlas equivalence (matrix equivalence). Let us note that this equivalence needs not to be a congruence. Nevertheless, the following holds. 

\begin{prop}[\cite{brown-suszko1973}, theorem 8] \label{pr-greqv}
	Each $\Lan$-atlas has the greatest atlas congruence $($matrix congruence$)$.
\end{prop}
\begin{proof}
	Suppose that $\M = \lr{\alg{A},\D}$. Recall that the class of all congruences on $\alg{A}$ is a complete lattice in which the join $\bigvee_{i \in I}\theta_i$ of congruences is defined as follows (cf. \cite{burris-sankapp1981}, theorem 4.7):
	\[
	\begin{aligned}
		& (a,b) \in \bigvee_{i \in I}\theta_i \iff\text{ there are } a_0,\dots,a_{n} \in \alg{A} \text{ such that } \\
		& a = a_0, a_n = b \text{ and } (a_0,a_1) \in \theta_{i_1}, \dots, (a_{n-1},a_n) \in \theta_{n}.  
	\end{aligned}
	\]
	It is not hard to see that if $\theta_i$ are atlas (or matrix) congruences, then their join is an atlas (or matrix) congruence. Thus, the join of all atlas (or matrix) congruences is the greatest atlas (or matrix) equivalence. 
\end{proof}

\begin{prop}\label{pr-matcong} Let {\em$\matM = \lr{\alg{A},D}$} be a logical matrix. Then relation {\em$\meqvA$} is a matrix congruence on the Lindenbaum matrix {\em$\LinM = \lr{\FormsL, \textbf{Cn}_{\aLog_{\mat{M}}}(\varnothing)}$}. 
\end{prop}
\begin{proof}
	It should be clear that $\meqvA$ is an equivalence relation. Moreover, if $\alpha \in \ThmM=\textbf{Cn}_{\aLog_{\mat{M}}}(\varnothing)$, then  $v[\alpha] \in D$, for all valuations $v$. If $\beta \in \FormsL$ and $\beta \meqvA \alpha$, then, by definition of $\meqvA$, $v[\beta] = v[\alpha]$ and, hence, $v[\beta] \in D$, for all valuations $v$. Thus, $\beta \in \ThmM$ and $\meqvA$ is a matrix congruence.
	
	Suppose that $F$ is an $n$-ary connective, $\alpha_1,\beta_1\dots,\alpha_n,\beta_n, \in \FormsL$ and  $\alpha_i \meqvA \beta_i$ for $i \in \{1,\ldots,n\}$. Then, by definition of $\meqvA$, for each valuation $v$, $v[\alpha_i] = v[\beta_i]$ for all $i \in\{1,\ldots,n\}$ and, hence,
	\[
	v[F\alpha_1,\dots,\alpha_n] = F(v[\alpha_1],\dots,v[\alpha_n]) = F(v[\beta_1],\dots,v[\beta_n]) = v[F\beta_1,\dots,\beta_n].   
	\]
	and hence, $F\alpha_1,\dots,\alpha_n \meqvA F\beta_1,\dots,\beta_n$.
\end{proof}

\begin{cor}\label{pr-atlcong} Let {\em$\M = \lr{\alg{A},\D}$} be an $\Lan$-atlas. Then {\em$\meqvA$} is an atlas congruence on the Lindenbaum atlas of the abstract logic defined by $\M$. 
\end{cor}

\subsection{Restricted Lindenbaum matrices and atlases}\label{sec-restr}

Let $\FormsLn$ denote the subset of all formulas from $\FormsL$ containing  only the variables $p_1,\dots,p_n$. 

Suppose that $\aLog$ is an abstract logic and $\LinS = \lr{\FormsL,\theory}$ is its Lindenbaum atlas. Then, for any $n > 0$, one can define a \textbf{\textit{restricted Lindenbaum atlas}}:\index{matrix!Lindenbaum atlas!restricted}
\[
\LinnS := \lr{\FormsLn, \Sgn}, \text{ where } \Sgn = \set{D \cap \FormsLn}{T \in \theory}.
\]

Let us observe that if $T$ is an $\aLog$-theory, that is $T \in \theory$, then the set $\rton{T}$ is closed in $\FormsLn$ relative to $\vdashS$ in the following sense:
\[
\text{for any } \alpha \cup X \subseteq \FormsLn,\,  X \vdash_S \alpha \text{ and } X \subseteq \rton{T} \implies \alpha \in \rton{T}.
\]
Moreover, if $\aLog$ is structural and $\sigma$ is a substitution mapping all variables into $\FormsLn$, then for any $\aLog$-theory $T$, $\sigma(T) \subseteq \FormsLn$ and is closed in $\FormsLn$ relative to $\vdashS$.\\

In this chapter, dealing chiefly with abstract logics determined by $\Lan$-atlases and matrices, we use the following notation:
\begin{itemize}
	\item $\theoryMM: = \Sigma_{\aLog_{\mathcal{M}}}$ and $\theoryM:=\Sigma_{\aLog_{\mat{M}}}$, where $\M$ is an $\Lan$-atlas and, respectively, $\matM$ is an $\Lan$-matrix;
	\item $\LinMM:=\lr{\FormsL,\theoryMM}$ and $\LinM:=\lr{\FormsL,\theoryM}$;
	\item $\theorynMM:=\set{T\cap\FormsLn}{T\in\theoryMM}$ and $\theorynM:=\set{T\cap\FormsLn}{T\in\theoryM}$, where $\M$ is an $\Lan$-atlas and, respectively, $\matM$ is an $\Lan$-matrix; 
	\item $\LinnMM:=\lr{\FormsLn,\theorynMM}$ and $\LinnM:=\lr{\FormsLn,\theorynM}$;
	\item $\Thm[\matM]:=\textbf{Cn}_{\aLog_{\mat{M}}}(\varnothing)$ and 
	$\Thmn[\matM]:=\Thm[\matM]\cap\FormsLn$, where $\matM$ is an $\Lan$-matrix;\footnote{Although $\Thm[\matM]=L\matM$ (the logic of $\matM$), however, it will be convenient to use the notations $\Thm[\matM]$ and $\Thmn[\matM]$ in the same context.}
	\item given two $\Lan$-atlases $\M_1$ and $\M_2$,
	\[
	\aLog[\M_1]\preccurlyeq\aLog[\M_2]~\define~\models_{\M_1}\subseteq\,\models_{\M_2};
	\]
	accordingly,
	\[
	\aLog[\M_1]=\aLog[\M_2]~\define~\aLog[\M_1]\preccurlyeq\aLog[\M_2]~\text{and}~\aLog[\M_2]\preccurlyeq\aLog[\M_1];
	\]
	\item given two $\Lan$-matrices $\matM_1$ and $\matM_2$,
	\[
	\aLog[\matM_1]\preccurlyeq\aLog[\matM_2]~\define~\models_{\matM_1}\subseteq\,\models_{\matM_2};
	\]
	accordingly,
	\[
	\aLog[\matM_1]=\aLog[\matM_2]~\define~\aLog[\matM_1]\preccurlyeq\aLog[\matM_2]~\text{and}~\aLog[\matM_2]\preccurlyeq\aLog[\matM_1].
	\]
\end{itemize}

The notion of restricted Lindenbaum atlas is just an example of a more general concept.
\begin{defn}\label{D:subatlas-submartix}\index{subatlas}\index{submatrix}
	Let {\em$\M_1=\lr{\alg{A}_1,\D_1}$} and {\em$\M_2=\lr{\alg{A}_2,\D_2}$} be $\Lan$-atlases. $\M_1$ is called a \textbf{subatlas} of $\M_2$ if {\em$\alg{A}_1$} is a subalgebra of {\em$\alg{A}_2$} and {\em$\D_1=\set{D\cap|\alg{A}_1|}{D\in\D_2}$}. In particular, given $\Lan$-matrices {\em$\mat{M}_1=\lr{\alg{A}_1,D_1}$} and {\em$\mat{M}_2=\lr{\alg{A}_2,D_2}$}, {\em$\mat{M}_1$} is a \textbf{submatrix} of {\em$\mat{M}_2$} if {\em$\alg{A}_1$} is a subalgebra of {\em$\alg{A}_2$} and {\em$D_1=D_2\cap|\alg{A}_1|$}.
\end{defn}

Thus, in the terminology of Definition~\ref{D:subatlas-submartix}, $\LinnMM$ is a subatlas of $\LinMM$ and $\LinnM$ is a submatrix of $\LinM$.

\begin{prop}\label{P:subatlas-property}
	Let {\em$\M_1=\lr{\alg{A}_1,\D_1}$} be a subatlas of {\em$\M_2=\lr{\alg{A}_2,\D_2}$}. Then $\aLog[\M_2]\subseteq\,\aLog[\M_1]$. In particular, if an $\Lan$-matrix $\mat{M}_1$ is a submatrix of a matrix $\mat{M}_2$, then $\Thm[\mat{M}_2]\subseteq\Thm[\mat{M}_1]$.
\end{prop}
\begin{proof}
	Suppose $X\models_{\M_2}\alpha$. Let $v$ be a valuation in $\alg{A}_1$ such that $v[X]\subseteq D$, where $D\in\D_1$. By definition, there is a logical filter $D^\prime\in\D_2$ such that $D=D^{\prime}\cap\alg{A}_1$. Thus, $v[X]\subseteq D^{\prime}$. Hence, by premise, $v[\alpha]\in D^\prime$, which implies that $v[\alpha]\in D$.
	
	The second statement is obvious.
\end{proof}

An obvious consequence of Proposition~\ref{P:subatlas-property} is the following.
\begin{cor} \label{pr-atlincl}
	Let $\M$ be an $\Lan$-atlas. Then {\em$\aLog[\M]\,=\,\aLog[\LinMM]\preccurlyeq\,\aLog[\LinnMM]$}. In particular, {\em$\Thm[\LinM]\subseteq\Thm[\LinnM]$}.
\end{cor}
\begin{proof}
	The first equality follows from Corollary~\ref{C:lindenbaum-matrix}, the first inclusion from Proposition~\ref{P:subatlas-property}.
\end{proof}

We note that every atlas congruence on $\LinMM$ (matrix congruence on $\LinM$) can be restricted to $\LinnMM$ (respectively, to $\LinnM$). 

Let us observe that for any $\Lan$-matrix $\matM = \lr{\alg{A},D}$, relation $\meqvA$ depends only on $\alg{A}$ and does not depend on $D$ and hence, Proposition \ref{pr-matcong} yields that for any $D \subseteq |\alg{A}|$,  $\meqvA$ is a matrix congruence on Lindenbaum matrix $\lr{\FormsL,\Thm[\matM]}$. Thus, the following holds.

\begin{cor}\label{cor-mateqvincl}
	Let {\em$\matM_1 = \lr{\alg{A},D_1}$} and {\em$\matM_1 = \lr{\alg{A},D_2}$} be $\Lan$-matrices. Then
	\[
	D_1 \subseteq D_2 \implies \Thm[\matM_1]/\mmeqvA \ \subseteq \ \Thm[\matM_2]/\mmeqvA.
	\]
\end{cor}
\begin{proof}
	Indeed, if $D_1 \subseteq D_2$, then $\Thm[\matM_1] \subseteq \Thm[\matM_2]$ and because $\meqvA$ is a matrix equivalence and the $\meqvA$-cosets do not cross boundaries of $\Thm[\matM_1]$ and $\Thm[\matM_2]$ (cf. Fig. \ref{fig-meqv2}),   
	$\Thm[\matM_1]/\mmeqvA \ \subseteq \ \Thm[\matM_2]/\mmeqvA$.
	
	\begin{figure}[ht]
		\[
		\begin{array}{lcr}
			\ctdiagram{
				\ctnohead
				\cten 0,0,0,90:{}
				\cten 0,90,90,90:{}
				\cten 90,90,90,0:{}
				\cten 90,0,0,0:{}	
				\cten 0,70,90,70:{}	
				\cten 0,50,90,50:{}	
				\ctnohead\ctdash
				\def\zzctdrawdashedge{\drawdashedge{1pt}{1pt}11}
				\cten 0,70,10,90:{}	
				\cten 10,70,20,90:{}	
				\cten 20,70,30,90:{}
				\cten 30,70,40,90:{}
				\cten 40,70,50,90:{}
				\cten 50,70,60,90:{}
				\cten 60,70,70,90:{}
				\cten 70,70,80,90:{}	
				\cten 80,70,90,90:{}	
				\cten 0,90,10,50:{} 	
				\cten 10,90,20,50:{} 
				\cten 20,90,30,50:{} 
				\cten 30,90,40,50:{} 
				\cten 40,90,50,50:{} 
				\cten 50,90,60,50:{} 	
				\cten 60,90,70,50:{} 
				\cten 70,90,80,50:{} 
				\cten 80,90,90,50:{} 	
				\ctv -30,79:{\Thmn[\matM_2]}
				\ctv -30,59:{\Thmn[\matM_1]}
				\ctv 38,-12:{\FormsLn}
			}
			&
			&
			\ctdiagram{
				\ctnohead
				\cten 0,0,0,90:{}
				\cten 0,90,90,90:{}
				\cten 90,90,90,0:{}
				\cten 90,0,0,0:{}	
				\cten 0,70,90,70:{}	
				\cten 0,50,90,50:{}	
				\ctnohead\ctdash
				\def\zzctdrawdashedge{\drawdashedge{1pt}{1pt}11}
				\cten 0,70,10,90:{}	
				\cten 10,70,20,90:{}	
				\cten 20,70,30,90:{}
				\cten 30,70,40,90:{}
				\cten 40,70,50,90:{}
				\cten 50,70,60,90:{}
				\cten 60,70,70,90:{}
				\cten 70,70,80,90:{}	
				\cten 80,70,90,90:{}	
				\cten 0,90,10,50:{} 	
				\cten 10,90,20,50:{} 
				\cten 20,90,30,50:{} 
				\cten 30,90,40,50:{} 
				\cten 40,90,50,50:{} 
				\cten 50,90,60,50:{} 	
				\cten 60,90,70,50:{} 
				\cten 70,90,80,50:{} 
				\cten 80,90,90,50:{} 	
				\ctnohead\ctdash
				\def\zzctdrawdashedge{\drawdashedge{3pt}{2.5pt}11}
				\cten 0,10,90,10:{} 	 	
				\cten 0,20,90,20:{} 	 	
				\cten 0,30,90,30:{} 	 	
				\cten 0,40,90,40:{} 	 	
				\cten 0,60,90,60:{} 	 	
				\cten 0,80,90,80:{} 	 	
				\cten 10,0,10,90:{} 	 	
				\cten 20,0,20,90:{} 	 	
				\cten 30,0,30,90:{} 	 	
				\cten 40,0,40,90:{} 	 	
				\cten 50,0,50,90:{} 	 	
				\cten 60,0,60,90:{} 	 	
				\cten 70,0,70,90:{} 	 	
				\cten 80,0,80,90:{} 	 	
				\ctv -40,79:{\Thmn[\matM_2]/\!{\meqvA}}
				\ctv -40,59:{\Thmn[\matM_1]/\!{\meqvA}}
				\ctv 38,-12:{\FormsLn/\!{\meqvA}}
			}
		\end{array}
		\]
		\caption{$\meqv{_\alg{A}}$-Equivalence partitioning}\label{fig-meqv2}
	\end{figure}
\end{proof}

In some cases, we the equality $\aLog[\LinMM]=\aLog[\LinnMM]$ may hold.

\begin{figure}[ht]
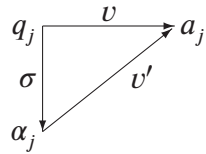

	\[
	\ctdiagram{
		\ctet 0,40,50,40:{v}
		\ctel 0,40,0,0:{\sigma}
		\cter 0,0,49,39:{v^{\prime}}		
		\ctv -7,40:{q_j}
		\ctv 57,40:{a_j}
		\ctv -7,0:{\alpha_j}
	}
	\]
	\caption{Refuting valuation}\label{fig-pr-redton}
\end{figure}

\begin{prop} \label{pr-atleq}
	Let {\em$\M = \lr{\alg{A},\D}$} be an $\Lan$-atlas and {\em$\alg{A}$} be generated by elements $a_1,\dots,a_n$. Then {\em$\aLog[\LinMM]=\aLog[\LinnMM]$}. 
\end{prop}
\begin{proof}
	In view of Corollary~\ref{pr-atlincl}, it remains to prove that $\aLog[\LinnMM]\preccurlyeq\aLog[\LinMM]$, that is, for any set $X\cup\{\alpha\}\subseteq\FormsL$, 
	\[
	X \models_{\LinnMM} \alpha \implies X \models_{\LinMM} \alpha; 
	\] 
	and, in virtue of of the first equality of Corollary~\ref{pr-atlincl}, it suffices to show that
	\[
	X \models_{\LinnMM} \alpha \implies X \models_{\M} \alpha.
	\]
	
	For contradiction, assume that there is a valuation $v$ in $\alg{A}$ such that
	\[
	v[X] \subseteq D \text{ and } v[\alpha] \notin D, \text{ for some } D \in \D.
	\]
	We aim to prove that there is a substitution $\sigma$ in $\LinnMM$ such that
	\[
	\sigma(X) \subseteq T^{(n)} \text{ while } \sigma(\alpha) \notin T^{(n)}, \text{ for some theory } T^{(n)} \in \Sgn.
	\]
	
	Indeed, valuations in $\LinnMM$ are substitutions mapping the variables of $\VarL$ into $\FormsLn$. Assume that $\alpha$ contains variables only from a set $\{q_1,\dots,q_m\}$. Then, because $\alg{A}$ is generated by elements $a_1,\dots,a_n$, there are formulas $\alpha_j$ in variables $p_1,\dots,p_n$ and such that for every $j=1,\dots,m$,
	\[
	v[q_i] = \alpha_j[a_1,\dots,a_n].
	\]	
	Let us consider a substitution $\sigma$:
	\[
	\sigma(q) = \begin{cases}
		\alpha_j \text{ if } q = q_j,\\
		\alpha_1 \text{ otherwise}.
	\end{cases}
	\]
	
	It should be clear that $\sigma$ is a valuation in $\LinnMM = \lr{\FormsLn, \Sgn}$. 
	
	Let 
	\[
	T := \set{\beta \in \FormsL}{\sigma[X] \vdash_{S[\M]} \beta}.
	\]
	Then, $\rton{T} \in  \Sgn$ and $\sigma(X) \subseteq \rton{T}$. We aim to show that $\sigma(\alpha) \notin \rton{T}$.
	
	Let $v^\prime$ be a valuation in $\alg{A}$ such that $v^{\prime}: p_k \mapsto a_k$, where $k=1,\dots,n$. We observe that $v = v^{\prime}\circ\sigma$ (cf. Fig. \ref{fig-pr-redton}). 	Therefore, 
	\[
	v^{\prime}[\sigma(\alpha)] = v[\alpha] \notin D,
	\]
	while for each $\beta \in X$,
	\[
	v^{\prime}[\sigma(\beta)] = v[\beta] \in D.	
	\]
	Hence, $\sigma (X) \models_{\M} \sigma(\alpha)$ and consequently, $\sigma(\alpha) \notin \rton{T}$.
\end{proof}

We observe that it follows from Proposition \ref{pr-atleq} that if a algebraic carrier of an $\Lan$-atlas $\M = \lr{\alg{A},\D}$ is generated by $n$ elements, then the logic $\aLog_{\mathcal{M}}$ coincides with the abstract logic determined by the Lindenbaum atlas $\LinnMM$. Thus, if an algebra $\alg{A}$ is finite, say of $\card{\alg{A}} = n$, the logic $\aLog_{\mathcal{M}}$ coincides with $\aLog_{\LinnMM}$. This conclusion is important because, as we will see in Proposition~\ref{th-mefini}, in contrast to the cardinality of $\LinMM$, the cardinality of $\LinnMM$ is finite.

\begin{cor}\label{cor-logn}
	Suppose that {\em$\M = \lr{\alg{A},\D}$} is a finite $\Lan$-atlas of cardinality $n$. Then for every $m \ge n$, {\em$\aLog[\M] = \aLog[\Lin^{(m)}_{\M}]$}.
\end{cor}
\noindent\emph{Proof}~is left to the reader. (Exercise~\ref{section:logics-defined-by-atlases}.\ref{EX:cor-logn})\\

\begin{cor}\label{cor-theorn}
	Suppose that {\em$\matM = \lr{\alg{A},D}$} is a finite $\Lan$-matrix of cardinality $n$. Then
	{\em$\Thm[\matM] = \Thm[\LinnM]$}.
\end{cor}

We also note the following corollaries.
\begin{cor}\label{cor-redtongen}
	Suppose that {\em$\matM_1 = \lr{\alg{A},D_1}$} and {\em$\matM_2 = \lr{\alg{A},D_2}$} are  $\Lan$-matrices and algebra {\em$\alg{A}$} is generated by $n$ elements $a_1,\dots,a_n$. Then
	\[
	\Thm[\matM_1] = \Thm[\matM_2] \iff \Thmn[\matM_1] = \Thmn[\matM_2]. 
	\] 
\end{cor}
\noindent\emph{Proof}~is left to the reader. (Exercise~\ref{section:logics-defined-by-atlases}.\ref{EX:cor-redtongen})

\begin{cor}\label{cor-redtn}
	Suppose that {\em$\matM_1 = \lr{\alg{A},D_1}$} and {\em$\matM_2 = \lr{\alg{A},D_2}$} are finite $\Lan$-matrices of cardinality $n$. Then
	\[
	\Thm[\matM_1] = \Thm[\matM_2] \iff \Thmn[\matM_1] = \Thmn[\matM_2]. 
	\] 
\end{cor}
\noindent\emph{Proof}~is left to the reader. (Exercise~\ref{section:logics-defined-by-atlases}.\ref{EX:cor-redtn})\\


A very important property of equivalence $\meqvA$, which will be used for constructing an effective procedure, is formulated in the following proposition.

\begin{prop}\label{th-mefini}
	For any finite $\Lan$-atlas {\em$\M = \lr{\alg{A}, \D}$} and any natural number $n >0$, the set {\em$\FormsLn/\mmeqvA$} is finite.
\end{prop}
\begin{proof}
	We need to prove that $\meqvA$ partitions set $\FormsLn$ in a finite number of equivalence classes. The proof follows from the observation that each formula $\alpha \in \FormsLn$ can be regarded as a function on $\alg{A}$ in $n$ variables:
	\[
	f_\alpha(a_1,\dots,a_n) = v[\alpha],
	\]
	where $v$ is a valuation with $v: p_i \mapsto a_i$. From the definition of $\meqvA$, it follows that for any $\alpha,\beta \in \FormsLn$,
	\[
	\alpha \meqvA \beta \iff f_\alpha = f_\beta.
	\] 
	Hence, the cardinality of $\FormsLn/\mmeqvA$ does not exceed the cardinality of the set of all $n$-ary functions on $\alg{A}$, that is $\card{\FormsLn/\mmeqvA} \leq k^{k^n}$, where $k = \card{\alg{A}}$.	
\end{proof}

Suppose that $\sim$ is an equivalence relation on a set $A$. A subset of $A$ containing exactly one element from each $\sim$-equivalence class is called a  \textit{\textbf{complete set of representatives}} of $A/\!\!\sim$. Proposition \ref{th-mefini} entails that if $\matM = \lr{\alg{A},D}$ is a finite $\Lan$-matrix, then a complete set of representatives of $\FormsLn/\mmeqvA$ is finite. Moreover, there is an effective procedure that, given a finite $\Lan$-matrix or $\Lan$-atlas, delivers a complete set of representatives of $\FormsLn/\mmeqvA$.

Indeed, let us observe that we can enumerate $\FormsLn$; see Section~\ref{section:word-sets} and Section~\ref{section:enumerable-sets}. To that end, we enumerate the formulas in order of growing depth: we take all formulas of depth 0, and then we construct all formulas of depth $n$ based on formulas of lesser depth.   

Because algebra $\alg{A}$ is finite, for any two formulas $\alpha$ and $\beta$, we can easily determine whether $\alpha \meqvA \beta$ or not. Thus, after obtaining a new formula during enumeration of $\FormsLn$, we decide whether this formula is $\meqvA$-equivalent to a preceding formula, and if it is, we omit this formula, otherwise we retain this formula as a representative of the respective $\meqvA$-coset. And we stop when we determine that such a process does not produce any new representatives. Proposition \ref{pr-complsr} gives us a way to determine when to stop: if we did not retain any of formulas of depth $n+1$, we have obtain representatives of all $\meqvA$-cosets.   

\paragraph{Exercises~\ref{section:logics-defined-by-atlases}}
\begin{enumerate}
	\item \label{EX:eq-def-indist}Prove that that the relation $\meqvA$ of $\alg{A}$-indistinguishability is an equivalence on $\FormsL$.
	\item\label{EX:pr-indtheor}Prove Proposition~\ref{pr-indtheor}.
	\item\label{EX:lemma-pr-complsr} Suppose that $\alpha,\beta \in \FormsL$ such that $\alpha \meqvA \beta $. Also, let $\gamma[\alpha]$ be a formula having an occurrence of $\alpha$, while $\gamma[\beta]$ be a formula obtained from $\gamma[\alpha]$ by replacing this occurrence of $\alpha$ with $\beta$. Then $\gamma[\alpha] \meqvA \gamma[\beta]$ (Lemma~\ref{L:pr-complsr}). 
	\item\label{EX:pr-repr}Prove Proposition~\ref{pr-repr}.
	\item\label{EX:greates-atlas-equivalence}Prove that if $\M = \lr{\alg{A},\D}$ is an $\Lan$-atlas and $\sim_D$ is a greatest matrix equivalence on $\lr{\alg{A},D}$, where $D \in \D$, then  $\bigcap\set{\sim_D}{D \in \D}$ is the greatest atlas equivalence on $\M$.
	\item\label{EX:cor-logn} Prove Corollary~\ref{cor-logn}.
	\item\label{EX:cor-redtongen}Prove Corollary~\ref{cor-redtongen}.
	\item\label{EX:cor-redtn}Prove Corollary~\ref{cor-redtn}.
\end{enumerate}

\section{Enumerating procedure}\label{section:enumerating-procedure}
We consider the formulas of $\FormsLn$ as words of a finite alphabet (Section~\ref{section:word-sets}). Actually, the set $\FormsLn$ is an effectively decidable word set. (Exercises~\ref{section:enumerating-procedure}.\ref{EX:one})

Now we aim to show that, given a finite algebra $\alg{A}$, a complete set 
of representatives of $\FormsLn/\mmeqvA$ can be effectively generated. 

If \alg{A} has one element, then all formulas of $\FormsLn$ are equivalent to each other with respect to ~$\mmeqvA$. Thus a set consisting of any formula of $\FormsLn$ is a complete set of representatives.

Now we assume that \alg{A} has more than one element. We describe a generating procedure inductively.

\textsf{Step} $1$. We include all variables of $\FormsLn$ as representatives of $\FormsLn/\mmeqvA$.  Since $\card{\alg{A}} > 1$, for any two distinct variables $p_i,p_j \in \FormsLn$, there is a valuation that maps $p_i$ in one element of $\alg{A}$ and $p_j$ in a distinct element of $\alg{A}$, that is, $p_i \meqvA p_j$ does not hold. 

If  two $\Lan$-constants coincide in $\alg{A}$, we include only one of them as a representative of $\FormsLn/\mmeqvA$. 

\textsf{Step} $n$. For each $k$-ary connective $F \in \Lan$, we construct all formulas $F\alpha_1,\dots,\alpha_n$, where each $\alpha_i$ has previously been selected as a representative and at least one of $\alpha_i$ has depth $n-1$, i.e. $\alpha_i$ has been selected on the Step $n$. Then we compare the formula $F\alpha_1,\dots,\alpha_n$ with all previously selected representatives. If $F\alpha_1,\dots,\alpha_n$ is $\meqvA$-equivalent to one of them, we skip it; otherwise, we add $F\alpha_1,\dots,\alpha_n$ as a representative and thus extending the list.

If we have not selected a new representative on the Step $n$, then, by Proposition \ref{pr-complsr}, we have already selected representatives of all $\meqvA$-cosets and thus obtained a complete set of representatives.  

\begin{example}\label{ex-tr}
	{\em	Let $\alg{A} = \lr{\{\zero,\tb,\one\}; \lor, \zero}$, where $\lor$ is defined as in $\godelThree$ (Section~\ref{section:goedel}), and let $\matM = \lr{\alg{A},\{\one\}}$.
		
		To obtain a complete set of representatives of $\FormsL^{(1)}/\mmeqvA$, we do the following:
		
		\[
		\begin{array}{l|l|l}
			\text{depth} & \text{formulas} & \text{representatives}\\
			0 & p, \zero & p, \zero\\
			1 & p \lor p, p \lor \zero,\zero \lor p, \zero \lor \zero  & \\
		\end{array}
		\]
		
		It is not hard to see that neither formula $p$ or $\zero$ is valid in $\matM$ and  thus,the logic of $\matM$ is trivial.}
\end{example}

\paragraph{Exercises~\ref{section:enumerating-procedure}}
\begin{enumerate}
	\item \label{EX:one}Show that the set $\FormsLn$ is effectively decidable.
	\item Using Proposition~\ref{P:nishimura}, prove that the formulas defined in~\eqref{E:powers_p^n} form a complete set of representatives of $\FormsL^{(1)}/\mmeqvA$, where $\alg{A}$ is a Rieger-Nishimura algebra (Section~\ref{section:rieger-nichimura-algebra}).
\end{enumerate}

\section{Solution to Problem (m.1)} \label{sec-procetriv}

We need to demonstrate that there is an effective procedure that, given a finite matrix $\matM$, determines whether $\Thm[\matM] \neq\varnothing$.

The idea behind such a procedure is this: given a finite $\Lan$-matrix $\matM = \lr{\alg{A},D}$, we enumerate the formulas of $\FormsL$, and for each formula, we test whether it belongs to $\Thm[\matM]$, that is, we check whether the formula is valid in $\matM$. If we have found a valid formula, $\Thm[\matM] \neq \varnothing$ and we stop. The problem arises if $\Thm[\matM] = \varnothing$, because in this case, the procedure does not stop. To alleviate this problem, we do the following.

\begin{itemize}
	\item[(a)] Firstly, we reduce the problem to formulas in one variable: because set $\Thm[\matM]$ is closed relative to substitutions, if $\alpha \in \Thm[\matM]$, then $\alpha^{(p)} \in \Thm[\matM]$, where $\alpha^{(p)}$ is a formula obtained from $\alpha$ by substituting $p$ for each variable. That is, $\Thm[\matM] = \varnothing$ if and only if $\Thm^{(1)}[\matM] = \varnothing$, and we need to list only formulas from $\FormsL^{(1)}$, although, if $\Lan$ contains connectives, $\FormsL^{(1)}$ is still an infinite set and we need to make further reductions.	
	\item[(b)] Since the algebra $\alg{A}$ is finite, by Proposition \ref{th-mefini}, the set $\FormsL^{(1)}/\meqvA$ is finite. In addition, if one of the representatives of $\meqvA$-cosets is in $\Thm[\matM]$, then all of them are in $\Thm[\matM]$. Hence, it suffices to test just one of the representatives from each $\meqvA$-coset. 	
	\item[(c)] All we need is to be able to enumerate the representatives of  $\meqvA$-cosets: one for each class, and we will construct a procedure that does precisely this. Then, for each representative $\alpha$ we can check whether $\alpha$ is valid in $\matM$; and if so, clearly $\Thm[\matM] \neq \varnothing$, or we will exhaust the list of representatives and conclude that $\Thm[\matM] = \varnothing$.  
\end{itemize}

\begin{figure}[ht]
	\[
	\begin{array}{lcr}
		\ctdiagram{
			\ctnohead
			\cten 0,0,0,90:{}
			\cten 0,90,90,90:{}
			\cten 90,90,90,0:{}
			\cten 90,0,0,0:{}	
			\cten 0,70,90,70:{}	
			\ctnohead\ctdash
			\def\zzctdrawdashedge{\drawdashedge{1pt}{1pt}11}
			\cten 0,70,10,90:{}	
			\cten 10,70,20,90:{}	
			\cten 20,70,30,90:{}
			\cten 30,70,40,90:{}
			\cten 40,70,50,90:{}
			\cten 50,70,60,90:{}
			\cten 60,70,70,90:{}
			\cten 70,70,80,90:{}	
			\cten 80,70,90,90:{}	
			\ctv -30,79:{\Thm^{(1)}[\matM]}
			\ctv 38,-12:{\FormsL^{(1)}}
		}
		&
		&
		\ctdiagram{
			\ctnohead
			\cten 0,0,0,90:{}
			\cten 0,90,90,90:{}
			\cten 90,90,90,0:{}
			\cten 90,0,0,0:{}	
			\cten 0,70,90,70:{}	
			\ctnohead\ctdash
			\def\zzctdrawdashedge{\drawdashedge{1pt}{1pt}11}
			\cten 0,70,10,90:{}	
			\cten 10,70,20,90:{}	
			\cten 20,70,30,90:{}
			\cten 30,70,40,90:{}
			\cten 40,70,50,90:{}
			\cten 50,70,60,90:{}
			\cten 60,70,70,90:{}
			\cten 70,70,80,90:{}	
			\cten 80,70,90,90:{}	
			\ctnohead\ctdash
			\def\zzctdrawdashedge{\drawdashedge{3pt}{2.5pt}11}
			\cten 0,10,90,10:{} 	 	
			\cten 0,20,90,20:{} 	 	
			\cten 0,30,90,30:{} 	 
			\cten 0,50,90,50:{} 		
			\cten 0,40,90,40:{} 	 	
			\cten 0,60,90,60:{} 	 	
			\cten 0,80,90,80:{} 	 	
			\cten 10,0,10,90:{} 	 	
			\cten 20,0,20,90:{} 	 	
			\cten 30,0,30,90:{} 	 	
			\cten 40,0,40,90:{} 	 	
			\cten 50,0,50,90:{} 	 	
			\cten 60,0,60,90:{} 	 	
			\cten 70,0,70,90:{} 	 	
			\cten 80,0,80,90:{} 	 	
			\ctv -40,79:{\Thm^{(1)}[\matM]/\!{\meqvA}}
			\ctv 38,-12:{\FormsL^{(1)}/\!{\meqvA}}
		}
	\end{array}
	\]
	\caption{$\meqv{_\alg{A}}$-Equivalence partitioning}\label{fig-meqv}
\end{figure}

\begin{example} \label{ex-nontr}
	{\em	Let $\matM = \lr{\alg{A},\{\one\}}$, where $\alg{A} = \lr{\{\zero,\tb,\one\}; \impl, \zero}$ and $\impl$ is defined as in $\godelThree$ (Section~\ref{section:goedel}). We construct a complete set of representatives of $\FormsL^{(1)}/\meqvA$ as follows.
		\[
		\begin{array}{l|l|l}
			\text{depth} & \text{formulas} & \text{representatives}\\
			0 & p, \zero& p, \zero\\
			1 & p \impl p, p \impl \zero, \zero \impl p, \zero \impl \zero & p \impl p, p \impl \zero\\
			2 & (p \impl p) \impl p, p \impl (p \impl p), (p \impl p) \impl \zero, \zero \impl (p \impl p), &\\
			2 & (p \impl p) \impl (p \impl p),&\\
			2 & (p \impl p) \impl (p \impl \zero), (p \impl \zero) \impl (p \impl p),  &\\
			2 & (p \impl \zero) \impl p, p \impl (p \impl \zero), (p \impl \zero) \impl \zero, \zero \impl (p \impl \zero),& (p \impl \zero) \impl p\\
			2 & (p \impl \zero) \impl (p \impl p),(p \impl \zero) \impl (p \impl \zero),  &\\ 
			3 & ((p \impl \zero) \impl p) \impl p, p \impl ((p \impl \zero) \impl p), &((p \impl \zero) \impl p) \impl p\\
			3 & ((p \impl \zero) \impl p) \impl \zero, \zero \impl ((p \impl \zero) \impl p) &\\
			3 & ((p \impl \zero) \impl p) \impl (p \impl p), (p \impl p) \impl ((p \impl \zero) \impl p), &\\
			3 & ((p \impl \zero) \impl p) \impl (p \impl \zero), (p \impl \zero) \impl ((p \impl \zero) \impl p), &\\
			3 & ((p \impl \zero) \impl p) \impl ((p \impl \zero) \impl p)&
		\end{array}	
		\]
		We leave to the reader to verify that every formula of depth $4$ obtained from the selected representatives is $\meqvA$-equal to one of the selected representatives. (Exercise~\ref{sec-procetriv}.\ref{EX:ex-nontr})
		
		It is clear that abstract logic defined by $\matM$ is not trivial: formula $p \impl p$ is a theorem. }
\end{example}

\begin{rem}
	{\em In Example \ref{ex-nontr}, we could have stopped after selecting $p \impl p$ as a representative. We did not do this to illustrate how the enumeration procedure works.  }
\end{rem}

\paragraph{Exercises~\ref{sec-procetriv}}
\begin{enumerate}
	\item \label{EX:ex-nontr}Complete Example~\ref{ex-nontr} by proving that every formula of depth $4$ obtained from the selected representatives is $\meqvA$-equal to one of the selected representatives.
	\item Let $\mat{M}=\lr{\{\zero,\one\};\land,\lor,\{\one\}}$, where $\land$ and $\lor$ are defined as in the matrix $\booleTwo$ (Section~\ref{S:two-valued}). Prove that $\Thm[\matM]$ is empty.
\end{enumerate}
\section{Solution to Problem (m.2)} \label{sec-m2}

\begin{defn}\index{matrix!weak equivalence}\index{matrix!weakly adequate}
	Two $\Lan$-matrices $\matM_1$ and $\matM_2$ are said to be \textbf{weakly equivalent} if they are weakly adequate for the same abstract logic.
\end{defn}

One can rephrase the above definition: $\Lan$-matrices $\matM_1$ and $\matM_2$ are weakly equivalent if $\Thm[\matM_1] = \Thm[\matM_2]$. In terms of Lindenbaum matrices, two matrices $\matM_1$ and $\matM_2$ are weakly equivalent if and only if the Lindenbaum matrices of their respective abstract logic coincide.

We consider a bit more general problem: Is there an effective procedure that, given two finite $\Lan$-matrices $\matM_1$ and $\matM_2$, determines whether $\Thm[\matM_1] \subseteq \Thm[\matM_2]$? Clearly, if there is such a procedure, we also have a procedure that determines whether $\Thm[\matM_1] = \Thm[\matM_2]$.

To solve this problem we will do the following.
\begin{itemize}
	\item[(a)] We will reduce the problem to a case when the matrices $\matM_1$ and $\matM_2$ have the same carrier, more precisely, when $\matM_i = \lr{\alg{A}; D_i}$, where $i\in\{1,2\}$, and $D_2 \subseteq D_1$ and we need to determine whether $\Thm[\matM_1] = \Thm[\matM_2]$; cf. Section \ref{sec-Kal-a}.
	
	\item[(b)] We can reduce the problem even further to formulas on $n$ variables: by Corollary \ref{cor-redtn},
	if $n$ is cardinality of $\alg{A}$, then $\Thm[\matM_1] = \Thm[\matM_2]$ if, and only if, $\Thm^{(n)}[\matM_1] = \Thm^{(n)}[\matM_2]$.  
	
	\item[(c)] Furthermore, by Corollary \ref{cor-mateqvincl}, $D_2 \subseteq D_1$ entails $\Thm[\matM_1] \subseteq \Thm[\matM_2]$ and hence, 	$\Thmn[\matM_1]/\mmeqvA \ \subseteq \ \Thmn[\matM_2]/\mmeqvA$.
	Thus, we need only to determine whether $\Thmn[\matM_2]/\mmeqvA \ \subseteq \ \Thmn[\matM_1]/\mmeqvA$.
	
	By Proposition \ref{th-mefini}, $\FormsLn/\mmeqvA$ is a finite set. Thus, using the enumerating procedure from Section \ref{sec-procetriv}, for each selected representative $\alpha$ we can check whether $\alpha \in \Thmn[\matM_2]$ and if it is, we can check whether $\alpha \in \Thmn[\matM_1]$. It is clear that $\alpha \in \Thmn[\matM_i]$ if and only if $\alpha$ is valid in $\matM_i$ and because matrices $\matM_i$ are finite, we can effectively test whether $\alpha \in \Thmn[\matM_i]$ or not. During such a procedure, either we found a formula $\alpha$ such that $\alpha \in \Thmn[\matM_2]$ while $\alpha$ such that $\alpha \notin \Thmn[\matM_1]$, that is, $\Thmn[\matM_1] \neq \Thmn[\matM_2]$, or after we have obtain a complete set of representatives, we can conclude that $\Thmn[\matM_1] \neq \Thmn[\matM_2]$. 
\end{itemize}

\subsubsection{Reduction to the matrices with the same algebraic carrier}  \label{sec-Kal-a}

Let us consider the following operations on matrices.

\begin{defn}\label{def-oper}\index{$\matM_1 \lsum \matM_2$}\index{$\matM_1 \rsum \matM_2$}\index{$\matM_1 \times \matM_2$}\index{$\matM_1 + \matM_2$}
	Suppose that {\em$\matM_i = \lr{\alg{A}_i,D_i}$}, where $i=1,2$, are $\Lan$-matrices. Then, let
	{\em\[
		\begin{aligned}
			&(\text{a}) \quad \matM_1 \lsum \matM_2 = \lr{\alg{A}_1 \times \alg{A}_2, D_1 \times |\alg{A}_2|}; \\
			&(\text{b}) \quad \matM_1 \rsum \matM_2 = \lr{\alg{A}_1 \times \alg{A}_2, |\alg{A}_1| \times D_2};\\
			&(\text{c}) \quad \matM_1 \times \matM_2 = \lr{\alg{A}_1 \times \alg{A}_2, D_1 \times D_2};\\
			&(\text{d}) \quad \matM_1 + \matM_2 = \lr{\alg{A}_1 \times \alg{A}_2, |D_1| \times |\alg{A}| \cup |\alg{A}_1| \times D_2}.
		\end{aligned}
		\]}
\end{defn}

As we can see, all matrices on the right hand side have the same carrier.

\begin{prop}\label{pr-sum} Let {\em$\matM_1 =  \lr{\alg{A}_1,D_1}$} and {\em$\matM_2 =  \lr{\alg{A}_2,D_2}$} be $\Lan$-matrices. Then 
	{\em\[
		\begin{aligned}
			&(\text{a}) \quad  \aLog[\matM_1 \lsum \matM_2] = \aLog[\matM_1]; \\
			&(\text{b}) \quad  \aLog[\matM_1 \rsum \matM_2] = \aLog[\matM_2];\\
			&(\text{c}) \quad  \Thm[\matM_1 \times \matM_2] = \Thm[\matM_1] \cap \Thm[\matM_2];\\
			&(\text{d}) \quad  \Thm[\matM_1 + \matM_2] = \Thm[\matM_1] \cup \Thm[\matM_2].
		\end{aligned}
		\]}
\end{prop}
\noindent\emph{Proof}~is left to the reader. (Exercise~\ref{sec-m2}.\ref{EX:pr-sum})\\

\begin{cor}\label{cor-redone}  Suppose that {\em$\matM_1 =  \lr{\alg{A}_1,D_1}$} and {\em$\matM_2 =  \lr{\alg{A}_2,D_2}$}. Then 
	{\em\begin{itemize}
			\item[(a)] $\aLog[\matM_1] \preccurlyeq \aLog[\matM_2] \iff \aLog[\matM_1 \lsum \matM_2] \preccurlyeq \aLog[\matM_1 \rsum \matM_2]$;
			\item[(b)] $\Thm[\matM_1] \subseteq \Thm[\matM_2] \iff \Thm[\matM_1 \times \matM_2] = \Thm[\matM_1 \lsum \matM_2]$;
			\item[(c)] $\Thm[\matM_1] = \Thm[\matM_2] \iff \Thm[\matM_1 \times \matM_2] = \Thm[\matM_1 + \matM_2]$.	
	\end{itemize}}
\end{cor} 
\begin{proof}
	We note that if matrices $\matM_1$ and $\matM_2$ are finite, there is an effective procedure for constructing matrices $\matM_1 \lsum \matM_2$,  $\matM_1 \rsum \matM_2$, $\matM_1 \times \matM_2$ and $\matM_1 + \matM_2$. Moreover, Corollary \ref{cor-redone}(b) entails that the problems to determine by given finite matrices $\matM_1$ and $\matM_2$ whether (i) $\Thm[\matM_1] \subseteq \Thm[\matM_2]$ and whether (ii) $\Thm[\matM_1] = \Thm[\matM_2]$ are equivalent.
\end{proof}

\paragraph{Exercises~\ref{sec-m2}}
\begin{enumerate}
	\item \label{EX:pr-sum}Prove Proposition~\ref{pr-sum}.
	\item Let $\mathbf{M}_1 = \langle(\{\zero,a,\one\},\land,\lor,\to,\neg),\{\one\} \rangle$ and $\mathbf{M}_2 = \langle(\{\zero,a,\one\},\land,\lor,\to,-),\{\one\} \rangle$, where $\land,\lor,\to$ and $\neg$ are defined as in Heyting algebra, and $-$
	is~{\L}ukasiewicz's involutive negation: $-\zero = \one, -a = a, -\one = \zero$. For $\to$ is a Heyting implication, the sets $\Thm[\mathbf{M}_1]$ and $\Thm[\mathbf{M}_2]$ are closed under modus ponens. Prove that $\Thm[\mathbf{M}_1 + \mathbf{M}_2]$ is not closed under modus ponens.
	(Hint: Observe that  $\neg\neg (p \lor \neg p) \in\Thm[\mathbf{M}_1]$ and $\neg\neg (p \lor \neg p) \to (p \lor \neg p) \in\Thm[\mathbf{M}_2]$, while $(p \lor \neg p) \notin\Thm[\mathbf{M}_1] \cup\Thm[\mathbf{M}_2] = \Thm[\mathbf{M}_1 + \mathbf{M}_2]$.)
\end{enumerate}

\section{Solution to Problem (m.3)}\label{section:solution-of-Problem-m3}

Given two finite $\Lan$-atlases $\M_1$ and $\M_2$, we need to determine whether $\aLog[\M_1]$ and $\aLog[\M_2]$, where $\aLog[\M_1]$ and $\aLog[\M_2]$ are  abstract logics defined by these atlases.


To solve this problem we take a path similar to the path we took to solve the problem of weak equivalence of matrices, and we will do the following.

\begin{itemize}
	\item[(a)] We can reduce the problem to a case when the atlases $\M_1$ and $\M_2$ have the same carrier, let say $\alg{A}$ (cf. Proposition \ref{pr-atlsame}).
	
	\item[(b)] We can reduce the problem further: if $\M_1$ and $\M_2$ are $\Lan$-atlases of cardinalities $n_1$ and $n_2$,
	by Corollary \ref{cor-logn}, for $ m = \max(n_1,n_2)$,
	\[
	\aLog[\M_1] = \aLog[\Lin^{(m)}[\M_1]] \text{ and } \aLog[\M_2] = \aLog[\Lin^{(m)}[\M_2]] 
	\]
	and thus,
	\[
	\aLog[\M_1] \preccurlyeq \aLog[\M_2] \iff  \aLog[\Lin^{(m)}[\M_1]] \preccurlyeq \aLog[\Lin^{(m)}[\M_2]]. 
	\] 
	
	\item[(c)] Similarly to Step (c) from Section \ref{sec-m2}, we can obtain a complete representatives of $\Lin^{(m)}/\mmeqvA$, let say $X$. Because by Proposition \ref{th-mefini} the set $X$ is finite, we can check whether for every $\alpha \cup Y \subseteq X$,
	\[
	Y \models_{\Lin^{(m)}[\M_1]}\alpha  \quad  \text{ yields }  \quad  Y \models_{\Lin^{(m)}[\M_2]} \alpha 
	\]   
	and in such a way we can  to determine whether $S[\Lin^{(m)}[\M_1]] \preccurlyeq S[\Lin^{(m)}[\M_2]]$, or not.
\end{itemize}

\subsubsection{Reduction to the atlases with the same algebraic carrier}  \label{sec-Kal-atl}

Suppose that $\M_i = \lr{\alg{A}_i,\D_i}, i=1,2$ are $\Lan$-atlases. Then, let

\begin{equation}\label{E:10.3}
	\begin{array}{c}
		\M_1 \lsum \M_2 := \lr{\alg{A}_1 \times \alg{A}_2, \set{D \times |\alg{A}_2|}{D \in \D_1}}, \\
		\M_1 \rsum \M_2 := \lr{\alg{A}_1 \times \alg{A}_2, \set{|\alg{A}_1| \times D}{D \in \D_2}}.\\
	\end{array}
\end{equation}

As we can see,  atlases $\M_1$ and $\M_2$ have the same carrier.

\begin{prop}\label{pr-atlsame}
	Let {\em$\M_1 = \lr{\alg{A}_1,\D_1}$} and {\em$\M_2 = \lr{\alg{A}_2,\D_2}$} be $\Lan$-atlases. Then
	$	\aLog[\M_1] = \aLog[\M_1 \lsum \M_2]$ and $\aLog[\M_2] = \aLog[\M_1 \rsum \M_2]$.
\end{prop}
\noindent\emph{Proof}~is left to the reader. (Exercise~\ref{section:solution-of-Problem-m3}.\ref{EX:pr-atlsame})\\

Thus, there is an effective procedure to determine whether two finite $\Lan$-atlases are equivalent. Because there is a procedure that  

\begin{rem}
	{\em	Let us observe that for any $\Lan$-atlases $\M_1$ and $\M_2$, by Proposition \ref{pr-atlsame}, $S[\M_1 \lsum \M_2] = S[\M_1]$ and because $S[\M_1 \times \M_2] = S[\M_1] \cap S[\M_2]$,
		\[
		S[\M_1] \preccurlyeq S[\M_2] \iff S[\M_1 \lsum \M_2] = S[\M_1 \times \M_2]. 
		\]}
\end{rem}

Thus, the problems to determine by two finite $\Lan$-atlases $\M_1$ and $\M_2$ whether $\aLog{\M_1}\preccurlyeq \aLog[\M_2]$ and whether $\aLog[\M_1] = \aLog[\M_2]$ are equivalent.

By Proposition 4.4.10, for any finite set $\matM_i$, where $i\in \{1,\ldots, m\}$, of finite $\Lan$-matrices, there is a finite $\Lan$-atlas $\M$ defining the same abstract logic. The proof of Proposition 4.4.10 gives an effective procedure how to construct $\M$. Therefore, there is an effective procedure which, given two finite sets of finite $\Lan$-matrices, determines whether these sets define the same abstract logic.

We note that Corollary \ref{cor-logn} holds for the abstract logics understood as multiple conclusion consequence relations. Hence, the procedure that solves the problem of equivalence of two finite sets of finite matrices can be easily modified fo the multiple conclusion consequence relations.

\paragraph{Exercises~\ref{section:solution-of-Problem-m3}}
\begin{enumerate}
	\item \label{EX:pr-atlsame}Prove Proposition~\ref{pr-atlsame}.
	\item Derive the equalities $\aLog[\M \lsum \M]= \aLog[\M \rsum \M]=\aLog[\M]$ from the definition~\eqref{E:10.3}, without using Proposition~\ref{pr-atlsame}.
\end{enumerate}

\section{Historical notes}

In \cite{los1949} {\L}o{\'s} observed that if $S$ is an abstract logic defined by a logical matrix $\matM = \lr{\alg{A},D}$, for any $n > 0$, one can consider a functional matrix: $\matM^{(n)} := \lr{\alg{A}^{(n)},D^{(n)}}$, where $\alg{A}^{(n)}$ is a set of all functions $f:|\alg{A}|^n \longrightarrow |\alg{A}|$, $D' \subseteq \alg{A}^{(n)}$ consisting of all functions from $\alg{A}^{(n)}$ taking values only from $D$, and connectives are defined pointwise. Moreover, one can consider a submatrix $\matM_e^{(n)}$ of $\matM^{(n)}$ consisting of all functions expressible via superposition of fundamental operations of $\alg{A}$. Such submatrices are closely related to restricted Lindenbaum matrices, and their interrelations were studied in \cite[\S 15]{los1949}. In particular, if $\matM$ is a finite matrix of cardinality $m$, then $\Thm[\matM] = \Thm[\matM_e^{(n)}]$ for any $n \ge m$.

An effective procedure for enumerating logical functions on $n$ variables defined on a finite matrix was suggested in \cite{wernick1939}. Let us observe that given a finite matrix $\matM$, one can effectively construct matrix $\matM_e^{(n)}$ using approach similar to enumerating procedure described in Section \ref{sec-restr}: recall that two formulas $\alpha$ and $\beta$ are $\meqvA$-equivalent if and only if $\alpha$ and $\beta$ define the same function.

Proposition 25 from \cite{los1949} states that abstract logic $\aLog$ of a finite matrix $\matM$ has no theorem if and only if the set of designated elements of $\matM_e^{(1)}$ is empty and hence, there is an effective procedure to determine whether abstract logic defined by a finite matrix has theorems or not. 

Two matrices $\matM_i = \lr{\alg{A}_i,D_i}, i=1,2$ are \textbf{matrix isomorphic} if there is an isomorphism $\varphi$ between algebras $\alg{A}_i$ such that $\varphi(D_1) = D_2$. Proposition 27* states that for any two finite matrices $\matM'$ and $\matM''$ of cardinalities $n_1$ and $n_2$,  $\Thm[\matM'] = \Thm[\matM'']$ if and only if matrices $\matM_e'^{(n)}$ and $\matM_e''^{(n)}$, where $n = \max(n_1,n_2)$, are matrix isomorphic. Thus, given two finite matrices $\matM'$ and $\matM''$, one can effectively construct matrices $\matM_e'^{(n)}$ and $\matM_e''^{(n)}$ and check whether these matrices are matrix isomorphic and in such a way, to determine whether $\Thm[\matM'] = \Thm[\matM'']$.

An effective procedure to determine whether an abstract logic defined by a given finite matrix has theorems was suggested in \cite{kalicki1950}. In \cite{kalicki1952}, Kalicki developed further his approach from his 1950 paper and constructed an effective procedure to determine whether two finite logical matrices define the abstract logics with the same sets of theorems, and he applied the same ideas to finite algebras: to determine whether two finite algebras having the same finite similarity type, have the same sets of valid equalities (in other words, whether two given finite algebras generate the same variety).   

In \cite{kalicki1950_1} Kalicki introduced and studied the operations on matrices: the direct product and sum of matrices (Definition \ref{def-oper} (c) and (d)). 

Existence of the greatest matrix congruence was observed in \cite{porte1965}. We used a proof suggested in \cite{brown-suszko1973}.

In \cite{citkin1975}, Citkin, using an approach similar to \cite{los1949}, described an effective procedure which given two finite $\Lan$-atlases determined whether thee atlases define the same abstract logic.\footnote{More precisely, the procedure in \cite{citkin1975} is for Smiley matrices, which can be viewed as atlases.}

Later, in \cite{zygmunt1983} suggested an effective procedure that given two finite $\Lan$-atlases  $\M_1$ and $\M_2$ (ramified matrices in Zygmunt's terms), determines whether $\aLog[\M_1] \preccurlyeq \aLog[\M_2]$.

In \cite{devyatkin2013}, the Kalicki's approach from \cite{kalicki1950} was used to construct an effective procedure which by two finite $\Lan$-matrices determines whether these matrices define the same abstract logic.

Even though there are effective procedures to determine equivalence of two finite $\Lan$-atlases, these procedures have a high complexity: it follows from \cite{bergman-slutski2000}  that even problem to determine by two finite algebras whether they generate the same quasivariety, it is NP-complete. And for algebraziable logics, the latter problem is equivalent to the problem of determining equivalence of two finite matrices.

\cleardoublepage
\phantomsection
\addcontentsline{toc}{chapter}{Bibliography}

\bibliographystyle{apalike}
\bibliography{Book}

\printindex

\end{document}